\pdfoutput=1
\documentclass[12pt]{amsart}
\usepackage{lmodern}
\usepackage{ifthen}
\usepackage{amsfonts}
\usepackage{amsxtra}
\usepackage{amssymb}
\usepackage{mathdots}
\usepackage{array}
\usepackage[margin=0.8in]{geometry}
\usepackage{xcolor}
\definecolor{cite}{rgb}{0.30,0.60,1.00}
\definecolor{url}{rgb}{0.00,0.00,0.80}
\definecolor{link}{rgb}{0.40,0.10,0.20}
\usepackage[colorlinks,linkcolor=link,urlcolor=url,citecolor=cite,pagebackref,breaklinks]{hyperref}
\usepackage{bbm}
\usepackage{mathtools}
\usepackage{mathrsfs}
\usepackage{appendix}
\usepackage[all]{xy}
\usepackage[lite,abbrev,msc-links,alphabetic]{amsrefs}

\usepackage{graphicx}
\usepackage{multirow}
\usepackage{pstricks}
\usepackage{pst-pdf}
\usepackage{enumitem}
\usepackage{pifont}
\usepackage{marvosym}

\usepackage[OT2,T1]{fontenc}
\DeclareSymbolFont{cyrletters}{OT2}{wncyr}{m}{n}
\DeclareMathSymbol{\Sha}{\mathalpha}{cyrletters}{"58}

\usepackage{graphicx}

\makeatletter
\providecommand*{\Dashv}{%
  \mathrel{%
    \mathpalette\@Dashv\vDash
  }%
}
\newcommand*{\@Dashv}[2]{%
  \reflectbox{$\m@th#1#2$}%
}
\makeatother

\numberwithin{equation}{section}

\theoremstyle{plain}
\newtheorem{proposition}{Proposition}[subsection]
\newtheorem{conjecture}[proposition]{Conjecture}
\newtheorem{corollary}[proposition]{Corollary}
\newtheorem{lem}[proposition]{Lemma}
\newtheorem{theorem}[proposition]{Theorem}

\theoremstyle{definition}
\newtheorem{definition}[proposition]{Definition}
\newtheorem{construction}[proposition]{Construction}
\newtheorem{notation}[proposition]{Notation}
\newtheorem{assumption}[proposition]{Assumption}
\newtheorem{hypothesis}[proposition]{Hypothesis}

\theoremstyle{remark}
\newtheorem{remark}[proposition]{Remark}

\renewcommand{\b}[1]{\mathbf{#1}}
\renewcommand{\c}[1]{\mathcal{#1}}
\renewcommand{\d}[1]{\mathbb{#1}}
\newcommand{\f}[1]{\mathfrak{#1}}
\renewcommand{\r}[1]{\mathrm{#1}}
\newcommand{\s}[1]{\mathscr{#1}}
\renewcommand{\sf}[1]{\mathsf{#1}}
\renewcommand{\(}{\left(}
\renewcommand{\)}{\right)}
\newcommand{\res}{\mathbin{|}}
\newcommand{\ol}[1]{\overline{#1}{}}
\newcommand{\wt}[1]{\widetilde{#1}{}}
\newcommand{\ul}{\underline}
\renewcommand{\leq}{\leqslant}
\renewcommand{\geq}{\geqslant}

\newcommand{\bG}{\b G}

\newcommand{\bM}{\b M}

\newcommand{\bP}{\b P}
\newcommand{\bQ}{\b Q}

\newcommand{\bT}{\b T}

\newcommand{\bX}{\b X}

\newcommand{\bi}{\b i}

\newcommand{\bm}{\b m}

\newcommand{\cA}{\c A}

\newcommand{\cD}{\c D}
\newcommand{\cE}{\c E}
\newcommand{\cF}{\c F}
\newcommand{\cG}{\c G}
\newcommand{\cH}{\c H}
\newcommand{\cI}{\c I}

\newcommand{\cL}{\c L}

\newcommand{\cN}{\c N}
\newcommand{\cO}{\c O}
\newcommand{\cP}{\c P}

\newcommand{\cS}{\c S}
\newcommand{\cT}{\c T}

\newcommand{\cV}{\c V}

\newcommand{\dA}{\d A}
\newcommand{\dB}{\d B}
\newcommand{\dC}{\d C}

\newcommand{\dE}{\d E}
\newcommand{\dF}{\d F}
\newcommand{\dG}{\d G}

\newcommand{\dN}{\d N}

\newcommand{\dP}{\d P}
\newcommand{\dQ}{\d Q}
\newcommand{\dR}{\d R}

\newcommand{\dT}{\d T}

\newcommand{\dX}{\d X}

\newcommand{\dZ}{\d Z}

\newcommand{\fC}{\f C}
\newcommand{\fD}{\f D}

\newcommand{\fK}{\f K}

\newcommand{\fP}{\f P}

\newcommand{\fS}{\f S}
\newcommand{\fT}{\f T}

\newcommand{\fa}{\f a}

\newcommand{\fc}{\f c}
\newcommand{\fd}{\f d}

\newcommand{\ff}{\f f}

\newcommand{\fm}{\f m}
\newcommand{\fn}{\f n}

\newcommand{\fp}{\f p}
\newcommand{\fq}{\f q}
\newcommand{\fr}{\f r}
\newcommand{\fs}{\f s}

\newcommand{\rA}{\r A}
\newcommand{\rB}{\r B}

\newcommand{\rE}{\r E}
\newcommand{\rF}{\r F}
\newcommand{\rG}{\r G}
\newcommand{\rH}{\r H}
\newcommand{\rI}{\r I}
\newcommand{\rJ}{\r J}
\newcommand{\rK}{\r K}
\newcommand{\rL}{\r L}
\newcommand{\rM}{\r M}
\newcommand{\rN}{\r N}

\newcommand{\rP}{\r P}
\newcommand{\rQ}{\r Q}
\newcommand{\rR}{\r R}
\newcommand{\rS}{\r S}
\newcommand{\rT}{\r T}
\newcommand{\rU}{\r U}
\newcommand{\rV}{\r V}
\newcommand{\rW}{\r W}
\newcommand{\rX}{\r X}
\newcommand{\rY}{\r Y}
\newcommand{\rZ}{\r Z}

\newcommand{\rb}{\r b}
\newcommand{\rc}{\r c}
\newcommand{\rd}{\r d}

\newcommand{\rh}{\r h}

\renewcommand{\rm}{\r m}

\newcommand{\rs}{\r s}
\newcommand{\rt}{\r t}

\newcommand{\sD}{\s D}

\newcommand{\sG}{\s G}

\newcommand{\sI}{\s I}

\newcommand{\sP}{\s P}

\newcommand{\sS}{\s S}

\newcommand{\sV}{\s V}
\newcommand{\sW}{\s W}

\newcommand{\sfH}{\sf H}

\newcommand{\sfP}{\sf P}

\newcommand{\sfR}{\sf R}

\newcommand{\sfT}{\sf T}

\newcommand{\sfs}{\sf s}

\newcommand{\sfv}{\sf v}

\newcommand{\sfx}{\sf x}

\newcommand{\tF}{\mathtt{F}}

\newcommand{\tI}{\mathtt{I}}

\newcommand{\tR}{\mathtt{R}}

\newcommand{\tT}{\mathtt{T}}

\newcommand{\tV}{\mathtt{V}}

\newcommand{\tc}{\mathtt{c}}
\newcommand{\td}{\mathtt{d}}

\newcommand{\ti}{\mathtt{i}}
\newcommand{\tj}{\mathtt{j}}

\newcommand{\tw}{\mathtt{w}}

\newcommand{\balpha}{\boldsymbol{\alpha}}

\newcommand{\bmu}{\boldsymbol{\mu}}
\newcommand{\bbt}{\boldsymbol{t}}

\newcommand{\obj}{\text{\Flatsteel}}
\newcommand{\qbinom}[2]{\genfrac{[}{]}{0pt}{}{#1}{#2}}

\newcommand{\tp}[1]{\prescript{\rt\!}{}{#1}}
\newcommand{\pres}[2]{\prescript{#1}{}{#2}}
\newcommand{\floor}[1]{\lfloor{#1}\rfloor}
\newcommand{\ceil}[1]{\lceil{#1}\rceil}
\newcommand{\ab}{\r{ab}}

\newcommand{\bad}{\r{bad}}

\newcommand{\CF}{\mathbbm{1}}
\newcommand{\cInd}{\rc\text{-}\r{Ind}}
\newcommand{\cl}{\r{cl}}

\newcommand{\cris}{\r{cris}}

\newcommand{\dr}{\r{dR}}
\newcommand{\et}{{\acute{\r{e}}\r{t}}}

\newcommand{\ext}{\r{ext}}

\newcommand{\FL}{\r{FL}}

\newcommand{\free}{\r{fr}}
\newcommand{\gr}{\r{gr}}
\renewcommand{\graph}{\triangle}

\newcommand{\Hk}{\r{Hk}}
\newcommand{\HOM}{\c{H}om}

\newcommand{\id}{\r{id}}
\newcommand{\inc}{\r{inc}}
\newcommand{\Inc}{\r{Inc}}

\newcommand{\loc}{\r{loc}}
\newcommand{\lr}{\r{lr}}

\newcommand{\mix}{\r{mix}}
\newcommand{\mnm}{\r{min}}

\newcommand{\ns}{\r{ns}}
\newcommand{\op}{\r{op}}
\newcommand{\prim}{\r{prim}}
\newcommand{\ram}{\r{ram}}
\newcommand{\rat}{\r{rat}}

\newcommand{\Sat}{\r{Sat}}
\newcommand{\sh}{\r{sh}}
\newcommand{\sing}{\r{sing}}
\renewcommand{\sp}{\r{sp}}

\newcommand{\std}{\r{std}}
\newcommand{\St}{\r{St}}
\newcommand{\sym}{\r{sym}}
\newcommand{\tor}{\r{tor}}
\newcommand{\univ}{\r{univ}}
\newcommand{\unr}{\r{unr}}
\newcommand{\ur}{\r{ur}}

\newcommand{\hra}{\hookrightarrow}

\newcommand{\Mod}{\sf{Mod}}
\newcommand{\Set}{\sf{Set}}
\newcommand{\Fun}{\sf{Fun}}
\newcommand{\Sch}{\sf{Sch}}

\DeclareMathOperator{\ad}{ad}
\DeclareMathOperator{\AJ}{AJ}

\DeclareMathOperator{\As}{As}
\DeclareMathOperator{\Aut}{Aut}

\DeclareMathOperator{\BC}{BC}
\DeclareMathOperator{\CH}{CH}

\DeclareMathOperator{\coker}{coker}

\DeclareMathOperator{\Def}{Def}
\DeclareMathOperator{\diag}{diag}
\DeclareMathOperator{\disc}{disc}
\DeclareMathOperator{\Disc}{Disc}

\DeclareMathOperator{\DL}{DL}
\DeclareMathOperator{\End}{End}

\DeclareMathOperator{\Fr}{Fr}

\DeclareMathOperator{\Gal}{Gal}

\DeclareMathOperator{\GL}{GL}

\DeclareMathOperator{\GU}{GU}
\DeclareMathOperator{\Hom}{Hom}
\DeclareMathOperator{\IP}{Im}
\DeclareMathOperator{\IM}{im}

\DeclareMathOperator{\Ind}{Ind}

\DeclareMathOperator{\Iso}{Iso}

\DeclareMathOperator{\Ker}{ker}
\DeclareMathOperator{\Lat}{Lat}
\DeclareMathOperator{\Lie}{Lie}
\DeclareMathOperator{\Map}{Map}

\DeclareMathOperator{\modulo}{mod}
\DeclareMathOperator{\Nm}{Nm}

\DeclareMathOperator{\ord}{ord}

\DeclareMathOperator{\rank}{rank}

\DeclareMathOperator{\Res}{Res}

\DeclareMathOperator{\Sh}{Sh}

\DeclareMathOperator{\Sp}{Sp}

\DeclareMathOperator{\Spec}{Spec}

\DeclareMathOperator{\Sym}{Sym}

\DeclareMathOperator{\Tr}{Tr}
\DeclareMathOperator{\tr}{tr}

\DeclareMathOperator{\val}{val}
\DeclareMathOperator{\Ver}{Ver}

\DeclareMathOperator{\WD}{WD}

\begin{document}

\title{On the Beilinson--Bloch--Kato conjecture for Rankin--Selberg motives}

\author{Yifeng Liu}
\address{Institute for Advanced Study in Mathematics, Zhejiang University, Hangzhou 310058, China}
\email{liuyf0719@zju.edu.cn}

\author{Yichao Tian}
\address{Morningside Center of Mathematics, AMSS, Chinese Academy of Sciences, Beijing 100190, China}
\email{yichaot@math.ac.cn}

\author{Liang Xiao}
\address{Beijing International Center for Mathematical Research, Peking University, Beijing 100871, China}
\email{lxiao@bicmr.pku.edu.cn}

\author{Wei Zhang}
\address{Department of Mathematics, Massachusetts Institute of Technology, Cambridge MA 02139, United States}
\email{weizhang@mit.edu}

\author{Xinwen Zhu}
\address{Division of Physics, Mathematics and Astronomy, California Institute of Technology, Pasadena CA 91125, United States}
\email{xzhu@caltech.edu}

\date{\today}

\begin{abstract}
  In this article, we study the Beilinson--Bloch--Kato conjecture for motives associated to Rankin--Selberg products of conjugate self-dual automorphic representations, within the framework of the Gan--Gross--Prasad conjecture. We show that if the central critical value of the Rankin--Selberg $L$-function does not vanish, then the Bloch--Kato Selmer group with coefficients in a favorable field of the corresponding motive vanishes. We also show that if the class in the Bloch--Kato Selmer group constructed from a certain diagonal cycle does not vanish, which is conjecturally equivalent to the nonvanishing of the central critical first derivative of the Rankin--Selberg $L$-function, then the Bloch--Kato Selmer group is of rank one.
\end{abstract}

\maketitle

\tableofcontents

\section{Introduction}
\label{ss:1}

In this article, we study the Beilinson--Bloch--Kato conjecture for motives associated to Rankin--Selberg products of conjugate self-dual automorphic representations of $\GL_n(\dA_F)\times\GL_{n+1}(\dA_F)$ for a CM number field $F$, within the framework of the Gan--Gross--Prasad conjecture \cite{GGP12} for the pair of unitary groups $\rU(n)\times\rU(n+1)$. For background on the Beilinson--Bloch--Kato conjecture, which is a generalization of the Birch and Swinnerton-Dyer conjecture from elliptic curves to higher dimensional algebraic varieties, we refer to \cite{BK90} (see also the introduction of \cite{Liu1}).

\subsection{Main results}
\label{ss:main}

Let $F/F^+$ be a totally imaginary quadratic extension of a totally real number field. We first state one of our main results that is least technical to understand.

\begin{theorem}[Corollary~\ref{co:sym}]\label{th:sym}
Let $n\geq 2$ be an integer. Let $A$ and $A'$ be two modular elliptic curves over $F^+$ such that $\End(A_{\ol{F}})=\End(A'_{\ol{F}})=\dZ$. Suppose that
\begin{enumerate}[label=(\alph*)]
  \item $A_{\ol{F}}$ and $A'_{\ol{F}}$ are not isogenous to each other;

  \item both $\Sym^{n-1}A$ and $\Sym^nA'$ are modular; and

  \item $F^+\neq\dQ$.
\end{enumerate}
If the (central critical) $L$-value $L(n,\Sym^{n-1}A_F\times\Sym^nA'_F)$ does not vanish, then the Bloch--Kato Selmer group
\[
\rH^1_f(F,\Sym^{n-1}\rH^1_\et(A_{\ol{F}},\dQ_\ell)\otimes_{\dQ_\ell}\Sym^n\rH^1_\et(A'_{\ol{F}},\dQ_\ell)(n))
\]
vanishes for all but finitely many rational primes $\ell$.
\end{theorem}

\begin{remark}
The finite set of rational primes $\ell$ that are excluded in Theorem \ref{th:sym} can be effectively bounded. We now explain the three conditions in Theorem \ref{th:sym}.
\begin{enumerate}[label=(\alph*)]
  \item is necessary since otherwise (L3) and (L5) in Definition \ref{de:lambda} fail for all rational primes $\ell$.

  \item is necessary since our approach only applies to Galois representations arising from automorphic representations. We summarise the current knowledge on the modularity of symmetric powers of elliptic curves in Remark \ref{re:sym}.

  \item is necessary only for technical reasons. First, we do not know Hypothesis \ref{hy:unitary_cohomology}, which concerns cohomology of unitary Shimura varieties, yet for $N\geq 4$ if $F^+=\dQ$. Second, we do not have (an appropriate replacement for) Proposition \ref{th:generic}, a result generalizing \cite{CS17}, when $F^+=\dQ$. Indeed, as long as we have these results as expected, (c) can be lifted.
\end{enumerate}
\end{remark}

Theorem \ref{th:sym} is a special case of a more general result concerning the Bloch--Kato Selmer groups of Galois representations associated to conjugate self-dual automorphic representations. To reduce the burden of long and technical terminology in the future, we first introduce the following definition, which will serve for the entire article.

\begin{definition}\label{de:relevant}
We say that a complex representation $\Pi$ of $\GL_N(\dA_F)$ with $N\geq 1$ is \emph{relevant} if
\begin{enumerate}
  \item $\Pi$ is an irreducible cuspidal automorphic representation;

  \item $\Pi\circ\tc\simeq\Pi^\vee$, where $\tc\in\Gal(F/F^+)$ is the complex conjugation;

  \item for every archimedean place $w$ of $F$, $\Pi_w$ is isomorphic to the (irreducible) principal series representation induced by the characters $(\arg^{1-N},\arg^{3-N},\dots,\arg^{N-3},\arg^{N-1})$, where $\arg\colon\dC^\times\to\dC^\times$ is the \emph{argument character} defined by the formula $\arg(z)\coloneqq z/\sqrt{z\ol{z}}$.
\end{enumerate}
\end{definition}

\begin{remark}\label{re:relevant}
If $\Pi$ is relevant, then it is regular algebraic in the sense of \cite{Clo90}*{Definition~3.12}. Moreover, it is well-known that $L(s,\Pi,\As^{(-1)^N})$ is regular at $s=1$ (see, for example, \cite{GHL}*{\S6.1}).
\end{remark}

Now we can state our main result in the context of automorphic representations, of which Theorem \ref{th:sym} is a special case. Till the end of the next subsection, we will take an integer $n\geq 2$, and denote by $n_0$ and $n_1$ the unique even and odd numbers in $\{n,n+1\}$, respectively.

\begin{theorem}[Theorem \ref{th:selmer0}]\label{th:selmer0_pre}
Let $\Pi_0$ and $\Pi_1$ be relevant representations of $\GL_{n_0}(\dA_F)$ and $\GL_{n_1}(\dA_F)$, respectively. Let $E\subseteq\dC$ be a strong coefficient field of both $\Pi_0$ and $\Pi_1$ (Definition \ref{de:weak_field}). Suppose that $F^+\neq\dQ$. If $L(\frac{1}{2},\Pi_0\times\Pi_1)\neq 0$, then for all admissible primes $\lambda$ of $E$ with respect to $(\Pi_0,\Pi_1)$, the Bloch--Kato Selmer group $\rH^1_f(F,\rho_{\Pi_0,\lambda}\otimes_{E_\lambda}\rho_{\Pi_1,\lambda}(n))$ vanishes. Here, $\rho_{\Pi_\alpha,\lambda}$ is the Galois representation of $F$ with coefficients in $E_\lambda$ associated to $\Pi_\alpha$ for $\alpha=0,1$, as described in Proposition \ref{pr:galois} and Definition \ref{de:weak_field}.
\end{theorem}

\begin{remark}
The notion of admissible primes appearing in Theorem \ref{th:selmer0_pre} is introduced in Definition \ref{de:lambda}, which consists of a long list of assumptions, some of which are rather technical. Here, we would like to comment on the essence of these assumptions.
\begin{enumerate}
  \item[(L1,2)] are elementary and exclude only finitely many primes $\lambda$.

  \item[(L3)] is expected to hold for every prime $\lambda$ if and only if the (conjectural) automorphic product $\Pi_0\boxtimes\Pi_1$, as an irreducible admissible representation of $\GL_{n(n+1)}(\dA_F)$, remains cuspidal.

  \item[(L4)] is expected to hold for all but finitely many primes $\lambda$.

  \item[(L5)] is basically saying that, under (L4), the image of the pair of residual Galois representations $(\bar\rho_{\Pi_0,\lambda},\bar\rho_{\Pi_1,\lambda})$ contains an element of a particular form. It is expected to hold for all but finitely many primes $\lambda$ if the two automorphic representations $\Pi_0$ and $\Pi_1$ are not correlated in some manner. For example, when $n=2$, we expect that as long as $\Pi_1$ is not an automorphic twist of $\Sym^2\Pi_0$ after any base change, then (L5) holds for all but finitely many primes $\lambda$.

  \item[(L6)] is a technical assumption that is only used in the argument of an R=T theorem concerning Galois deformations in \cite{LTXZZ}. It is expected to hold for all but finitely many primes $\lambda$ (see \cite{LTXZZ}*{\S4.2}).

  \item[(L7)] is a technical assumption for the vanishing of certain Hecke localized cohomology of unitary Shimura varieties off middle degree. In fact, when $F^+\neq\dQ$, (L7) holds for all but finitely many primes $\lambda$ by Corollary \ref{co:generic}.
\end{enumerate}
\end{remark}

In fact, we have dedicated ourselves to obtaining the following family of abstract examples in which all but finitely many primes are admissible. Note that neither the following theorem nor Theorem \ref{th:sym} implies each other.

\begin{theorem}[Corollary \ref{co:abstract}]\label{th:abstract_pre}
Let $\Pi_0$, $\Pi_1$, and $E$ be as in Theorem \ref{th:selmer0_pre}. Suppose that
\begin{enumerate}[label=(\alph*)]
  \item there exists a very special inert prime $\fp$ of $F^+$ (Definition \ref{de:special_inert}) such that $\Pi_{0,\fp}$ is Steinberg, and $\Pi_{1,\fp}$ is unramified whose Satake parameter contains $1$ exactly once;\footnote{Note that the Satake parameter of $\Pi_{1,\fp}$ has to contain $1$ at least once by Definition \ref{de:relevant}(2).}

  \item for $\alpha=0,1$, there exists a nonarchimedean place $w_\alpha$ of $F$ such that $\Pi_{\alpha,w_\alpha}$ is supercuspidal; and

  \item $F^+\neq\dQ$.
\end{enumerate}
If $L(\frac{1}{2},\Pi_0\times\Pi_1)\neq 0$, then for all but finitely many primes $\lambda$ of $E$, the Bloch--Kato Selmer group
$\rH^1_f(F,\rho_{\Pi_0,\lambda}\otimes_{E_\lambda}\rho_{\Pi_1,\lambda}(n))$ vanishes.
\end{theorem}

\begin{remark}
In (a) of Theorem \ref{th:abstract_pre}, if the CM field $F$ is Galois or contains an imaginary quadratic field, then a very special inert prime of $F^+$ is simply a prime of $F^+$ that is inert in $F$, of degree $1$ over $\dQ$, whose underlying rational prime is odd and unramified in $F$.
\end{remark}

Now we state our result in the (Selmer) rank $1$ case. Let $\Pi_0$ and $\Pi_1$ be relevant representations of $\GL_{n_0}(\dA_F)$ and $\GL_{n_1}(\dA_F)$, respectively. Let $E\subseteq\dC$ be a strong coefficient field of both $\Pi_0$ and $\Pi_1$ (Definition \ref{de:weak_field}). Suppose that the global epsilon factor of $\Pi_0\times\Pi_1$ is $-1$. Then the Beilinson--Bloch--Kato conjecture predicts that if $L'(\frac{1}{2},\Pi_0\times\Pi_1)\neq 0$, then the Bloch--Kato Selmer group $\rH^1_f(F,\rho_{\Pi_0,\lambda}\otimes_{E_\lambda}\rho_{\Pi_1,\lambda}(n))$ has rank $1$. However, what we can prove now is half of this implication. Namely, for every prime $\lambda$ of $E$, we will construct explicitly an element $\graph_\lambda$ in (the direct sum of finitely many copies of) $\rH^1_f(F,\rho_{\Pi_0,\lambda}\otimes_{E_\lambda}\rho_{\Pi_1,\lambda}(n))$ in \S\ref{ss:main_1} as the image of the Abel--Jacobi map of the diagonal cycle of the product unitary Shimura variety (see \eqref{eq:selmer1} for the precise definition). In fact, by Conjecture \ref{co:aggp} and Beilinson's conjecture on the injectivity of the $\ell$-adic Abel--Jacobi map, the nonvanishing of $\graph_\lambda$ is equivalent to the nonvanishing of $L'(\frac{1}{2},\Pi_0\times\Pi_1)$. Our theorem in the rank $1$ case reads as follows.

\begin{theorem}[Theorem \ref{th:selmer1}]\label{th:selmer1_pre}
Let $\Pi_0$ and $\Pi_1$ be relevant representations of $\GL_{n_0}(\dA_F)$ and $\GL_{n_1}(\dA_F)$, respectively. Let $E\subseteq\dC$ be a strong coefficient field of both $\Pi_0$ and $\Pi_1$ (Definition \ref{de:weak_field}). Suppose that $F^+\neq\dQ$. For all admissible primes $\lambda$ of $E$ with respect to $(\Pi_0,\Pi_1)$, if $\graph_\lambda\neq 0$, then $\rH^1_f(F,\rho_{\Pi_0,\lambda}\otimes_{E_\lambda}\rho_{\Pi_1,\lambda}(n))$ has dimension $1$ over $E_\lambda$.
\end{theorem}

We also have an analogue of Theorem \ref{th:abstract_pre} in the rank $1$ case, whose statement we omit here.

\begin{remark}
In both Theorem \ref{th:selmer0_pre} and Theorem \ref{th:selmer1_pre}, the assumption that $F^+\neq\dQ$ if $n\geq 3$ can be lifted once Hypothesis \ref{hy:unitary_cohomology} is known for $N\geq 4$ when $F^+=\dQ$.

In fact, when $n=2$, we have a slightly different argument that can lift the restriction $F^+\neq\dQ$, and the assumptions (L6) and (L7) in Definition \ref{de:lambda} in all the results above.
\end{remark}

\if false

\subsection{Main results}
\label{ss:main}

Let $F/F^+$ be a totally imaginary quadratic extension of a totally real number field. We first state one of our main results that is least technical to understand.

\begin{theorem}[Corollary~\ref{co:sym}]\label{th:sym}
Let $n\geq 2$ be an integer. Let $A$ and $A'$ be two modular elliptic curves over $F^+$ such that $\End(A_{\ol{F}})=\End(A'_{\ol{F}})=\dZ$. Suppose that
\begin{enumerate}[label=(\alph*)]
  \item $A_{\ol{F}}$ and $A'_{\ol{F}}$ are not isogenous to each other;

  \item both $\Sym^{n-1}A$ and $\Sym^nA'$ are modular; and

  \item $F^+\neq\dQ$ if $n\geq 3$.
\end{enumerate}
If the (central critical) $L$-value $L(n,\Sym^{n-1}A_F\times\Sym^nA'_F)$ does not vanish, then the Bloch--Kato Selmer group
\[
\rH^1_f(F,\Sym^{n-1}\rH^1_\et(A_{\ol{F}},\dQ_\ell)\otimes_{\dQ_\ell}\Sym^n\rH^1_\et(A'_{\ol{F}},\dQ_\ell)(n))
\]
vanishes for all but finitely many rational primes $\ell$.
\end{theorem}

\begin{remark}
The finite set of rational primes $\ell$ that are excluded in Theorem \ref{th:sym} can be effectively bounded. We now explain the three conditions in Theorem \ref{th:sym}.
\begin{enumerate}[label=(\alph*)]
  \item is necessary since otherwise (L3) and (L5) in Definition \ref{de:lambda} fail for all rational primes $\ell$.

  \item is necessary since our approach only applies to Galois representations arising from automorphic representations. We summarise the current knowledge on the modularity of symmetric powers of elliptic curves in Remark \ref{re:sym}.

  \item is necessary only for technical reasons. First, we do not know Hypothesis \ref{hy:unitary_cohomology}, which concerns cohomology of unitary Shimura varieties, yet for $N\geq 4$ if $F^+=\dQ$. Second, we do not have (an appropriate replacement for) Proposition \ref{th:generic}, a result generalizing \cite{CS17}, when $F^+=\dQ$. Indeed, as long as we have these results as expected, (c) can be lifted.
\end{enumerate}
\end{remark}

Theorem \ref{th:sym} is a special case of a more general result concerning the Bloch--Kato Selmer groups of Galois representations associated to conjugate self-dual automorphic representations. To reduce the burden of long and technical terminology in the future, we first introduce the following definition, which will serve for the entire article.

\begin{definition}\label{de:relevant}
We say that a complex representation $\Pi$ of $\GL_N(\dA_F)$ with $N\geq 1$ is \emph{relevant} if
\begin{enumerate}
  \item $\Pi$ is an irreducible cuspidal automorphic representation;

  \item $\Pi\circ\tc\simeq\Pi^\vee$, where $\tc\in\Gal(F/F^+)$ is the complex conjugation;

  \item for every archimedean place $\tau$ of $F$, $\Pi_\tau$ is isomorphic to the (irreducible) principal series representation induced by the characters $(\arg^{N-1},\arg^{N-3},\dots,\arg^{3-N},\arg^{1-N})$, where $\arg\colon\dC^\times\to\dC^\times$ is the \emph{argument character} defined by the formula $\arg(z)\coloneqq z/\sqrt{z\ol{z}}$.
\end{enumerate}
If $\Pi$ is relevant, then it is regular algebraic in the sense of \cite{Clo90}*{Definition~3.12}.
\end{definition}

Now we can state our main result in the context of automorphic representations, of which Theorem \ref{th:sym} is a special case. Till the end of the next subsection, we will take an integer $n\geq 2$, and denote by $n_0$ and $n_1$ the unique even and odd numbers in $\{n,n+1\}$, respectively.

\begin{theorem}[Theorem \ref{th:selmer0}]\label{th:selmer0_pre}
Let $\Pi_0$ and $\Pi_1$ be relevant representations of $\GL_{n_0}(\dA_F)$ and $\GL_{n_1}(\dA_F)$, respectively. Let $E\subseteq\dC$ be a strong coefficient field of both $\Pi_0$ and $\Pi_1$ (Definition \ref{de:weak_field}). Suppose that $F^+\neq\dQ$ if $n\geq 3$. If $L(\frac{1}{2},\Pi_0\times\Pi_1)\neq 0$, then for all admissible primes $\lambda$ of $E$ with respect to $(\Pi_0,\Pi_1)$, the Bloch--Kato Selmer group $\rH^1_f(F,\rho_{\Pi_0,\lambda}\otimes_{E_\lambda}\rho_{\Pi_1,\lambda}(n))$ vanishes. Here, $\rho_{\Pi_\alpha,\lambda}$ is the Galois representation of $F$ with coefficients in $E_\lambda$ associated to $\Pi_\alpha$ for $\alpha=0,1$, as described in Proposition \ref{pr:galois} and Definition \ref{de:weak_field}.
\end{theorem}

In fact, Theorem \ref{th:selmer0} is slightly stronger than the one stated here.

\begin{remark}
The notion of admissible primes appearing in Theorem \ref{th:selmer0_pre} is introduced in Definition \ref{de:lambda}, which consists of a long list of assumptions, some of which are rather technical. Here, we would like to comment on the essence of these assumptions.
\begin{enumerate}
  \item[(L1,2)] are elementary and exclude only finitely many primes $\lambda$.

  \item[(L3)] is expected to hold for every prime $\lambda$ if and only if the (conjectural) automorphic product $\Pi_0\boxtimes\Pi_1$, as an irreducible admissible representation of $\GL_{n(n+1)}(\dA_F)$, remains cuspidal.

  \item[(L4)] is expected to hold for all but finitely many primes $\lambda$.

  \item[(L5)] is basically saying that, under (L4), the image of the pair of residual Galois representations $(\bar\rho_{\Pi_0,\lambda},\bar\rho_{\Pi_1,\lambda})$ contains an element of a particular form. It is expected to hold for all but finitely many primes $\lambda$ if the two automorphic representations $\Pi_0$ and $\Pi_1$ are not correlated in some manner. For example, when $n=2$, we expect that as long as $\Pi_1$ is not an automorphic twist of $\Sym^2\Pi_0$ after any base change, then (L5) holds for all but finitely many primes $\lambda$.

  \item[(L6)] is a technical assumption that is only used in the argument of an R=T theorem concerning Galois deformations in \cite{LTXZZ}. It is expected to hold for all but finitely many primes $\lambda$ (see \cite{LTXZZ}*{\S4.2}).

  \item[(L7)] is a technical assumption for the vanishing of certain Hecke localized cohomology of unitary Shimura varieties off middle degree. In fact, when $F^+\neq\dQ$, (L7) holds for all but finitely many primes $\lambda$ by Corollary \ref{co:generic}.
\end{enumerate}
\end{remark}

In fact, we have dedicated ourselves to obtaining the following family of abstract examples in which all but finitely many primes are admissible. Note that neither the following theorem nor Theorem \ref{th:sym} implies each other.

\begin{theorem}[Corollary \ref{co:abstract}]\label{th:abstract_pre}
Let $\Pi_0$, $\Pi_1$, and $E$ be as in Theorem \ref{th:selmer0_pre}. Suppose that
\begin{enumerate}[label=(\alph*)]
  \item there exists a very special inert prime $\fp$ of $F^+$ (Definition \ref{de:special_inert}) such that $\Pi_{0,\fp}$ is Steinberg, and $\Pi_{1,\fp}$ is unramified whose Satake parameter contains $1$ exactly once;\footnote{Note that the Satake parameter of $\Pi_{1,\fp}$ has to contain $1$ at least once by Definition \ref{de:relevant}(2).}

  \item for $\alpha=0,1$, there exists a nonarchimedean place $w_\alpha$ of $F$ such that $\Pi_{\alpha,w_\alpha}$ is supercuspidal; and

  \item $F^+\neq\dQ$ if $n\geq 3$.
\end{enumerate}
If $L(\frac{1}{2},\Pi_0\times\Pi_1)\neq 0$, then for all but finitely many primes $\lambda$ of $E$, the Bloch--Kato Selmer group
$\rH^1_f(F,\rho_{\Pi_0,\lambda}\otimes_{E_\lambda}\rho_{\Pi_1,\lambda}(n))$ vanishes.
\end{theorem}

\begin{remark}
In (a) of Theorem \ref{th:abstract_pre}, if the CM field $F$ is Galois or contains an imaginary quadratic field, then a very special inert prime of $F^+$ is simply a prime of $F^+$ that is inert in $F$, of degree $1$ over $\dQ$, whose underlying rational prime is odd and unramified in $F$.
\end{remark}

Now we state our result in the (Selmer) rank $1$ case. Let $\Pi_0$ and $\Pi_1$ be relevant representations of $\GL_{n_0}(\dA_F)$ and $\GL_{n_1}(\dA_F)$, respectively. Let $E\subseteq\dC$ be a strong coefficient field of both $\Pi_0$ and $\Pi_1$ (Definition \ref{de:weak_field}). Suppose that the global epsilon factor of $\Pi_0\times\Pi_1$ is $-1$. Then the Beilinson--Bloch--Kato conjecture predicts that if $L'(\frac{1}{2},\Pi_0\times\Pi_1)\neq 0$, then the Bloch--Kato Selmer group $\rH^1_f(F,\rho_{\Pi_0,\lambda}\otimes_{E_\lambda}\rho_{\Pi_1,\lambda}(n))$ has rank $1$. However, what we can prove now is half of this implication. Namely, for every prime $\lambda$ of $E$, we will construct explicitly an element $\graph_\lambda$ in (the direct sum of finitely many copies of) $\rH^1_f(F,\rho_{\Pi_0,\lambda}\otimes_{E_\lambda}\rho_{\Pi_1,\lambda}(n))$ in \S\ref{ss:main_1} as the image of the Abel--Jacobi map of the diagonal cycle of the product unitary Shimura variety (see \eqref{eq:selmer1} for the precise definition). In fact, by Conjecture \ref{co:aggp} and Beilinson's conjecture on the injectivity of the $\ell$-adic Abel--Jacobi map, the nonvanishing of $\graph_\lambda$ is equivalent to the nonvanishing of $L'(\frac{1}{2},\Pi_0\times\Pi_1)$. Our theorem in the rank $1$ case reads as follows.

\begin{theorem}[Theorem \ref{th:selmer1}]\label{th:selmer1_pre}
Let $\Pi_0$ and $\Pi_1$ be relevant representations of $\GL_{n_0}(\dA_F)$ and $\GL_{n_1}(\dA_F)$, respectively. Let $E\subseteq\dC$ be a strong coefficient field of both $\Pi_0$ and $\Pi_1$ (Definition \ref{de:weak_field}). Suppose that $F^+\neq\dQ$ if $n\geq 3$. For all admissible primes $\lambda$ of $E$ with respect to $(\Pi_0,\Pi_1)$, if $\graph_\lambda\neq 0$, then $\rH^1_f(F,\rho_{\Pi_0,\lambda}\otimes_{E_\lambda}\rho_{\Pi_1,\lambda}(n))$ has dimension $1$ over $E_\lambda$.
\end{theorem}

In fact, Theorem \ref{th:selmer1} is slightly stronger than the one stated here. We also have an analogue of Theorem \ref{th:abstract_pre} in the rank $1$ case, whose statement we omit here.

\begin{remark}
In both Theorem \ref{th:selmer0_pre} and Theorem \ref{th:selmer1_pre}, the assumption that $F^+\neq\dQ$ if $n\geq 3$ can be lifted once Hypothesis \ref{hy:unitary_cohomology} is known for $N\geq 4$ when $F^+=\dQ$.
\end{remark}

\fi

\subsection{Road map for the article}

The very basic idea of bounding Selmer groups as in our main theorems follows from Kolyvagin \cite{Kol90}, namely, we construct a system of torsion Galois cohomology classes serving as annihilators of (reduction of) Selmer groups. However, our system is not a generalization of the Euler--Kolyvagin system originally constructed by Kolyvagin. Instead, our system is constructed via level-raising congruences,\footnote{What we need from level-raising congruences is much more than merely the existence part. In fact, we have to identify the level-raising explicitly through the geometry of the special fiber of some Shimura variety, for which we call \emph{arithmetic level-raising}.} which was first introduced by Bertolini and Darmon in the case of Heegner points in the study of certain Iwasawa main conjecture of elliptic curves \cite{BD05}. The first example where such level-raising system was used to bound Selmer groups beyond the Heegner point case was performed by one of us in \cite{Liu1}, for the so-called twisted triple product automorphic motives. In the sequels \cite{Liu2} and \cite{LT}, the case of the so-called cubic triple product automorphic motives was also studied. From this point of view, our current article is a vast generalization of the previous results mentioned above.


The following is a road map for reading the main part of the article, where we indicate the need from the five appendices in the parentheses.
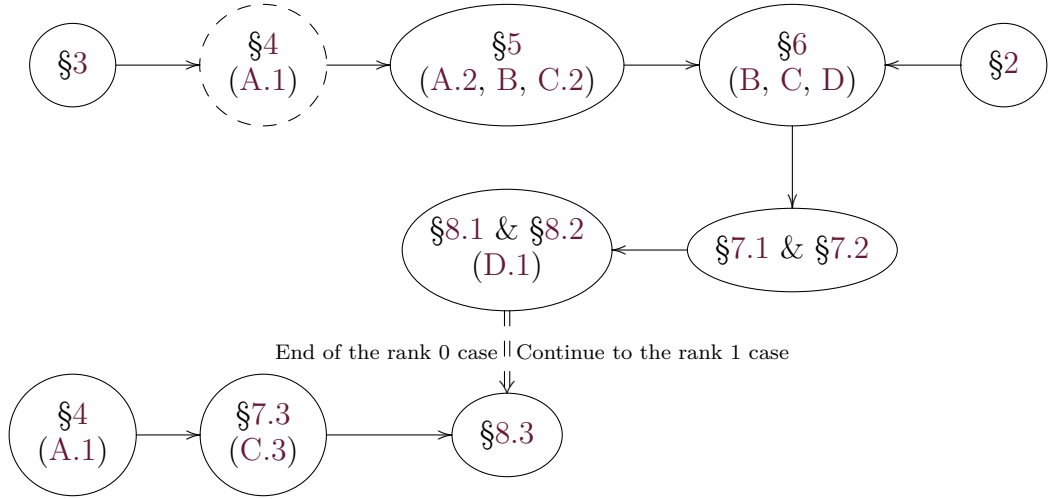
\[
\xymatrix{
*+++[o][F]{\txt{\S\ref{ss:3} }}*
\ar[r]&
*+++[o][F--]{\txt{\S\ref{ss:qs} \\ (\ref{ss:dl_smooth}) }}*
\ar[r]&
*+++[o][F]{\txt{\S\ref{ss:ns} \\ (\ref{ss:dl_semistable}, \ref{ss:a}, \ref{ss:omega}) }}*
\ar[r]&
*+++[o][F]{\txt{\S\ref{ss:6} \\ (\ref{ss:a}, \ref{ss:b}, \ref{ss:c}) }}* \ar[d]
& *+++[o][F]{\txt{\S\ref{ss:2} }}* \ar[l]
\\
&&*+++[o][F]{\txt{\S\ref{ss:admissible_prime} \& \S\ref{ss:main_0} \\ (\ref{ss:vanishing}) }}*
\ar@{==>}[d]_-{\text{End of the rank 0 case}}^-{\text{Continue to the rank 1 case}}
&*+++[o][F]{\txt{\S\ref{ss:setup} \& \S\ref{ss:first_reciprocity} }}* \ar[l]
\\
*+++[o][F]{\txt{\S\ref{ss:qs} \\ (\ref{ss:dl_smooth}) }}*
\ar[r]&
*+++[o][F]{\txt{\S\ref{ss:second_reciprocity} \\ (\ref{ss:endoscopic}) }}*
\ar[r]&
*+++[o][F]{\txt{\S\ref{ss:main_1} }}*
}
\]
The proof of Theorem \ref{th:selmer1_pre} is based on the proof of Theorem \ref{th:selmer0_pre}. We may regard the transition from the rank $0$ case to the rank $1$ case as an induction step. As seen from the road map, for the rank $0$ case alone, \S\ref{ss:qs}, \S\ref{ss:dl_smooth}, \S\ref{ss:second_reciprocity}, and, of course, \S\ref{ss:main_1} are not needed. However, we strongly recommend the readers to go through \S\ref{ss:qs} even if they are only interested in the rank $0$ case, as \S\ref{ss:qs} is an appropriate warm-up for reading \S\ref{ss:ns}, which is parallel but much more complicated.

In what follows, we explain the main steps in the proof of Theorem \ref{th:selmer0_pre}. Some of the notations in the rest of this subsection are \emph{ad hoc}, only for the purpose of explaining ideas, hence will be obsolete or differ from the main text.

The initial step (which although will not appear until \S\ref{ss:main_0}) is to translate the condition that $L(\frac{1}{2},\Pi_0\times\Pi_1)\neq 0$ into a more straightforward statement. This is exactly the content of the global Gan--Gross--Prasad conjecture \cite{GGP12}. In fact, as stated in Lemma \ref{le:ggp}, we may construct a pair of hermitian spaces $(\rV^\circ_n,\rV^\circ_{n+1})$ over $F$ (with respect to $F/F^+$) in which $\rV^\circ_n$ is totally positive definite of rank $n$, and $\rV^\circ_{n+1}=\rV^\circ_n\oplus F\cdot 1$ where $1$ has norm $1$. For $\alpha=0,1$, put $\Sh(\rV^\circ_{n_\alpha})\coloneqq\rU(\rV^\circ_{n_\alpha})(F^+)\backslash\rU(\rV^\circ_{n_\alpha})(\dA_{F^+}^\infty)$ as a Shimura (pro-)set. We may further find cuspidal automorphic representations $\pi_0$ and $\pi_1$ contained in the space of locally constant functions on $\Sh(\rV^\circ_{n_0})$ and $\Sh(\rV^\circ_{n_1})$ satisfying $\BC(\pi_0)\simeq\Pi_0$ and $\BC(\pi_1)\simeq\Pi_1$, respectively, such that
\begin{align}\label{eq:ggp}
\cP(f_0,f_1)\coloneqq\int_{\Sh(\rV^\circ_n)}f_0(h)f_1(h)\rd h\neq 0
\end{align}
for some $f_0\in\pi_0$ and $f_1\in\pi_1$ valued in $O_E$. Such result was first obtained by one of us \cite{Zha} under some local restrictions. Those restrictions are all lifted till very recently through some new techniques in the study of trace formulae \cite{BPLZZ}. In what follows, we will fix open compact subgroups of $\rU(\rV^\circ_{n_0})(\dA_{F^+}^\infty)$ and $\rU(\rV^\circ_{n_1})(\dA_{F^+}^\infty)$ that fix $f_0$ and $f_1$, respectively, and will carry them implicitly in the notation.

The next step is to bring the set $\Sh(\rV^\circ_{n_\alpha})$ into arithmetic geometry so that the period \eqref{eq:ggp} can be related to certain Galois cohomology classes. Now we choose a special inert prime $\fp$ of $F^+$ (see Definition \ref{de:special_inert}) with sufficiently large underlying rational prime $p$, such that all data appearing so far are unramified above $p$. For $\alpha=0,1$, we attach to $\rV^\circ_{n_\alpha}$ \emph{canonically} a strictly semistable scheme $\bM_\fp(\rV^\circ_{n_\alpha})$ over $\Spec\dZ_{p^2}$ of relative dimension $n_\alpha-1$, whose complex generic fiber is \emph{non-canonically} isomorphic to the disjoint union of finitely many Shimura varieties attached to the \emph{nearby} hermitian space of $\rV^\circ_{n_\alpha}$ by changing local components at $\fp$ and one archimedean place. Moreover, we can write its special fiber $\rM_\fp(\rV^\circ_{n_\alpha})$ over $\Spec\dF_{p^2}$ as the union of $\rM^\circ_\fp(\rV^\circ_{n_\alpha})$ and $\rM^\bullet_\fp(\rV^\circ_{n_\alpha})$, in which $\rM^\circ_\fp(\rV^\circ_{n_\alpha})$ is geometrically a $\dP^{n_\alpha-1}$-fibration over the Shimura set $\Sh(\rV^\circ_{n_\alpha})$. The other stratum $\rM^\bullet_\fp(\rV^\circ_{n_\alpha})$, which is rather mysterious, will also be involved in the later computation. In fact, one key effort we make is to show that only the basic locus of the stratum $\rM^\bullet_\fp(\rV^\circ_{n_\alpha})$ will play a role in the computation. For the basic locus, we show that its normalization is geometrically a fibration over the Shimura set $\Sh(\rV^\circ_{n_\alpha})$ (but with a slightly different level structure at $\fp$) by certain Deligne--Lusztig varieties of dimension $r_\alpha\coloneqq\floor{\tfrac{n_\alpha}{2}}$, introduced in \S\ref{ss:dl_semistable}. The study of various geometric aspects of the scheme $\bM_\fp(\rV^\circ_{n_\alpha})$, including its associated Rapoport--Zink spectral sequence and its functorial behavior from $n$ to $n+1$, will be carried out in \S\ref{ss:ns}.

The automorphic input will be thrown into the scheme $\bM_\fp(\rV^\circ_{n_\alpha})$ from the third step, in \S\ref{ss:6}, where we study the local Galois cohomology of certain cohomology of $\bM_\fp(\rV^\circ_{n_\alpha})$ localized at some Hecke ideals. More precisely, we fix an admissible prime $\lambda$ of $E$ with respect to $(\Pi_0,\Pi_1)$, and denote by $O_\lambda$ and $k_\lambda$ the ring of integers and the residue field of $E_\lambda$, respectively. For $\alpha=0,1$, the Satake parameters of $\Pi_\alpha$ induce a homomorphism $\phi_\alpha\colon\dT_{n_\alpha}\to k_\lambda$ with kernel $\fm_\alpha$, where $\dT_{n_\alpha}$ is a certain abstract spherical Hecke algebra for unitary groups of rank $n_\alpha$. When $\alpha=0$ (resp.\ $\alpha=1$), we need to study the singular (resp.\ unramified) part of the local Galois cohomology
\begin{align}\label{eq:introduction_1}
\rH^1(\dQ_{p^2},\rH^{n_\alpha-1}_\fT(\ol\rM_\fp(\rV^\circ_{n_\alpha}),\rR\Psi O_\lambda(r_\alpha))_{\fm_\alpha}),
\end{align}
where $\ol\rM_\fp(\rV^\circ_{n_\alpha})\coloneqq\rM_\fp(\rV^\circ_{n_\alpha})\otimes_{\dF_{p^2}}\ol\dF_p$, and $\rH_\fT$ denotes the certain invariant part of the \'{e}tale cohomology (a subtlety that can be ignored at this moment). The question boils down to the arithmetic level-raising phenomenon (resp.\ existence of Tate cycles) when $\alpha=0$ (resp.\ $\alpha=1$). However, in both cases, we have to rely on the recent progress on the Tate conjecture for Shimura varieties achieved by two of us \cite{XZ17}. Now we would like to continue the discussion on the case where $\alpha=0$, since it is more interesting and more involved, and omit the case where $\alpha=1$. The first key point is to figure out the correct condition such that the level-raising phenomenon (namely, from unramified to mildly ramified at the place $\fp$) happens on the cohomology \eqref{eq:introduction_1} in a way that can be understood: we say that $\fp$ is a \emph{level-raising prime} with respect to $\lambda$ if $\ell\nmid p(p^2-1)$ where $\ell$ is the underlying rational prime of $\lambda$, and the $\modulo\lambda$ Satake parameter of $\Pi_{0,\fp}$ contains the pair $\{p,p^{-1}\}$ exactly once and does not contain the pair $\{-1,-1\}$. Suppose that $\fp$ is such a prime, we show that there is a canonical isomorphism
\begin{align}\label{eq:introduction_2}
\rH^1_\sing(\dQ_{p^2},\rH^{n_0-1}_\fT(\ol\rM_\fp(\rV^\circ_{n_0}),\rR\Psi O_\lambda(r_0))/\fm_0)
\simeq O_\lambda[\Sh(\rV^\circ_{n_0})]/\fm_0
\end{align}
of $k_\lambda$-vector spaces of finite dimension. Note that by our condition on $\fp$, the right-hand side of \eqref{eq:introduction_2} is nonvanishing, which implies that the left-hand side is also nonvanishing; in other words, we see the level-raising phenomenon in $\rH^{n_0-1}_\fT(\ol\rM_\fp(\rV^\circ_{n_0}),\rR\Psi O_\lambda(r_0))$. The proof of \eqref{eq:introduction_2} is the technical heart of this article (for example, it uses materials from all of the five appendices). Through studying the geometry and intersection theory on the special fiber $\rM_\fp(\rV^\circ_{n_0})$ in \S\ref{ss:ns} and some of the appendices, we can conclude that $O_\lambda[\Sh(\rV^\circ_{n_0})]/\fm_0$ is canonically a subquotient of $\rH^1_\sing(\dQ_{p^2},\rH^{n_0-1}_\fT(\ol\rM_\fp(\rV^\circ_{n_0}),\rR\Psi O_\lambda(r_0))/\fm_0)$. Thus, it remains to show that the two sides of \eqref{eq:introduction_2} have the same cardinality. For this, we use the theory of Galois deformations. We construct a global Galois deformation ring $\sfR^\mix$ over $O_\lambda$ with two quotient rings $\sfR^\unr$ and $\sfR^\ram$, together with a natural $\sfR^\unr$-module $\sfH^\unr$ and a natural $\sfR^\ram$-module $\sfH^\ram$. They satisfy the following relation: if we put $\sfR^{\r{cong}}\coloneqq\sfR^\unr\otimes_{\sfR^\mix}\sfR^\ram$, which is an Artinian ring over $O_\lambda$, then we have natural isomorphisms
\begin{align*}
\sfH^\unr\otimes_{\sfR^\unr}\sfR^{\r{cong}}\otimes_{O_\lambda}k_\lambda &\simeq O_\lambda[\Sh(\rV^\circ_{n_0})]/\fm_0, \\
\sfH^\ram\otimes_{\sfR^\ram}\sfR^{\r{cong}}\otimes_{O_\lambda}k_\lambda &\simeq \rH^1_\sing(\dQ_{p^2},\rH^{n_0-1}_\fT(\ol\rM_\fp(\rV^\circ_{n_0}),\rR\Psi O_\lambda(r_0))/\fm_0).
\end{align*}
Thus, we only need to show that $\sfH^\unr$ and $\sfH^\ram$ are both finite free over $\sfR^\unr$ and $\sfR^\ram$, respectively, of the same rank. The finite-freeness follows from an R=T theorem, proved in \cite{LTXZZ}. The comparison of ranks can be performed over $E_\lambda$, which turns out to be an automorphic problem and is solved in \S\ref{ss:raising2} based on \S\ref{ss:dimension}. Summarizing the discussion above, we obtain \eqref{eq:introduction_2}. In practice, we also need a $\modulo\lambda^m$ version of \eqref{eq:introduction_2}.

The fourth step is to merge the study of \eqref{eq:introduction_1} for $n_0$ and $n_1$ together, to obtain the so-called \emph{first explicit reciprocity law} for the Rankin--Selberg product of Galois representations. As an application, we construct a system of torsion Galois cohomology classes whose image in the singular part of the local Galois cohomology at $\fp$ of the product Galois representation is controlled by the period integral \eqref{eq:ggp}. This step is sort of routine, once we have enough knowledge on \eqref{eq:introduction_1}; it is completed in \S\ref{ss:first_reciprocity}.

The final step of the proof of Theorem \ref{th:selmer0_pre} will be performed in \S\ref{ss:main_0}, where we use the system of torsion Galois cohomology classes constructed in the previous step, together with some Galois theoretical facts from \S\ref{ss:2}, to bound the Selmer group, which is possible due to the nonvanishing of \eqref{eq:ggp}.

\subsection{Notations and conventions}
\label{ss:notation}

In this subsection, we setup some common notations and conventions for the entire article, including appendices, unless otherwise specified. The notations in the previous two subsections will not be relied on from this moment, and should not be kept for further reading.

\subsubsection*{Generalities:}

\begin{itemize}[label={\ding{109}}]
  \item Denote by $\dN=\{0,1,2,3,\dots\}$ the monoid of nonnegative integers.

  \item We only apply the operation $\sqrt{\phantom{a}}$ to positive real numbers, which takes values in positive real numbers as well.

  \item For a set $S$, we denote by $\CF_S$ the characteristic function of $S$.

  \item The eigenvalues or generalized eigenvalues of a matrix over a field $k$ are counted with multiplicity (namely, dimension of the corresponding eigenspace or generalized eigenspace); in other words, they form a multi-subset of an algebraic extension of $k$.

  \item For every rational prime $p$, we fix an algebraic closure $\ol\dQ_p$ of $\dQ_p$ with the residue field $\ol\dF_p$. For every integer $r\geq 1$, we denote by $\dQ_{p^r}$ the subfield of $\ol\dQ_p$ that is an unramified extension of $\dQ_p$ of degree $r$, by $\dZ_{p^r}$ its ring of integers, and by $\dF_{p^r}$ its residue field.

  \item For a nonarchimedean place $v$ of a number field $K$, we write $\|v\|$ for the cardinality of the residue field of $K_v$.

  \item We use standard notations from the category theory. The category of sets is denoted by $\Set$. For a category $\fC$, we denote by $\fC^{\r{op}}$ its opposite category, and denote by $\fC_{/A}$ the category of morphisms to $A$ for an object $A$ of $\fC$. For another category $\fD$, we denote by $\Fun(\fC,\fD)$ the category of functors from $\fC$ to $\fD$. In particular, we denote by $\sfP\fC\coloneqq\Fun(\fC^{\r{op}},\Set)$ the category of presheaves on $\fC$, which contains $\fC$ as a full subcategory by the Yoneda embedding. Isomorphisms in a category will be indicated by $\simeq$. We also use the symbol $\obj$ to indicate a virtual object.

  \item For an algebra $A$, we denote by $\Mod(A)$ the category of left $A$-modules.

  \item All rings are commutative and unital; and ring homomorphisms preserve units. For a (topological) ring $L$, a (topological) $L$-ring is a (topological) ring $R$ together with a (continuous) ring homomorphism from $L$ to $R$. However, we use the word \emph{algebra} in the general sense, which is not necessarily commutative or unital.

  \item If a base ring is not specified in the tensor operation $\otimes$, then it is $\dZ$.

  \item For a ring $L$ and a set $S$, denote by $L[S]$ the $L$-module of $L$-valued functions on $S$ of finite support.
\end{itemize}


\subsubsection*{Algebraic geometry:}

\begin{itemize}[label={\ding{109}}]
  \item We denote by the category of schemes by $\Sch$ and its full subcategory of locally Noetherian schemes by $\Sch'$. For a scheme $S$ (resp. Noetherian scheme $S$), we denote by $\Sch_{/S}$ (resp.\ $\Sch'_{/S}$) the category of $S$-schemes (resp.\ locally Noetherian $S$-schemes). If $S=\Spec R$ is affine, we also write $\Sch_{/R}$ (resp.\ $\Sch'_{/R}$) for $\Sch_{/S}$ (resp.\ $\Sch'_{/S}$).

  \item The structure sheaf of a scheme $X$ is denoted by $\cO_X$.

  \item For a scheme $X$ over an affine scheme $\Spec R$ and an $R$-ring $S$, we write $X\otimes_RS$ or even $X_S$ for $X\times_{\Spec R}\Spec S$.

  \item For a scheme $S$ in characteristic $p$ for some rational prime $p$, we denote by $\sigma\colon S\to S$ the absolute $p$-power Frobenius morphism. For a perfect field $\kappa$ of characteristic $p$, we denote by $W(\kappa)$ its Witt ring, and by abuse of notation, $\sigma\colon W(\kappa)\to W(\kappa)$ the canonical lifting of the $p$-power Frobenius map.

  \item For a smooth morphism $S\to T$ of schemes, we denote by $\cT_{S/T}$ the relative tangent sheaf, which is a locally free $\cO_S$-module.

  \item For a scheme $S$ and a locally free $\cO_S$-module $\cV$ of finite rank, we denote by $\dP(\cV)\to S$ the moduli scheme of quotient line bundles of $\cV$ over $S$, known as the \emph{projective fibration} associated to $\cV$.

  \item For a scheme $S$ and (sheaves of) $\cO_S$-modules $\cF$ and $\cG$, we denote by $\HOM(\cF,\cG)$ the quasi-coherent sheaf of $\cO_S$-linear homomorphisms from $\cF$ to $\cG$.

  \item For two positive integers $r,s$, we denote by $\rM_{r,s}$ the scheme over $\dZ$ of $r$-by-$s$ matrices, and put $\rM_r\coloneqq\rM_{r,r}$ for short; we also denote by $\GL_r\subseteq\rM_r$ the subscheme of invertible $r$-by-$r$ matrices. Then $\GL_1$ is simply the multiplicative group $\bG_m\coloneqq\dZ[T,T^{-1}]$; but we will distinguish between $\GL_1$ and $\bG_m$ according to the context.

  \item For a number field $K$, a commutative group scheme $G\to S$ equipped with an action by $O_K$ over some base scheme $S$, and an ideal $\fa\subset O_K$, we denote by $G[\fa]$ the maximal closed subgroup scheme of $G$ annihilated by all elements in $\fa$.

  \item By a \emph{coefficient ring} for \'{e}tale cohomology, we mean either a finite ring, or a finite extension of $\dQ_\ell$, or the ring of integers of a finite extension of $\dQ_\ell$. In the latter two cases, we regard the \'{e}tale cohomology as the continuous one. We say that a coefficient ring $L$ is \emph{$n$-coprime} for a positive integer $n$ if $n$ is invertible in $L$ in the first case, and $\ell\nmid n$ in the latter two cases.
\end{itemize}

\subsubsection*{Group theory:}

Let $G$ and $\tilde\Gamma$ be groups, and $\Gamma$ a subgroup of $\tilde\Gamma$. Let $L$ be a ring.

\begin{itemize}[label={\ding{109}}]
  \item Denote by $\Gamma^\ab$ the maximal abelian quotient of $\Gamma$.

  \item For a homomorphism $\rho\colon\Gamma\to\GL_r(L)$ for some $r\geq 1$, we denote by $\rho^\vee\colon\Gamma\to\GL_r(L)$ the contragredient homomorphism, which is defined by the formula $\rho^\vee(x)=\tp\rho(x)^{-1}$ for every $x\in\Gamma$.

  \item For a homomorphism $\rho\colon\Gamma\to G$ and an element $\gamma\in\tilde\Gamma$ that normalizes $\Gamma$, we let $\rho^\gamma\colon\Gamma\to G$ be the homomorphism defined by $\rho^\gamma(x)=\rho(\gamma x\gamma^{-1})$ for every $x\in\Gamma$.

  \item We say that two homomorphisms $\rho_1,\rho_2\colon\Gamma\to G$ are conjugate if there exists an element $g\in G$ such that $\rho_1=g\circ\rho_2\circ g^{-1}$, that is, $\rho_1(x)=g\rho_2(x)g^{-1}$ for every $x\in\Gamma$.

  \item The $L$-module $L[G]$ is naturally an $L$-algebra, namely, the group algebra of $G$ with coefficients in $L$.

  \item Suppose that $G$ is a locally compact and totally disconnected topological group. For an open compact subgroup $K$ of $G$, the $L$-module $L[K\backslash G/K]$ (of bi-$K$-invariant compactly supported $L$-valued functions on $G$) is naturally an $L$-algebra, where the algebra structure is given by the composition of cosets. In particular, the unit element of $L[K\backslash G/K]$ is always $\CF_K$.
\end{itemize}

\subsubsection*{Combinatorics:}

\begin{notation}\label{no:numerical}
We recall the $q$-analogues of binomial coefficients:
\[
[0]_q=1,\quad [n]_q=\frac{q^n-1}{q-1},\quad [n]_q!=[n]_q\cdot [n-1]_q\cdots[1]_q,
\quad\qbinom{n}{m}_q=\frac{[n]_q!}{[n-m]_q!\cdot [m]_q!}
\]
for integers $0\leq m\leq n$. For $r\geq 0$ and $q\in\dN$, we put
\begin{align*}
\begin{dcases}
\td_{r,q}\coloneqq\sum_{\delta=0}^r(-1)^\delta(2\delta+1)q^{\delta(\delta+1)}\qbinom{2r+1}{r-\delta}_{-q},\\
\td^\bullet_{r,q}\coloneqq\frac{1}{q+1}\(\td_{r,q}+\frac{(-q)^{r+1}-1}{q+1}(q+1)(q^3+1)\cdots(q^{2r-1}+1)\).
\end{dcases}
\end{align*}
\end{notation}

\subsubsection*{Ground fields:}

\begin{itemize}[label={\ding{109}}]
  \item Let $\tc\in\Aut(\dC/\dQ)$ be the complex conjugation.

  \item Throughout the article, we fix a \emph{subfield} $F\subseteq\dC$ that is a number field and is stable under $\tc$; it is assumed to be a CM field except in \S\ref{ss:2}.

  \item Let $F^+\subseteq F$ be the maximal subfield on which $\tc$ acts by the identity.

  \item Let $\ol{F}$ be \emph{the} Galois closure of $F$ in $\dC$. Put $\Gamma_F\coloneqq\Gal(\ol{F}/F)$ and $\Gamma_{F^+}\coloneqq\Gal(\ol{F}/F^+)$.

  \item Denote by $\Sigma_\infty$ (resp.\ $\Sigma^+_\infty$) the set of complex embeddings of $F$ (resp.\ $F^+$) with $\tau_\infty\in\Sigma_\infty$ (resp.\ $\ul\tau_\infty\in\Sigma^+_\infty$) the default one. For $\tau\in\Sigma_\infty$, we denote by $\tau^\tc$ the its complex conjugation.

  \item For every rational prime $p$, denote by $\Sigma^+_p$ the set of all $p$-adic places of $F^+$.

  \item Denote by $\Sigma^+_\bad$ the union of $\Sigma^+_p$ for all $p$ that ramifies in $F$.

  \item Denote by $\eta_{F/F^+}\colon\Gamma_{F^+}\to\{\pm1\}$ the character associated to the extension $F/F^+$.

  \item For every prime $\ell$, denote by $\epsilon_\ell\colon\Gamma_{F^+}\to\dZ_\ell^\times$ the $\ell$-adic cyclotomic character.
\end{itemize}

For every place $v$ of $F^+$, we
\begin{itemize}[label={\ding{109}}]
  \item put $F_v\coloneqq F\otimes_{F^+}F^+_v$; and define $\delta(v)$ to be $1$ (resp.\ $2$) if $v$ splits (resp.\ does not split) in $F$;

  \item fix an algebraic closure $\ol{F}^+_v$ of $F^+_v$ containing $\ol{F}$; and put $\Gamma_{F^+_v}\coloneqq\Gal(\ol{F}^+_v/F^+_v)$ as a subgroup of $\Gamma_{F^+}$;

  \item for a homomorphism $r$ from $\Gamma_{F^+}$ to another group, denote by $r_v$ the restriction of $r$ to the subgroup $\Gamma_{F^+_v}$.
\end{itemize}

For every nonarchimedean place $w$ of $F$, we
\begin{itemize}[label={\ding{109}}]
  \item identify the Galois group $\Gamma_{F_w}$ with $\Gamma_{F^+_v}\cap\Gamma_F$ (resp.\ $\tc(\Gamma_{F^+_v}\cap\Gamma_F)\tc$), where $v$ is the underlying place of $F^+$, if the embedding $F\hookrightarrow\ol{F}^+_v$ induces (resp.\ does not induce) the place $w$;

  \item let $\rI_{F_w}\subseteq\Gamma_{F_w}$ be the inertia subgroup;

  \item let $\kappa_w$ be the residue field of $F_w$, and identify its Galois group $\Gamma_{\kappa_w}$ with $\Gamma_{F_w}/\rI_{F_w}$;

  \item denote by $\phi_w\in\Gamma_{F_w}$ a lifting of the \emph{arithmetic} Frobenius element in $\Gamma_{\kappa_w}$.
\end{itemize}

\begin{definition}\label{de:strongly_disjoint}
We say that two subsets $\Sigma^+_1$ and $\Sigma^+_2$ of nonarchimedean places of $F^+$ are \emph{strongly disjoint} if there is no common rational prime underlying the places from both sets.
\end{definition}

\subsection*{Acknowledgements}

This article is the main outcome of the AIM SQuaREs project \emph{Geometry of Shimura varieties and arithmetic application to L-functions} conducted by the five authors from 2017 to 2019. We would like to express our sincere gratitude and appreciation to the American Institute of Mathematics for their constant and generous support of the project, and to the staff members at the AIM facility in San Jose, California for their excellent coordination and hospitality.

We would like to thank Sug~Woo~Shin and Yihang~Zhu for the discussion concerning Hypothesis \ref{hy:unitary_cohomology} and the endoscopic classification for unitary groups, Zipei~Nie and Jun~Su for suggesting a proof of a combinatorial lemma (Lemma \ref{le:enumeration_even_3}), Ana~Caraiani and Peter~Scholze for the discussion concerning \S\ref{ss:vanishing}, Kai-Wen Lan for the discussion concerning the reference \cite{LS18}, and Ruiqi~Bai and Murilo~Corato~Zanarella for correcting some errors in early drafts. Finally, we thank the anonymous referees for careful reading and many valuable suggestions.

The research of Y.~L. is partially supported by the NSF grant DMS--1702019 and a Sloan Research Fellowship. The research of L.~X. is partially supported by the NSF grant DMS--1502147 and DMS--1752703, the Chinese NSF grant under agreement No. NSFC--12071004, and a grant from the Chinese Ministry of Education. The research of W.~Z. is partially supported by the NSF grant DMS--1838118 and DMS--1901642. The research of X.~Z. is partially supported by the NSF grant DMS--1902239 and a Simons Fellowship.

\section{Galois cohomology and Selmer groups}
\label{ss:2}

In this section, we make the Galois theoretical preparation for the proof of the main theorems. Most discussions in this section are generalizations from \cites{Liu1,Liu2}. The material of this section will not be used until \S\ref{ss:6}. In \S\ref{ss:adic_modules}, we collect some lemmas on $\ell$-adic modules with certain group actions. In \S\ref{ss:local_galois}, we study local Galois cohomology. In \S\ref{ss:galois_lemma}, we perform the discussion that is typical for Kolyvagin's type of argument. The Selmer group and its variant will be introduced in \S\ref{ss:bloch_kato}. In \S\ref{ss:conjugate_selfdual}, we discuss extension of essentially conjugate self-dual representations. In \S\ref{ss:localization}, we study localization of Selmer groups. In \S\ref{ss:rankin_selberg}, we study an example related to the Rankin--Selberg product.

We will start from a more general setup in order to make the discussion applicable to the orthogonal case as well, which may be studied in the future. Thus, we fix a \emph{subfield} $F\subseteq\dC$ that is a number field, \emph{not} necessarily CM.

We fix an odd rational prime $\ell$ that is unramified in $F$, and consider a finite extension $E_\lambda/\dQ_\ell$, with the ring of integers $O_\lambda$ and the maximal ideal $\lambda$ of $O_\lambda$. We denote by $\dB_\cris$ Fontaine's crystalline period ring for $\dQ_\ell$, and recall from \S\ref{ss:notation} that $\epsilon_\ell\colon\Gamma_{F^+}\to\dZ_\ell^\times$ is the $\ell$-adic cyclotomic character.

\subsection{Preliminaries on $\ell$-adic modules with group actions}
\label{ss:adic_modules}

Let $\Gamma$ be a topological group and $L$ a $\dZ_\ell$-ring that is finite over either $\dZ_\ell$ or $\dQ_\ell$. Note that in this case, every finitely generated $L$-module is equipped with the natural $\ell$-adic topology.

\begin{notation}
We denote by $\Mod(\Gamma,L)$ the category of finitely generated $L$-modules equipped with a continuous action of $\Gamma$, and by $\Mod(\Gamma,L)_\tor$ (resp.\ $\Mod(\Gamma,L)_\free$) the full subcategory of $\Mod(\Gamma,L)$ consisting of those objects whose underlying $L$-modules are torsion (resp.\ free).
\end{notation}


\begin{definition}\label{de:weakly_semisimple}
We say that an $L[\Gamma]$-module $M$ is \emph{weakly semisimple} if
\begin{enumerate}
  \item $M$ is an object of $\Mod(\Gamma,L)$; and

  \item the natural map $M^\Gamma\to M_\Gamma$ is an isomorphism.
\end{enumerate}
\end{definition}

\begin{lem}\label{le:weakly_semisimple0}
Suppose that $\Gamma$ is isomorphic to $\widehat\dZ$. Let $M$ be an object of $\Mod(\Gamma,L)$. Then
\begin{enumerate}
  \item $M_\Gamma=0$ implies $M^\Gamma=0$;

  \item if the natural map $M^\Gamma\to M_\Gamma$ is surjective, then $M$ is weakly semisimple.
\end{enumerate}
\end{lem}

\begin{proof}
Take a topological generator $\gamma$ of $\Gamma$.

For (1), we have the exact sequence
\[
0 \to M^\Gamma \to M \xrightarrow{\gamma-1} M \to M_\Gamma \to 0.
\]
Since $M_\Gamma=0$, $\gamma-1\colon M\to M$ is surjective. As $M$ is Noetherian, it follows that $M^\Gamma=0$.

For (2), taking (continuous) $\Gamma$-cohomology of the short exact sequence
\[
0 \to M^\Gamma \to M \to M/M^\Gamma \to 0,
\]
we obtain the sequence
\[
\(M/M^\Gamma\)^\Gamma \to M^\Gamma \to M_\Gamma \to \(M/M^\Gamma\)_\Gamma \to 0.
\]
Since $M^\Gamma\to M_\Gamma$ is surjective, it follows that $\(M/M^\Gamma\)_\Gamma=0$. By (1), we have $\(M/M^\Gamma\)^\Gamma=0$, hence the map $M^\Gamma\to M_\Gamma$ is injective as well.

The lemma is proved.
\end{proof}

\begin{lem}\label{le:weakly_semisimple1}
Suppose that $\Gamma$ is isomorphic to $\widehat\dZ$.
\begin{enumerate}
  \item A finite direct sum of weakly semisimple $L[\Gamma]$-modules is weakly semisimple.

  \item A subquotient $L[\Gamma]$-module of a weakly semisimple $L[\Gamma]$-module is weakly semisimple.
\end{enumerate}
\end{lem}

\begin{proof}
Part (1) is obvious.

For (2), let $M$ be a weakly semisimple $L[\Gamma]$-module and consider a short exact sequence
\[
0 \to N \to M \to Q \to 0
\]
of $L[\Gamma]$-module. We obtain the diagram
\begin{align}\label{eq:weakly_semisimple}
\xymatrix{
0 \ar[r] & N^\Gamma \ar[r]\ar[d] & M^\Gamma \ar[r]\ar[d]^-\simeq & Q^\Gamma \ar[d] \\
& N_\Gamma \ar[r] & M_\Gamma \ar[r] & Q_\Gamma \ar[r] & 0
}
\end{align}
in which the middle vertical arrow is an isomorphism. It follows that $Q^\Gamma\to Q_\Gamma$ is surjective, which implies that $Q$ is weakly semisimple by Lemma \ref{le:weakly_semisimple0}(2). It also follows that $M^\Gamma\to Q^\Gamma$ is surjective, which implies that $N_\Gamma\to M_\Gamma$ is injective. Thus, \eqref{eq:weakly_semisimple} is an isomorphism of exact sequences. Part (2) is proved.
\end{proof}

\begin{lem}\label{le:weakly_semisimple2}
Suppose that $\Gamma$ is isomorphic to $\widehat\dZ$. Let $M$ be an object of $\Mod(\Gamma,O_\lambda)_\free$. Suppose that $M\otimes_{O_\lambda}O_\lambda/\lambda$ is weakly semisimple, and $\dim_{E_\lambda}(M\otimes_{O_\lambda}E_\lambda)^\Gamma\geq\dim_{O_\lambda/\lambda}(M\otimes_{O_\lambda}O_\lambda/\lambda)^\Gamma$. Then $M$ is weakly semisimple as well, and $\dim_{E_\lambda}(M\otimes_{O_\lambda}E_\lambda)^\Gamma=\dim_{O_\lambda/\lambda}(M\otimes_{O_\lambda}O_\lambda/\lambda)^\Gamma$.
\end{lem}

\begin{proof}
Since $M$ is a finitely generated free $O_\lambda$-module, both $M^\Gamma$ and $M/M^\Gamma$ are finitely generated free $O_\lambda$-modules. In particular, the map $M^\Gamma\otimes_{O_\lambda}O_\lambda/\lambda\to(M\otimes_{O_\lambda}O_\lambda/\lambda)^\Gamma$ is injective. As we have
\[
\dim_{O_\lambda/\lambda}M^\Gamma\otimes_{O_\lambda}O_\lambda/\lambda
=\rank_{O_\lambda}M^\Gamma=\dim_{E_\lambda}(M\otimes_{O_\lambda}E_\lambda)^\Gamma,
\]
the map $M^\Gamma\otimes_{O_\lambda}O_\lambda/\lambda\to(M\otimes_{O_\lambda}O_\lambda/\lambda)^\Gamma$ is an isomorphism. It follows that
\[
\dim_{E_\lambda}(M\otimes_{O_\lambda}E_\lambda)^\Gamma=\dim_{O_\lambda/\lambda}(M\otimes_{O_\lambda}O_\lambda/\lambda)^\Gamma.
\]
It also follows that the maps
\[
M^\Gamma\otimes_{O_\lambda}O_\lambda/\lambda\to(M\otimes_{O_\lambda}O_\lambda/\lambda)^\Gamma
\to(M\otimes_{O_\lambda}O_\lambda/\lambda)_\Gamma\simeq M_\Gamma\otimes_{O_\lambda}O_\lambda/\lambda
\]
are isomorphisms since $M\otimes_{O_\lambda}O_\lambda/\lambda$ is weakly semisimple. By Nakayama's lemma, the map $M^\Gamma\to M_\Gamma$ is surjective. By Lemma \ref{le:weakly_semisimple0}(2), $M$ is weakly semisimple. The lemma is proved.
\end{proof}

To end this subsection, we record the following definition which slightly generalizes \cite{Liu1}*{Definition~5.1}, and will be used in later sections.

\begin{definition}\label{de:divisibility}
Consider an $O_\lambda$-module $M$ and an element $x\in M$. We define the \emph{exponent} and the \emph{order} of $x$ to be
\begin{align*}
\exp_\lambda(x,M)&\coloneqq\min\{d\in\dZ_{\geq 0}\cup\{\infty\}\res \lambda^d x=0\},\\
\ord_\lambda(x,M)&\coloneqq\sup\{d\in\dZ_{\geq 0}\res x\in\lambda^d M\},
\end{align*}
respectively.
\end{definition}

\subsection{Local Galois cohomology}
\label{ss:local_galois}

In this subsection, we study Galois cohomology locally at nonarchimedean places of $F$. Let $w$ be a nonarchimedean place of $F$. We recall from \S\ref{ss:notation} various notations concerning $F_w$.

\begin{notation}\label{no:galois_module}
For a $\dZ_\ell$-ring $L$ that is finite over either $\dZ_\ell$ or $\dQ_\ell$ and $?\in\{\;,\tor,\free\}$, we
\begin{enumerate}
  \item put $\Mod(F_w,L)_?\coloneqq\Mod(\Gamma_{F_w},L)_?$;

  \item denote by $\obj(j)\colon\Mod(F_w,L)_?\to\Mod(F_w,L)_?$ the functor of $j$-th Tate twist for $j\in\dZ$; and

  \item denote by $\obj^\vee\colon\Mod(F_w,L)_?^\op\to\Mod(F_w,L)_?$ the functor sending $M$ to $\Hom_L(M,L)$.
\end{enumerate}
We also denote
\[
\obj_\dQ\colon\Mod(F_w,O_\lambda)\to\Mod(F_w,E_\lambda)
\]
the base change functor sending $M$ to $M\otimes_{O_\lambda}E_\lambda$, and
\[
\obj^*\colon\Mod(F_w,O_\lambda)_\tor^\op\to\Mod(F_w,O_\lambda)
\]
the $E_\lambda$-Pontryagin duality functor sending $M$ to $\Hom_{O_\lambda}(M,E_\lambda/O_\lambda)$. For every pair $m,m'\in\{1,2,\dots,\infty\}$ with $m'\geq m$, we have a ``reduction modulo $\lambda^m$'' functor
\[
\bar\obj^{(m)}\coloneqq\obj\otimes_{O_\lambda}O_\lambda/\lambda^m\colon\Mod(F_w,O_\lambda/\lambda^{m'})\to\Mod(F_w,O_\lambda/\lambda^m)
\]
(that is, it sends $\rR$ to $\bar\rR^{(m)}$).\footnote{Here, $O_\lambda/\lambda^\infty$ is understood as $O_\lambda$.} We usually write $\bar\obj$ for $\bar\obj^{(1)}$.
\end{notation}

For every object $\rR\in\Mod(F_w,O_\lambda)$, we have the local Tate pairing
\begin{align}\label{eq:local_tate}
\langle\;,\;\rangle_w\colon \rH^1(F_w,\rR)\times \rH^1(F_w,\rR^*(1))\xrightarrow{\cup}\rH^2(F_w,E_\lambda/O_\lambda(1))\simeq E_\lambda/O_\lambda,
\end{align}
which we will study in the following. We will define a submodule functor $\rH^1_\ns(F_w,\obj)$ of $\rH^1(F_w,\obj)$ for every nonarchimedean place $w$ of $F$, which is usually denoted as $\rH^1_\ur(F_w,\obj)$ and $\rH^1_f(F_w,\obj)$ when $\ell\nmid w$ and $\ell\mid w$, respectively. We choose this unconventional notation only to uniformize the two cases.

First, we study the case where $\ell\nmid w$.

\begin{definition}\label{de:local_finite}
For every object $\rR$ in either $\Mod(F_w,E_\lambda)$ or $\Mod(F_w,O_\lambda)$, we put
\[
\rH^1_\sing(F_w,\rR)\coloneqq\rH^1(\rI_{F_w},\rR)^{\Gamma_{\kappa_w}};
\]
and denote by $\rH^1_\ns(F_w,\rR)$ the kernel of the canonical map
\[
\partial_w\colon\rH^1(F_w,\rR)\to\rH^1_\sing(F_w,\rR).
\]
\end{definition}

By the inflation-restriction exact sequence (see, for example, \cite{Liu2}*{Lemma~2.6}), we know that $\partial_w$ is surjective, and that $\rH^1_\ns(F_w,\rR)$ is canonically isomorphic to $\rH^1(\kappa_w,\rR^{\rI_{F_w}})$.

\begin{lem}\label{le:tate}
For $\rR\in\Mod(F_w,O_\lambda)_\tor$, the restriction of the local Tate pairing $\langle\;,\;\rangle_w$ \eqref{eq:local_tate} to $\rH^1_\ns(F_w,\rR)\times\rH^1_\ns(F_w,\rR^*(1))$ vanishes.
\end{lem}

\begin{proof}
This is well-known. In fact, the cup product of $\rH^1_\ns(F_w,\rR)$ and $\rH^1_\ns(F_w,\rR^*(1))$ factors through $\rH^2(\kappa_w,\rR^{\rI_{F_w}}\otimes\rR^*(1)^{\rI_{F_w}})$, which is the zero group.
\end{proof}

Second, we study the case where $\ell\mid w$. In particular, $F_w$ is a finite unramified extension of $\dQ_\ell$. Denote by $\obj_0\colon\Mod(F_w,O_\lambda)\to\Mod(F_w,\dZ_\ell)$ the obvious forgetful functor.

\begin{definition}\label{de:crystalline}
Let $a\leq b$ be two integers.
\begin{enumerate}
  \item For an object $\rR\in\Mod(F_w,\dZ_\ell)_\tor$, we say that $\rR$ is \emph{crystalline (with Hodge--Tate weights in $[a,b]$)} if $\rR=\rR''/\rR'$ where $\rR'\subseteq\rR''$ are two $\Gamma_{F_w}$-stable $\dZ_\ell$-lattices in a crystalline $\dQ_\ell$-representation of $\Gamma_{F_w}$ (with Hodge--Tate weights in $[a,b]$).\footnote{We adopt the convention that $\dQ_\ell(1)$ has Hodge--Tate weight $-1$.}

  \item For an object $\rR\in\Mod(F_w,\dZ_\ell)$, we say that $\rR$ is \emph{crystalline (with Hodge--Tate weights in $[a,b]$)} if  $\rR/\ell^m\rR$ is a torsion crystalline module (with Hodge--Tate weights in $[a,b]$) for every integer $m\geq 1$.\footnote{In fact, by Lemma \ref{le:crystalline} below, when $a\leq 0\leq b$ and $b-a\leq\ell-2$, an object $\rR\in\Mod(F_w,\dZ_\ell)_\free$ is crystalline with Hodge--Tate weights in $[a,b]$ if and only if $\rR_\dQ$ is.}

  \item For an object $\rR\in\Mod(F_w,O_\lambda)$, we say that $\rR$ is \emph{crystalline (with Hodge--Tate weights in $[a,b]$)} if $\rR_0$ is.
\end{enumerate}
\end{definition}

\begin{definition}[\cite{Niz93}*{\S4}]\label{de:crystalline_finite}
For an object $\rR\in\Mod(F_w,O_\lambda)$ that is crystalline, we define $\rH^1_\ns(F_w,\rR)$ to be the subset of $\rH^1(F_w,\rR)=\rH^1(F_w,\rR_0)$ consisting of elements $s$ represented by an extension
\[
0\to \rR_0 \to \rR_s \to \dZ_\ell \to 0
\]
in the category $\Mod(F_w,\dZ_\ell)$ such that $\rR_s$ is crystalline.\footnote{It is clear that if $\rR$ is crystalline with Hodge--Tate weights in $[a,b]$ for $a\leq 0\leq b$, then $\rR_s$ in the extension representing an element in $\rH^1_\ns(F_w,\rR)$ is also crystalline with Hodge--Tate weights in $[a,b]$.}
\end{definition}

It follows that $\rH^1_\ns(F_w,\rR)$ is an $O_\lambda$-submodule of $\rH^1(F_w,\rR)$.

\begin{lem}\label{le:crystalline}
Let $\rR$ be an object of $\Mod(F_w,O_\lambda)_\free$ such that $\rR_\dQ$ is crystalline with Hodge--Tate weights in $[a,b]$ with $a\leq 0\leq b$ and $b-a\leq\ell-2$. Then $\rH^1_\ns(F_w,\rR)$ coincides with the preimage of
\[
\Ker\(\rH^1(F_w,\rR_\dQ)\to\rH^1(F_w,\rR_\dQ\otimes_{\dQ_\ell}\dB_\cris)\)
\]
under the natural map $\rH^1(F_w,\rR)\to\rH^1(F_w,\rR_\dQ)$.
\end{lem}

\begin{proof}
This is proved in \cite{Bre99}*{Proposition~6}.
\end{proof}

\begin{lem}\label{le:crystalline_tate}
Suppose that the integers $a,b$ satisfy $a<0\leq b$ and $b-a\leq \tfrac{\ell-2}{2}$. Then for every $\rR\in\Mod(F_w,O_\lambda)_\tor$ that is crystalline with Hodge--Tate weights in $[a,b]$, the restriction of the local Tate pairing $\langle\;,\;\rangle_w$ \eqref{eq:local_tate} to $\rH^1_\ns(F_w,\rR)\times\rH^1_\ns(F_w,\rR^*(1))$ takes values in $\fd_\lambda^{-1}/O_\lambda$, where $\fd_\lambda\subseteq O_\lambda$ is the different ideal of $E_\lambda$ over $\dQ_\ell$.
\end{lem}

\begin{proof}
We have a canonical map $\Tr\colon(\rR^*)_0\to(\rR_0)^*$ in the category $\Mod(F_w,\dZ_\ell)$ induced by the trace map $\Tr_{E_\lambda/\dQ_\ell}$, which induces a map $\rH^1(F_w,\rR^*(1))\to\rH^1(F_w,(\rR_0)^*(1))$ under which the image of $\rH^1_\ns(F_w,\rR^*(1))$ is contained in $\rH^1_\ns(F_w,(\rR_0)^*(1))$. Take arbitrary elements $x\in\rH^1_\ns(F_w,\rR)$ and $y\in\rH^1_\ns(F_w,\rR^*(1))$. Then we have
\[
\Tr_{E_\lambda/\dQ_\ell}(\langle x,y\rangle_w)=\Tr_{E_\lambda/\dQ_\ell}\langle x,y\rangle_w=\langle x,\Tr(y)\rangle_w\in\dQ_\ell/\dZ_\ell.
\]
However, $\langle x,\Tr(y)\rangle_w=0$ by \cite{Niz93}*{Proposition~6.2}. The lemma follows.
\end{proof}

\subsection{Some Galois-theoretical lemmas}
\label{ss:galois_lemma}

In this subsection, we generalize some lemmas from \cite{Liu1}. For a finite set $\Sigma$ of places of $F$, we denote by $\Gamma_{F,\Sigma}$ the Galois group of the maximal subextension of $\ol{F}/F$ that is unramified outside $\Sigma$.

\begin{notation}
For a $\dZ_\ell$-ring $L$ that is finite over either $\dZ_\ell$ or $\dQ_\ell$ and $?\in\{\;,\tor,\free\}$, we put
\[
\Mod(F,L)_?\coloneqq\varinjlim_\Sigma\Mod(\Gamma_{F,\Sigma},L)_?,
\]
where the colimit is taken over all finite sets $\Sigma$ of places of $F$ with inflation as transition functors. We have functors $\obj(j)$, $\obj^\vee$, $\obj_\dQ$, $\obj^*$, and $\bar\obj^{(m)}$ similar to those in Notation \ref{no:galois_module}. For an object $\rR\in\Mod(F,L)$ and $i\in\dZ$, we put
\[
\rH^i(F,\rR)\coloneqq\varinjlim_\Sigma\rH^i(\Gamma_{F,\Sigma},\rR).
\]
Moreover, for every place $w$ of $F$, we have the restriction functor $\Mod(F,L)\to\Mod(F_w,L)$; and denote
\[
\loc_w\colon \rH^i(F,\rR)\to\rH^i(F_w,\rR)
\]
the localization map.
\end{notation}

\begin{definition}[\cite{Liu1}*{Definition~5.1}]
Let $G$ be a profinite group. For an object $\rR\in\Mod(G,O_\lambda)_\tor$, we define its \emph{reducibility depth} to be the smallest integer $\fr_\rR\geq 0$ such that
\begin{enumerate}
  \item if $\rR'$ is a $G$-stable $O_\lambda$-submodule that is not contained in $\lambda\rR$, then $\rR'$ contains $\lambda^{\fr_\rR}\rR$;

  \item for every positive integer $m$, the group $\End_{O_\lambda[G]}(\bar\rR^{(m)})/O_\lambda\cdot\id$ is annihilated by $\lambda^{\fr_\rR}$.
\end{enumerate}
\end{definition}

Note that if $\rR/\lambda\rR$ is absolutely irreducible, then $\fr_\rR=0$.

\begin{lem}\label{le:reducibility}
Let $\rR\in\Mod(F,O_\lambda)$ be an object such that $\rR_\dQ$ is absolutely irreducible. Then there exists an integer $\fr_\rR$ depending on $\rR$ only, such that $\bar\rR^{(m)}$ has reducibility depth at most $\fr_\rR$ for every positive integer $m$.
\end{lem}

\begin{proof}
The same argument in \cite{Liu1}*{Lemma~5.2} applies to our case as well, with $\dZ/p^n$ replaced by $O_\lambda/\lambda^m$.
\end{proof}

Now we fix a positive integer $m$. Consider an object $\rR\in\Mod(F,O_\lambda/\lambda^m)_\free$. We denote by $\rho\colon\Gamma_F\to\GL(\rR)$ the associated homomorphism. Let $F_\rho/F$ be the Galois extension fixed by the kernel of $\rho$, and $G\coloneqq\Gal(F_\rho/F)$ the image of $\rho$. we have the restriction map
\begin{align}\label{eq:restriction}
\Res_\rho\colon\rH^1(F,\rR)\to\rH^1(F_\rho,\rR)^G=\Hom_G(\Gamma_{F_\rho}^\ab,\rR),
\end{align}
where $\Gamma_{F_\rho}^\ab\coloneqq\Gal(F_\rho^\ab/F_\rho)$ with $F_\rho^\ab\subseteq\ol{F}$ the maximal abelian extension of $F_\rho$, which is equipped with the natural conjugation action by $G=\Gal(F_\rho/F)$.

The map $\Res_\rho$ induces an $O_\lambda$-linear pairing
\[
[\;,\;]\colon\rH^1(F,\rR)\times\Gamma_{F_\rho}^\ab\to \rR,
\]
such that the action of $G$ on $\Gamma_{F_\rho}^\ab$ is compatible with that on $\rR$. Let $S$ be a finitely generated $O_\lambda/\lambda^m$-submodule of $\rH^1(F,\rR)$, and let $F_S/F_\rho$ be the finite abelian extension such that $\Gal(F_\rho^\ab/F_S)$ is the subgroup of $\Gamma_{F_\rho}^\ab$ consisting of $\gamma$ satisfying $[s,\gamma]=0$ for every $s\in S$. Then the above pairing induces an injective map
\begin{align}\label{eq:theta}
\theta_S\colon\Gal(F_S/F_\rho)\to\Hom_{O_\lambda}(S,\rR)
\end{align}
of abelian groups that is compatible with $G$-actions.

As in \cite{Liu1}*{\S5.1}, we introduce a sequence $\ff$ that is given by $\ff(0)=1$, $\ff(1)=1$, $\ff(2)=4$, $\ff(r+1)=2(\ff(r)+1)$ for $r\geq 2$.

\begin{lem}\label{le:image}
Let the notation be as above. Suppose that the map $\Res_\rho$ is injective. If $S$ is a free $O_\lambda/\lambda^m$-module of rank $r_S$ for some positive integer $m$, then the $O_\lambda$-submodule of $\Hom_{O_\lambda}(S,\rR)$ generated by the image of $\theta_S$ contains $\lambda^{\ff(r_S)\fr_\rR}\Hom_{O_\lambda}(S,\rR)$, where $\fr_\rR$ is the reducibility depth of $\rR$.
\end{lem}

\begin{proof}
The same argument in \cite{Liu1}*{Lemma~5.4} applies to our case as well, with $\dZ/p^n$ replaced by $O_\lambda/\lambda^mp$. Note that the proof only uses the injectivity, not the surjectivity, of the map $\Res_\rho$.
\end{proof}

Concerning the injectivity of the map $\Res_\rho$ \eqref{eq:restriction}, we have the following lemma.

\begin{lem}\label{le:restriction}
Suppose that either one of the following two assumptions holds:
\begin{enumerate}[label=(\alph*)]
  \item the image of $\Gamma_F$ in $\GL(\bar\rR)$ contains a nontrivial scalar element;

  \item $\dim_{O_\lambda/\lambda}\bar\rR\leq\min\{\tfrac{\ell+1}{2},\ell-3\}$, $\bar\rR$ is a semisimple $(O_\lambda/\lambda)[\Gamma_F]$-module, and moreover $\Hom_{(O_\lambda/\lambda)[\Gamma_F]}(\End(\bar\rR),\bar\rR)=0$.
\end{enumerate}
Then $\Res_\rho$ is injective.
\end{lem}

\begin{proof}
By the inflation-restriction exact sequence, it suffices to show that $\rH^1(G,\rR)=0$.

In the situation (a), it follows that $G$ contains a nontrivial scalar element of order coprime to $\ell$. Then by the same argument in \cite{Gro89}*{Proposition~9.1}, we have $\rH^1(G,\rR)=0$. More precisely, let $\gamma\in G$ be a nontrivial scalar element of order coprime to $\ell$. Then we have $\rH^1(G/\langle\gamma\rangle,\rR^\gamma)=0$ and $\rH^1(\langle\gamma\rangle,\rR)=0$, which imply $\rH^1(G,\rR)=0$.

Now we consider the situation (b). We prove by induction that $\rH^1(G,\bar\rR^{(i)})=0$ for $1\leq i\leq m$. Suppose that $\rH^1(G,\bar\rR^{(j)})=0$ for $1\leq j\leq i<m$. By the short exact sequence
\[
0 \to \bar\rR^{(i+1)}\otimes_{O_\lambda/\lambda^{i+1}}\lambda^i/\lambda^{i+1} \to \bar\rR^{(i+1)} \to \bar\rR^{(i)} \to 0
\]
of $O_\lambda[G]$-modules, in which $\bar\rR^{(i+1)}\otimes_{O_\lambda/\lambda^{i+1}}\lambda^i/\lambda^{i+1}$ is isomorphic to $\bar\rR$, we know that $\rH^1(G,\bar\rR^{(i+1)})=0$. Therefore, it remains to check the initial step that $\rH^1(G,\bar\rR)=0$.

Let $G^i\subseteq G$ be the kernel of the composite homomorphism $G\to\GL(\rR)\to\GL(\bar\rR^{(i)})$ for $1\leq i\leq m$, so we obtain a filtration $0=G^m\subseteq G^{m-1}\subseteq\cdots\subseteq G^1\subseteq G$ of normal subgroups of $G$. We prove by induction that $\rH^1(G/G^i,\bar\rR)=0$. For $i=1$, since $\bar\rR$ is a faithful semisimple $(O_\lambda/\lambda)[G/G^1]$-module, $G/G^1$ has no nontrivial normal $\ell$-subgroup. As $\dim_{O_\lambda/\lambda}\bar\rR\leq\ell-3$, we have $\rH^1(G/G^1,\bar\rR)=0$ by \cite{Gur99}*{Theorem~A}. Suppose that $\rH^1(G/G^j,\bar\rR)=0$ for $1\leq j\leq i<m$. By the inflation-restriction exact sequence
\[
0 \to \rH^1(G/G^i,\bar\rR) \to \rH^1(G/G^{i+1},\bar\rR) \to \Hom_G(G^i/G^{i+1},\bar\rR),
\]
it suffices to show that $\Hom_G(G^i/G^{i+1},\bar\rR)=0$, or equivalently, $\Hom_{(O_\lambda/\lambda)[G]}(G^i/G^{i+1}\otimes O_\lambda/\lambda,\bar\rR)=0$. Note that $G^i/G^{i+1}$ is an $\dF_\ell[G]$-submodule of $\End(\bar\rR)$, hence $(G^i/G^{i+1})\otimes O_\lambda/\lambda$ is an $(O_\lambda/\lambda)[G]$-submodule of $\End(\bar\rR)\otimes(O_\lambda/\lambda)\simeq\End(\bar\rR)^d$, where $d\coloneqq[O_\lambda/\lambda:\dF_\ell]$ is the degree. Since $\bar\rR$ is a semisimple $(O_\lambda/\lambda)[G]$-module and $2\dim_{O_\lambda/\lambda}\bar\rR<\ell+2$, by \cite{Ser94}*{Corollaire~1}, we know that $\End(\bar\rR)$ is a semisimple $(O_\lambda/\lambda)[G]$-module. In particular, we have $\Hom_{(O_\lambda/\lambda)[G]}(G^i/G^{i+1}\otimes O_\lambda/\lambda,\bar\rR)=0$ as $\Hom_G(\End(\bar\rR),\bar\rR)=0$.

The lemma is proved.
\end{proof}

\subsection{Reduction of Selmer groups}
\label{ss:bloch_kato}

We recall the following definition of the Bloch--Kato Selmer group from \cite{BK90}.

\begin{definition}[Bloch--Kato Selmer group]\label{de:bk_rational}
For an object $\rR\in\Mod(F,E_\lambda)$, we define the \emph{Bloch--Kato Selmer group} $\rH^1_f(F,\rR)$ of $\rR$ to be the $E_\lambda$-subspace of $\rH^1(F,\rR)$ consisting of elements $s$ such that
\begin{enumerate}
  \item $\loc_w(s)\in\rH^1_\ns(F_w,\rR)$ (Definition \ref{de:local_finite}) for every nonarchimedean place $w$ of $F$ not above $\ell$; and

  \item $\loc_w(s)\in\Ker\(\rH^1(F_w,\rR)\to\rH^1(F_w,\rR\otimes_{\dQ_\ell}\dB_\cris)\)$ for every place $w$ of $F$ above $\ell$.
\end{enumerate}
\end{definition}

\begin{definition}\label{de:bk_integral}
Consider an object $\rR\in\Mod(F,O_\lambda)_\free$.
\begin{enumerate}
  \item We define the \emph{(integral) Bloch--Kato Selmer group} $\rH^1_f(F,\rR)$ of $\rR$ to be inverse image of $\rH^1_f(F,\rR_\dQ)$ under the obvious map $\rH^1(F,\rR)\to\rH^1(F,\rR_\dQ)$.

  \item For $m\in\{1,2,\dots,\infty\}$, we define $\rH^1_{f,\rR}(F,\bar\rR^{(m)})$ to be the image of $\rH^1_f(F,\rR)$ under the obvious map $\rH^1(F,\rR)\to\rH^1(F,\bar\rR^{(m)})$.
\end{enumerate}
\end{definition}

\begin{lem}\label{le:global_local}
Consider an object $\rR\in\Mod(F,O_\lambda)_\free$. Suppose that we are in one of the two following cases
\begin{enumerate}
  \item $w$ is a nonarchimedean place of $F$ not above $\ell$ at which $\rR$ is unramified.

  \item $w$ is a place of $F$ above $\ell$ at which $\rR_\dQ$ is crystalline with Hodge--Tate weights in $[a,b]$ with $a\leq 0\leq b$ and $b-a\leq\ell-2$.
\end{enumerate}
Then for every positive integer $m$, the image of $\rH^1_{f,\rR}(F,\bar\rR^{(m)})$ under the localization map $\loc_w\colon\rH^1(F,\bar\rR^{(m)})\to\rH^1(F_w,\bar\rR^{(m)})$ is contained in $\rH^1_\ns(F_w,\bar\rR^{(m)})$.
\end{lem}

\begin{proof}
Case (1) follows from \cite{Rub00}*{Lemma~1.3.5 \& Lemma~1.3.8}. Case (2) follows from Lemma \ref{le:crystalline}.
\end{proof}

We recall the notion of purity for a local Galois representation.

\begin{definition}\label{de:purity}
Let $w$ be a nonarchimedean place of $F$ not above $\ell$. Consider an object $\rR\in\Mod(F_w,E_\lambda)$. Let $\WD(\rR)$ be the attached Weil--Deligne representation, and $\gr_n\WD(\rR)$ the $n$-th graded piece of the monodromy filtration on $\WD(\rR)$. For $\mu\in\dZ$, we say that $\rR$ is \emph{pure of weight $\mu$} if $\gr_n\WD(\rR)$ is strictly pure of weight $\mu+n$ for each $n$, that is, all eigenvalues of $\phi_w$ on $\gr_n\WD(\rR)$ are Weil $\|w\|^{-(\mu+n)}$-numbers.\footnote{In particular, $E_\lambda(1)$ is (strictly) pure of weight $-2$.}
\end{definition}

From now to the end of this section, we suppose that the complex conjugation $\tc$ restricts to an automorphism of $F$ (of order at most two). We adopt the notation concerning ground fields in \S\ref{ss:notation}; in particular, we put $F^+\coloneqq F^{\tc=1}$. We also have a functor
\[
\obj^\tc\colon\Mod(F,L)\to\Mod(F,L)
\]
induced by the conjugation by $\tc$.

\begin{lem}\label{le:conjugate}
For every object $\rR\in\Mod(F,E_\lambda)$, the functor $\obj^\tc$ induces an isomorphism
\[
\rH^1_f(F,\rR)\simeq\rH^1_f(F,\rR^\tc)
\]
of Selmer groups.
\end{lem}

\begin{proof}
Regard elements in $\rH^1(F,\obj)$ as extensions. Then applying $\obj^\tc$ to extensions induces maps
\[
\rH^1(F,\rR)\to\rH^1(F,\rR^\tc),\quad
\rH^1(F,\rR^\tc)\to\rH^1(F,\rR)
\]
which are inverses to each other. It is clear that conditions (1) and (2) in Definition \ref{de:bk_rational} are preserved under such maps. The lemma follows.
\end{proof}

\begin{proposition}\label{pr:selmer_reduction}
Let $\rR$ be an object in $\Mod(F,O_\lambda)_\free$.
\begin{enumerate}
  \item Let $S$ be a free $O_\lambda$-submodule of $\rH^1_f(F,\rR)$ whose image in $\rH^1_f(F,\rR)/\rH^1_f(F,\rR)_\tor$ is saturated. For every positive integer $m$, if we denote by $S^{(m)}$ the image of $S$ in $\rH^1_{f,\rR}(F,\bar\rR^{(m)})$, then it is a free $O_\lambda/\lambda^m$-module of the same rank as $S$.

  \item Suppose that $\rR$ satisfies $\rR_\dQ^\tc\simeq\rR_\dQ^\vee(1)$ and such that $\rR_\dQ$ is pure of weight $-1$ at every nonarchimedean place $w$ of $F$ not above $\ell$. For every finite set $\Sigma$ of places of $F$, there exists a positive integer $m_\Sigma$, depending on $\rR$ and $\Sigma$, such that for every $S$ as in (1) and every integer $m>m_\Sigma$, we have $\loc_w(\lambda^{m_\Sigma}S^{(m)})=0$ for every nonarchimedean place $w\in\Sigma$ not above $\ell$.
\end{enumerate}
\end{proposition}

\begin{proof}
For (1), let $\rT$ be the image of $\rH^1_f(F,\rR)_\tor$ in $\rH^1(F,\bar\rR^{(m)})$, which is contained in $\rH^1_{f,\rR}(F,\bar\rR^{(m)})$. Then we have a natural injective map
\[
\frac{\rH^1_f(F,\rR)/\rH^1_f(F,\rR)_\tor}{\lambda^m(\rH^1_f(F,\rR)/\rH^1_f(F,\rR)_\tor)}\to\rH^1_{f,\rR}(F,\bar\rR^{(m)})/\rT.
\]
Since the image of $S$ in $\rH^1_f(F,\rR)/\rH^1_f(F,\rR)_\tor$ is saturated, (1) follows immediately.

For (2), we look at the map
\[
\loc_\Sigma^{\infty\ell}\colon\rH^1_{f,\rR}(F,\bar\rR^{(m)})\to\bigoplus_{w\in\Sigma,w\nmid\infty\ell}\rH^1(F_w,\bar\rR^{(m)}).
\]
For every $w\nmid\infty\ell$, since $\rR_\dQ$ is of pure weight $-1$ at $w$, $\rR_\dQ^\tc$ and $\rR_\dQ^\vee(1)$ are of pure weight of $-1$ at $w$ as well. Thus, we have $\rH^0(F_w,\rR_\dQ)=0$ and $\rH^2(F_w,\rR_\dQ)\simeq\rH^0(F_w,\rR_\dQ^\vee(1))^\vee=0$, hence $\rH^1(F_w,\rR_\dQ)=0$ by the Euler characteristic formula (see also the proof of \cite{Nek07}*{Proposition~4.2.2(1)}). Thus, $\rH^1(F_w,\rR)$ is annihilated by $\lambda^{m_w}$ for some integer $m_w\geq 0$. We may enlarge $m_w$ such that $\lambda^{m_w}$ also annihilates $\rH^2(F_w,\rR)_\tor$. Then it follows that $\rH^1(F_w,\bar\rR^{(m)})$ is annihilated by $\lambda^{2m_w}$. Now if we put $m_\Sigma\coloneqq\max\{2m_w\res w\in\Sigma,w\nmid\infty\ell\}$, then (2) follows. This completes the proof of the proposition.
\end{proof}

\subsection{Extension of essentially conjugate self-dual representations}
\label{ss:conjugate_selfdual}

In this subsection, we collect some notion and facts on the extension of essentially conjugate self-dual representations.

\begin{notation}\label{no:sg}
When $[F:F^+]=2$, we introduce the group scheme $\sG_N$ from \cite{CHT08}*{\S1} as
\[
\sG_N\coloneqq(\GL_N\times\GL_1)\rtimes\{1,\fc\}
\]
with $\fc^2=1$ and
\[
\fc(g,\mu)\fc=(\mu\tp{g}^{-1},\mu)
\]
for $(g,\mu)\in\GL_N\times\GL_1$. Denote by $\nu\colon\sG_N\to\GL_1$ the homomorphism such that $\nu\res_{\GL_N\times\GL_1}$ is the projection to the factor $\GL_1$ and that $\nu(\fc)=-1$.

When $[F:F^+]=1$, we put $\sG_N\coloneqq\GL_N\times\GL_1$ and regard the symbol $\fc$ as the identity element. Denote by $\nu\colon\sG_N\to\GL_1$ the projection to the second factor.
\end{notation}

\begin{notation}\label{no:sg_extension}
Let $R$ be a topological ring. For a continuous homomorphism
\[
r\colon\Gamma_{F^+}\to\sG_N(R)
\]
such that the image of $r\res_{\Gamma_F}$ lies in $\GL_N(R)\times R^\times$, we denote
\[
r^\natural\colon\Gamma_F\to\GL_N(R)\times R^\times\to\GL_N(R)
\]
the composition of $r\res_{\Gamma_F}$ with the projection to $\GL_N(R)$.
\end{notation}

To end this subsection, we recall the extension along $j$-polarization. This has been introduced in \cite{CHT08}*{\S1} when $[F:F^+]=2$.

\begin{definition}\label{de:polarization}
For a $\dZ_\ell$-ring $L$ that is finite over either $\dZ_\ell$ or $\dQ_\ell$, an integer $j$, and an object $\rR$ in $\Mod(F,L)$, a \emph{$j$-polarization} of $\rR$ is an isomorphism
\[
\Xi\colon\rR^\tc\xrightarrow{\sim}\rR^\vee(j)
\]
in $\Mod(F,L)$, such that $\Xi^{\tc,\vee}(j)=(-1)^{\mu_\Xi+j+1}\cdot\Xi$ for some $\mu_\Xi\in\dZ/2\dZ$. We say that $\rR$ is \emph{$j$-polarizable} if there exists a $j$-polarization.
\end{definition}

\begin{construction}\label{cs:extension}
Let $\rR$ be a nonzero object in $\Mod(F,L)_\free$ with the associated continuous homomorphism $\rho\colon\Gamma_F\to\GL(\rR)$, equipped with a $j$-polarization $\Xi\colon\rR^\tc\xrightarrow{\sim}\rR^\vee(j)$. Choose an isomorphism $\rR\simeq L^{\oplus N}$ of the underlying $L$-modules for a unique integer $N\geq 1$.
\begin{enumerate}
  \item When $[F:F^+]=1$, we let
     \[
     \rho_+\colon\Gamma_{F^+}\to\sG_N(L)
     \]
     be the continuous homomorphism sending $g\in\Gamma_{F^+}=\Gamma_F$ to $(\rho(g),\epsilon_\ell^j(g))$.

  \item When $[F:F^+]=2$, the $j$-polarization $\Xi$ gives rise to an element $B\in\GL_N(L)$ satisfying $\rho^\tc=B\circ\epsilon_\ell^j\rho^\vee\circ B^{-1}$ and $B\tp{B}^{-1}=(-1)^{\mu_\Xi+j+1}$. We let
      \[
      \rho_+\colon\Gamma_{F^+}\to\sG_N(L)
      \]
      be the continuous homomorphism given by the formula $\rho_+\res_{\Gamma_F}=(\rho,\epsilon_\ell^j\res_{\Gamma_F})1$ and $\rho_+(\tc)=(B,(-1)^{\mu_\Xi+j+1})\fc$.
\end{enumerate}
In both cases, we call $\rho_+$ an \emph{extension} of $\rho$.
\end{construction}

\subsection{Localization of Selmer groups}
\label{ss:localization}

In this subsection, we study the behavior of Selmer groups under localization maps.

\begin{notation}\label{no:reduction}
We take a nonzero object $\rR\in\Mod(F,O_\lambda)_\free$ with the associated homomorphism $\rho\colon\Gamma_F\to\GL(\rR)$, together with a $j$-polarization $\Xi\colon\rR^\tc\xrightarrow{\sim}\rR^\vee(j)$. We fix an isomorphism $\rR\simeq O_\lambda^{\oplus N}$. Let
\[
\rho_+\colon\Gamma_{F^+}\to\sG_N(O_\lambda)
\]
be the extension of $\rho$ from Construction \ref{cs:extension}. For every integer $m\geq 1$, we have the induced homomorphisms
\begin{align*}
\bar\rho^{(m)}&\colon\Gamma_F\to\GL(\bar\rR^{(m)})\simeq\GL_N(O_\lambda/\lambda^m),\\
\bar\rho^{(m)}_+&\colon\Gamma_{F^+}\to\sG_N(O_\lambda/\lambda^m),
\end{align*}
and we omit the superscript $(m)$ when $m=1$.
\end{notation}

We denote by $F^{(m)}\coloneqq F_{\bar\rho^{(m)}}$ and $F^{(m)}_+$ the subfields of $\ol{F}$ fixed by $\Ker\bar\rho^{(m)}$ and $\Ker\bar\rho^{(m)}_+$, respectively. In particular, we have $F\subseteq F^{(m)}\subseteq F^{(m)}_+\subseteq F^{(m)}(\zeta_{\ell^m})$.

\begin{notation}\label{no:h_gamma}
For a positive integer $m$ and an element
\[
\gamma\in(\GL_N(O_\lambda/\lambda^m)\times(O_\lambda/\lambda^m)^\times)\fc\subseteq\sG_N(O_\lambda/\lambda^m),
\]
we denote by $h_\gamma\in\GL_N(O_\lambda/\lambda^m)$ the first component of $\gamma^{[F:F^+]}\in\GL_N(O_\lambda/\lambda^m)\times(O_\lambda/\lambda^m)^\times$.
\end{notation}

Now we fix a positive integer $m$ and a finitely generated $O_\lambda$-submodule $S$ of $\rH^1_{f,\rR}(F,\bar\rR^{(m)})$. We have the finite abelian extension $F_S/F^{(m)}$ from \S\ref{ss:galois_lemma}. Consider an element $\gamma$ as in Notation \ref{no:h_gamma} that belongs to the image of $\bar\rho^{(m)}_+$. The following definition is essentially \cite{Liu1}*{Definition~5.6}.

\begin{definition}\label{de:associated}
We say that a place $w^{(m)}_+$ of $F^{(m)}_+$ is \emph{$\gamma$-associated} if
\begin{itemize}[label={\ding{109}}]
  \item $w^{(m)}_+$ is not above $\infty$ or $\ell$;

  \item $w^{(m)}_+$ is unramified over $F^+$;

  \item its underlying place of $F^{(m)}$ is unramified in $F_S$; and

  \item its arithmetic Frobenius substitution in $\Gal(F^{(m)}_+/F^+)\simeq\IM\bar\rho^{(m)}_+$ coincides with $\gamma$.
\end{itemize}
\end{definition}

Recall the injective map
\[
\theta_S\colon\Gal(F_S/F^{(m)})\to\Hom_{O_\lambda}(S,\bar\rR^{(m)})
\]
of abelian groups from \eqref{eq:theta} with $\rho=\bar\rho^{(m)}$, which is equivariant under the action of $\Gal(F^{(m)}/F)$. Take a $\gamma$-associated place $w^{(m)}_+$ of $F^{(m)}_+$, and denote by its underlying places of $F^{(m)}$ and $F$ by $w^{(m)}$ and $w$, respectively. Since $F_S/F^{(m)}$ is abelian, $w^{(m)}$ has a well-defined arithmetic Frobenius substitution $\Psi_{w^{(m)}}\in\Gal(F_S/F^{(m)})$. Denote by $G_{S,\gamma}$ the subset of $\Gal(F_S/F^{(m)})$ of elements $\Psi_{w^{(m)}}$ for all $\gamma$-associated places $w^{(m)}_+$.

\begin{lem}\label{le:associated}
Suppose that the order of $\gamma$ is coprime to $\ell$. Then we have
\[
G_{S,\gamma}=\theta_S^{-1}\Hom_{O_\lambda}(S,(\bar\rR^{(m)})^{h_\gamma}).
\]
\end{lem}

\begin{proof}
Note that the arithmetic Frobenius substitution of $w^{(m)}$ in $\Gal(F^{(m)}/F)$ coincides with $h_\gamma$, which implies that the action of $h_\gamma$ on $\Gal(F_S/F^{(m)})$ fixes $\Psi_{w^{(m)}}$. Thus, the image of $G_{S,\gamma}$ under $\theta_S$ is contained in $\Hom_{O_\lambda}(S,(\bar\rR^{(m)})^{h_\gamma})$.

Conversely, suppose that $\Psi\in\Gal(F_S/F^{(m)})$ satisfies $\theta_S(\Psi)\in\Hom_{O_\lambda}(S,(\bar\rR^{(m)})^{h_\gamma})$. We need to find a $\gamma$-associated place $w^{(m)}_+$ such that $\Psi=\Psi_{w^{(m)}}$. We regard $\gamma$ as an element in $\Gal(F^{(m)}_+/F^+)$ and $h_\gamma$ as an element in $\Gal(F^{(m)}/F)$. Let $g$ be the order of $h_\gamma$, which is coprime to $\ell$. Consider the element $(g^{-1}\Psi)h_\gamma\in\Gal(F_S/F)=\Gal(F_S/F^{(m)})\rtimes\Gal(F^{(m)}/F)$. Let $\tilde{F}_S$ be the smallest subfield of $\dC$ that is Galois over $F^+$ and contains $F_S$ and $F^{(m)}_+$. Since $\gamma$ has order prime to $\ell$, it is easy to see that there is an element $\tilde\gamma\in\Gal(\tilde{F}_S/F^+)$ lifting $\gamma$ such that the image of $\tilde\gamma^{[F:F^+]}\in\Gal(\tilde{F}_S/F)$ in $\Gal(F_S/F)$ coincides with $(g^{-1}\Psi)h_\gamma$. By the Chebotarev density theorem, we can find a place $\tilde{w}$ of $\tilde{F}_S$ whose arithmetic Frobenius substitution coincides with $\gamma$ and whose underlying place $w^{(m)}_+$ of $F^{(m)}_+$ is $\gamma$-associated. Then it is clear that $\Psi=\Psi_{w^{(m)}}$.
\end{proof}

By the above lemma, for every $r\in\dN$, we have a map
\[
\theta_{S,\gamma}^r\colon G_{S,\gamma}^r\to\Hom_{O_\lambda}(S,((\bar\rR^{(m)})^{h_\gamma})^{\oplus r})
\]
of abelian groups induced by $\theta_S$.

\begin{definition}\label{de:abundant}
Suppose that $S$ is a free $O_\lambda/\lambda^{m-m_0}$-module of rank $r_S$ for some $m_0\in\dN$ and $r_S\in\dN$. We say that an $r_S$-tuple $(\Psi_1,\dots,\Psi_{r_S})\in G_{S,\gamma}^{r_S}$ is \emph{$(S,\gamma)$-abundant} if the image of the homomorphism $\theta_{S,\gamma}^{r_S}(\Psi_1,\dots,\Psi_{r_S})$ contains $\lambda^{m_0+\ff(r_S)\fr_\rR}((\bar\rR^{(m)})^{h_\gamma})^{\oplus r_S}$, where $\fr_\rR$ and $\ff(r_S)$ are the integers appearing in Lemma \ref{le:reducibility} and Lemma \ref{le:image}, respectively.
\end{definition}

The following proposition provides $(S,\gamma)$-abundant tuples under certain conditions.

\begin{proposition}\label{pr:selmer_localization}
Suppose that $S$ is a free $O_\lambda/\lambda^{m-m_0}$-module of rank $r_S$ for some $m_0\in\dN$ and $r_S\in\dN$. Assume that the following are satisfied:
\begin{itemize}[label={\ding{109}}]
  \item $\rR_\dQ$ is absolutely irreducible;

  \item either one of the two assumptions in Lemma \ref{le:restriction} is satisfied;

  \item the order of $\gamma$ is coprime to $\ell$; and

  \item $(\bar\rR^{(m)})^{h_\gamma}$ is free over $O_\lambda/\lambda^m$ of rank $1$.
\end{itemize}
Then $(S,\gamma)$-abundant $r_S$-tuple exists.
\end{proposition}

\begin{proof}
By Lemma \ref{le:restriction}, $\Res_{\bar\rho^{(m)}}$ is injective. By Lemma \ref{le:reducibility} and Lemma \ref{le:image}, the $O_{\lambda}$-submodule generated by the image of $\theta_S$ contains $\lambda^{\ff(r_S)\fr_\rR}\Hom_{O_\lambda}(S,\bar\rR^{(m)})$.
Since $h_{\gamma}$ has order coprime to $\ell$, $\Hom_{O_\lambda}(S,(\bar\rR^{(m)})^{h_\gamma})$  is a direct summand of $\Hom_{O_\lambda}(S,\bar\rR^{(m)})$. It follows from Lemma \ref{le:associated} that the $O_\lambda$-submodule generated by $\theta_S(G_{S,\gamma})$ contains $\lambda^{\ff(r_S)\fr_\rR}\Hom_{O_\lambda}(S,(\bar\rR^{(m)})^{h_\gamma})$. As $(\bar\rR^{(m)})^{h_\gamma}$ is free $O_\lambda/\lambda^m$-module of rank $1$ and $S$ is a free $O_\lambda/\lambda^{m-m_0}$-module of rank $r_S$, the proposition follows immediately.
\end{proof}

\begin{proposition}\label{co:selmer_localization}
Let the assumptions be as in Proposition \ref{pr:selmer_localization} and put $r\coloneqq r_S$ for short. For every $(S,\gamma)$-abundant $r$-tuple $(\Psi_1,\dots,\Psi_r)$, one can choose a basis $\{s_1,\dots,s_r\}$ of $S$ such that $\theta_S(\Psi_i)(s_j)=0$ if $i\neq j$ and
\[
\exp_\lambda\(\theta_S(\Psi_j)(s_j),(\bar\rR^{(m)})^{h_\gamma}\)\geq m-m_0-\ff(r)\fr_\rR.
\]
Moreover, if we write $\Psi_i=\Psi_{w^{(m)}_i}$ with a $\gamma$-associated place $w^{(m)}_i$ of $F^{(m)}_+$ for $1\leq i\leq r$, then we have
$\loc_{w_i}(s_j)=0$ if $i\neq j$ and
\begin{align*}
\exp_\lambda\(\loc_{w_i}(s_i),\rH^1_\ns(F_{w_i},\bar\rR^{(m)})\)\geq  m-m_0-\ff(r)\fr_\rR.
\end{align*}
\end{proposition}

Note that by Definition \ref{de:associated} and Lemma \ref{le:global_local}, the image of $\loc_{w_i}\colon S\to\rH^1(F_{w_i},\bar\rR^{(m)})$ is contained in $\rH^1_\ns(F_{w_i},\bar\rR^{(m)})$.

\begin{proof}
The first part is obvious from Definition \ref{de:abundant}.

For the second part, note that $\rH^1_\ns(F^{(m)}_{w^{(m)}_i},\bar\rR^{(m)})$ is canonically isomorphic to $\bar\rR^{(m)}$ by evaluating on the element $\Psi_i=\Psi_{w^{(m)}_i}$. By the definition of $\theta_S$, the map $\theta_S(\Psi_i)\colon S\to\bar\rR^{(m)}$ coincides with the composite map
\[
S\xrightarrow{\loc_{w_i}}\rH^1_\ns(F_{w_i},\bar\rR^{(m)})\to\rH^1_\ns(F^{(m)}_{w^{(m)}_i},\bar\rR^{(m)})\simeq\bar\rR^{(m)}.
\]
The second part follows immediately.

The proposition is proved.
\end{proof}

\subsection{Case of Rankin--Selberg product}
\label{ss:rankin_selberg}

In this subsection, we discuss Galois modules that are related to Rankin--Selberg products. We take objects $\rR_\alpha\in\Mod(F,O_\lambda)_\free$ for $\alpha=0,1$ of rank $n_\alpha>0$ with the associated homomorphism $\rho_\alpha\colon\Gamma_F\to\GL(\rR_\alpha)$, together with a $(1-\alpha)$-polarization $\Xi_\alpha\colon\rR_\alpha^\tc\xrightarrow{\sim}\rR_\alpha^\vee(1-\alpha)$. We fix isomorphisms $\rR_\alpha\simeq O_\lambda^{\oplus n_\alpha}$ for $\alpha=0,1$.

We assume that $n_0=2r_0$ is even and $n_1=2r_1+1$ is odd. Put
\[
\rR\coloneqq\rR_0\otimes_{O_\lambda}\rR_1,\qquad\rho\coloneqq\rho_0\otimes\rho_1\colon\Gamma_F\to\GL(\rR),
\]
and $\Xi\coloneqq\Xi_0\otimes\Xi_1\colon\rR^\tc\xrightarrow{\sim}\rR^\vee(1)$ which is a $1$-polarization of $\rR$.

For a homomorphism $\rho$ from $\Gamma_F$ and a place $w$ of $F$, we write $\rho_w$ for the restriction of $\rho$ to the subgroup $\Gamma_{F_w}$. Moreover, for clarity, we denote by $\bar\epsilon_\ell^{(m)}\colon\Gamma_{F^+}\to(O_\lambda/\lambda^m)^\times$ the reduction of $\epsilon_\ell$ modulo $\lambda^m$ for a positive integer $m$, and put $\bar\epsilon_\ell\coloneqq\bar\epsilon_\ell^{(1)}$ for simplicity.

\begin{lem}\label{le:galois_element}
Let the notation be as above. Take a totally real finite Galois extension $F'/F^+$ contained in $\dC$ and a polynomial $\sP(T)\in\dZ[T]$. For every positive integer $m$, consider the following statement
\begin{description}
  \item[$(\r{GI}_{F',\sP}^m)$] The image of the restriction of the homomorphism
     \[
     (\bar\rho^{(m)}_{0+},\bar\rho^{(m)}_{1+},\bar\epsilon_\ell^{(m)})\colon\Gamma_{F^+}\to
     \sG_{n_0}(O_\lambda/\lambda^m)\times\sG_{n_1}(O_\lambda/\lambda^m)\times(O_\lambda/\lambda^m)^\times
     \]
     (see Notation \ref{no:reduction} for the notation) to $\Gal(\ol{F}/F')$ contains an element $(\gamma_0,\gamma_1,\xi)$ satisfying
     \begin{enumerate}[label=(\alph*)]
       \item $\sP(\xi)$ is invertible in $O_\lambda/\lambda^m$;

       \item for $\alpha=0,1$, $\gamma_\alpha$ belongs to $(\GL_{n_\alpha}(O_\lambda/\lambda^m)\times(O_\lambda/\lambda^m)^\times)\fc$ with order coprime to $\ell$;

       \item the kernels of $h_{\gamma_0}-1$, $h_{\gamma_1}-1$, and $h_{\gamma_0}\otimes h_{\gamma_1}-1$ (Notation \ref{no:h_gamma}) are all free over $O_\lambda/\lambda^m$ of rank $1$;

       \item if $[F:F^+]=2$, then $h_{\gamma_0}$ does not have an eigenvalue that is equal to $-1$ in $O_\lambda/\lambda$;

       \item if $[F:F^+]=2$, then $h_{\gamma_1}$ does not have an eigenvalue that is equal to $-\xi$ in $O_\lambda/\lambda$.
     \end{enumerate}
\end{description}
Then $(\r{GI}_{F',\sP}^1)$ implies $(\r{GI}_{F',\sP}^m)$ for every $m\geq 1$.
\end{lem}

\begin{proof}
Take an element $(\gamma_0,\gamma_1,\xi)$ obtained from $(\r{GI}_{F',\sP}^1)$. For every integer $m\geq 2$, we need to construct an element $(\gamma'_0,\gamma'_1,\xi')$ in the image of $(\bar\rho^{(m)}_{0+},\bar\rho^{(m)}_{1+},\bar\epsilon_\ell^{(m)})$ satisfying (a--e). First, we take $(\gamma'_0,\gamma'_1,\xi')$ to be an arbitrary lifting of $(\gamma_0,\gamma_1,\xi)$ in the image of $(\bar\rho^{(m)}_{0+},\bar\rho^{(m)}_{1+},\bar\epsilon_\ell^{(m)})$. Since the order of $\gamma_\alpha$ is coprime to $\ell$, there exists a positive integer $d_\alpha$ such that $\gamma_\alpha^{\ell^{d_\alpha}}=\gamma_\alpha$. On the other hand, we can find a positive integer $e_\alpha$ such that $(\gamma'_\alpha)^{\ell^{e_\alpha}}$ has order coprime to $\ell$ and that $1$ is an eigenvalue of $h_{(\gamma'_\alpha)^{\ell^{e_\alpha}}}$. Replacing $\gamma'_\alpha$ by $(\gamma'_\alpha)^{\ell^{d_\alpha e_\alpha}}$, we obtain the desired element $(\gamma'_0,\gamma'_1,\xi')$. The lemma follows.
\end{proof}

At the end of this section, we discuss an example using elliptic curves. Let $A_0$ and $A_1$ be two elliptic curves over $F^+$. For a rational prime $\ell$ (that is odd and unramified in $F$), we put
\[
\rR_\alpha\coloneqq(\Sym^{n_\alpha-1}_{\dZ_\ell}\rH^1_\et({A_\alpha}_{\ol{F}},\dZ_\ell))(r_\alpha)
\]
as a $\dZ_\ell[\Gamma_F]$-module for $\alpha=0,1$. Then $\rR_\alpha$ is an object in $\Mod(F,\dZ_\ell)_\free$ of rank $n_\alpha$ with a canonical $(1-\alpha)$-polarization $\Xi_\alpha\colon\rR_\alpha^\tc\xrightarrow{\sim}\rR_\alpha^\vee(1-\alpha)$. Put $\rR\coloneqq\rR_0\otimes_{\dZ_\ell}\rR_1$ and $\Xi\coloneqq\Xi_0\otimes\Xi_1$ as above.

\begin{proposition}\label{pr:elliptic_curve}
Suppose that ${A_0}_{\ol{F}}$ and ${A_1}_{\ol{F}}$ are not isogenous to each other and $\End({A_0}_{\ol{F}})=\End({A_1}_{\ol{F}})=\dZ$. Take a totally real finite Galois extension $F'/F^+$ contained in $\dC$ and a polynomial $\sP(T)\in\dZ[T]$. Then for sufficiently large $\ell$, we have that
\begin{enumerate}
  \item the image of $\bar\rho\colon\Gamma_F\to\GL(\rR\otimes\dF_\ell)$ contains a nontrivial scalar element;

  \item all of $\bar\rho_0$, $\bar\rho_1$, and $\bar\rho_0\otimes\bar\rho_1$ are absolutely irreducible; and

  \item $(\r{GI}_{F',\sP}^1)$ from Lemma \ref{le:galois_element} holds (with the coefficient field $E_\lambda=\dQ_\ell$).
\end{enumerate}
\end{proposition}

\begin{proof}
For $\alpha=0,1$ and every $\ell$, we have the homomorphism
\[
\bar\rho_{A_\alpha,\ell}\colon\Gamma_F\to\GL(\rH^1_\et({A_\alpha}_{\ol{F}},\dF_\ell))\simeq\GL_2(\dF_\ell).
\]
Then we have $\bar\rho_\alpha=(\Sym^{n_\alpha-1}\bar\rho_{A_\alpha,\ell})(r_\alpha)$ for $\alpha=0,1$. By our assumption on ${A_0}_{\ol{F}}$ and ${A_1}_{\ol{F}}$, and \cite{Ser72}*{Th\'{e}or\`{e}me~6}, for sufficiently large $\ell$, the image of the homomorphism
\[
(\bar\rho_{A_0,\ell},\bar\rho_{A_1,\ell},\bar\epsilon_\ell)\colon\Gamma_F\to\GL_2(\dF_\ell)\times\GL_2(\dF_\ell)\times\dF_\ell^\times
\]
consists exactly of the elements $(g_0,g_1,\xi)$ satisfying $\det g_0=\det g_1=\xi^{-1}$. Then both (1) and (2) follow immediately.

For (3), take an element $g\in\Gamma_F$ such that its image under $(\bar\rho_{A_0,\ell},\bar\rho_{A_1,\ell},\bar\epsilon_\ell)$ is in the conjugacy class of
\[
\(
\begin{pmatrix}
a & 0 \\
0 & 1
\end{pmatrix},
\begin{pmatrix}
ab & 0 \\
0 & b^{-1}
\end{pmatrix},
a^{-1}
\)
\]
for $a,b\in\dF_\ell^\times$ satisfying
\begin{itemize}[label={\ding{109}}]
  \item $\sP(a^{-1})\neq 0$,

  \item $(a^{2i}(ab^2)^{2j})^{[F':F^+]}\neq 1$ for $(i,j)\in\{r_0,r_0-1,\dots,1-r_0\}\times\{r_1,r_1-1,\dots,-r_1\}$ except for $(0,0)$,

  \item $(a^{2i-1})^{[F':F^+]}\neq -1$ for $i\in\{r_0,r_0-1,\dots,1-r_0\}$, and

  \item $(a(ab^2)^{2j})^{[F':F^+]}\neq -1$ for $j\in\{r_1,r_1-1,\dots,-r_1\}$.
\end{itemize}
Such pair $(a,b)$ always exists for sufficiently large $\ell$. Then it is straightforward to check that the image of $g^{[F':F^+]}\tc$ under $(\bar\rho_{0+},\bar\rho_{1+},\bar\epsilon_\ell)$ (under the notation of Lemma \ref{le:galois_element}) satisfies (a--e) of Lemma \ref{le:galois_element}. In particular, (3) follows.
\end{proof}

\section{Preliminaries on hermitian structures}
\label{ss:3}

In this section, we collect some constructions and results concerning objects carrying certain hermitian structures. In \S\ref{ss:hermitian_space}, we introduce hermitian spaces, their associated unitary groups and unitary Hecke algebras. In \S\ref{ss:unitary_shimura}, we introduce unitary Shimura varieties and unitary Shimura sets. In \S\ref{ss:cm}, we review the notion of (generalized) CM types. In \S\ref{ss:unitary_abelian_scheme}, we collect some facts about abelian schemes with hermitian structure, which will be parameterized by our unitary Shimura varieties. In \S\ref{ss:cm_moduli}, we introduce a moduli scheme parameterizing CM abelian varieties, which is an auxiliary moduli space in order to equip our unitary Shimura variety a moduli interpretation.

Let $N\geq 1$ be an integer.

\subsection{Unitary Satake parameters and unitary Hecke algebras}
\label{ss:hermitian_space}

We start by recalling the notion of the coefficient field for an automorphic representation of $\GL_N(\dA_F)$. Let $\Pi$ be an irreducible cuspidal automorphic (complex) representation of $\GL_N(\dA_F)$.

\begin{definition}[see \cite{Clo90}*{\S3.1}]\label{de:field_of_definition}
The \emph{coefficient field} of $\Pi$ is defined to be the smallest subfield of $\dC$, denoted by $\dQ(\Pi)$, such that for every $\rho\in\Aut(\dC/\dQ(\Pi))$, $\Pi^\infty$ and $\Pi^\infty\otimes_{\dC,\rho}\dC$ are isomorphic.
\end{definition}

For a nonarchimedean place $w$ of $F$ such that $\Pi_w$ is unramified, let
\[
\balpha(\Pi_w)\coloneqq\{\alpha(\Pi_w)_1,\dots,\alpha(\Pi_w)_N\}\subseteq\dC
\]
be the Satake parameter of $\Pi_w$ and $\dQ(\Pi_w)\subseteq\dC$ be the subfield generated by the coefficients of the polynomial
\[
\prod_{i=1}^N\(T-\alpha(\Pi_w)_i\cdot \sqrt{\|w\|}^{N-1}\)\in\dC[T].
\]

\begin{lem}\label{le:field_of_definition}
Suppose that $\Pi$ is regular algebraic \cite{Clo90}*{Definition~3.12}. Then the coefficient field $\dQ(\Pi)$ is a number field, and is the composition of $\dQ(\Pi_w)$ for all nonarchimedean places $w$ of $F$ such that $\Pi_w$ is unramified.
\end{lem}

\begin{proof}
By \cite{Clo90}*{Th\'{e}or\`{e}me~3.13}, $\dQ(\Pi)$ is a number field. Let $\dQ(\Pi)'$ be the composition of $\dQ(\Pi_w)$ for such $w$.

By the construction of unramified principal series, it is clear that for every $\gamma\in\Aut(\dC/\dQ(\Pi)')$ and every $w$ such that $\Pi_w$ is unramified, $\Pi_w$ and $\Pi_w\otimes_{\dC,\gamma}\dC$ have the same Satake parameter, hence are isomorphic. Since $\Pi$ is regular algebraic, by \cite{Clo90}*{Th\'{e}or\`{e}me~3.13}, there exists a cuspidal automorphic representation $\pres{\gamma}\Pi$ of $\GL_N(\dA_F)$ such that $\pres{\gamma}\Pi^\infty\simeq\Pi^\infty\otimes_{\dC,\gamma}\dC$. By the strong multiplicity one property for $\GL_N$ \cite{PS79}, we know that for $\gamma\in\Aut(\dC/\dQ(\Pi)')$, $\pres{\gamma}\Pi\simeq\Pi$, hence $\Pi^\infty\otimes_{\dC,\gamma}\dC\simeq\Pi^\infty$. It follows that $\dQ(\Pi)$ is contained in $\dQ(\Pi)'$.

Conversely, for $\gamma\in\Aut(\dC/\dQ(\Pi))$, $\Pi_w$ and $\Pi_w\otimes_{\dC,\gamma}\dC$ are isomorphic for every $w$. When $\Pi_w$ is unramified, $\dQ(\Pi_w)$ is simply the field of definition of $\Pi_w$, which implies that $\gamma$ fixes $\dQ(\Pi_w)$. It follows that $\dQ(\Pi')$ is contained in $\dQ(\Pi)$.

The lemma follows.
\end{proof}

\begin{definition}[Abstract Satake parameter]\label{de:satake_parameter}
Let $L$ be a ring. For a multi-subset $\balpha\coloneqq\{\alpha_1,\dots,\alpha_N\}\subseteq L$, we put
\[
P_{\balpha}(T)\coloneqq\prod_{i=1}^N(T-\alpha_i)\in L[T].
\]
Consider a nonarchimedean place $v$ of $F^+$ not in $\Sigma^+_\bad$.
\begin{enumerate}
  \item Suppose that $v$ is inert in $F$. We define an \emph{(abstract) Satake parameter} in $L$ at $v$ of rank $N$ to be a multi-subset $\balpha\subseteq L$ of cardinality $N$. We say that $\balpha$ is \emph{unitary} if $P_{\balpha}(T)=(-T)^N\cdot P_{\balpha}(T^{-1})$.

  \item Suppose that $v$ splits in $F$. We define an \emph{(abstract) Satake parameter} in $L$ at $v$ of rank $N$ to be a pair $\balpha\coloneqq(\balpha_1;\balpha_2)$ of multi-subsets $\balpha_1,\balpha_2\subseteq L$ of cardinality $N$, indexed by the two places $w_1,w_2$ of $F$ above $v$. We say that $\balpha$ is \emph{unitary} if $P_{\balpha_1}(T)=c\cdot T^N\cdot P_{\balpha_2}(T^{-1})$ for some constant $c\in L^\times$.
\end{enumerate}
For two Satake parameters $\balpha_0$ and $\balpha_1$ in $L$ at $v$ of rank $n_0$ and $n_1$, respectively, we may form their tensor product $\balpha_0\otimes\balpha_1$ which is of rank $n_0n_1$ in the obvious way. If $\balpha_0$ and $\balpha_1$ are both unitary, then so is $\balpha_0\otimes\balpha_1$.
\end{definition}

\begin{notation}\label{no:satake}
We denote by $\Sigma^+_\Pi$ the smallest (finite) set of nonarchimedean places of $F^+$ containing $\Sigma^+_\bad$ such that $\Pi_w$ is unramified for every nonarchimedean place $w$ of $F$ not above $\Sigma^+_\Pi$.

Take a nonarchimedean place $v$ of $F^+$ not in $\Sigma^+_\Pi$.
\begin{enumerate}
  \item If $v$ is inert in $F$, then we put $\balpha(\Pi_v)\coloneqq\balpha(\Pi_w)$ for the unique place $w$ of $F$ above $v$.

  \item If $v$ splits in $F$ into two places $w_1$ and $w_2$, then we put $\balpha(\Pi_v)\coloneqq(\balpha(\Pi_{w_1});\balpha(\Pi_{w_2}))$.
\end{enumerate}
Thus, $\balpha(\Pi_v)$ is a Satake parameter in $\dC$ at $v$ of rank $N$.
\end{notation}

\begin{definition}\label{de:satake_condition}
Let $v$ be a nonarchimedean place of $F^+$ inert in $F$, and $L$ a ring in which $\|v\|$ is invertible. Let $P\in L[T]$ be a monic polynomial of degree $N$ satisfying $P(T)=(-T)^N\cdot P(T^{-1})$.
\begin{enumerate}
  \item When $N$ is odd, we say that $P$ is \emph{Tate generic at $v$} if $P'(1)$ is invertible in $L$.

  \item When $N$ is odd, we say that $P$ is \emph{intertwining generic at $v$} if $P(-\|v\|)$ is invertible in $L$.

  \item When $N$ is even, we say that $P$ is \emph{level-raising special at $v$} if $P(\|v\|)=0$ and $P'(\|v\|)$ is invertible in $L$.

  \item When $N$ is even, we say that $P$ is \emph{intertwining generic at $v$} if $P(-1)$ is invertible in $L$.
\end{enumerate}
\end{definition}

\begin{remark}\label{re:satake_condition}
Suppose that $L$ is a field in Definition \ref{de:satake_condition}. It is easy to see that in Definition \ref{de:satake_condition}, if $P=P_{\balpha}$ for a unitary Satake parameter $\balpha$ in $L$ at $v$, then
\begin{enumerate}
  \item means that $1$ appears exactly once in $\balpha$;

  \item means that the pair $\{-\|v\|,-\|v\|^{-1}\}$ does not appear in $\balpha$;

  \item means that the pair $\{\|v\|,\|v\|^{-1}\}$ appears exactly once in $\balpha$;

  \item means that the pair $\{-1,-1\}$ does not appear in $\balpha$.
\end{enumerate}
Here, we note that when $N$ is odd, $1$ appears in $\balpha$ and all other elements appear in pairs of the form $\{\alpha,\alpha^{-1}\}$; when $N$ is even, elements in $\balpha$ appear in pairs of the form $\{\alpha,\alpha^{-1}\}$.
\end{remark}

We now introduce hermitian spaces.

\begin{definition}[Hermitian space]\label{de:hermitian_space}
Let $R$ be an $O_{F^+}[(\Sigma^+_\bad)^{-1}]$-ring. A \emph{hermitian space} over $O_F\otimes_{O_{F^+}}R$ of rank $N$ is a projective $O_F\otimes_{O_{F^+}}R$-module $\rV$ of rank $N$ together with a perfect pairing
\[
(\;,\;)_\rV\colon\rV\times\rV\to O_F\otimes_{O_{F^+}}R
\]
that is $O_F\otimes_{O_{F^+}}R$-linear in the first variable and $(O_F\otimes_{O_{F^+}}R,\tc\otimes\id_R)$-linear in the second variable, and satisfies $(x,y)_\rV=(y,x)_\rV^\tc$ for $x,y\in\rV$. We denote by $\rU(\rV)$ the group of $O_F\otimes_{O_{F^+}}R$-linear isometries of $\rV$, which is a reductive group over $R$.

Moreover, we denote by $\rV_\sharp$ the hermitian space $\rV\oplus O_F\otimes_{O_{F^+}}R\cdot 1$ where $1$ has norm $1$. For an $O_F\otimes_{O_{F^+}}R$-linear isometry $f\colon \rV\to\rV'$, we have the induced isometry $f_\sharp\colon\rV_\sharp\to\rV'_\sharp$.
\end{definition}

Let $v$ be a nonarchimedean place of $F^+$ not in $\Sigma^+_\bad$. Let $\Lambda_{N,v}$ be the unique up to isomorphism hermitian space over $O_{F_v}=O_F\otimes_{O_{F^+}}O_{F^+_v}$ of rank $N$, and $\rU_{N,v}$ its unitary group over $O_{F^+_v}$. Under a suitable basis, the associated  hermitian form of $\Lambda_{N,v}$ is given by the matrix
\[
\begin{pmatrix}
0 &\cdots & 0  & 1\\
0 & \cdots & 1 & 0\\
\vdots & \iddots & \vdots &\vdots\\
1 & \cdots & 0 &0
\end{pmatrix}.
\]
Consider the \emph{local spherical Hecke algebra}
\[
\dT_{N,v}\coloneqq\dZ[\rU_{N,v}(O_{F^+_v})\backslash\rU_{N,v}(F^+_v)/\rU_{N,v}(O_{F^+_v})].
\]
According to our convention, the unit element of $\dT_{N,v}$ is $\CF_{\rU_{N,v}(O_{F^+_v})}$. Let $\rA_{N,v}$ be the maximal split diagonal subtorus of $\rU_{N,v}$, and $\rX_*(\rA_{N,v})$ be its cocharacter group. Then there is a well-known Satake transform
\begin{align}\label{eq:satake_transform}
\dT_{N,v}\to\dZ[\|v\|^{\pm\delta(v)/2}][\rA_{N,v}(F^+_v)/\rA_{N,v}(O_{F^+_v})]\simeq \dZ[\|v\|^{\pm\delta(v)/2}][\rX_*(\rA_{N,v})]
\end{align}
as a homomorphism of algebras. Choose a uniformizer $\varpi_v$ of $F^+_v$.

\begin{construction}\label{cs:satake_hecke_pre}
Let $L$ be a $\dZ[\|v\|^{\pm\delta(v)/2}]$-ring. Let $\balpha$ be a unitary Satake parameter in $L$ at $v$ of rank $N$. There are two cases.
\begin{enumerate}
  \item Suppose that $v$ is inert in $F$. Then a set of representatives of $\rA_{N,v}(F^+_v)/\rA_{N,v}(O_{F^+_v})$ can be taken as
      \[
      \{\diag(\varpi_v^{t_1},\dots,\varpi_v^{t_N})\res t_1,\dots,t_N\in\dZ\text{ satisfying }t_i+t_{N+1-i}=0\text{ for all }1\leq i\leq N\}.
      \]
      Choose an ordering of $\balpha$ as $(\alpha_1,\dots,\alpha_N)$ satisfying $\alpha_i\alpha_{N+1-i}=1$; we have a unique homomorphism
      \[
      \dZ[\|v\|^{\pm\delta(v)/2}][\rA_{N,v}(F^+_v)/\rA_{N,v}(O_{F^+_v})]\to L
      \]
      of $\dZ[\|v\|^{\pm\delta(v)/2}]$-rings sending the class of $\diag(\varpi_v^{t_1},\dots,\varpi_v^{t_N})$ to
      $\prod_{i=1}^{\floor{\tfrac{N}{2}}}\alpha_i^{t_i}$. Composing with the Satake transform \eqref{eq:satake_transform}, we obtain a ring homomorphism
      \[
      \phi_{\balpha}\colon\dT_{N,v}\to L.
      \]
      It is independent of the choices of the uniformizer $\varpi_v$ and the ordering of $\balpha$.

  \item Suppose that $v$ splits in $F$ into two places $w_1$ and $w_2$. Then a set of representatives of $\rA_{N,v}(F^+_v)/\rA_{N,v}(O_{F^+_v})$ can be taken as
      \[
      \left\{
      \left.
      \(
      \diag(\varpi_v^{t_1},\dots,\varpi_v^{t_N}),
      \diag(\varpi_v^{-t_N},\dots,\varpi_v^{-t_1})
      \)
      \right| t_1,\dots,t_N\in\dZ
      \right\},
      \]
      where the first diagonal matrix (resp. the second diagonal matrix) is regarded as an element in $\rA_{N,v}(F_{w_1})$ (resp.\ $\rA_{N,v}(F_{w_2})$). Choose orders in $\balpha_1$ and $\balpha_2$ as $(\alpha_{1,1},\dots,\alpha_{1,N})$ and $(\alpha_{2,1},\dots,\alpha_{2,N})$ satisfying $\alpha_{1,i}\alpha_{2,N+1-i}=1$; we have a unique homomorphism
      \[
      \dZ[\|v\|^{\pm\delta(v)/2}][\rA_{N,v}(F^+_v)/\rA_{N,v}(O_{F^+_v})]\to L
      \]
      of $\dZ[\|v\|^{\pm\delta(v)/2}]$-rings sending the class of $\(\diag(\varpi_v^{t_1},\dots,\varpi_v^{t_N}),\diag(\varpi_v^{-t_N},\dots,\varpi_v^{-t_1})\)$ to $\prod_{i=1}^{N}\alpha_{1,i}^{t_i}$. Composing with the Satake transform \eqref{eq:satake_transform}, we obtain a ring homomorphism
      \[
      \phi_{\balpha}\colon\dT_{N,v}\to L.
      \]
      It is independent of the choices of the uniformizer $\varpi_v$, the order of the two places of $F$ above $v$, and the orders in $\balpha_1$ and $\balpha_2$.
\end{enumerate}
\end{construction}

\begin{definition}[Abstract unitary Hecke algebra]\label{de:abstract_hecke}
For a finite set $\Sigma^+$ of nonarchimedean places of $F^+$ containing $\Sigma^+_\bad$, we define the \emph{abstract unitary Hecke algebra away from $\Sigma^+$} to be the restricted tensor product
\[
\dT_N^{\Sigma^+}\coloneqq{\bigotimes_v}'\dT_{N,v}
\]
over all $v\not\in\Sigma^+_\infty\cup\Sigma^+$ with respect to unit elements. It is a ring.
\end{definition}

\begin{construction}\label{cs:satake_hecke}
Suppose that $\Pi$ satisfies $\Pi\circ\tc\simeq\Pi^\vee$. For $v\not\in\Sigma^+_\Pi$, the Satake parameter $\balpha(\Pi_v)$ is unitary. Thus by Construction \ref{cs:satake_hecke_pre}, we have a homomorphism
\[
\phi_\Pi\coloneqq\bigotimes_{v\not\in\Sigma^+_\infty\cup\Sigma^+_\Pi}\phi_{\balpha(\Pi_v)}\colon\dT_N^{\Sigma^+_\Pi}\to\dC,
\]
where we regard $\dC$ as a $\dZ[\|v\|^{\pm\delta(v)/2}]$-ring by sending $\|v\|^{\pm\delta(v)/2}$ to $\sqrt{\|v\|}^{\pm\delta(v)}$. If $\Pi$ is regular algebraic, then $\phi_\Pi$ takes values in $\dQ(\Pi)$ by Lemma \ref{le:field_of_definition}. Furthermore, by \cite{ST14}*{Proposition~4.1 \& Remark~4.2}, when $\Pi$ is relevant (Definition \ref{de:relevant}), $\phi_\Pi$ takes values in $O_{\dQ(\Pi)}$. In particular, we obtain a homomorphism
\[
\phi_\Pi\colon\dT_N^{\Sigma^+_\Pi}\to O_{\dQ(\Pi)}.
\]
\end{construction}

At last, we introduce some categories of open compact subgroups, which will be used later.

\begin{definition}\label{de:neat_category}
Let $\rV$ be a hermitian space over $F$ of rank $N$. Let $\Box$ be a finite set of nonarchimedean places of $F^+$.
\begin{enumerate}
  \item (Neat subgroups) For a nonarchimedean place $v$ of $F^+$ and an element $g_v\in\rU(\rV)({F^+_v})$, let $\Gamma(g_v)$ be the subgroup of $(\ol{F}^+_v)^\times$ generated by the eigenvalues of $g_v$ (regarded as an element in $\GL(\rV)(F_v)$), whose torsion subgroup $\Gamma(g_v)_{\r{tors}}$ lies in $\ol\dQ^\times$. We say an element $g=(g_v)\in\rU(\rV)(\dA_{F^+}^{\infty,\Box})$ is \emph{neat} if $\bigcap_{v\notin\Box}\Gamma(g_v)_{\mathrm{tors}}=\{1\}$, and a subgroup $\rK\subseteq\rU(\rV)(\dA_{F^+}^{\infty,\Box})$ is \emph{neat} if all its elements are neat.

  \item We define a category $\fK(\rV)^\Box$ whose objects are neat open compact subgroups $\rK$ of $\rU(\rV)(\dA_{F^+}^{\infty,\Box})$, and a morphism from $\rK$ to $\rK'$ is an element $g\in\rK\backslash\rU(\rV)(\dA_{F^+}^{\infty,\Box})/\rK'$ satisfying $g^{-1}\rK g\subseteq\rK'$. Denote by $\fK'(\rV)^\Box$ the subcategory of $\fK(\rV)^\Box$ that allows only identity double cosets as morphisms.

  \item We define a category $\fK(\rV)_\sp^\Box$ whose objects are pairs $\rK=(\rK_\flat,\rK_\sharp)$ where $\rK_\flat$ is an object of $\fK(\rV)^\Box$ and $\rK_\sharp$ is an object of $\fK(\rV_\sharp)^\Box$ such that $\rK_\flat\subseteq\rK_\sharp\cap\rU(\rV)(\dA_{F^+}^{\infty,\Box})$, and a morphism from $\rK=(\rK_\flat,\rK_\sharp)$ to $\rK'=(\rK'_\flat,\rK'_\sharp)$ is an element $g\in\rK_\flat\backslash\rU(\rV)(\dA_{F^+}^{\infty,\Box})/\rK'_\flat$ such that $g^{-1}\rK_\flat g\subseteq\rK'_\flat$ and $g^{-1}\rK_\sharp g\subseteq\rK'_\sharp$.\footnote{The subscript ``sp'' indicates that this notation will be related the special homomorphism of Shimura varieties later.} We have the obvious functors
      \begin{align*}
      \obj_\flat\colon\fK(\rV)_\sp^\Box\to\fK(\rV)^\Box,\quad
      \obj_\sharp\colon\fK(\rV)_\sp^\Box\to\fK(\rV_\sharp)^\Box
      \end{align*}
      sending $\rK=(\rK_\flat,\rK_\sharp)$ to $\rK_\flat$ and $\rK_\sharp$, respectively. Note that $\fK(\rV)_\sp^\Box$ is a non-full subcategory of $\fK(\rV)^\Box\times\fK(\rV_\sharp)^\Box$.
\end{enumerate}
When $\Box$ is the empty set, we suppress it from all the notations above.
\end{definition}

\subsection{Unitary Shimura varieties and sets}
\label{ss:unitary_shimura}

We introduce hermitian spaces over $F$ that will be used in this article.

\begin{definition}\label{de:standard_hermitian_space}
Let $\rV$ be a hermitian space over $F$ of rank $N$.
\begin{enumerate}
  \item We say that $\rV$ is \emph{standard definite} if it has signature $(N,0)$ at every place in $\Sigma^+_\infty$.

  \item We say that $\rV$ is \emph{standard indefinite} if it has signature $(N-1,1)$ at $\ul\tau_\infty$ and $(N,0)$ at other places in $\Sigma^+_\infty$.
\end{enumerate}
\end{definition}

First, we introduce unitary Shimura varieties. Take a standard indefinite hermitian space $\rV$ over $F$ of rank $N$. We have a functor
\begin{align*}
\Sh(\rV,\obj)\colon \fK(\rV)&\to\Sch_{/F} \\
\rK&\mapsto\Sh(\rV,\rK)
\end{align*}
of Shimura varieties associated to the reductive group $\Res_{F^+/\dQ}\rU(\rV)$ and the Deligne homomorphism
\begin{align*}
\rh\colon\Res_{\dC/\dR}\bG_m & \to(\Res_{F^+/\dQ}\rU(\rV))\otimes_\dQ\dR
=\prod_{\ul\tau\in\Sigma^+_\infty}\rU(\rV_{\ul\tau}) \\
z & \mapsto \(
\begin{pmatrix}
1_{N-1} &  \\
&  z^\tc/z \\
\end{pmatrix}
,
1_N,\dots,1_N
\)
\in\rU(\rV)(F^+_{\ul\tau_\infty})\times\prod_{\ul\tau\in\Sigma^+_\infty,\ul\tau\neq \ul\tau_\infty}\rU(\rV)(F^+_{\ul\tau}),\notag
\end{align*}
where we have identified $\rU(\rV)(F^+_{\ul\tau_\infty})$ with a subgroup of $\GL_N(\dC)$ via the complex embedding $\tau_\infty$ of $F$.

Second, we introduce unitary Shimura sets. Take a standard definite hermitian space $\rV$ over $F$ of rank $N$. We have a functor
\begin{align*}
\Sh(\rV,\obj)\colon \fK(\rV)&\to\Set \\
\rK&\mapsto\Sh(\rV,\rK)\coloneqq\rU(\rV)(F^+)\backslash\rU(\rV)(\dA_{F^+}^\infty)/\rK.
\end{align*}

\begin{remark}
Whether the notion $\Sh(\rV,\obj)$ stands for a scheme or a set depends on whether $\rV$ is standard indefinite or standard definite; so there will be no confusion about notation. Of course, one can equip $\Sh(\rV,\obj)$ with a natural scheme structure when $\rV$ is standard definite; but we will not take this point of view in this article.
\end{remark}

We now recall the notion of automorphic base change.

\begin{definition}[Automorphic base change]\label{de:bc_global}
Let $\rV$ be a hermitian space over $F$ of rank $N$, and $\pi$ an irreducible admissible representation of $\rU(\rV)(\dA_{F^+})$. An \emph{automorphic base change} of $\pi$ is defined to be an automorphic representation $\BC(\pi)$ of $\GL_N(\dA_F)$ that is a finite isobaric sum of discrete automorphic representations such that $\BC(\pi)_v\simeq\BC(\pi_v)$ holds for all but finitely many nonarchimedean places $v$ of $F^+$ such that $\pi_v$ is unramified. By the strong multiplicity one property for $\GL_N$ \cite{PS79}, if $\BC(\pi)$ exists, then it is unique up to isomorphism.
\end{definition}

\begin{proposition}\label{pr:galois}
Let $\Pi$ be a relevant representation of $\GL_N(\dA_F)$ (Definition \ref{de:relevant}).
\begin{enumerate}
  \item For every nonarchimedean place $w$ of $F$, $\Pi_w$ is tempered.

  \item For every rational prime $\ell$ and every isomorphism $\iota_\ell\colon\dC\xrightarrow{\sim}\ol\dQ_\ell$, there is a semisimple continuous homomorphism
      \[
      \rho_{\Pi,\iota_\ell}\colon\Gamma_F\to\GL_N(\ol\dQ_\ell),
      \]
      unique up to conjugation, satisfying that for every nonarchimedean place $w$ of $F$, the Frobenius semisimplification of the associated Weil--Deligne representation of $\rho_{\Pi,\iota_\ell}\res_{\Gamma_{F_w}}$ corresponds to the irreducible admissible representation $\iota_\ell\Pi_w|\det|_w^{\frac{1-N}{2}}$ of $\GL_N(F_w)$ under the local Langlands correspondence. Moreover, $\rho_{\Pi,\iota_\ell}^\tc$ and $\rho_{\Pi,\iota_\ell}^\vee(1-N)$ are conjugate.
\end{enumerate}
\end{proposition}

\begin{proof}
Part (1) is \cite{Car12}*{Theorem~1.2}. For (2), the Galois representation $\rho_{\Pi,\iota_\ell}$ is constructed in \cite{CH13}*{Theorem~3.2.3}, and the local-global compatibility is obtained in \cite{Car12}*{Theorem~1.1} and \cite{Car14}*{Theorem~1.1}. The last property in (2) follows from the previous one and the Chebotarev density theorem.
\end{proof}

\begin{definition}\label{de:weak_field}
Let $\Pi$ be a relevant representation of $\GL_N(\dA_F)$. We say that a subfield $E\subseteq\dC$ is a \emph{strong coefficient field} of $\Pi$ if $E$ is a number field containing $\dQ(\Pi)$ (Definition \ref{de:field_of_definition}); and for every prime $\lambda$ of $E$, there exists a continuous homomorphism
\[
\rho_{\Pi,\lambda}\colon\Gamma_F\to\GL_N(E_\lambda),
\]
necessarily unique up to conjugation, such that for every isomorphism $\iota_\ell\colon\dC\xrightarrow{\sim}\ol\dQ_\ell$ inducing the prime $\lambda$, $\rho_{\Pi,\lambda}\otimes_{E_\lambda}\ol\dQ_\ell$ and $\rho_{\Pi,\iota_\ell}$ are conjugate, where $\rho_{\Pi,\iota_\ell}$ is the homomorphism from Proposition \ref{pr:galois}(2).
\end{definition}

\begin{remark}\label{re:galois}
By \cite{CH13}*{Proposition~3.2.5}, a strong coefficient field of $\Pi$ exists for $\Pi$ relevant. Moreover, under Hypothesis \ref{hy:unitary_cohomology} below, $\dQ(\Pi)$ is already a strong coefficient field of $\Pi$ if $\Pi\simeq\BC(\pi)$ for a standard pair $(\rV,\pi)$ (see Definition \ref{de:standard_pair} below) in which $\rV$ is standard \emph{indefinite}.
\end{remark}

\begin{definition}\label{de:standard_pair}
Consider a pair $(\rV,\pi)$ where $\rV$ is a hermitian space over $F$ and $\pi$ is a discrete automorphic representation of $\rU(\rV)(\dA_{F^+})$. We say that $(\rV,\pi)$ is a \emph{standard pair} if either one of the following two situations happens:
\begin{enumerate}
  \item $\rV$ is standard definite, and $\pi^\infty$ appears in
     \[
     \varinjlim_{\rK\in\fK'(\rV)}\dC[\Sh(\rV,\rK)];
     \]

  \item $\rV$ is standard indefinite, and $\pi^\infty$ appears in
      \[
      \varinjlim_{\rK\in\fK'(\rV)}\iota_\ell^{-1}\rH^i_\et(\Sh(\rV,\rK)_{\ol{F}},\ol\dQ_\ell)
      \]
      for some isomorphism $\iota_\ell\colon\dC\xrightarrow{\sim}\ol\dQ_\ell$ and some $i\in\dZ$.
\end{enumerate}
\end{definition}

\begin{proposition}\label{pr:bc_global}
Let $(\rV,\pi)$ be a standard pair. Then $\BC(\pi)$ exists.
\end{proposition}

\begin{proof}
This is proved in \cite{Shi}*{Theorem~1.1}.\footnote{In fact, in \cite{Shi}, the author considers the case for unitary similitude group and assumes that $F$ contains an imaginary quadratic field. However, we can obtain the result in our setup by modifying the argument as in the proof of Proposition \ref{th:generic}.} When $\rV$ is standard definite, this is also proved in \cite{Lab}*{Corollaire~5.3}.
\end{proof}

\begin{remark}\label{re:bc_global}
In fact, in view of \cite{Shi}*{Theorem~1.1}, for a standard pair $(\rV,\pi)$, we have the associated Galois representation $\rho_{\BC(\pi),\iota_\ell}$ similar to the one in Proposition \ref{pr:galois} as well, with $N=\dim_F\rV$.
\end{remark}

\begin{hypothesis}\label{hy:unitary_cohomology}
Consider an integer $N\geq 1$. For every standard indefinite hermitian space $\rV$ over $F$ of rank $N$, every discrete automorphic representation $\pi$ of $\rU(\rV)(\dA_{F^+})$ such that $\BC(\pi)$ exists and is a relevant representation of $\GL_N(\dA_F)$, and every isomorphism $\iota_\ell\colon\dC\xrightarrow{\sim}\ol\dQ_\ell$, if $\rho_{\BC(\pi),\iota_\ell}$ is irreducible, then
\[
W^{N-1}(\pi)\coloneqq\Hom_{\ol\dQ_\ell[\rU(\rV)(\dA_{F^+}^\infty)]}
\(\iota_\ell\pi^\infty,\varinjlim_{\fK'(\rV)}\rH^{N-1}_\et(\Sh(\rV,\rK)_{\ol{F}},\ol\dQ_\ell)\)
\]
is isomorphic to the underlying $\ol\dQ_\ell[\Gamma_F]$-module of $\rho_{\BC(\pi),\iota_\ell}^\tc$.
\end{hypothesis}

\begin{proposition}\label{pr:unitary_cohomology}
Hypothesis \ref{hy:unitary_cohomology} holds for $N\leq 3$, and for $N>3$ if $F^+\neq\dQ$.
\end{proposition}

\begin{proof}
The case for $N=1$ follows directly from the definition of the canonical model of Shimura varieties over reflex fields. The case for $N=2$ is proved in \cite{Liu3}*{Theorem~D.6(2)}.\footnote{Note that our Deligne homomorphism is conjugate to the one in \cite{Liu3}*{\S C.1}, which is responsible for the $\tc$-conjugation in $\rho_{\BC(\pi),\iota_\ell}^\tc$.} The case for $N=3$ when $F^+=\dQ$ follows from the main result of \cite{Rog92}. The case for $N\geq 3$ when $F^+\neq\dQ$ will be proved in \cite{KSZ}.
\end{proof}

\subsection{Generalized CM type and reflexive closure}
\label{ss:cm}

We denote by $\dN[\Sigma_\infty]$ the commutative monoid freely generated by the set $\Sigma_\infty$, which admits an action of $\Aut(\dC/\dQ)$ via the set $\Sigma_\infty$.

\begin{definition}\label{de:cm_type}
A \emph{generalized CM type of rank $N$} is an element
\[
\Psi=\sum_{\tau\in\Sigma_\infty}r_\tau\tau\in\dN[\Sigma_\infty]
\]
satisfying $r_\tau+r_{\tau^\tc}=N$ for every $\tau\in\Sigma_\infty$. For such $\Psi$, we define its \emph{reflex field} $F_\Psi\subseteq\dC$ to be the fixed subfield of the stabilizer of $\Psi$ in $\Aut(\dC/\dQ)$. A \emph{CM type} is simply a generalized CM type of rank $1$. For a CM type $\Phi$, we say that $\Phi$ contains $\tau$ if its coefficient $r_\tau$ equals $1$.
\end{definition}

\begin{definition}\label{de:reflexive}
We define the \emph{reflexive closure} of $F$, denoted by $F_{\r{rflx}}$, to be the subfield of $\dC$ generated by $F$ and $F_\Phi$ for every CM type $\Phi$ of $F$. Put $F^+_{\r{rflx}}\coloneqq(F_{\r{rflx}})^{\tc=1}$.
\end{definition}

\begin{remark}
It is clear that $F_{\r{rflx}}$ is a CM field finite Galois over $F$; $F_{\r{rflx}}^+$ is the maximal totally real subfield of $F_{\r{rflx}}$ and is finite Galois over $F^+$. In many cases, we have $F_{\r{rflx}}=F$ and hence $F_{\r{rflx}}^+=F^+$, for example, when $F$ is Galois or contains an imaginary quadratic field.
\end{remark}

\begin{definition}\label{de:special_inert}
We say that a prime $\fp$ of $F^+$ is \emph{special inert} if the following are satisfied:
\begin{enumerate}
  \item $\fp$ is inert in $F$;

  \item the underlying rational prime $p$ of $\fp$ is odd and is unramified in $F$;

  \item $\fp$ is of degree one over $\dQ$, that is, $F^+_\fp=\dQ_p$.
\end{enumerate}
By abuse of notation, we also denote by $\fp$ for its induced prime of $F$.

We say that a special inert prime $\fp$ of $F^+$ is \emph{very special inert} if there exists a prime $\fp'$ of $F^+_{\r{rflx}}$ above $\fp$ satisfying $(F^+_{\r{rflx}})_{\fp'}=F^+_\fp(=\dQ_p)$.\footnote{This is equivalent to that for every prime $\fq$ of $F^+$ above $p$ that is inert in $F$, $[F^+_{\fq}:\dQ_p]$ is odd.}
\end{definition}

\begin{remark}
In Definition \ref{de:special_inert}, (3) is proposed only for the purpose of simplifying computations on Dieudonn\'{e} modules in \S\ref{ss:qs} and \S\ref{ss:ns}; it is not really necessary as results in these two sections should remain valid without (3). However, dropping (3) will vastly increase the burden of notations and computations in those two sections, where the technicality is already heavy.
\end{remark}

In what follows in this article, we will often take a rational prime $p$ that is unramified in $F$, and an isomorphism $\iota_p\colon\dC\xrightarrow{\sim}\ol\dQ_p$. By composing with $\iota_p$, we regard $\Sigma_\infty$ also as the set of $p$-adic embeddings of $F$. We also regard $\dQ_p$ as a subfield of $\dC$ via $\iota_p^{-1}$.

\begin{notation}\label{no:p_notation}
We introduce the following important notations.
\begin{enumerate}
  \item In what follows, whenever we introduce some finite unramified extension $\dQ_?^?$ of $\dQ_p$, we denote by $\dZ_?^?$ its ring of integers and put $\dF_?^?\coloneqq\dZ_?^?/p\dZ_?^?$.

  \item For every $\tau\in\Sigma_\infty$, we denote by $\dQ_p^\tau\subseteq\dC$ the composition of $\tau(F)$ and $\dQ_p$, which is unramified over $\dQ_p$. For a scheme $S\in\Sch_{/\dZ_p^\tau}$ and an $\cO_S$-module $\cF$ with an action $O_F\to\End_{\cO_S}(\cF)$, we denote by $\cF_\tau$ the maximal $\cO_S$-submodule of $\cF$ on which $O_F$ acts via the homomorphism $\tau\colon O_F\to\dZ_p^\tau\to\cO_S$.

  \item For every $\tau\in\Sigma_\infty$, we denote by $\dQ_p^\diamondsuit\subseteq\dC$ the composition of $\dQ_p^\tau$, which is unramified over $\dQ_p$. We can identify $\Sigma_\infty$ with $\Hom(O_F,\dZ_p^\diamondsuit)=\Hom(O_F,\dF_p^\diamondsuit)$. In particular, the $p$-power Frobenius map $\sigma$ acts on $\Sigma_\infty$.

  \item For a generalized CM type $\Psi$ of rank $N$, we denote by $\dQ_p^\Psi\subseteq\dC$ the composition of $\dQ_p$, $F$, and $F_\Psi$, which is contained in $\dQ_p^\diamondsuit$.

  \item For a (functor in) scheme over $\dZ^?_?$ written like $\bX_?(\cdot\cdot\cdot)$, we put $\rX_?(\cdot\cdot\cdot)\coloneqq\bX_?(\cdot\cdot\cdot)\otimes_{\dZ^?_?}\dF^?_?$ and $\bX^\eta_?(\cdot\cdot\cdot)\coloneqq\bX_?(\cdot\cdot\cdot)\otimes_{\dZ^?_?}\dQ^?_?$. For a (functor in) scheme over $\dF^?_?$ written like $\rX^?_?(\cdot\cdot\cdot)$, we put $\ol\rX^?_?(\cdot\cdot\cdot)\coloneqq\rX^?_?(\cdot\cdot\cdot)\otimes_{\dF^?_?}\ol\dF_p$. Similar conventions are applied to morphisms as well.
\end{enumerate}
\end{notation}

\subsection{Unitary abelian schemes}
\label{ss:unitary_abelian_scheme}

We first introduce some general notations about abelian schemes.

\begin{notation}
Let $A$ be an abelian scheme over a scheme $S$. We denote by $A^\vee$ the \emph{dual abelian variety} of $A$ over $S$. We denote by $\rH^\dr_1(A/S)$ (resp.\ $\Lie_{A/S}$, and $\omega_{A/S}$) for the \emph{relative de Rham homology} (resp.\ \emph{Lie algebra}, and \emph{dual Lie algebra}) of $A/S$, all regarded as locally free $\cO_S$-modules. We have the following \emph{Hodge exact sequence}
\begin{align}\label{eq:hodge_sequence}
0 \to \omega_{A^\vee/S} \to \rH^\dr_1(A/S) \to \Lie_{A/S} \to 0
\end{align}
of sheaves on $S$. When the base $S$ is clear from the context, we sometimes suppress it from the notation.
\end{notation}

\begin{definition}[Unitary abelian scheme]\label{de:unitary_abelian_scheme}
We prescribe a subring $\dP\subseteq\dQ$. Let $S$ be a scheme in $\Sch_{/\dP}$.
\begin{enumerate}
  \item An \emph{$O_F$-abelian scheme} over $S$ is a pair $(A,i)$ in which $A$ is an abelian scheme over $S$ and $i\colon O_F\to\End_S(A)\otimes\dP$ is a homomorphism of algebras sending $1$ to the identity endomorphism.

  \item A \emph{unitary $O_F$-abelian scheme} over $S$ is a triple $(A,i,\lambda)$ in which $(A,i)$ is an $O_F$-abelian scheme over $S$, and $\lambda\colon A\to A^\vee$ is a quasi-polarization such that $i(a^\tc)^\vee\circ\lambda=\lambda\circ i(a)$ for every $a\in O_F$, and there exists $c\in\dP^\times$ making $c\lambda$ a polarization.

  \item For two $O_F$-abelian schemes $(A,i)$ and $(A',i')$ over $S$, a (quasi-)homomorphism from $(A,i)$ to $(A',i')$ is a (quasi-)homomorphism $\varphi\colon A\to A'$ such that $\varphi\circ i(a)=i'(a)\circ\varphi$ for every $a\in O_F$. We will usually refer to such $\varphi$ as an $O_F$-linear (quasi-)homomorphism.
\end{enumerate}
Moreover, we will usually suppress the notion $i$ if the argument is insensitive to it.
\end{definition}

\begin{definition}[Signature type]\label{de:signature}
Let $\Psi$ be a generalized CM type of rank $N$ (Definition \ref{de:cm_type}). Consider a scheme $S\in\Sch_{/O_{F_\Psi}\otimes\dP}$. We say that an $O_F$-abelian scheme $(A,i)$ over $S$ has \emph{signature type} $\Psi$ if for every $a\in O_F$, the characteristic polynomial of $i(a)$ on $\Lie_{A/S}$ is given by
\[
\prod_{\tau\in\Sigma_\infty}(T-\tau(a))^{r_\tau}\in\cO_S[T].
\]
\end{definition}

\begin{construction}\label{cs:hermitian_structure}
Let $K$ be an $O_{F_\Psi}\otimes\dP$-ring that is an algebraically closed field. Suppose that we are given a unitary $O_F$-abelian scheme $(A_0,i_0,\lambda_0)$ over $K$ of signature type $\Phi$ that is a CM type, and a unitary $O_F$-abelian scheme $(A,i,\lambda)$ over $K$ of signature type $\Psi$. For every set $\Box$ of places of $\dQ$ containing $\infty$ and the characteristic of $K$, if not zero, we construct a hermitian space
\[
\Hom_{F\otimes_\dQ\dA^\Box}^{\lambda_0,\lambda}(\rH^\et_1(A_0,\dA^\Box),\rH^\et_1(A,\dA^\Box))
\]
over $F\otimes_\dQ\dA^\Box=F\otimes_{F^+}(F^+\otimes_\dQ\dA^\Box)$, with the underlying $F\otimes_\dQ\dA^\Box$-module
\[
\Hom_{F\otimes_\dQ\dA^\Box}(\rH^\et_1(A_0,\dA^\Box),\rH^\et_1(A,\dA^\Box))
\]
equipped with the pairing
\[
(x,y)\coloneqq i_0^{-1}\((\lambda_{0*})^{-1}\circ y^\vee\circ\lambda_*\circ x\)\in i_0^{-1}\End_{F\otimes_\dQ\dA^\Box}(\rH^\et_1(A_0,\dA^\Box))=F\otimes_\dQ\dA^\Box.
\]
\end{construction}

Now we take a rational prime $p$ that is unramified in $F$, and take the prescribed subring $\dP$ in Definition \ref{de:unitary_abelian_scheme} to be $\dZ_{(p)}$. We also choose an isomorphism $\iota_p\colon\dC\simeq \ol\dQ_p$, and adopt Notation \ref{no:p_notation}.

\begin{definition}\label{de:p_quasi}
Let $A$ and $B$ be two abelian schemes over a scheme $S\in\Sch_{/\dZ_{(p)}}$. We say that a quasi-homomorphism (resp.\ quasi-isogeny) $\varphi\colon A\to B$ is a \emph{quasi-$p$-homomorphism} (resp.\ \emph{quasi-$p$-isogeny}) if there exists some $c\in\dZ_{(p)}^\times$ such that $c\varphi$ is a homomorphism (resp.\ isogeny). A quasi-isogeny $\varphi$ is \emph{prime-to-$p$} if both $\varphi$ and $\varphi^{-1}$ are quasi-$p$-isogenies. We say that a quasi-polarization $\lambda$ of $A$ is \emph{$p$-principal} if $\lambda$ is a prime-to-$p$ quasi-isogeny.
\end{definition}

Note that for a unitary $O_F$-abelian scheme $(A,i,\lambda)$, the quasi-polarization $\lambda$ is a quasi-$p$-isogeny. To continue, take a generalized CM type $\Psi=\sum_{\tau\in\Sigma_\infty}r_\tau\tau$ of rank $N$.

\begin{remark}\label{re:hodge_sequence}
Let $A$ be an $O_F$-abelian scheme of signature type $\Psi$ over a scheme $S\in\Sch_{/\dZ_p^\tau}$ for some $\tau\in\Sigma_\infty$. Then \eqref{eq:hodge_sequence} induces a short exact sequence
\[
0 \to \omega_{A^\vee/S,\tau} \to \rH^\dr_1(A/S)_\tau \to \Lie_{A/S,\tau} \to 0
\]
of locally free $\cO_S$-modules of ranks $N-r_\tau$, $N$, and $r_\tau$, respectively. If $S$ belongs to $\Sch_{/\dZ_p^\diamondsuit}$, then we have decompositions
\begin{align*}
\rH^\dr_1(A/S)&=\bigoplus_{\tau\in\Sigma_\infty}\rH^\dr_1(A/S)_\tau,\\
\Lie_{A/S}&=\bigoplus_{\tau\in\Sigma_\infty}\Lie_{A/S,\tau},\\
\omega_{A/S}&=\bigoplus_{\tau\in\Sigma_\infty}\omega_{A/S,\tau}
\end{align*}
of locally free $\cO_S$-modules.
\end{remark}

\begin{notation}\label{no:weil_pairing}
Take $\tau\in\Sigma_\infty$. Let $(A,\lambda)$ be a unitary $O_F$-abelian scheme of signature type $\Psi$ over a scheme $S\in\Sch_{/\dZ_p^\tau}$. We denote
\[
\langle\;,\;\rangle_{\lambda,\tau}\colon\rH^\dr_1(A/S)_\tau\times\rH^\dr_1(A/S)_{\tau^\tc}\to\cO_S
\]
the $\cO_S$-bilinear pairing induced by the quasi-polarization $\lambda$, which is perfect if and only if $\lambda$ is $p$-principal. Moreover, for an $\cO_S$-submodule $\cF\subseteq\rH^\dr_1(A/S)_\tau$, we denote by $\cF^\perp\subseteq\rH^\dr_1(A/S)_{\tau^\tc}$ its (right) orthogonal complement under the above pairing, if $\lambda$ is clear from the context.
\end{notation}

Next we review some facts from the Serre--Tate theory \cite{Kat81} and the Grothendieck--Messing theory \cite{Mes72}, tailored to our application. Let $\Psi$ be a generalized CM type of rank $N$ such that $\min\{r_\tau,r_{\tau^\tc}\}=0$ for every $\tau$ not above $\ul\tau_\infty$. Consider a closed immersion $S\hookrightarrow\hat{S}$ in $\Sch_{/\dZ_p^\Psi}$ on which $p$ is locally nilpotent, with its ideal sheaf equipped with a PD structure, and a unitary $O_F$-abelian scheme $(A,\lambda)$ of signature type $\Psi$ over $S$. We let $\rH^\cris_1(A/\hat{S})$ be the evaluation of the first relative crystalline homology of $A/S$ at the PD-thickening $S\hookrightarrow\hat{S}$, which is a locally free $\cO_{\hat{S}}\otimes O_F$-module. The polarization $\lambda$ induces a pairing
\begin{align}\label{eq:weil_pairing_cris}
\langle\;,\;\rangle_{\lambda,\tau_\infty}^\cris\colon
\rH^\cris_1(A/\hat{S})_{\tau_\infty}\times\rH^\cris_1(A/\hat{S})_{\tau_\infty^\tc}\to\cO_{\hat{S}}.
\end{align}
We define two groupoids
\begin{itemize}[label={\ding{109}}]
  \item $\Def(S,\hat{S};A,\lambda)$, whose objects are unitary $O_F$-abelian schemes $(\hat{A},\hat{\lambda})$ of signature type $\Psi$ over $\hat{S}$ that lift $(A,\lambda)$;

  \item $\Def'(S,\hat{S};A,\lambda)$, whose objects are pairs $(\hat\omega_{\tau_\infty},\hat\omega_{\tau_\infty^\tc})$ where for each $\tau=\tau_\infty,\tau_\infty^\tc$, $\hat\omega_\tau\subseteq\rH_1^\cris(A/\hat{S})_\tau$ is a subbundle that lifts $\omega_{A^\vee/S,\tau}\subseteq\rH^\dr_1(A/S)_\tau$, such that $\langle\hat\omega_{\tau_\infty},\hat\omega_{\tau_\infty^\tc}\rangle_{\lambda,\tau_\infty}^\cris=0$.
\end{itemize}

\begin{proposition}\label{pr:deformation}
The functor from $\Def(S,\hat{S};A,\lambda)$ to $\Def'(S,\hat{S};A,\lambda)$ sending $(\hat{A},\hat{\lambda})$ to $(\omega_{\hat{A}^\vee/\hat{S},\tau_\infty},\omega_{\hat{A}^\vee/\hat{S},\tau_\infty^\tc})$ is a natural equivalence.
\end{proposition}

\begin{proof}
By \'{e}tale descent, we may replace $S\hookrightarrow\hat{S}$ by $S\otimes_{\dZ_p^\Psi}\dZ_p^\diamondsuit\hookrightarrow\hat{S}\otimes_{\dZ_p^\Psi}\dZ_p^\diamondsuit$. Then we have a decomposition
\[
\rH^\cris_1(A/\hat{S})=\bigoplus_{\tau\in\Sigma_\infty}\rH^\cris_1(A/\hat{S})_\tau
\]
similar to the one in Notation \ref{no:p_notation}. Note that for $\tau\not\in\{\tau_\infty,\tau_\infty^\tc\}$, the subbundle $\omega_{A^\vee/S,\tau}$ has a unique lifting to either zero or the entire $\rH^\cris_1(A/\hat{S})_\tau$. Thus, the proposition follows from the Serre--Tate and Grothendieck--Messing theories.
\end{proof}

To end this subsection, we review some notions for abelian schemes in characteristic $p$.

\begin{notation}\label{no:frobenius_verschiebung}
Let $A$ be an abelian scheme over a scheme $S\in\Sch_{/\dF_p}$. Put
\[
A^{(p)}\coloneqq A\times_{S,\sigma}S,
\]
where $\sigma$ is the absolute Frobenius morphism of $S$.
Then we have
\begin{enumerate}
  \item a canonical isomorphism $\rH^\dr_1(A^{(p)}/S)\simeq\sigma^*\rH^\dr_1(A/S)$ of $\cO_S$-modules;

  \item the Frobenius homomorphism $\Fr_A\colon A\to A^{(p)}$ which induces the \emph{Verschiebung map}
      \[
      \tV_A\coloneqq(\Fr_A)_*\colon\rH^\dr_1(A/S)\to\rH^\dr_1(A^{(p)}/S)
      \]
      of $\cO_S$-modules;

  \item the Verschiebung homomorphism $\Ver_A\colon A^{(p)}\to A$ which induces the \emph{Frobenius map}
      \[
      \tF_A\coloneqq(\Ver_A)_*\colon\rH^\dr_1(A^{(p)}/S)\to\rH^\dr_1(A/S)
      \]
      of $\cO_S$-modules.
\end{enumerate}
For a subbundle $H$ of $\rH^\dr_1(A/S)$, we denote by $H^{(p)}$ the subbundle of $\rH^\dr_1(A^{(p)}/S)$ that corresponds to $\sigma^*H$ under the isomorphism in (1). In what follows, we will suppress $A$ in the notations $\tF_A$ and $\tV_A$ if the reference to $A$ is clear.
\end{notation}

In Notation \ref{no:frobenius_verschiebung}, we have $\Ker\tF=\IM\tV=\omega_{A^{(p)}/S}$ and $\Ker\tV=\IM\tF$. Take $\tau\in\Sigma_\infty$. For a scheme $S\in\Sch_{/\dF_p^\tau}$ and an $O_F$-abelian scheme $A$ over $S$, we have $(\rH^\dr_1(A/S)_\tau)^{(p)}=\rH^\dr_1(A^{(p)}/S)_{\sigma\tau}$ under Notations \ref{no:p_notation} and \ref{no:frobenius_verschiebung}.

\begin{notation}\label{re:frobenius_verschiebung}
Suppose that $S=\Spec\kappa$ for a field $\kappa$ of characteristic $p$. Then we have a canonical isomorphism $\rH^\dr_1(A^{(p)}/\kappa)\simeq\rH^\dr_1(A/\kappa)\otimes_{\kappa,\sigma}\kappa$.
\begin{enumerate}
  \item By abuse of notation, we have
     \begin{itemize}[label={\ding{109}}]
        \item the $(\kappa,\sigma)$-linear Frobenius map $\tF\colon\rH^\dr_1(A/\kappa)\to\rH^\dr_1(A/\kappa)$ and

        \item if $\kappa$ is perfect, the $(\kappa,\sigma^{-1})$-linear Verschiebung map $\tV\colon\rH^\dr_1(A/\kappa)\to\rH^\dr_1(A/\kappa)$.
     \end{itemize}

  \item When $\kappa$ is perfect, recall that we have the \emph{covariant} Dieudonn\'{e} module $\cD(A)$ associated to the $p$-divisible group $A[p^\infty]$, which is a free $W(\kappa)$-module, such that $\cD(A)/p\cD(A)$ is canonically isomorphic to $\rH^\dr_1(A/\kappa)$. Again by abuse of notation, we have
     \begin{itemize}[label={\ding{109}}]
       \item the $(W(\kappa),\sigma)$-linear Frobenius map $\tF\colon\cD(A)\to\cD(A)$ lifting the one above, and

       \item the $(W(\kappa),\sigma^{-1})$-linear Verschiebung map $\tV\colon\cD(A)\to\cD(A)$ lifting the one above,
     \end{itemize}
     respectively, satisfying $\tF\circ\tV=\tV\circ\tF=p$.

  \item When $\kappa$ is perfect and contains $\dF_p^\tau$ for some $\tau\in\Sigma_\infty$, applying Notation \ref{no:p_notation} to the $W(\kappa)$-module $\cD(A)$, we obtain $W(\kappa)$-submodules $\cD(A)_{\sigma^i\tau}\subseteq\cD(A)$ for every $i\in\dZ$. Thus, we obtain
     \begin{itemize}[label={\ding{109}}]
       \item the $(W(\kappa),\sigma)$-linear Frobenius map $\tF\colon\cD(A)_\tau\to\cD(A)_{\sigma\tau}$ and

       \item the $(W(\kappa),\sigma^{-1})$-linear Verschiebung map $\tV\colon\cD(A)_\tau\to\cD(A)_{\sigma^{-1}\tau}$
     \end{itemize}
     by restriction. We have canonical isomorphisms and inclusions:
     \[
     \tV\cD(A)_{\sigma\tau}/p\cD(A)_\tau\simeq\omega_{A^\vee,\tau}\subseteq
     \cD(A)_\tau/p\cD(A)_\tau\simeq\rH^\dr_1(A)_\tau.
     \]
\end{enumerate}
\end{notation}

\begin{notation}\label{no:weil_pairing_dieudonne}
Take $\tau\in\Sigma_\infty$. Let $(A,\lambda)$ be a unitary $O_F$-abelian scheme of signature type $\Psi$ over $\Spec\kappa$ for a perfect field $\kappa$ containing $\dF_p^\tau$. We have a pairing
\[
\langle\;,\;\rangle_{\lambda,\tau}\colon\cD(A)_\tau\times\cD(A)_{\tau^\tc}\to W(\kappa)
\]
lifting the one in Notation \ref{no:weil_pairing}. We denote by $\cD(A)_\tau^\vee$ the $W(\kappa)$-dual of $\cD(A)_\tau$, as a submodule of $\cD(A)_{\tau^\tc}\otimes\dQ$. In what follows, unless we specify, the dual is always with respect to the default quasi-polarization.
\end{notation}

The following lemma will be repeatedly used in later discussion.

\begin{lem}\label{le:inverse_isogeny}
Suppose that $F^+$ is contained in $\dQ_p$ (via the embedding $\tau\colon F^+\hra\dC\simeq\ol\dQ_p$) with $\fp$ the induced $p$-adic  prime. Let $\varpi\in O_{F^+}$ be an element such that $\val_\fp(\varpi)=1$. Consider two  $O_F$-abelian schemes $A$ and $B$ over a scheme $S\in\Sch_{/\dF_{p^2}}$. Let $\alpha\colon A\to B$ and $\beta\colon B\to A$ be two $O_F$-linear quasi-$p$-isogenies (Definition \ref{de:p_quasi}) such that $\beta\circ\alpha=\varpi\cdot\id_A$ (hence $\alpha\circ\beta=\varpi\cdot\id_B$). Then
\begin{enumerate}
  \item For $\tau\in\{\tau_\infty,\tau_\infty^\tc\}$, the induced maps
      \begin{align*}
      \alpha_{*,\tau}&\colon\rH^\dr_1(A/S)_\tau\to\rH^\dr_1(B/S)_\tau,\\
      \beta_{*,\tau}&\colon\rH^\dr_1(B/S)_\tau\to\rH^\dr_1(A/S)_\tau
      \end{align*}
      satisfy the relations $\Ker\alpha_{*,\tau}=\IM\beta_{*,\tau}$ and $\Ker\beta_{*,\tau}=\IM\alpha_{*,\tau}$; and these kernels and images are locally free $\cO_S$-modules.

  \item We have
      \[
      \rank_{\cO_S}\Lie_{B/S,\tau_\infty}-\rank_{\cO_S}\Lie_{A/S,\tau_\infty}
      =\rank_{\cO_S}(\Ker\alpha_{*,\tau_\infty})-\rank_{\cO_S}(\Ker\alpha_{*,\tau_\infty^\tc}).
      \]

  \item Let $\lambda_A$ and $\lambda_B$ be two quasi-polarizations on $A$ and $B$, respectively, such that $(A,\lambda_A)$ and $(B,\lambda_B)$ become unitary $O_F$-abelian schemes of dimension $N[F^+:\dQ]$ for some integer $N\geq 1$. Suppose that $\alpha^\vee\circ\lambda_B\circ\alpha=\varpi\lambda_A$.
      \begin{enumerate}
      \item If both $\lambda_A$ and $\lambda_B$ are $p$-principal, then we have
          \[
          \rank_{\cO_S}(\Ker\alpha_{*,\tau_\infty})+\rank_{\cO_S}(\Ker\alpha_{*,\tau_\infty^\tc})=N.
          \]

      \item If $\lambda_A$ is $p$-principal and $\Ker\lambda_B[\fp^\infty]$ is of rank $p^2$, then we have
          \[
          \rank_{\cO_S}(\Ker\alpha_{*,\tau_\infty})+\rank_{\cO_S}(\Ker\alpha_{*,\tau_\infty^\tc})=N-1.
          \]

      \item If $\Ker\lambda_A[\fp^\infty]$ is of rank $p^2$ and $\lambda_B$ is $p$-principal, then we have
          \[
          \rank_{\cO_S}(\Ker\alpha_{*,\tau_\infty})+\rank_{\cO_S}(\Ker\alpha_{*,\tau_\infty^\tc})=N+1.
          \]

      \item If both $\Ker\lambda_A[\fp^\infty]$ and $\Ker\lambda_B[\fp^\infty]$ are of rank $p^2$, respectively, then we have
          \[
          \rank_{\cO_S}(\Ker\alpha_{*,\tau_\infty})+\rank_{\cO_S}(\Ker\alpha_{*,\tau_\infty^\tc})=N.
          \]
      \end{enumerate}

  \item Let $\lambda_A$ and $\lambda_B$ be two quasi-polarizations on $A$ and $B$, respectively, such that $(A,\lambda_A)$ and $(B,\lambda_B)$ become unitary $O_F$-abelian schemes of dimension $N[F^+:\dQ]$ for some integer $N\geq 1$. Suppose that $\alpha^\vee\circ\lambda_B\circ\alpha=\lambda_A$. If $\Ker\lambda_A[\fp^\infty]$ is of rank $p^2$ and $\lambda_B$ is $p$-principal, then we have
      \[
      \rank_{\cO_S}(\Ker\alpha_{*,\tau_\infty})+\rank_{\cO_S}(\Ker\alpha_{*,\tau_\infty^\tc})=1.
      \]
\end{enumerate}

\end{lem}

\begin{proof}
We may assume $S$ connected. Up to replacing $\alpha$, $\beta$ and $\varpi$ by a common $\dZ_{(p)}^{\times}$-multiple, we may also assume that $\alpha$ and $\beta$ are genuine isogenies.

For (1), it suffices to show that the induced maps
\begin{align*}
\alpha_*&\colon\rH^\dr_1(A/S)\otimes_{O_{F^+}}\dZ_p\to\rH^\dr_1(B/S)\otimes_{O_{F^+}}\dZ_p,\\
\beta_*&\colon\rH^\dr_1(B/S)\otimes_{O_{F^+}}\dZ_p\to\rH^\dr_1(A/S)\otimes_{O_{F^+}}\dZ_p
\end{align*}
satisfy the relations $\Ker\alpha_*=\IM\beta_*$ and $\Ker\beta_*=\IM\alpha_*$; and these kernels and images are locally free $\cO_S$-modules.

Note that $A[\fp]$, $B[\fp]$, $\Ker\alpha[\fp]$, and $\Ker\beta[\fp]$ are all locally free finite group schemes over $S$ with an action by $O_F/\fp O_F$. By the relation among $\alpha,\beta,\varpi$, we may assume that $A[\fp]$ and $B[\fp]$ have degree $p^{2d}$; $\Ker\alpha[\fp]$ has degree $p^r$; and $\Ker\beta[\fp]$ has degree $p^{2d-r}$. As $\beta_*\circ\alpha_*=0$ and $\alpha_*\circ\beta_*=0$, it suffices to show that both $\Ker\alpha_*$ and $\IM\beta_*$ (resp.\ both $\Ker\beta_*$ and $\IM\alpha_*$) are locally direct factors of $\rH^\dr_1(A/S)\otimes_{O_{F^+}}\dZ_p$ (resp.\ $\rH^\dr_1(B/S)\otimes_{O_{F^+}}\dZ_p$) of rank $r$ (resp.\ $2d-r$), which will follow if we can show that $\coker\alpha_*$ and $\coker\beta_*$ are locally free $\cO_S$-modules of rank $r$ and $2d-r$, respectively.

We now prove that $\coker\alpha_*$ is a locally free $\cO_S$-modules of rank $r$; and the other case is similar. We follow the argument in \cite{dJ93}*{Lemma~2.3}. Consider the big crystalline site $(S/\dZ_p)_{\cris}$ with the structural sheaf $\cO_S^{\cris}$. Denote by $\cD(A[\fp^\infty])$ and $\cD(B[\fp^\infty])$ the covariant Dieudonn\'{e} crystals on $(S/\dZ_p)_{\cris}$ of $p$-divisible groups $A[\fp^\infty]$ and $B[\fp^\infty]$, respectively, which are locally free $\cO_S^{\cris}$-modules. We have a short exact sequence
\begin{align}\label{eq:inverse_isogeny1}
0 \to \alpha_*\cD(A[\fp^\infty])/\varpi\cD(B[\fp^\infty]) \to \cD(B[\fp^\infty])/\varpi\cD(B[\fp^\infty])
\to \cD(B[\fp^\infty])/\alpha_*\cD(A[\fp^\infty]) \to 0
\end{align}
and a surjective map
\begin{align}\label{eq:inverse_isogeny2}
\alpha_*\colon\cD(A[\fp^\infty])/\beta_*\cD(B[\fp^\infty])\to\alpha_*\cD(A[\fp^\infty])/\varpi\cD(B[\fp^\infty])
\end{align}
of $\cO_S^\cris$-modules. To show that $\coker\alpha_*$ is a locally free $\cO_S$-module of rank $r$, it suffices to show that $\cD(B[\fp^\infty])/\alpha_*\cD(A[\fp^\infty])$ is a locally free $\cO_S^\cris/p\cO_S^\cris$-module of rank $r$. By \cite{BBM82}*{Proposition~4.3.1}, $\cD(B[\fp^\infty])/\varpi\cD(B[\fp^\infty])$ is a locally free $\cO_S^\cris/p\cO_S^\cris$-module of rank $2d$. Thus, by \eqref{eq:inverse_isogeny1} and \eqref{eq:inverse_isogeny2}, it suffices to show that the $\cO_S^\cris/p\cO_S^\cris$-modules $\alpha_*\cD(A[\fp^\infty])/\varpi\cD(B[\fp^\infty])$ and $\cD(B[\fp^\infty])/\alpha_*\cD(A[\fp^\infty])$ are locally generated by $2d-r$ and $r$ sections, respectively. However, this can be easily checked using classical Dieudonn\'e modules after base change to geometric points of $S$. Thus, (1) is proved.

For (2), we know from (1) that both $\Ker\alpha_{*,\tau_\infty}$ and $\Ker\alpha_{*,\tau_\infty^\tc}$ are locally free $\cO_S$-modules. We may assume that $S=\Spec\kappa$ for a perfect field $\kappa$ containing $\dF_{p^2}$. Put $r\coloneqq\dim_\kappa\Lie_{A/\kappa,\tau_\infty}$ and $s\coloneqq\dim_\kappa\Lie_{B/\kappa,\tau_\infty}$. Then we have
\[
s=\dim_\kappa(\omega_{B^\vee/\kappa,\tau_\infty^\tc})=\dim_\kappa\frac{\tV\cD(B)_{\tau_\infty}}{p\cD(B)_{\tau_\infty^\tc}},\quad
r=\dim_\kappa(\omega_{A^\vee/\kappa,\tau_\infty^\tc})=\dim_\kappa\frac{\tV\cD(A)_{\tau_\infty}}{p\cD(A)_{\tau_\infty^\tc}}.
\]
Thus, we obtain
\begin{align}\label{eq:inverse_isogeny}
s-r=\dim_\kappa\frac{\tV\cD(B)_{\tau_\infty}}{p\cD(B)_{\tau_\infty^\tc}}-\dim_\kappa\frac{\tV\cD(A)_{\tau_\infty}}{p\cD(A)_{\tau_\infty^\tc}}.
\end{align}
Regarding $\cD(A)$ as a submodule of $\cD(B)$ via $\alpha_*$, it follows that
\begin{align*}
\eqref{eq:inverse_isogeny}&=\dim_\kappa\frac{\tV\cD(B)_{\tau_\infty}}{\tV\cD(A)_{\tau_\infty}}-
\dim_\kappa\frac{p\cD(B)_{\tau_\infty^\tc}}{p\cD(A)_{\tau_\infty^\tc}}
=\dim_\kappa\frac{\cD(B)_{\tau_\infty}}{\cD(A)_{\tau_\infty}}-
\dim_\kappa\frac{\cD(B)_{\tau_\infty^\tc}}{\cD(A)_{\tau_\infty^\tc}}\\
&=\dim_\kappa(\Ker\alpha_{*,\tau_\infty})-\dim_\kappa(\Ker\alpha_{*,\tau_\infty^\tc}).
\end{align*}
Thus, (2) is proved.

For (3) and (4), it suffices to show that $S=\Spec\kappa$ for an algebraically closed field $\kappa$ containing $\dF_{p^2}$. Then $\rank_{\cO_S}(\Ker\alpha_{*,\tau_\infty})+\rank_{\cO_S}(\Ker\alpha_{*,\tau_\infty^\tc})$ is half the length of the cokernel of the induced map $\alpha_*\colon\cD(A[\fp^\infty])\to\cD(B[\fp^\infty])$ regarded as $W(\kappa)$-modules. All cases follow immediately. In fact, for example, (3d) follows from the relation $2+2(\rank_{\cO_S}(\Ker\alpha_{*,\tau_\infty})+\rank_{\cO_S}(\Ker\alpha_{*,\tau_\infty^\tc}))=2N+2$; and others are similar.

\if false
For (3), consider the alternating pairing
\[
\langle\;,\;\rangle_{\lambda_A,\tau_\infty}\colon\rH^{\dr}_1(A/S)_{\tau_\infty}\times \rH^{\dr}_1(A/S)_{\tau^\tc_\infty}\to \cO_S
\]
induced by $\lambda_A$ as in Notation \ref{no:weil_pairing}.

In case (a), we have $\Ker\alpha_{*,\tau_\infty^\tc}=(\Ker\alpha_{*,\tau_\infty})^\perp$, since $\lambda_B$ is also $p$-principal.

In case (b), $\Ker\alpha_{*,\tau_\infty^\tc}$ is a subbundle of $(\Ker\alpha_{*,\tau_\infty})^\perp$ of corank $1$. The identity follows immediately from the relation $\rank_{\cO_S}(\Ker\alpha_{*,\tau_\infty})+\rank_{\cO_S}((\Ker\alpha_{*,\tau_\infty})^\perp)=N$.

In case (c), we apply case (b) to $\beta$ and use the relation
\[
\rank_{\cO_S}(\Ker\alpha_{*,\tau_\infty})+\rank_{\cO_S}(\Ker\beta_{*,\tau_\infty})
=\rank_{\cO_S}(\Ker\alpha_{*,\tau_\infty^\tc})+\rank_{\cO_S}(\Ker\beta_{*,\tau_\infty^\tc})
=N
\]
obtained from (1).

In case (d), we have both situations in (b) and (c), and the identity follows by a similar reason.

The proof for (4) is similar to (3). We leave the detail to the readers.
\fi
\end{proof}

\subsection{A CM moduli scheme}
\label{ss:cm_moduli}

In this subsection, we introduce an auxiliary moduli scheme parameterizing certain CM abelian varieties, which will be used in \S\ref{ss:qs} and \S\ref{ss:ns}.

\begin{definition}
Let $R$ be a $\dZ[(\disc F)^{-1}]$-ring.
\begin{enumerate}
  \item A \emph{rational skew-hermitian space} over $O_F\otimes R$ of rank $N$ is a free $O_F\otimes R$-module $\rW$ of rank $N$ together with an $R$-bilinear skew-symmetric perfect pairing
      \[
      \langle\;,\;\rangle_\rW\colon\rW\times\rW\to R
      \]
      satisfying $\langle ax,y\rangle_\rW=\langle x,a^\tc y\rangle_\rW$ for every $a\in O_F\otimes R$ and $x,y\in\rW$.

  \item Let $\rW$ and $\rW'$ be two rational skew-hermitian spaces over $O_F\otimes R$, a map $f\colon\rW\to\rW'$ is a \emph{similitude} if $f$ is an $O_F\otimes R$-linear isomorphism such that there exists some $c(f)\in R^\times$ satisfying $\langle f(x),f(y)\rangle_{\rW'}=c(f)\langle x,y\rangle_\rW$ for every $x,y\in\rW$.

  \item Two rational skew-hermitian spaces over $O_F\otimes R$ are \emph{similar} if there exists a similitude between them.

  \item For a rational skew-hermitian space $\rW$ over $O_F\otimes R$, we denote by $\GU(\rW)$ its \emph{group of similitude} as a reductive group over $R$; it satisfies that for every ring $R'$ over $R$, $\GU(\rW)(R')$ is the set of self-similitude of the rational skew-hermitian space $\rW\otimes_RR'$ over $O_F\otimes R'$.
\end{enumerate}
\end{definition}

We define a subtorus $\rT_0\subseteq(\Res_{O_F/\dZ}\bG_m)\otimes\dZ[(\disc F)^{-1}]$ such that for every $\dZ[(\disc F)^{-1}]$-ring $R$, we have
\[
\rT_0(R)=\{a\in O_F\otimes R\res \Nm_{F/{F^+}}a\in R^\times\}.
\]

Now we take a rational prime $p$ that is unramified in $F$. We take the prescribed subring $\dP$ in Definition \ref{de:unitary_abelian_scheme} to be $\dZ_{(p)}$.

\begin{remark}\label{re:hasse_principle}
Let $\rW_0$ be a rational skew-hermitian space over $O_F\otimes\dZ_{(p)}$ of rank $1$. Then $\GU(\rW_0)$ is canonically isomorphic to $\rT_0\otimes_{\dZ[(\disc F)^{-1}]}\dZ_{(p)}$. Moreover, the set of similarity classes of rational skew-hermitian spaces $\rW'_0$ over $O_F\otimes\dZ_{(p)}$ of rank $1$ such that $\rW'_0\otimes_{\dZ_{(p)}}\dA$ is similar to $\rW_0\otimes_{\dZ_{(p)}}\dA$ is canonically isomorphic to
\[
\Ker^1(\rT_0)\coloneqq\Ker\(\rH^1(\dQ,\rT_0)\to\prod_{v\leq\infty}\rH^1(\dQ_v,\rT_0)\),
\]
which is a finite abelian group.
\end{remark}

\begin{definition}\label{de:skew_hermitian_type}
Let $\Phi$ be a CM type. We say that a rational skew-hermitian space $\rW_0$ over $O_F\otimes\dZ_{(p)}$ of rank $1$ has \emph{type $\Phi$} if for every $x\in\rW_0$ and every totally imaginary element $a\in F^\times$ satisfying $\IP\tau(a)>0$ for all $\tau\in\Phi$, we have $\langle a x,x\rangle_{\rW_0}\geq 0$.
\end{definition}

\begin{definition}
For a rational skew-hermitian space $\rW_0$ over $O_F\otimes\dZ_{(p)}$ of rank $1$ and type $\Phi$ and an open compact subgroup $\rK_0^p\subseteq\rT_0(\dA^{\infty,p})$, we define a presheaf $\bT^1_p(\rW_0,\rK_0^p)$ on $\Sch'_{/O_{F_\Phi}\otimes\dZ_{(p)}}$ as follows: for every $S\in\Sch'_{/O_{F_\Phi}\otimes\dZ_{(p)}}$, we let $\bT^1_p(\rW_0,\rK_0^p)(S)$ be the set of equivalence classes of triples $(A_0,\lambda_0,\eta_0^p)$, where
\begin{itemize}[label={\ding{109}}]
  \item $(A_0,\lambda_0)$ is a unitary $O_F$-abelian scheme of signature type $\Phi$ over $S$ such that $\lambda_0$ is $p$-principal;

  \item $\eta_0^p$ is a $\rK_0^p$-level structure, that is, for a chosen geometric point $s$ on every connected component of $S$, a $\pi_1(S,s)$-invariant $\rK_0^p$-orbit of similitude
      \[
      \eta_0^p\colon\rW_0\otimes_\dQ\dA^{\infty,p}\to\rH^\et_1(A_{0s},\dA^{\infty,p})
      \]
      of rational skew-hermitian spaces over $F\otimes_\dQ\dA^{\infty,p}$, where $\rH^\et_1(A_{0s},\dA^{\infty,p})$ is equipped with the rational skew-hermitian form induced by $\lambda_0$.
\end{itemize}
Two triples $(A_0,\lambda_0,\eta_0^p)$ and $(A'_0,\lambda'_0,\eta_0^{p\prime})$ are equivalent if there exists a prime-to-$p$ $O_F$-linear quasi-isogeny $\varphi_0\colon A_0\to A'_0$ carrying $(\lambda_0,\eta_0^p)$ to $(c\lambda'_0,\eta_0^{p\prime})$ for some $c\in\dZ_{(p)}^\times$.
\end{definition}

For an object $(A_0,\lambda_0,\eta_0^p)\in\bT^1_p(\rW_0,\rK_0^p)(\dC)$, its first homology $\rH_1(A_0(\dC),\dZ_{(p)})$ is a rational skew-hermitian space over $O_F\otimes\dZ_{(p)}$ induced by $\lambda_0$, which is of rank $1$ and type $\Phi$, and is everywhere locally similar to $\rW_0$. Thus, by Remark \ref{re:hasse_principle}, we obtain a map
\[
\tw\colon\bT^1_p(\rW_0,\rK_0^p)(\dC)\to\Ker^1(\rT_0)
\]
sending $(A_0,\lambda_0,\eta_0^p)\in\bT^1_p(\rW_0,\rK_0^p)(\dC)$ to the similarity class of $\rH_1(A_0(\dC),\dZ_{(p)})$.

It is known that when $\rK_0^p$ is neat, $\bT^1_p(\rW_0,\rK_0^p)$ is represented by a scheme finite and \'{e}tale over $O_{F_\Phi}\otimes\dZ_{(p)}$. We define $\bT_p(\rW_0,\rK_0^p)$ to be the minimal open and closed subscheme of $\bT^1_p(\rW_0,\rK_0^p)$ containing $\tw^{-1}(\rW_0)$. The group $\rT_0(\dA^{\infty,p})$ acts on $\bT_p(\rW_0,\rK_0^p)$ via the formula
\[
a\cdot(A_0,\lambda_0,\eta_0^p)=(A_0,\lambda_0,\eta_0^p\circ a)
\]
whose stabilizer is $\rT_0(\dZ_{(p)})\rK_0^p$. In fact, $\rT_0(\dA^{\infty,p})/\rT_0(\dZ_{(p)})\rK_0^p$ is the Galois group of the Galois morphism
\[
\bT_p(\rW_0,\rK_0^p)\to\Spec (O_{F_\Phi}\otimes\dZ_{(p)}).
\]

\begin{definition}\label{de:groupoid}
We denote by $\fT$ the groupoid of $\rT_0(\dA^{\infty,p})/\rT_0(\dZ_{(p)})\rK_0^p$, that is, a category with a single object $\ast$ with $\Hom(\ast,\ast)=\rT_0(\dA^{\infty,p})/\rT_0(\dZ_{(p)})\rK_0^p$.
\end{definition}

\begin{remark}\label{re:groupoid}
As $\bT_p(\rW_0,\rK_0^p)$ is an object in $\Sch_{/O_{F_\Phi}\otimes\dZ_{(p)}}$ with an action by $\rT_0(\dA^{\infty,p})/\rT_0(\dZ_{(p)})\rK_0^p$, it induces a functor from $\fT$ to $\Sch_{/O_{F_\Phi}\otimes\dZ_{(p)}}$, which we still denote by $\bT_p(\rW_0,\rK_0^p)$. In what follows, we may often have another category $\fC$ and will regard $\bT_p(\rW_0,\rK_0^p)$ as a functor from $\fC\times\fT$ to $\Sch_{/O_{F_\Phi}\otimes\dZ_{(p)}}$ as the composition of the projection functor $\fC\times\fT\to\fT$ and the functor $\bT_p(\rW_0,\rK_0^p)\colon\fT\to\Sch_{/O_{F_\Phi}\otimes\dZ_{(p)}}$.
\end{remark}

\begin{notation}\label{no:groupoid}
For a functor $X\colon\fT\to\Sch$ and a coefficient ring $L$, we denote
\[
\rH^i_\fT(X,L(j))\subseteq\rH^i_\et(X(\ast),L(j)),\quad
\rH^i_{\fT,c}(X,L(j))\subseteq\rH^i_{\et,c}(X(\ast),L(j))
\]
the maximal $L$-submodules, respectively, on which $\rT_0(\dA^{\infty,p})/\rT_0(\dZ_{(p)})\rK_0^p$ acts trivially.
\end{notation}

\begin{definition}\label{de:trace}
Let $\kappa$ be an algebraically closed field of characteristic $p$, and $L$ a $p$-coprime coefficient ring. For a functor $X\colon\fT\to\Sch_{/\kappa}$ such that $X(\ast)$ is smooth of finite type of dimension $d$ and that $\fT$ acts freely on the set of connected components of $X(\ast)$, we define the \emph{$\fT$-trace map}
\[
\int_X^\fT\colon\rH^{2d}_{\fT,c}(X(\ast),L(d))\to L
\]
to be the composite map
\[
\rH^{2d}_{\fT,c}(X(\ast),L(d))\hookrightarrow\rH^{2d}_c(X(\ast),L(d))\to\bigoplus_Y\rH^{2d}_c(Y,L(d))\xrightarrow{\sum\tr_Y}L,
\]
where $\{Y\}$ is a set of representatives of $\fT$-orbits on the connected components of $X(\ast)$, and the second map is the natural projection. It is clear that the above composite map does not depend on the choice of $\{Y\}$.
\end{definition}

\section{Unitary moduli schemes: smooth case}
\label{ss:qs}

In this section, we define and study a certain smooth integral moduli scheme whose generic fiber is the product of a unitary Shimura variety and an auxiliary CM moduli. Since the materials in this section are strictly in the linear order, we will leave the summary of contents to each subsection.

\subsection{Initial setup}
\label{ss:qs_initial}

We fix a special inert prime (Definition \ref{de:special_inert}) $\fp$ of $F^+$ (with the underlying rational prime $p$). We take the prescribed subring $\dP$ in Definition \ref{de:unitary_abelian_scheme} to be $\dZ_{(p)}$. We choose the following data
\begin{itemize}[label={\ding{109}}]
  \item a CM type $\Phi$ containing $\tau_\infty$;

  \item a rational skew-hermitian space $\rW_0$ over $O_F\otimes\dZ_{(p)}$ of rank $1$ and type $\Phi$ (Definition \ref{de:skew_hermitian_type});

  \item a neat open compact subgroup $\rK_0^p\subseteq\rT_0(\dA^{\infty,p})$;

  \item an isomorphism $\iota_p\colon\dC\simeq \ol\dQ_p$ such  that $\iota_p\circ\ul\tau_\infty\colon F^+\hookrightarrow \ol\dQ_p$ induces the place $\fp$ of $F^+$;

  \item an element $\varpi\in O_{F^+}$ that is totally positive and satisfies $\val_\fp(\varpi)=1$, and $\val_\fq(\varpi)=0$ for every prime $\fq\neq\fp$ of $F^+$ above $p$.
\end{itemize}
We adopt Notation \ref{no:p_notation}. In particular, $\dF^\Phi_p$ contains $\dF_{p^2}$. Since the argument below is insensitive to the choices of $\rW_0$ and $\rK_0^p$, we will not include them in all notations. However, we will keep the prime $\fp$ in notations as, in later application, we need to choose different primes in a crucial step. Put $\bT_\fp\coloneqq\bT_p(\rW_0,\rK_0^p)\otimes_{O_{F_\Phi}\otimes\dZ_{(p)}}\dZ^\Phi_p$.

\subsection{Construction of moduli schemes}
\label{ss:qs_moduli_scheme}

In this subsection, we construct our initial moduli schemes. We start from the datum $(\rV,\{\Lambda_\fq\}_{\fq\mid p})$, where
\begin{itemize}[label={\ding{109}}]
  \item $\rV$ is a standard \emph{indefinite} hermitian space (Definition \ref{de:standard_hermitian_space}) over $F$ of rank $N\geq 1$, and

  \item $\Lambda_\fq$ is a self-dual $O_{F_\fq}$-lattice in $\rV\otimes_FF_\fq$ for every prime $\fq$ of $F^+$ above $p$.
\end{itemize}

Before defining the moduli functor, we need the following lemma to make sense of the later definition.

\begin{lem}
The field $\dQ^\Phi_p$ contains $F_\Psi$ with $\Psi=N\Phi-\tau_\infty+\tau_\infty^\tc$, which is a generalized CM type of rank $N$, for every $N\geq 1$.
\end{lem}

\begin{proof}
Take $\rho\in\Aut(\dC/\dQ^\Phi_p)\subseteq\Aut(\dC/F)$. Then we have $\rho\Phi=\Phi$ and $\rho\tau_\infty=\tau_\infty$. Thus, we have $\rho(N\Phi-\tau_\infty+\tau_\infty^\tc)=N\Phi-\tau_\infty+\tau_\infty^\tc$ for every $N\geq 1$. The lemma follows.
\end{proof}

Recall that we have the category $\Sch'_{/\dZ^\Phi_p}$ of locally Noetherian schemes over $\dZ^\Phi_p$, and $\sfP\Sch'_{/\dZ^\Phi_p}$ the category of presheaves on $\Sch'_{/\dZ^\Phi_p}$.

\begin{definition}\label{de:qs_moduli_scheme}
We define a functor
\begin{align*}
\bM_\fp(\rV,\obj)\colon\fK(\rV)^p\times\fT &\to\sfP\Sch'_{/\dZ^\Phi_p} \\
\rK^p &\mapsto \bM_\fp(\rV,\rK^p)
\end{align*}
such that for every $S\in\Sch'_{/\dZ^\Phi_p}$, $\bM_\fp(\rV,\rK^p)(S)$ is the set of equivalence classes of sextuples $(A_0,\lambda_0,\eta_0^p;A,\lambda,\eta^p)$, where
\begin{itemize}[label={\ding{109}}]
  \item $(A_0,\lambda_0,\eta_0^p)$ is an element in $\bT_\fp(S)$;

  \item $(A,\lambda)$ is a unitary $O_F$-abelian scheme of signature type $N\Phi-\tau_\infty+\tau_\infty^\tc$ over $S$ (Definitions \ref{de:unitary_abelian_scheme} and \ref{de:signature}) such that $\lambda$ is $p$-principal;

  \item $\eta^p$ is a $\rK^p$-level structure, that is, for a chosen geometric point $s$ on every connected component of $S$, a $\pi_1(S,s)$-invariant $\rK^p$-orbit of isomorphisms
      \[
      \eta^p\colon\rV\otimes_\dQ\dA^{\infty,p}\to
      \Hom_{F\otimes_\dQ\dA^{\infty,p}}^{\lambda_0,\lambda}(\rH^\et_1(A_{0s},\dA^{\infty,p}),\rH^\et_1(A_s,\dA^{\infty,p}))
      \]
      of hermitian spaces over $F\otimes_\dQ\dA^{\infty,p}=F\otimes_{F^+}\dA_{F^+}^{\infty,p}$. See Construction \ref{cs:hermitian_structure} (with $\Box=\{\infty,p\}$) for the right-hand side.
\end{itemize}
Two sextuples $(A_0,\lambda_0,\eta_0^p;A,\lambda,\eta^p)$ and $(A'_0,\lambda'_0,\eta_0^{p\prime};A',\lambda',\eta^{p\prime})$ are equivalent if there are prime-to-$p$ $O_F$-linear quasi-isogenies $\varphi_0\colon A_0\to A'_0$ and $\varphi\colon A\to A'$ such that
\begin{itemize}[label={\ding{109}}]
  \item $\varphi_0$ carries $\eta_0^p$ to $\eta_0^{p\prime}$;

  \item there exists $c\in\dZ_{(p)}^\times$ such that $\varphi_0^\vee\circ\lambda'_0\circ\varphi_0=c\lambda_0$ and $\varphi^\vee\circ\lambda'\circ\varphi=c\lambda$; and

  \item the $\rK^p$-orbit of maps $v\mapsto\varphi_*\circ\eta^p(v)\circ(\varphi_{0*})^{-1}$ for $v\in\rV\otimes_\dQ\dA^{\infty,p}$ coincides with $\eta^{p\prime}$.
\end{itemize}
On the level of morphisms,
\begin{itemize}[label={\ding{109}}]
  \item a morphism $g\in\rK^p\backslash\rU(\rV)(\dA_F^{\infty,p})/\rK^{p\prime}$ of $\fK(\rV)^p$ maps $\bM_\fp(\rV,\rK^p)(S)$ to $\bM_\fp(\rV,\rK^{p\prime})(S)$ by changing $\eta^p$ to $\eta^p\circ g$; and

  \item a morphism $a$ of $\fT$ acts on $\bM_\fp(\rV,\rK^p)(S)$ by changing $\eta_0^p$ to $\eta_0^p\circ a$.
\end{itemize}
\end{definition}

We clearly have the forgetful morphism
\begin{align}\label{eq:qs_moduli_scheme}
\bM_\fp(\rV,\obj)\to\bT_\fp
\end{align}
in $\Fun(\fK(\rV)^p\times\fT,\sfP\Sch'_{/\dZ_p^\Phi})$, the category of functors from $\fK(\rV)^p\times\fT$ to $\sfP\Sch'_{/\dZ_p^\Phi}$. Here, we regard $\bT_\fp$ as an object in $\Fun(\fK(\rV)^p\times\fT,\Sch'_{/\dZ_p^\Phi})$ as in Remark \ref{re:groupoid}. According to Notation \ref{no:p_notation}, we shall denote by the base change of \eqref{eq:qs_moduli_scheme} to $\dF_p^\Phi$ by $\rM_\fp(\rV,\obj)\to\rT_\fp$, which is a morphism in $\Fun(\fK(\rV)^p\times\fT,\sfP\Sch'_{/\dF_p^\Phi})$.

\begin{theorem}\label{th:qs_moduli_scheme}
The morphism \eqref{eq:qs_moduli_scheme} is represented by a quasi-projective smooth scheme over $\bT_{\fp}$ of relative dimension $N-1$. Moreover, for every $\rK^p\in\fK(\rV)^p$, we have a canonical isomorphism
\[
\cT_{\bM_\fp(\rV,\rK^p)/\bT_\fp}\simeq\HOM\(\omega_{\cA^\vee,\tau_\infty},\rH^\dr_1(\cA)_{\tau_\infty}/\omega_{\cA^\vee,\tau_\infty}\)
\]
of coherent sheaves on $\bM_\fp(\rV,\rK^p)$, where $(\cA_0,\lambda_0,\eta_0^p;\cA,\lambda,\eta^p)$ is the universal object over $\bM_\fp(\rV,\rK^p)$ and we recall that $\cT_{\bM_\fp(\rV,\rK^p)/\bT_\fp}$ is the relative tangent sheaf. Moreover, \eqref{eq:qs_moduli_scheme} is projective if and only if its base change to $\dQ_p^\Phi$ is.
\end{theorem}

\begin{proof}
The first claim is proved in \cite{RSZ}*{Theorem~4.4}. It remains to compute the tangent sheaf. Take an object $\rK^p\in\fK(\rV)^p$. Since both $\rK^p_0$ and $\rK^p$ are neat, $\bM_\fp(\rV,\rK^p)$ is an algebraic space. Thus, we have the universal object $(\cA_0,\lambda_0,\eta_0^p;\cA,\lambda,\eta^p)$ over $\bM_\fp(\rV,\rK^p)$. By a standard argument in deformation theory, using Proposition \ref{pr:deformation}, we know that the morphism $\bM_\fp(\rV,\rK^p)\to\bT_\fp$ is separated and smooth; and we have a canonical isomorphism for the tangent sheaf
\[
\cT_{\bM_\fp(\rV,\rK^p)/\bT_\fp}\simeq\HOM\(\omega_{\cA^\vee,\tau_\infty},\rH^\dr_1(\cA)_{\tau_\infty}/\omega_{\cA^\vee,\tau_\infty}\)
\]
which is locally free of rank $N-1$. The theorem is proved.
\end{proof}

Let $\rK_\fq$ be the stabilizer of $\Lambda_\fq$ for every $\fq\mid p$; and put $\rK_p\coloneqq\prod_{\fq\mid p}\rK_\fq$. As shown in \cite{RSZ}*{\S3.3}, there is a \emph{canonical} ``moduli interpretation'' isomorphism of varieties over $\dQ_p^{\Phi}$
\begin{align}\label{eq:qs_moduli_scheme_shimura}
\bM^\eta_\fp(\rV,\obj)\xrightarrow{\sim}\Sh(\rV,\obj\rK_p)\times_{\Spec{F}}\bT^\eta_\fp
\end{align}
(Notation \ref{no:p_notation}(5)) in $\Fun(\fK(\rV)^p\times\fT,\Sch_{/\dQ_p^\Phi})_{/\bT^\eta_\fp}$, where $\fT$ acts on $\Sh(\rV,\obj\rK_p)\times_{\Spec{F}}\bT^\eta_\fp$ through the second factor. See also Remark \ref{re:qs_moduli_scheme_shimura} below.

\begin{lem}\label{le:qs_pbc}
Let $L$ be a $p$-coprime coefficient ring. The two specialization maps
\begin{align*}
\rH^i_{\fT,c}(\bM_\fp(\rV,\obj)\otimes_{\dZ_p^\Phi}\ol\dQ_p,L)&\to\rH^i_{\fT,c}(\ol\rM_\fp(\rV,\obj),L), \\
\rH^i_\fT(\bM_\fp(\rV,\obj)\otimes_{\dZ_p^\Phi}\ol\dQ_p,L)&\to\rH^i_\fT(\ol\rM_\fp(\rV,\obj),L),
\end{align*}
are both isomorphisms. In particular, \eqref{eq:qs_moduli_scheme_shimura} induces isomorphisms
\begin{align*}
\rH^i_{\et,c}(\Sh(\rV,\obj\rK_p)_{\ol{F}},L)&\simeq\rH^i_{\fT,c}(\ol\rM_\fp(\rV,\obj),L), \\
\rH^i_\et(\Sh(\rV,\obj\rK_p)_{\ol{F}},L)&\simeq\rH^i_\fT(\ol\rM_\fp(\rV,\obj),L),
\end{align*}
in $\Fun(\fK(\rV)^p,\Mod(L[\Gal(\ol\dQ_p/\dQ^\Phi_p)]))$ for every $i\in\dZ$. Here, $\Gal(\ol\dQ_p/\dQ^\Phi_p)$ is regarded as a subgroup of $\Gal(\ol{F}/F)$ under our fixed isomorphism $\iota_p\colon\dC\simeq\ol\dQ_p$.
\end{lem}

\begin{proof}
Since $\bM_\fp(\rV,\obj)$ is smooth over $\dZ_p^\Phi$, we have a canonical isomorphism $L\simeq\rR\Psi L$. When $\bM_\fp(\rV,\obj)$ is proper, this is simply the proper base change. When $\bM_\fp(\rV,\obj)$ is not proper, this follows from \cite{LS18}*{Corollary~5.20}.
\end{proof}

\begin{remark}\label{re:qs_moduli_scheme_shimura}
For the readers' convenience, we describe the isomorphism \eqref{eq:qs_moduli_scheme_shimura} on complex points, which determines the isomorphism uniquely. It suffices to assign to every point
\[
x=(A_0,\lambda_0,\eta_0^p;A,\lambda,\eta^p)\in\bM_\fp(\rV,\rK^p)(\dC)
\]
 a point in
\[
\Sh(\rV,\rK^p\rK_p)(\dC)=\rU(\rV)(F^+)\backslash\(\rV(\dC)_-/\dC^\times\times\rU(\rV)(\dA_{F^+}^\infty)/\rK^p\rK_p\)
\]
where $\rV(\dC)_-/\dC^\times$ is the set of negative definite complex lines in $\rV\otimes_F\dC$. Put
\[
\rV_x\coloneqq\Hom_F(\rH_1(A_0(\dC),\dQ),\rH_1(A(\dC),\dQ))
\]
equipped with a pairing in the way similar to Construction \ref{cs:hermitian_structure}, which becomes a hermitian space over $F$ of rank $N$. Moreover, it is standard indefinite. By the comparison between singular homology and \'etale homology, we have a canonical isometry of hermitian spaces
\[
\rho\colon\rV_x\otimes_\dQ\dA^{\infty,p}\xrightarrow{\sim}
\Hom_{F\otimes_\dQ\dA^{\infty,p}}^{\lambda_0,\lambda}(\rH^\et_1(A_0,\dA^{\infty,p}),\rH^\et_1(A,\dA^{\infty,p})),
\]
which implies that $\rV_x\otimes_\dQ\dA^{\infty,p}\simeq\rV\otimes_\dQ\dA^{\infty,p}$ by the existence of the level structure $\eta^p$. On the other hand, we have a canonical decomposition
\[
\Hom_{O_F\otimes\dZ_p}(\rH^\et_1(A_0,\dZ_p),\rH^\et_1(A,\dZ_p))=\bigoplus_{\fq\mid p}\Lambda_{x,\fq}
\]
of $O_F\otimes\dZ_p$-modules in which $\Lambda_{x,\fq}$ is a self-dual lattice in $\rV\otimes_FF_\fq$ for every prime $\fq$ of $F^+$ above $p$. Thus, by the Hasse principle for hermitian spaces, this implies that hermitian spaces $\rV_x$ and $\rV$ are isomorphic. Choose an isometry $\eta_\rat\colon\rV_x\to\rV$. Thus, we obtain an isometry
\[
g^p\coloneqq \eta_\rat\circ\rho^{-1}\circ\eta^p\colon\rV\otimes_\dQ\dA^{\infty,p}\to\rV\otimes_\dQ\dA^{\infty,p}
\]
as an element in $\rU(\rV)(\dA_{F^+}^{\infty,p})$. For every $\fq$ above $p$, there exists an element $g_\fq\in\rU(\rV)(F^+_\fq)$ such that $g_\fq\Lambda_\fq=\eta_\rat\Lambda_{x,\fq}$. Together, we obtain an element $g_x\coloneqq(g^p,(g_\fq)_{\fq\mid p})\in\rU(\rV)(\dA_{F^+}^\infty)$. Finally,
\[
l_x\coloneqq\{\alpha\in\Hom_F(\rH_1^\dr(A_0/\dC),\rH_1^\dr(A/\dC))\res
\alpha(\omega_{A_0^\vee,\tau_\infty})\subseteq\omega_{A^\vee,\tau_\infty}\}
\]
is a line in $\rV_x(\dC)$ such that $\eta_\rat(l_x)$ is an element in $\rV(\dC)_-/\dC^\times$. It is easy to check that the coset
\[
\rU(\rV)(F^+)(\eta_\rat(l_x),g_x\rK^p\rK_p)
\]
does not depend on the choice of $\eta_\rat$, hence gives rise an element in $\Sh(\rV,\rK^p\rK_p)(\dC)$. It is clear that the action of a morphism $a$ of $\fT$ on $x$ does not change the above coset.
\end{remark}

\subsection{Basic correspondence for the special fiber}
\label{ss:qs_basic_correspondence}

In this subsection, we construct and study the basic correspondence for the special fiber $\rM_\fp(\rV,\obj)$. Recall that we have chosen an element $\varpi\in O_{F^+}$ that is totally positive and satisfies $\val_\fp(\varpi)=1$, and $\val_\fq(\varpi)=0$ for every prime $\fq\neq\fp$ of $F^+$ above $p$.

\begin{definition}\label{de:qs_definite}
We define a functor
\begin{align*}
\rS_\fp(\rV,\obj)\colon\fK(\rV)^p\times\fT &\to\sfP\Sch'_{/\dF^\Phi_p} \\
\rK^p &\mapsto \rS_\fp(\rV,\rK^p)
\end{align*}
such that for every $S\in\Sch'_{/\dF^\Phi_p}$, $\rS_\fp(\rV,\rK^p)(S)$ is the set of equivalence classes of sextuples $(A_0,\lambda_0,\eta_0^p;A^\star,\lambda^\star,\eta^{p\star})$, where
\begin{itemize}[label={\ding{109}}]
  \item $(A_0,\lambda_0,\eta_0^p)$ is an element in $\rT_\fp(S)$;

  \item $(A^\star,\lambda^\star)$ is a unitary $O_F$-abelian scheme of signature type $N\Phi$ over $S$ such that $\Ker\lambda^\star[p^\infty]$ is trivial (resp.\ contained in $A^\star[\fp]$ of rank $p^2$) if $N$ is odd (resp.\ even);

  \item $\eta^{p\star}$ is, for a chosen geometric point $s$ on every connected component of $S$, a $\pi_1(S,s)$-invariant $\rK^p$-orbit of isomorphisms
      \[
      \eta^{p\star}\colon\rV\otimes_\dQ\dA^{\infty,p}\to\Hom_{F\otimes_\dQ\dA^{\infty,p}}^{\varpi\lambda_0,\lambda^\star}
      (\rH^\et_1(A_{0s},\dA^{\infty,p}),\rH^\et_1(A^\star_s,\dA^{\infty,p}))
      \]
      of hermitian spaces over $F\otimes_\dQ\dA^{\infty,p}=F\otimes_{F^+}\dA_{F^+}^{\infty,p}$.\footnote{Note that here we are using $\varpi\lambda_0$ rather than $\lambda_0$ in order to be consistent with the compatibility condition for polarizations in the isogeny considered in Definition \ref{de:qs_basic_correspondence}.}
\end{itemize}
The equivalence relation and the action of morphisms in $\fK(\rV)^p\times\fT$ are defined similarly as in Definition \ref{de:qs_moduli_scheme}.
\end{definition}

We clearly have the forgetful morphism
\begin{align*}
\rS_\fp(\rV,\obj)\to\rT_\fp
\end{align*}
in $\Fun(\fK(\rV)^p\times\fT,\sfP\Sch'_{/\dF_p^\Phi})$, which is represented by finite and \'{e}tale schemes by \cite{RSZ}*{Theorem~4.4}.

Now we take a point $s^\star=(A_0,\lambda_0,\eta_0^p;A^\star,\lambda^\star,\eta^{p\star})\in\rS_\fp(\rV,\rK^p)(\kappa)$ where $\kappa$ is a field containing $\dF^\Phi_p$. Then $A^\star_{\ol\kappa}[\fp^\infty]$ is a supersingular $p$-divisible group by the signature condition and the fact that $\fp$ is inert in $F$. From Notation \ref{re:frobenius_verschiebung}, we have the $(\kappa,\sigma)$-linear Frobenius map
\[
\tF\colon\rH^\dr_1(A^\star/\kappa)_{\tau_\infty}\to\rH^\dr_1(A^\star/\kappa)_{\sigma\tau_\infty}
=\rH^\dr_1(A^\star/\kappa)_{\tau_\infty^\tc}.
\]
We define a pairing
\[
\{\;,\;\}_{s^\star}\colon\rH^\dr_1(A^\star/\kappa)_{\tau_\infty}\times\rH^\dr_1(A^\star/\kappa)_{\tau_\infty}\to\kappa
\]
by the formula $\{x,y\}_{s^\star}\coloneqq\langle\tF x,y\rangle_{\lambda^\star,\tau_\infty^\tc}$ (Notation \ref{no:weil_pairing}). To ease notation, we put
\[
\sV_{s^\star}\coloneqq\rH^\dr_1(A^\star/\kappa)_{\tau_\infty}.
\]

\begin{lem}\label{le:qs_dl}
The pair $(\sV_{s^\star},\{\;,\;\}_{s^\star})$ is admissible of rank $N$ (Definition \ref{de:dl_admissible}). In particular, the Deligne--Lusztig variety $\DL_{s^\star}\coloneqq\DL(\sV_{s^\star},\{\;,\;\}_{s^\star},\ceil{\tfrac{N+1}{2}})$ (Definition \ref{de:dl}) is a geometrically irreducible projective smooth scheme in $\Sch_{/\kappa}$ of dimension $\floor{\tfrac{N-1}{2}}$ with a canonical isomorphism for its tangent sheaf
\[
\cT_{\DL_{s^\star}/\kappa}\simeq\HOM\(\cH/\cH^\dashv,(\sV_{s^\star})_{\DL_{s^\star}}/\cH\)
\]
where $\cH\subseteq(\sV_{s^\star})_{\DL_{s^\star}}$ is the universal subbundle.
\end{lem}

\begin{proof}
It follows from the construction that $\{\;,\;\}_{s^\star}$ is $(\kappa,\sigma)$-linear in the first variable and $\kappa$-linear in the second variable. By the signature condition Definition \ref{de:qs_definite}(2), the map $\tF\colon\rH^\dr_1(A^\star/\kappa)_{\tau_\infty}\to
\rH^\dr_1(A^\star/\kappa)_{\tau_\infty^\tc}$ is an isomorphism, and the pairing $\langle\tF\;,\;\rangle_{\lambda^\star,\tau_\infty^\tc}$ has kernel of rank $0$ (resp.\ $1$) if $N$ is odd (resp.\ even). Thus, by Proposition \ref{pr:dl}, it suffices to show that $(\sV_{s^\star},\{\;,\;\}_{s^\star})$ is admissible.

Note that we have a canonical isomorphism $(\sV_{s^\star})_{\ol\kappa}=\rH^\dr_1(A^\star/\kappa)_{\tau_\infty}\otimes_\kappa\ol\kappa
\simeq\rH^\dr_1(A^\star_{\ol \kappa}/\ol\kappa)_{\tau_\infty}$, and that the $(\ol\kappa,\sigma)$-linear Frobenius map $\tF\colon\rH^\dr_1(A^\star_{\ol \kappa}/\ol\kappa)_{\tau_\infty}\to\rH^\dr_1(A^\star_{\ol \kappa}/\ol\kappa)_{\tau_\infty^\tc}$ and the $(\ol\kappa,\sigma^{-1})$-linear Verschiebung map $\tV\colon\rH^\dr_1(A^\star_{\ol \kappa}/\ol\kappa)_{\tau_\infty}\to\rH^\dr_1(A^\star_{\ol \kappa}/\ol\kappa)_{\tau_\infty^\tc}$ are both isomorphisms. Thus, we obtain a $(\ol\kappa,\sigma^2)$-linear isomorphism $\tV^{-1}\tF\colon\rH^\dr_1(A^\star_{\ol\kappa}/\ol\kappa)_{\tau_\infty}\to\rH^\dr_1(A^\star_{\ol\kappa}/\ol\kappa)_{\tau_\infty}$.
Denote by $\sV_0$ the subset of $\rH^\dr_1(A^\star_{\ol \kappa}/\ol\kappa)_{\tau_\infty}$ on which $\tV^{-1}\tF=\id$, which is an $\dF_{p^2}$-linear subspace. Since the $p$-divisible group $A^\star_{\ol\kappa}[\fp^\infty]$ is supersingular, by Dieudonn\'{e}'s classification of crystals, the canonical map $\sV_0\otimes_{\dF_{p^2}}\ol\kappa\to\rH^\dr_1(A^\star/\ol\kappa)_{\tau_\infty}=(\sV_{s^\star})_{\ol\kappa}$ is an isomorphism. For $x,y\in\sV_0$, we have
\[
\{x,y\}_{s^\star}=\langle\tF x,y\rangle_{\lambda^\star,\tau_\infty^\tc}=
\langle x,\tV y\rangle_{\lambda^\star,\tau_\infty}^\sigma
=\langle x,\tF y\rangle_{\lambda^\star,\tau_\infty}^\sigma
=-\langle\tF y,x\rangle_{\lambda^\star,\tau_\infty^\tc}^\sigma
=-\{y,x\}_{s^\star}^\sigma.
\]
Thus, $(\sV_{s^\star},\{\;,\;\}_{s^\star})$ is admissible. The lemma follows.
\end{proof}

\begin{definition}\label{de:qs_basic_correspondence}
We define a functor
\begin{align*}
\rB_\fp(\rV,\obj)\colon\fK(\rV)^p\times\fT &\to\sfP\Sch'_{/\dF^\Phi_p} \\
\rK^p &\mapsto \rB_\fp(\rV,\rK^p)
\end{align*}
such that for every $S\in\Sch'_{/\dF^\Phi_p}$, $\rB_\fp(\rV,\rK^p)(S)$ is the set of equivalence classes of decuples $(A_0,\lambda_0,\eta_0^p;A,\lambda,\eta^p;A^\star,\lambda^\star,\eta^{p\star};\alpha)$, where
\begin{itemize}[label={\ding{109}}]
  \item $(A_0,\lambda_0,\eta_0^p;A,\lambda,\eta^p)$ is an element of $\rM_\fp(\rV,\rK^p)(S)$;

  \item $(A_0,\lambda_0,\eta_0^p;A^\star,\lambda^\star,\eta^{p\star})$ is an element of $\rS_\fp(\rV,\rK^p)(S)$; and

  \item $\alpha\colon A\to A^\star$ is an $O_F$-linear quasi-$p$-isogeny (Definition \ref{de:p_quasi}) such that
  \begin{enumerate}[label=(\alph*)]
    \item $\Ker\alpha[p^\infty]$ is contained in $A[\fp]$;

    \item we have $\varpi\cdot\lambda=\alpha^\vee\circ\lambda^\star\circ\alpha$; and

    \item the $\rK^p$-orbit of maps $v\mapsto\alpha_*\circ\eta^p(v)$ for $v\in\rV\otimes_\dQ\dA^{\infty,p}$ coincides with $\eta^{p\star}$.
  \end{enumerate}
\end{itemize}
Two decuples $(A_0,\lambda_0,\eta_0^p;A,\lambda,\eta^p;A^\star,\lambda^\star,\eta^{p\star};\alpha)$ and $(A'_0,\lambda'_0,\eta_0^{p\prime};A',\lambda',\eta^{p\prime};A^{\star\prime},\lambda^{\star\prime},\eta^{p\star\prime};\alpha^\prime)$ are equivalent if there are prime-to-$p$ $O_F$-linear quasi-isogenies $\varphi_0\colon A_0\to A'_0$, $\varphi\colon A\to A'$, and $\varphi^\star\colon A^\star\to A^{\star\prime}$ such that
\begin{itemize}[label={\ding{109}}]
  \item $\varphi_0$ carries $\eta_0^p$ to $\eta_0^{p\prime}$;

  \item there exists $c\in\dZ_{(p)}^\times$ such that $\varphi_0^\vee\circ\lambda'_0\circ\varphi_0=c\lambda_0$, $\varphi^\vee\circ\lambda'\circ\varphi=c\lambda$, and $\varphi^{\star\vee}\circ\lambda^{\star\prime}\circ\varphi^\star=c\lambda^\star$;

  \item the $\rK^p$-orbit of maps $v\mapsto\varphi_*\circ\eta^p(v)\circ(\varphi_{0*})^{-1}$ for $v\in\rV\otimes_\dQ\dA^{\infty,p}$ coincides with $\eta^{p\prime}$;

  \item the $\rK^p$-orbit of maps $v\mapsto\varphi^\star_*\circ\eta^{p\star}(v)\circ(\varphi_{0*})^{-1}$ for $v\in\rV\otimes_\dQ\dA^{\infty,p}$ coincides with $\eta^{p\star\prime}$;

  \item $\varphi^\star\circ\alpha=\alpha'\circ\varphi$ holds.
\end{itemize}
On the level of morphisms,
\begin{itemize}[label={\ding{109}}]
  \item a morphism $g\in\rK^p\backslash\rU(\rV)(\dA_F^{\infty,p})/\rK^{p\prime}$ of $\fK(\rV)^p$ maps $\rB_\fp(\rV,\rK^p)(S)$ to $\rB_\fp(\rV,\rK^{p\prime})(S)$ by changing $\eta^p,\eta^{p\star}$ to $\eta^p\circ g,\eta^{p\star}\circ g$, respectively; and

  \item a morphism $a$ of $\fT$ acts on $\rM_\fp(\rV,\rK^p)(S)$ by changing $\eta_0^p$ to $\eta_0^p\circ a$.
\end{itemize}
\end{definition}

We obtain in the obvious way a correspondence
\begin{align}\label{eq:qs_basic_correspondence}
\xymatrix{
\rS_\fp(\rV,\obj)  &
\rB_\fp(\rV,\obj) \ar[r]^-{\iota}\ar[l]_-{\pi} &
\rM_\fp(\rV,\obj)
}
\end{align}
in $\Fun(\fK(\rV)^p\times\fT,\sfP\Sch'_{/\dF_p^\Phi})_{/\rT_\fp}$.

\begin{definition}[Basic correspondence]\label{de:basic_correspondence}
We refer to \eqref{eq:qs_basic_correspondence} as the \emph{basic correspondence} on $\rM_\fp(\rV,\obj)$,\footnote{We adopt this terminology since the image of $\iota$ is in fact the basic locus of $\rM_\fp(\rV,\obj)$.} with $\rS_\fp(\rV,\obj)$ being the \emph{source} of the basic correspondence.
\end{definition}

\begin{theorem}\label{th:qs_basic_correspondence}
In the diagram \eqref{eq:qs_basic_correspondence}, take a point
\[
s^\star=(A_0,\lambda_0,\eta_0^p;A^\star,\lambda^\star,\eta^{p\star})\in\rS_\fp(\rV,\rK^p)(\kappa)
\]
where $\kappa$ is a field containing $\dF^\Phi_p$. Put $\rB_{s^\star}\coloneqq\pi^{-1}(s^\star)$, and denote by $(\cA,\lambda,\eta^p;\alpha)$ the universal object over the fiber $\rB_{s^\star}$.
\begin{enumerate}
  \item The fiber $\rB_{s^\star}$ is a smooth scheme over $\kappa$, with a canonical isomorphism for its tangent bundle
        \[
        \cT_{\rB_{s^\star}/\kappa}\simeq\HOM\(\omega_{\cA^\vee,\tau_\infty},\Ker\alpha_{*,\tau_\infty}/\omega_{\cA^\vee,\tau_\infty}\).
        \]

  \item The restriction of $\iota$ to $\rB_{s^\star}$ is locally on $\rB_{s^{\star}}$ a  closed immersion, with a canonical isomorphism for its normal bundle
        \[
        \cN_{\iota\res\rB_{s^\star}}\simeq\HOM\(\omega_{\cA^\vee,\tau_\infty},\IM\alpha_{*,\tau_\infty}\).
        \]

  \item The assignment sending a point
        $(A_0,\lambda_0,\eta_0^p;A,\lambda,\eta^p;A^\star,\lambda^\star,\eta^{p\star};\alpha)\in\rB_{s^\star}(S)$  for every  $S\in \Sch'_{/\kappa}$ to the subbundle
        \[
        H\coloneqq(\breve\alpha_{*,\tau_\infty})^{-1}\omega_{A^\vee/S,\tau_\infty}\subseteq\rH^\dr_1(A^\star/S)_{\tau_\infty}=
        \rH^\dr_1(A^\star/\kappa)_{\tau_\infty}\otimes_\kappa\cO_S=(\sV_{s^\star})_S,
        \]
        where $\breve\alpha\colon A^\star\to A$ is the (unique) $O_F$-linear quasi-$p$-isogeny such that $\breve\alpha\circ\alpha=\varpi\cdot\id_A$, induces an isomorphism
        \[
        \zeta_{s^\star}\colon\rB_{s^\star}\xrightarrow{\sim}\DL_{s^\star}=\DL(\sV_{s^\star},\{\;,\;\}_{s^\star},\ceil{\tfrac{N+1}{2}}).
        \]
        In particular, $\rB_{s^\star}$ is a geometrically irreducible projective smooth scheme in $\Sch_{/\kappa}$ of dimension $\floor{\tfrac{N-1}{2}}$ by Lemma \ref{le:qs_dl}. In particular, $\iota$ is of pure codimension $\floor{\tfrac{N}{2}}$.
\end{enumerate}
\end{theorem}

\begin{proof}
For an object $(A_0,\lambda_0,\eta_0^p;A,\lambda,\eta^p;A^\star,\lambda^\star,\eta^{p\star};\alpha)\in\rB_\fp(\rV,\rK^p)(S)$, Definition \ref{de:qs_basic_correspondence}(a) implies that there is a (unique) $O_F$-linear quasi-$p$-isogeny $\breve\alpha\colon A^\star\to A$ such that $\breve\alpha\circ\alpha=\varpi\cdot\id_A$, hence $\alpha\circ\breve\alpha=\varpi\cdot\id_{A^\star}$. Moreover, we have the following properties from Definition \ref{de:qs_basic_correspondence}:
\begin{enumerate}[label=(\alph*')]
  \item $\Ker\breve\alpha[p^\infty]$ is contained in $A^\star[\fp]$;

  \item we have $\varpi\cdot\lambda^\star=\breve\alpha^\vee\circ\lambda\circ\breve\alpha$; and

  \item the $\rK^p$-orbit of maps $v\mapsto\varpi^{-1}\breve\alpha_*\circ\eta^{\star p}(v)$ for $v\in\rV\otimes_\dQ\dA^{\infty,p}$ coincides with $\eta^p$.
\end{enumerate}

First, we show (1). It is clear that $\rB_{s^\star}$ is a scheme of finite type over $\kappa$. Consider a closed immersion $S\hookrightarrow\hat{S}$ in $\Sch'_{/\kappa}$ defined by an ideal sheaf $\cI$ satisfying $\cI^2=0$. Take a point $x=(A_0,\lambda_0,\eta_0^p;A,\lambda,\eta^p;A^\star,\lambda^\star,\eta^{p\star};\alpha)\in\rB_{s^\star}(S)$. To compute lifting of $x$ to $\hat{S}$, we use the Serre--Tate and Grothendieck--Messing theories. Note that lifting $\alpha$ is equivalent to lifting both $\alpha$ and $\breve\alpha$, satisfying (b,c) in Definition \ref{de:qs_basic_correspondence} and (b',c') above, respectively. Thus, by Proposition \ref{pr:deformation}, to lift $x$ to an $\hat{S}$-point is equivalent to lifting
\begin{itemize}[label={\ding{109}}]
  \item $\omega_{A^\vee/S,\tau_\infty}$ to a subbundle $\hat\omega_{A^\vee,\tau_\infty}$ of $\rH^\cris_1(A/\hat{S})_{\tau_\infty}$ (of rank $1$),

  \item $\omega_{A^\vee/S,\tau_\infty^\tc}$ to a subbundle $\hat\omega_{A^\vee,\tau_\infty^\tc}$ of $\rH^\cris_1(A/\hat{S})_{\tau_\infty^\tc}$ (of rank $N-1$),
\end{itemize}
subject to the following requirements
\begin{enumerate}[label=(\alph*'')]
  \item $\hat\omega_{A^\vee,\tau_\infty}$ and $\hat\omega_{A^\vee,\tau_\infty^\tc}$ are orthogonal under $\langle\;,\;\rangle_{\lambda,\tau_\infty}^\cris$ \eqref{eq:weil_pairing_cris}; and

  \item $\breve\alpha_{*,\tau_\infty^\tc}\rH^\cris_1(A^\star/\hat{S})_{\tau_\infty^\tc}$ is contained in $\hat\omega_{A^\vee,\tau_\infty^\tc}$.
\end{enumerate}
Since $\langle\;,\;\rangle_{\lambda,\tau_\infty}^\cris$ is a perfect pairing, $\hat\omega_{A^\vee,\tau_\infty}$ uniquely determines $\hat\omega_{A^\vee,\tau_\infty^\tc}$ by (a''). Moreover, by Property (b') above, we know that $\Ker\alpha_{*,\tau_\infty}$ and $\IM\breve\alpha_{*,\tau_\infty^\tc}$ are orthogonal complements to each other under $\langle\;,\;\rangle_{\lambda,\tau_\infty}^\cris$. Thus, (b'') is equivalent to
\begin{enumerate}
  \item [(c'')] $\hat\omega_{A^\vee,\tau_\infty}$ is contained in the kernel of $\alpha_{*,\tau_\infty}\colon\rH^\cris_1(A/\hat{S})_{\tau_\infty}\to\rH^\cris_1(A^\star/\hat{S})_{\tau_\infty}$.
\end{enumerate}

To summarize, lifting $x$ to an $\hat{S}$-point is equivalent to lifting $\omega_{A^\vee/S,\tau_\infty}$ to a subbundle $\hat\omega_{A^\vee,\tau_\infty}$ of $\Ker\alpha_{*,\tau_\infty}$. In other words, the subset of $\rB_{s^\star}(\hat S)$ above $x$ is canonically a torsor over $\Hom_{\cO_S}(\omega_{A^\vee,\tau_\infty},(\Ker\alpha_{*,\tau_\infty}/\omega_{A^\vee,\tau_\infty})\otimes_{\cO_S}\cI)$. Thus, (1) follows.

Next, we show (2). By Theorem \ref{th:qs_moduli_scheme}, we have a canonical isomorphism
\[
\iota_\kappa^*\cT_{\rM_\fp(\rV,\rK^p)/\kappa}\res_{\rB_{s^\star}}\simeq
\HOM\(\omega_{\cA^\vee,\tau_\infty},\rH^\dr_1(\cA)_{\tau_\infty}/\omega_{\cA^\vee,\tau_\infty}\),
\]
and the induced map $\cT_{\rB_{s^\star}/\kappa}\to\iota_\kappa^*\cT_{\rM_\fp(\rV,\rK^p)/\kappa}\res_{\rB_{s^\star}}$ is identified with the canonical map
\[
\HOM\(\omega_{\cA^\vee,\tau_\infty},\Ker\alpha_{*,\tau_\infty}/\omega_{\cA^\vee,\tau_\infty}\)\to
\HOM\(\omega_{\cA^\vee,\tau_\infty},\rH^\dr_1(\cA)_{\tau_\infty}/\omega_{\cA^\vee,\tau_\infty}\).
\]
It is clearly injective, with cokernel canonically isomorphic to
\[
\HOM\(\omega_{\cA^\vee,\tau_\infty},\IM\alpha_{*,\tau_\infty}\).
\]
Thus, (2) follows.

Finally, we show (3). We first show that $\zeta_{s^\star}$ has the correct image, namely, $H$ is a locally free $\cO_S$-module of rank $\ceil{\tfrac{N+1}{2}}$, and satisfies $(\tF H^{(p)})^\perp\subseteq H$. Lemma \ref{le:inverse_isogeny}(1,2,3) implies that $H$ is locally free, and
\begin{align*}
\rank_{\cO_S}(\Ker\alpha_{*,\tau_\infty})-\rank_{\cO_S}(\Ker\alpha_{*,\tau_\infty^\tc})&=1,\\
\rank_{\cO_S}(\Ker\alpha_{*,\tau_\infty})+\rank_{\cO_S}(\Ker\alpha_{*,\tau_\infty^\tc})&=2\ceil{\tfrac{N}{2}}-1.
\end{align*}
Thus, we have $\rank_{\cO_S}(\Ker\alpha_{*,\tau_\infty})=\ceil{\tfrac{N}{2}}$ and
\[
\rank_{\cO_S}(\Ker\breve\alpha_{*,\tau_\infty})=N-\rank_{\cO_S}(\Ker\alpha_{*,\tau_\infty})=\ceil{\tfrac{N-1}{2}}.
\]
On the other hand, as $\omega_{A^\vee/S,\tau_\infty}$ has rank $1$ and $\omega_{A^{\star\vee}/S,\tau_\infty}$ has rank $0$, $\omega_{A^\vee/S,\tau_\infty}$ is contained in the kernel of $\alpha_{*,\tau_\infty}$, hence in the image of $\breve\alpha_{*,\tau_\infty}$. Together, we obtain $\rank_{\cO_S}H=\ceil{\tfrac{N+1}{2}}$. From the equalities
\begin{align*}
\breve\alpha_{*,\tau_\infty^\tc}(\tF H^{(p)})
&=\breve\alpha_{*,\tau_\infty^\tc}\tF_{A^\star}\((\breve\alpha_{*,\tau_\infty})^{-1}\omega_{A^\vee/S,\tau_\infty}\)^{(p)}
=\breve\alpha_{*,\tau_\infty^\tc}\tF_{A^\star}(\breve\alpha^{(p)}_{*,\tau_\infty^\tc})^{-1}\omega_{A^{(p)\vee}/S,\tau_\infty^\tc}\\
&=\tF_A\breve\alpha^{(p)}_{*,\tau_\infty^\tc}(\breve\alpha^{(p)}_{*,\tau_\infty^\tc})^{-1}\omega_{A^{(p)\vee}/S,\tau_\infty^\tc}
=\tF_A\omega_{A^{(p)\vee}/S,\tau_\infty^\tc}=0
\end{align*}
and the fact that $\tF H^{(p)}$ and $\Ker\breve\alpha_{*,\tau_\infty^\tc}$ are both subbundles of $\rH^\dr_1(A^\star/S)_{\tau_\infty^\tc}$ of rank $\ceil{\tfrac{N+1}{2}}$, we know $\tF H^{(p)}=\Ker\breve\alpha_{*,\tau_\infty^\tc}$. By Definition \ref{de:qs_basic_correspondence}(b) and the definition of $\breve\alpha$, we have
\[
\langle\Ker\breve\alpha_{*,\tau_\infty^\tc},\IM\alpha_{*,\tau_\infty}\rangle_{\lambda^\star,\tau_\infty^\tc}=
\langle\breve\alpha_{*,\tau_\infty^\tc}\Ker\breve\alpha_{*,\tau_\infty^\tc},\rH^\dr_1(A/S)_{\tau_\infty}\rangle_{\lambda,\tau_\infty^\tc}=0,
\]
which implies
\[
\Ker\breve\alpha_{*,\tau_\infty}=\IM\alpha_{*,\tau_\infty}\subseteq(\Ker\breve\alpha_{*,\tau_\infty^\tc})^\perp=(\tF H^{(p)})^\perp.
\]
As both sides are subbundles of $\rH^\dr_1(A^\star/S)_{\tau_\infty}$ of rank $\ceil{\tfrac{N-1}{2}}$, we must have $\Ker\breve\alpha_{*,\tau_\infty}=(\tF H^{(p)})^\perp$. In particular, we have $(\tF H^{(p)})^\perp\subseteq H$. Thus, $\zeta_{s^\star}$ is defined as we claim.

Since the target of $\zeta_{s^\star}$ is smooth over $\kappa$ by Lemma \ref{le:qs_dl}, to see that $\zeta_{s^\star}$ is an isomorphism, it suffices to check that for every algebraically closed field $\kappa'$ containing $\kappa$, the following statements hold:
\begin{enumerate}
  \item [(3--1)] $\zeta_{s^\star}$ induces a bijection on $\kappa'$-points; and

  \item [(3--2)] $\zeta_{s^\star}$ induces an isomorphism on the tangent spaces at every $\kappa'$-point.
\end{enumerate}
To ease notation, we may assume that $\kappa'=\kappa$, hence is perfect in particular.

For (3--1), we construct an inverse to the map $\zeta_{s^\star}(\kappa)$. Take a point $y\in\DL_{s^\star}(\kappa)$ represented by a $\kappa$-linear subspace $H\subseteq\sV_{s^\star}=\rH^\dr_1(A^\star/\kappa)_{\tau_\infty}$. We regard $\tF$ and $\tV$ as those sesquilinear maps in Notation \ref{re:frobenius_verschiebung}. In particular, we have $(\tF H)^\perp\subseteq H$. For every $\tau\in\Sigma_\infty$, we define a $W(\kappa)$-submodule $\cD_{A,\tau}\subseteq\cD(A^\star)_\tau$ as follows.
\begin{itemize}[label={\ding{109}}]
  \item If $\tau\not\in\{\tau_\infty,\tau_\infty^\tc\}$, then $\cD_{A,\tau}=\cD(A^\star)_\tau$.

  \item We set $\cD_{A,\tau_\infty}\coloneqq\tV^{-1}\tilde{H}^\tc$, where $\tilde{H}^\tc$ is the preimage of $H^\perp$ under the reduction map $\cD(A^\star)_{\tau_\infty^\tc}\to\cD(A^\star)_{\tau_\infty^\tc}/p\cD(A^\star)_{\tau_\infty^\tc}=\rH^\dr_1(A^\star)_{\tau_\infty^\tc}$.

  \item We set $\cD_{A,\tau_\infty^\tc}\coloneqq\tF\tilde{H}$, where $\tilde{H}$ is the preimage of $H$ under the reduction map $\cD(A^\star)_{\tau_\infty}\to\cD(A^\star)_{\tau_\infty}/p\cD(A^\star)_{\tau_\infty}=\rH^\dr_1(A^\star)_{\tau_\infty}$.
\end{itemize}
Finally, put $\cD_A\coloneqq\bigoplus_{\tau\in\Sigma_\infty}\cD_{A,\tau}$ as a $W(\kappa)$-submodule of $\cD(A^\star)$. We show that it is stable under $\tF$ and $\tV$. It suffices to show that both $\tF$ and $\tV$ stabilize $\cD_{A,\tau_\infty}\oplus\cD_{A,\tau_\infty^\tc}$, which breaks into checking that
\begin{itemize}[label={\ding{109}}]
  \item $\tF\cD_{A,\tau_\infty}\subseteq\cD_{A,\tau_\infty^\tc}$, that is, $\tF\tV^{-1}\tilde{H}^\tc\subseteq\tF\tilde{H}$. It suffices to show that $\tV^{-1}(H^\perp)$ (as a subspace of $\rH^\dr_1(A^\star)_{\tau_\infty}$) is contained in $H$. However, $\tV^{-1}(H^\perp)=(\tF H)^\perp$, which is contained in $H$.

  \item $\tF\cD_{A,\tau_\infty^\tc}\subseteq\cD_{A,\tau_\infty}$, that is, $\tF\tF\tilde{H}\subseteq\tV^{-1}\tilde{H}^\tc$. It suffices to show $p\tF\tilde{H}\subseteq\tilde{H}^\tc$, which obviously holds.

  \item $\tV\cD_{A,\tau_\infty}\subseteq\cD_{A,\tau_\infty^\tc}$, that is, $\tV\tV^{-1}\tilde{H}^\tc\subseteq\tF\tilde{H}$. it suffices to show $H^\perp\subseteq\tF H$ as subspaces of $\rH^\dr_1(A^\star)_{\tau_\infty^\tc}$, which follows from $(\tF H)^\perp\subseteq H$.

  \item $\tV\cD_{A,\tau_\infty^\tc}\subseteq\cD_{A,\tau_\infty}$, that is, $\tV\tF\tilde{H}\subseteq\tV^{-1}\tilde{H}^\tc$. It is obvious as $\tV^{-1}\tilde{H}^\tc$ contains $p\cD(A^\star)_{\tau_\infty}$.
\end{itemize}
Thus, $(\cD_A,\tF,\tV)$ is a Dieudonn\'{e} module over $W(\kappa)$. By the Dieudonn\'{e} theory, there is an $O_F$-abelian scheme $A$ over $\kappa$ with $\cD(A)_\tau=\cD_{A,\tau}$ for every $\tau\in\Sigma_\infty$, and an $O_F$-linear $p$-isogeny $\alpha\colon A\to A^\star$ inducing the inclusion of Dieudonn\'{e} modules $\cD(A)=\cD_A\subseteq\cD(A^\star)$. Moreover, since $p\cD(A^\star)\subseteq\cD(A)$, we have $\Ker\alpha[p^\infty]\subseteq A[\fp]$.

Let $\lambda\colon A\to A^\vee$ be the unique quasi-polarization such that $\varpi\lambda=\alpha^\vee\circ\lambda^\star\circ\alpha$. We claim that $\lambda$ is $p$-principal. It is enough to show the induced pairing
\[
p^{-1}\cdot\langle\;,\;\rangle_{\lambda^\star,\tau_\infty}\colon\cD(A)_{\tau_\infty}\times\cD(A)_{\tau_\infty^\tc}\to W(\kappa)
\]
(Notation \ref{no:weil_pairing_dieudonne}) is non-degenerate. Since $\tilde{H}$ is $W(\kappa)$-dual to $p^{-1}\tilde{H}^\tc$, hence $\cD(A)_{\tau_\infty^\tc}=\tF\tilde{H}$ is dual to $\tV^{-1}(p^{-1}\tilde{H}^\tc)=p^{-1}\tV^{-1}\tilde{H}^\tc=p^{-1}\cD(A)_{\tau_\infty}$, the above pairing is non-degenerate.

It is an easy consequence of Lemma \ref{le:inverse_isogeny}(2,3) that the $O_F$-abelian scheme $A$ has signature type $N\Phi-\tau_\infty+\tau_\infty^\tc$. Finally, let $\eta^p$ be the unique $\rK^p$-level structure such that Definition \ref{de:qs_basic_correspondence}(c) is satisfied. Putting together, we obtain a point
$x=(A_0,\lambda_0,\eta_0^p;A,\lambda,\eta^p;A^\star,\lambda^\star,\eta^{p\star};\alpha)\in\rB_{s^\star}(\kappa)$ such that $\zeta_{s^\star}(x)=y$. It is easy to see that such assignment gives rise to an inverse of $\zeta_{s^\star}(\kappa)$, hence (3--1) follows immediately.

For (3--2), let $\cT_x$ and $\cT_y$ be the tangent spaces at $x$ and $y$ as in (3--1), respectively. By (1) and Lemma \ref{le:qs_dl}, we have canonical isomorphisms
\[
\cT_x\simeq\Hom_\kappa(\omega_{A^\vee,\tau_\infty},\Ker\alpha_{*,\tau_\infty}/\omega_{A^\vee,\tau_\infty}),\quad
\cT_y\simeq\Hom_\kappa(H/(\tF H)^\perp,\rH^\dr_1(A^\star)_{\tau_\infty}/H).
\]
Moreover, by the definition of $\zeta_{s^\star}$, the map $(\zeta_{s^\star})_*\colon\cT_x\to\cT_y$ is induced by the following two maps
\begin{align*}
H/(\tF H)^\perp&=(\breve\alpha_{*,\tau_\infty})^{-1}\omega_{A^\vee,\tau_\infty}/\Ker\breve\alpha_{*,\tau_\infty}
\xrightarrow{\breve\alpha_{*,\tau_\infty}}\omega_{A^\vee,\tau_\infty},\\
\rH^\dr_1(A^\star)_{\tau_\infty}/H&=
\rH^\dr_1(A^\star)_{\tau_\infty}/(\breve\alpha_{*,\tau_\infty})^{-1}\omega_{A^\vee,\tau_\infty}
\xrightarrow{\breve\alpha_{*,\tau_\infty}}\Ker\alpha_{*,\tau_\infty}/\omega_{A^\vee,\tau_\infty},
\end{align*}
both being isomorphisms. Thus, (3--2) and hence (3) follow.
\end{proof}

\begin{remark}
In Theorem \ref{th:qs_basic_correspondence}, when $\rK^p$ is sufficiently small, the restriction of $\iota$ to $\rB_{s^\star}$ is a closed immersion for every point $s^\star\in\rS_\fp(\rV,\rK^p)(\kappa)$ and every field $\kappa$ containing $\dF^\Phi_p$.
\end{remark}

\subsection{Source of basic correspondence and Tate cycles}

In this subsection, we study the source $\rS_\fp(\rV,\obj)$ of the basic correspondence. We will describe the set $\rS_\fp(\rV,\obj)(\ol\dF_p)$ in terms of a certain Shimura set and study its Galois action. Such a description is not canonical, which depends on the choice of a definite uniformization datum defined as follows.

\begin{definition}\label{de:qs_uniformization_data}
We define a \emph{definite uniformization datum for $\rV$ (at $\fp$)} to be a collection of $(\rV^\star,\ti,\{\Lambda^\star_\fq\}_{\fq\mid p})$, where
\begin{itemize}[label={\ding{109}}]
  \item $\rV^\star$ is a standard definite hermitian space over $F$ of rank $N$;

  \item $\ti\colon\rV\otimes_\dQ\dA^{\infty,p}\to\rV^\star\otimes_\dQ\dA^{\infty,p}$ is an isometry;

  \item for every prime $\fq$ of $F^+$ above $p$ other than $\fp$, $\Lambda^\star_\fq$ is a self-dual $O_{F_\fq}$-lattice in $\rV^\star\otimes_FF_\fq$; and

  \item $\Lambda^\star_\fp$ is an $O_{F_\fp}$-lattice in $\rV^\star\otimes_FF_\fp$ satisfying $p\Lambda^\star_\fp\subseteq(\Lambda^\star_\fp)^\vee$ such that $(\Lambda^\star_\fp)^\vee/p\Lambda^\star_\fp$ has length $0$ (resp.\ $1$) if $N$ is odd (resp.\ even).
\end{itemize}
\end{definition}

By the Hasse principle for hermitian spaces, there exists a definite uniformization datum for which we fix one. Let $\rK^\star_\fq$ be the stabilizer of $\Lambda^\star_\fq$ for every $\fq$ over $p$; and put $\rK^\star_p\coloneqq\prod_{\fq\mid p}\rK^\star_\fq$. The isometry $\ti$ induces an equivalence of categories $\ti\colon\fK(\rV)^p\xrightarrow{\sim}\fK(\rV^\star)^p$.

\begin{construction}\label{cs:qs_uniformization}
We now construct a \emph{uniformization map}, denoted by the Greek letter \emph{upsilon}
\begin{align}\label{eq:qs_uniformization}
\upsilon\colon\rS_\fp(\rV,\obj)(\ol\dF_p)
\to\Sh(\rV^\star,(\ti\obj)\rK^\star_p)\times\rT_\fp(\ol\dF_p)
\end{align}
in $\Fun(\fK(\rV)^p\times\fT,\Set)_{/\rT_\fp(\ol\dF_p)}$, which turns out to be an isomorphism.

Take a point $s^\star=(A_0,\lambda_0,\eta_0^p;A^\star,\lambda^\star,\eta^{p\star})\in\rS_\fp(\rV,\rK^p)(\ol\dF_p)$. Let
\[
\rV_{s^\star}\coloneqq\Hom_{O_F}(A_0,A^\star)\otimes\dQ
\]
be the space of $O_F$-linear quasi-homomorphisms. We equip $\rV_{s^\star}$ with a pairing
\[
(x,y)=\varpi^{-1}\cdot\lambda_0^{-1}\circ y^\vee\circ\lambda^\star\circ x\in\End_{O_F}(A_0)\otimes\dQ=F,
\]
which becomes a hermitian space over $F$. Note that we have an extra factor $\varpi^{-1}$ in the above pairing. Moreover, for every prime $\fq$ of $F^+$ above $p$, put
\[
\Lambda_{s^\star,\fq}\coloneqq\Hom_{O_F}(A_0[\fq^\infty],A^\star[\fq^\infty])
\]
which is an $O_{F_\fq}$-lattice in $(\rV_{s^\star})_\fq$ since $A^\star$ is isogenous to $A_0^N$.

Now we construct $\upsilon$, whose process is very similar to Remark \ref{re:qs_moduli_scheme_shimura}. Note that we have an isometry
\[
\rho\colon\rV_{s^\star}\otimes_\dQ\dA^{\infty,p}\xrightarrow{\sim}
\Hom_{F\otimes_\dQ\dA^{\infty,p}}^{\varpi\lambda_0,\lambda^\star}(\rH^\et_1(A_0,\dA^{\infty,p}),\rH^\et_1(A^\star,\dA^{\infty,p})).
\]
By Lemma \ref{le:qs_uniformization} below, we can choose an isometry $\eta_\rat\colon\rV_{s^\star}\to\rV^\star$. Thus, we obtain an isometry
\[
g^p\coloneqq \eta_\rat\circ\rho^{-1}\circ\eta^{p\star}\circ\ti^{-1}\colon\rV^\star\otimes_\dQ\dA^{\infty,p}
\to\rV^\star\otimes_\dQ\dA^{\infty,p}
\]
as an element in $\rU(\rV^\star)(\dA_{F^+}^{\infty,p})$. By Lemma \ref{le:qs_uniformization}(1,2), for every $\fq$ above $p$, there exists an element $g_\fq\in\rU(\rV^\star)(F^+_\fq)$ such that $g_\fq\Lambda^\star_\fq=\eta_\rat\Lambda_{s^\star,\fq}$. Together, we obtain an element $g_{s^\star}\coloneqq(g^p,(g_\fq)_{\fq\mid p})\in\rU(\rV^\star)(\dA_{F^+}^\infty)$ such that the double coset $\rU(\rV^\star)(F)g(\ti\rK^p)\rK^\star_p$ depends only on the point $s^\star$. Thus, it allows us to define
\[
\upsilon(s^\star)\coloneqq\(\rU(\rV^\star)(F)g_{s^\star}(\ti\rK^p)\rK^\star_p,(A_0,\lambda_0,\eta_0^p)\)
\in\Sh(\rV^\star,(\ti\rK^p)\rK^\star_p)\times\rT_\fp(\ol\dF_p).
\]
\end{construction}

\begin{lem}\label{le:qs_uniformization}
The hermitian spaces $\rV_{s^\star}$ and $\rV^\star$ are isomorphic. Moreover,
\begin{enumerate}
  \item for every prime $\fq$ of $F^+$ above $p$ other than $\fp$, the lattice $\Lambda_{s^\star,\fq}$ is self-dual;

  \item the lattice $\Lambda_{s^\star,\fp}$ satisfies $p\Lambda_{s^\star,\fp}\subseteq(\Lambda_{s^\star,\fp})^\vee$ such that $(\Lambda_{s^\star,\fp})^\vee/p\Lambda_{s^\star,\fp}$ has length $0$ (resp.\ $1$) if $N$ is odd (resp.\ even).
\end{enumerate}
\end{lem}

\begin{proof}
We first prove (1) and (2).

For (1), note that $A^\star[\fq^\infty]$ is isomorphic to $(A_0[\fq^\infty])^N$, equipped with the polarization $\lambda^\star[\fq^\infty]$ that is principal. Thus, $\Lambda_{s^\star,\fq}$ is self-dual as $\lambda_0[\fq^\infty]$ is principal and $\val_\fq(\varpi)=0$.

For (2), note that $A^\star[\fp^\infty]$ is isomorphic to $(A_0[\fp^\infty])^N$, equipped with the polarization $\lambda^\star[\fp^\infty]$ satisfying such that $\Ker\lambda^\star[\fp^\infty]$ is trivial (resp.\ contained in $A^\star[\fp]$ of rank $p^2$) if $N$ is odd (resp.\ even). Thus, the statement follows as $\lambda_0[\fp^\infty]$ is principal and $\val_\fp(\varpi)=1$.

Now to prove the main statement, it suffices to show that
\begin{enumerate}[label=(\roman*)]
  \item $\rV_{s^\star}$ is totally positive definite; and

  \item the hermitian spaces $\rV_{s^\star}\otimes_\dQ\dA^{\infty,p}$ and $\rV\otimes_\dQ\dA^{\infty,p}$ are isomorphic.
\end{enumerate}

For (i), it follows from the same argument in \cite{KR14}*{Lemma~2.7}.

For (ii), we have a map
\[
\rV_{s^\star}\otimes_\dQ\dA^{\infty,p}\to
\Hom_{F\otimes_\dQ\dA^{\infty,p}}^{\varpi\lambda_0,\lambda^\star}(\rH^\et_1(A_0,\dA^{\infty,p}),\rH^\et_1(A^\star,\dA^{\infty,p}))
\]
of hermitian spaces, which is injective. As both sides have rank $N$ and the right-hand side is isomorphic to $\rV\otimes_\dQ\dA^{\infty,p}$, (ii) follows.
\end{proof}

\begin{proposition}\label{pr:qs_uniformization}
The uniformization map $\upsilon$ \eqref{eq:qs_uniformization} is an isomorphism. Moreover, the induced action of $\Gal(\ol\dF_p/\dF_p^\Phi)$ on the target of $\upsilon$ factors through the projection map
\[
\Sh(\rV^\star,(\ti\obj)\rK^\star_p)\times\rT_\fp(\ol\dF_p)\to\rT_\fp(\ol\dF_p).
\]
\end{proposition}

\begin{proof}
We first show that $\upsilon$ is an isomorphism. Take a point $t=(A_0,\lambda_0,\eta_0^p)\in\rT_\fp(\ol\dF_p)$. It suffices to show that, for every $\rK^p\in\fK(\rV)^p$, the restriction
\[
\upsilon\colon\rS_\fp(\rV,\rK^p)(\ol\dF_p)_{/t}\to\Sh(\rV^\star,(\ti\rK^p)\rK^\star_p)
\]
to the fiber over $t$ is an isomorphism. The injectivity follows directly from the definition. For the surjectivity, it suffices to show that for every $g\in\rU(\rV^\star)(\dA_{F^+}^{\infty,p})$, there is an object $s^\star=(A_0,\lambda_0,\eta_0^p;A^\star,\lambda^\star,\eta^{p\star})\in\rS_\fp(\rV,\rK^p)(\ol\dF_p)_{/t}$ whose image under $\upsilon$ is the image of $g$ in $\Sh(\rV^\star,(\ti\rK^p)\rK^\star_p)$. To construct $s^\star$, we take an $O_F$-lattice $\Lambda^\star$ in $\rV^\star$ satisfying $\Lambda^\star\otimes_FF_\fp=\Lambda^\star_\fp$. Put $A^\star\coloneqq A_0\otimes_{O_F}\Lambda^\star$, which is equipped with a unique quasi-polarization $\lambda^\star$ such that the canonical isomorphism
\[
\rV^\star\otimes_\dQ\dA^{\infty,p}\simeq\Hom_{F\otimes_\dQ\dA^{\infty,p}}(\rH^\et_1(A_0,\dA^{\infty,p}),\rH^\et_1(A^\star,\dA^{\infty,p}))
\]
of $F\otimes_\dQ\dA^{\infty,p}$-modules is an isometry of hermitian spaces. We let $\eta^{p\star}$ be the map
\[
\rV\otimes_\dQ\dA^{\infty,p}\xrightarrow{g\circ\ti}\rV^\star\otimes_\dQ\dA^{\infty,p}=
\Hom_{F\otimes_\dQ\dA^{\infty,p}}^{\varpi\lambda_0,\lambda^\star}(\rH^\et_1(A_0,\dA^{\infty,p}),\rH^\et_1(A^\star,\dA^{\infty,p})).
\]
Then $\upsilon(s^\star)=g$ in $\Sh(\rV^\star,(\ti\rK^p)\rK^\star_p)$. Thus, $\upsilon$ is an isomorphism.

Since $\upsilon$ is an isomorphism, the Galois group $\Gal(\ol\dF_p/\dF_p^\Phi)$ acts on the target of $\upsilon$. We show that it acts trivially on the first factor of the target of $\upsilon$. Take an element $\varsigma\in\Gal(\ol\dF_p/\dF_p^\Phi)$ and a point $s^\star=(A_0,\lambda_0,\eta_0^p;A^\star,\lambda^\star,\eta^{p\star})\in\rS_\fp(\rV,\rK^p)(\ol\dF_p)$. Then $\varsigma s^\star$ is simply represented by $(A_0^\varsigma,\lambda_0^\varsigma,\eta_0^{p\varsigma};A^{\star\varsigma},\lambda^{\star\varsigma},\eta^{p\star\varsigma})$, the $\varsigma$-twist of the previous object. We then have a canonical isomorphism
\[
\rV_{\varsigma s^\star}=\Hom_{O_F}(A_0^\varsigma,A^{\star\varsigma})\otimes\dQ\simeq\Hom_{O_F}(A_0,A^\star)\otimes\dQ=\rV_{s^\star}
\]
of hermitian spaces. Unraveling the definition, we see that $g_{s^\star}=g_{\varsigma s^\star}$. Thus, we have
\[
\upsilon(\varsigma s^\star)\coloneqq\(\rU(\rV^\star)(F)g_{s^\star}(\ti\rK^p)\rK^\star_p,(A_0^\varsigma,\lambda_0^\varsigma,\eta_0^{p\varsigma})\).
\]
The proposition follows.
\end{proof}

Next, we define an action of the Hecke algebra $\dZ[\rK^\star_\fp\backslash\rU(\rV^\star)(F^+_\fp)/\rK^\star_\fp]$ on $\rS_\fp(\rV,\obj)$ via finite \'{e}tale correspondences, that is compatible with the uniformization map \eqref{eq:qs_uniformization}.

\begin{construction}\label{cs:qs_definite_hecke}
For every element $g\in\rK^\star_\fp\backslash\rU(\rV^\star)(F^+_\fp)/\rK^\star_\fp$, we define a functor
\begin{align*}
\rS_\fp(\rV,\obj)_g\colon\fK(\rV)^p\times\fT &\to\sfP\Sch'_{/\dF^\Phi_p} \\
\rK^p &\mapsto \rS_\fp(\rV,\rK^p)_g
\end{align*}
such that for every $S\in\Sch'_{/\dF^\Phi_p}$, $\rS_\fp(\rV,\rK^p)_g(S)$ is the set of equivalence classes of decuples
$(A_0,\lambda_0,\eta_0^p;A^\star,\lambda^\star,\eta^{p\star};A^\star_g,\lambda^\star_g,\eta^{p\star}_g;\phi^\star)$, where
\begin{itemize}[label={\ding{109}}]
  \item $(A_0,\lambda_0,\eta_0^p;A^\star,\lambda^\star,\eta^{p\star})$ and $(A_0,\lambda_0,\eta_0^p;A^\star_g,\lambda^\star_g,\eta^{p\star}_g)$ are both elements in $\rS_\fp(\rV,\rK^p)(S)$; and

  \item $\phi^\star\colon A^\star\to A^\star_g$ is an $O_F$-linear quasi-isogeny such that
     \begin{enumerate}[label=(\alph*)]
       \item $\phi^{\star\vee}\circ\lambda^\star_g\circ\phi^\star=\lambda^\star$;

       \item $\phi^\star[\fp^\infty]\colon A^\star[\fp^\infty]\to A^\star_g[\fp^\infty]$ is a quasi-isogeny of height zero under which the two lattices $\Hom_{O_F}(A_{0s}[\fp^\infty],A^\star_s[\fp^\infty])$ and $\Hom_{O_F}(A_{0s}[\fp^\infty],A^\star_{gs}[\fp^\infty])$ are at the relative position determined by $g$ for every geometric point $s$ of $S$;

       \item $\phi^\star[\fq^\infty]$ is an isomorphism for every prime $\fq$ of $F^+$ above $p$ that is not $\fp$; and

       \item the $\rK^p$-orbit of maps $v\mapsto\phi^\star_*\circ\eta^{p\star}(v)$ for $v\in\rV\otimes_\dQ\dA^{\infty,p}$ coincides with $\eta^{p\star}_g$.
     \end{enumerate}
\end{itemize}
The equivalence relation and the action of morphisms in $\fK(\rV)^p\times\fT$ are defined similarly as in Definition \ref{de:qs_basic_correspondence}. Then we construct the \emph{Hecke correspondence} (of $g$) to be the morphism
\begin{align}\label{eq:qs_definite_hecke}
\Hk_g\colon\rS_\fp(\rV,\obj)_g\to\rS_\fp(\rV,\obj)\times\rS_\fp(\rV,\obj)
\end{align}
in $\Fun(\fK(\rV)^p\times\fT,\sfP\Sch'_{/\dF_p^\Phi})_{/\rT_\fp}$ induced by the assignment
\[
(A_0,\lambda_0,\eta_0^p;A^\star,\lambda^\star,\eta^{p\star};A^\star_g,\lambda^\star_g,\eta^{p\star}_g;\phi^\star)
\mapsto((A_0,\lambda_0,\eta_0^p;A^\star,\lambda^\star,\eta^{p\star}),(A_0,\lambda_0,\eta_0^p;A^\star_g,\lambda^\star_g,\eta^{p\star}_g)).
\]
Here, the product in \eqref{eq:qs_definite_hecke} is also taken in the category $\Fun(\fK(\rV)^p\times\fT,\sfP\Sch'_{/\dF_p^\Phi})_{/\rT_\fp}$, that is, $\rS_\fp(\rV,\obj)\times\rS_\fp(\rV,\obj)$ is a functor sending $\rK^p$ to $\rS_\fp(\rV,\rK^p)\times_{\rT_\fp}\rS_\fp(\rV,\rK^p)$ on which $\fT$ acts diagonally.
\end{construction}

\begin{proposition}\label{pr:qs_definite_hecke}
For every $g\in\rK^\star_\fp\backslash\rU(\rV^\star)(F^+_\fp)/\rK^\star_\fp$, we have
\begin{enumerate}
  \item The morphism $\Hk_g$ \eqref{eq:qs_definite_hecke} is finite \'{e}tale; in particular, it is a morphism in $\Fun(\fK(\rV)^p\times\fT,\Sch_{/\dF_p^\Phi})_{/\rT_\fp}$.

  \item The uniformization map $\upsilon$ \eqref{eq:qs_uniformization} lifts uniquely to an isomorphism making the diagram
    \[
    \xymatrix{
    \rS_\fp(\rV,\obj)_g(\ol\dF_p) \ar[r]^-{\upsilon}\ar[d]_-{\Hk_g(\ol\dF_p)} &
    \Sh(\rV^\star,(\ti\obj)(g\rK^\star_p g^{-1}\cap\rK^\star_p))\times\rT_\fp(\ol\dF_p) \ar[d] \\
    \rS_\fp(\rV,\obj)(\ol\dF_p)\times_{\rT_\fp(\ol\dF_p)}\rS_\fp(\rV,\obj)(\ol\dF_p) \ar[r]^-{\upsilon\times\upsilon} & \(\Sh(\rV^\star,(\ti\obj)\rK^\star_p)\times\Sh(\rV^\star,(\ti\obj)\rK^\star_p)\)\times\rT_\fp(\ol\dF_p)
    }
    \]
    in $\Fun(\fK(\rV)^p\times\fT,\Set)_{/\rT_\fp(\ol\dF_p)}$ commutative, where the right vertical map is induced by the set-theoretical Hecke correspondence of $g$.
\end{enumerate}
\end{proposition}

\begin{proof}
For (1), it suffices to consider those $\rK^p\in\fK(\rV)^p$ that are sufficiently small. Then the morphism
$\Hk_g\colon\rS_\fp(\rV,\rK^p)_g\to\rS_\fp(\rV,\rK^p)\times_{\rT_\fp}\rS_\fp(\rV,\rK^p)$ is closed, hence represented by a finite \'{e}tale scheme. Part (2) follows directly from the definition.
\end{proof}

\begin{remark}
In fact, the proof of Proposition \ref{pr:qs_definite_hecke}(1) together with Proposition \ref{pr:qs_uniformization} imply that $\Hk_g$ is a local isomorphism.
\end{remark}

\begin{remark}\label{re:qs_tate_hecke}
Note that since $\rK^\star_\fp$ is a special maximal open compact subgroup of $\rU(\rV^\star)(F^+_\fp)$, the algebra $\dZ[\rK^\star_\fp\backslash\rU(\rV^\star)(F^+_\fp)/\rK^\star_\fp]$ is commutative. Moreover, when $N$ is odd, $\Lambda_{s^\star,\fp}$ is a self-dual lattice under the pairing $\varpi\cdot(\;,\;)_{\rV^\star}$, hence $\dZ[\rK^\star_\fp\backslash\rU(\rV^\star)(F^+_\fp)/\rK^\star_\fp]$ is canonically isomorphic to $\dT_{N,\fp}$.
\end{remark}

Let $L$ be a $p$-coprime coefficient ring. The uniformization map \eqref{eq:qs_uniformization} induces an isomorphism
\begin{align*}
L[\Sh(\rV^\star,(\ti\obj)\rK^\star_p)]\simeq\rH^0_\fT(\ol\rS_\fp(\rV,\obj),L)=
\rH^0_\fT(\rS_\fp(\rV,\obj),L)
\end{align*}
in $\Fun(\fK(\rV)^p,\Mod(L[\rK^\star_\fp\backslash\rU(\rV^\star\otimes_FF_\fp)/\rK^\star_\fp]))$ by Proposition \ref{pr:qs_definite_hecke}. Recall from Theorem \ref{th:qs_basic_correspondence}(3) that the morphism $\iota$ in \eqref{eq:qs_basic_correspondence} is of pure codimension $\floor{\tfrac{N}{2}}$.

\begin{construction}\label{cs:qs_incidence}
Put $r\coloneqq\lfloor\tfrac{N}{2}\rfloor\geq 0$. We construct a pair of maps
\[
\begin{dcases}
\inc^\star_!\colon L[\Sh(\rV^\star,(\ti\obj)\rK^\star_p)]\xrightarrow{\sim}\rH^0_\fT(\rS_\fp(\rV,\obj),L) \\ \qquad\qquad\xrightarrow{\pi^*}\rH^0_\fT(\rB_\fp(\rV,\obj),L)\xrightarrow{\iota_!}\rH^{2r}_\fT(\rM_\fp(\rV,\obj),L(r)), \\
\inc_\star^*\colon \rH^{2(N-r-1)}_\fT(\rM_\fp(\rV,\obj),L(N-r-1))\xrightarrow{\iota^*}\rH^{2(N-r-1)}_\fT(\rB_\fp(\rV,\obj),L(N-r-1)) \\
\qquad\qquad\xrightarrow{\pi_!}\rH^0_\fT(\rS_\fp(\rV,\obj),L)\xrightarrow{\sim}L[\Sh(\rV^\star,(\ti\obj)\rK^\star_p)],
\end{dcases}
\]
in $\Fun(\fK(\rV)^p,\Mod(L))$. In fact, the two maps are essentially Poincar\'{e} dual to each other.
\end{construction}

\begin{theorem}\label{th:qs_tate}
Suppose that $N=2r+1$ is odd with $r\geq 0$. Then the composite map $\inc_\star^*\circ\inc^\star_!$ is equal to the Hecke operator
\[
\tT^\star_{N,\fp}\coloneqq\sum_{\delta=0}^r\td_{r-\delta,p}\cdot\tT_{N,\fp;\delta}\in\dT_{N,\fp}
\]
in which the numbers $\td_{r-\delta,p}$ are introduced in Notation \ref{no:numerical}, and the Hecke operators $\tT_{N,\fp;\delta}$ are introduced in Notation \ref{no:hecke} (as $\tT^\circ_{N;\delta}$).
\end{theorem}

Note that by Remark \ref{re:qs_tate_hecke}, $L[\Sh(\rV^\star,(\ti\obj)\rK^\star_p)]$ is a $\dT_{N,\fp}$-module when $N$ is odd.

\begin{proof}
This is \cite{XZ17}*{Theorem~9.3.5}.
\end{proof}

\subsection{Functoriality under special morphisms}
\label{ss:qs_functoriality}

In this subsection, we study the behavior of various moduli schemes under the special morphisms, which is closely related to the Rankin--Selberg motives for $\GL_n\times\GL_{n+1}$. We start from the datum $(\rV_n,\{\Lambda_{n,\fq}\}_{\fq\mid p})$ as in the beginning of \S\ref{ss:qs_moduli_scheme}, but with $\rV_n$ of rank $n\geq1$. We then have the induced datum
\[
(\rV_{n+1},\{\Lambda_{n+1,\fq}\}_{\fq\mid p})\coloneqq((\rV_n)_\sharp,\{(\Lambda_{n,\fq})_\sharp\}_{\fq\mid p})
\]
of rank $n+1$ by Definition \ref{de:hermitian_space}. For $N\in\{n,n+1\}$, we let $\rK_{N,\fq}$ be the stabilizer of $\Lambda_{N,\fq}$, and put $\rK_{N,p}\coloneqq\prod_{\fq\mid p}\rK_{N,\fq}$. Recall the category $\fK(\rV_n)_\sp^p$ and functors $\obj_\flat,\obj_\sharp$ from Definition \ref{de:neat_category}. To unify notation, we put $\obj_n\coloneqq\obj_\flat$ and $\obj_{n+1}\coloneqq\obj_\sharp$. There are five stages of functoriality we will consider.

The first stage concerns Shimura varieties. The canonical inclusions
\[
\rV_n\hookrightarrow\rV_{n+1},\quad\{\Lambda_{n,\fq}\hookrightarrow\Lambda_{n+1,\fq}\}_{\fq\mid p}
\]
induce a morphism
\begin{align}\label{eq:qs_functoriality_shimura}
\sh_\uparrow\colon\Sh(\rV_n,\obj_n\rK_{n,p})\to\Sh(\rV_{n+1},\obj_{n+1}\rK_{n+1,p})
\end{align}
in $\Fun(\fK(\rV_n)_\sp^p,\Sch_{/F})$, known as the \emph{special morphism}.

For the second stage of functoriality, we have a morphism
\begin{align}\label{eq:qs_functoriality_moduli_scheme}
\bm_\uparrow\colon \bM_\fp(\rV_n,\obj_n)\to \bM_{\fp}(\rV_{n+1},\obj_{n+1})
\end{align}
in $\Fun(\fK(\rV_n)_\sp^p\times\fT,\Sch_{/\dZ_p^\Phi})_{/\bT_\fp}$ sending an object $(A_0,\lambda_0,\eta_0^p;A,\lambda,\eta^p)\in\bM_\fp(\rV_n,\rK^p_n)(S)$ to the object $(A_0,\lambda_0,\eta_0^p;A\times A_0,\lambda\times\lambda_0,\eta^p\oplus(\id_{A_0})_*)\in\bM_\fp(\rV_{n+1},\rK^p_{n+1})(S)$. We then have the following commutative diagram
\begin{align}\label{eq:qs_functoriality_moduli_scheme1}
\xymatrix{
\bM^\eta_\fp(\rV_{n+1},\obj_{n+1}) \ar[rr]^-{\eqref{eq:qs_moduli_scheme_shimura}}
&&  \Sh(\rV_{n+1},\obj_{n+1}\rK_{n+1,p})\times_{\Spec{F}}\bT^\eta_\fp \\
\bM^\eta_\fp(\rV_n,\obj_n) \ar[rr]^-{\eqref{eq:qs_moduli_scheme_shimura}}\ar[u]^-{\bm^\eta_\uparrow}
&& \Sh(\rV_n,\obj_n\rK_{n,p})\times_{\Spec{F}}\bT^\eta_\fp \ar[u]_-{\sh_\uparrow\times\id}
}
\end{align}
in $\Fun(\fK(\rV_n)_\sp^p\times\fT,\Sch_{/\dQ_p^\Phi})_{/\bT^\eta_\fp}$.

At the third stage of functoriality, we study the basic correspondence \eqref{eq:qs_basic_correspondence} under the special morphisms. We will complete a commutative diagram in $\Fun(\fK(\rV_n)_\sp^p\times\fT,\Sch_{/\dF_p^\Phi})_{/\rT_\fp}$ as follows
\begin{align}\label{eq:qs_functoriality_basic_correspondence}
\xymatrix{
\rS_\fp(\rV_{n+1},\obj_{n+1})   &&
\rB_\fp(\rV_{n+1},\obj_{n+1})
\ar[rr]^-{\iota_{n+1}}\ar[ll]_-{\pi_{n+1}} &&
\rM_\fp(\rV_{n+1},\obj_{n+1})  \\
\rS_\fp(\rV_n,\obj)_\sp \ar@{}[drr]|{\Box}
\ar[u]_-{\rs_\uparrow}\ar[d]^-{\rs_\downarrow} &&
\rB_\fp(\rV_n,\obj)_\sp \ar[ll]_-{\pi_\sp}
\ar[u]_-{\rb_\uparrow}\ar[d]^-{\rb_\downarrow} \\
\rS_\fp(\rV_n,\obj_n)   &&
\rB_\fp(\rV_n,\obj_n) \ar[rr]^-{\iota_n}\ar[ll]_-{\pi_n}  &&
\rM_\fp(\rV_n,\obj_n) \ar[uu]_-{\rm_\uparrow}
}
\end{align}
in which the lower-left square is Cartesian; and the lower (resp. upper) line is the basic correspondences on $\rM_\fp(\rV_n,\obj_n)$ (resp.\ $\rM_\fp(\rV_{n+1},\obj_{n+1})$) as introduced in Definition \ref{de:basic_correspondence}.

\begin{definition}\label{de:qs_definite_special}
We define a functor
\begin{align*}
\rS_\fp(\rV_n,\obj)_\sp\colon\fK(\rV_n)_\sp^p\times\fT &\to\sfP\Sch'_{/\dF^\Phi_p} \\
\rK^p &\mapsto \rS_\fp(\rV_n,\rK^p)_\sp
\end{align*}
such that for every $S\in\Sch'_{/\dF^\Phi_p}$, $\rS_\fp(\rV_n,\rK^p)_\sp(S)$ is the set of equivalence classes of decuples $(A_0,\lambda_0,\eta_0^p;A^\star,\lambda^\star,\eta^{p\star};A^\star_\natural,\lambda^\star_\natural,\eta^{p\star}_\natural;\delta^\star)$, where
\begin{itemize}[label={\ding{109}}]
  \item $(A_0,\lambda_0,\eta_0^p;A^\star,\lambda^\star,\eta^{p\star})$ is an element in $\rS_\fp(\rV_n,\rK^p_n)(S)$;

  \item $(A_0,\lambda_0,\eta_0^p;A^\star_\natural,\lambda^\star_\natural,\eta^{p\star}_\natural)$ is an element in $\rS_\fp(\rV_{n+1},\rK^p_{n+1})(S)$; and

  \item $\delta^\star\colon A^\star\times A_0\to A^\star_\natural$ is an $O_F$-linear quasi-$p$-isogeny (Definition \ref{de:p_quasi}) such that
      \begin{enumerate}[label=(\alph*)]
          \item $\Ker\delta^\star[p^\infty]$ is contained in $(A^\star\times A_0)[\fp]$;

          \item we have $\lambda^\star\times\varpi\lambda_0=\delta^{\star\vee}\circ\lambda^\star_\natural\circ\delta^\star$; and

          \item the $\rK^p_{n+1}$-orbit of maps $v\mapsto\delta^\star_*\circ(\eta^{p\star}\oplus(\id_{A_0})_*)(v)$ for $v\in\rV_{n+1}\otimes_\dQ\dA^{\infty,p}$ coincides with $\eta^{p\star}_\natural$.
      \end{enumerate}
\end{itemize}
The equivalence relation and the action of morphisms in $\fK(\rV_n)_\sp^p\times\fT$ are defined similarly as in Definition \ref{de:qs_basic_correspondence}.
\end{definition}

We clearly have the forgetful morphism $\rS_\fp(\rV_n,\obj)_\sp\to\rT_\fp$ in $\Fun(\fK(\rV_n)_\sp^p\times\fT,\sfP\Sch'_{/\dF_p^\Phi})$, which is represented by finite and \'{e}tale schemes. By definition, we have the two forgetful morphisms
\begin{align*}
\rs_\downarrow&\colon\rS_\fp(\rV_n,\obj)_\sp\to\rS_\fp(\rV_n,\obj_n),\\
\rs_\uparrow&\colon\rS_\fp(\rV_n,\obj)_\sp\to\rS_\fp(\rV_{n+1},\obj_{n+1})
\end{align*}
in $\Fun(\fK(\rV_n)_\sp^p\times\fT,\Sch_{/\dF_p^\Phi})_{/\rT_\fp}$.

\begin{lem}\label{le:qs_functoriality_source}
We have the following properties concerning $\rs_{\downarrow}$.
\begin{enumerate}
  \item When $n$ is odd, $\rs_{\downarrow}$ is an isomorphism, and the morphism
      \[
      \rs_{\uparrow}\circ \rs_{\downarrow}^{-1}\colon\rS_{\fp}(\rV_n, \obj_{n})\to \rS_{\fp}(\rV_{n+1},\obj_{n+1})
      \]
      is given by the assignment
      \[
      (A_0, \lambda_0, \eta_0^p; A^\star, \lambda^\star, \eta^{p\star})\mapsto (A_0, \lambda_0, \eta_0^p; A^\star\times A_0, \lambda^\star\times\varpi\lambda_0, \eta^{p\star}\times (\id_{A_0})_*).
      \]

  \item When $n$ is even, $\rs_{\downarrow}$ is finite \'etale of degree $p+1$.
\end{enumerate}
\end{lem}

\begin{proof}
Take an object $\rK^p$ of $\fK(\rV_n)_\sp^p$, and a point $x=(A_0,\lambda_0,\eta_0^p;A^\star,\lambda^\star,\eta^{p\star})\in \rS_{\fp}(\rV_n,\rK^p_n)(\kappa)$ for some perfect field $\kappa$ containing $\dF_{p}^{\Phi}$.

For (1), it suffices to show that the fibre $\rs_{\downarrow}^{-1}(x)$ consists of the single point with the extra datum $(A^\star_\natural,\lambda^\star_\natural, \eta^{p\star}_\natural; \delta^\star)=(A^\star\times A_0,\lambda^\star\times\varpi\lambda_0,\eta^{p\star}\times\eta^p_0;\id)$. This follows from the fact that $\delta^\star$ as in Definition~\ref{de:qs_definite_special} induces an equivalence between $(A^\star_\natural, \lambda^\star_\natural, \eta^{p\star}_\natural)$ and $( A^\star\times A_0, \lambda^\star\times\varpi \lambda_0, \eta^{p\star}\times \eta^p_0)$.

For (2), we note first that a point in the fibre $\rs_{\downarrow}^{-1}(x)$ is determined by the quasi-$p$-isogeny $\delta^\star$, which is in turn determined, up to equivalence, by a totally isotropic $(O_F/\fp)$-subgroup of $\Ker(\lambda^\star\times\varpi\lambda_0)$ of order $p^2$. We classify such subgroups by using Dieudonn\'e theory. Let $\cD(A^\star\times A_0)_{\tau_\infty^{\tc}}^\vee$ be the dual lattice of $\cD(A^\star\times A_0)_{\tau_\infty^{\tc}}$ (Notation \ref{no:weil_pairing_dieudonne}) but with respect to the quasi-polarization $\lambda^\star\times\varpi\lambda_0$. The quotient $\sW_x\coloneqq\cD(A^\star\times A_0)_{\tau_\infty^{\tc}}^\vee/\cD(A^\star\times A_0)_{\tau_\infty}$ is $\kappa$-vector space  of dimension $2$ equipped with an induced \emph{nondegenerate} hermitian pairing. Then the hermitian space $\sW_x$ is admissible in the sense of Definition~\ref{de:dl_admissible} with underlying hermitian space over $\dF_{p^2}$ given by $\sW_{x,0}\coloneqq\sW_x^{\tV^{-1}\tF=1}$. Then $\sW_{x,0}$ is an $\dF_{p^2}$-vector space of dimension $2$. By the classical Dieudonn\'e theory for finite group schemes over $\kappa$, the set of totally isotropic $(O_F/\fp)$-subgroups of $\Ker(\lambda^\star\times \varpi\lambda_0)$ of order $p^2$ is in natural bijection with the set of isotropic $\dF_{p^2}$-lines in $\sW_{x,0}$, which has cardinality $p+1$.
\end{proof}

\begin{definition}\label{de:qs_base_special}
We define $\rB_\fp(\rV_n,\obj)_\sp$ to be the fiber product indicated in the following Cartesian diagram
\[
\xymatrix{
\rS_\fp(\rV_n,\obj)_\sp \ar[d]^-{\rs_\downarrow} && \rB_\fp(\rV_n,\obj)_\sp \ar[ll]_-{\pi_\sp}\ar[d]_-{\rb_\downarrow} &&  \\
\rS_\fp(\rV_n,\obj_n) &&  \rB_\fp(\rV_n,\obj_n) \ar[ll]_-{\pi_n}
}
\]
in $\Fun(\fK(\rV_n)_\sp^p\times\fT,\Sch_{/\dF_p^\Phi})_{/\rT_\fp}$.
\end{definition}

\begin{lem}\label{le:qs_functoriality}
The assignment sending an object
\[
((A_0,\lambda_0,\eta_0^p;A,\lambda,\eta^p;A^\star,\lambda^\star,\eta^{p\star};\alpha),
(A_0,\lambda_0,\eta_0^p;A^\star,\lambda^\star,\eta^{p\star};A^\star_\natural,\lambda^\star_\natural,\eta^{p\star}_\natural;\delta^\star))
\]
of $\rB_\fp(\rV_n,\rK^p)_\sp(S)$ to
\begin{align}\label{eq:qs_special}
(A_0,\lambda_0,\eta_0^p;A\times A_0,\lambda\times\lambda_0,\eta^p\oplus(\id_{A_0})_*;
A^\star_\natural,\lambda^\star_\natural,\eta^{p\star}_\natural;\delta^\star\circ(\alpha\times\id_{A_0}))
\end{align}
defines a morphism
\[
\rb_\uparrow\colon\rB_\fp(\rV_n,\obj)_\sp\to\rB_\fp(\rV_{n+1},\obj_{n+1})
\]
in $\Fun(\fK(\rV_n)_\sp^p\times\fT,\Sch_{/\dF_p^\Phi})_{/\rT_\fp}$.
\end{lem}

\begin{proof}
The lemma amounts to showing that \eqref{eq:qs_special} is an object of $\rB_\fp(\rV_{n+1},\rK^p_{n+1})(S)$. Put $\alpha_\natural\coloneqq\delta^\star\circ(\alpha\times\id_{A_0})\colon A\times A_0\to A^\star_\natural$. The only nontrivial condition in Definition \ref{de:qs_basic_correspondence} to check is that $\Ker\alpha_\natural[p^\infty]$ is contained in $(A\times A_0)[\fp]$. For this, we may assume $S=\Spec\kappa$ for a perfect field $\kappa$ containing $\dF_{p}^{\Phi}$.

Consider the following injective maps of Dieudonn\'{e} modules
\[
\cD(A)_\tau\oplus\cD(A_0)_\tau \xrightarrow{\alpha_{*,\tau}\oplus\id} \cD(A^\star)_\tau\oplus\cD(A_0)_\tau
\xrightarrow{\delta^\star_{*,\tau}}\cD(A^\star_\natural)_\tau
\]
for every $\tau\in\Sigma_\infty$. We have the inclusion $\cD(A^\star_\natural)_\tau\subseteq\cD(A^\star)_{\tau^\tc}^\vee\oplus\varpi^{-1}\cD(A_0)_\tau$ (Notation \ref{no:weil_pairing_dieudonne}). Thus, it suffices to show $p\cD(A^\star)^\vee_{\tau^\tc}\subseteq\cD(A)_\tau$ for every $\tau\in\Sigma_\infty$. For $\tau\not\in\{\tau_\infty,\tau_\infty^\tc\}$, we have $\cD(A^\star)^\vee_{\tau^\tc}=\cD(A)_\tau$. It remains to show $p\cD(A^\star)^\vee_{\tau^\tc}\subseteq\cD(A)_\tau$ for $\tau\in\{\tau_\infty,\tau_\infty^\tc\}$. Recall the subspace $H\coloneqq(\breve\alpha_{*,\tau_\infty})^{-1}\omega_{A^\vee/\kappa,\tau_\infty}\subseteq\rH^\dr_1(A^\star/\kappa)_{\tau_\infty}$ from Theorem \ref{th:qs_basic_correspondence}. Under the notation in proof of Theorem \ref{th:qs_basic_correspondence}, since $(\tF H)^\perp\subseteq H$, we have $p\cD(A^\star)_{\tau_\infty^\tc}^\vee\subseteq\tilde{H}$, hence $p\tilde\cD(A^\star)_{\tau_\infty}^\vee\subseteq\tilde{H}^\tc$. Thus, we have
\begin{align*}
p\cD(A^\star)^\vee_{\tau_\infty^\tc}&=p\tV^{-1}(\cD(A^\star)_{\tau_\infty}^\vee)\subseteq \tV^{-1}\tilde{H}^\tc=\cD(A)_{\tau_\infty},\\
p\cD(A^\star)^\vee_{\tau_\infty}&=p\tF(\cD(A^\star)_{\tau_\infty^\tc}^\vee)\subseteq\tF\tilde{H}=\cD(A)_{\tau_\infty^\tc}.
\end{align*}
The lemma follows.
\end{proof}

By the above lemma, we obtain our desired diagram \eqref{eq:qs_functoriality_basic_correspondence}. Moreover, we have the following result.

\begin{proposition}\label{pr:qs_functoriality}
When $n$ is even, the square
\[
\xymatrix{
\rB_\fp(\rV_{n+1},\obj_{n+1}) \ar[rr]^-{\iota_{n+1}} &&  \rM_\fp(\rV_{n+1},\obj_{n+1}) \\
\rB_\fp(\rV_n,\obj)_\sp \ar[rr]^-{\iota_n\circ\rb_\downarrow} \ar[u]^-{\rb_\uparrow}
&& \rM_\fp(\rV_n,\obj_n) \ar[u]_-{\rm_\uparrow}
}
\]
extracted from the diagram \eqref{eq:qs_functoriality_basic_correspondence} is Cartesian.
\end{proposition}

We remark that the above proposition is not correct on the nose when $n$ is odd and at least $3$.

\begin{proof}
The square in the proposition induces a morphism
\begin{align*}
\iota_\sp\colon\rB_\fp(\rV_n,\obj)_\sp\to\rB_\fp(\rV_{n+1},\obj_{n+1})\times_{\rM_\fp(\rV_{n+1},\obj_{n+1})}\rM_\fp(\rV_n,\obj_n).
\end{align*}
We need to prove that $\iota_\sp$ is an isomorphism. By Theorem \ref{th:qs_basic_correspondence}, we know that $\iota_\sp$ is locally for the Zariski topology on the source  a  closed immersion, such that both the  source  and the target are smooth. Thus, it suffices to show that for a given algebraically closed field $\kappa$ containing $\dF_p^\Phi$, we have that
\begin{enumerate}
  \item $\iota_\sp(\kappa)$ is an isomorphism in $\Fun(\fK(\rV_n)_\sp^p\times\fT,\Set)$; and

  \item for every $\rK^p\in\fK(\rV_n)^p_\sp$ and every $x\in\rB_\fp(\rV_n,\rK^p)_\sp(\kappa)$, the induced diagram
     \begin{align}\label{eq:qs_functoriality3}
     \xymatrix{
     \cT_{\rb_\uparrow(x)} \ar[rr]^-{\iota_{n+1*}} && \cT_{\iota_{n+1}(\rb_\uparrow(x))} \\
     \cT_x \ar[rr]^-{\iota_{n*}\circ\rb_{\downarrow*}}\ar[u]^-{\rb_{\uparrow*}}  && \cT_{\iota(\rb_\downarrow(x))} \ar[u]_-{\rm_{\uparrow*}}
     }
     \end{align}
     of tangent spaces is a Cartesian square of $\kappa$-modules.
\end{enumerate}

For (1), we take an object $\rK^p\in\fK(\rV_n)_\sp^p$ and construct an inverse of $\iota_\sp(\kappa)$. Take a point
\[
(A_0,\lambda_0,\eta_0^p;A,\lambda,\eta^p;A^\star_\natural,\lambda^\star_\natural,\eta^{p\star}_\natural;\alpha_\natural)
\]
in the target of $\iota_\sp(\kappa)$. Then $\alpha_\natural$ induces an inclusion
\[
\cD(A)_\tau\oplus \cD(A_0)_\tau\subseteq\cD(A^\star_\natural)_\tau
\]
of Dieudonn\'{e} modules, which is an equality if $\tau\not\in\{\tau_\infty,\tau_\infty^\tc\}$. We put
\[
\cD_{A^\star}\coloneqq\bigoplus_{\tau\in\Sigma_\infty}\cD_{A^\star,\tau}
\]
where $\cD_{A^\star,\tau}=\cD(A)_\tau$ for $\tau\not\in\{\tau_\infty,\tau_\infty^\tc\}$ and $\cD_{A^\star,\tau}=\cD(A^\star_\natural)_\tau\cap p^{-1}\cD(A)_\tau$ for $\tau\in\{\tau_\infty,\tau_\infty^\tc\}$. Then $\cD_{A^\star}$ is a Dieudonn\'{e} module containing $\cD(A)$. By the Dieudonn\'{e} theory, there is an $O_F$-abelian scheme $A^\star$ over $\kappa$ with $\cD(A^\star)_\tau=\cD_{A^\star,\tau}$ for every $\tau\in\Sigma_\infty$, and an $O_F$-linear isogeny $\alpha\colon A\to A^\star$ inducing the inclusion of Dieudonn\'{e} modules $\cD(A)\subseteq\cD(A^\star)$. We factors $\alpha_\natural$ as
\[
A\times A_0\xrightarrow{\alpha\times\id_{A_0}}A^\star\times A_0\xrightarrow{\delta^\star} A^\star_\natural.
\]
It is clear that there is a unique quasi-polarization $\lambda^\star$ of $A^\star$ such that $\lambda^\star\times\varpi\lambda_0=\delta^{\star\vee}\circ\lambda^\star_\natural\circ\delta^\star$. Let $\eta^{p\star}$ be the $\rK^p_n$-level structure induced from $\eta^p$ under $\alpha$. We claim that the datum
\[
((A_0,\lambda_0,\eta_0^p;A,\lambda,\eta^p;A^\star,\lambda^\star,\eta^{p\star};\alpha),
(A_0,\lambda_0,\eta_0^p;A^\star,\lambda^\star,\eta^{p\star};A^\star_\natural,\lambda^\star_\natural,\eta^{p\star}_\natural;\delta^\star))
\]
gives rise to an element in $\rB_\fp(\rV_n,\rK^p)_\sp(\kappa)$.
 It suffices to show that $(A_0,\lambda_0,\eta_0^p;A^\star,\lambda^\star,\eta^{p\star})$ is an element in $\rS_\fp(\rV_{n},\rK^p_{n})(\kappa)$. Moreover precisely, we need to show that
\begin{enumerate}
  \item [(1--1)] the $O_F$-abelian scheme $A^\star$ has signature type $n\Phi$; and

  \item [(1--2)] $\Ker\lambda^\star[p^\infty]$ is contained in $A^\star[\fp]$ of degree $p^2$.
\end{enumerate}
To prove these, we add two auxiliary properties
\begin{enumerate}
  \item [(1--3)] the composite map $\cD(A^\star_\natural)_\tau\subseteq p^{-1}\cD(A)_\tau\oplus p^{-1}\cD(A_0)_\tau\to p^{-1}\cD(A_0)_\tau$ is surjective for $\tau\in\{\tau_\infty,\tau_\infty^\tc\}$; and

  \item [(1--4)] the cokernel of the inclusion $\cD(A^\star)_\tau\oplus\cD(A_0)_\tau\subseteq\cD(A^\star_\natural)_\tau$ is isomorphic to $\kappa$ for $\tau\in\{\tau_\infty,\tau_\infty^\tc\}$.
\end{enumerate}

For (1--3), if not surjective, then we have $\cD(A^\star_\natural)_\tau\subseteq p^{-1}\cD(A)_\tau\oplus \cD(A_0)_\tau$ for both $\tau\in\{\tau_\infty,\tau_\infty^\tc\}$. As $\varpi\lambda\times\varpi\lambda_0=\alpha_\natural^\vee\circ\lambda^\star_\natural\circ\alpha_\natural$, this contradicts with the fact that $\lambda^\star_\natural$ is $p$-principal.

For (1--4), it follows (1--3) and the fact that the kernel of $\cD(A^\star_\natural)_\tau\to p^{-1}\cD(A_0)_\tau$ is $\cD(A^\star)_\tau$ for $\tau\in\{\tau_\infty,\tau_\infty^\tc\}$.

For (1--1), it amounts to showing that $\tF\colon\cD(A^\star)_\tau\to\cD(A^\star)_{\tau^\tc}$ is an isomorphism for every $\tau\in\Phi$. This is obvious for $\tau\neq\tau_\infty$. When $\tau=\tau_\infty$, this follows from (1.4) and the fact that both $\tF\colon\cD(A^\star_\natural)_\tau\to\cD(A^\star_\natural)_{\tau^\tc}$ and $\tF\colon\cD(A_0)_\tau\to\cD(A_0)_{\tau^\tc}$ are isomorphisms.

For (1--2), it follows from (1--4) and the fact that $\lambda^\star_\natural$ is $p$-principal.

Thus, (1) is proved.

For (2), the diagram \eqref{eq:qs_functoriality3} is identified with
\[
\xymatrix{
\Hom_\kappa\(\omega_{A^\vee,\tau_\infty},\Ker\alpha_{\natural*,\tau_\infty}/\omega_{A^\vee,\tau_\infty}\)  \ar[r] &
\Hom_\kappa\(\omega_{A^\vee,\tau_\infty},\rH^\dr_1(A\times A_0)_{\tau_\infty}/\omega_{A^\vee,\tau_\infty}\) \\
\Hom_\kappa\(\omega_{A^\vee,\tau_\infty},\Ker\alpha_{*,\tau_\infty}/\omega_{A^\vee,\tau_\infty}\)  \ar[r]\ar[u] &
\Hom_\kappa\(\omega_{A^\vee,\tau_\infty},\rH^\dr_1(A)_{\tau_\infty}/\omega_{A^\vee,\tau_\infty}\) \ar[u]
}
\]
by Theorem \ref{th:qs_moduli_scheme} and Theorem \ref{th:qs_basic_correspondence}. However, it is an easy consequence of (1--3) that $\Ker\alpha_{\natural*,\tau_\infty}\cap\rH^\dr_1(A)_{\tau_\infty}=\Ker\alpha_{*,\tau_\infty}$. Thus, the above diagram is Cartesian; and (2) follows.
\end{proof}

At the fourth stage of functoriality, we compare the special morphisms for basic correspondences and for Deligne--Lusztig varieties. Take a point
\[
s^\star=(A_0,\lambda_0,\eta_0^p;A^\star,\lambda^\star,\eta^{p\star};A^\star_\natural,\lambda^\star_\natural,\eta^{p\star}_\natural;\delta^\star)
\in\rS_\fp(\rV_n,\rK^p)_\sp(\kappa)
\]
for a field $\kappa$ containing $\dF_p^\Phi$. Put
\[
s^\star_n\coloneqq\rs_\downarrow(s^\star),\quad s^\star_{n+1}\coloneqq\rs_\uparrow(s^\star);
\]
and denote by $\rB_{s^\star}$, $\rB_{s^\star_n}$, and $\rB_{s^\star_{n+1}}$ their preimages under $\pi_\sp$, $\pi_n$, and $\pi_{n+1}$ in \eqref{eq:qs_functoriality_basic_correspondence}, respectively. By Lemma \ref{le:qs_dl}, we have admissible pairs $(\sV_{s^\star_n},\{\;,\;\}_{s^\star_n})$ and $(\sV_{s^\star_{n+1}},\{\;,\;\}_{s^\star_{n+1}})$. As in Construction \ref{cs:dl_special_morphism}, we extend the pair $(\sV_{s^\star_n},\{\;,\;\}_{s^\star_n})$ to $(\sV_{s^\star_n,\sharp},\{\;,\;\}_{s^\star_n,\sharp})$. Then the homomorphism $\delta^\star\colon A^\star\times A_0\to A^\star_\natural$ induces a $\kappa$-linear map
\[
\delta_{s^\star}\colon\sV_{s^\star_n,\sharp}\to\sV_{s^\star_{n+1}}
\]
satisfying $\{\delta_{s^\star}(x),\delta_{s^\star}(y)\}_{s^\star_{n+1}}=\{x,y\}_{s^\star_n,\sharp}$ for every $x,y\in\sV_{s^\star_n,\sharp}$. By Construction \ref{cs:dl_special_morphism}, we obtain a morphism
\[
\delta_{s^\star\uparrow}\colon
\DL_{s^\star_n}=\DL(\sV_{s^\star_n},\{\;,\;\}_{s^\star_n},\ceil{\tfrac{n+1}{2}})\to
\DL_{s^\star_{n+1}}=\DL(\sV_{s^\star_{n+1}},\{\;,\;\}_{s^\star_{n+1}},\ceil{\tfrac{n+2}{2}})
\]
of the corresponding Deligne--Lusztig varieties.

\begin{proposition}\label{pr:qs_functoriality_dl}
Let the notation be as above. The following diagram
\[
\xymatrix{
\rB_{s^\star_{n+1}} \ar[rr]^-{\zeta_{s^\star_{n+1}}}_-{\simeq} && \DL_{s^\star_{n+1}} \\
\rB_{s^\star} \ar[rr]^-{\zeta_{s^\star_n}\circ\rb_\downarrow}_-{\simeq} \ar[u]^-{\rb_\uparrow}
&& \DL_{s^\star_n} \ar[u]_-{\delta_{s^\star\uparrow}}
}
\]
in $\Sch_{/\kappa}$ commutes, where $\zeta_{s^\star_n}$ and $\zeta_{s^\star_{n+1}}$ are the isomorphisms in Theorem \ref{th:qs_basic_correspondence}(3). In particular, $\rb_\uparrow\colon\rB_{s^\star}\to\rB_{s^\star_{n+1}}$ is an isomorphism if $n$ is odd, and is a regular embedding of codimension one if $n$ is even.
\end{proposition}

\begin{proof}
Note that by Lemma \ref{le:qs_functoriality_source}, the restricted morphism $\rb_\downarrow\colon\rB_{s^\star}\to\rB_{s^\star_n}$ is an isomorphism. Thus, the last claim follows from the commutativity and Proposition \ref{pr:dl_special_morphism}.

When $n$ is odd, the commutativity is obvious. When $n$ is even, it suffices to show that for every point
\[
(A_0,\lambda_0,\eta_0^p;A,\lambda,\eta^p;A^\star,\lambda^\star,\eta^{p\star};\alpha)\in\rB_{s^\star}(S),
\]
we have
\begin{align}\label{eq:qs_functoriality_dl}
\delta^\star_{*,\tau_\infty}\((\breve\alpha_{*,\tau_\infty})^{-1}\omega_{A^\vee/S,\tau_\infty}\oplus\rH^\dr_1(A_0/S)_{\tau_\infty}\)
=(\breve\alpha_{\natural*,\tau_\infty})^{-1}\omega_{A^\vee\times A_0^\vee/S,\tau_\infty}
\end{align}
in view of the diagram
\[
\xymatrix{
A\times A_0 \ar@{=}[r] \ar[d]_-{\alpha\times\id_{A_0}} & A\times A_0
\ar[d]^-{\alpha_\natural\coloneqq\delta^\star\circ(\alpha\times\id_{A_0})} \\
A^\star\times A_0 \ar[r]^-{\delta^\star}\ar[d]_-{\breve\alpha\times\varpi\id_{A_0}} & A^\star_\natural \ar[d]^-{\breve\alpha_\natural} \\
A\times A_0 \ar@{=}[r] & A\times A_0
}
\]
in which $\breve\alpha\circ\alpha=\varpi\cdot\id_A$ and $\breve\alpha_\natural\circ\alpha_\natural=\varpi\cdot\id_{A\times A_0}$. Since both sides of \eqref{eq:qs_functoriality_dl} have the same rank, it suffices to show that
\[
\breve\alpha_{\natural*,\tau_\infty}
\(\delta^\star_{*,\tau_\infty}\((\breve\alpha_{*,\tau_\infty})^{-1}\omega_{A^\vee/S,\tau_\infty}\oplus\rH^\dr_1(A_0/S)_{\tau_\infty}\)\)
\subseteq\omega_{A^\vee\times A_0^\vee/S,\tau_\infty},
\]
which is obvious as $\varpi$ annihilates $\rH^\dr_1(A_0/S)_{\tau_\infty}$. The proposition is proved.
\end{proof}

At the final stage of functoriality, we relate the special morphisms for sources of basic correspondences to Shimura sets under the uniformization map $\upsilon$ \eqref{eq:qs_uniformization}.

\begin{notation}\label{no:qs_uniformization}
As in Definition \ref{de:qs_uniformization_data}, we choose a definite uniformization datum $(\rV^\star_n,\ti_n,\{\Lambda^\star_{n,\fq}\}_{\fq\mid p})$ for $\rV$. We also fix a definite uniformization datum $(\rV^\star_{n+1},\ti_{n+1},\{\Lambda^\star_{n+1,\fq}\}_{\fq\mid p})$ for $\rV_{n+1}$ satisfying
\begin{itemize}[label={\ding{109}}]
  \item $\rV^\star_{n+1}=(\rV^\star_n)_\sharp$ and $\ti_{n+1}=(\ti_n)_\sharp$;

  \item $\Lambda^\star_{n+1,\fq}=(\Lambda^\star_{n,\fq})_\sharp$ for $\fq\neq\fp$; and

  \item $(\Lambda^\star_{n,\fp})_\sharp\subseteq\Lambda^\star_{n+1,\fp}\subseteq p^{-1}(\Lambda^\star_{n,\fp})_\sharp^\vee$.
\end{itemize}
Let $\rK^\star_{n+1,\fq}$ be the stabilizer of $\Lambda^\star_{n+1,\fq}$ for every $\fq$ over $p$; and put $\rK^\star_{n+1,p}\coloneqq\prod_{\fq\mid p}\rK^\star_{n+1,\fq}$. Moreover, we put $\rK^\star_{\sp,\fp}\coloneqq\rK^\star_{n,\fp}\cap\rK^\star_{n+1,\fp}$ (as a subgroup of $\rK^\star_{n,\fp}$) and $\rK^\star_{\sp,p}\coloneqq\rK^\star_{\sp,\fp}\times\prod_{\fq\neq\fp}\rK^\star_{n,\fq}$.
\end{notation}

\begin{remark}\label{re:qs_uniformization}
When $n$ is odd, since $(\Lambda^\star_{n,\fp})^\vee=p\Lambda^\star_{n,\fp}$, we must have $\Lambda^\star_{n+1,\fp}=(\Lambda^\star_{n,\fp})_\sharp$ as well, hence $\rK^\star_{\sp,p}=\rK^\star_{n,p}$. When $n$ is even, the number of choices of $\Lambda^\star_{n+1,\fp}$ is $p+1$.
\end{remark}

Similar to Construction \ref{cs:qs_uniformization}, we may construct a uniformization map
\begin{align}\label{eq:qs_uniformization_basic}
\upsilon_\sp\colon\rS_\fp(\rV_n,\obj)_\sp(\ol\dF_p)
\to\Sh(\rV^\star_n,(\ti_n\obj_n)\rK^\star_{\sp,p})\times\rT_\fp(\ol\dF_p)
\end{align}
in $\Fun(\fK(\rV_n)_\sp^p\times\fT,\Set)_{/\rT_\fp(\ol\dF_p)}$ which is an isomorphism, whose details we leave to the readers.

\begin{proposition}\label{pr:qs_uniformization_basic}
The following diagram
\[
\xymatrix{
\rS_\fp(\rV_{n+1},\obj_{n+1})(\ol\dF_p) \ar[rr]^-{\upsilon_{n+1}}_-{\eqref{eq:qs_uniformization}} &&
\Sh(\rV^\star_{n+1},(\ti_{n+1}\obj_{n+1})\rK^\star_{n+1,p})\times\rT_\fp(\ol\dF_p) \\
\rS_\fp(\rV_n,\obj)_\sp(\ol\dF_p) \ar[rr]^-{\upsilon_\sp}_-{\eqref{eq:qs_uniformization_basic}}
\ar[u]^-{\rs_\uparrow(\ol\dF_p)}\ar[d]_-{\rs_\downarrow(\ol\dF_p)} &&
\Sh(\rV^\star_n,(\ti_n\obj_n)\rK^\star_{\sp,p})\times\rT_\fp(\ol\dF_p)
\ar[u]_-{\sh^\star_\uparrow\times\id}\ar[d]^-{\sh^\star_\downarrow\times\id} \\
\rS_\fp(\rV_n,\obj_n)(\ol\dF_p) \ar[rr]^-{\upsilon_n}_-{\eqref{eq:qs_uniformization}}  &&
\Sh(\rV^\star_n,(\ti_n\obj_n)\rK^\star_{n,p})\times\rT_\fp(\ol\dF_p)
}
\]
in $\Fun(\fK(\rV_n)_\sp^p\times\fT,\Set)_{/\rT_\fp(\ol\dF_p)}$ commutes, where $\sh^\star_\downarrow$ and $\sh^\star_\uparrow$ are obvious maps on Shimura sets. Moreover, the induced actions of $\Gal(\ol\dF_p/\dF_p^\Phi)$ on all terms on the right-hand side factor through the projection to the factor $\rT_\fp(\ol\dF_p)$.
\end{proposition}

\begin{proof}
The commutativity follows directly from definition. The proof of the last claim is same to Proposition \ref{pr:qs_uniformization}.
\end{proof}

\subsection{Second geometric reciprocity law}
\label{ss:qs_reciprocity}

In this subsection, we state and prove a theorem we call \emph{second geometric reciprocity law}, which can be regarded a geometric template for the second explicit reciprocity law studied in \S\ref{ss:second_reciprocity} once throw the automorphic input.

We keep the setup in \S\ref{ss:qs_functoriality}. However, we allow $\obj=(\obj_n,\obj_{n+1})$ to be an object of $\fK(\rV_n)^p\times\fK(\rV_{n+1})^p$, rather than $\fK(\rV_n)^p_\sp$. Denote by $n_0$ and $n_1$ the unique even and odd numbers in $\{n,n+1\}$, respectively. Write $n_0=2r_0$ and $n_1=2r_1+1$ for unique integers $r_0,r_1\geq 1$. In particular, we have $n=r_0+r_1$. Let $L$ be a $p$-coprime coefficient ring.

To ease notation, we put $\rX^?_{n_\alpha}\coloneqq\rX^?_\fp(\rV_{n_\alpha},\obj_{n_\alpha})$ for meaningful triples $(\rX,?,\alpha)\in\{\bM,\rM,\rB,\rS\}\times\{\;,\eta\}\times\{0,1\}$.

\begin{construction}\label{cs:qs_reciprocity}
We construct two maps and two graphs.
\begin{enumerate}
  \item For every integers $i,j$, we define
     \[
     \loc'_\fp\colon\rH^i_\et(\Sh(\rV_{n_0},\obj_{n_0}\rK_{n_0,p})\times_{\Spec{F}}\Sh(\rV_{n_1},\obj_{n_1}\rK_{n_1,p}),L(j))
     \to\rH^i_\fT(\rM_{n_0}\times_{\rT_\fp}\rM_{n_1},L(j))
     \]
     to be the composition of the localization map
     \begin{align*}
     \loc_\fp&\colon\rH^i_\et(\Sh(\rV_{n_0},\obj_{n_0}\rK_{n_0,p})\times_{\Spec{F}}\Sh(\rV_{n_1},\obj_{n_1}\rK_{n_1,p}),L(j)) \\
     &\to\rH^i_\et((\Sh(\rV_{n_0},\obj_{n_0}\rK_{n_0,p})\times_{\Spec{F}}\Sh(\rV_{n_1},\obj_{n_1}\rK_{n_1,p}))\otimes_F\dQ_p^\Phi,L(j)),
     \end{align*}
     the pullback map
     \[
     \rH^i_\et((\Sh(\rV_{n_0},\obj_{n_0}\rK_{n_0,p})\times_{\Spec{F}}\Sh(\rV_{n_1},\obj_{n_1}\rK_{n_1,p}))\otimes_F\dQ_p^\Phi,L(j))
     \to\rH^i_\fT(\bM^\eta_{n_0}\times_{\bT^\eta_\fp}\bM^\eta_{n_1},L(j))
     \]
     induced from \eqref{eq:qs_moduli_scheme_shimura}, and the isomorphism
     \[
     \rH^i_\fT(\rM_{n_0}\times_{\rT_\fp}\rM_{n_1},\rR\Psi L(j))
     \xrightarrow{\sim}\rH^i_\fT(\rM_{n_0}\times_{\rT_\fp}\rM_{n_1},L(j))
     \]
     due to the fact $L\simeq\rR\Psi L$ by Theorem \ref{th:qs_moduli_scheme}.

  \item Analogous to Construction \ref{cs:qs_incidence}, we define the map
     \begin{align*}
     \inc^{\star,\star}_!&\colon L[\Sh(\rV^\star_{n_0},(\ti_{n_0}\obj_{n_0})\rK^\star_{n_0,p})]\otimes_LL[\Sh(\rV^\star_{n_1},(\ti_{n_1}\obj_{n_1})\rK^\star_{n_1,p})] \\
     &\xrightarrow{\sim}\rH^0_\fT(\rS_{n_0},L)\otimes_L\rH^0_\fT(\rS_{n_1},L)
     =\rH^0_\fT(\rS_{n_0}\times_{\rT_\fp}\rS_{n_1},L) \\
     &\xrightarrow{(\pi_{n_0}\times\pi_{n_1})^*}\rH^0_\fT(\rB_{n_0}\times_{\rT_\fp}\rB_{n_1},L)
     \xrightarrow{(\iota_{n_0}\times\iota_{n_1})_!}\rH^{2n}_\fT(\rM_{n_0}\times_{\rT_\fp}\rM_{n_1},L(n))
     \end{align*}
     in $\Fun(\fK(\rV_n)^p\times\fK(\rV_{n+1})^p,\Mod(L))$.
\end{enumerate}
Suppose that $\obj$ is taken in the subcategory $\fK(\rV_n)^p_\sp$.
\begin{enumerate}\setcounter{enumi}{2}
  \item We define $\graph\Sh(\rV_n,\obj_n\rK_{n,p})$ to be the graph of the morphism $\sh_\uparrow$ \eqref{eq:qs_functoriality_shimura}, as a closed subscheme of $\Sh(\rV_{n_0},\obj_{n_0}\rK_{n_0,p})\times_{\Spec{F}}\Sh(\rV_{n_1},\obj_{n_1}\rK_{n_1,p})$, which gives rise to a class
      \begin{align*}
      [\graph\Sh(\rV_n,\obj_n\rK_{n,p})]\in
      \rH^{2n}_\et(\Sh(\rV_{n_0},\obj_{n_0}\rK_{n_0,p})\times_{\Spec{F}}\Sh(\rV_{n_1},\obj_{n_1}\rK_{n_1,p}),L(n))
      \end{align*}
      by the absolute cycle class map.

  \item We define $\graph\Sh(\rV^\star_n,(\ti_n\obj_n)\rK^\star_{\sp,p})$ to be the graph of the correspondence $(\sh^\star_\downarrow,\sh^\star_\uparrow)$, which is a subset of $\Sh(\rV^\star_{n_0},(\ti_{n_0}\obj_{n_0})\rK^\star_{n_0,p})\times\Sh(\rV^\star_{n_1},(\ti_{n_1}\obj_{n_1})\rK^\star_{n_1,p})$.
\end{enumerate}
\end{construction}

The following theorem, which we call the \emph{second geometric reciprocity law}, relates the class $[\graph\Sh(\rV_n,\obj_n\rK_{n,p})]$ with an explicit class coming from the Shimura set.

\begin{theorem}[Second geometric reciprocity law]\label{th:qs_reciprocity}
Suppose that $\obj$ is taken in the subcategory $\fK(\rV_n)^p_\sp$. We have
\[
\resizebox{\hsize}{!}{
\xymatrix{
\tT^\star_{n_1,\fp}.(\id\times\pi_{n_1})_!(\id\times\iota_{n_1})^*\loc'_\fp\([\graph\Sh(\rV_n,\obj_n\rK_{n,p})]\)
=(\id\times\pi_{n_1})_!(\id\times\iota_{n_1})^*\inc^{\star,\star}_!(\CF_{\graph\Sh(\rV^\star_n,(\ti_n\obj_n)\rK^\star_{\sp,p})})
}
}
\]
in $\rH^{2r_0}_\fT(\rM_{n_0}\times_{\rT_\fp}\rS_{n_1},L(r_0))$, where $\tT^\star_{n_1,\fp}\in\dT_{n_1,\fp}$ is the Hecke operator appearing in Theorem \ref{th:qs_tate}.
\end{theorem}

Note that by Proposition \ref{pr:qs_definite_hecke} and Remark \ref{re:qs_tate_hecke}, $\rH^{2r_0}_\fT(\rM_{n_0}\times_{\rT_\fp}\rS_{n_1},L(r_0))$ is a $\dT_{n_1,\fp}$-module. For the readers' convenience, we illustrate the identity in the above theorem through the following diagram
\[
\resizebox{\hsize}{!}{
\xymatrix{
\rH^{2n}_\et(\Sh(\rV_{n_0},\obj_{n_0}\rK_{n_0,p})\times_{\Spec{F}}\Sh(\rV_{n_1},\obj_{n_1}\rK_{n_1,p}),L(n)) \ar[r]^-{\loc'_\fp}
& \rH^{2n}_\fT(\rM_{n_0}\times_{\rT_\fp}\rM_{n_1},L(n)) \ar[d]^-{(\id\times\iota_{n_1})^*}
& L[\Sh(\rV^\star_{n_0},(\ti_{n_0}\obj_{n_0})\rK^\star_{n_0,p})]\otimes_LL[\Sh(\rV^\star_{n_1},(\ti_{n_1}\obj_{n_1})\rK^\star_{n_1,p})]
\ar[l]_-{\inc^{\star,\star}_!} \\
[\graph\Sh(\rV_n,\obj_n\rK_{n,p})] \ar@{}[u]|-{\rotatebox[origin=c]{90}{$\in$}} \ar@{|->}@/_2pc/[rdd]
& \rH^{2n}_\fT(\rM_{n_0}\times_{\rT_\fp}\rB_{n_1},L(n)) \ar[d]^-{(\id\times\pi_{n_1})_!}
& \CF_{\graph\Sh(\rV^\star_n,(\ti_n\obj_n)\rK^\star_{\sp,p})}  \ar@{}[u]|-{\rotatebox[origin=c]{90}{$\in$}} \ar@{|->}@/^2pc/[ldd] \\
& \rH^{2r_0}_\fT(\rM_{n_0}\times_{\rT_\fp}\rS_{n_1},L(r_0)) \\
& \tT^\star_{n_1,\fp}.\cdots=\cdots \ar@{}[u]|-{\rotatebox[origin=c]{90}{$\in$}}
}
}
\]

\begin{proof}
We denote
\[
\bm_\graph\colon\bM_n\to\bM_n\times_{\bT_\fp}\bM_{n+1}=\bM_{n_0}\times_{\bT_\fp}\bM_{n_1}
\]
the diagonal morphism of the correspondence $(\id,\bm_\uparrow)$ \eqref{eq:qs_functoriality_moduli_scheme} in $\Fun(\fK(\rV_n)_\sp^p\times\fT,\Sch_{/\dZ_p^\Phi})_{/\bT_\fp}$. Then we have the identity
\[
\loc'_\fp\([\graph\Sh(\rV_n,\obj_n\rK_{n,p})]\)=\rm_{\graph!}[\rM_n]\in\rH^{2n}_\fT(\rM_n\times_{\rT_\fp}\rM_{n+1},L(n))
\]
by the commutative diagram \eqref{eq:qs_functoriality_moduli_scheme1}.

Put $\rB_\sp\coloneqq\rB_\fp(\rV_n,\obj)_\sp$ for short, and denote
\[
\rb_\graph\coloneqq(\rb_\downarrow,\rb_\uparrow)\colon\rB_\sp\to\rB_n\times_{\rT_\fp}\rB_{n+1}=\rB_{n_0}\times_{\rT_\fp}\rB_{n_1}
\]
the diagonal morphism of the correspondence $(\rb_\downarrow,\rb_\uparrow)$. By Proposition \ref{pr:qs_functoriality} (resp.\ Lemma \ref{le:qs_functoriality_source}) when $n=n_0$ (resp.\ $n=n_1$), the following commutative diagram
\[
\xymatrix{
\rB_\sp \ar[rr]^-{(\iota_{n_0}\times\id)\circ\rb_\graph}\ar[d]_-{\iota_n\circ\rb_\downarrow}
&& \rM_{n_0}\times_{\rT_\fp}\rB_{n_1} \ar[d]^-{\id\times\iota_{n_1}} \\
\rM_n \ar[rr]^-{\rm_\graph}  && \rM_{n_0}\times_{\rT_\fp}\rM_{n_1}
}
\]
is Cartesian. Then by Proper Base Change, we have
\begin{align*}
\tT^\star_{n_1,\fp}.(\id\times\pi_{n_1})_!(\id\times\iota_{n_1})^*\rm_{\graph!}[\rM_n]
&=\tT^\star_{n_1,\fp}.(\id\times\pi_{n_1})_!((\iota_{n_0}\times\id)\circ\rb_\graph)_!(\iota_n\circ\rb_\downarrow)^*[\rM_n] \\
&=\tT^\star_{n_1,\fp}.(\id\times\pi_{n_1})_!((\iota_{n_0}\times\id)\circ\rb_\graph)_![\rB_\sp].
\end{align*}
The commutative diagram
\[
\xymatrix{
\rB_\sp \ar[rr]^-{(\iota_{n_0}\times\id)\circ\rb_\graph}\ar[d]_-{(\id\times\pi_{n_1})\circ\rb_\graph} && \rM_{n_0}\times_{\rT_\fp}\rB_{n_1}
\ar[d]^-{\id\times\pi_1} \\
\rB_{n_0}\times_{\rT_\fp}\rS_{n_1} \ar[rr]^-{\iota_{n_0}\times\id}
&& \rM_{n_0}\times_{\rT_\fp}\rS_{n_1}
}
\]
implies the identity
\[
\tT^\star_{n_1,\fp}.(\id\times\pi_{n_1})_!((\iota_{n_0}\times\id)\circ\rb_\graph)_![\rB_\sp]
=\tT^\star_{n_1,\fp}.(\iota_{n_0}\times\id)_!((\id\times\pi_{n_1})\circ\rb_\graph)_![\rB_\sp].
\]
Now by the definition of $\rB_\sp$ (Definition \ref{de:qs_base_special}), we have
\[
((\id\times\pi_{n_1})\circ\rb_\graph)_![\rB_\sp]=(\pi_{n_0}\times\id)^*(\CF_{\graph\Sh(\rV^\star_n,(\ti_n\obj_n)\rK^\star_{\sp,p})}).
\]
In all, we have
\[
\tT^\star_{n_1,\fp}.(\id\times\pi_{n_1})_!(\id\times\iota_{n_1})^*\rm_{\graph!}[\rM_n]
=(\iota_{n_0}\times\id)_!(\pi_{n_0}\times\id)^*(\tT^\star_{n_1,\fp}.\CF_{\graph\Sh(\rV^\star_n,(\ti_n\obj_n)\rK^\star_{\sp,p})}),
\]
which, by Theorem \ref{th:qs_tate}, is equal to
\begin{align*}
&\quad(\iota_{n_0}\times\id)_!(\pi_{n_0}\times\id)^*(\id\times\pi_{n_1})_!(\id\times\iota_{n_1})^*(\id\times\iota_{n_1})_!(\id\times\pi_{n_1})^*
(\CF_{\graph\Sh(\rV^\star_n,(\ti_n\obj_n)\rK^\star_{\sp,p})}) \\
&=(\id\times\pi_{n_1})_!(\id\times\iota_{n_1})^*\inc^{\star,\star}_!(\CF_{\graph\Sh(\rV^\star_n,(\ti_n\obj_n)\rK^\star_{\sp,p})}).
\end{align*}
The theorem follows.
\end{proof}

\section{Unitary moduli schemes: semistable case}
\label{ss:ns}

In this section, we define and study a certain semistable integral moduli scheme whose generic fiber is the product of a unitary Shimura variety and an auxiliary CM moduli. Since the materials in this section are strictly in the linear order, we will leave the summary of contents to each subsection.

\subsection{Initial setup}
\label{ss:ns_initial}

We fix a special inert prime (Definition \ref{de:special_inert}) $\fp$ of $F^+$ (with the underlying rational prime $p$). We take the prescribed subring $\dP$ in Definition \ref{de:unitary_abelian_scheme} to be $\dZ_{(p)}$. We choose following data
\begin{itemize}[label={\ding{109}}]
  \item a CM type $\Phi$ containing $\tau_\infty$;

  \item a rational skew-hermitian space $\rW_0$ over $O_F\otimes\dZ_{(p)}$ of rank $1$ and type $\Phi$ (Definition \ref{de:skew_hermitian_type});

  \item a neat open compact subgroup $\rK_0^p\subseteq\rT_0(\dA^{\infty,p})$;

  \item an isomorphism $\ol\dQ_p\simeq\dC$ that induces the place $\fp$ of $F^+$;

  \item an element $\varpi\in O_{F^+}$ that is totally positive and satisfies $\val_\fp(\varpi)=1$, and $\val_\fq(\varpi)=0$ for every prime $\fq\neq\fp$ of $F^+$ above $p$.
\end{itemize}
We adopt Notation \ref{no:p_notation}. In particular, $\dF^\Phi_p$ contains $\dF_{p^2}$. Since the argument below is insensitive to the choices of $\rW_0$ and $\rK_0^p$, we will not include them in all notations. However, we will keep the prime $\fp$ in notations as later in application, we need to choose different primes in a crucial step. Put $\bT_\fp\coloneqq\bT_p(\rW_0,\rK_0^p)\otimes_{O_{F_\Phi}\otimes\dZ_{(p)}}\dZ^\Phi_p$.

\subsection{Construction of moduli schemes}
\label{ss:ns_moduli_scheme}

In this subsection, we construct our initial moduli schemes. We start from the datum $(\rV^\circ,\{\Lambda^\circ_\fq\}_{\fq\mid p})$, where
\begin{itemize}[label={\ding{109}}]
  \item $\rV^\circ$ is a standard \emph{definite} hermitian space (Definition \ref{de:standard_hermitian_space}) over $F$ of rank $N\geq 1$, and

  \item for every prime $\fq$ of $F^+$ above $p$, a self-dual $O_{F_\fq}$-lattice $\Lambda^\circ_\fq$ in $\rV^\circ\otimes_FF_\fq$.
\end{itemize}

\begin{definition}\label{de:ns_moduli_scheme}
We define a functor
\begin{align*}
\bM_\fp(\rV^\circ,\obj)\colon\fK(\rV^\circ)^p\times\fT &\to\sfP\Sch'_{/\dZ^\Phi_p} \\
\rK^{p\circ} &\mapsto \bM_\fp(\rV^\circ,\rK^{p\circ})
\end{align*}
such that for every $S\in\Sch'_{/\dZ^\Phi_p}$, $\bM_\fp(\rV^\circ,\rK^{p\circ})(S)$ is the set of equivalence classes of sextuples $(A_0,\lambda_0,\eta_0^p;A,\lambda,\eta^p)$, where
\begin{itemize}[label={\ding{109}}]
  \item $(A_0,\lambda_0,\eta_0^p)$ is an element in $\bT_\fp(S)$;

  \item $(A,\lambda)$ is a unitary $O_F$-abelian scheme of signature type $N\Phi-\tau_\infty+\tau_\infty^\tc$ over $S$ (Definitions \ref{de:unitary_abelian_scheme} and \ref{de:signature}) such that $\Ker\lambda[p^\infty]$ is contained in $A[\fp]$ of rank $p^2$;

  \item $\eta^p$ is a $\rK^{p\circ}$-level structure, that is, for a chosen geometric point $s$ on every connected component of $S$, a $\pi_1(S,s)$-invariant $\rK^{p\circ}$-orbit of isomorphisms
      \[
      \eta^p\colon\rV^\circ\otimes_\dQ\dA^{\infty,p}\to
      \Hom_{F\otimes_\dQ\dA^{\infty,p}}^{\lambda_0,\lambda}(\rH^\et_1(A_{0s},\dA^{\infty,p}),\rH^\et_1(A_s,\dA^{\infty,p}))
      \]
      of hermitian spaces over $F\otimes_\dQ\dA^{\infty,p}=F\otimes_{F^+}\dA_{F^+}^{\infty,p}$. See Construction \ref{cs:hermitian_structure} (with $\Box=\{\infty,p\}$) for the right-hand side.
\end{itemize}
The equivalence relation and the action of morphisms in $\fK(\rV^\circ)^p\times\fT$ are defined similarly as in Definition \ref{de:qs_moduli_scheme}.
\end{definition}

\begin{remark}
In the definition of the moduli functor $\bM_\fp(\rV^\circ,\obj)$, we use the \emph{definite} hermitian space $\rV^\circ$ to define the tame level structure -- this is different from the usual treatment. The reason for doing this is to make the uniformization map \eqref{eq:ns_uniformization1} for a certain stratum in the special fiber of $\bM_\fp(\rV^\circ,\obj)$ \emph{canonical}, since our main interest is the Shimura set $\Sh(\rV^\circ,\obj\rK^\circ_p)$, while the trade-off is that the relation between the generic fiber of $\bM_\fp(\rV^\circ,\obj)$ and unitary Shimura varieties cannot be made canonical (see Definition \ref{de:ns_uniformization_data}).
\end{remark}

We clearly have the forgetful morphism
\begin{align}\label{eq:ns_moduli_scheme}
\bM_\fp(\rV^\circ,\obj)\to\bT_\fp
\end{align}
in $\Fun(\fK(\rV^\circ)^p\times\fT,\sfP\Sch'_{/\dZ_p^\Phi})$, which is representable by quasi-projective schemes. According to Notation \ref{no:p_notation}, we shall denote by the base change of \eqref{eq:ns_moduli_scheme} to $\dF_p^\Phi$ by $\rM_\fp(\rV^\circ,\obj)\to\rT_\fp$, which is a morphism in $\Fun(\fK(\rV^\circ)^p\times\fT,\Sch_{/\dF_p^\Phi})$.

\begin{definition}\label{de:ns_moduli_scheme1}
For every $\rK^{p\circ}\in\fK(\rV^\circ)^p$, let $(\cA_0,\lambda_0,\eta_0^p;\cA,\lambda,\eta^p)$ be the universal object over $\bM_\fp(\rV^\circ,\rK^{p\circ})$. We define
\begin{enumerate}
  \item $\rM^\circ_\fp(\rV^\circ,\rK^{p\circ})$ to be the locus of $\rM_\fp(\rV^\circ,\rK^{p\circ})$ on which $\omega_{\cA^\vee,\tau_\infty}$ coincides with $\rH^\dr_1(\cA)_{\tau_\infty^\tc}^\perp$, which we call the \emph{balloon stratum};\footnote{This terminology is borrowed from an unpublished note by Kudla and Rapoport, where they study the corresponding Rapoport--Zink space. The intuition becomes clear after Theorem \ref{th:ns_basic_correspondence1} where we show that this stratum is a projective space fibration over a zero-dimensional scheme.}

  \item $\rM^\bullet_\fp(\rV^\circ,\rK^{p\circ})$ to be the locus of $\rM_\fp(\rV^\circ,\rK^{p\circ})$ on which $\rH^\dr_1(\cA)_{\tau_\infty}^\perp$ is a line subbundle of $\omega_{\cA^\vee,\tau_\infty^\tc}$, which we call the \emph{ground stratum};

  \item $\rM^\dag_\fp(\rV^\circ,\rK^{p\circ})$ to be $\rM^\circ_\fp(\rV^\circ,\rK^{p\circ})\bigcap\rM^\bullet_\fp(\rV^\circ,\rK^{p\circ})$, which we call the \emph{link stratum}.\footnote{This is the stratum linking balloons to the ground.}
\end{enumerate}
We denote
\begin{align*}
\rm^{\dag\circ}&\colon\rM^\dag_\fp(\rV^\circ,\obj)\to\rM^\circ_\fp(\rV^\circ,\obj),\\
\rm^{\dag\bullet}&\colon\rM^\dag_\fp(\rV^\circ,\obj)\to\rM^\bullet_\fp(\rV^\circ,\obj),
\end{align*}
the obvious inclusion morphisms.
\end{definition}

\begin{remark}
When $N=1$, the ground stratum and the link stratum are both empty.
\end{remark}

\begin{theorem}\label{th:ns_moduli_scheme}
For every $\rK^{p\circ}\in\fK(\rV^\circ)^p$, we have
\begin{enumerate}
  \item The scheme $\bM_\fp(\rV^\circ,\rK^{p\circ})$ is quasi-projective and strictly semistable over $\bT_\fp$ of relative dimension $N-1$; and we have
      \[
      \rM_\fp(\rV^\circ,\rK^{p\circ})=\rM^\circ_\fp(\rV^\circ,\rK^{p\circ})\bigcup\rM^\bullet_\fp(\rV^\circ,\rK^{p\circ}).
      \]
      Moreover, \eqref{eq:ns_moduli_scheme} is projective if and only if its base change to $\dQ_p^\Phi$ is.

  \item The loci $\rM^\circ_\fp(\rV^\circ,\rK^{p\circ})$ and $\rM^\bullet_\fp(\rV^\circ,\rK^{p\circ})$ are both closed subsets of $\rM_\fp(\rV^\circ,\rK^{p\circ})$, smooth over $\rT_\fp$ if we endow them with the induced reduced scheme structure.

  \item We have a canonical isomorphism
      \[
      \cT_{\rM^\circ_\fp(\rV^\circ,\rK^{p\circ})/\rT_\fp}\simeq\HOM\(\omega_{\cA^\vee,\tau_\infty^\tc},\Lie_{\cA,\tau_\infty^\tc}\)
      \]
      of coherent sheaves over $\rM^\circ_\fp(\rV^\circ,\rK^{p\circ})$ for the relative tangent sheaf.

  \item When $N\geq 2$, the relative tangent sheaf $\cT_{\rM^\bullet_\fp(\rV^\circ,\rK^{p\circ})/\rT_\fp}$ fits canonically into an exact sequence
      \[
      \resizebox{\hsize}{!}{
      \xymatrix{
      0 \ar[r]& \HOM\(\omega_{\cA^\vee,\tau_\infty},\omega_{\cA^\vee,\tau_\infty^\tc}^\perp/\omega_{\cA^\vee,\tau_\infty}\)
      \ar[r]& \cT_{\rM^\bullet_\fp(\rV^\circ,\rK^{p\circ})/\rT_\fp}
      \ar[r]& \HOM\(\omega_{\cA^\vee,\tau_\infty^\tc}/\rH^\dr_1(\cA)_{\tau_\infty}^\perp,\Lie_{\cA,\tau_\infty^\tc}\) \ar[r]& 0
      }
      }
      \]
      of coherent sheaves over $\rM^\bullet_\fp(\rV^\circ,\rK^{p\circ})$.

  \item When $N\geq 2$, the natural map $\cT_{\rM^\dag_\fp(\rV^\circ,\rK^{p\circ})/\rT_\fp}
      \to\cT_{\rM^\bullet_\fp(\rV^\circ,\rK^{p\circ})/\rT_\fp}\res_{\rM^\dag_\fp(\rV^\circ,\rK^{p\circ})}$ between relative tangent sheaves induces an isomorphism
      \[
      \cT_{\rM^\dag_\fp(\rV^\circ,\rK^{p\circ})/\rT_\fp}\simeq
      \HOM\(\omega_{\cA^\vee,\tau_\infty^\tc}/\rH^\dr_1(\cA)_{\tau_\infty}^\perp,\Lie_{\cA,\tau_\infty^\tc}\)
      \]
      of coherent sheaves over $\rM^\dag_\fp(\rV^\circ,\rK^{p\circ})$ under the exact sequence in (4). In particular, the exact sequence in (4) splits over $\rM^\dag_\fp(\rV^\circ,\rK^{p\circ})$.
\end{enumerate}
\end{theorem}

\begin{proof}
For (1), the (quasi-)projectiveness part is well-known. We consider the remaining assertions. Take a point $x=(A_0,\lambda_0,\eta_0^p;A,\lambda,\eta^p)\in\bM_\fp(\rV^\circ,\rK^{p\circ})(\kappa)$ for a perfect field $\kappa$ containing $\dF_p^\Phi$, and denote by $\cO_x$ the completed local ring of $\bM_\fp(\rV^\circ,\rK^{p\circ})$ at $x$. We have a $W(\kappa)$-bilinear pairing $\langle\;,\;\rangle_{\lambda,\tau_\infty}\colon\cD(A)_{\tau_\infty}\times\cD(A)_{\tau_\infty^\tc}\to W(\kappa)$ as in Notation \ref{no:weil_pairing_dieudonne}. By Proposition \ref{pr:deformation}, we have for every Artinian $W(\kappa)$-ring $R$ that is a quotient of $\cO_x$, that $\Hom_{W(\kappa)}(\cO_x,R)$ is the set of pairs of $R$-subbundles
\[
M_{\tau_\infty}\subseteq\cD(A)_{\tau_\infty}\otimes_{W(\kappa)}R,\quad
M_{\tau_\infty^\tc}\subseteq\cD(A)_{\tau_\infty^\tc}\otimes_{W(\kappa)}R
\]
of ranks $1$ and $N-1$ lifting $\omega_{A^\vee/\kappa,\tau_\infty}$ and $\omega_{A^\vee/\kappa,\tau_\infty^\tc}$, respectively, such that $\langle M_{\tau_\infty},M_{\tau_\infty^\tc}\rangle_{\lambda,\tau_\infty}=0$. We choose isomorphisms $\cD(A)_{\tau_\infty}\simeq W(\kappa)^{\oplus N}$ and $\cD(A)_{\tau_\infty^\tc}\simeq W(\kappa)^{\oplus N}$ under which the pairing $\langle\;,\;\rangle_{\lambda,\tau_\infty}$ is given by
\[
\langle(x_1,\dots,x_N),(y_1,\dots,y_N)\rangle_{\lambda,\tau_\infty}=px_1y_1+x_2y_2+\cdots +x_Ny_N.
\]
There are four possible cases.
\begin{enumerate}[label=(\roman*)]
  \item If $\omega_{A^\vee/\kappa,\tau_\infty}$ is generated by $(1,0,\dots,0)$ and $\omega_{A^\vee/\kappa,\tau_\infty^\tc}$ contains $(1,0,\dots,0)$, then possibly after changing coordinates, we may assume that $\omega_{A^\vee/\kappa,\tau_\infty^\tc}=\{(y_1,\dots,y_{N-1},0)\}$. Then we have $\cO_x\simeq W(\kappa)[[x_1,\dots,x_{N-1},x_N]]/(x_1x_N-p)$. In this case, $x$ must belong to $\rM_\fp^\dag(\rV^\circ,\rK^{p\circ})(\kappa)$.

  \item If $\omega_{A^\vee/\kappa,\tau_\infty}$ is generated by $(1,0,\dots,0)$ and $\omega_{A^\vee/\kappa,\tau_\infty^\tc}$ does not contain $(1,0,\dots,0)$, then possibly after changing coordinates, we may assume that $\omega_{A^\vee/\kappa,\tau_\infty^\tc}=\{(0,y_2,\dots,y_N)\}$. It is clear that $M_{\tau_\infty}$ is determined by $M_{\tau_\infty^\tc}$; and $\cO_x\simeq W(\kappa)[[x_2,\dots,x_N]]$.

  \item If $\omega_{A^\vee/\kappa,\tau_\infty}$ is not generated by $(1,0,\dots,0)$ and $\omega_{A^\vee/\kappa,\tau_\infty^\tc}$ contains $(1,0,\dots,0)$, then possibly after changing coordinates, we may assume that $\omega_{A^\vee/\kappa,\tau_\infty}$ is generated by $(0,\dots,0,1)$. It is clear that $M_{\tau_\infty^\tc}$ is determined by $M_{\tau_\infty}$; and $\cO_x\simeq W(\kappa)[[x_1,\dots,x_{N-1}]]$.

  \item If $\omega_{A^\vee/\kappa,\tau_\infty}$ is not generated by $(1,0,\dots,0)$ and $\omega_{A^\vee/\kappa,\tau_\infty^\tc}$ does not contain $(1,0,\dots,0)$, then this would not happen.
\end{enumerate}
Together with the fact that $\bM_\fp(\rV^\circ,\rK^{p\circ})\otimes\dQ$ is smooth of dimension $N-1$, $\bM_\fp(\rV^\circ,\rK^{p\circ})$ is strictly semistable over $\bT_\fp$ of relative dimension $N-1$. Moreover, $\rM^\circ_\fp(\rV^\circ,\rK^{p\circ})$ is the locus where (i) or (ii) happens; and $\rM^\bullet_\fp(\rV^\circ,\rK^{p\circ})$ is the locus where (i) or (iii) happens. Thus, both (1) and (2) follow.

For (3--5), we will use deformation theory. For common use, we consider a closed immersion $S\hookrightarrow\hat{S}$ in $\Sch'_{/\rT_\fp}$ defined by an ideal sheaf $\cI$ with $\cI^2=0$. Take an $S$-point $(A_0,\lambda_0,\eta_0^p;A,\lambda,\eta^p)$ in various schemes we will consider. By Proposition \ref{pr:deformation}, we need to lift $\omega_{A^\vee,\tau_\infty}$ and $\omega_{A^\vee,\tau_\infty^\tc}$ to subbundles $\hat\omega_{A^\vee,\tau_\infty}\subseteq\rH^\cris_1(A/\hat{S})_{\tau_\infty}$ and $\hat\omega_{A^\vee,\tau_\infty^\tc}\subseteq\rH^\cris_1(A/\hat{S})_{\tau_\infty^\tc}$, respectively, that are orthogonal to each other under the pairing \eqref{eq:weil_pairing_cris}.

For (3), since we require $\langle\hat\omega_{A^\vee,\tau_\infty},\rH^\cris_1(A/\hat{S})_{\tau_\infty^\tc}\rangle_{\lambda,\tau_\infty}^\cris=0$, it remains to lift $\hat\omega_{A^\vee,\tau_\infty^\tc}$ without restriction. Thus, (3) follows by Remark \ref{re:hodge_sequence}.

For (4), we need to first find lifting $\hat\omega_{A^\vee,\tau_\infty^\tc}$ that contains $\rH^\cris_1(A/\hat{S})_{\tau_\infty}^\perp$; and then find lifting $\hat\omega_{A^\vee,\tau_\infty}$ satisfying $\langle\hat\omega_{A^\vee,\tau_\infty},\hat\omega_{A^\vee,\tau_\infty^\tc}\rangle_{\lambda,\tau_\infty}^\cris=0$. Thus, (4) follows by Remark \ref{re:hodge_sequence}.

For (5), we only need to find lifting $\hat\omega_{A^\vee,\tau_\infty^\tc}$ that contains $\rH^\cris_1(A/\hat{S})_{\tau_\infty}^\perp$, which implies (5).
\end{proof}

In the remaining part of this subsection, we discuss the relation between $\bM_\fp(\rV^\circ,\obj)$ and certain unitary Shimura varieties. Since we use a standard definite hermitian space to parameterize the level structures, such relation is not canonical, which depends on the choice of an indefinite uniformization datum defined as follows.

\begin{definition}\label{de:ns_uniformization_data}
We define an \emph{indefinite uniformization datum for $\rV^\circ$ (at $\fp$)} to be a collection of $(\rV',\tj,\{\Lambda'_\fq\}_{\fq\mid p})$, where
\begin{itemize}[label={\ding{109}}]
  \item $\rV'$ is a standard indefinite hermitian space over $F$ of rank $N$;

  \item $\tj\colon\rV^\circ\otimes_\dQ\dA^{\infty,p}\to\rV'\otimes_\dQ\dA^{\infty,p}$ is an isometry;

  \item for every prime $\fq$ of $F^+$ above $p$ other than $\fp$, $\Lambda'_\fq$ is a self-dual $O_{F_\fq}$-lattice in $\rV'\otimes_FF_\fq$; and

  \item $\Lambda'_\fp$ is an $O_{F_\fp}$-lattice in $\rV'\otimes_FF_\fp$ satisfying $\Lambda'_\fp\subseteq(\Lambda'_\fp)^\vee$ and $(\Lambda'_\fp)^\vee/\Lambda'_\fp$ has length $1$.
\end{itemize}
\end{definition}

By the Hasse principle for hermitian spaces, there exists an indefinite uniformization datum for which we fix one. Let $\rK'_\fq$ be the stabilizer of $\Lambda'_\fq$ for every $\fq$ over $p$; and put $\rK'_p\coloneqq\prod_{\fq\mid p}\rK'_\fq$. The isometry $\tj$ induces an equivalence of categories $\tj\colon\fK(\rV^\circ)^p\xrightarrow{\sim}\fK(\rV')^p$.

Then similar to Remark \ref{re:qs_moduli_scheme_shimura}, we obtain a ``moduli interpretation'' isomorphism
\begin{align}\label{eq:ns_moduli_scheme_shimura}
\bM^\eta_\fp(\rV^\circ,\obj)\xrightarrow{\sim}\Sh(\rV',\tj\obj\rK'_p)\times_{\Spec{F}}\bT^\eta_\fp
\end{align}
(Notation \ref{no:p_notation}(5)) in $\Fun(\fK(\rV^\circ)^p\times\fT,\Sch_{/\dQ_p^\Phi})_{/\bT^\eta_\fp}$, where $\fT$ acts on $\Sh(\rV',\tj\obj\rK'_p)\times_{\Spec{F}}\bT^\eta_\fp$ via the second factor.

\begin{lem}\label{le:ns_pbc}
Let $L$ be a $p$-coprime coefficient ring. The two specialization maps
\begin{align*}
\rH^i_{\fT,c}(\bM_\fp(\rV^\circ,\obj)\otimes_{\dZ_p^\Phi}\ol\dQ_p,L)&\to\rH^i_{\fT,c}(\ol\rM_\fp(\rV^\circ,\obj),\rR\Psi L), \\
\rH^i_\fT(\bM_\fp(\rV^\circ,\obj)\otimes_{\dZ_p^\Phi}\ol\dQ_p,L)&\to\rH^i_\fT(\ol\rM_\fp(\rV^\circ,\obj),\rR\Psi L),
\end{align*}
are both isomorphisms. In particular, \eqref{eq:ns_moduli_scheme_shimura} induces isomorphisms
\begin{align*}
\rH^i_{\et,c}(\Sh(\rV',\tj\obj\rK'_p)_{\ol{F}},L)&\simeq\rH^i_{\fT,c}(\ol\rM_\fp(\rV^\circ,\obj),\rR\Psi L),\\
\rH^i_\et(\Sh(\rV',\tj\obj\rK'_p)_{\ol{F}},L)&\simeq\rH^i_\fT(\ol\rM_\fp(\rV^\circ,\obj),\rR\Psi L),
\end{align*}
in $\Fun(\fK(\rV^\circ)^p,\Mod(L[\Gal(\ol\dQ_p/\dQ^\Phi_p)]))$ for every $i\in\dZ$. Here, $\Gal(\ol\dQ_p/\dQ^\Phi_p)$ is regarded as a subgroup of $\Gal(\ol{F}/F)$ under our fixed isomorphism $\iota_p\colon\dC\simeq\ol\dQ_p$.
\end{lem}

\begin{proof}
When $\bM_\fp(\rV,\obj)$ is proper, this is simply the proper base change. When $\bM_\fp(\rV,\obj)$ is not proper, this follows from \cite{LS18}*{Corollary~5.20}.
\end{proof}

\begin{remark}\label{re:ns_smooth}
When $F^+\neq\dQ$, the Shimura variety $\Sh(\rV',\rK^{p\prime}\rK'_p)$ is proper over $F$ for $\rK^{p\prime}\in\fK(\rV')^p$. We explain that $\Sh(\rV',\rK^{p\prime}\rK'_p)$ has proper smooth reduction at every place $w$ of $F$ above $\Sigma^+_p\setminus\{\fp\}$.

Take a place $w$ of $F$ above $\Sigma^+_p\setminus\{\fp\}$. Choose a CM type $\Phi$ containing $\tau_\infty$ and an isomorphism $\dC\simeq\ol\dQ_p$ that induces $w$ (not the unique place above $\fp$!). Put $\bT_w\coloneqq\bT_p(\rW_0,\rK_0^p)\otimes_{O_{F_\Phi}\otimes\dZ_{(p)}}\dZ^\Phi_p$. We define a functor $\bM_w(\rV',\rK^{p\prime})$ on $\Sch'_{/\dZ^\Phi_p}$ such that for every $S\in\Sch'_{/\dZ^\Phi_p}$, $\bM_w(\rV',\rK^{p\prime})(S)$ is the set of equivalence classes of sextuples $(A_0,\lambda_0,\eta_0^p;A,\lambda,\eta^p)$, where
\begin{itemize}[label={\ding{109}}]
  \item $(A_0,\lambda_0,\eta_0^p)$ is an element in $\bT_w(S)$;

  \item $(A,\lambda)$ is a unitary $O_F$-abelian scheme of signature type $N\Phi-\tau_\infty+\tau_\infty^\tc$ over $S$ (Definitions \ref{de:unitary_abelian_scheme} and \ref{de:signature}) such that $\Ker\lambda[p^\infty]$ is contained in $A[\fp]$ of rank $p^2$;

  \item $\eta^p$ is a $\rK^{p\prime}$-level structure, similarly defined as in Definition \ref{de:ns_moduli_scheme}.
\end{itemize}
Then $\bM_w(\rV',\rK^{p\prime})$ is represented by a projective scheme over $\dZ_p^\Phi$. An easy computation of the tangent sheaf as in Theorem \ref{th:qs_moduli_scheme} shows that $\bM_w(\rV',\rK^{p\prime})$ is smooth of relative dimension $N-1$. Moreover, we have a canonical isomorphism
\[
\bM^\eta_w(\rV',\rK^{p\prime})\simeq\Sh(\rV',\rK^{p\prime}\rK'_p)\times_{\Spec{F}}\bT^\eta_w
\]
over $\bT^\eta_w$. Thus, $\Sh(\rV',\rK^{p\prime}\rK'_p)$ has proper smooth reduction at $w$ as $\bT_w$ is finite \'{e}tale over $O_{F_w}$.
\end{remark}

\subsection{Basic correspondence for the balloon stratum}

In this subsection, we construct and study the basic correspondence for the balloon stratum $\rM^\circ_\fp(\rV^\circ,\obj)$.

\begin{definition}\label{de:qs_definite1}
We define a functor
\begin{align*}
\rS^\circ_\fp(\rV^\circ,\obj)\colon\fK(\rV^\circ)^p\times\fT &\to\sfP\Sch'_{/\dF^\Phi_p} \\
\rK^{p\circ} &\mapsto \rS^\circ_\fp(\rV^\circ,\rK^{p\circ})
\end{align*}
such that for every $S\in\Sch'_{/\dF^\Phi_p}$, $\rS^\circ_\fp(\rV^\circ,\rK^{p\circ})(S)$ is the set of equivalence classes of sextuples $(A_0,\lambda_0,\eta_0^p;A^\circ,\lambda^\circ,\eta^{p\circ})$, where
\begin{itemize}[label={\ding{109}}]
  \item $(A_0,\lambda_0,\eta_0^p)$ is an element in $\rT_\fp(S)$;

  \item $(A^\circ,\lambda^\circ)$ is a unitary $O_F$-abelian scheme of signature type $N\Phi$ over $S$ such that $\lambda^\circ$ is $p$-principal;

  \item $\eta^{p\circ}$ is, for a chosen geometric point $s$ on every connected component of $S$, a $\pi_1(S,s)$-invariant $\rK^{p\circ}$-orbit of isomorphisms
      \[
      \eta^{p\circ}\colon\rV^\circ\otimes_\dQ\dA^{\infty,p}\to
      \Hom_{F\otimes_\dQ\dA^{\infty,p}}^{\lambda_0,\lambda^\circ}(\rH^\et_1(A_{0s},\dA^{\infty,p}),\rH^\et_1(A^\circ_s,\dA^{\infty,p}))
      \]
      of hermitian spaces over $F\otimes_\dQ\dA^{\infty,p}=F\otimes_{F^+}\dA_{F^+}^{\infty,p}$.
\end{itemize}
The equivalence relation and the action of morphisms in $\fK(\rV^\circ)^p\times\fT$ are defined similarly as in Definition \ref{de:qs_moduli_scheme}.
\end{definition}

We clearly have the forgetful morphism
\begin{align*}
\rS^\circ_\fp(\rV^\circ,\obj)\to\rT_\fp
\end{align*}
in $\Fun(\fK(\rV^\circ)^p\times\fT,\sfP\Sch'_{/\dF_p^\Phi})$, which is represented by finite and \'{e}tale schemes by \cite{RSZ}*{Theorem~4.4}.

Now we take a point $s^\circ=(A_0,\lambda_0,\eta_0^p;A^\circ,\lambda^\circ,\eta^{p\circ})\in\rS^\circ_\fp(\rV^\circ,\rK^{p\circ})(\kappa)$ where $\kappa$ is a perfect field containing $\dF^\Phi_p$. Then $A^\circ_{\ol\kappa}[\fp^\infty]$ is a supersingular $p$-divisible by the signature condition and the fact that $\fp$ is inert in $F$. The $(\kappa,\sigma^{-1})$-linear Verschiebung map
\[
\tV\colon\rH^\dr_1(A^\circ/\kappa)_{\tau_\infty}\to\rH^\dr_1(A^\circ/\kappa)_{\sigma^{-1}\tau_\infty}
=\rH^\dr_1(A^\circ/\kappa)_{\tau_\infty^\tc}
\]
(Notation \ref{re:frobenius_verschiebung}) is an isomorphism. Thus, we obtain a $(\kappa,\sigma)$-linear isomorphism
\[
\tV^{-1}\colon\rH^\dr_1(A^\circ/\kappa)_{\tau_\infty^\tc}\to\rH^\dr_1(A^\circ/\kappa)_{\tau_\infty}.
\]
We define a non-degenerate pairing
\[
\{\;,\;\}_{s^\circ}\colon\rH^\dr_1(A^\circ/\kappa)_{\tau_\infty^\tc}\times\rH^\dr_1(A^\circ/\kappa)_{\tau_\infty^\tc}\to\kappa
\]
by the formula $\{x,y\}_{s^\circ}\coloneqq\langle\tV^{-1} x,y\rangle_{\lambda^\circ,\tau_\infty}$ (Notation \ref{no:weil_pairing}). To ease notation, we put
\[
\sV_{s^\circ}\coloneqq\rH^\dr_1(A^\circ/\kappa)_{\tau_\infty^\tc}.
\]
By the same proof of Lemma \ref{le:qs_dl}, we know that $(\sV_{s^\circ},\{\;,\;\}_{s^\circ})$ is admissible. Thus, we have the Deligne--Lusztig variety $\DL_{s^\circ}\coloneqq\DL(\sV_{s^\circ},\{\;,\;\}_{s^\circ},N-1)$ (Definition \ref{de:dl}).

\begin{definition}\label{de:ns_basic_correspondence1}
We define a functor
\begin{align*}
\rB^\circ_\fp(\rV^\circ,\obj)\colon\fK(\rV^\circ)^p\times\fT &\to\sfP\Sch'_{/\dF^\Phi_p} \\
\rK^{p\circ} &\mapsto \rB^\circ_\fp(\rV^\circ,\rK^{p\circ})
\end{align*}
such that for every $S\in\Sch'_{/\dF^\Phi_p}$, $\rB^\circ_\fp(\rV^\circ,\rK^{p\circ})(S)$ is the set of equivalence classes of decuples $(A_0,\lambda_0,\eta_0^p;A,\lambda,\eta^p;A^\circ,\lambda^\circ,\eta^{p\circ};\beta)$, where
\begin{itemize}[label={\ding{109}}]
  \item $(A_0,\lambda_0,\eta_0^p;A,\lambda,\eta^p)$ is an element of $\rM^\circ_\fp(\rV^\circ,\rK^{p\circ})(S)$;

  \item $(A_0,\lambda_0,\eta_0^p;A^\circ,\lambda^\circ,\eta^{p\circ})$ is an element of $\rS^\circ_\fp(\rV^\circ,\rK^{p\circ})(S)$; and

  \item $\beta\colon A\to A^\circ$ is an $O_F$-linear quasi-$p$-isogeny (Definition \ref{de:p_quasi}) such that
      \begin{enumerate}[label=(\alph*)]
        \item $\Ker\beta[p^\infty]$ is contained in $A[\fp]$;

        \item we have $\lambda=\beta^\vee\circ\lambda^\circ\circ\beta$; and

        \item the $\rK^{p\circ}$-orbit of maps $v\mapsto\beta_*\circ\eta^p(v)$ for $v\in\rV^\circ\otimes_\dQ\dA^{\infty,p}$ coincides with $\eta^{p\circ}$.
      \end{enumerate}
\end{itemize}
The equivalence relation and the action of morphisms in $\fK(\rV^\circ)^p\times\fT$ are defined similarly as in Definition \ref{de:qs_basic_correspondence}.
\end{definition}

We obtain in the obvious way a correspondence
\begin{align}\label{eq:ns_basic_correspondence1}
\xymatrix{
\rS^\circ_\fp(\rV^\circ,\obj)  &
\rB^\circ_\fp(\rV^\circ,\obj) \ar[r]^-{\iota^\circ}\ar[l]_-{\pi^\circ} &
\rM^\circ_\fp(\rV^\circ,\obj)
}
\end{align}
in $\Fun(\fK(\rV^\circ)^p\times\fT,\sfP\Sch'_{/\dF_p^\Phi})_{/\rT_\fp}$.

\begin{definition}[Basic correspondence]
We refer to \eqref{eq:ns_basic_correspondence1} as the \emph{basic correspondence} on the balloon stratum $\rM^\circ_\fp(\rV^\circ,\obj)$, with $\rS^\circ_\fp(\rV^\circ,\obj)$ being the \emph{source} of the basic correspondence.
\end{definition}

\begin{theorem}\label{th:ns_basic_correspondence1}
In the diagram \eqref{eq:ns_basic_correspondence1}, $\iota^\circ$ is an isomorphism. Moreover, for every point
$s^\circ=(A_0,\lambda_0,\eta_0^p;A^\circ,\lambda^\circ,\eta^{p\circ})\in\rS^\circ_\fp(\rV^\circ,\rK^{p\circ})(\kappa)$ where $\kappa$ is a perfect field containing $\dF^\Phi_p$, if we put $\rB^\circ_{s^\circ}\coloneqq\pi^{\circ-1}(s^\circ)$, then the assignment sending
$(A_0,\lambda_0,\eta_0^p;A,\lambda,\eta^p;A^\circ,\lambda^\circ,\eta^{p\circ};\beta)\in\rB^\circ_{s^\circ}(S)$ to the subbundle
\[
H\coloneqq \beta_{*,\tau_\infty^\tc}\omega_{A^\vee/S,\tau_\infty^\tc}\subseteq\rH^\dr_1(A^\circ/S)_{\tau_\infty^\tc}
=\rH^\dr_1(A^\circ/\kappa)_{\tau_\infty^\tc}\otimes_\kappa\cO_S=(\sV_{s^\circ})_S
\]
induces an isomorphism $\zeta^\circ_{s^\circ}\colon\rB^\circ_{s^\circ}\xrightarrow{\sim}\dP(\sV_{s^\circ})$ satisfying that
\begin{enumerate}
  \item $\zeta^\circ_{s^\circ}$ restricts to an isomorphism
     \[
     \zeta^\circ_{s^\circ}\colon\rB^\circ_{s^\circ}\bigcap\iota^{\circ-1}\rM^\dag_\fp(\rV^\circ,\rK^{p\circ})
     \xrightarrow{\sim}\DL_{s^\circ}=\DL(\sV_{s^\circ},\{\;,\;\}_{s^\circ},N-1);
     \]

  \item we have an isomorphism
     \[
     \HOM\(\omega_{\cA^\vee,\tau_\infty},\omega_{\cA^\vee,\tau_\infty^\tc}^\perp/\omega_{\cA^\vee,\tau_\infty}\)
     \simeq(\zeta^\circ_{s^\circ})^*\cO_{\dP(\sV_{s^\circ})}(-(p+1)).
     \]
\end{enumerate}
In particular, $\rB^\circ_{s^\circ}\bigcap\iota^{\circ-1}\rM^\dag_\fp(\rV^\circ,\rK^{p\circ})$ is a Fermat hypersurface in $\rB^\circ_{s^\circ}\simeq\dP(\sV_{s^\circ})$.
\end{theorem}

\begin{proof}
Take an object $\rK^{p\circ}\in\fK(\rV^\circ)^p$. It is clear that $\rB^\circ_\fp(\rV^\circ,\obj)$ is a scheme. We denote by $(\cA_0,\lambda_0,\eta_0^p;\cA,\lambda,\eta^p;\cA^\circ,\lambda^\circ,\eta^{p\circ};\beta)$ the universal object over $\rB^\circ_\fp(\rV^\circ,\rK^{p\circ})$.

First, we show that $\iota^\circ$ is an isomorphism. It is an easy exercise from Grothendieck--Messing theory that the canonical map $\cT_{\rB^\circ_\fp(\rV^\circ,\rK^{p\circ})/\rT_\fp}\to\iota^{\circ*}\cT_{\rM^\circ_\fp(\rV^\circ,\rK^{p\circ})/\rT_\fp}$ is an isomorphism. Thus, it suffices to show that $\iota^\circ(\kappa')$ is a bijection for every algebraically closed field $\kappa'$ containing $\kappa$. To ease notation, we may assume $\kappa'=\kappa$. We construct an inverse of $\iota^\circ(\kappa)$. Take a point $(A_0,\lambda_0,\eta_0^p;A,\lambda,\eta^p)\in\rM^\circ_\fp(\rV^\circ,\rK^{p\circ})(\kappa)$. Write $\tilde\omega_{A^\vee,\tau_\infty}$ the preimage of $\omega_{A^\vee,\tau_\infty}$ under the reduction map $\cD(A)_{\tau_\infty}\to\rH^\dr_1(A/\kappa)_{\tau_\infty}$. As $\langle\omega_{A^\vee,\tau_\infty},\rH^\dr_1(A/\kappa)_{\tau_\infty^\tc}\rangle_{\lambda,\tau_\infty}=0$, we have $\cD(A)_{\tau_\infty^\tc}^\vee=p^{-1}\tilde\omega_{A^\vee,\tau_\infty}$. Now we put $\cD_{A^\circ,\tau}\coloneqq\cD(A)_{\tau}$ for $\tau\neq\tau_\infty$, and $\cD_{A^\circ,\tau_\infty}\coloneqq p^{-1}\tilde\omega_{A^\vee,\tau_\infty}$. We claim that $\cD_{A^\circ}\coloneqq\bigoplus_{\tau\in\Sigma_\infty}\cD_{A^\circ,\tau}$ is a Dieudonn\'{e} module, which amounts to the inclusions $\tF\cD_{A^\circ,\tau_\infty}\subseteq\cD_{A^\circ,\tau_\infty^\tc}$ and $\tV\cD_{A^\circ,\tau_\infty}\subseteq\cD_{A^\circ,\tau_\infty^\tc}$. The first one is obvious; and the second one is equivalent to the first one as $\cD_{A^\circ,\tau_\infty}$ and $\cD_{A^\circ,\tau_\infty^\tc}$ are integrally dual under $\langle\;,\;\rangle_{\lambda,\tau_\infty}^\cris$. Then by the Dieudonn\'{e} theory, there is an $O_F$-abelian scheme $A^\circ$ over $\kappa$ with $\cD(A^\circ)_\tau=\cD_{A^\circ,\tau}$ for every $\tau\in\Sigma_\infty$, and an $O_F$-linear isogeny $\beta\colon A\to A^\circ$ inducing the inclusion of Dieudonn\'{e} modules $\cD(A)\subseteq\cD(A^\circ)$. By Lemma \ref{le:inverse_isogeny}(2,4), the $O_F$-abelian scheme $A^\circ$ has signature type $N\Phi$. Let $\lambda^\circ$ be the unique quasi-polarization of $A^\circ$ satisfying $\lambda=\beta^\vee\circ\lambda^\circ\circ\beta$, which is $p$-principal as $\cD_{A^\circ,\tau_\infty^\tc}=\cD_{A^\circ,\tau_\infty}^\vee$. Finally, we let $\eta^{p\circ}$ be the map sending $v\in\rV^\circ\otimes_\dQ\dA^{\infty,p}$ to $\beta_*\circ\eta^p(v)$. Thus, we obtain an object $(A_0,\lambda_0,\eta_0^p;A,\lambda,\eta^p;A^\circ,\lambda^\circ,\eta^{p\circ};\beta)\in\rS^\circ_\fp(\rV^\circ,\rK^{p\circ})(\kappa)$. It is straightforward to check that such assignment gives rise to an inverse of $\iota^\circ(\kappa)$.

Second, we show that $\zeta^\circ_{s^\circ}$ is well-defined, namely, $H$ is a subbundle of rank $N-1$. By Lemma \ref{le:inverse_isogeny}(2,4) and Definition \ref{de:ns_basic_correspondence1}(b), we have $\rank_{\cO_S}(\Ker\beta_{*,\tau_\infty})-\rank_{\cO_S}(\Ker\beta_{*,\tau_\infty^\tc})=1$ and $\rank_{\cO_S}(\Ker\beta_{*,\tau_\infty})+\rank_{\cO_S}(\Ker\beta_{*,\tau_\infty^\tc})=1$. Thus, $\beta_{*,\tau_\infty^\tc}$ is an isomorphism, hence $H$ is a subbundle of rank $N-1$.

Third, we show that $\zeta^\circ_{s^\circ}$ is an isomorphism. Denote by $\cH\subseteq(\sV_{s^\circ})_{\dP(\sV_{s^\circ})}$ the universal subbundle (of rank $N-1$). Then we have a canonical isomorphism
\[
\cT_{\dP(\sV_{s^\circ})/\kappa}\simeq\Hom_{\cO_{\dP(\sV_{s^\circ})}}\(\cH,\rH^\dr_1(A^\circ/\kappa)_{\tau_\infty^\tc}/\cH\).
\]
By Theorem \ref{th:ns_moduli_scheme}(1) and the fact that $\beta_{*,\tau_\infty^\tc}$ is an isomorphism, we obtain an isomorphism
\[
\(\iota^{\circ*}\cT_{\rM^\circ_\fp(\rV^\circ,\rK^{p\circ})/\rT_\fp}\)\res_{\rB^\circ_{s^\circ}}
\xrightarrow{\sim}\zeta^{\circ*}_{s^\circ}\cT_{\dP(\sV_{s^\circ})/\kappa}.
\]
Thus, to show that $\zeta^\circ_{s^\circ}\colon\rB^\circ_{s^\circ}\to\dP(\sV_{s^\circ})$ is an isomorphism, it suffices to construct an inverse of $\zeta^\circ_{s^\circ}(\kappa')$ for every algebraically closed field $\kappa'$ containing $\kappa$. To ease notation, we may assume $\kappa'=\kappa$. Take a $\kappa$-linear subspace $H\subseteq\sV_{s^\circ}=\rH^\dr_1(A^\circ)_{\tau_\infty^\tc}$ of rank $N-1$. Let $\tilde{H}$ denote by its preimage under the reduction map $\cD(A^\circ)_{\tau_\infty^\tc}\to\rH^\dr_1(A^\circ)_{\tau_\infty^\tc}$. We put $\cD_{A,\tau}\coloneqq\cD(A^\circ)_\tau$ for $\tau\neq\tau_\infty$, and $\cD_{A,\tau_\infty}\coloneqq\tV^{-1}\tilde{H}\subseteq\cD(A^\circ)_{\tau_\infty}$. It is clear that $\cD_A\coloneqq\bigoplus_{\tau\in\Sigma_\infty}\cD_{A,\tau}$ is a Dieudonn\'{e} module. By the Dieudonn\'{e} theory, there is an $O_F$-abelian scheme $A$ over $\kappa$ with $\cD(A)_\tau=\cD_{A,\tau}$ for every $\tau\in\Sigma_\infty$, and an $O_F$-linear isogeny $\beta\colon A\to A^\circ$ inducing the inclusion of Dieudonn\'{e} modules $\cD(A)\subseteq\cD(A^\circ)$. By a similar argument as for $\iota^\circ$, we obtain a point $(A,\lambda,\eta^p;\beta)\in\rB^\circ_{s^\circ}(\kappa)$; and it follows that such assignment is an inverse of $\zeta^\circ_{s^\circ}(\kappa)$.

Finally, we check the two properties of $\zeta^\circ_{s^\circ}$.

For (1), we check that the closed subscheme $\zeta^\circ_{s^\circ}(\rB^\circ_{s^\circ}\cap\iota^{\circ-1}\rM^\dag_\fp(\rV^\circ,\rK^{p\circ}))$ coincides with $\DL(\sV_{s^\circ},\{\;,\;\}_{s^\circ},N-1)$. Recall that $\rM^\dag_\fp(\rV^\circ,\rK^{p\circ})$ is define by the condition
\[
\rH^1_\dr(A/S)_{\tau_\infty}^\perp\subseteq\omega_{A^\vee/S,\tau_\infty^\tc}.
\]
Note that we have $H=\beta_{*,\tau_\infty^\tc}\omega_{A^\vee/S,\tau_\infty^\tc}$ and $\tV^{-1}H^{(p)}=\beta_{*,\tau_\infty}\rH^\dr_1(A/S)_{\tau_\infty}$, which implies $(\tV^{-1}H^{(p)})^\perp=(\beta_{*,\tau_\infty}\rH^\dr_1(A/S)_{\tau_\infty})^\perp
=\beta_{*,\tau_\infty^\tc}(\rH^1_\dr(A/S)_{\tau_\infty}^\perp)$. Applying the isomorphism $\beta_{*,\tau_\infty^\tc}$, the above condition is equivalent to
\[
(\tV^{-1}H^{(p)})^\perp\subseteq H,
\]
which is the condition defining $\DL(\sV_{s^\circ},\{\;,\;\}_{s^\circ},N-1)$.

For (2), we have
\[
\omega_{A^\vee,\tau_\infty}=\Ker\beta_{*,\tau_\infty}\simeq\rH^\dr_1(A^\circ/S)_{\tau_\infty}/\beta_{*,\tau_\infty}\rH^\dr_1(A/S)_{\tau_\infty}
=\rH^\dr_1(A^\circ/S)_{\tau_\infty}/\tV^{-1}H^{(p)}
\]
and
\[
\omega_{A^\vee,\tau_\infty^\tc}^\perp/\omega_{A^\vee,\tau_\infty}\simeq\beta_{*,\tau_\infty}\omega_{A^\vee,\tau_\infty^\tc}^\perp
=(\beta_{*,\tau_\infty^\tc}\omega_{A^\vee/S,\tau_\infty^\tc})^\perp= H^\perp.
\]
Thus, we have
\[
\omega_{\cA^\vee,\tau_\infty}\simeq\zeta^{\circ*}_{s^\circ}\cO_{\dP(\sV_{s^\circ})}(p),\quad
\omega_{\cA^\vee,\tau_\infty^\tc}^\perp/\omega_{\cA^\vee,\tau_\infty}\simeq\zeta^{\circ*}_{s^\circ}\cO_{\dP(\sV_{s^\circ})}(-1)
\]
from which (2) follows.

The theorem is all proved.
\end{proof}

\begin{corollary}\label{co:normal}
When $N\geq 2$, the normal bundle of the closed immersion
\[
\rm^{\dag\bullet}\colon\rM^\dag_\fp(\rV^\circ,\rK^{p\circ})\to\rM^\bullet_\fp(\rV^\circ,\rK^{p\circ})
\]
is isomorphic to $(\rm^{\dag\circ})^*\cO_{\rM^\circ_\fp(\rV^\circ,\rK^{p\circ})}(-(p+1))$.
\end{corollary}

\begin{proof}
By Theorem \ref{th:ns_moduli_scheme}(4,5), we have that the normal bundle is isomorphic to
\[
\HOM\(\omega_{\cA^\vee,\tau_\infty},\omega_{\cA^\vee,\tau_\infty^\tc}^\perp/\omega_{\cA^\vee,\tau_\infty}\).
\]
Thus, the claim follows from Theorem \ref{th:ns_basic_correspondence1}. We can also argue that the normal bundle of $\rm^{\dag\bullet}$ is dual to the normal bundle of $\rm^{\dag\circ}$ which is isomorphic to $(\rm^{\dag\circ})^*\cO_{\rM^\circ_\fp(\rV^\circ,\rK^{p\circ})}(p+1)$ by Theorem \ref{th:ns_basic_correspondence1}.
\end{proof}

\begin{construction}\label{cs:ns_uniformization1}
Let $\rK^\circ_\fq$ be the stabilizer of $\Lambda^\circ_\fq$ for every $\fq\mid p$; and put $\rK^\circ_p\coloneqq\prod_{\fq\mid p}\rK^\circ_\fq$. Similar to Construction \ref{cs:qs_uniformization}, we may construct a \emph{uniformization map}, \emph{canonical} this time,
\begin{align}\label{eq:ns_uniformization1}
\upsilon^\circ\colon\rS^\circ_\fp(\rV^\circ,\obj)(\ol\dF_p)
\xrightarrow{\sim}\Sh(\rV^\circ,\obj\rK^\circ_p)\times\rT_\fp(\ol\dF_p)
\end{align}
in $\Fun(\fK(\rV^\circ)^p\times\fT,\Set)_{/\rT_\fp(\ol\dF_p)}$ which is an isomorphism, under which the induced action of $\Gal(\ol\dF_p/\dF_p^\Phi)$ on the target is trivial on $\Sh(\rV^\circ,\obj\rK^\circ_p)$.

Moreover, similar to Construction \ref{cs:qs_definite_hecke} and Proposition \ref{pr:qs_definite_hecke}, for every $g\in\rK^\circ_\fp\backslash\rU(\rV^\circ)(F^+_\fp)/\rK^\circ_\fp$, we may construct the Hecke correspondence
\[
\Hk_g\colon\rS^\circ_\fp(\rV^\circ,\obj)_g\to\rS^\circ_\fp(\rV^\circ,\obj)\times\rS^\circ_\fp(\rV^\circ,\obj)
\]
as a morphism in $\Fun(\fK(\rV^\circ)^p\times\fT,\Sch_{/\dF_p^\Phi})_{/\rT_\fp}$ that is finite \'{e}tale and compatible with the uniformization map.
\end{construction}

\subsection{Basic correspondence for the ground stratum}
\label{ss:ns_ground}

In this subsection, we construct and study the basic correspondence for the ground stratum $\rM^\bullet_\fp(\rV^\circ,\obj)$. We assume $N\geq 2$.

\begin{definition}\label{de:ns_definite2}
We define a functor
\begin{align*}
\rS^\bullet_\fp(\rV^\circ,\obj)\colon\fK(\rV^\circ)^p\times\fT &\to\sfP\Sch'_{/\dF^\Phi_p} \\
\rK^{p\circ} &\mapsto \rS^\bullet_\fp(\rV^\circ,\rK^{p\circ})
\end{align*}
such that for every $S\in\Sch'_{/\dF^\Phi_p}$, $\rS^\bullet_\fp(\rV^\circ,\rK^{p\circ})(S)$ is the set of equivalence classes of sextuples $(A_0,\lambda_0,\eta_0^p;A^\bullet,\lambda^\bullet,\eta^{p\bullet})$, where
\begin{itemize}[label={\ding{109}}]
  \item $(A_0,\lambda_0,\eta_0^p)$ is an element in $\rT_\fp(S)$;

  \item $(A^\bullet,\lambda^\bullet)$ is a unitary $O_F$-abelian scheme of signature type $N\Phi$ over $S$ such that $\Ker\lambda^\bullet[p^\infty]$ is trivial (resp.\ contained in $A^\bullet[\fp]$ of rank $p^2$) if $N$ is even (resp.\ odd);

  \item $\eta^{p\bullet}$ is, for a chosen geometric point $s$ on every connected component of $S$, a $\pi_1(S,s)$-invariant $\rK^{p\bullet}$-orbit of isomorphisms
      \[
      \eta^{p\bullet}\colon\rV^\circ\otimes_\dQ\dA^{\infty,p}\to
      \Hom_{F\otimes_\dQ\dA^{\infty,p}}^{\varpi\lambda_0,\lambda^\bullet}
      (\rH^\et_1(A_{0s},\dA^{\infty,p}),\rH^\et_1(A^\bullet_s,\dA^{\infty,p}))
      \]
      of hermitian spaces over $F\otimes_\dQ\dA^{\infty,p}=F\otimes_{F^+}\dA_{F^+}^{\infty,p}$.\footnote{Note that here we are using $\varpi\lambda_0$ rather than $\lambda_0$ in order to be consistent with the compatibility condition for polarizations in the isogeny considered in Definition \ref{de:ns_basic_correspondence2}.}
\end{itemize}
The equivalence relation and the action of morphisms in $\fK(\rV^\circ)^p\times\fT$ are defined similarly as in Definition \ref{de:qs_moduli_scheme}.
\end{definition}

We clearly have the forgetful morphism
\begin{align*}
\rS^\bullet_\fp(\rV^\circ,\obj)\to\rT_\fp
\end{align*}
in $\Fun(\fK(\rV^\circ)^p\times\fT,\sfP\Sch'_{/\dF_p^\Phi})$, which is represented by finite and \'{e}tale schemes by \cite{RSZ}*{Theorem~4.4}.\footnote{In fact, \cite{RSZ}*{Theorem~4.4} only considers the case where the polarization is $p$-principal (namely, $\Ker\lambda^\bullet[p^\infty]$ is trivial), but its proof works in the case where $\Ker\lambda^\bullet[p^\infty]$ is contained in $A^\bullet[\fp]$ of rank $p^2$ as well since the computation of the tangent space is the same.}

Now we take a point $s^\bullet=(A_0,\lambda_0,\eta_0^p;A^\bullet,\lambda^\bullet,\eta^{p\bullet})\in\rS^\bullet_\fp(\rV^\circ,\rK^{p\circ})(\kappa)$ where $\kappa$ is a perfect field containing $\dF^\Phi_p$. Then $A^\bullet_{\ol\kappa}[\fp^\infty]$ is a supersingular $p$-divisible by the signature condition and the fact that $\fp$ is inert in $F$. The $(\kappa,\sigma^{-1})$-linear Verschiebung map
\[
\tV\colon\rH^\dr_1(A^\bullet/\kappa)_{\tau_\infty}\to\rH^\dr_1(A^\bullet/\kappa)_{\sigma^{-1}\tau_\infty}
=\rH^\dr_1(A^\bullet/\kappa)_{\tau_\infty^\tc}
\]
(Notation \ref{re:frobenius_verschiebung}) is an isomorphism. Thus, we obtain a $(\kappa,\sigma)$-linear isomorphism
\[
\tV^{-1}\colon\rH^\dr_1(A^\bullet/\kappa)_{\tau_\infty^\tc}\to\rH^\dr_1(A^\bullet/\kappa)_{\tau_\infty}.
\]
We define a pairing
\[
\{\;,\;\}_{s^\bullet}\colon\rH^\dr_1(A^\bullet/\kappa)_{\tau_\infty^\tc}\times\rH^\dr_1(A^\bullet/\kappa)_{\tau_\infty^\tc}\to\kappa
\]
by the formula $\{x,y\}_{s^\bullet}\coloneqq\langle\tV^{-1} x,y\rangle_{\lambda^\bullet,\tau_\infty}$ (Notation \ref{no:weil_pairing}). To ease notation, we put
\[
\sV_{s^\bullet}\coloneqq\rH^\dr_1(A^\bullet/\kappa)_{\tau_\infty^\tc}.
\]
By the same proof of Lemma \ref{le:qs_dl}, we know that $(\sV_{s^\bullet},\{\;,\;\}_{s^\bullet})$ is admissible. Thus, we have the Deligne--Lusztig variety $\DL^\bullet_{s^\bullet}\coloneqq\DL^\bullet(\sV_{s^\bullet},\{\;,\;\}_{s^\bullet})$ (Definition \ref{de:dl_bullet}). Moreover, $\dim_\kappa\sV_{s^\bullet}^\perp$ is equal to $0$ (resp.\ $1$) when $N$ is even (resp.\ odd).

\begin{definition}\label{de:ns_basic_correspondence2}
We define a functor
\begin{align*}
\rB^\bullet_\fp(\rV^\circ,\obj)\colon\fK(\rV^\circ)^p\times\fT &\to\sfP\Sch'_{/\dF^\Phi_p} \\
\rK^{p\circ} &\mapsto \rB^\bullet_\fp(\rV^\circ,\rK^{p\circ})
\end{align*}
such that for every $S\in\Sch'_{/\dF^\Phi_p}$, $\rB^\bullet_\fp(\rV^\circ,\rK^{p\circ})(S)$ is the set of equivalence classes of decuples $(A_0,\lambda_0,\eta_0^p;A,\lambda,\eta^p;A^\bullet,\lambda^\bullet,\eta^{p\bullet};\gamma)$, where
\begin{itemize}[label={\ding{109}}]
  \item $(A_0,\lambda_0,\eta_0^p;A,\lambda,\eta^p)$ is an element of $\rM^\bullet_\fp(\rV^\circ,\rK^{p\circ})(S)$;

  \item $(A_0,\lambda_0,\eta_0^p;A^\bullet,\lambda^\bullet,\eta^{p\bullet})$ is an element of $\rS^\bullet_\fp(\rV^\circ,\rK^{p\circ})(S)$; and

  \item $\gamma\colon A\to A^\bullet$ is an $O_F$-linear quasi-$p$-isogeny (Definition \ref{de:p_quasi}) such that
    \begin{enumerate}[label=(\alph*)]
      \item $\Ker\gamma[p^\infty]$ is contained in $A[\fp]$;

      \item $(\Ker\gamma_{*,\tau_\infty})^\perp$ is contained in $\omega_{A^\vee/S,\tau_\infty^\tc}$;

      \item $\Ker\gamma_{*,\tau_\infty}$ contains $\rH^\dr_1(A/S)_{\tau_\infty^\tc}^\perp$;\footnote{This condition is implied by the others when $N$ is even.}

      \item we have $\varpi\cdot\lambda=\gamma^\vee\circ\lambda^\bullet\circ\gamma$; and

      \item the $\rK^{p\circ}$-orbit of maps $v\mapsto\gamma_*\circ\eta^p(v)$ for $v\in\rV^\circ\otimes_\dQ\dA^{\infty,p}$ coincides with $\eta^{p\bullet}$.
    \end{enumerate}
\end{itemize}
The equivalence relation and the action of morphisms in $\fK(\rV^\circ)^p\times\fT$ are defined similarly as in Definition \ref{de:qs_basic_correspondence}.
\end{definition}

We obtain in the obvious way a correspondence
\begin{align}\label{eq:ns_basic_correspondence2}
\xymatrix{
\rS^\bullet_\fp(\rV^\circ,\obj)  &
\rB^\bullet_\fp(\rV^\circ,\obj) \ar[r]^-{\iota^\bullet}\ar[l]_-{\pi^\bullet} &
\rM^\bullet_\fp(\rV^\circ,\obj)
}
\end{align}
in $\Fun(\fK(\rV^\circ)^p\times\fT,\sfP\Sch'_{/\dF_p^\Phi})_{/\rT_\fp}$.

\begin{definition}[Basic correspondence]
We refer to \eqref{eq:ns_basic_correspondence2} as the \emph{basic correspondence} on the ground stratum $\rM^\bullet_\fp(\rV^\circ,\obj)$, with $\rS^\bullet_\fp(\rV^\circ,\obj)$ being the \emph{source} of the basic correspondence.
\end{definition}

\begin{theorem}\label{th:ns_basic_correspondence2}
In the diagram \eqref{eq:ns_basic_correspondence2}, take a point
\[
s^\bullet=(A_0,\lambda_0,\eta_0^p;A^\bullet,\lambda^\bullet,\eta^{p\bullet})\in\rS^\bullet_\fp(\rV^\circ,\rK^{p\circ})(\kappa)
\]
where $\kappa$ is a perfect field containing $\dF^\Phi_p$. Put $\rB^\bullet_{s^\bullet}\coloneqq\pi^{\bullet-1}(s^\bullet)$, and denote by $(\cA,\lambda,\eta^p;\gamma)$ the universal object over the fiber $\rB^\bullet_{s^\bullet}$.
\begin{enumerate}
  \item The fiber $\rB^\bullet_{s^\bullet}$ is a smooth scheme over $\kappa$, whose tangent sheaf $\cT_{\rB^\bullet_{s^\bullet}/\kappa}$ fits canonically into an exact sequence
      \[
      0 \to \HOM\(\omega_{\cA^\vee,\tau_\infty},\omega_{\cA^\vee,\tau_\infty^\tc}^\perp/\omega_{\cA^\vee,\tau_\infty}\)
      \to \cT_{\rB^\bullet_{s^\bullet}/\kappa}
      \to\HOM\(\omega_{\cA^\vee,\tau_\infty^\tc}/(\Ker\gamma_{*,\tau_\infty})^\perp,\Lie_{\cA^\vee,\tau_\infty^\tc}\)  \to 0.
      \]

  \item The restriction of $\iota^\bullet_\kappa$ to $\rB^\bullet_{s^\bullet}$ is locally on $\rB^\bullet_{s^\bullet}$ a  closed immersion, with a canonical isomorphism for its normal sheaf
      \[
      \cN_{\iota^\bullet_\kappa\res\rB^\bullet_{s^\bullet}}\simeq
      \HOM\((\Ker\gamma_{*,\tau_\infty})^\perp/\rH^\dr_1(\cA)_{\tau_\infty}^\perp,\Lie_{\cA^\vee,\tau_\infty^\tc}\)\simeq
      \(\IM\gamma_{*,\tau_\infty}\)\otimes_{\cO_{\rB^\bullet_{s^\bullet}}}\Lie_{\cA^\vee,\tau_\infty^\tc}.
      \]

  \item We have $\gamma_{*,\tau_\infty^\tc}(\Ker\gamma_{*,\tau_\infty})^\perp=\rH^\dr_1(A^\bullet/S)_{\tau_\infty}^\perp$.

  \item The assignment sending
        $(A_0,\lambda_0,\eta_0^p;A,\lambda,\eta^p;A^\bullet,\lambda^\bullet,\eta^{p\bullet};\gamma)\in\rB^\bullet_{s^\bullet}(S)$ to the subbundles
        \begin{align*}
        H_1\coloneqq((\breve\gamma_{*,\tau_\infty})^{-1}\omega_{A^\vee/S,\tau_\infty})^\perp&\subseteq
        \rH^\dr_1(A^\bullet/S)_{\tau_\infty^\tc}=\rH^\dr_1(A^\bullet/\kappa)_{\tau_\infty^\tc}\otimes_\kappa\cO_S=(\sV_{s^\bullet})_S,\\
        H_2\coloneqq\gamma_{*,\tau_\infty^\tc}\omega_{A^\vee/S,\tau_\infty^\tc}&\subseteq
        \rH^\dr_1(A^\bullet/S)_{\tau_\infty^\tc}=\rH^\dr_1(A^\bullet/\kappa)_{\tau_\infty^\tc}\otimes_\kappa\cO_S=(\sV_{s^\bullet})_S,
        \end{align*}
        where $\breve\gamma\colon A^\bullet\to A$ is the (unique) $O_F$-linear quasi-$p$-isogeny such that $\breve\gamma\circ\gamma=\varpi\cdot\id_A$, induces an isomorphism
        \[
        \zeta^\bullet_{s^\bullet}\colon\rB^\bullet_{s^\bullet}\xrightarrow{\sim}\DL^\bullet_{s^\bullet}
        =\DL^\bullet(\sV_{s^\bullet},\{\;,\;\}_{s^\bullet}).
        \]
        In particular, $\rB^\bullet_{s^\bullet}$ is a geometrically irreducible projective smooth scheme in $\Sch_{/\kappa}$ of dimension $\floor{\tfrac{N}{2}}$.

  \item If we denote by $(\cH_{s^\bullet 1},\cH_{s^\bullet 2})$ the universal object over $\DL^\bullet_{s^\bullet}$, then there is a canonical isomorphism
      \[
      \zeta^{\bullet*}_{s^\bullet}\(\cH_{s^\bullet 1}^\dashv/\cH_{s^\bullet 2}\)\simeq\iota^{\bullet*}\Lie_{\cA,\tau_\infty^\tc}
      \]
      of line bundles on $\rB^\bullet_{s^\bullet}$.
\end{enumerate}
\end{theorem}

\begin{proof}
By Lemma \ref{le:inverse_isogeny}(2,3) and Definition \ref{de:ns_basic_correspondence2}, we have
\begin{align*}
\rank_{\cO_S}(\Ker\gamma_{*,\tau_\infty})+\rank_{\cO_S}(\Ker\gamma_{*,\tau_\infty^\tc})&=2\lfloor\tfrac{N}{2}\rfloor+1,\\
\rank_{\cO_S}(\Ker\gamma_{*,\tau_\infty})-\rank_{\cO_S}(\Ker\gamma_{*,\tau_\infty^\tc})&=1,
\end{align*}
which imply
\begin{align}\label{eq:ns_basic_correspondence21}
\rank_{\cO_S}(\Ker\gamma_{*,\tau_\infty})=\ceil{\tfrac{N+1}{2}},\quad\rank_{\cO_S}(\Ker\gamma_{*,\tau_\infty^\tc})=\ceil{\tfrac{N-1}{2}}.
\end{align}
Note that under Definition \ref{de:ns_basic_correspondence2}(a,b,d), Definition \ref{de:ns_basic_correspondence2}(c) is equivalent to that $(\Ker\gamma_{*,\tau_\infty})^\perp$ is a subbundle of $\rH^\dr_1(A/S)_{\tau_\infty^\tc}$ of rank $\ceil{\tfrac{N}{2}}$.

For an object $(A_0,\lambda_0,\eta_0^p;A,\lambda,\eta^p;A^\bullet,\lambda^\bullet,\eta^{p\bullet};\gamma)\in\rB^\bullet_\fp(\rV^\circ,\rK^{p\circ})(S)$, Definition \ref{de:ns_basic_correspondence2}(a) implies that there is a (unique) $O_F$-linear quasi-$p$-isogeny $\breve\gamma\colon A^\bullet\to A$ such that $\breve\gamma\circ\gamma=\varpi\cdot\id_A$, hence $\gamma\circ\breve\gamma=\varpi\cdot\id_{A^\bullet}$. Moreover, we have the following properties from Definition \ref{de:ns_basic_correspondence2}:
\begin{enumerate}[label=(\alph*')]
  \item $\Ker\breve\gamma[p^\infty]$ is contained in $A^\bullet[\fp]$;

  \item $(\IM\breve\gamma_{*,\tau_\infty})^\perp$ is contained in $\omega_{A^\vee,\tau_\infty^\tc}$;

  \item $\IM\breve\gamma_{*,\tau_\infty}$ contains $\rH^\dr_1(A/S)_{\tau_\infty^\tc}^\perp$;

  \item we have $\varpi\cdot\lambda^\bullet=\breve\gamma^\vee\circ\lambda\circ\breve\gamma$; and

  \item the $\rK^p$-orbit of maps $v\mapsto\varpi^{-1}\breve\gamma_*\circ\eta^{\bullet p}(v)$ for $v\in\rV^\circ\otimes_\dQ\dA^{\infty,p}$ coincides with $\eta^p$.
\end{enumerate}

First, we show (1). It is clear that $\rB^\bullet_{s^\bullet}$ is a scheme of finite type over $\kappa$. Consider a closed immersion $S\hookrightarrow\hat{S}$ in $\Sch'_{/\kappa}$ defined by an ideal sheaf $\cI$ satisfying $\cI^2=0$. Take a point $x=(A_0,\lambda_0,\eta_0^p;A,\lambda,\eta^p;A^\bullet,\lambda^\bullet,\eta^{p\bullet};\gamma)\in\rB^\bullet_{s^\bullet}(S)$. To compute lifting of $x$ to $\hat{S}$, we use the Serre--Tate and Grothendieck--Messing theories. Note that lifting $\gamma$ is equivalent to lifting both $\gamma$ and $\breve\gamma$, satisfying (b--e) in Definition \ref{de:ns_basic_correspondence2} and (b'--e') above, respectively. Thus, by Proposition \ref{pr:deformation}, to lift $x$ to an $\hat{S}$-point is equivalent to lifting
\begin{itemize}[label={\ding{109}}]
  \item $\omega_{A^\vee/S,\tau_\infty}$ to a subbundle $\hat\omega_{A^\vee,\tau_\infty}$ of $\rH^\cris_1(A/\hat{S})_{\tau_\infty}$ (of rank $1$),

  \item $\omega_{A^\vee/S,\tau_\infty^\tc}$ to a subbundle $\hat\omega_{A^\vee,\tau_\infty^\tc}$ of $\rH^\cris_1(A/\hat{S})_{\tau_\infty^\tc}$ (of rank $N-1$),
\end{itemize}
subject to the following requirements
\begin{enumerate}[label=(\alph*'')]
  \item $\hat\omega_{A^\vee,\tau_\infty}$ and $\hat\omega_{A^\vee,\tau_\infty^\tc}$ are orthogonal under $\langle\;,\;\rangle_{\lambda,\tau_\infty}^\cris$ \eqref{eq:weil_pairing_cris};

  \item $(\breve\gamma_{*,\tau_\infty}\rH^\cris_1(A^\bullet/\hat{S})_{\tau_\infty})^\perp$ is contained in $\hat\omega_{A^\vee,\tau_\infty^\tc}$.
\end{enumerate}
As $\breve\gamma_{*,\tau_\infty}\rH^\cris_1(A^\bullet/\hat{S})_{\tau_\infty}=\Ker\gamma_{*,\tau_\infty}
\subseteq\rH^\cris_1(A/\hat{S})_{\tau_\infty}$, (b'') is equivalent to
\begin{enumerate}
  \item [(c'')] $(\Ker\gamma_{*,\tau_\infty})^\perp$ is contained in $\hat\omega_{A^\vee,\tau_\infty^\tc}$.
\end{enumerate}
To summarize, lifting  $x$ to an $\hat{S}$-point is equivalent to lifting $\omega_{A^\vee/S,\tau_\infty^\tc}$ to a subbundle $\hat\omega_{A^\vee,\tau_\infty^\tc}$ of $\rH^\cris_1(A/\hat{S})_{\tau_\infty^\tc}$ containing $(\Ker\gamma_{*,\tau_\infty})^\perp$, and then lifting $\omega_{A^\vee/S,\tau_\infty}$ to a subbundle $\hat\omega_{A^\vee,\tau_\infty}$ of $\hat\omega_{A^\vee,\tau_\infty^\tc}^\perp$. Thus, (1) follows.

Next, we show (2). By Theorem \ref{th:ns_moduli_scheme}(4), the map $\cT_{\rB^\bullet_{s^\bullet}/\kappa}\to\iota^{\bullet*}\cT_{\rM^\bullet_\fp(\rV^\circ,\rK^{p\circ})/\kappa}\res_{\rB^\bullet_{s^\bullet}}$
is induced by the canonical map
\[
\HOM\(\omega_{\cA^\vee,\tau_\infty^\tc}/(\Ker\gamma_{*,\tau_\infty})^\perp,\Lie_{\cA^\vee,\tau_\infty^\tc}\)
\to\HOM\(\omega_{\cA^\vee,\tau_\infty^\tc}/\rH^\dr_1(\cA)_{\tau_\infty}^\perp,\Lie_{\cA^\vee,\tau_\infty^\tc}\).
\]
It is clearly injective, whose cokernel is canonically isomorphic to
\begin{align*}
&\HOM\((\Ker\gamma_{*,\tau_\infty})^\perp/\rH^\dr_1(\cA)_{\tau_\infty}^\perp,\Lie_{\cA^\vee,\tau_\infty^\tc}\)\\
&\simeq\HOM\((\IM\gamma_{*,\tau_\infty})^\vee,\Lie_{\cA^\vee,\tau_\infty^\tc}\)
\simeq\(\IM\gamma_{*,\tau_\infty}\)\otimes_{\cO_{\rB^\bullet_{s^\bullet}}}\Lie_{\cA^\vee,\tau_\infty^\tc}.
\end{align*}
We obtain (2).

Next, we show (3). By Definition \ref{de:ns_basic_correspondence2}(d) and the definition of $\breve\gamma$, we have $\lambda\circ\breve\gamma=\gamma^\vee\circ\lambda^\bullet$, which implies
\begin{align}\label{eq:ns_basic_correspondence22}
(\Ker\gamma_{*,\tau_\infty})^\perp=\gamma_{*,\tau_\infty^\tc}^{-1}(\rH^\dr_1(A^\bullet/S)_{\tau_\infty}^\perp).
\end{align}
It remains to show that $\rH^\dr_1(A^\bullet/S)_{\tau_\infty}^\perp$ is contained in $\IM\gamma_{*,\tau_\infty^\tc}=\Ker\breve\gamma_{*,\tau_\infty^\tc}$. By Definition \ref{de:ns_basic_correspondence2}(c), we know that $\breve\gamma_{*,\tau_\infty}^{-1}(\rH^\dr_1(A/S)_{\tau_\infty^\tc}^\perp)$ is a subbundle of $\rH^\dr_1(A^\bullet/S)_{\tau_\infty}$ of rank $\ceil{\tfrac{N}{2}}$. Similarly to \eqref{eq:ns_basic_correspondence22}, we have $(\Ker\breve\gamma_{*,\tau_\infty^\tc})^\perp=\breve\gamma_{*,\tau_\infty}^{-1}(\rH^\dr_1(A/S)_{\tau_\infty^\tc}^\perp)$, which is also a subbundle of $\rH^\dr_1(A^\bullet/S)_{\tau_\infty}$ of rank $\ceil{\tfrac{N}{2}}$. Thus, $\Ker\breve\gamma_{*,\tau_\infty^\tc}$ contains $\rH^\dr_1(A^\bullet/S)_{\tau_\infty}^\perp$.

Next, we show (4). We first show that $\zeta^\bullet_{s^\bullet}$ has the correct image, namely, we check
\begin{itemize}[label={\ding{109}}]
  \item $\rank_{\cO_S}H_1=\ceil{\tfrac{N}{2}}$ and $\rank_{\cO_S}H_2=\ceil{\tfrac{N}{2}}-1$: By \ref{eq:ns_basic_correspondence21}, we obtain $\rank_{\cO_S}H_1=\ceil{\tfrac{N}{2}}$. Since $\Ker\gamma_{*,\tau_\infty^\tc}\subseteq(\Ker\gamma_{*,\tau_\infty})^\perp\subseteq\omega_{A^\vee/S,\tau_\infty^\tc}$, we have $H_2=\gamma_{*,\tau_\infty^\tc}\omega_{A^\vee/S,\tau_\infty^\tc}\simeq\omega_{A^\vee/S,\tau_\infty^\tc}/\Ker\gamma_{*,\tau_\infty^\tc}$. Thus, we obtain $\rank_{\cO_S}H_2=\ceil{\tfrac{N}{2}}-1$.

  \item $\rH^\dr_1(A^\bullet/S)_{\tau_\infty}^\perp\subseteq H_2$: By Definition \ref{de:ns_basic_correspondence2}(b), $H_2$ contains $\gamma_{*,\tau_\infty^\tc}(\Ker\gamma_{*,\tau_\infty})^\perp$ in which the latter coincides with $\rH^\dr_1(A^\bullet/S)_{\tau_\infty}^\perp$ by (3).

  \item $H_2\subseteq H_1$: As $\lambda\circ\breve\gamma=\gamma^\vee\circ\lambda^\bullet$, we have
      \[
      \langle(\breve\gamma_{*,\tau_\infty})^{-1}\omega_{A^\vee/S,\tau_\infty},
      \gamma_{*,\tau_\infty^\tc}\omega_{A^\vee/S,\tau_\infty^\tc}\rangle_{\lambda^\bullet,\tau_\infty}
      =\langle\breve\gamma_{*,\tau_\infty}(\breve\gamma_{*,\tau_\infty})^{-1}\omega_{A^\vee/S,\tau_\infty},
      \omega_{A^\vee/S,\tau_\infty^\tc}\rangle_{\lambda,\tau_\infty}=0.
      \]
      Thus, we have $H_2\subseteq H_1$.

  \item $H_2\subseteq H_1^\dashv$: Note that we have
      \[
      \IM\gamma_{*,\tau_\infty^\tc}=
      \Ker\breve\gamma_{*,\tau_\infty^\tc}=(\breve\gamma_{*,\tau_\infty^\tc})^{-1}(\tF\omega_{A^\vee/S,\tau_\infty}^{(p)})\subseteq
      \tF((\breve\gamma_{*,\tau_\infty})^{-1}\omega_{A^\vee/S,\tau_\infty})=\tF((H_1^{(p)})^\perp).
      \]
      Thus, $(\tF((H_1^{(p)})^\perp))^\perp\subseteq(\IM\gamma_{*,\tau_\infty^\tc})^\perp$, which in turn implies $H_1^{(p)}\subseteq\tV((\IM\gamma_{*,\tau_\infty^\tc})^\perp)$, which further implies $\tV^{-1}H_1^{(p)}\subseteq(\IM\gamma_{*,\tau_\infty^\tc})^\perp$, which implies $\IM\gamma_{*,\tau_\infty^\tc}\subseteq H_1^\dashv$. By comparing ranks via \eqref{eq:ns_basic_correspondence21}, we obtain
      \begin{align}\label{eq:ns_basic_correspondence23}
      \IM\gamma_{*,\tau_\infty^\tc}= H_1^\dashv.
      \end{align}
      In particular, $H_1^\dashv$ contains $H_2$ as $\IM\gamma_{*,\tau_\infty^\tc}$ does.

  \item $H_1\subseteq H_2^\dashv$: Note that $H_2^{(p)}=\gamma_{*,\tau_\infty}(\tV\rH^\dr_1(A/S)_{\tau_\infty^\tc})=\tV(\IM\gamma_{*,\tau_\infty})=
      \tV(\Ker\breve\gamma_{*,\tau_\infty})\subseteq\tV(H_1^\perp)$. Thus, $\tV^{-1}H_2^{(p)}\subseteq H_1^\perp$, which implies $H_1\subseteq(\tV^{-1}H_2^{(p)})^\perp=H_2^\dashv$.

  \item $H_1^\dashv\subseteq H_2^\dashv$: This follows from $H_2\subseteq H_1$.
\end{itemize}

Since the target of $\zeta^\bullet_{s^\bullet}$ is smooth over $\kappa$ by Proposition \ref{pr:dl_bullet}, to see that $\zeta^\bullet_{s^\bullet}$ is an isomorphism, it suffices to check that for every algebraically closed field $\kappa'$ containing $\kappa$, the following statements hold:
\begin{enumerate}
  \item [(4--1)] $\zeta^\bullet_{s^\bullet}$ induces a bijection on $\kappa'$-points; and

  \item [(4--2)] $\zeta^\bullet_{s^\bullet}$ induces an isomorphism on the tangent spaces at every $\kappa$-point.
\end{enumerate}
To ease notation, we may assume $\kappa'=\kappa$.

For (4--1), we construct an inverse to the map $\zeta^\bullet_{s^\bullet}(\kappa)$. Take a point $y\in\DL^\bullet_{s^\bullet}(\kappa)$ represented by $\kappa$-linear subspaces
\[
\rH^\dr_1(A^\bullet)_{\tau_\infty}^\perp\subseteq H_2\subseteq H_1\subseteq\sV_{s^\bullet}=\rH^\dr_1(A^\bullet)_{\tau_\infty^\tc}.
\]
We regard $\tF$ and $\tV$ as those sesquilinear maps in Notation \ref{re:frobenius_verschiebung}. For every $\tau\in\Sigma_\infty$, we define a $W(\kappa)$-submodule $\cD_{A,\tau}\subseteq\cD(A^\bullet)_\tau$ as follows.
\begin{itemize}[label={\ding{109}}]
  \item If $\tau\not\in\{\tau_\infty,\tau_\infty^\tc\}$, then $\cD_{A,\tau}=\cD(A^\bullet)_\tau$.

  \item We set $\cD_{A,\tau_\infty}\coloneqq\tV^{-1}\tilde{H}_2$, where $\tilde{H}_2$ is the preimage of $H_2$ under the reduction map $\cD(A^\bullet)_{\tau_\infty^\tc}\to\cD(A^\bullet)_{\tau_\infty^\tc}/p\cD(A^\bullet)_{\tau_\infty^\tc}
      =\rH^\dr_1(A^\bullet)_{\tau_\infty^\tc}$.

  \item We set $\cD_{A,\tau_\infty^\tc}\coloneqq\tF\tilde{H}_1^\tc$, where $\tilde{H}_1^\tc$ is the preimage of $H_1^\perp$ under the reduction map $\cD(A^\bullet)_{\tau_\infty}\to\cD(A^\bullet)_{\tau_\infty}/p\cD(A^\bullet)_{\tau_\infty}=\rH^\dr_1(A^\bullet)_{\tau_\infty}$.
\end{itemize}
Finally, put $\cD_A\coloneqq\bigoplus_{\tau\in\Sigma_\infty}\cD_{A,\tau}$ as a $W(\kappa)$-submodule of $\cD(A^\bullet)$. We show that it is stable under $\tF$ and $\tV$. It suffices to show that both $\tF$ and $\tV$ stabilize $\cD_{A,\tau_\infty}\oplus\cD_{A,\tau_\infty^\tc}$, which breaks into checking that
\begin{itemize}[label={\ding{109}}]
  \item $\tF\cD_{A,\tau_\infty}\subseteq\cD_{A,\tau_\infty^\tc}$, that is, $\tF\tV^{-1}\tilde{H}_2\subseteq\tF\tilde{H}_1^\tc$. It suffices to show that $\tV^{-1}H_2$ (as a subspace of $\rH^\dr_1(A^\bullet)_{\tau_\infty}$) is contained in $H_1^\perp$, which follow from the relation $H_1\subseteq H_2^\dashv$.

  \item $\tF\cD_{A,\tau_\infty^\tc}\subseteq\cD_{A,\tau_\infty}$, that is, $\tF\tF\tilde{H}_1^\tc\subseteq\tV^{-1}\tilde{H}_2$. It suffices to show $p\tF\tilde{H}_1^\tc\subseteq\tilde{H}_2$, which obviously holds.

  \item $\tV\cD_{A,\tau_\infty}\subseteq\cD_{A,\tau_\infty^\tc}$, that is, $\tV\tV^{-1}\tilde{H}_2\subseteq\tF\tilde{H}_1^\tc$. it suffices to show $H_2\subseteq\tF H_1^\perp$, which follows from the identity $\tF H_1^\perp=(\tV^{-1}H_1)^\perp$ and the relation $H_2\subseteq H_1^\dashv$.

  \item $\tV\cD_{A,\tau_\infty^\tc}\subseteq\cD_{A,\tau_\infty}$, that is, $\tV\tF\tilde{H}_1^\tc\subseteq\tV^{-1}\tilde{H}_2$. It is obvious as $\tV^{-1}\tilde{H}_2$ contains $p\cD(A^\bullet)_{\tau_\infty}$.
\end{itemize}
Thus, $(\cD_A,\tF,\tV)$ is a Dieudonn\'{e} module over $W(\kappa)$. By the Dieudonn\'{e} theory, there is an $O_F$-abelian scheme $A$ over $\kappa$ with $\cD(A)_\tau=\cD_{A,\tau}$ for every $\tau\in\Sigma_\infty$, and an $O_F$-linear isogeny $\gamma\colon A\to A^\bullet$ inducing the inclusion of Dieudonn\'{e} modules $\cD(A)=\cD_A\subseteq\cD(A^\bullet)$. Moreover, since $\fp\cD(A^\bullet)\subseteq\cD(A)$, we have $\Ker\gamma[p^\infty]\subseteq A[\fp]$. Now we check that $(\Ker\gamma_{*,\tau_\infty})^\perp$ is contained in $\omega_{A^\vee/S,\tau_\infty^\tc}$, which is equivalent to that $p\cD(A^\bullet)_{\tau_\infty}^\vee\cap\cD(A)_{\tau_\infty^\tc}\subseteq\tV\cD(A)_{\tau_\infty}$. However, as $H_2$ contains $\rH^\dr_1(A^\bullet)_{\tau_\infty}^\perp$, we have $p\cD(A^\bullet)_{\tau_\infty}^\vee\subseteq\tilde{H}_2=\tV\cD(A)_{\tau_\infty}$.

Let $\lambda\colon A\to A^\vee$ be the unique quasi-polarization such that $\varpi\lambda=\gamma^\vee\circ\lambda^\bullet\circ\gamma$. We claim that $\lambda[p^\infty]$ is a polarization whose kernel is contained in $A[\fp]$ of rank $p^2$. Since $H_2\subseteq H_1$, we have $\langle\tilde{H}_1^\tc,\tilde{H}_2\rangle_{\lambda^\bullet,\tau_\infty}\subseteq pW(\kappa)$, which implies $\langle\cD(A)_{\tau_\infty},\cD(A)_{\tau_\infty^\tc}\rangle_{\lambda^\bullet,\tau_\infty}\subseteq pW(\kappa)$. It is enough to show
that the inclusion $\cD(A)_{\tau_\infty^\tc}\to\cD(A)_{\tau_\infty}^\vee$ induced from $\langle\;,\;\rangle_{\lambda^\bullet,\tau_\infty}$ has cokernel of length $N+1$. This follows from the facts that the cokernel of $\cD(A^\bullet)_{\tau_\infty^\tc}\hookrightarrow\cD(A^\bullet)_{\tau_\infty}^\vee$ has length $N-2\floor{\tfrac{N}{2}}$, and the cokernel of $\cD(A)_{\tau_\infty}\oplus\cD(A)_{\tau_\infty^\tc}\hookrightarrow
\cD(A^\bullet)_{\tau_\infty}\oplus\cD(A^\bullet)_{\tau_\infty^\tc}$ has length $2\floor{\tfrac{N}{2}}+1$.

It is an easy consequence of Lemma \ref{le:inverse_isogeny}(2) that the $O_F$-abelian scheme $A$ has signature type $N\Phi-\tau_\infty+\tau_\infty^\tc$. Finally, let $\eta^p$ be the unique $\rK^p$-level structure such that Definition \ref{de:qs_basic_correspondence}(d) is satisfied. Putting together, we obtain a point
$x=(A_0,\lambda_0,\eta_0^p;A,\lambda,\eta^p;A^\bullet,\lambda^\bullet,\eta^{p\bullet};\gamma)\in\rB^\bullet_{s^\bullet}(\kappa)$ such that $\zeta^\bullet_{s^\bullet}(x)=y$. It is easy to see that such assignment gives rise to an inverse of $\zeta^\bullet_{s^\bullet}(\kappa)$, hence (4--1) follows immediately.

For (4--2), let $\cT_x$ and $\cT_y$ be the tangent spaces at $x$ and $y$ as in (4--1), respectively. By Proposition \ref{pr:dl_bullet} and the construction, the induced map $(\zeta^\bullet_{s^\bullet})_*\colon\cT_x\to\cT_y$ fits into a commutative diagram
\[
\xymatrix{
\Hom_\kappa\(\omega_{A^\vee,\tau_\infty},\omega_{A^\vee,\tau_\infty^\tc}^\perp/\omega_{A^\vee,\tau_\infty}\) \ar[r]\ar[d]
& \cT_x \ar[r]\ar[d]^-{(\zeta^\bullet_{s^\bullet})_*} & \Hom_\kappa\(\omega_{A^\vee,\tau_\infty^\tc}/(\Ker\gamma_{*,\tau_\infty})^\perp,\Lie_{A^\vee,\tau_\infty^\tc}\) \ar[d] \\
\Hom_\kappa\(H_1/H_2,H_2^\dashv/H_1\)  \ar[r]& \cT_y  \ar[r]& \Hom_\kappa(H_2/\sV_{s^\bullet}^\dashv,H_1^\dashv/H_2)
}
\]
in $\Mod(\kappa)$. The right vertical arrow is induced by maps
\[
\omega_{A^\vee,\tau_\infty^\tc}/(\Ker\gamma_{*,\tau_\infty})^\perp\xrightarrow{\gamma_{*,\tau_\infty^\tc}}H_2/\sV_{s^\bullet}^\dashv,\quad
\Lie_{A^\vee,\tau_\infty^\tc}\simeq\rH^\dr_1(A)_{\tau_\infty^\tc}/\omega_{A^\vee,\tau_\infty^\tc}
\xrightarrow{\gamma_{*,\tau_\infty^\tc}} H_1^\dashv/H_2
\]
which are both isomorphisms by \eqref{eq:ns_basic_correspondence22} and \eqref{eq:ns_basic_correspondence23}, respectively. The left vertical arrow is the composition
\[
\Hom_\kappa\(\omega_{A^\vee,\tau_\infty},\omega_{A^\vee,\tau_\infty^\tc}^\perp/\omega_{A^\vee,\tau_\infty}\)
\to\Hom_\kappa\(H_1^\perp/\tV^{-1} H_2,H_2^\perp/H_1^\perp\)\xrightarrow{\sim}\Hom_\kappa\(H_1/H_2,H_2^\dashv/H_1\)
\]
in which the first arrow is induced by maps
\[
H_1^\perp/\tV^{-1} H_2\xrightarrow{\breve\gamma_{*,\tau_\infty}}\omega_{A^\vee,\tau_\infty},\quad
H_2^\perp/H_1^\perp\xrightarrow{\breve\gamma_{*,\tau_\infty}}\omega_{A^\vee,\tau_\infty^\tc}^\perp/\omega_{A^\vee,\tau_\infty}
\]
which are both isomorphisms as $\breve\gamma_{*,\tau_\infty}(H_1^\perp)=\omega_{A^\vee,\tau_\infty}$, $\breve\gamma_{*,\tau_\infty}(\tV^{-1} H_2)=0$, and $\breve\gamma_{*,\tau_\infty}(H_2^\perp)=\omega_{A^\vee,\tau_\infty^\tc}^\perp$. Thus, $(\zeta^\bullet_{s^\bullet})_*\colon\cT_x\to\cT_y$ is an isomorphism by the Five Lemma, hence (4--2) and (4) follow.

Finally, (5) is a consequence of \eqref{eq:ns_basic_correspondence23}.
\end{proof}

\begin{remark}
We have the following remarks concerning Theorem \ref{th:ns_basic_correspondence2}.
\begin{enumerate}
  \item When $\rK^{p\circ}$ is sufficiently small, the restriction of $\iota^\bullet_\kappa$ to $\rB^\bullet_{s^\bullet}$ is a closed immersion for every point $s^\bullet\in\rS^\bullet_\fp(\rV^\circ,\rK^{p\circ})(\kappa)$ and every perfect field $\kappa$ containing $\dF^\Phi_p$.

  \item In fact, one can show that the union of $\rM^\dag_\fp(\rV^\circ,\rK^{p\circ})$ and the image of $\iota^\bullet\colon\rB^\bullet_\fp(\rV^\circ,\rK^{p\circ})\to\rM^\bullet_\fp(\rV^\circ,\rK^{p\circ})$ is exactly the basic locus of $\rM^\bullet_\fp(\rV^\circ,\rK^{p\circ})$. In particular, as long as $N\geq 5$, the basic locus of $\rM^\bullet_\fp(\rV^\circ,\rK^{p\circ})$ is \emph{not} equidimensional.
\end{enumerate}
\end{remark}

\begin{construction}\label{cs:ns_uniformization2}
To construct a uniformization map for $\rS^\bullet_\fp(\rV^\circ,\obj)$, we need to choose an $O_{F_\fp}$-lattice $\Lambda^\bullet_\fp$ in $\rV^\circ\otimes_FF_\fp$ satisfying
\begin{itemize}[label={\ding{109}}]
  \item $\Lambda^\circ_\fp\subseteq\Lambda^\bullet_\fp\subseteq p^{-1}\Lambda^\circ_\fp$, and

  \item $p\Lambda^\bullet_\fp\subseteq(\Lambda^\bullet_\fp)^\vee$ such that $(\Lambda^\bullet_\fp)^\vee/p\Lambda^\bullet_\fp$ has length $0$ (resp.\ $1$) if $N$ is even (resp.\ odd).
\end{itemize}
Let $\rK^\bullet_\fp$ be the stabilizer of $\Lambda^\bullet_\fp$; and put $\rK^\bullet_p\coloneqq\rK^\bullet_\fp\times\prod_{\fq\mid p,\fq\neq\fp}\rK^\circ_\fq$. Similar to Construction \ref{cs:qs_uniformization}, we may construct a \emph{uniformization map}
\begin{align}\label{eq:ns_uniformization2}
\upsilon^\bullet\colon\rS^\bullet_\fp(\rV^\circ,\obj)(\ol\dF_p)
\xrightarrow{\sim}\Sh(\rV^\circ,\obj\rK^\bullet_p)\times\rT_\fp(\ol\dF_p)
\end{align}
in $\Fun(\fK(\rV^\circ)^p\times\fT,\Set)_{/\rT_\fp(\ol\dF_p)}$ which is an isomorphism, under which the induced action of $\Gal(\ol\dF_p/\dF_p^\Phi)$ on the target is trivial on $\Sh(\rV^\circ,\obj\rK^\bullet_p)$.

Moreover, similar to Construction \ref{cs:qs_definite_hecke} and Proposition \ref{pr:qs_definite_hecke}, for every $g\in\rK^\bullet_\fp\backslash\rU(\rV^\circ)(F^+_\fp)/\rK^\bullet_\fp$, we may construct the Hecke correspondence
\[
\Hk_g\colon\rS^\bullet_\fp(\rV^\circ,\obj)_g\to\rS^\bullet_\fp(\rV^\circ,\obj)\times\rS^\bullet_\fp(\rV^\circ,\obj)
\]
as a morphism in $\Fun(\fK(\rV^\circ)^p\times\fT,\Sch_{/\dF_p^\Phi})_{/\rT_\fp}$ that is finite \'{e}tale and compatible with the uniformization map.
\end{construction}

\subsection{Basic correspondence for the link stratum}

In this subsection, we construct and study the basic correspondence for the link stratum $\rM^\dag_\fp(\rV^\circ,\obj)$. We also discuss its relation with the two previously constructed basic correspondences. We assume $N\geq 2$.

\begin{definition}\label{de:ns_definite3}
We define a functor
\begin{align*}
\rS^\dag_\fp(\rV^\circ,\obj)\colon\fK(\rV^\circ)^p\times\fT &\to\sfP\Sch'_{/\dF^\Phi_p} \\
\rK^{p\circ} &\mapsto \rS^\dag_\fp(\rV^\circ,\rK^{p\circ})
\end{align*}
such that for every $S\in\Sch'_{/\dF^\Phi_p}$, $\rS^\dag_\fp(\rV^\circ,\rK^{p\circ})(S)$ is the set of equivalence classes of decuples $(A_0,\lambda_0,\eta_0^p;A^\circ,\lambda^\circ,\eta^{p\circ};A^\bullet,\lambda^\bullet,\eta^{p\bullet};\psi)$, where
\begin{itemize}[label={\ding{109}}]
  \item $(A_0,\lambda_0,\eta_0^p;A^\circ,\lambda^\circ,\eta^{p\circ})$ is an element in $\rS^\circ_\fp(\rV^\circ,\rK^{p\circ})(S)$;

  \item $(A_0,\lambda_0,\eta_0^p;A^\bullet,\lambda^\bullet,\eta^{p\bullet})$ is an element in $\rS^\bullet_\fp(\rV^\circ,\rK^{p\circ})(S)$; and

  \item $\psi\colon A^\circ\to A^\bullet$ is an $O_F$-linear quasi-$p$-isogeny (Definition \ref{de:p_quasi}) such that
  \begin{enumerate}[label=(\alph*)]
    \item $\Ker\psi[p^\infty]$ is contained in $A^\circ[\fp]$;

    \item we have $\varpi\cdot\lambda^\circ=\psi^\vee\circ\lambda^\bullet\circ\psi$; and

    \item the $\rK^{p\circ}$-orbit of maps $v\mapsto\psi_*\circ\eta^{p\circ}(v)$ for $v\in\rV^\circ\otimes_\dQ\dA^{\infty,p}$ coincides with $\eta^{p\bullet}$.
  \end{enumerate}
\end{itemize}
The equivalence relation and the action of morphisms in $\fK(\rV^\circ)^p\times\fT$ are defined similarly as in Definition \ref{de:qs_basic_correspondence}.
\end{definition}

We clearly have the forgetful morphism
\begin{align*}
\rS^\dag_\fp(\rV^\circ,\obj)\to\rT_\fp
\end{align*}
in $\Fun(\fK(\rV^\circ)^p\times\fT,\sfP\Sch'_{/\dF_p^\Phi})$, which is represented by finite and \'{e}tale schemes.

By definition, we have the two forgetful morphisms
\[
\rs^{\dag\circ}\colon\rS^\dag_\fp(\rV^\circ,\obj)\to\rS^\circ_\fp(\rV^\circ,\obj),\quad
\rs^{\dag\bullet}\colon\rS^\dag_\fp(\rV^\circ,\obj)\to\rS^\bullet_\fp(\rV^\circ,\obj)
\]
in $\Fun(\fK(\rV^\circ)^p\times\fT,\Sch_{/\dF_p^\Phi})_{/\rT_\fp}$.

\begin{definition}\label{de:ns_basic_correspondence3}
We define $\rB^\dag_\fp(\rV^\circ,\obj)$ to be the limit of the following diagram
\[
\xymatrix{
\rS^\circ_\fp(\rV^\circ,\obj)  &
\rB^\circ_\fp(\rV^\circ,\obj) \ar[r]^-{\iota^\circ}\ar[l]_-{\pi^\circ} &
\rM^\circ_\fp(\rV^\circ,\obj) \\
\rS^\dag_\fp(\rV^\circ,\obj) \ar[u]_-{\rs^{\dag\circ}}\ar[d]^-{\rs^{\dag\bullet}}
&& \rM^\dag_\fp(\rV^\circ,\obj) \ar[u]_-{\rm^{\dag\circ}}\ar[d]^-{\rm^{\dag\bullet}} \\
\rS^\bullet_\fp(\rV^\circ,\obj)  &
\rB^\bullet_\fp(\rV^\circ,\obj) \ar[r]^-{\iota^\bullet}\ar[l]_-{\pi^\bullet} &
\rM^\bullet_\fp(\rV^\circ,\obj)
}
\]
in the category $\Fun(\fK(\rV^\circ)^p\times\fT,\Sch_{/\dF_p^\Phi})_{/\rT_\fp}$.
\end{definition}

From the definition above, we have the following commutative diagram
\begin{align}\label{eq:ns_link}
\xymatrix{
\rS^\circ_\fp(\rV^\circ,\obj)  &
\rB^\circ_\fp(\rV^\circ,\obj) \ar[r]^-{\iota^\circ}_-{\sim}\ar[l]_-{\pi^\circ} &
\rM^\circ_\fp(\rV^\circ,\obj) \\
& \rS^\dag_\fp(\rV^\circ,\obj) \ar[lu]_-{\rs^{\dag\circ}}\ar[rd]^-{\rs^{\dag\bullet}} &
\rB^\dag_\fp(\rV^\circ,\obj) \ar[r]^-{\iota^\dag}\ar[l]_-{\pi^\dag} \ar[lu]_-{\rb^{\dag\circ}}\ar[rd]^-{\rb^{\dag\bullet}} &
\rM^\dag_\fp(\rV^\circ,\obj) \ar[lu]_-{\rm^{\dag\circ}}\ar[rd]^-{\rm^{\dag\bullet}} \\
&& \rS^\bullet_\fp(\rV^\circ,\obj)  &
\rB^\bullet_\fp(\rV^\circ,\obj) \ar[r]^-{\iota^\bullet}\ar[l]_-{\pi^\bullet} &
\rM^\bullet_\fp(\rV^\circ,\obj)
}
\end{align}
in $\Fun(\fK(\rV^\circ)^p\times\fT,\Sch_{/\dF_p^\Phi})_{/\rT_\fp}$, together with the four new morphisms from $\rB^\dag_\fp(\rV^\circ,\obj)$ as indicated. It will be clear in \S\ref{ss:ns_functoriality} why we draw the diagram oblique.

\begin{theorem}\label{th:ns_basic_correspondence3}
In the diagram \eqref{eq:ns_link}, we have
\begin{enumerate}
  \item The square
    \[
    \xymatrix{
    \rB^\dag_\fp(\rV^\circ,\obj) \ar[r]^-{\iota^\dag} \ar[d]^-{\rb^{\dag\bullet}} & \rM^\dag_\fp(\rV^\circ,\obj) \ar[d]^-{\rm^{\dag\bullet}} \\
    \rB^\bullet_\fp(\rV^\circ,\obj) \ar[r]^-{\iota^\bullet} & \rM^\bullet_\fp(\rV^\circ,\obj)
    }
    \]
    is a Cartesian diagram.

  \item Take a point $s^\dag=(A_0,\lambda_0,\eta_0^p;A^\circ,\lambda^\circ,\eta^{p\circ};A^\bullet,\lambda^\bullet,\eta^{p\bullet};\psi)
    \in\rS^\dag_\fp(\rV^\circ,\rK^{p\circ})(\kappa)$ where $\kappa$ is a perfect field containing $\dF^\Phi_p$.
    Put $\rB^\dag_{s^\dag}\coloneqq\pi^{\dag-1}(s^\dag)$ and  $\sV_{s^\dag}\coloneqq(\IM\psi_{*,\tau_\infty^\tc})/\rH^\dr_1(A^\bullet/\kappa)_{\tau_\infty}^\perp$ which has dimension $\floor{\tfrac{N}{2}}$. Then the assignment sending
    \[
    ((A_0,\lambda_0,\eta_0^p;A,\lambda,\eta^p;A^\circ,\lambda^\circ,\eta^{p\circ};\beta),
    (A_0,\lambda_0,\eta_0^p;A,\lambda,\eta^p;A^\bullet,\lambda^\bullet,\eta^{p\bullet};\gamma))
    \in\rB^\dag_{s^\dag}(S)
    \]
    (with $\gamma=\psi\circ\beta$) to $
    (\gamma_{*,\tau_\infty^\tc}\omega_{A^\vee/S,\tau_\infty^\tc})/\rH^\dr_1(A^\bullet/S)_{\tau_\infty}^\perp$ induces an isomorphism
    \[
    \zeta^\dag_{s^\dag}\colon\rB^\dag_{s^\dag}\xrightarrow{\sim}\dP(\sV_{s^\dag}).
    \]

\end{enumerate}
\end{theorem}

\begin{proof}
For (1), unravelling all the definitions, it suffices to show that for every object
\[
((A_0,\lambda_0,\eta_0^p;A,\lambda,\eta^p;A^\circ,\lambda^\circ,\eta^{p\circ};\beta),
(A_0,\lambda_0,\eta_0^p;A,\lambda,\eta^p;A^\bullet,\lambda^\bullet,\eta^{p\bullet};\gamma))
\]
of $\rM^\dag_\fp(\rV^\circ,\rK^{p\circ})(S)\times_{\rM^\bullet_\fp(\rV^\circ,\rK^{p\circ})(S)}\rB^\bullet_\fp(\rV^\circ,\rK^{p\circ})(S)=
\rB^\circ_\fp(\rV^\circ,\rK^{p\circ})(S)\times_{\rM_\fp(\rV^\circ,\rK^{p\circ})(S)}\rB^\bullet_\fp(\rV^\circ,\rK^{p\circ})(S)$, the quasi-isogeny $\psi\coloneqq\gamma\circ\beta^{-1}\colon A^\circ\to A^\bullet$ is a quasi-$p$-isogeny. However, since $\beta_{*,\tau_\infty^\tc}\colon\rH^\dr_1(A)_{\tau_\infty^\tc}\to\rH^\dr_1(A^\circ)_{\tau_\infty^\tc}$ is an isomorphism and $\Ker\beta_{*,\tau_\infty}=\omega_{A^\vee,\tau_\infty}$, it suffices to show that $\omega_{A^\vee,\tau_\infty}$ is contained in $\Ker\gamma$, which is clear as $\omega_{A^{\bullet\vee},\tau_\infty}=0$.

For (2), we first show that for a point
\[
x^\bullet=(A_0,\lambda_0,\eta_0^p;A,\lambda,\eta^p;A^\bullet,\lambda^\bullet,\eta^{p\bullet};\gamma)
\in\rB^\bullet_\fp(\rV^\circ,\rK^{p\circ})(S),
\]
$\iota^\bullet(x^\bullet)$ belongs to $\rM^\dag_\fp(\rV^\circ,\rK^{p\circ})(S)$ if and only if $H_1=H_1^\dashv$, where we recall from Theorem \ref{th:ns_basic_correspondence2} that $H_1\coloneqq((\breve\gamma_{*,\tau_\infty})^{-1}\omega_{A^\vee,\tau_\infty})^\perp$. In fact, by Definition \ref{de:ns_moduli_scheme1}, $\iota^\bullet(x^\bullet)\in\rM^\dag_\fp(\rV^\circ,\rK^{p\circ})(S)$ if and only if $\omega_{A^\vee,\tau_\infty}=\rH^1_\dr(A)_{\tau_\infty^\tc}^\perp$. In the proof of Theorem \ref{th:ns_basic_correspondence2}, we see $\IM\gamma_{*,\tau_\infty^\tc}= H_1^\dashv$ \eqref{eq:ns_basic_correspondence23}. As $\lambda\circ\breve\gamma=\gamma^\vee\circ\lambda^\bullet$, we have $(\IM\gamma_{*,\tau_\infty^\tc})^\perp=(\breve\gamma_{*,\tau_\infty})^{-1}\rH^1_\dr(A)_{\tau_\infty^\tc}^\perp$. Thus, if $\omega_{A^\vee,\tau_\infty}=\rH^1_\dr(A)_{\tau_\infty^\tc}^\perp$, then $H_1=((\IM\gamma_{*,\tau_\infty^\tc})^\perp)^\perp$ which equals $\IM\gamma_{*,\tau_\infty^\tc}=H_1^\dashv$, as $\IM\gamma_{*,\tau_\infty^\tc}$ contains $\rH^\dr_1(A^\bullet)_{\tau_\infty}^\perp$. On the other hand, if $H_1=H_1^\dashv$, then $(\breve\gamma_{*,\tau_\infty})^{-1}\omega_{A^\vee,\tau_\infty}=(\IM\gamma_{*,\tau_\infty^\tc})^\perp
=(\breve\gamma_{*,\tau_\infty})^{-1}\rH^1_\dr(A)_{\tau_\infty^\tc}^\perp$, which implies easily that $\omega_{A^\vee,\tau_\infty}=\rH^1_\dr(A)_{\tau_\infty^\tc}^\perp$.

Second, we show $H_1=\IM\psi_{*,\tau_\infty^\tc}$ if $x^\bullet\in\rB^\dag_{s^\dag}(S)$. Since $\gamma=\psi\circ\beta$, we have $\IM\gamma_{*,\tau_\infty^\tc}\subseteq\IM\psi_{*,\tau_\infty^\tc}$. As $\IM\gamma_{*,\tau_\infty^\tc}=H_1^\dashv=H_1$, we have $H_1\subseteq\IM\psi_{*,\tau_\infty^\tc}$. On the other hand, it follows easily from Lemma \ref{le:inverse_isogeny}(2,3) that $\IM\psi_{*,\tau_\infty^\tc}$ has rank $\ceil{\tfrac{N}{2}}$. Thus, we must have $H_1=\IM\psi_{*,\tau_\infty^\tc}$.

The above two claims together with Theorem \ref{th:ns_basic_correspondence2}(4) imply (2).
\end{proof}

\begin{remark}\label{re:ns_basic_correspondence3}
It follows from the proof of Theorem \ref{th:ns_basic_correspondence3} that for every $s^\dag\in\rS^\dag_\fp(\rV^\circ,\rK^{p\circ})(\kappa)$, if we put $s^\circ\coloneqq\rs^{\dag\circ}(s^\dag)$ and $s^\bullet\coloneqq\rs^{\dag\bullet}(s^\dag)$, then
\begin{enumerate}
  \item the morphism $\zeta^\circ_{s^\circ}\circ\rb^{\dag\circ}\circ(\zeta^\dag_{s^\dag})^{-1}$ identifies $\dP(\sV_{s^\dag})$ as a closed subscheme of $\dP(\sV_{s^\circ})$ induced by the obvious $\kappa$-linear (surjective) map $\sV_{s^\circ}\to\sV_{s^\dag}$; and

  \item the morphism $\zeta^\bullet_{s^\bullet}\circ\rb^{\dag\bullet}\circ(\zeta^\dag_{s^\dag})^{-1}$ identifies $\dP(\sV_{s^\dag})$ as a closed subscheme (of codimension one) of $\DL^\bullet(\sV_{s^\bullet},\{\;,\;\}_{s^\bullet})$ defined by the condition $H_1=H_1^\dashv$.
\end{enumerate}
\end{remark}

\begin{construction}\label{cs:ns_uniformization3}
Put $\rK^\dag_p\coloneqq\rK^\circ_p\cap\rK^\bullet_p$. Similar to Construction \ref{cs:qs_uniformization}, we construct a \emph{uniformization map}
\begin{align}\label{eq:ns_uniformization3}
\upsilon^\dag\colon\rS^\dag_\fp(\rV^\circ,\obj)(\ol\dF_p)
\xrightarrow{\sim}\Sh(\rV^\circ,\obj\rK^\dag_p)\times\rT_\fp(\ol\dF_p)
\end{align}
in $\Fun(\fK(\rV^\circ)^p\times\fT,\Set)_{/\rT_\fp(\ol\dF_p)}$ which is an isomorphism, under which the induced action of $\Gal(\ol\dF_p/\dF_p^\Phi)$ on the target is trivial on $\Sh(\rV^\circ,\obj\rK^\dag_p)$.
\end{construction}

\subsection{Cohomology of the link stratum}
\label{ss:ns_cohomology}

In this subsection, we study the cohomology of the link stratum. We assume $N\geq 2$.

We first construct certain Hecke correspondences for $\rB^\circ_\fp(\rV^\circ,\obj)$ extending Construction \ref{cs:ns_uniformization1}. Unlike the functor $\rS^\circ_\fp(\rV^\circ,\obj)$, the natural action of $\rK^\circ_\fp=\rU(\Lambda^\circ_\fp)(O_{F^+_\fp})$ on the functor $\rB^\circ_\fp(\rV^\circ,\obj)$ is nontrivial. However, as we will see, such action factors through the quotient $\rU(\Lambda^\circ_\fp)(O_{F^+_\fp})\to\rU(\Lambda^\circ_\fp)(\dF_p)$. Let $\rK^\circ_{\fp1}$ be the kernel of the reduction map $\rK^\circ_\fp=\rU(\Lambda^\circ_\fp)(O_{F^+_\fp})\to\rU(\Lambda^\circ_\fp)(\dF_p)$.

\begin{construction}\label{cs:ns_balloon_hecke}
We first define a functor
\begin{align*}
\rS^\circ_{\fp1}(\rV^\circ,\obj)\colon\fK(\rV^\circ)^p\times\fT &\to\sfP\Sch'_{/\dF^\Phi_p} \\
\rK^{p\circ} &\mapsto \rS^\circ_\fp(\rV^\circ,\rK^{p\circ})
\end{align*}
such that for every $S\in\Sch'_{/\dF^\Phi_p}$, $\rS^\circ_{\fp1}(\rV^\circ,\rK^{p\circ})(S)$ is the set of equivalence classes of septuples $(A_0,\lambda_0,\eta_0^p;A^\circ,\lambda^\circ,\eta^{p\circ};\eta^\circ_\fp)$, where
\begin{itemize}[label={\ding{109}}]
  \item $(A_0,\lambda_0,\eta_0^p;A^\circ,\lambda^\circ,\eta^{p\circ})$ is an element in $\rS^\circ_\fp(\rV^\circ,\rK^{p\circ})(S)$;

  \item $\eta^\circ_\fp$ is, for a chosen geometric point $s$ on every connected component of $S$, an isomorphism
      \[
      \eta^\circ_\fp\colon\Lambda^\circ_\fp\otimes\dF_p\to\Hom_{O_F}(A_{0s}[\fp],A^\circ_s[\fp])
      \]
      of hermitian spaces over $O_{F_\fp}\otimes\dF_p$, where $\Hom_{O_F}(A_{0s}[\fp],A^\circ_s[\fp])$ is equipped with the hermitian form constructed similarly as in Construction \ref{cs:hermitian_structure} with respect to $(\lambda_0,\lambda^\circ)$.
\end{itemize}
The equivalence relation and the action of morphisms in $\fK(\rV^\circ)^p\times\fT$ are defined similarly as in Definition \ref{de:qs_moduli_scheme}. In fact, we have a further action of $\rU(\Lambda^\circ_\fp)(\dF_p)$ on $\rS^\circ_{\fp1}(\rV^\circ,\obj)$. Moreover, similar to Construction \ref{cs:qs_definite_hecke} and Proposition \ref{pr:qs_definite_hecke}, for every $g\in\rK^\circ_{\fp1}\backslash\rU(\rV^\circ)(F^+_\fp)/\rK^\circ_{\fp1}$, we may construct the Hecke correspondence
\begin{align}\label{eq:ns_balloon_hecke_pre}
\Hk_g\colon\rS^\circ_{\fp1}(\rV^\circ,\obj)_g\to\rS^\circ_{\fp1}(\rV^\circ,\obj)\times\rS^\circ_{\fp1}(\rV^\circ,\obj)
\end{align}
as a morphism in $\Fun(\fK(\rV^\circ)^p\times\fT,\Sch_{/\dF_p^\Phi})_{/\rT_\fp}$ that is finite \'{e}tale.

On the other hand, Theorem \ref{th:ns_basic_correspondence1} implies that we have a canonical isomorphism
\[
\rB^\circ_\fp(\rV^\circ,\obj)\simeq
\rS^\circ_{\fp1}(\rV^\circ,\obj)\overset{\rU(\Lambda^\circ_\fp)(\dF_p)}{\times}\dP(\Lambda^\circ_\fp\otimes\dF_p)
\]
in the category $\Fun(\fK(\rV^\circ)^p\times\fT,\Sch_{/\dF_p^\Phi})_{/\rT_\fp}$. Thus, for every $g\in\rK^\circ_{\fp1}\backslash\rU(\rV^\circ)(F^+_\fp)/\rK^\circ_{\fp1}$, we obtain from \eqref{eq:ns_balloon_hecke_pre} the Hecke correspondence
\begin{align*}
\Hk_g\colon\rB^\circ_\fp(\rV^\circ,\obj)_g\to\rB^\circ_\fp(\rV^\circ,\obj)\times\rB^\circ_\fp(\rV^\circ,\obj)
\end{align*}
as a morphism in $\Fun(\fK(\rV^\circ)^p\times\fT,\Sch_{/\dF_p^\Phi})_{/\rT_\fp}$ that is finite \'{e}tale.
\end{construction}

Now we study cohomology.

\begin{lem}\label{le:ns_cohomology}
Consider a $p$-coprime coefficient ring $L$.
\begin{enumerate}
  \item If $p+1$ is invertible in $L$, then the restriction map
     \[
     (\rm^{\dag\circ})^*\colon\rH^i_\fT(\ol\rM^\circ_\fp(\rV^\circ,\obj),L)\to\rH^i_\fT(\ol\rM^\dag_\fp(\rV^\circ,\obj),L)
     \]
     is an isomorphism for every integer $i\not\in\{N-2,2N-2\}$. In particular, $\rH^i_\fT(\ol\rM^\circ_\fp(\rV^\circ,\obj),L)$ and $\rH^i_\fT(\ol\rM^\dag_\fp(\rV^\circ,\obj),L)$ vanish if $i$ is odd and different from $N-2$.

  \item For every $i\in\dZ$, both $\rH^i_\fT(\ol\rM^\circ_\fp(\rV^\circ,\obj),L)$ and $\rH^i_\fT(\ol\rM^\dag_\fp(\rV^\circ,\obj),L)$ are free $L$-modules.

  \item When $N$ is even, the action of $\Gal(\ol\dF_p/\dF_p^\Phi)$ on $\rH^{N-2}_\fT(\ol\rM^\dag_\fp(\rV^\circ,\obj),L(\tfrac{N-2}{2}))$ is trivial.
\end{enumerate}
\end{lem}

\begin{proof}
By Theorem \ref{th:ns_basic_correspondence1}, for every $\rK^{p\circ}\in\fK(\rV^\circ)^p$ and every $s^\circ\in\rS^\circ_\fp(\rV^\circ,\rK^{p\circ})(\ol\dF_p)$, the restriction of $(\rm^{\dag\circ})^*$ to the fibers over $s^\circ$ is a morphism appearing in Lemma \ref{le:dl_cohomology}.

Part (1) then follows from Lemma \ref{le:dl_cohomology}(2). Part (2) follows from Lemma \ref{le:dl_cohomology}(3). Part (3) follows from Lemma \ref{le:dl_cohomology}(4) and Construction \ref{cs:ns_uniformization1}.
\end{proof}

\begin{definition}\label{de:ns_primitive}
Let $\xi\in\rH^2_{\fT}(\ol\rB^{\circ}_{\fp}(\rV^{\circ},\obj),L(1))$ be the first Chern class of the tautological quotient line bundle on $\ol\rB^{\circ}_{\fp}(\rV^{\circ},\obj)$ (that is, in the situation of Theorem \ref{th:ns_basic_correspondence1}, the restriction of $\xi$ to $\rB^\circ_{s^\circ}$ is isomorphic to $\zeta^{\circ*}_{s^\circ}\cO_{\dP(\sV_{s^\circ})}(1)$ for every $\rK^{p\circ}\in\fK(\rV^\circ)^p$ and every $s^\circ\in\rS^\circ_\fp(\rV^\circ,\rK^{p\circ})(\ol\dF_p)$). We define the \emph{primitive cohomology} $\rH^\prim(\ol\rM^\dag_\fp(\rV^\circ,\obj),L(i))$ to be the kernel of the map
\[
\cup(\rm^{\dag\circ*}\iota^\circ_!\xi)\colon
\rH^{N-2}_\fT(\ol\rM^\dag_\fp(\rV^\circ,\obj),L(i)) \to \rH^N_\fT(\ol\rM^\dag_\fp(\rV^\circ,\obj),L(i+1)),
\]
which is canonically a direct summand of $\rH^{N-2}_\fT(\ol\rM^\dag_\fp(\rV^\circ,\obj),L(i))$.
\end{definition}

\begin{proposition}\label{pr:ns_link_cohomology}
Take an object $\rK^{p\circ}\in\fK(\rV^\circ)^p$, a rational prime $\ell\neq p$, and an isomorphism $\iota_\ell\colon\dC\simeq\ol\dQ_\ell$. Then we have an isomorphism
\begin{align}\label{eq:tate}
\iota_\ell^{-1}\rH^\prim(\ol\rM^\dag_\fp(\rV^\circ,\rK^{p\circ}),\ol\dQ_\ell)\simeq
\Map_{\rK^\circ_\fp}\(\rU(\rV^\circ)(F^+)\backslash\rU(\rV^\circ)(\dA^\infty_{F^+})/\rK^{p\circ}\prod_{\fq\mid p,\fq\neq\fp}\rK^\circ_\fq,\Omega_N\)
\end{align}
of $\dC[\rK^{p\circ}\rK^\circ_{p1}\backslash\rU(\rV^\circ)(\dA^\infty_{F^+})/\rK^{p\circ}\rK^\circ_{p1}]$-modules, where $\Omega_N$ is the Tate--Thompson representation of $\rK^\circ_\fp$ introduced in \S\ref{ss:omega}. Moreover, let $\pi^{\infty,p}$ be an irreducible admissible representation of $\rU(\rV^\circ)(\dA^{\infty,p}_{F^+})$ such that $(\pi^{\infty,p})^{\rK^{p\circ}}$ is a constituent of $\iota_\ell^{-1}\rH^\prim(\ol\rM^\dag_\fp(\rV^\circ,\rK^{p\circ}),\ol\dQ_\ell)$. Then one can complete $\pi^{\infty,p}$ to an automorphic representation $\pi=\pi^{\infty,p}\otimes\pi_\infty\otimes\prod_{\fq\mid p}\pi_\fq$ of $\rU(\rV^\circ)(\dA_{F^+})$ such that $\pi_\infty$ is trivial; $\pi_\fq$ is unramified for $\fq\neq\fp$; and
\begin{enumerate}
  \item when $N$ is even, $\pi_\fp$ is a constituent of an unramified principal series;

  \item when $N$ is odd, $\BC(\pi_\fp)$ is a constituent of an unramified principal series of $\GL_N(F_\fp)$ whose Satake parameter contains $\{-p,-p^{-1}\}$.
\end{enumerate}
\end{proposition}

\begin{proof}
Put $\rK^\circ_{p1}\coloneqq\rK^\circ_{\fp 1}\times\prod_{\fq\mid p,\fq\neq\fp}\rK^\circ_\fq$. By Construction \ref{cs:ns_balloon_hecke}, the cohomology $\rH^{N-2}_\fT(\ol\rM^\dag_\fp(\rV^\circ,\rK^{p\circ}),\ol\dQ_\ell)$ is an $\ol\dQ_\ell[\rK^{p\circ}\rK^\circ_{p1}\backslash\rU(\rV^\circ)(\dA^\infty_{F^+})/\rK^{p\circ}\rK^\circ_{p1}]$-module for which $\rH^\prim(\ol\rM^\dag_\fp(\rV^\circ,\rK^{p\circ}),\ol\dQ_\ell)$ is a submodule.

In the uniformization map \eqref{eq:ns_uniformization1}, we let $s_0\in\rS^\circ_\fp(\rV^\circ,\rK^{p\circ})(\ol\dF_p)$ be the point corresponding to the unit element on the right-hand side. Put
\[
\rH^\prim_{s_0}(\ol\rM^\dag_\fp(\rV^\circ,\rK^{p\circ}),\ol\dQ_\ell)
\coloneqq\rH^\prim(\ol\rM^\dag_\fp(\rV^\circ,\rK^{p\circ}),\ol\dQ_\ell)\bigcap
\rH^{N-2}(\ol\rM^\dag_\fp(\rV^\circ,\rK^{p\circ})\cap\pi^{\circ-1}(s_0),\ol\dQ_\ell).
\]
Then $\rH^\prim_{s_0}(\ol\rM^\dag_\fp(\rV^\circ,\rK^{p\circ}),\ol\dQ_\ell)$ is a representation of $\rU(\Lambda^\circ_\fp)(\dF_p)=\rK^\circ_\fp/\rK^\circ_{\fp1}$, which is (isomorphic to) $\iota_\ell\Omega_N$. Thus, we obtain \eqref{eq:tate}.

For the remaining part, note that the right-hand side of \eqref{eq:tate} is a $\dC[\rK^{p\circ}\rK^\circ_{p1}\backslash\rU(\rV^\circ)(\dA^\infty_{F^+})/\rK^{p\circ}\rK^\circ_{p1}]$-submodule of $\Map(\rU(\rV^\circ)(F^+)\backslash\rU(\rV^\circ)(\dA^\infty_{F^+})/\rK^{p\circ}\rK^\circ_{p1},\dC)$. In particular, we can complete $\pi^{\infty,p}$ to an automorphic representation $\pi=\pi^{\infty,p}\otimes\pi_\infty\otimes\prod_{\fq\mid p}\pi_\fq$ of $\rU(\rV^\circ)(\dA_{F^+})$ such that $\pi_\infty$ is trivial; $\pi_\fq$ is unramified for $\fq\neq\fp$; and $\pi_\fp\res_{\rK^\circ_\fp}$ contains $\Omega_N$.

In case (1), by Proposition \ref{pr:omega}(2), we know that $\Omega_N$ has nonzero Borel fixed vectors. Thus, $\pi_\fp$ is a constituent of an unramified principal series.

In case (2), we first consider the case where $N=3$. As $\pi_\fp\res_{\rK^\circ_\fp}$ contains $\Omega_3$, it has to be $\cInd_{\rK_3}^{\rU_3}\Omega_3$ by Proposition \ref{pr:omega}(3) and \cite{MP96}*{Theorem~6.11(2)}. Thus, by \cite{MP96}*{Proposition~6.6}, $\pi_\fp\res_{\rK^\circ_\fp}$ is irreducible supercuspidal, which is actually the unique supercuspidal unipotent representation of $\rU(\rV^\circ)(F^+_\fp)$. In fact, $\cInd_{\rK_3}^{\rU_3}\Omega_3$ is the representation $\pi^s(\mathbf{1})$ appearing in \cite{Rog90}*{Proposition~13.1.3(d)}, after identifying $\ol\dQ_\ell$ with $\dC$. By \cite{Rog90}*{Proposition~13.2.2(c)}, $\BC(\pi^s(\mathbf{1}))$ is the tempered constituent of the unramified principal series of $\GL_3(F_\fp)$ with the Satake parameter $\{-p,1,-p^{-1}\}$. Now for general $N=2r+1$, as $\pi_\fp\res_{\rK^\circ_\fp}$ contains $\Omega_N$, by Proposition \ref{pr:omega}(4) and \cite{MP96}*{Theorem~6.11(2)}, $\pi_\fp$ is a constituent the normalized parabolic induction of $\pi^s(\mathbf{1})\boxtimes\chi_1\boxtimes\cdots\boxtimes\chi_{r-1}$ for some unramified characters $\chi_1,\dots,\chi_{r-1}$ of $F^\times$. Therefore, by the compatibility of local base change and induction, $\BC(\pi_\fp)$ is a constituent of an unramified principal series of $\GL_N(F_\fp)$ whose Satake parameter contains $\{-p,-p^{-1}\}$.

The proposition is proved.
\end{proof}

\subsection{Intersection on the ground stratum}
\label{ss:ns_intersection}

In this subsection, we describe a certain scheme-theoretical intersection on the ground stratum, which will be used in the next subsection. We assume $N\geq 2$.

Take an object $\rK^{p\circ}\in\fK(\rV^\circ)^p$. Given two (possibly same) points $s^\bullet_1,s^\bullet_2\in\rS^\bullet_\fp(\rV^\circ,\rK^{p\circ})(\kappa)$ for a perfect field $\kappa$ containing $\dF_p^\Phi$, we put
\[
\rB^\bullet_{s^\bullet_1,s^\bullet_2}\coloneqq
\rB^\bullet_{s^\bullet_1}\times_{\rM^\bullet_\fp(\rV^\circ,\rK^{p\circ})_\kappa}\rB^\bullet_{s^\bullet_2}
\]
as the (possibly empty) fiber product of $\iota^\bullet_\kappa\res\rB^\bullet_{s^\bullet_1}$ and $\iota^\bullet_\kappa\res\rB^\bullet_{s^\bullet_2}$. To describe $\rB^\bullet_{s^\bullet_1,s^\bullet_2}$, we need to use some particular cases of the Hecke correspondences introduced in Construction \ref{cs:ns_uniformization2}. We now give more details.

\begin{definition}\label{de:ns_definite_hecke}
For every integer $0\leq j\leq N$, we define a functor
\begin{align*}
\rS^\bullet_\fp(\rV^\circ,\obj)_j\colon\fK(\rV^\circ)^p\times\fT &\to\sfP\Sch'_{/\dF^\Phi_p} \\
\rK^{p\circ} &\mapsto \rS^\bullet_\fp(\rV^\circ,\rK^{p\circ})_j
\end{align*}
such that for every $S\in\Sch'_{/\dF^\Phi_p}$, $\rS^\bullet_\fp(\rV^\circ,\rK^{p\circ})_j(S)$ is the set of equivalence classes of decuples
$(A_0,\lambda_0,\eta_0^p;A^\bullet_1,\lambda^\bullet_1,\eta^{p\bullet}_1;A^\bullet_2,\lambda^\bullet_2,\eta^{p\bullet}_2;\phi^\bullet)$, where
\begin{itemize}[label={\ding{109}}]
  \item $(A_0,\lambda_0,\eta_0^p;A^\bullet_i,\lambda^\bullet_i,\eta^{p\bullet}_i)$ for $i=1,2$ are two elements in $\rS^\bullet_\fp(\rV^\circ,\rK^{p\circ})(S)$; and

  \item $\phi^\bullet\colon A^\bullet_1\to A^\bullet_2$ is an $O_F$-linear quasi-isogeny such that
     \begin{enumerate}[label=(\alph*)]
       \item $p\phi^\bullet\circ\lambda^{\bullet-1}_1$ is a quasi-$p$-isogeny; and $\Ker(p\phi^\bullet)[\fp]$ has rank $p^{2(N-j)}$;

       \item $\phi^\bullet[\fq^\infty]$ is an isomorphism for every prime $\fq$ of $F^+$ above $p$ that is not $\fp$;

       \item we have $\phi^{\bullet\vee}\circ\lambda^\bullet_2\circ\phi^\bullet=\lambda^\bullet_1$; and

       \item the $\rK^{p\circ}$-orbit of maps $v\mapsto\phi^\bullet_*\circ\eta^{p\bullet}_1(v)$ for $v\in\rV^\circ\otimes_\dQ\dA^{\infty,p}$ coincides with $\eta^{p\bullet}_2$.
     \end{enumerate}
\end{itemize}
The equivalence relation and the action of morphisms in $\fK(\rV^\circ)^p\times\fT$ are defined similarly as in Definition \ref{de:qs_basic_correspondence}. Finally, we denote
\begin{align*}
\Hk_j\colon\rS^\bullet_\fp(\rV^\circ,\obj)_j\to\rS^\bullet_\fp(\rV^\circ,\obj)\times\rS^\bullet_\fp(\rV^\circ,\obj)
\end{align*}
the morphism in $\Fun(\fK(\rV^\circ)^p\times\fT,\Sch_{/\dF_p^\Phi})_{/\rT_\fp}$ induced by the assignment
\[
(A_0,\lambda_0,\eta_0^p;A^\bullet_1,\lambda^\bullet_1,\eta^{p\bullet}_1;A^\bullet_2,\lambda^\bullet_2,\eta^{p\bullet}_2;\phi^\bullet)
\mapsto((A_0,\lambda_0,\eta_0^p;A^\bullet_1,\lambda^\bullet_1,\eta^{p\bullet}_1),
(A_0,\lambda_0,\eta_0^p;A^\bullet_2,\lambda^\bullet_2,\eta^{p\bullet}_2)).
\]
\end{definition}

\begin{remark}
When $\rK^{p\circ}$ is sufficiently small, the morphism
\[
\Hk_j\colon\rS^\bullet_\fp(\rV^\circ,\rK^{p\circ})_j\to\rS^\bullet_\fp(\rV^\circ,\rK^{p\circ})\times\rS^\bullet_\fp(\rV^\circ,\rK^{p\circ})
\]
is a closed immersion for every $j$; and the images of $\Hk_j$ for all $j$ are mutually disjoint.
\end{remark}

Now we take a point $s^\bullet=(A_0,\lambda_0,\eta_0^p;A^\bullet_1,\lambda^\bullet_1,\eta^{p\bullet}_1;A^\bullet_2,\lambda^\bullet_2,\eta^{p\bullet}_2;\phi^\bullet)
\in\rS^\bullet_\fp(\rV^\circ,\rK^{p\circ})_j(\kappa)$ where $\kappa$ is a perfect field containing $\dF^\Phi_p$. By Definition \ref{de:ns_definite_hecke}(c), we have $(p\phi^\bullet\circ\lambda^{\bullet-1}_1)^\vee=p\phi^{\bullet-1}\circ\lambda^{\bullet-1}_2$. Thus, $p\phi^{\bullet-1}\circ\lambda^{\bullet-1}_2$, hence $p\phi^{\bullet-1}$ are quasi-$p$-isogenies as well. In particular, for every $\tau\in\Sigma_\infty$, we may consider
\begin{align*}
\Ker(p\phi^\bullet)_{*,\tau}&\coloneqq
\Ker\((p\phi^\bullet)_{*,\tau}\colon\rH^\dr_1(A^\bullet_1/\kappa)_\tau\to\rH^\dr_1(A^\bullet_2/\kappa)_\tau\),\\
\IM(p\phi^{\bullet-1})_{*,\tau}&\coloneqq
\IM\((p\phi^{\bullet-1})_{*,\tau}\colon\rH^\dr_1(A^\bullet_2/\kappa)_\tau\to\rH^\dr_1(A^\bullet_1/\kappa)_\tau\).
\end{align*}

\begin{lem}\label{le:ns_definite_hecke1}
We have
\begin{enumerate}
  \item $\IM(p\phi^{\bullet-1})_{*,\tau}\subseteq\Ker(p\phi^\bullet)_{*,\tau}$ for every $\tau\in\Sigma_\infty$;

  \item $\dim_\kappa\Ker(p\phi^\bullet)_{*,\tau}=N-j$ for $\tau\in\{\tau_\infty,\tau_\infty^\tc\}$;

  \item $\IM(p\phi^{\bullet-1})_{*,\tau}\cap\rH^\dr_1(A^\bullet_1/\kappa)_{\tau^\tc}^\perp=0$ for $\tau\in\{\tau_\infty,\tau_\infty^\tc\}$;

  \item $(\IM(p\phi^{\bullet-1})_{*,\tau})^\perp=\Ker(p\phi^\bullet)_{*,\tau^\tc}$ for $\tau\in\{\tau_\infty,\tau_\infty^\tc\}$; and

  \item $\dim_\kappa\IM(p\phi^{\bullet-1})_{*,\tau}=j$ for $\tau\in\{\tau_\infty,\tau_\infty^\tc\}$.
\end{enumerate}
In particular, $\rS^\bullet_\fp(\rV^\circ,\rK^{p\circ})_j$ is empty if $j>\floor{\tfrac{N}{2}}$.
\end{lem}

\begin{proof}
For (1), it is obvious since $(p\phi^\bullet)\circ(p\phi^{\bullet-1})=p^2$.

For (2), by Definition \ref{de:ns_definite_hecke}(a), we have $\dim_\kappa\Ker(p\phi^\bullet)_{*,\tau_\infty}+\dim_\kappa\Ker(p\phi^\bullet)_{*,\tau_\infty^\tc}=2(N-j)$. Using the isomorphisms $\tV\colon\rH^\dr_1(A^\bullet_1/\kappa)_{\tau_\infty}\to\rH^\dr_1(A^\bullet_1/\kappa)_{\tau_\infty^\tc}$ and $\tV\colon\rH^\dr_1(A^\bullet_2/\kappa)_{\tau_\infty}\to\rH^\dr_1(A^\bullet_2/\kappa)_{\tau_\infty^\tc}$, we have $\dim_\kappa\Ker(p\phi^\bullet)_{*,\tau_\infty}=\dim_\kappa\Ker(p\phi^\bullet)_{*,\tau_\infty^\tc}$, hence both are equal to $N-j$.

For (3), it suffices to consider $\tau=\tau_\infty$ due to the isomorphism $\tV$. Via $\phi^\bullet$, we regard $\cD(A^\bullet_2)$ as a lattice in $\cD(A^\bullet_1)_\dQ$. By Definition \ref{de:ns_definite_hecke}(a), we have $p\cD(A^\bullet_2)_{\tau_\infty}\subseteq\cD(A^\bullet_1)_{\tau_\infty}\subseteq \cD(A^\bullet_2)_{\tau_\infty^\tc}^\vee$ (Notation \ref{no:weil_pairing_dieudonne}). Suppose that $\rH^\dr_1(A^\bullet_1/\kappa)_{\tau_\infty^\tc}^\perp\cap\IM(p\phi^{\bullet-1})_{*,\tau_\infty}\neq 0$. Then one can find $x_2\in\cD(A^\bullet_2)_{\tau_\infty}$ and $x_1\in\cD(A^\bullet_1)_{\tau_\infty^\tc}^\vee\setminus\cD(A^\bullet_1)_{\tau_\infty}$ such that $px_1=px_2$. It follows that $\langle x_2,\tV x_2\rangle_{\lambda^\bullet_2,\tau_\infty}=\langle x_1,\tV x_1\rangle_{\lambda^\bullet_1,\tau_\infty}$ does not belong to $W(\kappa)$, which is a contradiction. Here, we regard $\tV$ as Verschiebung maps on for Dieudonn\'{e} modules of $A^\bullet_1$ and $A^\bullet_2$, which are isomorphisms.

For (4), as $\lambda^\bullet_1\circ\phi^{\bullet-1}=\phi^{\bullet\vee}\circ\lambda^\bullet_2$, we have for $\tau\in\{\tau_\infty,\tau_\infty^\tc\}$ that
\[
(\IM(p\phi^{\bullet-1})_{*,\tau})^\perp=((p\phi^\bullet)_{*,\tau^\tc})^{-1}\rH^\dr_1(A^\bullet_2/\kappa)_\tau^\perp,
\]
which equals $\Ker(p\phi^\bullet)_{*,\tau^\tc}$ by (3).

For (5), by (2,3,4), we have $\dim_\kappa\IM(p\phi^{\bullet-1})_{*,\tau}=j$ for $\tau\in\{\tau_\infty,\tau_\infty^\tc\}$.

The last claim follows from (1,2,5).
\end{proof}

By Lemma \ref{le:ns_definite_hecke1}(1,4), for $\tau\in\{\tau_\infty,\tau_\infty^\tc\}$, we may put
\[
\rH^\dr_1(\phi^\bullet)_\tau\coloneqq
\frac{\Ker(p\phi^\bullet)_{*,\tau}}{\IM(p\phi^{\bullet-1})_{*,\tau}};
\]
and we have the induced $\kappa$-bilinear pairing
\[
\langle\;,\;\rangle_{\lambda^\bullet_1,\tau_\infty}\colon\rH^\dr_1(\phi^\bullet)_{\tau_\infty}
\times\rH^\dr_1(\phi^\bullet)_{\tau_\infty^\tc}\to\kappa.
\]
On the other hand, the $(\kappa,\sigma^{-1})$-linear Verschiebung map $\tV\colon\rH^\dr_1(A^\bullet_1/\kappa)_{\tau_\infty}\to\rH^\dr_1(A^\bullet_1/\kappa)_{\tau_\infty^\tc}$ induces a $(\kappa,\sigma^{-1})$-linear isomorphism $\tV\colon\rH^\dr_1(\phi^\bullet)_{\tau_\infty}\to\rH^\dr_1(\phi^\bullet)_{\tau_\infty^\tc}$. We define a pairing
\[
\{\;,\;\}_{s^\bullet}\colon\rH^\dr_1(\phi^\bullet)_{\tau_\infty^\tc}\times\rH^\dr_1(\phi^\bullet)_{\tau_\infty^\tc}\to\kappa
\]
by the formula $\{x,y\}_{s^\bullet}\coloneqq\langle\tV^{-1}x,y\rangle_{\lambda^\bullet_1,\tau_\infty}$. To ease notation, we put
\[
\sV_{s^\bullet}\coloneqq\rH^\dr_1(\phi^\bullet)_{\tau_\infty^\tc}.
\]

\begin{lem}\label{le:ns_definite_hecke2}
Suppose that $j\leq\floor{\tfrac{N}{2}}-1$. The pair $(\sV_{s^\bullet},\{\;,\;\}_{s^\bullet})$ is admissible of rank $N-2j$ (Definition \ref{de:dl_admissible}) satisfying $\dim_\kappa\sV_{s^\bullet}^\dashv=N-2\floor{\tfrac{N}{2}}$. In particular, we have the geometrically irreducible smooth projective scheme $\DL^\bullet(\sV_{s^\bullet},\{\;,\;\}_{s^\bullet})\in\Sch_{/\kappa}$ of dimension $\floor{\tfrac{N}{2}}-j$ as introduced in Definition \ref{de:dl_bullet}.
\end{lem}

\begin{proof}
By Lemma \ref{le:ns_definite_hecke1}(2,5), we have $\dim_\kappa\sV_{s^\bullet}=N-2j$. By Lemma \ref{le:ns_definite_hecke1}(3,4), we have $\dim_\kappa\sV_{s^\bullet}^\dashv=N-2\floor{\tfrac{N}{2}}$. The lemma follows by Proposition \ref{pr:dl_bullet}.
\end{proof}

Now consider a connected scheme $S\in\Sch'_{/\kappa}$ and a point $x\in\rB^\bullet_{s^\bullet_1,s^\bullet_2}(S)$ represented by a quattuordecuple $(A_0,\lambda_0,\eta_0^p;A,\lambda,\eta^p;A^\bullet_1,\lambda^\bullet_1,\eta^{p\bullet}_1;\gamma_1;
A^\bullet_2,\lambda^\bullet_2,\eta^{p\bullet}_2;\gamma_2)$.

\begin{lem}\label{le:ns_definite_hecke3}
There exists a unique integer $j$ satisfying $0\leq j\leq \floor{\tfrac{N}{2}}-1$ such that $s^\bullet\coloneqq(A_0,\lambda_0,\eta_0^p;
A^\bullet_1,\lambda^\bullet_1,\eta^{p\bullet}_1;A^\bullet_2,\lambda^\bullet_2,\eta^{p\bullet}_2;\phi^\bullet)$ is an element in $\rS^\bullet_\fp(\rV^\circ,\rK^{p\circ})_j(S)$, where $\phi^\bullet\coloneqq\gamma_2\circ\gamma_1^{-1}$. Moreover, we have
\begin{align}\label{eq:ns_definite_hecke1}
\IM(p\phi^{\bullet-1})_{*,\tau_\infty^\tc}\subseteq H_2\subseteq H_1 \subseteq\Ker(p\phi^\bullet)_{*,\tau_\infty^\tc},
\end{align}
where $H_2\subseteq H_1\subseteq \rH^\dr_1(A_1^\bullet/S)_{\tau_\infty^\tc}$ are subbundles in Theorem \ref{th:ns_basic_correspondence2} for the image of $x$ in $\rB^\bullet_{s^\bullet_1}(S)$.
\end{lem}

\begin{proof}
First, by definition, we have $\Ker(p\phi^\bullet)[\fp]=\Ker(\gamma_2\circ\breve\gamma_1)[\fp]$, which is an $O_F$-stable finite flat subgroup of $A^\bullet_1[\fp]$. Thus, as $S$ is connected, there is a unique integer $j$ satisfying $0\leq j\leq N$ such that $\Ker(p\phi^\bullet)[\fp]$ has rank $p^{2(N-j)}$.

Second, we show that $p\phi^\bullet\circ\lambda^{\bullet-1}_1$ is a quasi-$p$-isogeny, that is, $\gamma_2\circ\breve\gamma_1\circ\lambda^{\bullet-1}_1$ is a quasi-$p$-isogeny. By Theorem \ref{th:ns_basic_correspondence2}(4), $\gamma_{1*,\tau_\infty^\tc}\omega_{A^\vee/S,\tau_\infty^\tc}$ contains $\rH^\dr_1(A^\bullet_1)_{\tau_\infty}^\perp$, which implies $\breve\gamma_{1*,\tau_\infty^\tc}\rH^\dr_1(A^\bullet_1)_{\tau_\infty}^\perp=0$ hence $(\gamma_2\circ\breve\gamma_1)_{*,\tau_\infty^\tc}\rH^\dr_1(A^\bullet_1)_{\tau_\infty}^\perp=0$. On the other hand, as $\breve\gamma_{1*,\tau_\infty}\rH^\dr_1(A^\bullet_1)_{\tau_\infty^\tc}^\perp\subseteq\rH^\dr_1(A)_{\tau_\infty^\tc}^\perp$, we have $(\gamma_2\circ\breve\gamma_1)_{*,\tau_\infty}\rH^\dr_1(A^\bullet_1)_{\tau_\infty^\tc}^\perp=0$ by Definition \ref{de:ns_basic_correspondence2}(c). In other words, $\Ker\lambda^\bullet_1[\fp^\infty]$ is contained in $\Ker\breve\gamma_1[\fp^\infty]$. Thus, $p\phi^\bullet\circ\lambda^{\bullet-1}_1=\gamma_2\circ\breve\gamma_1\circ\lambda^{\bullet-1}_1$ a quasi-$p$-isogeny.

Third, we show that $j$ is at most $\floor{\tfrac{N}{2}}-1$. (Note that Lemma \ref{le:ns_definite_hecke1} already implies that $j\leq\floor{\tfrac{N}{2}}$.) Theorem \ref{th:ns_basic_correspondence2}(4) implies $\rank_{\cO_S}H_2+1=\rank_{\cO_S}H_1$ and $\rH^\dr_1(A^\bullet_1/S)_{\tau_\infty}^\perp\subseteq H_2$. Lemma \ref{le:ns_definite_hecke1}(3) implies $\rank_{\cO_S}H_2\geq\rank_{\cO_S}\IM(p\phi^{\bullet-1})_{*,\tau_\infty^\tc}+1$. Thus, by Lemma \ref{le:ns_definite_hecke1}(2,5) and \eqref{eq:ns_definite_hecke1}, we have $(N-j)-j\geq 2$, that is, $j\leq\floor{\tfrac{N}{2}}-1$.

Definition \ref{de:ns_definite_hecke}(b,c,d) are obvious. Thus, it remains to check \eqref{eq:ns_definite_hecke1}. On one hand, we have
\begin{align*}
\IM(p\phi^{\bullet-1})_{*,\tau_\infty^\tc}&=\IM(\gamma_1\circ\breve\gamma_2)_{*,\tau_\infty^\tc}
=\gamma_{1*,\tau_\infty^\tc}\breve\gamma_{2*,\tau_\infty^\tc}\rH^\dr_1(A^\bullet_2/S)_{\tau_\infty^\tc}\\
&=\gamma_{1*,\tau_\infty^\tc}\breve\gamma_{2*,\tau_\infty^\tc}\omega_{A^{\bullet\vee}_2/S,\tau_\infty^\tc}
\subseteq\gamma_{1*,\tau_\infty^\tc}\omega_{A^{\bullet\vee}_1/S,\tau_\infty^\tc}=H_2.
\end{align*}
On the other hand, since $\breve\gamma_{1*,\tau_\infty}\IM(p\phi^{\bullet-1})_{*,\tau_\infty}=
\breve\gamma_{1*,\tau_\infty}\IM(\gamma_1\circ\breve\gamma_2)_{*,\tau_\infty}=0$, we have the inclusion $\IM(p\phi^{\bullet-1})_{*,\tau_\infty}\subseteq(\breve\gamma_{1*,\tau_\infty})^{-1}\omega_{A^\vee,\tau_\infty}$. Thus, $H_1=((\breve\gamma_{1*,\tau_\infty})^{-1}\omega_{A^\vee,\tau_\infty})^\perp$ is contained in $(\IM(p\phi^{\bullet-1})_{*,\tau_\infty})^\perp$, which is $\Ker(p\phi^\bullet)_{*,\tau_\infty^\tc}$ by Lemma \ref{le:ns_definite_hecke1}(4). The lemma is proved.
\end{proof}

\begin{definition}\label{de:ns_definite_hecke_1}
By Lemma \ref{le:ns_definite_hecke3}, we have a morphism
\[
\rB^\bullet_{s^\bullet_1,s^\bullet_2}\to\coprod_{j=0}^{\floor{\tfrac{N}{2}}-1}\Hk_j^{-1}(s^\bullet_1,s^\bullet_2).
\]
For a point $s^\bullet\in\Hk_j^{-1}(s^\bullet_1,s^\bullet_2)(\kappa)$ for some $0\leq j\leq\floor{\tfrac{N}{2}}-1$, we denote by $\rB^\bullet_{s^\bullet}$ the inverse image under the above morphism, which is an open and closed subscheme of $\rB^\bullet_{s^\bullet_1,s^\bullet_2}$.
\end{definition}

\begin{theorem}\label{th:ns_incidence}
Let $s^\bullet_1,s^\bullet_2\in\rS^\bullet_\fp(\rV^\circ,\rK^{p\circ})(\kappa)$ be two points for a perfect field $\kappa$ containing $\dF_p^\Phi$. We have
\[
\rB^\bullet_{s^\bullet_1,s^\bullet_2}=\coprod_{j=0}^{\floor{\tfrac{N}{2}}-1}\coprod_{s^\bullet\in\Hk_j^{-1}(s^\bullet_1,s^\bullet_2)(\kappa)}
\rB^\bullet_{s^\bullet}.
\]
Take $s^\bullet=(A_0,\lambda_0,\eta_0^p;A^\bullet_1,\lambda^\bullet_1,\eta^{p\bullet}_1;A^\bullet_2,\lambda^\bullet_2,\eta^{p\bullet}_2;\phi^\bullet)
\in\Hk_j^{-1}(s^\bullet_1,s^\bullet_2)(\kappa)$ for some $0\leq j\leq\floor{\tfrac{N}{2}}-1$.
\begin{enumerate}
  \item Denote by $\bar{H}_i$ the image of $H_i$ in $\rH^\dr_1(\phi^\bullet)_{\tau_\infty^\tc}\otimes_\kappa\cO_S=(\sV_{s^\bullet})_S$ for $i=1,2$. Then the assignment sending $(A_0,\lambda_0,\eta_0^p;A,\lambda,\eta^p;A^\bullet_1,\lambda^\bullet_1,\eta^{p\bullet}_1;\gamma_1; A^\bullet_2,\lambda^\bullet_2,\eta^{p\bullet}_2;\gamma_2)\in\rB^\bullet_{s^\bullet}(S)$ to $(\bar{H}_1,\bar{H_2})$ induces an isomorphism
      \[
      \zeta^\bullet_{s^\bullet}\colon\rB^\bullet_{s^\bullet}\to
      \DL^\bullet(\sV_{s^\bullet},\{\;,\;\}_{s^\bullet})
      \]
      (Definition \ref{de:dl_bullet}) in $\Sch_{/\kappa}$.

  \item The cokernel of the map
      \[
      \cT_{\rB^\bullet_{s^\bullet_1}/\kappa}\res_{\rB^\bullet_{s^\bullet}}
      \bigoplus\cT_{\rB^\bullet_{s^\bullet_2}/\kappa}\res_{\rB^\bullet_{s^\bullet}}
      \to\iota^{\bullet*}\cT_{\rM^\bullet_\fp(\rV^\circ,\rK^{p\circ})/\kappa}\res_{\rB^\bullet_{s^\bullet}}
      \]
      is canonically isomorphic to
      \[
      \zeta^{\bullet*}_{s^\bullet}\(\(\sigma^*\bar\cH_{s^\bullet2}\)
      \otimes_{\cO_{\DL^\bullet(\sV_{s^\bullet},\{\;,\;\}_{s^\bullet})}}\(\bar\cH_{s^\bullet1}^\dashv/\bar\cH_{s^\bullet2}\)\)
      \]
      where $(\bar\cH_{s^\bullet1},\bar\cH_{s^\bullet2})$ is the universal object over $\DL^\bullet(\sV_{s^\bullet},\{\;,\;\}_{s^\bullet})$.
\end{enumerate}
\end{theorem}

\begin{proof}
The decomposition of $\rB^\bullet_{s^\bullet_1,s^\bullet_2}$ follows directly from the definition and the fact that $\Hk_j^{-1}(s^\bullet_1,s^\bullet_2)$ is isomorphic to a finite disjoint union of $\Spec\kappa$.

First, we show (1). We first notice that Lemma \ref{le:ns_definite_hecke1} implies that $(\bar{H}_1,\bar{H}_2)$ is an element in $\DL^\bullet(\sV_{s^\bullet},\{\;,\;\}_{s^\bullet})(S)$.

Since the target of $\zeta^\bullet_{s^\bullet}$ is smooth over $\kappa$ by Lemma \ref{le:ns_definite_hecke2}, to see that $\zeta^\bullet_{s^\bullet}$ is an isomorphism, it suffices to check that for every algebraically closed field $\kappa'$ containing $\kappa$
\begin{enumerate}
  \item[(1--1)] $\zeta^\bullet_{s^\bullet}$ induces a bijection on $\kappa'$-points; and

  \item[(1--2)] $\zeta^\bullet_{s^\bullet}$ induces an isomorphism on the tangent spaces at every $\kappa'$-point.
\end{enumerate}
To ease notation, we may assume $\kappa'=\kappa$.

For (1--1), we construct an inverse to the map $\zeta^\bullet_{s^\bullet}(\kappa)$. Take a point $y\in\DL^\bullet(\sV_{s^\bullet},\{\;,\;\}_{s^\bullet})(\kappa)$ represented by $\kappa$-linear subspaces $\sV_{s^\bullet}^\dashv\subseteq\bar{H}_2\subseteq\bar{H}_1\subseteq\sV_{s^\bullet}$, or equivalently, subspaces
\[
\IM(p\phi^{\bullet-1})_{*,\tau_\infty^\tc}\oplus\rH^\dr_1(A^\bullet_1/\kappa)_{\tau_\infty}^\perp\subseteq H_2\subseteq H_1 \subseteq\Ker(p\phi^\bullet)_{*,\tau_\infty^\tc}\subseteq\rH^\dr_1(A^\bullet_1/\kappa)_{\tau_\infty^\tc}.
\]
These give rise to a point $y_1\in\DL^\bullet(\sV_{s^\bullet_1},\{\;,\;\}_{s^\bullet_1})(\kappa)$. By Theorem \ref{th:ns_basic_correspondence2}(4), we obtain a unique point $x_1=(A_0,\lambda_0,\eta_0^p;A,\lambda,\eta^p;A^\bullet_1,\lambda^\bullet_1,\eta^{p\bullet}_1;\gamma_1)\in\rB^\bullet_{s^\bullet_1}(\kappa)$ such that $\zeta^\bullet_{s^\bullet_1}(x_1)=y_1$. Put $\gamma_2\coloneqq\phi^\bullet\circ\gamma_1\colon A\to A^\bullet_2$. We claim that $\gamma_2$ is a quasi-$p$-isogeny. In fact, as $\lambda\circ\breve\gamma_1=\gamma_1^\vee\circ\lambda^\bullet_1$, $\langle\IM\gamma_{1*,\tau_\infty},\IM\gamma_{1*,\tau_\infty^\tc}\rangle_{\lambda^\bullet_1,\tau_\infty}=0$. Thus, we have
\[
\IM\gamma_{1*,\tau_\infty^\tc}\subseteq(\IM\gamma_{1*,\tau_\infty})^\perp
=(\tV^{-1}\gamma_{1*,\tau_\infty^\tc}\omega_{A^\vee,\tau_\infty^\tc})^\perp=H_2^\dashv\subseteq\Ker(p\phi^\bullet)_{*,\tau_\infty^\tc}.
\]
By the isomorphisms $\tV\colon\rH^\dr_1(A^\bullet_1/\kappa)_{\tau_\infty}\to\rH^\dr_1(A^\bullet_1/\kappa)_{\tau_\infty^\tc}$ and $\tV\colon\rH^\dr_1(A^\bullet_2/\kappa)_{\tau_\infty}\to\rH^\dr_1(A^\bullet_2/\kappa)_{\tau_\infty^\tc}$, we obtain $\IM\gamma_{1*,\tau_\infty}\subseteq\Ker(p\phi^\bullet)_{*,\tau_\infty}$. In particular, $\IM(p\phi^\bullet\circ\gamma_1)_{*,\tau}=0$ for every $\tau\in\Sigma_\infty$; in other words, $\gamma_2$ is a quasi-$p$-isogeny. Now we show that $x_2\coloneqq(A_0,\lambda_0,\eta_0^p;A,\lambda,\eta^p;A^\bullet_2,\lambda^\bullet_2,\eta^{p\bullet}_2;\gamma_2)$ satisfies Definition \ref{de:ns_basic_correspondence2}(a--e).

For (a), it suffices to show that $p\gamma_2^{-1}$ is a quasi-$p$-isogeny, equivalently, $\gamma_1^{-1}\circ(p\phi^{\bullet-1})$ is a quasi-$p$-isogeny. However, we have $\IM(p\phi^{\bullet-1})_{*,\tau_\infty}=\tV^{-1}\IM(p\phi^{\bullet-1})_{*,\tau_\infty^\tc}\subseteq\tV^{-1}H_2=\IM\gamma_{1*,\tau_\infty}$, hence $\IM(p\phi^{\bullet-1})_{*,\tau_\infty^\tc}\subseteq\IM\gamma_{1*,\tau_\infty^\tc}$ using the action of $\tV$, which together imply that $\gamma_1^{-1}\circ(p\phi^{\bullet-1})$ is a quasi-$p$-isogeny.

For (b), we identify $\cD(A)$ as submodules of both $\cD(A^\bullet_1)$ and $\cD(A^\bullet_2)$ via $\gamma_1$ and $\gamma_2$, respectively. Then we need to show that $p\cD(A^\bullet_2)_{\tau_\infty}^\vee\cap\cD(A)_{\tau_\infty^\tc}\subseteq\tV\cD(A)_{\tau_\infty}$. As $p\phi^{\bullet-1}\circ\lambda^{\bullet-1}_2$ is a quasi-$p$-isogeny, we have $p\cD(A^\bullet_2)_{\tau_\infty}^\vee\subseteq\cD(A^\bullet_1)_{\tau_\infty^\tc}$. Moreover, the image of $p\cD(A^\bullet_2)_{\tau_\infty}^\vee$ in $\cD(A^\bullet_1)_{\tau_\infty^\tc}/p\cD(A^\bullet_1)_{\tau_\infty^\tc}=\rH^\dr_1(A^\bullet_1)_{\tau_\infty^\tc}$ is contained in $\IM(p\phi^{\bullet-1})_{*,\tau_\infty^\tc}\oplus\rH^\dr_1(A^\bullet_1/\kappa)_{\tau_\infty}^\perp$, which is further contained in $H_2$. Thus, $p\cD(A^\bullet_2)_{\tau_\infty}^\vee\cap\cD(A)_{\tau_\infty^\tc}\subseteq\tV\cD(A)_{\tau_\infty}$ as $\tV\cD(A)_{\tau_\infty}$ is the inverse image of $H_2$ in $\cD(A^\bullet_1)_{\tau_\infty^\tc}$.

For (c), suppose that $\rH^\dr_1(A)_{\tau_\infty^\tc}^\perp$ is not contained in $\Ker\gamma_{2*,\tau_\infty}$. Since $\gamma_{2*,\tau_\infty}$ maps $\rH^\dr_1(A)_{\tau_\infty^\tc}^\perp$ into $\rH^\dr_1(A^\bullet_2)_{\tau_\infty^\tc}^\perp$, we have $\gamma_{2*,\tau_\infty}\rH^\dr_1(A)_{\tau_\infty^\tc}^\perp\cap\rH^\dr_1(A^\bullet_2)_{\tau_\infty^\tc}^\perp\neq 0$. On the other hand, since $\rH^\dr_1(A)_{\tau_\infty^\tc}^\perp$ is contained in $\Ker\gamma_{1*,\tau_\infty}=\IM\breve\gamma_{1*,\tau_\infty}$, we have $\gamma_{2*,\tau_\infty}\rH^\dr_1(A)_{\tau_\infty^\tc}^\perp\subseteq\IM(\gamma_2\circ\breve\gamma_1)_{*,\tau_\infty}
=\IM(p\phi^\bullet)_{*,\tau_\infty}$. Thus, $\IM(\gamma_2\circ\breve\gamma_1)_{*,\tau_\infty}\cap\rH^\dr_1(A^\bullet_2)_{\tau_\infty^\tc}^\perp\neq 0$, which contradicts with \ref{le:ns_definite_hecke1}(3) (with $\phi^\bullet$ replaced by $\phi^{\bullet-1}$).

For (d) and (e), they follow obviously.

To summarize, $x_2$ belongs to $\rB^\bullet_{s^\bullet_2}(\kappa)$; and $x\coloneqq(x_1,x_2)$ is an element in $\rB^\bullet_{s^\bullet}(\kappa)$ such that $\zeta^\bullet_{s^\bullet}(x)=y$. It is easy to see that such assignment gives rise to an inverse of $\zeta^\bullet_{s^\bullet}(\kappa)$, hence (1--1) follows immediately.

For (1--2), let $\cT_x$ and $\cT_y$ be the tangent spaces at $x$ and $y$ as in (1--1), respectively. By Theorem \ref{th:ns_basic_correspondence2}(1), we have a canonical short exact sequence
\[
0 \to \Hom_\kappa\(\omega_{A^\vee,\tau_\infty},\frac{\omega_{A^\vee,\tau_\infty^\tc}^\perp}{\omega_{A^\vee,\tau_\infty}}\)
\to \cT_x\to\Hom_\kappa\(\frac{\omega_{A^\vee,\tau_\infty^\tc}}
{(\Ker\gamma_{1*,\tau_\infty})^\perp+(\Ker\gamma_{2*,\tau_\infty})^\perp},\Lie_{A^\vee,\tau_\infty^\tc}\)  \to 0.
\]
Then by Proposition \ref{pr:dl_bullet} and the construction, the induced map $(\zeta^\bullet_{s^\bullet})_*\colon\cT_x\to\cT_y$ fits into a commutative diagram
\[
\xymatrix{
\Hom_\kappa\(\omega_{A^\vee,\tau_\infty},\dfrac{\omega_{A^\vee,\tau_\infty^\tc}^\perp}{\omega_{A^\vee,\tau_\infty}}\) \ar[r]\ar[d]
& \cT_x \ar[r]\ar[d]^-{(\zeta^\bullet_{s^\bullet})_*} & \Hom_\kappa\(\dfrac{\omega_{A^\vee,\tau_\infty^\tc}}
{(\Ker\gamma_{1*,\tau_\infty})^\perp+(\Ker\gamma_{2*,\tau_\infty})^\perp},\Lie_{A^\vee,\tau_\infty^\tc}\) \ar[d] \\
\Hom_\kappa\(\bar{H}_1/\bar{H}_2,\bar{H}_2^\dashv/\bar{H}_1\)  \ar[r]& \cT_y  \ar[r]& \Hom_\kappa(\bar{H}_2/\sV_{s^\bullet}^\dashv,\bar{H}_1^\dashv/\bar{H}_2)
}
\]
in $\Mod(\kappa)$. The left vertical arrow is the composition
\begin{align*}
&\quad\Hom_\kappa\(\omega_{A^\vee,\tau_\infty},\omega_{A^\vee,\tau_\infty^\tc}^\perp/\omega_{A^\vee,\tau_\infty}\) \\
&\to\Hom_\kappa\(H_1^\perp/\tV^{-1} H_2,H_2^\perp/H_1^\perp\)\\
&\xrightarrow{\sim}\Hom_\kappa\(H_1/H_2,H_2^\dashv/H_1\)
\simeq\Hom_\kappa\(\bar{H}_1/\bar{H}_2,\bar{H}_2^\dashv/\bar{H}_1\),
\end{align*}
which is an isomorphism. The right vertical arrow is induced by maps
\begin{align}\label{eq:ns_incidence3}
\frac{\omega_{A^\vee,\tau_\infty^\tc}}{(\Ker\gamma_{1*,\tau_\infty})^\perp+(\Ker\gamma_{2*,\tau_\infty})^\perp}
&\xrightarrow{\gamma_{1*,\tau_\infty^\tc}}
\frac{H_2}{\IM(p\phi^{\bullet-1})_{*,\tau_\infty^\tc}\oplus\rH^\dr_1(A^\bullet_1/\kappa)_{\tau_\infty}^\perp}
\simeq\bar{H}_2/\sV_{s^\bullet}^\dashv,
\end{align}
\begin{align}\label{eq:ns_incidence4}
\Lie_{A^\vee,\tau_\infty^\tc}\simeq\rH^\dr_1(A)_{\tau_\infty^\tc}/\omega_{A^\vee,\tau_\infty^\tc}
&\xrightarrow{\gamma_{1*,\tau_\infty^\tc}}H_1^\dashv/H_2\simeq\bar{H}_1^\dashv/\bar{H}_2.
\end{align}
Note that in \eqref{eq:ns_incidence3}, we have used Lemma \ref{le:ns_definite_hecke1}(3) to write the direct sum.

We show that \eqref{eq:ns_incidence3} is well-defined and is an isomorphism. It is clear that $\Ker\gamma_{1*,\tau_\infty^\tc}$ is contained in $(\Ker\gamma_{1*,\tau_\infty})^\perp$. Thus, it suffices to show that the image of $(\Ker\gamma_{1*,\tau_\infty})^\perp+(\Ker\gamma_{2*,\tau_\infty})^\perp$ under $\gamma_{1*,\tau_\infty^\tc}$ is $\IM(p\phi^{\bullet-1})_{*,\tau_\infty^\tc}\oplus\rH^\dr_1(A^\bullet_1/\kappa)_{\tau_\infty}^\perp$. By Theorem \ref{th:ns_basic_correspondence2}(3), we have $\gamma_{1*,\tau_\infty^\tc}(\Ker\gamma_{1*,\tau_\infty})^\perp=\rH^\dr_1(A^\bullet_1/\kappa)_{\tau_\infty}^\perp$. It is easy to see that $\gamma_{1*,\tau_\infty^\tc}(\Ker\gamma_{2*,\tau_\infty})^\perp$ is contained in $\Ker(\gamma_2\circ\breve\gamma_1)_{*,\tau_\infty}^\perp=\Ker(p\phi^\bullet)_{*,\tau_\infty}^\perp$, which coincides with $\IM(p\phi^{\bullet-1})_{*,\tau_\infty^\tc}\oplus\rH^\dr_1(A^\bullet_1/\kappa)_{\tau_\infty}^\perp$ by Lemma \ref{le:ns_definite_hecke1}(3,4). On the other hand, $\gamma_{1*,\tau_\infty^\tc}(\Ker\gamma_{2*,\tau_\infty})^\perp$ contains $\gamma_{1*,\tau_\infty^\tc}(\Ker\gamma_{2*,\tau_\infty^\tc})=\IM(\gamma_1\circ\breve\gamma_2)_{*,\tau_\infty^\tc}$, which is $\IM(p\phi^{\bullet-1})_{*,\tau_\infty^\tc}$. It follows that \eqref{eq:ns_incidence3} is an isomorphism.

By Theorem \ref{th:ns_basic_correspondence2}(5), \eqref{eq:ns_incidence4} is an isomorphism as well. Thus, $(\zeta^\bullet_{s^\bullet})_*\colon\cT_x\to\cT_y$ is an isomorphism by the Five Lemma, hence (1--2) and (1) follow.

Next, we show (2). Theorem \ref{th:ns_basic_correspondence2}(2) implies that the cokernel of the map
\[
\cT_{\rB^\bullet_{s^\bullet_1}/\kappa}\res_{\rB^\bullet_{s^\bullet}}
\bigoplus\cT_{\rB^\bullet_{s^\bullet_2}/\kappa}\res_{\rB^\bullet_{s^\bullet}}
\to\iota^{\bullet*}\cT_{\rM^\bullet_\fp(\rV^\circ,\rK^{p\circ})/\kappa}\res_{\rB^\bullet_{s^\bullet}}
\]
is canonically isomorphic to
\begin{align}\label{eq:ns_incidence1}
\HOM\((\Ker\gamma_{1*,\tau_\infty}+\Ker\gamma_{2*,\tau_\infty})^\perp/\rH^\dr_1(\cA)_{\tau_\infty}^\perp,\Lie_{\cA^\vee,\tau_\infty^\tc}\).
\end{align}
As $\Ker\gamma_{2*,\tau_\infty}=\IM\breve\gamma_{2*,\tau_\infty}$, we have
\begin{align}\label{eq:ns_incidence2}
\frac{\rH^\dr_1(\cA)_{\tau_\infty}}{\Ker\gamma_{1*,\tau_\infty}+\Ker\gamma_{2*,\tau_\infty}}
\simeq\frac{\IM\gamma_{1*,\tau_\infty}}{\IM(\gamma_1\circ\breve\gamma_2)_{*,\tau_\infty}}
=\frac{\IM\gamma_{1*,\tau_\infty}}{\IM(p\phi^{\bullet-1})_{*,\tau_\infty}}
\simeq\frac{\tV\IM\gamma_{1*,\tau_\infty}}{\tV\IM(p\phi^{\bullet-1})_{*,\tau_\infty}}.
\end{align}
However, we have $\tV\IM\gamma_{1*,\tau_\infty}=(\gamma_{1*,\tau_\infty^\tc}\omega_{\cA,\tau_\infty^\tc})^{(p)}$ and $\tV\IM(p\phi^{\bullet-1})_{*,\tau_\infty}=(\IM(p\phi^{\bullet-1})_{*,\tau_\infty^\tc})^{(p)}$. Thus, \eqref{eq:ns_incidence2} is isomorphic to $\sigma^*\bar\cH_{s^\bullet2}$, hence
\[
\eqref{eq:ns_incidence1}\simeq\HOM\((\sigma^*\bar\cH_{s^\bullet2})^\vee,\Lie_{\cA^\vee,\tau_\infty^\tc}\)
\simeq\(\sigma^*\bar\cH_{s^\bullet2}\)
\otimes_{\cO_{\DL^\bullet(\sV_{s^\bullet},\{\;,\;\}_{s^\bullet})}}\(\bar\cH_{s^\bullet1}^\dashv/\bar\cH_{s^\bullet2}\),
\]
where we use Theorem \ref{th:ns_basic_correspondence2}(5) for the last isomorphism. We have proved (2) and the theorem.
\end{proof}

We also need a description for
\[
\rB^\dag_{s^\bullet}\coloneqq\rB^\bullet_{s^\bullet}\times_{\rM^\bullet_\fp(\rV^\circ,\rK^{p\circ})}\rM^\dag_\fp(\rV^\circ,\rK^{p\circ})
\]
for $s^\bullet\in\Hk_j^{-1}(s^\bullet_1,s^\bullet_2)(\kappa)$. It is clear that if we put
\[
\rB^\dag_{s^\bullet_i}\coloneqq\rB^\bullet_{s^\bullet_i}\times_{\rM^\bullet_\fp(\rV^\circ,\rK^{p\circ})}\rM^\dag_\fp(\rV^\circ,\rK^{p\circ})
\]
for $i=1,2$, then
\[
\rB^\dag_{s^\bullet}=\rB^\dag_{s^\bullet_1}\times_{\rM^\dag_\fp(\rV^\circ,\rK^{p\circ})}\rB^\dag_{s^\bullet_2}.
\]
By definition, for every $S\in\Sch_{/\kappa}$, $\rB^\dag_{s^\bullet}(S)$ is the set of equivalence classes of unvigintuples
\[
(A_0,\lambda_0,\eta_0^p;A,\lambda,\eta^p;A^\circ,\lambda^\circ,\eta^{p\circ};
A^\bullet_1,\lambda^\bullet_1,\eta^{p\bullet}_1;A^\bullet_2,\lambda^\bullet_2,\eta^{p\bullet}_2;
\beta,\gamma_1,\gamma_2,\psi_1,\psi_2,\phi^\bullet)
\]
rendering the diagram
\[
\xymatrix{
A^\bullet_1 \ar[rr]^-{\phi^\bullet} && A^\bullet_2 \\
& A \ar[lu]_-{\gamma_1}\ar[ru]^-{\gamma_2}\ar[d]^-\beta \\
& A^\circ \ar@/^1pc/[luu]_{\psi_1}\ar@/_1pc/[ruu]^{\psi_2}
}
\]
commute. Here, the letters remain the same meaning as in our previous moduli problems. Put
\[
\rS^\dag_{s^\bullet}\coloneqq\{s^\bullet\}\times_{\rS^\bullet_\fp(\rV^\circ,\rK^{p\circ})\times\rS^\bullet_\fp(\rV^\circ,\rK^{p\circ})}
\(\rS^\dag_\fp(\rV^\circ,\rK^{p\circ})\times\rS^\dag_\fp(\rV^\circ,\rK^{p\circ})\)
\times_{\rS^\circ_\fp(\rV^\circ,\rK^{p\circ})\times\rS^\circ_\fp(\rV^\circ,\rK^{p\circ})}
\rS^\circ_\fp(\rV^\circ,\rK^{p\circ})
\]
where $\rS^\circ_\fp(\rV^\circ,\rK^{p\circ})\to\rS^\circ_\fp(\rV^\circ,\rK^{p\circ})\times\rS^\circ_\fp(\rV^\circ,\rK^{p\circ})$ is the diagonal morphism. Then we have a canonical map
\[
\pi^\dag_{s^\bullet}\colon\rB^\dag_{s^\bullet}\to\rS^\dag_{s^\bullet}
\]
of $\kappa$-schemes by forgetting $(A,\lambda,\eta^p)$ and related morphisms.

\begin{theorem}\label{th:ns_incidence_dag}
Let $s^\bullet_1,s^\bullet_2\in\rS^\bullet_\fp(\rV^\circ,\rK^{p\circ})(\kappa)$ be two points for a perfect field $\kappa$ containing $\dF_p^\Phi$. Take $s^\bullet\in\Hk_j^{-1}(s^\bullet_1,s^\bullet_2)(\kappa)$ for some $0\leq j\leq\floor{\tfrac{N}{2}}-1$. Then
the scheme $\rS^\dag_{s^\bullet}$ is a disjoint of $(p+1)(p^3+1)\cdots(p^{2\lfloor\tfrac{N}{2}\rfloor-2j-1}+1)$ copies of $\Spec\kappa$.

Take a point $t^\dag=(A_0,\lambda_0,\eta_0^p;A^\circ,\lambda^\circ,\eta^{p\circ};A^\bullet_1,\lambda^\bullet_1,\eta^{p\bullet}_1;
A^\bullet_2,\lambda^\bullet_2,\eta^{p\bullet}_2;\psi_1,\psi_2,\phi^\bullet)\in\rS^\dag_{s^\bullet}(\kappa)$.
\begin{enumerate}
  \item The assignment sending
     \[
     (A_0,\lambda_0,\eta_0^p;A,\lambda,\eta^p;A^\circ,\lambda^\circ,\eta^{p\circ};
     A^\bullet_1,\lambda^\bullet_1,\eta^{p\bullet}_1;A^\bullet_2,\lambda^\bullet_2,\eta^{p\bullet}_2;
     \beta,\gamma_1,\gamma_2,\psi_1,\psi_2,\phi^\bullet)\in\rB^\dag_{s^\bullet}(S)
     \]
     to $H_2/(\IM(p\phi^{\bullet-1})_{*,\tau_\infty^\tc}+\rH^\dr_1(A_1^\bullet/S)_{\tau_\infty}^\perp)$ induces an isomorphism
     \[
     \zeta^\dag_{t^\dag}\colon(\pi^{\dag}_{s^\bullet})^{-1}(t^\dag)\xrightarrow{\sim}\dP(\sV_{t^\dag})
     \]
     where we put
     \[
     \sV_{t^\dag}\coloneqq\frac{\IM(\psi_1)_{*,\tau_\infty^\tc}}
     {\IM(p\phi^{\bullet-1})_{*,\tau_\infty^\tc}+\rH^\dr_1(A_1^\bullet/S)_{\tau_\infty}^\perp}
     \]
     which has dimension $\floor{\tfrac{N}{2}}-j$.

  \item The cokernel of the map
      \[
      \cT_{\rB^\dag_{s^\bullet_1}/\kappa}\res_{(\pi^{\dag}_{s^\bullet})^{-1}(t^\dag)}
      \bigoplus\cT_{\rB^\dag_{s^\bullet_2}/\kappa}\res_{(\pi^{\dag}_{s^\bullet})^{-1}(t^\dag)}
      \to\iota^{\bullet*}\cT_{\rM^\dag_\fp(\rV^\circ,\rK^{p\circ})/\kappa}\res_{(\pi^{\dag}_{s^\bullet})^{-1}(t^\dag)}
      \]
      is canonically isomorphic to
      \[
      \zeta^{\dag*}_{t^\dag}\(\(\sigma^*\cH_{t^\dag}\)\otimes_{\cO_{\dP(\sV_{t^\dag})}}\cO_{\dP(\sV_{t^\dag})}(1)\)
      \]
      where $\cH_{t^\dag}$ is the universal object, namely, the tautological bundle on $\dP(\sV_{t^\dag})$.
\end{enumerate}
\end{theorem}

\begin{proof}
In fact, the assignment sending $(A_0,\lambda_0,\eta_0^p;A^\circ,\lambda^\circ,\eta^{p\circ};A^\bullet_1,\lambda^\bullet_1,\eta^{p\bullet}_1;
A^\bullet_2,\lambda^\bullet_2,\eta^{p\bullet}_2;\psi_1,\psi_2,\phi^\bullet)\in\rS^\dag_{s^\bullet}(S)$ to $\IM(\psi_1)_{*,\tau_\tc^\infty}$ induces a bijection from $\rS^\dag_{s^\bullet}(S)$ to the subbundles $H\subseteq\rH^\dr_1(A_1^\bullet/S)_{\tau_\infty^\tc}$ of rank $\lceil\tfrac{N}{2}\rceil$ satisfying $\IM(p\phi^{\bullet-1})_{*,\tau_\infty^\tc}\otimes_\kappa\cO_S\subseteq H\subseteq \Ker(p\phi^\bullet)_{*,\tau_\infty^\tc}\otimes_\kappa\cO_S$ and $\langle\tV^{-1}H,H\rangle_{\tau_\infty^\tc}=0$. Thus, we know that $\rS^\dag_{s^\bullet}$ is a disjoint of $(p+1)(p^3+1)\cdots(p^{2\lfloor\frac{N}{2}\rfloor-2j-1}+1)$ copies of $\Spec\kappa$.

For (1), we denote by $s^\dag_1$ the image of $t^\dag$ in $\rS^\dag_\fp(\rV^\circ,\rK^{p\circ})(\kappa)$ in the first factor. Then a point $(A_0,\lambda_0,\eta_0^p;A,\lambda,\eta^p;A^\circ,\lambda^\circ,\eta^{p\circ};A^\bullet_1,\lambda^\bullet_1,\eta^{p\bullet}_1;\beta,\gamma_1)
\in\rB^\dag_{s^\dag_1}(S)$ belongs to $\rB^\dag_{s^\bullet}(S)$ if and only if $H_2$ contains $\IM(p\phi^{\bullet-1})_{*,\tau_\infty^\tc}\otimes_\kappa\cO_S$. Thus, (1) follows from Theorem \ref{th:ns_basic_correspondence3}(2).

For (2), it follows from Theorem \ref{th:ns_incidence}(2) and the isomorphism
\[
\(\bar\cH_{s^\bullet1}^\dashv/\bar\cH_{s^\bullet2}\)\res_{\dP(\sV_{t^\dag})}
=\(\bar\cH_{s^\bullet1}/\bar\cH_{s^\bullet2}\)\res_{\dP(\sV_{t^\dag})}
\simeq\cO_{\dP(\sV_{t^\dag})}(1).
\]
\end{proof}

\subsection{Incidence maps on the ground stratum}
\label{ss:ns_incidence}

In this subsection, we define and study the incidence maps on ground stratum. We assume $N\geq 2$. In order to have a uniformization map for $\rS^\bullet_\fp(\rV^\circ,\obj)$, we also choose data as in Construction \ref{cs:ns_uniformization2}.

\begin{definition}\label{de:ns_hecke}
We denote
\begin{itemize}[label={\ding{109}}]
  \item $\dT^\circ_{N,\fp}$ the Hecke algebra $\dZ[\rK^\circ_\fp\backslash\rU(\rV^\circ)(F^+_\fp)/\rK^\circ_\fp]$;

  \item $\dT^\bullet_{N,\fp}$ the Hecke algebra $\dZ[\rK^\bullet_\fp\backslash\rU(\rV^\circ)(F^+_\fp)/\rK^\bullet_\fp]$;

  \item $\tT^{\bullet\circ}_{N,\fp}\in\dZ[\rK^\bullet_\fp\backslash\rU(\rV^\circ)(F^+_\fp)/\rK^\circ_\fp]$ the characteristic function of $\rK^\bullet_\fp\rK^\circ_\fp$; and

  \item $\tT^{\circ\bullet}_{N,\fp}\in\dZ[\rK^\circ_\fp\backslash\rU(\rV^\circ)(F^+_\fp)/\rK^\bullet_\fp]$ the characteristic function of $\rK^\circ_\fp\rK^\bullet_\fp$.
\end{itemize}
Moreover, we define the \emph{intertwining Hecke operator} to be
\[
\tI^\circ_{N,\fp}\coloneqq\tT^{\circ\bullet}_{N,\fp}\circ\tT^{\bullet\circ}_{N,\fp}\in\dT^\circ_{N,\fp}
\]
where the composition is taken as composition of cosets.
\end{definition}

\begin{remark}
We remind the readers that according to our convention, the unit elements of $\dZ[\rK^\circ_\fp\backslash\rU(\rV^\circ)(F^+_\fp)/\rK^\circ_\fp]$ and $\dZ[\rK^\bullet_\fp\backslash\rU(\rV^\circ)(F^+_\fp)/\rK^\bullet_\fp]$ are $\CF_{\rK^\circ_\fp}$ and $\CF_{\rK^\bullet_\fp}$, respectively. However, when $N$ is odd, $\rK^\circ_\fp$ and $\rK^\bullet_\fp$ have different volumes under a common Haar measure on $\rU(\rV^\circ)(F^+_\fp)$; in other words, the convolution products on the two Hecke algebras are not induced by the same Haar measure on $\rU(\rV^\circ)(F^+_\fp)$.
\end{remark}

Let $L$ be a $p$-coprime coefficient ring. By Construction \ref{cs:ns_uniformization1} and Construction \ref{cs:ns_uniformization2}, we have canonical isomorphisms
\begin{align*}
L[\Sh(\rV^\circ,\obj\rK^\circ_p)]&\simeq \rH^0_\fT(\ol\rS^\circ_\fp(\rV^\circ,\obj),L), \\
L[\Sh(\rV^\circ,\obj\rK^\bullet_p)]&\simeq \rH^0_\fT(\ol\rS^\bullet_\fp(\rV^\circ,\obj),L),
\end{align*}
in $\Fun(\fK(\rV^\circ)^p,\Mod(L[\rK^\circ_\fp\backslash\rU(\rV^\circ)(F^+_\fp)/\rK^\circ_\fp]))$ and in $\Fun(\fK(\rV^\circ)^p,\Mod(L[\rK^\bullet_\fp\backslash\rU(\rV^\circ)(F^+_\fp)/\rK^\bullet_\fp]))$, induced by $\upsilon^\circ$ \eqref{eq:ns_uniformization1} and $\upsilon^\bullet$ \eqref{eq:ns_uniformization2}, respectively.

\begin{construction}\label{cs:ns_incidence}
Recall from Definition \ref{de:ns_primitive} the class $\xi\in \rH^2_{\fT}(\ol\rB^{\circ}_{\fp}(\rV^{\circ},\obj),L(1))$, which is the first Chern class of the tautological quotient line bundle on $\ol\rB^{\circ}_{\fp}(\rV^{\circ},\obj)$. Put $r\coloneqq\lfloor\tfrac{N}{2}\rfloor\geq 1$. We construct three pairs of maps in $\Fun(\fK(\rV^\circ)^p,\Mod(L))$ as follows:
\begin{align*}
\begin{dcases}
\inc^\circ_!\colon L[\Sh(\rV^\circ,\obj\rK^\circ_p)]
\xrightarrow{\sim}\rH^0_\fT(\ol\rS^\circ_\fp(\rV^\circ,\obj),L)
\xrightarrow{\pi^{\circ*}}\rH^0_\fT(\ol\rB^\circ_\fp(\rV^\circ,\obj),L) \\
\qquad\qquad\xrightarrow{\cup \xi^{N-r-1}}\rH^{2(N-r-1)}_\fT(\ol\rB^\circ_\fp(\rV^\circ,\obj),L(N-r-1)) \\
\qquad\qquad\xrightarrow{\iota^\circ_!}\rH^{2(N-r-1)}_\fT(\ol\rM^\circ_\fp(\rV^\circ,\obj),L(N-r-1)), \\     \inc_\circ^*\colon\rH^{2r}_\fT(\ol\rM^\circ_\fp(\rV^\circ,\obj),L(r))
\xrightarrow{\iota^{\circ*}}\rH^{2r}_\fT(\ol\rB^\circ_\fp(\rV^\circ,\obj),L(r))  \\
\qquad\qquad\xrightarrow{\cup \xi^{N-r-1}}\rH^{2(N-1)}_\fT(\ol\rB^\circ_\fp(\rV^\circ,\obj),L(N-1)) \\
\qquad\qquad\xrightarrow{\pi^\circ_!}\rH^0_\fT(\ol\rS^\circ_\fp(\rV^\circ,\obj),L)
\xrightarrow{\sim}L[\Sh(\rV^\circ,\obj\rK^\circ_p)];
\end{dcases}
\end{align*}
\begin{align*}
\begin{dcases}
\inc^\dag_!\colon L[\Sh(\rV^\circ,\obj\rK^\circ_p)]
\xrightarrow{\sim}\rH^0_\fT(\ol\rS^\circ_\fp(\rV^\circ,\obj),L)
\xrightarrow{\pi^{\circ*}}\rH^0_\fT(\ol\rB^\circ_\fp(\rV^\circ,\obj),L) \\
\qquad\qquad\xrightarrow{\cup\xi^{N-r-2}}\rH^{2(N-r-2)}_\fT(\ol\rB^\circ_\fp(\rV^\circ,\obj),L(N-r-2)) \\
\qquad\qquad\xrightarrow{\iota^\circ_!}\rH^{2(N-r-2)}_\fT(\ol\rM^\circ_\fp(\rV^\circ,\obj),L(N-r-2)) \\
\qquad\qquad\xrightarrow{\rm^{\dag\circ*}}\rH^{2(N-r-2)}_\fT(\ol\rM^\dag_\fp(\rV^\circ,\obj),L(N-r-2)) \\
\qquad\qquad\xrightarrow{\rm^{\dag\bullet}_!}\rH^{2(N-r-1)}_\fT(\ol\rM^\bullet_\fp(\rV^\circ,\obj),L(N-r-1)), \\
\inc_\dag^*\colon\rH^{2r}_\fT(\ol\rM^\bullet_\fp(\rV^\circ,\obj),L(r))
\xrightarrow{\rm^{\dag\bullet*}}\rH^{2r}_\fT(\ol\rM^\dag_\fp(\rV^\circ,\obj),L(r)) \\
\qquad\qquad\xrightarrow{\rm^{\dag\circ}_!}\rH^{2(r+1)}_\fT(\ol\rM^\circ_\fp(\rV^\circ,\obj),L(r+1)) \\
\qquad\qquad\xrightarrow{\iota^{\circ*}}\rH^{2(r+1)}_\fT(\ol\rB^\circ_\fp(\rV^\circ,\obj),L(r+1)  \\
\qquad\qquad\xrightarrow{\cup \xi^{N-r-2}}\rH^{2(N-1)}_\fT(\ol\rB^\circ_\fp(\rV^\circ,\obj),L(N-1)) \\
\qquad\qquad\xrightarrow{\pi^\circ_!}\rH^0_\fT(\ol\rS^\circ_\fp(\rV^\circ,\obj),L)
\xrightarrow{\sim}L[\Sh(\rV^\circ,\obj\rK^\circ_p)];
\end{dcases}
\end{align*}
\begin{align*}
\begin{dcases}
\inc^\bullet_!\colon L[\Sh(\rV^\circ,\obj\rK^\bullet_p)]
\xrightarrow{\sim}\rH^0_\fT(\ol\rS^\bullet_\fp(\rV^\circ,\obj),L)
\xrightarrow{\pi^{\bullet*}}\rH^0_\fT(\ol\rB^\bullet_\fp(\rV^\circ,\obj),L) \\
\qquad\qquad\xrightarrow{\iota^\bullet_!}\rH^{2(N-r-1)}_\fT(\ol\rM^\bullet_\fp(\rV^\circ,\obj),L(N-r-1)), \\
\inc_\bullet^*\colon\rH^{2r}_\fT(\ol\rM^\bullet_\fp(\rV^\circ,\obj),L(r))
\xrightarrow{\iota^{\bullet*}}\rH^{2r}_\fT(\ol\rB^\bullet_\fp(\rV^\circ,\obj),L(r))  \\ \qquad\qquad\xrightarrow{\pi^\bullet_!}\rH^0_\fT(\ol\rS^\bullet_\fp(\rV^\circ,\obj),L)
\xrightarrow{\sim}L[\Sh(\rV^\circ,\obj\rK^\bullet_p)].
\end{dcases}
\end{align*}
Note that the construction of the second pair only makes sense when $N\geq 3$; and when $N=2$, we regard $\inc^\dag_!$ and $\inc_\dag^*$ as zero maps. In fact, the two maps in each pair are essentially Poincar\'{e} dual to each other.
\end{construction}

\begin{definition}\label{de:ns_incidence_odd}
Suppose that $N=2r+1$ is odd with $r\geq 1$.  We define the \emph{incidence map (on the ground stratum)} to be the map
\[
\inc\colon L[\Sh(\rV^\circ,\obj\rK^\circ_p)]\bigoplus L[\Sh(\rV^\circ,\obj\rK^\bullet_p)]
\to L[\Sh(\rV^\circ,\obj\rK^\circ_p)]\bigoplus L[\Sh(\rV^\circ,\obj\rK^\bullet_p)]
\]
in $\Fun(\fK(\rV^\circ)^p,\Mod(L))$ given by the matrix
\[
\begin{pmatrix}
\inc_\dag^*\circ\inc^\dag_! & \inc_\dag^*\circ\inc^\bullet_! \\
\inc_\bullet^*\circ\inc^\dag_! & \inc_\bullet^*\circ\inc^\bullet_!
\end{pmatrix}
\]
if we write elements in the column form.
\end{definition}

\begin{remark}\label{re:ns_incidence_odd}
The construction of the incidence map can be encoded in the following diagram
\[
\xymatrix{
L[\Sh(\rV^\circ,\obj\rK^\circ_p)] \ar[d] && L[\Sh(\rV^\circ,\obj\rK^\bullet_p)] \ar[d] \\
\rH^{2r-2}_\fT(\ol\rM^\dag_\fp(\rV^\circ,\obj),L(r-1)) \ar[rd]^{\rm^{\dag\bullet}_!} &&
\rH^0_\fT(\ol\rB^\bullet_\fp(\rV^\circ,\obj),L) \ar[ld]_-{\iota^\bullet_!} \\
& \rH^{2r}_\fT(\ol\rM^\bullet_\fp(\rV^\circ,\obj),L(r)) \ar[ld]^-{\rm^{\dag\bullet*}} \ar[rd]_-{\iota^{\bullet*}} \\
\rH^{2r}_\fT(\ol\rM^\dag_\fp(\rV^\circ,\obj),L(r)) \ar[d] && \rH^{2r}_\fT(\ol\rB^\bullet_\fp(\rV^\circ,\obj),L(r)) \ar[d]\\
L[\Sh(\rV^\circ,\obj\rK^\circ_p)]  && L[\Sh(\rV^\circ,\obj\rK^\bullet_p)]
}
\]
in $\Fun(\fK(\rV^\circ)^p,\Mod(L))$.
\end{remark}

\begin{proposition}\label{pr:ns_incidence_odd}
Suppose that $N=2r+1$ is odd with $r\geq 1$. Then the incidence map $\inc$ is given by the matrix
\[
\begin{pmatrix}
-(p+1)^2 & \tT^{\circ\bullet}_{N,\fp} \\
\tT^{\bullet\circ}_{N,\fp} & \tT^\bullet_{N,\fp}
\end{pmatrix}
\]
where
\[
\tT^\bullet_{N,\fp}\coloneqq\sum_{\delta=0}^{r-1}\td^\bullet_{r-\delta,p}\cdot\tT^\bullet_{N,\fp;\delta}
\]
in which the numbers $\td^\bullet_{r-\delta,p}$ are introduced in Notation \ref{no:numerical}, and the Hecke operators $\tT^\bullet_{N,\fp;\delta}$ are introduced in Notation \ref{no:hecke} (as $\tT^\bullet_{N;\delta}$).
\end{proposition}

\begin{proof}
Take an object $\rK^{p\circ}\in\fK(\rV^\circ)^p$.

First, we show $\inc_\dag^*\circ\inc^\dag_!=-(p+1)^2$. Since $\rm^{\dag\circ*}\cO_{\ol\rM^\circ_\fp(\rV^\circ,\rK^{p\circ})}(1)$ has degree $p+1$, it follows from Corollary \ref{co:normal}.

Second, we show $\inc_\dag^*\circ\inc^\bullet_!=\tT^{\circ\bullet}_{N,\fp}$ and $\inc_\bullet^*\circ\inc^\dag_!=\tT^{\bullet\circ}_{N,\fp}$. However, these are consequences of Theorem \ref{th:ns_basic_correspondence3} and Construction \ref{cs:ns_uniformization3}.

Finally, we show $\inc_\bullet^*\circ\inc^\bullet_!=\tT^\bullet_{N,\fp}$. By Theorem \ref{th:ns_incidence}(1), it suffices to show that for every $s^\bullet_1,s^\bullet_2\in\rS^\bullet_\fp(\rV^\circ,\rK^{p\circ})(\ol\dF_p)$ and every $s^\bullet\in\Hk_j^{-1}(s^\bullet_1,s^\bullet_2)$, the intersection multiplicity of $\rB^\bullet_{s^\bullet_1}$ and $\rB^\bullet_{s^\bullet_2}$ at the component $\rB^\bullet_{s^\bullet}$ equals $\td^\bullet_{r-j,p}$. This is true by Theorem \ref{th:ns_incidence}(2), Proposition \ref{pr:dl_excess}(1), and the excess intersection formula.

The proposition is proved.
\end{proof}

Now we assume that $N=2r$ is even with $r\geq 2$. The readers may have noticed that the situation is different from Definition \ref{de:ns_incidence_odd} since now $\rM^\bullet_\fp(\rV^\circ,\obj)$ has dimension $2r-1$ while $\rB^\bullet_\fp(\rV^\circ,\obj)$ still has dimension $r$. Thus to obtain a similar diagram as in Remark \ref{re:ns_incidence_odd}, we have to insert a map
\[
\Theta\colon\rH^{2r-2}_\fT(\ol\rM^\bullet_\fp(\rV^\circ,\obj),L(r-1))\to\rH^{2r}_\fT(\ol\rM^\bullet_\fp(\rV^\circ,\obj),L(r))
\]
to obtain a diagram like
\[
\xymatrix{
L[\Sh(\rV^\circ,\obj\rK^\circ_p)] \ar[d] && L[\Sh(\rV^\circ,\obj\rK^\bullet_p)] \ar[d] \\
\rH^{2r-4}_\fT(\ol\rM^\dag_\fp(\rV^\circ,\obj),L(r-2)) \ar[rd]^{\rm^{\dag\bullet}_!} &&
\rH^0_\fT(\ol\rB^\bullet_\fp(\rV^\circ,\obj),L) \ar[ld]_-{\iota^\bullet_!} \\
& \rH^{2r-2}_\fT(\ol\rM^\bullet_\fp(\rV^\circ,\obj),L(r-1)) \ar[d]^-{\Theta} \\
& \rH^{2r}_\fT(\ol\rM^\bullet_\fp(\rV^\circ,\obj),L(r))
\ar[ld]^-{\rm^{\dag\bullet*}} \ar[rd]_-{\iota^{\bullet*}} \\
\rH^{2r}_\fT(\ol\rM^\dag_\fp(\rV^\circ,\obj),L(r)) \ar[d] && \rH^{2r}_\fT(\ol\rB^\bullet_\fp(\rV^\circ,\obj),L(r)) \ar[d] \\
L[\Sh(\rV^\circ,\obj\rK^\circ_p)]  && L[\Sh(\rV^\circ,\obj\rK^\bullet_p)].
}
\]

\begin{definition}\label{de:ns_incidence_even}
For every line bundle $\cL$ on $\rM^\bullet_\fp(\rV^\circ,\obj)$,\footnote{A line bundle $\cL$ on $\rM^\bullet_\fp(\rV^\circ,\obj)$ is a collection of a line bundle $\cL(\rK^{p\circ})$ on every $\rM^\bullet_\fp(\rV^\circ,\rK^{p\circ})$, compatible with respect to pullbacks.} we denote
\[
\Theta_\cL\colon\rH^{2r-2}_\fT(\ol\rM^\bullet_\fp(\rV^\circ,\obj),L(r-1))\to\rH^{2r}_\fT(\ol\rM^\bullet_\fp(\rV^\circ,\obj),L(r))
\]
the map by taking cup product with $c_1(\cL)$, and define the \emph{$\cL$-incidence map (on the ground stratum)} to be the map
\[
\inc_\cL\colon L[\Sh(\rV^\circ,\obj\rK^\circ_p)]\bigoplus L[\Sh(\rV^\circ,\obj\rK^\bullet_p)]
\to L[\Sh(\rV^\circ,\obj\rK^\circ_p)]\bigoplus L[\Sh(\rV^\circ,\obj\rK^\bullet_p)]
\]
in $\Fun(\fK(\rV^\circ)^p,\Mod(L))$ given by the matrix
\[
\begin{pmatrix}
\inc_\dag^*\circ\Theta_\cL\circ\inc^\dag_! & \inc_\dag^*\circ\Theta_\cL\circ\inc^\bullet_! \\
\inc_\bullet^*\circ\Theta_\cL\circ\inc^\dag_! & \inc_\bullet^*\circ\Theta_\cL\circ\inc^\bullet_!
\end{pmatrix},
\]
if we write elements in the column form.
\end{definition}

We now compute $\Theta_\cL$ for two natural choices of $\cL$, namely, $\cO(\rM^\dag_\fp(\rV^\circ,\obj))$ and $\Lie_{\cA,\tau_\infty^\tc}$.

\begin{proposition}\label{pr:ns_incidence_even1}
Suppose that $N=2r$ is even with $r\geq 2$. Let $L$ be a $p$-coprime coefficient ring. For $\cL=\cO(\rM^\dag_\fp(\rV^\circ,\obj))$, the incidence map $\inc_\cL$ is given by
\[
\begin{pmatrix}
(p+1)^3 & -(p+1)\tT^{\circ\bullet}_{N,\fp} \\
-(p+1)\tT^{\bullet\circ}_{N,\fp} & \tR^\bullet_{N,\fp}
\end{pmatrix},
\]
where
\[
\tR^\bullet_{N,\fp}\coloneqq\sum_{\delta=0}^{r-1}
\frac{1-(-p)^{r-\delta}}{p+1} (p+1)(p+3)\cdots(p^{2(r-\delta)-1}+1)\cdot\tT^\bullet_{N,\fp;\delta}
\]
in which the Hecke operators $\tT^\bullet_{N,\fp;\delta}$ are introduced in Notation \ref{no:hecke} (as $\tT^\bullet_{N;\delta}$).
\end{proposition}

\begin{proof}
Take an object $\rK^{p\circ}\in\fK(\rV^\circ)^p$.

First, we show $\inc_\dag^*\circ\Theta_\cL\circ\inc^\dag_!=(p+1)^3$. Since $\rm^{\dag\circ*}\cO_{\ol\rM^\circ_\fp(\rV^\circ,\rK^{p\circ})}(1)$ has degree $p+1$, it follows from Corollary \ref{co:normal}.

Second, we show $\inc_\dag^*\circ\Theta_\cL\circ\inc^\bullet_!=-(p+1)\tT^{\circ\bullet}_{N,\fp}$ and $\inc_\bullet^*\circ\Theta_\cL\circ\inc^\dag_!=-(p+1)\tT^{\bullet\circ}_{N,\fp}$. However, these are consequences of Corollary \ref{co:normal}, Theorem \ref{th:ns_basic_correspondence3}, and Construction \ref{cs:ns_uniformization3}.

It remains to compute $\inc_\bullet^*\circ\Theta_\cL\circ\inc^\bullet_!$. By Theorem \ref{th:ns_incidence}(1), it suffices to show that for every $s^\bullet_1,s^\bullet_2\in\rS^\bullet_\fp(\rV^\circ,\rK^{p\circ})(\ol\dF_p)$ and every $s^\bullet\in\Hk_j^{-1}(s^\bullet_1,s^\bullet_2)$, the intersection multiplicity of $\rB^\dag_{s^\bullet_1}$ and $\rB^\dag_{s^\bullet_2}$ at the component $\rB^\dag_{s^\bullet}$ equals
\[
\frac{1-(-p)^{r-j}}{p+1} (p+1)(p+3)\cdots(p^{2(r-j)-1}+1).
\]
By Theorem \ref{th:ns_incidence_dag} and the excess intersection formula, such intersection multiplicity equals
\[
\sum_{t^\dag\in\rS^\dag_{s^\bullet}(\ol\dF_p)}\int_{\dP(\sV_{t^\dag})}
c_{r-j-1}\(\(\sigma^*\cH_{t^\dag}\)\otimes_{\cO_{\dP(\sV_{t^\dag})}}\cO_{\dP(\sV_{t^\dag})}(1)\).
\]
A simple exercise shows that
\[
\int_{\dP(\sV_{t^\dag})}
c_{r-j-1}\(\(\sigma^*\cH_{t^\dag}\)\otimes_{\cO_{\dP(\sV_{t^\dag})}}\cO_{\dP(\sV_{t^\dag})}(1)\)
=\frac{1-(-p)^{r-j}}{p+1}
\]
for every $t^\dag\in\rS^\dag_{s^\bullet}(\ol\dF_p)$. Thus, the claim follows from Theorem \ref{th:ns_incidence_dag}.
\end{proof}

\begin{proposition}\label{pr:ns_incidence_even2}
Suppose that $N=2r$ is even with $r\geq 2$. Let $L$ be a $p$-coprime coefficient ring. For $\cL=\Lie_{\cA,\tau_\infty^\tc}$, the incidence map $\inc_\cL$ is given by
\[
\begin{pmatrix}
-(p+1)^2 & \tT^{\circ\bullet}_{N,\fp} \\
\tT^{\bullet\circ}_{N,\fp} & \tT^\bullet_{N,\fp}
\end{pmatrix},
\]
where
\[
\tT^\bullet_{N,\fp}\coloneqq\sum_{\delta=0}^{r-1}\rd^\bullet_{r-\delta,p}\cdot\tT^\bullet_{N,\fp;\delta}
\]
in which the numbers $\td^\bullet_{r-\delta,p}$ are introduced in Notation \ref{no:numerical}, and the Hecke operators $\tT^\bullet_{N,\fp;\delta}$ are introduced in Notation \ref{no:hecke} (as $\tT^\bullet_{N;\delta}$).
\end{proposition}

\begin{proof}
Take an object $\rK^{p\circ}\in\fK(\rV^\circ)^p$. By Theorem \ref{th:ns_basic_correspondence1}, we have an isomorphism
\begin{align}\label{eq:ns_incidence_even21}
\iota^{\bullet*}\Lie_{\cA,\tau_\infty^\tc}\simeq\rm^{\dag\circ*}\cO_{\rM^\circ_\fp(\rV^\circ,\rK^{p\circ})}(1)
\end{align}
of line bundles on $\rM^\dag_\fp(\rV^\circ,\rK^{p\circ})$.

First, we show $\inc_\dag^*\circ\Theta_\cL\circ\inc^\dag_!=-(p+1)^2$. This is a consequence of \eqref{eq:ns_incidence_even21}, Corollary \ref{co:normal} and the fact that $\rm^{\dag\circ*}\cO_{\ol\rM^\circ_\fp(\rV^\circ,\rK^{p\circ})}(1)$ has degree $p+1$.

Second, we show $\inc_\dag^*\circ\Theta_\cL\circ\inc^\bullet_!=\tT^{\circ\bullet}_{N,\fp}$ and $\inc_\bullet^*\circ\Theta_\cL\circ\inc^\dag_!=\tT^{\bullet\circ}_{N,\fp}$. These are consequences of \eqref{eq:ns_incidence_even21} and Corollary \ref{co:normal}, Theorem \ref{th:ns_basic_correspondence3}, and Construction \ref{cs:ns_uniformization3}.

It remains to compute $\inc_\bullet^*\circ\Theta_\cL\circ\inc^\bullet_!$. By Theorem \ref{th:ns_incidence} and the excess intersection formula, it suffices to show that for every $s^\bullet_1,s^\bullet_2\in\rS^\bullet_\fp(\rV^\circ,\rK^{p\circ})(\ol\dF_p)$ and every $s^\bullet\in\Hk_j^{-1}(s^\bullet_1,s^\bullet_2)$, we have
\begin{align}\label{eq:ns_incidence_even22}
\int_{\DL^\bullet(\sV_{s^\bullet},\{\;,\;\}_{s^\bullet})}
c_{r-1}\(\(\sigma^*\bar\cH_{s^\bullet 2}\)\otimes_{\cO_{\DL^\bullet(\sV_{s^\bullet},\{\;,\;\}_{s^\bullet})}}\(\bar\cH_{s^\bullet 1}^\dashv/\bar\cH_{s^\bullet 2}\)\)\cdot c_1\((\zeta^\bullet_{s^\bullet})_*\Lie_{\cA,\tau_\infty^\tc}\)
=\td^\bullet_{r-j,p},
\end{align}
where $(\bar\cH_{s^\bullet 1},\bar\cH_{s^\bullet 2})$ is the universal object over $\DL^\bullet(\sV_{s^\bullet},\{\;,\;\}_{s^\bullet})$. However, by Theorem \ref{th:ns_basic_correspondence2}(5), we have $(\zeta^\bullet_{s^\bullet})_*\Lie_{\cA,\tau_\infty^\tc}\simeq\bar\cH_{s^\bullet 1}^\dashv/\bar\cH_{s^\bullet 2}$. Thus, \eqref{eq:ns_incidence_even22} follows from Proposition \ref{pr:dl_excess}(2). The proposition is proved.
\end{proof}

\subsection{Weight spectral sequence}
\label{ss:ns_weight}

In this subsection, we study the weight spectral sequence associated to $\bM_\fp(\rV^\circ,\obj)$. Our goal is to express certain important terms of the weight spectral sequence in terms of $\Sh(\rV^\circ,\obj\rK^\circ_\fp)$ and $\Sh(\rV^\circ,\obj\rK^\bullet_\fp)$. We keep the setup in \S\ref{ss:ns_incidence}. In particular, $N$ is an integer at least $2$ with $r\coloneqq\lfloor\tfrac{N}{2}\rfloor\geq 1$, and $L$ is a $p$-coprime coefficient ring. To ease notation, we put $\rX^?_N\coloneqq\rX^?_\fp(\rV^\circ,\obj)$ for meaningful pairs $(\rX,?)\in\{\bM,\rM,\rB,\rS\}\times\{\;,\circ,\bullet,\dag\}$.

\begin{construction}\label{cs:ns_weight}
By Theorem \ref{th:ns_moduli_scheme}(1), we have the weight spectral sequence $(\rE^{p,q}_s,\rd^{p,q}_s)$, with terms in the category $L[\Gal(\ol\dF_p/\dF_p^\Phi)]$, abutting to the cohomology $\rH^{p+q}_\fT(\ol\rM_N,\rR\Psi L(r))$. In particular, we have
\[
\rE^{0,2d}_1=\rH^{2d}_\fT(\ol\rM^\circ_N,L(r))\bigoplus\rH^{2d}_\fT(\ol\rM^\bullet_N,L(r)).
\]
Thus, the six maps in Construction \ref{cs:ns_incidence} give rise to another six maps
\begin{align*}
\begin{dcases}
\Inc^\circ_!\colon L[\Sh(\rV^\circ,\obj\rK^\circ_\fp)]\to\rE^{0,2(N-r-1)}_1(N-2r-1),\\
\Inc^\dag_!\colon L[\Sh(\rV^\circ,\obj\rK^\circ_\fp)]\to\rE^{0,2(N-r-1)}_1(N-2r-1),\\
\Inc^\bullet_!\colon L[\Sh(\rV^\circ,\obj\rK^\bullet_\fp)]\to\rE^{0,2(N-r-1)}_1(N-2r-1),\\
\Inc_\circ^*\colon \rE^{0,2r}_1\to L[\Sh(\rV^\circ,\obj\rK^\circ_\fp)], \\
\Inc_\dag^*\colon \rE^{0,2r}_1\to L[\Sh(\rV^\circ,\obj\rK^\circ_\fp)], \\
\Inc_\bullet^*\colon \rE^{0,2r}_1\to L[\Sh(\rV^\circ,\obj\rK^\bullet_\fp)],
\end{dcases}
\end{align*}
in $\Fun(\fK(\rV^\circ)^p,\Mod(L))$.
\end{construction}

In the future, we will have to study the composite maps
\[
\begin{pmatrix}
\Inc_\circ^* \\
\Inc_\dag^* \\
\Inc_\bullet^*
\end{pmatrix}
\begin{pmatrix}
\Inc^\circ_! & \Inc^\dag_! & \Inc^\bullet_!
\end{pmatrix},\quad
\begin{pmatrix}
\Inc_\circ^* \\
\Inc_\dag^* \\
\Inc_\bullet^*
\end{pmatrix}
\circ\rd^{-1,2r}_1\circ\rd^{0,2r-2}_1(-1)\circ
\begin{pmatrix}
\Inc^\circ_! & \Inc^\dag_! & \Inc^\bullet_!
\end{pmatrix}
\]
when $N$ is odd and even, respectively. In the next two lemmas, we will study the spectral sequence and prove two formulae related to the above maps, according to the parity of $N$.

\begin{lem}\label{le:ns_weight_odd}
Suppose that $N=2r+1$ is odd with $r\geq 1$.
\begin{enumerate}
  \item The first page of $\rE^{p,q}_s$ is as follows:
      \[
      \resizebox{\hsize}{!}{
      \boxed{
      \xymatrix{
      q\geq 2r+2 & \cdots \ar[r] & \cdots \ar[r] & \cdots \\
      q=2r+1 & \rH^{2r-1}_\fT(\ol\rM^\dag_N,L(r-1)) \ar[r]^-{\rd^{-1,2r+1}_1} & \rH^{2r+1}_\fT(\ol\rM^\circ_N,L(r))\oplus
      \rH^{2r+1}_\fT(\ol\rM^\bullet_N,L(r)) \ar[r]^-{\rd^{0,2r+1}_1} &
      \rH^{2r+1}_\fT(\ol\rM^\dag_N,L(r)) \\
      q=2r & \rH^{2r-2}_\fT(\ol\rM^\dag_N,L(r-1)) \ar[r]^-{\rd^{-1,2r}_1} & \rH^{2r}_\fT(\ol\rM^\circ_N,L(r))\oplus
      \rH^{2r}_\fT(\ol\rM^\bullet_N,L(r)) \ar[r]^-{\rd^{0,2r}_1} &
      \rH^{2r}_\fT(\ol\rM^\dag_N,L(r)) \\
      q=2r-1 & \rH^{2r-3}_\fT(\ol\rM^\dag_N,L(r-1)) \ar[r]^-{\rd^{-1,2r-1}_1} & \rH^{2r-1}_\fT(\ol\rM^\circ_N,L(r))\oplus
      \rH^{2r-1}_\fT(\ol\rM^\bullet_N,L(r)) \ar[r]^-{\rd^{0,2r-1}_1} &
      \rH^{2r-1}_\fT(\ol\rM^\dag_N,L(r)) \\
      q\leq 2r-2 & \cdots \ar[r] & \cdots \ar[r] & \cdots \\
      \rE^{p,q}_1 & p=-1 & p=0 & p=1
      }
      }
      }
      \]
      with $\rd^{-1,i}_1=(\rm^{\dag\circ}_!,-\rm^{\dag\bullet}_!)$, $\rd^{0,i}_1=(\rm^{\dag\circ})^*-(\rm^{\dag\bullet})^*$ for every $i\in\dZ$; and $\rE^{p,q}_1=0$ if $|p|>1$.

  \item We have
      \[
      \begin{pmatrix}
      \Inc_\circ^* \\
      \Inc_\dag^* \\
      \Inc_\bullet^*
      \end{pmatrix}
      \begin{pmatrix}
      \Inc^\circ_! & \Inc^\dag_! & \Inc^\bullet_!
      \end{pmatrix}
      =
      \begin{pmatrix}
      1 & 0 & 0 \\
      0 & -(p+1)^2 & \tT^{\circ\bullet}_{N,\fp} \\
      0 & \tT^{\bullet\circ}_{N,\fp} & \tT^\bullet_{N,\fp}
      \end{pmatrix}.
      \]

  \item We have $(\tT^{\bullet\circ}_{N,\fp}\circ\Inc_\dag^*+(p+1)^2\Inc_\bullet^*)\circ\rd^{-1,2r}_1=0$.
\end{enumerate}
\end{lem}

\begin{proof}
Part (1) is immediate. Part (2) is a consequence of Proposition \ref{pr:ns_incidence_odd}.

For (3), note that under the composite isomorphism
\begin{align*}
i&\colon L[\Sh(\rV^\circ,\obj\rK^\circ_\fp)]
\xrightarrow{\sim}\rH^0_\fT(\ol\rS^\circ_N,L)
\xrightarrow{\pi^{\circ*}}\rH^0_\fT(\ol\rB^\circ_N,L)
\xrightarrow{\cup\xi^{r-1}}\rH^{2r-2}_\fT(\ol\rB^\circ_N,L(r-1)) \\
&\xrightarrow{\iota^\circ_!}\rH^{2r-2}_\fT(\ol\rM^\circ_N,L(r-1))
\xrightarrow{\rm^{\dag\circ*}}\rH^{2r-2}_\fT(\ol\rM^\dag_N,L(r-1))
=\rE^{-1,2r}_1,
\end{align*}
the map $\rd^{-1,2r}_1\circ i\colon L[\Sh(\rV^\circ,\obj\rK^\circ_p)]\to\rE^{0,2r}_1$ coincides with $(p+1)\Inc^\circ_!-\Inc^\dag_!$. Thus, (3) follows by (2) as we have
\[
\begin{pmatrix}
0 & \tT^{\bullet\circ}_{N,\fp} & (p+1)^2
\end{pmatrix}
\begin{pmatrix}
1 & 0 & 0 \\
0 & -(p+1)^2 & \tT^{\circ\bullet}_{N,\fp} \\
0 & \tT^{\bullet\circ}_{N,\fp} & \tT^\bullet_{N,\fp}
\end{pmatrix}
\begin{pmatrix}
p+1 \\
-1 \\
0
\end{pmatrix}
=0.
\]
The lemma is proved.
\end{proof}

For $N$ even, we first recall that there is an (increasing) monodromy filtration $\rF_\bullet\rR\Psi L(r)$ of $\rR\Psi L(r)$. Such filtration induces a filtration $\rF_\bullet\rH^i_\fT(\ol\rM_N,\rR\Psi L(r))$ of $\rH^i_\fT(\ol\rM_N,\rR\Psi L(r))$, and a corresponding filtration $\rF_\bullet\rH^1(\rI_{\dQ_p^\Phi},\rH^i_\fT(\ol\rM_N,\rR\Psi L(r)))$ of the quotient module $\rH^1(\rI_{\dQ_p^\Phi},\rH^i_\fT(\ol\rM_N,\rR\Psi L(r)))$.

\begin{lem}\label{le:ns_weight_even}
Suppose that $N=2r$ is even with $r\geq 1$.
\begin{enumerate}
  \item The first page of $\rE^{p,q}_s$ is as follows:
     \[
     \resizebox{\hsize}{!}{
     \boxed{
     \xymatrix{
     q\geq 2r+1 & \cdots \ar[r] & \cdots \ar[r] & \cdots \\
     q=2r & \rH^{2r-2}_\fT(\ol\rM^\dag_N,L(r-1)) \ar[r]^-{\rd^{-1,2r}_1} & \rH^{2r}_\fT(\ol\rM^\circ_N,L(r))\oplus
     \rH^{2r}_\fT(\ol\rM^\bullet_N,L(r)) \ar[r]^-{\rd^{0,2r}_1} &
     \rH^{2r}_\fT(\ol\rM^\dag_N,L(r)) \\
     q=2r-1 & 0\ar[r] & \rH^{2r-1}_\fT(\ol\rM^\bullet_N,L(r)) \ar[r]  & 0 \\
     q=2r-2 & \rH^{2r-4}_\fT(\ol\rM^\dag_N,L(r-1)) \ar[r]^-{\rd^{-1,2r-2}_1} & \rH^{2r-2}_\fT(\ol\rM^\circ_N,L(r))\oplus
     \rH^{2r-2}_\fT(\ol\rM^\bullet_N,L(r)) \ar[r]^-{\rd^{0,2r-2}_1} &
     \rH^{2r-2}_\fT(\ol\rM^\dag_N,L(r)) \\
     q\leq 2r-3 & \cdots \ar[r] & \cdots \ar[r] & \cdots \\
     \rE^{p,q}_1 & p=-1 & p=0 & p=1
     }
     }
     }
     \]
     with $\rd^{-1,i}_1=(\rm^{\dag\circ}_!,-\rm^{\dag\bullet}_!)$, $\rd^{0,i}_1=(\rm^{\dag\circ})^*-(\rm^{\dag\bullet})^*$ for every $i\in\dZ$; and $\rE^{p,q}_1=0$ if $|p|>1$.

  \item The spectral sequence $\rE^{p,q}_s$ degenerates at the second page.

  \item In the (three-step) filtration $\rF_\bullet\rH^{2r-1}_\fT(\ol\rM_N,\rR\Psi L(r))$, we have canonical isomorphisms
      \begin{align*}
      \begin{dcases}
      \rF_{-1}\rH^{2r-1}_\fT(\ol\rM_N,\rR\Psi L(r)) \simeq\rE^{1,2r-2}_2=\coker\rd^{0,2r-2}_1,\\
      \frac{\rF_0\rH^{2r-1}_\fT(\ol\rM_N,\rR\Psi L(r))}{\rF_{-1}\rH^{2r-1}_\fT(\ol\rM_N,\rR\Psi L(r))}
      \simeq\rE^{0,2r-1}_2=\rH^{2r-1}_\fT(\ol\rM^\bullet_N,L(r)), \\
      \frac{\rH^{2r-1}_\fT(\ol\rM_N,\rR\Psi L(r))}{\rF_0\rH^{2r-1}_\fT(\ol\rM_N,\rR\Psi L(r))} \simeq\rE^{-1,2r}_2=\Ker\rd^{-1,2r}_1,
      \end{dcases}
      \end{align*}
      in $\Fun(\fK(\rV^\circ)^p,\Mod(L[\Gal(\ol\dF_p/\dF_p^\Phi)]))$.

  \item The monodromy map on $\rH^{2r-1}_\fT(\ol\rM_N,\rR\Psi L(r))$ is trivial on $\rF_0\rH^{2r-1}_\fT(\ol\rM_N,\rR\Psi L(r))$ and is given by the composite map
      \[
      \rE^{-1,2r}_2\xrightarrow{\mu}\rE^{1,2r-2}_2\hookrightarrow\rH^{2r-1}_\fT(\ol\rM_N,\rR\Psi L(r))
      \]
      in view of (3), where $\mu$ is the map induced from the identity map on $\rH^{2r-2}_\fT(\ol\rM^\dag_N,L(r-1))$.

  \item We have a canonical isomorphism
      \[
      \rF_{-1}\rH^1(\rI_{\dQ_p^\Phi},\rH^{2r-1}_\fT(\ol\rM_N,\rR\Psi L(r)))
      \simeq\(\frac{\rE^{1,2r-2}_2}{\mu\rE^{-1,2r}_2}\)(-1);
      \]
      in $\Fun(\fK(\rV^\circ)^p,\Mod(L[\Gal(\ol\dF_p/\dF_p^\Phi)]))$; and the map $\rd^{-1,2r}_1$ induces an isomorphism
      \[
      \(\frac{\rE^{1,2r-2}_2}{\mu\rE^{-1,2r}_2}\)(-1)\simeq\frac{\IM\rd^{-1,2r}_1}{\IM(\rd^{-1,2r}_1\circ\rd^{0,2r-2}_1(-1))}
      \]
      in $\Fun(\fK(\rV^\circ)^p,\Mod(L[\Gal(\ol\dF_p/\dF_p^\Phi)]))$.

  \item If $p^2-1$ is invertible in $L$, then we have a canonical short exact sequence
      \[
      \resizebox{\hsize}{!}{
      \xymatrix{
      0 \ar[r] & \rF_{-1}\rH^1(\rI_{\dQ_p^\Phi},\rH^{2r-1}_\fT(\ol\rM_N,\rR\Psi L(r)))
      \ar[r] & \rH^1_\sing(\dQ_p^\Phi,\rH^{2r-1}_\fT(\ol\rM_N,\rR\Psi L(r)))
      \ar[r] & \rH^{2r-1}_\fT(\ol\rM^\bullet_N,L(r-1))^{\Gal(\ol\dF_p/\dF_p^\Phi)} \to 0
      }
      }
      \]
      in $\Fun(\fK(\rV^\circ)^p,\Mod(L))$.

  \item The composite map
      \[
      \begin{pmatrix}
      \Inc_\circ^* \\
      \Inc_\dag^* \\
      \Inc_\bullet^*
      \end{pmatrix}
      \circ\rd^{-1,2r}_1\circ\rd^{0,2r-2}_1(-1)\circ
      \begin{pmatrix}
      \Inc^\circ_! & \Inc^\dag_! & \Inc^\bullet_!
      \end{pmatrix}
      \]
      coincides with
      \[
      \begin{pmatrix}
      p+1 & (p+1)^2 & -\tT^{\circ\bullet}_{N,\fp} \\
      (p+1)^2 & (p+1)^3 & -(p+1)\tT^{\circ\bullet}_{N,\fp} \\
      -\tT^{\bullet\circ}_{N,\fp} &
      -(p+1)\tT^{\bullet\circ}_{N,\fp} & \tR^\bullet_{N,\fp}
      \end{pmatrix},\quad
      \begin{pmatrix}
      p+1 & 0 & -\tT^{\circ\bullet}_{N,\fp} \\
      0 & 0 & 0 \\
      -\tT^{\bullet\circ}_{N,\fp} &
      0 & \tR^\bullet_{N,\fp}
      \end{pmatrix}
      \]
      when $N\geq 4$ and when $N=2$, respectively.

  \item The image of the map
      \begin{align*}
      &(\tT^{\bullet\circ}_{N,\fp}\circ\Inc_\circ^*+(p+1)\Inc_\bullet^*)
      \circ\rd^{-1,2r}_1\circ\rd^{0,2r-2}_1(-1)\circ(\Inc^\circ_!+\Inc^\dag_!+\Inc^\bullet_!)\colon \\
      &\quad L[\Sh(\rV^\circ,\obj\rK^\circ_\fp)]^{\oplus 2}\bigoplus L[\Sh(\rV^\circ,\obj\rK^\bullet_\fp)]
      \to L[\Sh(\rV^\circ,\obj\rK^\bullet_\fp)]
      \end{align*}
      is exactly $((p+1)\tR^\bullet_{N,\fp}-\tT^{\bullet\circ}_{N,\fp}\circ\tT^{\circ\bullet}_{N,\fp})L[\Sh(\rV^\circ,\obj\rK^\bullet_\fp)]$, where $\tR^\bullet_{N,\fp}$ is introduced in Proposition \ref{pr:ns_incidence_even1}.
\end{enumerate}
\end{lem}

\begin{proof}
For (1), note that by Lemma \ref{le:ns_cohomology}(1), both $\rH^i_\fT(\ol\rM^\dag_N,L)$ and $\rH^i_\fT(\ol\rM^\circ_N,L)$ vanish for $i$ odd. Thus, (1) follows.

Parts (2--4) follow directly from the description of $\rE^{p,q}_1$ and \cite{Sai03}*{Corollary~2.8(2)} for the description of the monodromy map (which does not require the scheme to be proper over the base). Part (5) follows from (1--4).

For (6), by Lemma \ref{le:ns_cohomology}(3), we know that the action of $\Gal(\ol\dF_p/\dF_p^\Phi)$ on $\rE^{1,2r-2}_2(-1)$ is trivial. As $p^2-1$ is invertible in $L$, we further have $\rE^{-1,2r}_2(-1)^{\Gal(\ol\dF_p/\dF_p^\Phi)}=0$ and
\[
\rH^1(\Gal(\ol\dF_p/\dF_p^\Phi),\rF_{-1}\rH^1(\rI_{\dQ_p^\Phi},\rH^{2r-1}_\fT(\ol\rM_N,\rR\Psi L(r))))=0.
\]
In particular, we have the isomorphism
\begin{align*}
\resizebox{\hsize}{!}{
\xymatrix{
\rH^1_\sing(\dQ_p^\Phi,\rH^{2r-1}_\fT(\ol\rM_N,\rR\Psi L(r)))\simeq
\rH^1(\rI_{\dQ_p^\Phi},\rH^{2r-1}_\fT(\ol\rM_N,\rR\Psi L(r)))^{\Gal(\ol\dF_p/\dF_p^\Phi)}
\simeq\rF_0\rH^1(\rI_{\dQ_p^\Phi},\rH^{2r-1}_\fT(\ol\rM_N,\rR\Psi L(r)))^{\Gal(\ol\dF_p/\dF_p^\Phi)}
}
}
\end{align*}
and that (6) follows from the induced long exact sequence.

For (7), when $N\geq 4$ (that is, $r\geq 2$), it follows from Theorem \ref{th:ns_basic_correspondence1}(2) and Proposition \ref{pr:ns_incidence_even1}; when $N=2$, it follows from a direct computation.

For (8), we have the identity
\begin{align*}
&\quad
\begin{pmatrix}
\tT^{\bullet\circ}_{N,\fp} & 0 & p+1
\end{pmatrix}
\begin{pmatrix}
\Inc_\circ^* \\
\Inc_\dag^* \\
\Inc_\bullet^*
\end{pmatrix}
\circ\rd^{-1,2r}_1\circ\rd^{0,2r-2}_1(-1)\circ
\begin{pmatrix}
\Inc^\circ_! & \Inc^\dag_! & \Inc^\bullet_!
\end{pmatrix}\\
&=
\begin{pmatrix}
0 & 0 & (p+1)\tR^\bullet_{N,\fp}-\tT^{\bullet\circ}_{N,\fp}\circ\tT^{\circ\bullet}_{N,\fp}
\end{pmatrix}
\end{align*}
by (7), which implies (8).
\end{proof}

\begin{construction}\label{cs:ns_nabla}
We construct
\begin{enumerate}
  \item when $N=2r+1$ is odd, the map
     \[
     \nabla^1\colon\rE^{0,2r}_2\to L[\Sh(\rV^\circ_N,\rK^\circ_N)]
     \]
      to be the restriction of the map $\tT^{\bullet\circ}_{N,\fp}\circ\Inc_\dag^*+(p+1)^2\Inc_\bullet^*\colon\rE^{0,2r}_1\to L[\Sh(\rV^\circ_N,\rK^\bullet_N)]$ to $\Ker\rd^{0,2r}_1$, which factors through $\rE^{0,2r}_2$ by Lemma \ref{le:ns_weight_odd}(3), composed with the map $\tT^{\circ\bullet}_{N,\fp}\colon L[\Sh(\rV^\circ_N,\rK^\bullet_N)]\to L[\Sh(\rV^\circ_N,\rK^\circ_N)]$;

  \item when $N=2r$ is even, the map
     \[
     \nabla^0\colon\Ker\rd^{0,2r}_1\to L[\Sh(\rV^\circ_N,\rK^\circ_N)]
     \]
     to be the restriction of the map $\tT^{\bullet\circ}_{N,\fp}\circ\Inc_\circ^*+(p+1)\Inc_\bullet^*\colon\rE^{0,2r}_1\to L[\Sh(\rV^\circ_N,\rK^\bullet_N)]$ in Lemma \ref{le:ns_weight_even}(8) to $\Ker\rd^{0,2r}_1$, composed with the map $\tT^{\circ\bullet}_{N,\fp}\colon L[\Sh(\rV^\circ_N,\rK^\bullet_N)]\to L[\Sh(\rV^\circ_N,\rK^\circ_N)]$.
\end{enumerate}
\end{construction}

\begin{remark}\label{re:ns_nabla}
By the descriptions of the Galois actions in Construction \ref{cs:ns_uniformization1} and Construction \ref{cs:ns_uniformization2}, the map $\nabla^1$ factors through the quotient map $\rE^{0,2r}_2\to(\rE^{0,2r}_2)_{\Gal(\ol\dF_p/\dF_p^\Phi)}$.
\end{remark}

To temporarily end the discussion on weight spectral sequences, we record the following easy lemma, which will be used later.

\begin{lem}\label{le:ns_weight_pre}
Suppose that $N\geq 3$. The following diagram
\[
\xymatrix{
\rE^{0,2r}_1 \ar[rrr]^-{(\Inc^*_\circ,\Inc^*_\dag,\Inc^*_\bullet)}\ar[d]_-{\rd^{0,2r}_1}
&&& L[\Sh(\rV^\circ,\obj\rK^\circ_\fp)]^{\oplus 2}\bigoplus L[\Sh(\rV^\circ,\obj\rK^\bullet_\fp)] \ar[d]^-{(p+1,-1,0)} \\
\rE^{1,2r}_1 \ar[rrr] &&& L[\Sh(\rV^\circ,\obj\rK^\circ_\fp)]
}
\]
is commutative, where the lower horizontal arrow is the composite map
\begin{align*}
&\rH^{2r}_\fT(\ol\rM^\dag_\fp(\rV^\circ,\obj),L(r))
\xrightarrow{\rm^{\dag\circ}_!}\rH^{2(r+1)}_\fT(\ol\rM^\circ_\fp(\rV^\circ,\obj),L(r+1))
\xrightarrow{\iota^{\circ*}}\rH^{2(r+1)}_\fT(\ol\rB^\circ_\fp(\rV^\circ,\obj),L(r+1) \\
&\qquad\qquad\xrightarrow{\cup \xi^{N-r-2}}\rH^{2(N-1)}_\fT(\ol\rB^\circ_\fp(\rV^\circ,\obj),L(N-1))
\xrightarrow{\pi^\circ_!}\rH^0_\fT(\ol\rS^\circ_\fp(\rV^\circ,\obj),L)
\xrightarrow{\sim}L[\Sh(\rV^\circ,\obj\rK^\circ_p)],
\end{align*}
which is an isomorphism.
\end{lem}

\begin{proof}
The commutativity of the diagram follows from the formula $\rd^{0,2r}_1=(\rm^{\dag\circ})^*-(\rm^{\dag\bullet})^*$, and the fact that $\rM^\dag_\fp(\rV^\circ,\obj)$ is a hypersurface in $\rM^\circ_\fp(\rV^\circ,\obj)$ of degree $p+1$ by Theorem \ref{th:ns_basic_correspondence1} and Lemma \ref{le:dl_cohomology}(1). By Lemma \ref{le:ns_cohomology} and the Poincar\'{e} duality theorem, the lower horizontal arrow is an isomorphism.
\end{proof}

\if false

\subsection{Special results in the rank 3 case}
\label{ss:ns_three}

In this subsection, we study some special properties of the ground stratum $\rM^\bullet_\fp(\rV^\circ,\obj)$ when $N=3$. The results here will only be used in the situation (b) of Lemma \ref{le:intertwining} and are only necessary for the main theorems in \S\ref{ss:main} in the case where $n=2$ and $F^+=\dQ$, so the readers may skip this subsection at this moment.

To begin with, we recall the following definition.

\begin{definition}[\cite{FK}*{Chapter~I.~Definition~3.7~\&~Note~3.10}]
A proper morphism $f\colon X\to Y$ of schemes of characteristic $p$ is \emph{purely inseparable} if the following two equivalent conditions hold:
\begin{enumerate}
  \item For every (scheme-theoretical) point $y$ of $Y$, there lies exactly one point $x$ of $X$, and the residue field extension is purely inseparable.

  \item For every algebraically field $\kappa$ of characteristic $p$, the induced map $f(\kappa)\colon X(\kappa)\to Y(\kappa)$ is a bijection.
\end{enumerate}
\end{definition}

We now assume $\dim_F\rV^\circ=N=3$.

\begin{definition}\label{de:ns_moduli_scheme_three}
We define a functor
\begin{align*}
\bM^\prime_\fp(\rV^\circ,\obj)\colon\fK(\rV^\circ)^p\times\fT &\to\sfP\Sch'_{/\dZ^\Phi_p} \\
\rK^{p\circ} &\mapsto \bM^\prime_\fp(\rV^\circ,\rK^{p\circ})
\end{align*}
such that for every $S\in\Sch'_{/\dZ^\Phi_p}$, $\bM^\prime_\fp(\rV^\circ,\rK^{p\circ})(S)$ is the set of equivalence classes of sextuples $(A_0,\lambda_0,\eta_0^p;A^\prime,\lambda^\prime,\eta^{p\prime})$, where
\begin{itemize}[label={\ding{109}}]
  \item $(A_0,\lambda_0,\eta_0^p)$ is an element in $\bT_\fp(S)$;

  \item $(A^\prime,\lambda^\prime)$ is a unitary $O_F$-abelian scheme of signature type $3\Phi-2\tau_\infty+2\tau_\infty^\tc$ over $S$ (Definitions \ref{de:unitary_abelian_scheme} and \ref{de:signature}) such that $\lambda^\prime$ is $p$-principal;

  \item $\eta^{p\prime}$ is a $\rK^{p\circ}$-level structure, that is, for a chosen geometric point $s$ on every connected component of $S$, a $\pi_1(S,s)$-invariant $\rK^{p\circ}$-orbit of isomorphisms
      \[
      \eta^{p\prime}\colon\rV^\circ\otimes_\dQ\dA^{\infty,p}\to
      \Hom_{F\otimes_\dQ\dA^{\infty,p}}^{\lambda_0,\lambda^\prime}
      (\rH^\et_1(A^\prime_{0s},\dA^{\infty,p}),\rH^\et_1(A^\prime_s,\dA^{\infty,p}))
      \]
      of hermitian spaces over $F\otimes_\dQ\dA^{\infty,p}=F\otimes_{F^+}\dA_{F^+}^{\infty,p}$. See Construction \ref{cs:hermitian_structure} (with $\Box=\{\infty,p\}$) for the right-hand side.
\end{itemize}
The equivalence relation and the action of morphisms in $\fK(\rV^\circ)^p\times\fT$ are defined similarly as in Definition \ref{de:qs_moduli_scheme}.
\end{definition}

We clearly have the forgetful morphism
\begin{align}\label{eq:ns_moduli_scheme_three}
\bM^\prime_\fp(\rV^\circ,\obj)\to\bT_\fp
\end{align}
in $\Fun(\fK(\rV^\circ)^p\times\fT,\sfP\Sch'_{/\dZ_p^\Phi})$. By a similar proof of Theorem \ref{th:qs_moduli_scheme}, the morphism \eqref{eq:ns_moduli_scheme_three} is represented by quasi-projective smooth schemes of relative dimension $2$. We denote by the base change of \eqref{eq:ns_moduli_scheme_three} to $\dF_p^\Phi$ by $\rM^\prime_\fp(\rV^\circ,\obj)\to\rT_\fp$, which is a morphism in $\Fun(\fK(\rV^\circ)^p\times\fT,\Sch_{/\dF_p^\Phi})$.

\begin{definition}\label{de:ns_three_correspondence1}
We define a functor
\begin{align*}
\rN_\fp(\rV^\circ,\obj)\colon\fK(\rV^\circ)^p\times\fT &\to\sfP\Sch'_{/\dF^\Phi_p} \\
\rK^{p\circ} &\mapsto \rN_\fp(\rV^\circ,\rK^{p\circ})
\end{align*}
such that for every $S\in\Sch'_{/\dF^\Phi_p}$, $\rN_\fp(\rV^\circ,\rK^{p\circ})(S)$ is the set of equivalence classes of decuples $(A_0,\lambda_0,\eta_0^p;A,\lambda,\eta^p;A^\prime,\lambda^\prime,\eta^{p\prime};\delta)$, where
\begin{itemize}[label={\ding{109}}]
  \item $(A_0,\lambda_0,\eta_0^p;A,\lambda,\eta^p)$ is an element of $\rM_\fp(\rV^\circ,\rK^{p\circ})(S)$;

  \item $(A_0,\lambda_0,\eta_0^p;A^\prime,\lambda^\prime,\eta^{p\prime})$ is an element of $\rM^\prime_\fp(\rV^\circ,\rK^{p\circ})(S)$; and

  \item $\delta\colon A\to A^\prime$ is an $O_F$-linear quasi-$p$-isogeny (Definition \ref{de:p_quasi}) such that
  \begin{enumerate}[label=(\alph*)]
    \item $\Ker\delta[p^\infty]$ is contained in $A[\fp]$;

    \item we have $\lambda=\delta^\vee\circ\lambda^\prime\circ\delta$; and

    \item the $\rK^{p\circ}$-orbit of maps $v\mapsto\delta_*\circ\eta^p(v)$ for $v\in\rV^\circ\otimes_\dQ\dA^{\infty,p}$ coincides with $\eta^{p\prime}$.
  \end{enumerate}
\end{itemize}
The equivalence relation and the action of morphisms in $\fK(\rV^\circ)^p\times\fT$ are defined similarly as in Definition \ref{de:qs_basic_correspondence}.
\end{definition}

By definition, we have the following two obvious forgetful morphisms.
\begin{align*}
\xymatrix{
&  \rN_\fp(\rV^\circ,\obj) \ar[ld]_-{\mu} \ar[rd]^-{\mu'} \\
\rM_\fp(\rV^\circ,\obj) && \rM^\prime_\fp(\rV^\circ,\obj)
}
\end{align*}
in $\Fun(\fK(\rV^\circ)^p\times\fT,\Sch_{/\dF_p^\Phi})$. By the extension property of isogeny, it is clear that both $\mu$ and $\mu'$ are proper. We apply the Stein factorization to the morphism $\mu'$ and obtain the following diagram
\begin{align}\label{eq:three}
\xymatrix{
&  \rN_\fp(\rV^\circ,\obj) \ar[ld]_-{\mu}\ar[r]^-{\nu}  &  \rN^\prime_\fp(\rV^\circ,\obj)  \ar[rd]^-{\nu'} \\
\rM_\fp(\rV^\circ,\obj) &&& \rM^\prime_\fp(\rV^\circ,\obj)
}
\end{align}
in $\Fun(\fK(\rV^\circ)^p\times\fT,\Sch_{/\dF_p^\Phi})$. For every $\rK^{p\circ}\in\fK(\rV^\circ)^p$ a perfect field $\kappa$ containing $\dF_p^\Phi$, we say that a point $(A_0,\lambda_0,\eta_0^p;A^\prime,\lambda^\prime,\eta^{p\prime})\in\rM^\prime_\fp(\rV^\circ,\rK^{p\circ})(\kappa)$ is \emph{special} if we have $\tF\rH^\dr_1(A^\prime/\kappa)_{\tau_\infty}=\tV\rH^\dr_1(A^\prime/\kappa)_{\tau_\infty}$. We denote by $\rM^\prime_\fp(\rV^\circ,\rK^{p\circ})_{\sp}$ the locus of special points in $\rM^\prime_\fp(\rV^\circ,\rK^{p\circ})$, regarded as a Zariski closed subset, and by $\rN^\prime_\fp(\rV^\circ,\rK^{p\circ})_{\sp}$ the (set-theoretical) inverse image of $\rM^\prime_\fp(\rV^\circ,\rK^{p\circ})_{\sp}$ under $\nu'$. An easy deformation argument shows that $\rM^\prime_\fp(\rV^\circ,\rK^{p\circ})_{\sp}$ is of dimension zero.

\begin{proposition}\label{pr:three}
In \eqref{eq:three}, for every $\rK^{p\circ}\in\fK(\rV^\circ)^p$, we have
\begin{enumerate}
  \item The morphism $\mu\colon\rN_\fp(\rV^\circ,\rK^{p\circ})\to\rM_\fp(\rV^\circ,\rK^{p\circ})$ induces a purely inseparable morphism onto its image which is $\rM^\bullet_\fp(\rV^\circ,\rK^{p\circ})$.

  \item The morphism $\nu'\colon\rN^\prime_\fp(\rV^\circ,\rK^{p\circ})\to\rM^\prime_\fp(\rV^\circ,\rK^{p\circ})$ is purely inseparable.

  \item The morphism $\nu\colon\rN_\fp(\rV^\circ,\rK^{p\circ})\to\rN^\prime_\fp(\rV^\circ,\rK^{p\circ})$ is the blow-up along $\rN^\prime_\fp(\rV^\circ,\rK^{p\circ})_{\sp}$.\footnote{Note that blow-up along a zero-dimensional closed subscheme $Z$ of a regular scheme depends only on the underlying closed subset of $Z$.}
\end{enumerate}
\end{proposition}

\begin{proof}
For (1), it suffices to show that for every algebraically closed field $\kappa$ containing $\dF_p^\Phi$, $\mu(\kappa)$ is an isomorphism from $\rN_\fp(\rV^\circ,\rK^{p\circ})(\kappa)$ to $\rM^\bullet_\fp(\rV^\circ,\rK^{p\circ})(\kappa)$.

We first show that the image of $\mu(\kappa)$ is contained in $\rM^\bullet_\fp(\rV^\circ,\rK^{p\circ})(\kappa)$. Take a point $y=(A_0,\lambda_0,\eta_0^p;A,\lambda,\eta^p;A^\prime,\lambda^\prime,\eta^{p\prime};\delta)\in\rN_\fp(\rV^\circ,\rK^{p\circ})(\kappa)$. By Lemma~\ref{le:inverse_isogeny}(2,4) and the relation $\lambda=\delta^\vee\circ\lambda^\prime\circ\delta$, we know that $\delta_{*,\tau}\colon\rH_1^\dr(A/\kappa)_\tau\to\rH_1^\dr(A^\prime/\kappa)_\tau$ is an isomorphism if $\tau\neq\tau_\infty^\tc$; and $\Ker\delta_{*,\tau_\infty^\tc}$ has dimension $1$. Moreover, since $\lambda'$ is $p$-principal, we have $\Ker\delta_{*,\tau_\infty^\tc}=\rH_1^\dr(A/\kappa)_{\tau_\infty}^\perp$. By the signature condition again, we also have $\Ker\delta_{*,\tau_\infty^\tc}\subseteq\omega_{A^\vee,\tau_\infty^\tc}$. Thus, $\mu(y)$ belongs to $\rM^\bullet_\fp(\rV^\circ,\rK^{p\circ})(\kappa)$ by Definition~\ref{de:ns_moduli_scheme1}.

It remains to construct an inverse to $\mu(\kappa)\colon\rN_\fp(\rV^\circ,\rK^{p\circ})(\kappa)\to\rM^\bullet_\fp(\rV^\circ,\rK^{p\circ})(\kappa)$. Take a point $x=(A_0,\lambda_0,\eta_0^p;A,\lambda,\eta^p)\in\rM^\bullet_\fp(\rV^\circ,\rK^{p\circ})(\kappa)$. Write $\tilde\rH_1^\dr(A/\kappa)_{\tau_\infty}^\perp$ the preimage of $\rH_1^\dr(A/\kappa)_{\tau_\infty}^\perp$ under the reduction map $\cD(A)_{\tau_\infty^\tc}\to\rH^\dr_1(A/\kappa)_{\tau_\infty^\tc}$. As $\langle\rH_1^\dr(A/\kappa)_{\tau_\infty}^\perp,\rH^\dr_1(A/\kappa)_{\tau_\infty}\rangle_{\lambda,\tau_\infty^\tc}=0$, we have $\cD(A)_{\tau_\infty}^\vee=p^{-1}\tilde\rH_1^\dr(A/\kappa)_{\tau_\infty}^\perp$. Now we put $\cD_{A^\prime,\tau}\coloneqq\cD(A)_\tau$ for $\tau\neq\tau_\infty^\tc$, and $\cD_{A^\prime,\tau_\infty^\tc}\coloneqq p^{-1}\tilde\rH_1^\dr(A/\kappa)_{\tau_\infty}^\perp$. We claim that $\cD_{A^\prime}\coloneqq\bigoplus_{\tau\in\Sigma_\infty}\cD_{A^\prime,\tau}$ is a Dieudonn\'{e} module, which amounts to the inclusions $\tF\cD_{A^\prime,\tau_\infty^\tc}\subseteq\cD_{A^\prime,\tau_\infty}$ and $\tV\cD_{A^\prime,\tau_\infty^\tc}\subseteq\cD_{A^\prime,\tau_\infty}$. The first one follows from the relation $\tF(\rH^1_\dr(A/\kappa)^\perp_{\tau^\infty})\subseteq\tF\omega_{A^\vee,\tau_\infty^\tc}=0$ in which the first inclusion is due to Definition \ref{de:ns_moduli_scheme1}; and the second one is equivalent to the first one as $\cD_{A^\prime,\tau_\infty^\tc}$ and $\cD_{A^\prime,\tau_\infty}$ are integrally dual under $\langle\;,\;\rangle_{\lambda,\tau_\infty^\tc}^\cris$. Then by the Dieudonn\'{e} theory, there is an $O_F$-abelian scheme $A^\prime$ over $\kappa$ with $\cD(A^\prime)_\tau=\cD_{A^\prime,\tau}$ for every $\tau\in\Sigma_\infty$, and an $O_F$-linear isogeny $\delta\colon A\to A^\prime$ inducing the inclusion of Dieudonn\'{e} modules $\cD(A)\subseteq\cD(A^\prime)$. By Lemma \ref{le:inverse_isogeny}(2,4), the $O_F$-abelian scheme $A^\prime$ has signature type $3\Phi-2\tau_\infty+2\tau_\infty^\tc$. Let $\lambda^\prime$ be the unique quasi-polarization of $A^\prime$ satisfying $\lambda=\delta^\vee\circ\lambda^\prime\circ\delta$, which is $p$-principal as $\cD_{A^\prime,\tau_\infty}=\cD_{A^\prime,\tau_\infty^\tc}^\vee$. Finally, we let $\eta^{p\prime}$ be the map sending $v\in\rV^\circ\otimes_\dQ\dA^{\infty,p}$ to $\delta_*\circ\eta^p(v)$. Thus, we obtain an object $(A_0,\lambda_0,\eta_0^p;A,\lambda,\eta^p;A^\prime,\lambda^\prime,\eta^{p\prime};\delta)\in\rN_\fp(\rV^\circ,\rK^{p\circ})(\kappa)$. It is straightforward to check that such assignment gives rise to an inverse of $\mu(\kappa)$.

We now consider (2) and (3) simultaneously. Let $\rN_\fp(\rV^\circ,\rK^{p\circ})_{\sp}$ be the inverse image of $\rM^\prime_\fp(\rV^\circ,\rK^{p\circ})_{\sp}$ under $\mu'$. By Lemma \ref{le:three}(1) below, the induced morphism
\[
\mu'\colon\rN_\fp(\rV^\circ,\rK^{p\circ})\setminus\rN_\fp(\rV^\circ,\rK^{p\circ})_{\sp}\to
\rM^\prime_\fp(\rV^\circ,\rK^{p\circ})\setminus\rM^\prime_\fp(\rV^\circ,\rK^{p\circ})_{\sp}
\]
is purely inseparable. Thus, the induced morphism
\[
\nu\colon\rN_\fp(\rV^\circ,\rK^{p\circ})\setminus\rN_\fp(\rV^\circ,\rK^{p\circ})_{\sp}\to
\rN^\prime_\fp(\rV^\circ,\rK^{p\circ})\setminus\rN^\prime_\fp(\rV^\circ,\rK^{p\circ})_{\sp}
\]
is an isomorphism, and the induced morphism
\[
\nu'\colon\rN^\prime_\fp(\rV^\circ,\rK^{p\circ})\setminus\rN^\prime_\fp(\rV^\circ,\rK^{p\circ})_{\sp}\to
\rM^\prime_\fp(\rV^\circ,\rK^{p\circ})\setminus\rM^\prime_\fp(\rV^\circ,\rK^{p\circ})_{\sp}
\]
is purely inseparable. Since $\rN_\fp(\rV^\circ,\rK^{p\circ})$ is quasi-projective, $\nu$ is projective. Thus, $\nu$ is a projective birational morphism, which has to be the blow-up along a subset $Z$ of $\rN^\prime_\fp(\rV^\circ,\rK^{p\circ})_{\sp}$ (see, for example, \cite{Liu02}*{Theorem~8.1.24}). Now we take a point $x'$ of $\rM^\prime_\fp(\rV^\circ,\rK^{p\circ})_{\sp}$ with the residue field $\kappa$, which is a finite extension of $\dF_p^\Phi$. Since $\nu'$ is a finite morphism, the inverse image of $x'$ consists of finitely many points $y'_1,\dots,y'_n$ of $\rN^\prime_\fp(\rV^\circ,\rK^{p\circ})_{\sp}$ with residue fields $\kappa_1,\dots,\kappa_n$, respectively. By Lemma \ref{le:three}(2) below, the residue field extension $\kappa_i/\kappa$ is trivial for every $1\leq i\leq n$; and moreover, $Z$ has nonempty intersection with $\{y'_1,\dots,y'_n\}$. Thus, $\mu^{\prime-1}(x')$ has cardinality at least $|\kappa|+n$. But we know that $\mu^{\prime-1}(x')$ has cardinality exactly $|\kappa|+1$. Therefore, we must have $n=1$. We immediately have both (2) and (3).
\end{proof}

\begin{lem}\label{le:three}
Consider an element $x'\in\rM^\prime_\fp(\rV^\circ,\rK^{p\circ})(\kappa)$ for some $\rK^{p\circ}\in\fK(\rV^\circ)^p$ and a perfect field $\kappa$ containing $\dF_p^\Phi$. We have
\begin{enumerate}
  \item If the image of $x'$ is not special, then $\mu^{\prime-1}(x')$ is a singleton.

  \item If the image of $x'$ is special, then $\mu^{\prime-1}(x')$ is isomorphic to $\dP^1(\kappa)$.
\end{enumerate}
\end{lem}

\begin{proof}
Write $x'=(A_0,\lambda_0,\eta_0^p;A^\prime,\lambda^\prime,\eta^{p\prime})$. By the Dieudonn\'{e} theory and Lemma~\ref{le:inverse_isogeny}(2,4), we see that $\mu^{\prime-1}(x')$ is bijective to Dieudonn\'{e} submodules $\cD_A\subseteq\cD(A^\prime)$ satisfying $\cD_{A,\tau}=\cD(A^\prime)_\tau$ for $\tau\neq\tau_\infty^\tc$, and that $\cD(A^\prime)_{\tau_\infty^\tc}/\cD_{A,\tau_\infty^\tc}$ is a vector space over $\kappa$ of dimension $1$. This amounts to the subspaces of $\rH^\dr_1(A^\prime/\kappa)_{\tau_\infty^\tc}$ of dimension $2$ containing $\tF\rH^\dr_1(A^\prime/\kappa)_{\tau_\infty}+\tV\rH^\dr_1(A^\prime/\kappa)_{\tau_\infty}$. Since both $\tF\rH^\dr_1(A^\prime/\kappa)_{\tau_\infty}$ and $\tV\rH^\dr_1(A^\prime/\kappa)_{\tau_\infty}$ have dimension $1$, the lemma follows by the definition of special points.
\end{proof}

\begin{remark}
In fact, one can show that $\mu$ induces an isomorphism from $\rN_\fp(\rV^\circ,\obj)$ to $\rM^\bullet_\fp(\rV^\circ,\obj)$; and $\nu'$ is purely inseparable of degree $p$. But we do not need these facts.
\end{remark}

\fi

\subsection{Functoriality under special morphisms}
\label{ss:ns_functoriality}

In this subsection, we study the behavior of various moduli schemes under the special morphisms, which is closely related to the Rankin--Selberg motives for $\GL_n\times\GL_{n+1}$. We start from the datum $(\rV^\circ_n,\{\Lambda^\circ_{n,\fq}\}_{\fq\mid p})$ as in the beginning of \S\ref{ss:ns_moduli_scheme}, but with $\rV^\circ_n$ of rank $n\geq2$. (See Remark \ref{re:ns_functoriality} below for the case $n=1$.) We then have the induced datum
\[
(\rV^\circ_{n+1},\{\Lambda^\circ_{n+1,\fq}\}_{\fq\mid p})\coloneqq((\rV^\circ_n)_\sharp,\{(\Lambda^\circ_{n,\fq})_\sharp\}_{\fq\mid p})
\]
of rank $n+1$ by Definition \ref{de:hermitian_space}. For $N\in\{n,n+1\}$, we let $\rK^\circ_{N,\fq}$ be the stabilizer of $\Lambda^\circ_{N,\fq}$, and put $\rK^\circ_{N,p}\coloneqq\prod_{\fq\mid p}\rK^\circ_{N,\fq}$. Recall the category $\fK(\rV^\circ_n)_\sp^p$ and functors $\obj_\flat,\obj_\sharp$ from Definition \ref{de:neat_category}. To unify notation, we put $\obj_n\coloneqq\obj_\flat$ and $\obj_{n+1}\coloneqq\obj_\sharp$. Similar to the case of smooth moduli schemes considered in \S\ref{ss:qs_functoriality}, there are five stages of functoriality we will consider.

The first stage concerns Shimura varieties.

\begin{notation}\label{no:ns_uniformization_indefinite}
We choose an indefinite uniformization datum $(\rV'_n,\tj_n,\{\Lambda'_{n,\fq}\}_{\fq\mid p})$ for $\rV^\circ_n$ as in Definition \ref{de:ns_uniformization_data}. Put $\rV'_{n+1}\coloneqq(\rV'_n)_\sharp$, $\tj_{n+1}\coloneqq(\tj_n)_\sharp$, and $\Lambda'_{n+1,\fq}\coloneqq(\Lambda'_{n,\fq})_\sharp$. Then $(\rV'_{n+1},\tj_{n+1},\{\Lambda'_{n+1,\fq}\}_{\fq\mid p})$ is an indefinite uniformization datum for $\rV^\circ_{n+1}$. For $N\in\{n,n+1\}$, we let $\rK'_{N,\fq}$ be the stabilizer of $\Lambda'_{N,\fq}$, and put $\rK'_{N,p}\coloneqq\prod_{\fq\mid p}\rK'_{N,\fq}$.
\end{notation}

We obtain a morphism
\begin{align*}
\sh'_\uparrow\colon\Sh(\rV'_n,\tj_n\obj_n\rK'_{n,p})\to\Sh(\rV'_{n+1},\tj_{n+1}\obj_{n+1}\rK'_{n+1,p})
\end{align*}
in $\Fun(\fK(\rV^\circ_n)_\sp^p,\Sch_{/F})$.

For the second stage of functoriality, we have a morphism
\begin{align}\label{eq:ns_functoriality_moduli_scheme}
\bm_\uparrow\colon \bM_\fp(\rV^\circ_n,\obj_n)\to \bM_\fp(\rV^\circ_{n+1},\obj_{n+1})
\end{align}
in $\Fun(\fK(\rV^\circ_n)_\sp^p\times\fT,\Sch_{/\dZ_p^\Phi})_{/\bT_\fp}$ sending an object
$(A_0,\lambda_0,\eta_0^p;A,\lambda,\eta^p)\in\bM_\fp(\rV^\circ_n,\rK^{p\circ}_n)(S)$ to the object
$(A_0,\lambda_0,\eta_0^p;A\times A_0,\lambda\times\lambda_0,\eta^p\oplus(\id_{A_0})_*)\in\bM_\fp(\rV^\circ_{n+1},\rK^{p\circ}_{n+1})(S)$. It is clear that $\bm_\uparrow$ restricts to three morphisms
\begin{align}\label{eq:ns_functoriality_moduli_scheme1}
\begin{dcases}
\rm^\circ_\uparrow \colon \rM^\circ_\fp(\rV^\circ_n,\obj_n)\to \rM^\circ_\fp(\rV^\circ_{n+1},\obj_{n+1}),\\
\rm^\dag_\uparrow \colon \rM^\dag_\fp(\rV^\circ_n,\obj_n)\to \rM^\dag_\fp(\rV^\circ_{n+1},\obj_{n+1}),\\
\rm^\bullet_\uparrow \colon \rM^\bullet_\fp(\rV^\circ_n,\obj_n)\to \rM^\bullet_\fp(\rV^\circ_{n+1},\obj_{n+1}).
\end{dcases}
\end{align}

Moreover, we have the following commutative diagram
\begin{align}\label{eq:ns_functoriality_shimura1}
\xymatrix{
\bM^\eta_\fp(\rV^\circ_{n+1},\obj_{n+1}) \ar[rr]^-{\eqref{eq:ns_moduli_scheme_shimura}}
&&  \Sh(\rV'_{n+1},\tj_{n+1}\obj_{n+1}\rK'_{n+1,p})\times_{\Spec{F}}\bT^\eta_\fp \\
\bM^\eta_\fp(\rV^\circ_n,\obj_n) \ar[rr]^-{\eqref{eq:ns_moduli_scheme_shimura}}\ar[u]^-{\bm^\eta_\uparrow}
&& \Sh(\rV'_n,\tj_n\obj_n\rK'_{n,p})\times_{\Spec{F}}\bT^\eta_\fp \ar[u]_-{\sh'_\uparrow\times\id}
}
\end{align}
in $\Fun(\fK(\rV^\circ_n)_\sp^p\times\fT,\Sch_{/\dQ_p^\Phi})_{/\bT^\eta_\fp}$.

At the third stage of functoriality, we study the basic correspondence diagram \eqref{eq:ns_link} for $N=n,n+1$ under the special morphisms. We will complete a commutative diagram in $\Fun(\fK(\rV^\circ_n)_\sp^p\times\fT,\Sch_{/\dF_p^\Phi})_{/\rT_\fp}$ as follows
\begin{align}\label{eq:ns_functoriality}
\resizebox{0.45\hsize}{!}{
\rotatebox{-90}{
\xymatrix{
\rS^\circ_\fp(\rV^\circ_{n+1},\obj_{n+1}) & & & \rB^\circ_\fp(\rV^\circ_{n+1},\obj_{n+1}) \ar[lll]_-{\pi^\circ_{n+1}}\ar[rrr]^-{\iota^\circ_{n+1}} & & & \rM^\circ_\fp(\rV^\circ_{n+1},\obj_{n+1}) & & \\
& \rS^\dag_\fp(\rV^\circ_{n+1},\obj_{n+1}) \ar[ul]_-{\rs^{\dag\circ}_{n+1}}\ar[rd]^-{\rs^{\dag\bullet}_{n+1}} & & & \rB^\dag_\fp(\rV^\circ_{n+1},\obj_{n+1}) \ar[ul]_-{\rb^{\dag\circ}_{n+1}}\ar[rd]^-{\rb^{\dag\bullet}_{n+1}}\ar[lll]_-{\pi^\dag_{n+1}}\ar[rrr]^-{\iota^\dag_{n+1}} & & & \rM^\dag_\fp(\rV^\circ_{n+1},\obj_{n+1})  \ar[ul]_-{\rm^{\dag\circ}_{n+1}}\ar[rd]^-{\rm^{\dag\bullet}_{n+1}} &  \\
&& \rS^\bullet_\fp(\rV^\circ_{n+1},\obj_{n+1}) & & & \rB^\bullet_\fp(\rV^\circ_{n+1},\obj_{n+1}) \ar[lll]_-{\pi^\bullet_{n+1}}\ar[rrr]^-{\iota^\bullet_{n+1}} & & & \rM^\bullet_\fp(\rV^\circ_{n+1},\obj_{n+1})  \\
\\
& \rS^\dag_\fp(\rV^\circ_n,\obj)_\sp \ar[rd]^-{\rs^{\dag\bullet}_\sp}\ar[uuu]_-{\rs^\dag_\uparrow}\ar[ddd]^-{\rs^\dag_\downarrow} & & & \rB^\dag_\fp(\rV^\circ_n,\obj)_\sp \ar[rd]^-{\rb^{\dag\bullet}_\sp} \ar[lll]|!{[lluu];[lld]}\hole_-{\pi^\dag_\sp} \ar[ddd]|!{[dr];[dll]}\hole^-{\rb^\dag_\downarrow} \ar[uuu]|!{[uur];[uull]}\hole_-{\rb^\dag_\uparrow} & & & & \\
&& \rS^\bullet_\fp(\rV^\circ_n,\obj)_\sp \ar[uuu]_-{\rs^\bullet_\uparrow}\ar[ddd]^-{\rs^\bullet_\downarrow} & & & \rB^\bullet_\fp(\rV^\circ_n,\obj)_\sp \ar[uuu]_-{\rb^\bullet_\uparrow}\ar[ddd]^-{\rb^\bullet_\downarrow}\ar[lll]_-{\pi^\bullet_\sp} & & & \\
\rS^\circ_\fp(\rV^\circ_n,\obj_n) \ar[uuuuuu]_-{\rs^\circ_\uparrow} & & & \rB^\circ_\fp(\rV^\circ_n,\obj_n) \ar[rrr]|!{[ruu];[rd]}\hole|!{[rru];[rrdd]}\hole^-{\iota^\circ_n} \ar[lll]|!{[lu];[ldd]}\hole|!{[lluu];[lld]}\hole_-{\pi^\circ_n} \ar[uuuuuu]|!{[ul];[urr]}\hole|!{[uull];[uur]}\hole|!{[uuuul];[uuuurr]}\hole|!{[uuuuull];[uuuuur]}\hole_-{\rb^\circ_\uparrow} & & & \rM^\circ_\fp(\rV^\circ_n,\obj_n)
\ar[uuuuuu]|!{[uuuul];[uuuurr]}\hole|!{[uuuuull];[uuuuur]}\hole_-{\rm^\circ_\uparrow} & & \\
& \rS^\dag_\fp(\rV^\circ_n,\obj_n) \ar[ul]_-{\rs^{\dag\circ}_n}\ar[rd]^-{\rs^{\dag\bullet}_n} & & & \rB^\dag_\fp(\rV^\circ_n,\obj_n) \ar[ul]_-{\rb^{\dag\circ}_n}\ar[rd]^-{\rb^{\dag\bullet}_n}\ar[rrr]|!{[ruu];[rd]}\hole^-{\iota^\dag_n} \ar[lll]|!{[lluu];[lld]}\hole_-{\pi^\dag_n} & & & \rM^\dag_\fp(\rV^\circ_n,\obj_n) \ar[ul]_-{\rm^{\dag\circ}_n}\ar[rd]^-{\rm^{\dag\bullet}_n} \ar[uuuuuu]|!{[uuuuull];[uuuuur]}\hole_-{\rm^\dag_\uparrow} &  \\
&& \rS^\bullet_\fp(\rV^\circ_n,\obj_n) & & & \rB^\bullet_\fp(\rV^\circ_n,\obj_n) \ar[lll]_-{\pi^\bullet_n}\ar[rrr]^-{\iota^\bullet_n} & & & \rM^\bullet_\fp(\rV^\circ_n,\obj_n) \ar[uuuuuu]_-{\rm^\bullet_\uparrow}
}
}
}
\end{align}
in which the bottom (resp. top) layer is the basic correspondence diagram \eqref{eq:ns_link} for $\rM_\fp(\rV^\circ_n,\obj_n)$ (resp.\ $\rM_\fp(\rV^\circ_{n+1},\obj_{n+1})$).

First, we consider the basic correspondences on the balloon strata, that is, the back layer of the diagram \eqref{eq:ns_functoriality}.

We define $\rs^\circ_\uparrow\colon\rS^\circ_\fp(\rV^\circ_n,\obj_n)\to\rS^\circ_\fp(\rV^\circ_{n+1},\obj_{n+1})$ to be the morphism sending an object
\[
(A_0,\lambda_0,\eta_0^p;A^\circ,\lambda^\circ,\eta^{p\circ})\in\rS^\circ_\fp(\rV^\circ_n,\rK^{p\circ}_n)(S)
\]
to the object
\[
(A_0,\lambda_0,\eta_0^p;A^\circ\times A_0,\lambda^\circ\times\lambda_0,\eta^{p\circ}\oplus(\id_{A_0})_*)\in\rS^\circ_\fp(\rV^\circ_{n+1},\rK^{p\circ}_{n+1})(S).
\]

\begin{remark}\label{re:ns_functoriality_set}
The canonical inclusions
\[
\rV^\circ_n\hookrightarrow\rV^\circ_{n+1},\quad\{\Lambda^\circ_{n,\fq}\hookrightarrow\Lambda^\circ_{n+1,\fq}\}_{\fq\mid p}
\]
induce a morphism
\begin{align*}
\sh^\circ_\uparrow\colon\Sh(\rV^\circ_n,\obj_n\rK^\circ_{n,p})\to\Sh(\rV^\circ_{n+1},\obj_{n+1}\rK^\circ_{n+1,p})
\end{align*}
in $\Fun(\fK(\rV^\circ_n)_\sp^p,\Set)$. It is clear that the following diagram
\[
\xymatrix{
\rS^\circ_\fp(\rV^\circ_{n+1},\obj_{n+1})(\ol\dF_p) \ar[r]^-{\upsilon^\circ_{n+1}} & \Sh(\rV^\circ_{n+1},\obj_{n+1}\rK^\circ_{n+1,p})\times\rT_\fp(\ol\dF_p) \\
\rS^\circ_\fp(\rV^\circ_n,\obj_n)(\ol\dF_p) \ar[r]^-{\upsilon^\circ_n}\ar[u]^-{\rs^\circ_\uparrow(\ol\dF_p)} & \Sh(\rV^\circ_n,\obj_n\rK^\circ_{n,p})\times\rT_\fp(\ol\dF_p) \ar[u]_-{\sh^\circ_\uparrow\times\id_{\rT_\fp(\ol\dF_p)}}
}
\]
in $\Fun(\fK(\rV^\circ_n)_\sp^p,\Set)_{/\rT_\fp(\ol\dF_p)}$ commutes, where $\upsilon^\circ_{n+1}$ and $\upsilon^\circ_n$ are uniformization maps in Construction \ref{cs:ns_uniformization1}.
\end{remark}

We define $\rb^\circ_\uparrow\colon\rB^\circ_\fp(\rV^\circ_n,\obj_n)\to\rB^\circ_\fp(\rV^\circ_{n+1},\obj_{n+1})$ to be the morphism sending an object
\[
(A_0,\lambda_0,\eta_0^p;A,\lambda,\eta^p;A^\circ,\lambda^\circ,\eta^{p\circ};\beta)\in\rB^\circ_\fp(\rV^\circ_n,\rK^{p\circ}_n)(S)
\]
to the object
\[
(A_0,\lambda_0,\eta_0^p;A\times A_0,\lambda\times\lambda_0,\eta^p\oplus(\id_{A_0})_*;A^\circ\times A_0,\lambda^\circ\times\lambda_0,\eta^{p\circ}\oplus(\id_{A_0})_*;\beta\times\id_{A_0})
\in\rB^\circ_\fp(\rV^\circ_{n+1},\rK^{p\circ}_{n+1})(S).
\]

Second, we consider the basic correspondences on the ground strata, that is, the front layer of the diagram \eqref{eq:ns_functoriality}.

\begin{definition}\label{de:ns_definite_special2}
We define a functor
\begin{align*}
\rS^\bullet_\fp(\rV^\circ_n,\obj)_\sp\colon\fK(\rV^\circ_n)_\sp^p\times\fT &\to\sfP\Sch'_{/\dF^\Phi_p} \\
\rK^{p\circ} &\mapsto \rS^\bullet_\fp(\rV^\circ_n,\rK^{p\circ})_\sp
\end{align*}
such that for every $S\in\Sch'_{/\dF^\Phi_p}$, $\rS^\bullet_\fp(\rV^\circ_n,\rK^{p\circ})_\sp(S)$ is the set of equivalence classes of decuples $(A_0,\lambda_0,\eta_0^p;A^\bullet,\lambda^\bullet,\eta^{p\bullet};A^\bullet_\natural,\lambda^\bullet_\natural,\eta^{p\bullet}_\natural;
\delta^\bullet)$, where
\begin{itemize}[label={\ding{109}}]
  \item $(A_0,\lambda_0,\eta_0^p;A^\bullet,\lambda^\bullet,\eta^{p\bullet})$ is an element in $\rS^\bullet_\fp(\rV^\circ_n,\rK^{p\circ}_n)(S)$;

  \item $(A_0,\lambda_0,\eta_0^p;A^\bullet_\natural,\lambda^\bullet_\natural,\eta^{p\bullet}_\natural)$ is an element in $\rS^\bullet_\fp(\rV^\circ_{n+1},\rK^{p\circ}_{n+1})(S)$; and

  \item $\delta^\bullet\colon A^\bullet\times A_0\to A^\bullet_\natural$ is an $O_F$-linear quasi-$p$-isogeny (Definition \ref{de:p_quasi}) such that
  \begin{enumerate}[label=(\alph*)]
    \item $\Ker\delta^\bullet[p^\infty]$ is contained in $(A^\bullet\times A_0)[\fp]$;

    \item we have $\lambda^\bullet\times\varpi\lambda_0=\delta^{\bullet\vee}\circ\lambda^\bullet_\natural\circ\delta^\bullet$; and

    \item the $\rK^p_{n+1}$-orbit of maps $v\mapsto\delta^\bullet_*\circ(\eta^{p\bullet}\oplus(\id_{A_0})_*)(v)$ for $v\in\rV^\circ_\sharp\otimes_\dQ\dA^{\infty,p}$ coincides with $\eta^{p\bullet}_\natural$.
  \end{enumerate}
\end{itemize}
The equivalence relation and the action of morphisms in $\fK(\rV^\circ_n)_\sp^p\times\fT$ are defined similarly as in Definition \ref{de:qs_basic_correspondence}.
\end{definition}

We clearly have the forgetful morphism
\begin{align*}
\rS^\bullet_\fp(\rV^\circ_n,\obj)_\sp\to\rT_\fp
\end{align*}
in $\Fun(\fK(\rV^\circ_n)_\sp^p\times\fT,\sfP\Sch'_{/\dF_p^\Phi})$, which is represented by finite and \'{e}tale schemes. By definition, we have the two forgetful morphisms
\[
\rs^\bullet_\downarrow\colon\rS^\bullet_\fp(\rV^\circ_n,\obj)_\sp\to\rS^\bullet_\fp(\rV^\circ_n,\obj_n),\quad
\rs^\bullet_\uparrow\colon\rS^\bullet_\fp(\rV^\circ_n,\obj)_\sp\to\rS^\bullet_\fp(\rV^\circ_{n+1},\obj_{n+1})
\]
in $\Fun(\fK(\rV^\circ_n)_\sp^p\times\fT,\Sch_{/\dF_p^\Phi})_{/\rT_\fp}$.

\begin{lem}\label{le:ns_functoriality_source}
We have the following properties concerning $\rs^\bullet_{\downarrow}$.
\begin{enumerate}
  \item When $n$ is even, $\rs^\bullet_{\downarrow}$ is an isomorphism, and the morphism
      \[
      \rs^\bullet_{\uparrow}\circ \rs_{\downarrow}^{\bullet-1}\colon\rS^\bullet_{\fp}(\rV^\circ_n, \obj_n)\to\rS^\bullet_{\fp}(\rV^\circ_{n+1},\obj_{n+1})
      \]
      is given by the assignment
      \[
      (A_0,\lambda_0,\eta_0^p;A^\bullet,\lambda^\bullet,\eta^{p\bullet})\mapsto
      (A_0,\lambda_0,\eta_0^p;A^\bullet\times A_0,\lambda^\bullet\times\varpi\lambda_0,\eta^{p\bullet}\times(\id_{A_0})_*).
      \]

  \item When $n$ is odd, $\rs^\bullet_{\downarrow}$ is finite \'etale of degree $p+1$.
\end{enumerate}
\end{lem}

\begin{proof}
The proof is very similar to Lemma \ref{le:qs_functoriality_source}, which we leave to the readers.
\end{proof}

\begin{definition}\label{de:ns_base_special2}
We define $\rB^\bullet_\fp(\rV^\circ_n,\obj)_\sp$ to be the fiber product indicated in the following Cartesian diagram
\[
\xymatrix{
\rS^\bullet_\fp(\rV^\circ_n,\obj)_\sp \ar[d]^-{\rs^\bullet_\downarrow} &&
\rB^\bullet_\fp(\rV^\circ_n,\obj)_\sp \ar[ll]_-{\pi^\bullet_\sp}\ar[d]_-{\rb^\bullet_\downarrow}  \\
\rS^\bullet_\fp(\rV^\circ_n,\obj_n) &&
\rB^\bullet_\fp(\rV^\circ_n,\obj_n) \ar[ll]_-{\pi^\bullet_n}
}
\]
in $\Fun(\fK(\rV^\circ_n)_\sp^p\times\fT,\Sch_{/\dF_p^\Phi})_{/\rT_\fp}$. We define $\rb^\bullet_\uparrow\colon\rB^\bullet_\fp(\rV^\circ_n,\obj)_\sp\to\rB^\bullet_\fp(\rV^\circ_{n+1},\obj_{n+1})$ to be the morphism sending an object
\[
((A_0,\lambda_0,\eta_0^p;A,\lambda,\eta^p;A^\bullet,\lambda^\bullet,\eta^{p\bullet};\gamma),(A_0,\lambda_0,\eta_0^p;
A^\bullet,\lambda^\bullet,\eta^{p\bullet};A^\bullet_\natural,\lambda^\bullet_\natural,\eta^{p\bullet}_\natural;\delta^\bullet))
\in\rB^\bullet_\fp(\rV^\circ_n,\rK^{p\circ})_\sp(S)
\]
to $(A_0,\lambda_0,\eta_0^p;A\times A_0,\lambda\times\lambda_0,\eta^p\oplus(\id_{A_0})_*;
A^\bullet_\natural,\lambda^\bullet_\natural,\eta^{p\bullet}_\natural;\delta^\bullet\circ(\gamma\times\id_{A_0}))$, which is an object of $\rB^\bullet_\fp(\rV^\circ_{n+1},\rK^{p\circ}_{n+1})(S)$ by a similar argument of Lemma \ref{le:qs_functoriality}.
\end{definition}

We have the following result.

\begin{proposition}\label{pr:ns_functoriality2}
When $n$ is odd, the square
\[
\xymatrix{
\rB^\bullet_\fp(\rV^\circ_{n+1},\obj_{n+1}) \ar[rr]^-{\iota^\bullet_{n+1}} &&  \rM^\bullet_\fp(\rV^\circ_{n+1},\obj_{n+1}) \\
\rB^\bullet_\fp(\rV^\circ_n,\obj)_\sp \ar[rr]^-{\iota^\bullet_n\circ\rb^\bullet_\downarrow} \ar[u]^-{\rb^\bullet_\uparrow}
&& \rM^\bullet_\fp(\rV^\circ_n,\obj_n) \ar[u]_-{\rm^\bullet_\uparrow}
}
\]
extracted from the diagram \eqref{eq:ns_functoriality} is Cartesian.
\end{proposition}

\begin{proof}
The proof is very similar to Proposition \ref{pr:qs_functoriality}, which we leave to the readers.
\end{proof}

Third, we consider the basic correspondences on the link strata, that is, the middle (vertical) layer of the diagram \eqref{eq:ns_functoriality}.

\begin{definition}\label{de:ns_definite_special3}
We define $\rS^\dag_\fp(\rV^\circ_n,\obj)_\sp$ to be the fiber product indicated in the following Cartesian diagram
\[
\xymatrix{
\rS^\dag_\fp(\rV^\circ_n,\obj)_\sp \ar[rr]^-{\rs^{\dag\bullet}_\sp}\ar[d]_-{\rs^\dag_\downarrow} && \rS^\bullet_\fp(\rV^\circ_n,\obj)_\sp \ar[d]^-{\rs^\bullet_\downarrow} \\
\rS^\dag_\fp(\rV^\circ_n,\obj_n) \ar[rr]^-{\rs^{\dag\bullet}_n} && \rS^\bullet_\fp(\rV^\circ_n,\obj_n)
}
\]
in $\Fun(\fK(\rV^\circ_n)_\sp^p\times\fT,\Sch_{/\dF_p^\Phi})_{/\rT_\fp}$. By Lemma \ref{le:ns_functoriality_source}, we know that $\rs^\dag_\downarrow$ is an isomorphism (resp. finite \'etale of degree $p+1$) when $n$ is even (resp.\ odd). We define $\rs^\dag_\uparrow\colon\rS^\dag_\fp(\rV^\circ_n,\obj)_\sp\to\rS^\dag_\fp(\rV^\circ_{n+1},\obj_{n+1})$ to be the morphism sending an object
\[
((A_0,\lambda_0,\eta_0^p;A^\circ,\lambda^\circ,\eta^{p\circ};A^\bullet,\lambda^\bullet,\eta^{p\bullet};\psi),
(A_0,\lambda_0,\eta_0^p;A^\bullet,\lambda^\bullet,\eta^{p\bullet};A^\bullet_\natural,\lambda^\bullet_\natural,\eta^{p\bullet}_\natural;
\delta^\bullet))
\in\rS^\dag_\fp(\rV^\circ_n,\rK^{p\circ})_\sp(S)
\]
to the object
\[
(A_0,\lambda_0,\eta_0^p;A^\circ\times A_0,\lambda^\circ\times\lambda_0,\eta^{p\circ}\oplus(\id_{A_0})_*;
A^\bullet_\natural,\lambda^\bullet_\natural,\eta^{p\bullet}_\natural;\delta^\bullet\circ(\psi\times\id_{A_0}))
\in\rS^\dag_\fp(\rV^\circ_{n+1},\rK^{p\circ}_{n+1})(S).
\]
\end{definition}

\begin{lem}\label{le:ns_functoriality_source_dagger}
We have
\begin{enumerate}
  \item When $n$ is even, the square
     \[
     \xymatrix{
     \rS^\dag_\fp(\rV^\circ_{n+1},\obj_{n+1})  \ar[r]^-{\rs^{\dag\bullet}_{n+1}} & \rS^\bullet_\fp(\rV^\circ_{n+1},\obj_{n+1}) \\
     \rS^\dag_\fp(\rV^\circ_n,\obj)_\sp \ar[r]^-{\rs^{\dag\bullet}_\sp} \ar[u]_-{\rs^\dag_\uparrow}
     & \rS^\bullet_\fp(\rV^\circ_n,\obj)_\sp \ar[u]_-{\rs^\bullet_\uparrow}
     }
     \]
     extracted from \eqref{eq:ns_functoriality} is a Cartesian diagram.

  \item When $n$ is odd, the square
     \[
     \xymatrix{
     \rS^\circ_\fp(\rV^\circ_{n+1},\obj_{n+1})  & \rS^\dag_\fp(\rV^\circ_{n+1},\obj_{n+1}) \ar[l]_-{\rs^{\dag\circ}_{n+1}} \\
     \rS^\circ_\fp(\rV^\circ_n,\obj_n)  \ar[u]_-{\rs^\circ_\uparrow}
     & \rS^\dag_\fp(\rV^\circ_n,\obj)_\sp  \ar[l]_-{\rs^{\dag\circ}_n\circ\rs^\dag_\downarrow} \ar[u]_-{\rs^\dag_\uparrow}
     }
     \]
     extracted from \eqref{eq:ns_functoriality} is a Cartesian diagram.
\end{enumerate}
\end{lem}

\begin{proof}
Let $\rS^\ddag_\fp(\rV^\circ_n,\obj)_\sp$ be the actual fiber product in both cases. Take an object $\rK^{p\circ}\in\fK(\rV^\circ_n)_\sp^p$. We have to show that the natural morphism $\rs^\ddag\colon\rS^\dag_\fp(\rV^\circ_n,\rK^{p\circ})_\sp\to\rS^\ddag_\fp(\rV^\circ_n,\rK^{p\circ})_\sp$ is an isomorphism. Since $\rs^\ddag$ is a morphism of \'{e}tale schemes over $\dF_p^\Phi$, it suffices to show that $\rs^\ddag(\kappa)$ is an isomorphism for every perfect field $\kappa$ containing $\dF_{p}^{\Phi}$.

For (1), by Lemma \ref{le:ns_functoriality_source}(1), an object in $\rS^\ddag_\fp(\rV^\circ_n,\rK^{p\circ})_\sp(S)$ is given by a pair of objects:
\begin{align*}
(A_0,\lambda_0,\eta_0^p;A^\bullet,\lambda^\bullet,\eta^{p\bullet};
A^\bullet\times A_0,\lambda^\bullet\times\varpi\lambda_0,\eta^{p\bullet}\times(\id_{A_0})_*)
&\in\rS^\bullet_\fp(\rV^\circ_n,\rK^{p\circ})_\sp(\kappa), \\
(A_0,\lambda_0,\eta_0^p;A^\circ_\natural,\lambda^\circ_\natural,\eta^{p\circ}_\natural;
A^\bullet\times A_0,\lambda^\bullet\times\varpi\lambda_0,\eta^{p\bullet}\times(\id_{A_0})_*;\psi_\natural)
&\in\rS^\dag_\fp(\rV^\circ_{n+1},\rK^{p\circ}_{n+1})(\kappa).
\end{align*}
Let $A^\circ$ be the cokernel of the kernel of the composite map $A^\circ_\natural\xrightarrow{\psi_\natural}A^\bullet\times A_0\to A^\bullet$, and $\psi\colon A^\circ\to A^\bullet$ the induced map. Let $\lambda^\circ$ be the unique quasi-polarization of $A^\circ$ satisfying $\varpi\cdot\lambda^\circ=\psi^\vee\circ\lambda^\bullet\circ\psi$. Since $\lambda^\circ_\natural$ is $p$-principal and we have $\varpi\cdot\lambda^\circ_\natural=\psi_\natural^\vee\circ(\lambda^\bullet\times\varpi\cdot\lambda_0)\circ\psi_\natural$, the composite map $A^\circ_\natural\xrightarrow{\psi_\natural}A^\bullet\times A_0\to A_0$ splits. Thus, the natural map $A^\circ_\natural\to A^\circ\times A_0$ is an isomorphism. Then $\lambda^\circ$ is $p$-principal, and we obtain an object
\[
(A_0,\lambda_0,\eta_0^p;A^\circ,\lambda^\circ,\eta^{p\circ};A^\bullet,\lambda^\bullet,\eta^{p\bullet};\psi)
\in\rS^\dag_\fp(\rV^\circ_n,\rK^{p\circ}_n)(\kappa)=\rS^\dag_\fp(\rV^\circ_n,\rK^{p\circ})_\sp(\kappa),
\]
where $\eta^{p\circ}$ is chosen such that Definition \ref{de:ns_definite3}(c) is satisfied. In other words, we obtain a morphism from $\rS^\ddag_\fp(\rV^\circ_n,\rK^{p\circ})_\sp(\kappa)$ to $\rS^\dag_\fp(\rV^\circ_n,\rK^{p\circ})_\sp(\kappa)$. It is straightforward to check that it is an inverse to the morphism $\rs^\ddag(\kappa)$.

For (2), an object in $\rS^\ddag_\fp(\rV^\circ_n,\rK^{p\circ})_\sp(\kappa)$ is given by a pair of objects:
\begin{align*}
(A_0,\lambda_0,\eta_0^p;A^\circ,\lambda^\circ,\eta^{p\circ})&\in\rS^\circ_\fp(\rV^\circ_n,\rK^{p\circ}_n)(\kappa), \\
(A_0,\lambda_0,\eta_0^p;A^\circ\times A_0,\lambda^\circ\times\lambda_0,\eta^{p\circ}\times(\id_{A_0})_*;
A^\bullet_\natural,\lambda^\bullet_\natural,\eta^{p\bullet}_\natural;\psi_\natural)
&\in\rS^\dag_\fp(\rV^\circ_{n+1},\rK^{p\circ}_{n+1})(\kappa).
\end{align*}
Let $A^{\bullet\vee}$ be the cokernel of the kernel of the composite map $A^{\bullet\vee}_\natural\xrightarrow{\psi_\natural^\vee}A^{\circ\vee}\times A_0^\vee\to A^{\circ\vee}$, and $\psi^\vee\colon A^{\circ\vee}\to A^{\bullet\vee}$ the induced map. Taking dual, we obtain a map $\psi\colon A^\circ\to A^\bullet$ and an induced map $\delta^\bullet\colon A^\bullet\times A_0\to A^\bullet_\natural$. Let $\lambda^\bullet$ be the unique quasi-polarization of $A^\bullet$ satisfying $\varpi\cdot\lambda^\circ=\psi^\vee\circ\lambda^\bullet\circ\psi$. Since $\lambda^\bullet_\natural$ is $p$-principal and we have $\lambda^\bullet\times\varpi\cdot\lambda_0=\delta^{\bullet\vee}\circ\lambda^\bullet_\natural\circ\delta^\bullet$, we know that $\Ker\lambda^\bullet[p^\infty]$ is contained in $A^\bullet[\fp]$ of rank $p^2$, and we obtain an object
\[
\((A_0,\lambda_0,\eta_0^p;A^\circ,\lambda^\circ,\eta^{p\circ};A^\bullet,\lambda^\bullet,\eta^{p\bullet};\psi),
(A_0,\lambda_0,\eta_0^p;A^\bullet,\lambda^\bullet,\eta^{p\bullet};
A^\bullet_\natural,\lambda^\bullet_\natural,\eta^{p\bullet}_\natural;\delta^\bullet)\)
\in\rS^\dag_\fp(\rV^\circ_n,\rK^{p\circ})_\sp(\kappa),
\]
where $\eta^{p\bullet}$ is chosen such that Definition \ref{de:ns_definite3}(c) is satisfied. In other words, we obtain a morphism from $\rS^\ddag_\fp(\rV^\circ_n,\rK^{p\circ})_\sp(\kappa)$ to $\rS^\dag_\fp(\rV^\circ_n,\rK^{p\circ})_\sp(\kappa)$. It is straightforward to check that it is an inverse to the morphism $\rs^\ddag(\kappa)$.
\end{proof}

\begin{definition}\label{de:ns_base_special3}
We define $\rB^\dag_\fp(\rV^\circ_n,\obj)_\sp$ to be the fiber product indicated in the following Cartesian diagram
\[
\xymatrix{
\rB^\dag_\fp(\rV^\circ_n,\obj)_\sp \ar[rr]^-{\pi^\dag_\sp}\ar[d]_-{\rb^\dag_\downarrow} && \rS^\dag_\fp(\rV^\circ_n,\obj)_\sp \ar[d]^-{\rs^\dag_\downarrow} \\
\rB^\dag_\fp(\rV^\circ_n,\obj_n) \ar[rr]^-{\pi^\dag_n} && \rS^\dag_\fp(\rV^\circ_n,\obj_n)
}
\]
in $\Fun(\fK(\rV^\circ_n)_\sp^p\times\fT,\Sch_{/\dF_p^\Phi})_{/\rT_\fp}$.
\end{definition}

By the universal property of Cartesian diagrams, we obtain a unique morphism
\[
\rb^{\dag\bullet}_\sp\colon\rB^\dag_\fp(\rV^\circ_n,\obj)_\sp\to\rB^\bullet_\fp(\rV^\circ_n,\obj)_\sp
\]
rendering the front lower-left cube of \eqref{eq:ns_functoriality} commute. Finally, an easy diagram chasing indicates that we have a unique morphism
\[
\rb^\dag_\uparrow\colon\rB^\dag_\fp(\rV^\circ_n,\obj)_\sp\to\rB^\dag_\fp(\rV^\circ_{n+1},\obj_{n+1})
\]
rendering the entire diagram \eqref{eq:ns_functoriality} commute. Thus, we obtain our desired diagram \eqref{eq:ns_functoriality}.

\begin{remark}\label{pr:ns_functoriality3}
By Proposition \ref{pr:ns_functoriality2} and Theorem \ref{th:ns_basic_correspondence3}(1), one can show that when $n$ is odd, the square
\[
\xymatrix{
\rB^\dag_\fp(\rV^\circ_{n+1},\obj_{n+1}) \ar[rr]^-{\iota^\dag_{n+1}} &&  \rM^\dag_\fp(\rV^\circ_{n+1},\obj_{n+1}) \\
\rB^\dag_\fp(\rV^\circ_n,\obj)_\sp \ar[rr]^-{\iota^\dag_n\circ\rb^\dag_\downarrow} \ar[u]^-{\rb^\dag_\uparrow}
&& \rM^\dag_\fp(\rV^\circ_n,\obj_n) \ar[u]_-{\rm^\dag_\uparrow}
}
\]
extracted from the diagram \eqref{eq:ns_functoriality} is Cartesian.
\end{remark}

\begin{remark}
By Lemma \ref{le:ns_functoriality_source}(1), Definition \ref{de:ns_base_special2}, Definition \ref{de:ns_definite_special3}, and Definition \ref{de:ns_base_special3}, the four downward arrows in the diagram \eqref{eq:ns_functoriality} are isomorphisms when $n$ is even.
\end{remark}

At the fourth stage of functoriality, we compare the special morphisms for basic correspondences and for Deligne--Lusztig varieties. Take a point
$s^\dag\in\rS^\dag_\fp(\rV^\circ_n,\rK^{p\circ})_\sp(\kappa)$ for a perfect field $\kappa$ containing $\dF_p^\Phi$. Put
\begin{align*}
s^\dag_n&\coloneqq\rs^\dag_\downarrow(s^\dag),\quad s^\dag_{n+1}\coloneqq\rs^\dag_\uparrow(s^\dag);\\
s^\circ_n&\coloneqq\rs^{\dag\circ}_n(s^\dag_n),\quad s^\circ_{n+1}\coloneqq\rs^{\dag\circ}_{n+1}(s^\dag_{n+1});\\
s^\bullet&\coloneqq\rs^{\dag\bullet}_\sp(s^\dag),\quad s^\bullet_n\coloneqq\rs^{\dag\bullet}_n(s^\dag_n),\quad
s^\bullet_{n+1}\coloneqq\rs^{\dag\bullet}_{n+1}(s^\dag_{n+1}).
\end{align*}
Denote by $\rB^\dag_{s^\dag}$, $\rB^\dag_{s^\dag_n}$, $\rB^\dag_{s^\dag_{n+1}}$, $\rB^\circ_{s^\circ_n}$, $\rB^\circ_{s^\circ_{n+1}}$, $\rB^\bullet_{s^\bullet}$, $\rB^\bullet_{s^\bullet_n}$, and $\rB^\bullet_{s^\bullet_{n+1}}$ their preimages under $\pi^\dag_\sp$, $\pi^\dag_n$, $\pi^\dag_{n+1}$, $\pi^\circ_n$, $\pi^\circ_{n+1}$, $\pi^\bullet_\sp$, $\pi^\bullet_n$, and $\pi^\bullet_{n+1}$, respectively.

\begin{proposition}\label{pr:ns_functoriality_dl}
Let the notation be as above. The following diagram
\begin{align*}
\xymatrix{
\rB^\circ_{s^\circ_{n+1}} \ar[rrrr]^-{\zeta^\circ_{s^\circ_{n+1}}} & & && \dP(\sV_{s^\circ_{n+1}}) & & \\
& \rB^\dag_{s^\dag_{n+1}} \ar[ul]_-{\rb^{\dag\circ}_{n+1}}\ar[rd]^-{\rb^{\dag\bullet}_{n+1}}\ar[rrrr]^-{\zeta^\dag_{s^\dag_{n+1}}} & & && \dP(\sV_{s^\dag_{n+1}}) \ar[ul]\ar[rd] & \\
& & \rB^\bullet_{s^\bullet_{n+1}} \ar[rrrr]^-{\zeta^\bullet_{s^\bullet_{n+1}}} & & && \DL^\bullet_{s^\bullet_{n+1}} \\
\rB^\circ_{s^\circ_n} \ar[uuu]^-{\rb^\circ_\uparrow} \ar[rrrr]|!{[rd];[ruu]}\hole|!{[rrdd];[rru]}\hole^-{\zeta^\circ_{s^\circ_n}} & & && \dP(\sV_{s^\circ_n}) \ar[uuu]|!{[ull];[urr]}\hole|!{[uulll];[uur]}\hole & & \\
& \rB^\dag_{s^\dag} \ar[uuu]^-{\rb^\dag_\uparrow} \ar[ul]^-{\rb^{\dag\circ}_n\circ\rb^\dag_\downarrow}\ar[rd]_-{\rb^{\dag\bullet}_\sp} \ar[rrrr]|!{[rd];[ruu]}\hole^-{\zeta^\dag_{s^\dag_n}\circ\rb^\dag_\downarrow} & & && \dP(\sV_{s^\dag_n})  \ar[ul]\ar[rd] \ar[uuu]|!{[uulll];[uur]}\hole & \\
& & \rB^\bullet_{s^\bullet} \ar[uuu]^-{\rb^\bullet_\uparrow} \ar[rrrr]^-{\zeta^\bullet_{s^\bullet_n}\circ\rb^\bullet_\downarrow} & & &&  \DL^\bullet_{s^\bullet_n} \ar[uuu]_-{\delta_{s^\bullet\uparrow}}
}
\end{align*}
in $\Sch_\kappa$ commutes, where
\begin{itemize}[label={\ding{109}}]
  \item $\zeta^\circ_{s^\circ_n}$ and $\zeta^\circ_{s^\circ_{n+1}}$ are the isomorphisms in Theorem \ref{th:ns_basic_correspondence1};

  \item $\zeta^\bullet_{s^\bullet_n}$ and $\zeta^\bullet_{s^\bullet_{n+1}}$ are the isomorphisms in Theorem \ref{th:ns_basic_correspondence2}(4);

  \item $\zeta^\dag_{s^\dag_n}$ and $\zeta^\dag_{s^\dag_{n+1}}$ are the isomorphisms in Theorem \ref{th:ns_basic_correspondence3}(2);

  \item $\dP(\sV_{s^\dag_n})\to\dP(\sV_{s^\circ_n})$ and $\dP(\sV_{s^\dag_{n+1}})\to\dP(\sV_{s^\circ_{n+1}})$ are closed embeddings in Remark \ref{re:ns_basic_correspondence3}(1);

  \item $\dP(\sV_{s^\dag_n})\to\DL^\bullet_{s^\bullet_n}=\DL^\bullet(\sV_{s^\bullet_n},\{\;,\;\}_{s^\bullet_n})$ and $\dP(\sV_{s^\dag_{n+1}})\to\DL^\bullet_{s^\bullet_{n+1}}=\DL^\bullet(\sV_{s^\bullet_{n+1}},\{\;,\;\}_{s^\bullet_{n+1}})$ are closed embeddings in Remark \ref{re:ns_basic_correspondence3}(2);

  \item $\dP(\sV_{s^\circ_n})\to\dP(\sV_{s^\circ_{n+1}})$ is the morphism induced by the obvious $\kappa$-linear (surjective) map $\sV_{s^\circ_{n+1}}\to\sV_{s^\circ_n}$;

  \item $\delta_{s^\bullet\uparrow}$ is the morphism in Construction \ref{cs:dl_bullet_special_morphism} with respect to the map $\delta_{s^\bullet}\colon\sV_{s^\bullet_n,\sharp}\to\sV_{s^\bullet_{n+1}}$ induced by $\delta^\bullet\colon A^\bullet\times A_0\to A^\bullet_\natural$; and

  \item $\dP(\sV_{s^\dag_n})\to\dP(\sV_{s^\dag_{n+1}})$ is the restriction of $\delta_{s^\bullet\uparrow}$, in view of Remark \ref{re:ns_basic_correspondence3}(2).
\end{itemize}
In particular, $\rb^\bullet_\uparrow\colon\rB^\bullet_{s^\bullet}\to\rB^\bullet_{s^\bullet_{n+1}}$ is an isomorphism when $n$ is even.
\end{proposition}

\begin{proof}
The proof is very similar to Proposition \ref{pr:qs_functoriality_dl}, which we leave to the readers. The last assertion follows as $\rb^\bullet_\downarrow\colon\rB^\bullet_{s^\bullet}\to\rB^\bullet_{s^\bullet_n}$ is always an isomorphism, and $\delta_{s^\bullet\uparrow}$ is an isomorphism when $n$ is even.
\end{proof}

At the final stage of functoriality, we relate the special morphisms for sources of basic correspondences to Shimura sets under the uniformization maps $\upsilon^\circ$ \eqref{eq:ns_uniformization1}, $\upsilon^\bullet$ \eqref{eq:ns_uniformization2}, and $\upsilon^\dag$ \eqref{eq:ns_uniformization3}. Recall that we have data $(\rV^\circ_n,\{\Lambda^\circ_{n,\fq}\}_{\fq\mid p})$ and $(\rV^\circ_{n+1},\{\Lambda^\circ_{n+1,\fq}\}_{\fq\mid p})$.

\begin{notation}\label{no:ns_uniformization}
We choose a lattice chain $\Lambda^\circ_{n,\fp}\subseteq\Lambda^\bullet_{n,\fp}\subseteq p^{-1}\Lambda^\circ_{n,\fp}$ of $\rV^\circ_n\otimes_FF_\fp$ and a lattice chain $\Lambda^\circ_{n+1,\fp}\subseteq\Lambda^\bullet_{n+1,\fp}\subseteq p^{-1}\Lambda^\circ_{n+1,\fp}$ of $\rV^\circ_{n+1}\otimes_FF_\fp$ satisfying the requirements in Construction \ref{cs:ns_uniformization2} for $N=n,n+1$, for which we assume that $(\Lambda^\bullet_{n,\fp})_\sharp\subseteq\Lambda^\bullet_{n+1,\fp}\subseteq p^{-1}(\Lambda^\bullet_{n,\fp})_\sharp^\vee$ holds. We now introduce various open compact subgroups at $p$.
\begin{itemize}[label={\ding{109}}]
  \item For $N\in\{n,n+1\}$, we have $\rK^\circ_{N,p}$ from Construction \ref{cs:ns_uniformization1}, $\rK^\bullet_{N,p}$ from Construction \ref{cs:ns_uniformization2}, and $\rK^\dag_{N,p}=\rK^\circ_{N,p}\cap\rK^\bullet_{N,p}$ from Construction \ref{cs:ns_uniformization3}.

  \item Put $\rK^\bullet_{\sp,\fp}\coloneqq\rK^\bullet_{n,\fp}\cap\rK^\bullet_{n+1,\fp}$ (as a subgroup of $\rK^\bullet_{n,\fp}$) and $\rK^\bullet_{\sp,p}\coloneqq\rK^\bullet_{\sp,\fp}\times\prod_{\fq\mid p,\fq\neq\fp}\rK^\circ_{n,\fq}$.

  \item Put $\rK^\dag_{\sp,p}\coloneqq\rK^\bullet_{\sp,p}\cap\rK^\circ_{n,p}$.
\end{itemize}
For later use, we also introduce natural maps
\begin{align*}
\begin{dcases}
\sh^\circ_\uparrow\colon\Sh(\rV^\circ_n,\obj_n\rK^\circ_{n,p})\to\Sh(\rV^\circ_{n+1},\obj_{n+1}\rK^\circ_{n+1,p}),\\
\sh^\bullet_\uparrow\colon\Sh(\rV^\circ_n,\obj_n\rK^\bullet_{\sp,p})\to\Sh(\rV^\circ_{n+1},\obj_{n+1}\rK^\bullet_{n+1,p}),\\
\sh^\bullet_\downarrow\colon\Sh(\rV^\circ_n,\obj_n\rK^\bullet_{n,p})\to\Sh(\rV^\circ_n,\obj_n\rK^\bullet_{\sp,p}),\\
\sh^\dag_\uparrow\colon\Sh(\rV^\circ_n,\obj_n\rK^\dag_{\sp,p})\to\Sh(\rV^\circ_{n+1},\obj_{n+1}\rK^\dag_{n+1,p}),\\
\sh^\dag_\downarrow\colon\Sh(\rV^\circ_n,\obj_n\rK^\dag_{n,p})\to\Sh(\rV^\circ_n,\obj_n\rK^\dag_{\sp,p}),\\
\sh^{\dag\circ}_n\colon\Sh(\rV^\circ_n,\obj_n\rK^\dag_{n,p})\to\Sh(\rV^\circ_n,\obj_n\rK^\circ_{n,p}), \\
\sh^{\dag\bullet}_n\colon\Sh(\rV^\circ_n,\obj_n\rK^\dag_{n,p})\to\Sh(\rV^\circ_n,\obj_n\rK^\bullet_{n,p}), \\
\sh^{\dag\circ}_{n+1}\colon\Sh(\rV^\circ_{n+1},\obj_{n+1}\rK^\dag_{n+1,p})\to\Sh(\rV^\circ_{n+1},\obj_{n+1}\rK^\circ_{n+1,p}), \\
\sh^{\dag\bullet}_{n+1}\colon\Sh(\rV^\circ_{n+1},\obj_{n+1}\rK^\dag_{n+1,p})\to\Sh(\rV^\circ_{n+1},\obj_{n+1}\rK^\bullet_{n+1,p}), \\
\sh^{\dag\bullet}_\sp\colon\Sh(\rV^\circ_n,\obj_n\rK^\dag_{\sp,p})\to\Sh(\rV^\circ_n,\obj_n\rK^\bullet_{\sp,p}),
\end{dcases}
\end{align*}
in $\Fun(\fK(\rV^\circ)^p_\sp,\Set)$. Note that $\sh^\circ_\uparrow$ has already appeared in Remark \ref{re:ns_functoriality_set}.
\end{notation}

Similar to Construction \ref{cs:qs_uniformization}, we may construct two uniformization maps
\begin{align}\label{eq:ns_uniformization4}
\upsilon^\bullet_\sp&\colon\rS^\bullet_\fp(\rV^\circ_n,\obj)_\sp(\ol\dF_p)
\to\Sh(\rV^\circ_n,\obj_n\rK^\bullet_{\sp,p})\times\rT_\fp(\ol\dF_p)
\end{align}
\begin{align}\label{eq:ns_uniformization5}
\upsilon^\dag_\sp&\colon\rS^\dag_\fp(\rV^\circ_n,\obj)_\sp(\ol\dF_p)
\to\Sh(\rV^\circ_n,\obj_n\rK^\dag_{\sp,p})\times\rT_\fp(\ol\dF_p)
\end{align}
in $\Fun(\fK(\rV^\circ_n)_\sp^p\times\fT,\Set)_{/\rT_\fp(\ol\dF_p)}$, which are isomorphisms. We leave the details to the readers.

\begin{proposition}\label{pr:ns_uniformization_basic}
The following diagram
\begin{align*}
\resizebox{0.45\hsize}{!}{
\rotatebox{-90}{
\xymatrix{
\rS^\circ_\fp(\rV^\circ_{n+1},\obj_{n+1})(\ol\dF_p) \ar[rrrr]^-{\upsilon^\circ_{n+1}}_-{\eqref{eq:ns_uniformization1}} &&&& \Sh(\rV^\circ_{n+1},\obj_{n+1}\rK^\circ_{n+1,p})\times\rT_\fp(\ol\dF_p)\\
& \rS^\dag_\fp(\rV^\circ_{n+1},\obj_{n+1})(\ol\dF_p) \ar[ul]_-{\rs^{\dag\circ}_{n+1}(\ol\dF_p)}\ar[rd]^-{\rs^{\dag\bullet}_{n+1}(\ol\dF_p)}
\ar[rrrr]^-{\upsilon^\dag_{n+1}}_-{\eqref{eq:ns_uniformization3}} &&&& \Sh(\rV^\circ_{n+1},\obj_{n+1}\rK^\dag_{n+1,p})\times\rT_\fp(\ol\dF_p) \ar[ul]_-{\sh^{\dag\circ}_{n+1}\times\id} \ar[rd]^-{\sh^{\dag\bullet}_{n+1}\times\id} \\
&& \rS^\bullet_\fp(\rV^\circ_{n+1},\obj_{n+1})(\ol\dF_p) \ar[rrrr]^-{\upsilon^\bullet_{n+1}}_-{\eqref{eq:ns_uniformization2}} &&&& \Sh(\rV^\circ_{n+1},\obj_{n+1}\rK^\bullet_{n+1,p})\times\rT_\fp(\ol\dF_p) \\
\\
& \rS^\dag_\fp(\rV^\circ_n,\obj)_\sp(\ol\dF_p)  \ar[rd]^-{\rs^{\dag\bullet}_\sp(\ol\dF_p)}\ar[uuu]_-{\rs^\dag_\uparrow(\ol\dF_p)}\ar[ddd]^-{\rs^\dag_\downarrow(\ol\dF_p)}
\ar[rrrr]|!{[rd];[ruu]}\hole^-{\upsilon^\dag_\sp}_-{\eqref{eq:ns_uniformization5}} &&&& \Sh(\rV^\circ_n,\obj_n\rK^\dag_{\sp,p})\times\rT_\fp(\ol\dF_p) \ar[dr]^-{\sh^{\dag\bullet}_\sp\times\id} \ar[ddd]|!{[dl];[dr]}\hole^-{\sh^\dag_\downarrow\times\id} \ar[uuu]|!{[uul];[uur]}\hole_-{\sh^\dag_\uparrow\times\id}  \\
&& \rS^\bullet_\fp(\rV^\circ_n,\obj)_\sp(\ol\dF_p) \ar[uuu]_-{\rs^\bullet_\uparrow(\ol\dF_p)}\ar[ddd]^-{\rs^\bullet_\downarrow(\ol\dF_p)}
\ar[rrrr]^-{\upsilon^\bullet_\sp}_-{\eqref{eq:ns_uniformization4}} &&&& \Sh(\rV^\circ_n,\obj_n\rK^\bullet_{\sp,p})\times\rT_\fp(\ol\dF_p)
\ar[uuu]_-{\sh^\bullet_\uparrow\times\id}\ar[ddd]^-{\sh^\bullet_\downarrow\times\id} \\
\rS^\circ_\fp(\rV^\circ_n,\obj_n)(\ol\dF_p) \ar[uuuuuu]_-{\rs^\circ_\uparrow(\ol\dF_p)}
\ar[rrrr]|!{[rd];[ruu]}\hole|!{[rrdd];[rru]}\hole^-{\upsilon^\circ_n}_-{\eqref{eq:ns_uniformization1}} &&&& \Sh(\rV^\circ_n,\obj_n\rK^\circ_{n,p})\times\rT_\fp(\ol\dF_p)
\ar[uuuuuu]|!{[ul];[ur]}\hole|!{[uul];[uur]}\hole|!{[uuuul];[uuuur]}\hole|!{[uuuuul];[uuuuur]}\hole_-{\sh^\circ_\uparrow\times\id}\\
& \rS^\dag_\fp(\rV^\circ_n,\obj_n)(\ol\dF_p) \ar[ul]_-{\rs^{\dag\circ}_n(\ol\dF_p)}\ar[rd]^-{\rs^{\dag\bullet}_n(\ol\dF_p)}
\ar[rrrr]|!{[rd];[ruu]}\hole^-{\upsilon^\dag_n}_-{\eqref{eq:ns_uniformization3}} &&&& \Sh(\rV^\circ_n,\obj_n\rK^\dag_{n,p})\times\rT_\fp(\ol\dF_p) \ar[ul]_-{\sh^{\dag\circ}_n\times\id} \ar[dr]^-{\sh^{\dag\bullet}_n\times\id} \\
&& \rS^\bullet_\fp(\rV^\circ_n,\obj_n)(\ol\dF_p)
\ar[rrrr]^-{\upsilon^\bullet_n}_-{\eqref{eq:ns_uniformization2}} &&&& \Sh(\rV^\circ_n,\obj_n\rK^\bullet_{n,p})\times\rT_\fp(\ol\dF_p)
}
}
}
\end{align*}
in $\Fun(\fK(\rV^\circ_n)_\sp^p\times\fT,\Set)_{/\rT_\fp(\ol\dF_p)}$ commutes (in which all uniformization maps are isomorphisms). Moreover, the induced actions of $\Gal(\ol\dF_p/\dF_p^\Phi)$ on all terms on the right-hand side factor through the projection to the factor $\rT_\fp(\ol\dF_p)$.
\end{proposition}

\begin{proof}
This follows from Constructions \ref{cs:ns_uniformization1}, \ref{cs:ns_uniformization2}, and \ref{cs:ns_uniformization3}.
\end{proof}

\begin{remark}\label{re:ns_functoriality}
When $n=1$, we have the diagram \eqref{eq:ns_functoriality} in which all terms not in the top or back layers are empty. Propositions \ref{pr:ns_functoriality_dl} and \ref{pr:ns_uniformization_basic} can be modified in the obvious way.
\end{remark}

\subsection{First geometric reciprocity law}
\label{ss:ns_reciprocity}

In this subsection, we state and prove a theorem we call \emph{first geometric reciprocity law}, which can be regarded a geometric template for the first explicit reciprocity law studied in \S\ref{ss:first_reciprocity} once we plug in the automorphic input.

We maintain the setup in \S\ref{ss:ns_functoriality}. However, we allow $\obj=(\obj_n,\obj_{n+1})$ to be an object of $\fK(\rV^\circ_n)^p\times\fK(\rV^\circ_{n+1})^p$, rather than $\fK(\rV^\circ_n)^p_\sp$. Denote by $n_0$ and $n_1$ the unique even and odd numbers in $\{n,n+1\}$, respectively. Write $n_0=2r_0$ and $n_1=2r_1+1$ for unique integers $r_0,r_1\geq 1$. In particular, we have $n=r_0+r_1$. Let $L$ be a $p$-coprime coefficient ring.

To ease notation, we put $\rX^?_{n_\alpha}\coloneqq\rX^?_\fp(\rV^\circ_{n_\alpha},\obj_{n_\alpha})$ for meaningful triples $(\rX,?,\alpha)\in\{\bM,\rM,\rB,\rS\}\times\{\;,\eta,\circ,\bullet,\dag\}\times\{0,1\}$.

\begin{notation}\label{no:ns_product}
We introduce following objects.
\begin{enumerate}
  \item Put $\bP\coloneqq\bM_{n_0}\times_{\bT_\fp}\bM_{n_1}$.

  \item For $(?_0,?_1)\in\{\circ,\bullet,\dag\}^2$, put $\rP^{?_0,?_1}\coloneqq\rM^{?_0}_{n_0}\times_{\rT_\fp}\rM^{?_1}_{n_1}$, which is a closed subscheme of $\rP$.\footnote{Recall from Notation \ref{no:p_notation}(5) that $\rP$ is $\bP\otimes_{\dZ_p^\Phi}\dF_p^\Phi$.}

  \item Let $\sigma\colon\bQ\to\bP$ be the blow-up along the subscheme $\rP^{\circ,\circ}$, which is a morphism in $\Fun(\fK(\rV^\circ_n)^p\times\fK(\rV^\circ_{n+1})^p\times\fT,\Sch_{/\dZ_p^\Phi})_{/\bT_\fp}$.

  \item For $(?_0,?_1)\in\{\circ,\bullet,\dag\}^2$, let $\rQ^{?_0,?_1}$ be the strict transform of $\rP^{?_0,?_1}$ under $\sigma$, which is a closed subscheme of $\rQ$.

  \item Let $\gamma^{?_0,?_1}_{?'_0,?'_1}\colon\rP^{?_0,?_1}\to\rP^{?'_0,?'_1}$ be the closed embedding if $\rP^{?_0,?_1}$ is contained in $\rP^{?'_0,?'_1}$, and $\delta^{?_0,?_1}_{?'_0,?'_1}\colon\rQ^{?_0,?_1}\to\rQ^{?'_0,?'_1}$ the closed embedding if $\rQ^{?_0,?_1}$ is contained in $\rQ^{?'_0,?'_1}$.
\end{enumerate}
Suppose that $\obj$ is taken in the subcategory $\fK(\rV^\circ_n)^p_\sp$.
\begin{enumerate}\setcounter{enumi}{5}
  \item Let $\bP_\graph$ be the graph of $\bm_\uparrow\colon\bM_n\to\bM_{n+1}$ \eqref{eq:ns_functoriality_moduli_scheme} over $\bT_\fp$ in $\Fun(\fK(\rV^\circ_n)_\sp^p\times\fT,\Sch_{/\dZ_p^\Phi})_{/\bT_\fp}$, as a closed subscheme of $\bP$.

  \item For $?=\bullet,\circ$, let $\rP^?_\graph$ be the graph of $\rm^?_\uparrow\colon\rM^?_n\to\rM^?_{n+1}$ \eqref{eq:ns_functoriality_moduli_scheme1} over $\rT_\fp$ in $\Fun(\fK(\rV^\circ_n)_\sp^p\times\fT,\Sch_{/\dF_p^\Phi})_{/\rT_\fp}$, as a closed subscheme of $\rP^{?,?}$.

  \item Let $\bQ_\graph$ be the strict transform of $\bP_\graph$ under $\sigma$, which is a closed subscheme of $\bQ$.
\end{enumerate}
\end{notation}

\begin{lem}\label{le:ns_product_pbc}
The two specialization maps
\begin{align*}
\rH^i_{\fT,c}(\bQ\otimes_{\dZ_{p^2}}\ol\dQ_p,L)&\to\rH^i_{\fT,c}(\ol\rQ,\rR\Psi L),\\
\rH^i_\fT(\bQ\otimes_{\dZ_{p^2}}\ol\dQ_p,L)&\to\rH^i_\fT(\ol\rQ,\rR\Psi L),
\end{align*}
are both isomorphisms.
\end{lem}

\begin{proof}
When $\bQ$ is proper, this is simply the proper base change. When $\bQ$ is not proper, this again follows from \cite{LS18}*{Corollary~5.20}.
\end{proof}

\begin{lem}\label{le:ns_product_weight}
The scheme $\bQ$ (valued at any object of $\fK(\rV^\circ_n)_\sp^p$) is strictly semistable over $\dZ_p^\Phi$ of relative dimension $2n-1$. Moreover, we have
\begin{enumerate}
  \item The reduction graph of $\bQ$ is as follows
    
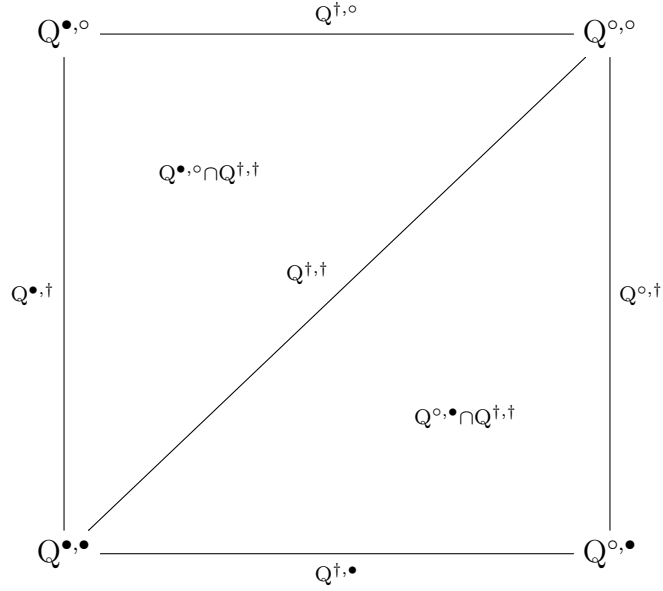
\[
    \xymatrix{
    \rQ^{\bullet,\circ} \ar@{-}[rrrrrr]^-{\rQ^{\dag,\circ}}\ar@{-}[dddddd]_-{\rQ^{\bullet,\dag}}
    \ar@{}[dddrrr]|-{\rQ^{\bullet,\circ}\cap\rQ^{\dag,\dag}}
    &&&&&& \rQ^{\circ,\circ} \ar@{-}[dddddd]^-{\rQ^{\circ,\dag}} \\
    \\
    \\
    &&&&&&\\
    \\
    \\
    \rQ^{\bullet,\bullet} \ar@{-}[rrrrrr]_-{\rQ^{\dag,\bullet}}\ar@{-}[rrrrrruuuuuu]^-{\rQ^{\dag,\dag}}&&&&&&
    \rQ^{\circ,\bullet} \ar@{}[llluuu]|-{\rQ^{\circ,\bullet}\cap\rQ^{\dag,\dag}}
    }
    \]
    so that
    \[
    \begin{dcases}
    \rQ^{(0)}=\rQ^{\circ,\circ}\coprod\rQ^{\circ,\bullet}\coprod\rQ^{\bullet,\bullet}\coprod\rQ^{\bullet,\circ},\\
    \rQ^{(1)}=\rQ^{\circ,\dag}\coprod\rQ^{\dag,\bullet}\coprod\rQ^{\bullet,\dag}\coprod\rQ^{\dag,\circ}\coprod\rQ^{\dag,\dag},\\
    \rQ^{(2)}=(\rQ^{\bullet,\circ}\cap\rQ^{\dag,\dag})\coprod(\rQ^{\circ,\bullet}\cap\rQ^{\dag,\dag}),\\
    \rQ^{(c)}=\emptyset,\text{ for $c\geq 3$.}
    \end{dcases}
    \]
    Here, $\rQ^{(c)}$ denotes the disjoint union of the strata of $\rQ$ of codimension $c$.

  \item For the morphism $\sigma$, we have that
    \begin{itemize}[label={\ding{109}}]
      \item the induced morphism $\sigma\colon\rQ^{?_0,?_1}\to\rP^{?_0,?_1}$ is an isomorphism if $?_0\neq ?_1$;

      \item the induced morphism $\sigma\colon\rQ^{?_0,?_1}\to\rP^{?_0,?_1}$ is the blow-up along $\rP^{\dag,\dag}$ if $(?_0,?_1)\in\{(\circ,\circ),(\bullet,\bullet)\}$;

      \item the induced morphism $\sigma\colon\rQ^{\dag,\dag}\to\rP^{\dag,\dag}$ is a trivial $\dP^1$-bundle;

      \item the induced morphisms $\sigma\colon\rQ^{\bullet,\circ}\cap\rQ^{\dag,\dag}\to\rP^{\dag,\dag}$ and $\sigma\colon\rQ^{\circ,\bullet}\cap\rQ^{\dag,\dag}\to\rP^{\dag,\dag}$ are both isomorphisms.

    \end{itemize}

  \item The natural map
     \[
     \sigma^*\colon\rH^i_\fT(\ol\rP^{?_0,?_1},O_\lambda)\to\rH^i_\fT(\ol\rQ^{?_0,?_1},O_\lambda)
     \]
     is injective, and moreover an isomorphism if $?_0\neq ?_1$.

  \item For $(?_0,?_1)\in\{(\circ,\circ),(\bullet,\bullet)\}$, the map
      \[
      (\delta^{\dag,\dag}_{?_0,?_1})_!\circ\sigma^*\colon\rH^{i-2}_\fT(\ol\rP^{\dag,\dag},O_\lambda(-1))\to\rH^i_\fT(\ol\rQ^{?_0,?_1},O_\lambda)
      \]
      is injective; and we have
      \[
      \rH^i_\fT(\ol\rQ^{?_0,?_1},O_\lambda)=\sigma^*\rH^i_\fT(\ol\rP^{?_0,?_1},O_\lambda)\bigoplus
      (\delta^{\dag,\dag}_{?_0,?_1})_!\sigma^*\rH^{i-2}_\fT(\ol\rP^{\dag,\dag},O_\lambda(-1)).
      \]

  \item If we denote by $\ff\in\rH^2_\fT(\ol\rQ^{\dag,\dag},O_\lambda(1))$ the cycle class of an arbitrary $\fT$-orbit of fibers of the trivial $\dP^1$-fibration $\ol\sigma\colon\ol\rQ^{\dag,\dag}\to\ol\rP^{\dag,\dag}$, then the map
      \[
      (\ff\cup)\circ\sigma^*\colon\rH^{i-2}_\fT(\ol\rP^{\dag,\dag},O_\lambda(-1))\to\rH^i_\fT(\ol\rQ^{\dag,\dag},O_\lambda)
      \]
      is injective; and we have
      \[
      \rH^i_\fT(\ol\rQ^{\dag,\dag},O_\lambda)=
      \sigma^*\rH^i_\fT(\ol\rP^{\dag,\dag},O_\lambda)\bigoplus\ff\cup\sigma^*\rH^{i-2}_\fT(\ol\rP^{\dag,\dag},O_\lambda(-1)).
      \]
\end{enumerate}
\end{lem}

\begin{proof}
Parts (1,2) follow from a standard computation of blow-up. Parts (3--5) follow from (2).
\end{proof}

Let $(\dE^{p,q}_s,\rd^{p,q}_s)$ be the weight spectral sequence abutting to the cohomology $\rH^{p+q}_\fT(\ol\rQ,\rR\Psi O_\lambda(n))$,\footnote{Strictly speaking, the differential maps $\rd^{p,q}_s$ depend on the choice of the ordering of (types of) irreducible components of $\rQ$, which we choose to be the clockwise order $\rQ^{\circ,\circ}<\rQ^{\circ,\bullet}<\rQ^{\bullet,\bullet}<\rQ^{\bullet,\circ}$.} whose first page is as follows:
\begin{align}\label{eq:ns_product_weight}
\resizebox{0.45\hsize}{!}{
\rotatebox{-90}{
\boxed{
\xymatrix{
q\geq 2n+1 & \cdots \ar[r] & \cdots \ar[r] & \cdots \ar[r] & \cdots \ar[r] & \cdots \\
q=2n & \rH^{2n-4}_\fT(\ol\rQ^{(2)},O_\lambda(n-2)) \ar[r]^-{\rd^{-2,2n}_1}
& \rH^{2n-2}_\fT(\ol\rQ^{(1)},O_\lambda(n-1)) \ar[r]^-{\rd^{-1,2n}_1}
& \txt{$\rH^{2n}_\fT(\ol\rQ^{(0)},O_\lambda(n))$\\$\bigoplus$\\$\rH^{2n-2}_\fT(\ol\rQ^{(2)},O_\lambda(n-1))$} \ar[r]^-{\rd^{0,2n}_1}
& \rH^{2n}_\fT(\ol\rQ^{(1)},O_\lambda(n)) \ar[r]^-{\rd^{1,2n}_1}
& \rH^{2n}_\fT(\ol\rQ^{(2)},O_\lambda(n)) \\
q=2n-1 & \rH^{2n-5}_\fT(\ol\rQ^{(2)},O_\lambda(n-2)) \ar[r]^-{\rd^{-2,2n-1}_1}
& \rH^{2n-3}_\fT(\ol\rQ^{(1)},O_\lambda(n-1)) \ar[r]^-{\rd^{-1,2n-1}_1}
& \txt{$\rH^{2n-1}_\fT(\ol\rQ^{(0)},O_\lambda(n))$\\$\bigoplus$\\$\rH^{2n-3}_\fT(\ol\rQ^{(2)},O_\lambda(n-1))$} \ar[r]^-{\rd^{0,2n-1}_1}
& \rH^{2n-1}_\fT(\ol\rQ^{(1)},O_\lambda(n)) \ar[r]^-{\rd^{1,2n-1}_1}
& \rH^{2n-1}_\fT(\ol\rQ^{(2)},O_\lambda(n)) \\
q=2n-2 & \rH^{2n-6}_\fT(\ol\rQ^{(2)},O_\lambda(n-2)) \ar[r]^-{\rd^{-2,2n-2}_1}
& \rH^{2n-4}_\fT(\ol\rQ^{(1)},O_\lambda(n-1)) \ar[r]^-{\rd^{-1,2n-2}_1}
& \txt{$\rH^{2n-2}_\fT(\ol\rQ^{(0)},O_\lambda(n))$\\$\bigoplus$\\$\rH^{2n-4}_\fT(\ol\rQ^{(2)},O_\lambda(n-1))$} \ar[r]^-{\rd^{0,2n-2}_1}
& \rH^{2n-2}_\fT(\ol\rQ^{(1)},O_\lambda(n)) \ar[r]^-{\rd^{1,2n-2}_1}
& \rH^{2n-2}_\fT(\ol\rQ^{(2)},O_\lambda(n)) \\
q\leq 2n-3 & \cdots \ar[r] & \cdots \ar[r] & \cdots \ar[r] & \cdots \ar[r] & \cdots \\
\dE^{p,q}_1 & p=-2 & p=-1 & p=0 & p=1 & p=2
}
}
}
}
\end{align}
with $\dE^{p,q}_1=0$ if $|p|>2$.

\begin{construction}\label{cs:ns_incidence_product}
For $\alpha=0,1$, let $\xi_\alpha\in\rH^2_\fT(\ol\rB^\circ_{n_\alpha},L(1))$ be the first Chern class of the tautological quotient line bundle on $\ol\rB^\circ_{n_\alpha}$. We construct four new pairs of maps in $\Fun(\fK(\rV^\circ_n)^p\times\fK(\rV^\circ_{n+1})^p,\Mod(L))$ as follows:
\begin{align*}
&\begin{dcases}
\inc^{\circ,\dag}_!\colon L[\Sh(\rV^\circ_{n_0},\obj_{n_0}\rK^\circ_{n_0,p})]\otimes_L L[\Sh(\rV^\circ_{n_1},\obj_{n_1}\rK^\circ_{n_1,p})] \\
\qquad\qquad\xrightarrow{\sim}
\rH^0_\fT(\ol\rS^\circ_{n_0},L)\otimes_L\rH^0_\fT(\ol\rS^\circ_{n_1},L)=\rH^0_\fT(\ol\rS^\circ_{n_0}\times_{\ol\rT_\fp}\ol\rS^\circ_{n_1},L) \\
\qquad\qquad\xrightarrow{(\pi^\circ_{n_0}\times\pi^\circ_{n_1})^*}
\rH^0_\fT(\ol\rB^\circ_{n_0}\times_{\ol\rT_\fp}\ol\rB^\circ_{n_1},L) \\
\qquad\qquad\xrightarrow{\cup\xi_0^{r_0-1}\cup\xi_1^{r_1-1}}
\rH^{2(n-2)}_\fT(\ol\rB^\circ_{n_0}\times_{\ol\rT_\fp}\ol\rB^\circ_{n_1},L(n-2)) \\
\qquad\qquad\xrightarrow{(\iota^\circ_{n_0}\times\iota^\circ_{n_1})_!}
\rH^{2(n-2)}_\fT(\ol\rM^\circ_{n_0}\times_{\ol\rT_\fp}\ol\rM^\circ_{n_1},L(n-2)) \\
\qquad\qquad\xrightarrow{(\id\times\rm^{\dag\circ}_{n_1})^*}
\rH^{2(n-2)}_\fT(\ol\rM^\circ_{n_0}\times_{\ol\rT_\fp}\ol\rM^\dag_{n_1},L(n-2)) \\
\qquad\qquad\xrightarrow{(\id\times\rm^{\dag\bullet}_{n_1})_!}
\rH^{2(n-1)}_{\fT,c}(\ol\rM^\circ_{n_0}\times_{\ol\rT_\fp}\ol\rM^\bullet_{n_1},L(n-1))=\rH^{2(n-1)}_{\fT,c}(\ol\rP^{\circ,\bullet},L(n-1)), \\
\inc_{\circ,\dag}^*\colon\rH^{2n}_\fT(\ol\rP^{\circ,\bullet},L(n))
=\rH^{2n}_\fT(\ol\rM^\circ_{n_0}\times_{\ol\rT_\fp}\ol\rM^\bullet_{n_1},L(n)) \\
\qquad\qquad\xrightarrow{(id\times\rm^{\dag\bullet}_{n_1})^*}
\rH^{2n}_\fT(\ol\rM^\circ_{n_0}\times_{\ol\rT_\fp}\ol\rM^\dag_{n_1},L(n)) \\
\qquad\qquad\xrightarrow{(\id\times\rm^{\dag\circ}_{n_1})_!}
\rH^{2n+2}_\fT(\ol\rM^\circ_{n_0}\times_{\ol\rT_\fp}\ol\rM^\circ_{n_1},L(n+1)) \\
\qquad\qquad\xrightarrow{(\iota^\circ_{n_0}\times\iota^\circ_{n_1})^*}
\rH^{2n+2}_\fT(\ol\rB^\circ_{n_0}\times_{\ol\rT_\fp}\ol\rB^\circ_{n_1},L(n+1)) \\
\qquad\qquad\xrightarrow{\cup\xi_0^{r_0-1}\cup\xi_1^{r_1-1}}
\rH^{4n-2}_\fT(\ol\rB^\circ_{n_0}\times_{\ol\rT_\fp}\ol\rB^\circ_{n_1},L(2n-1)) \\
\qquad\qquad\xrightarrow{(\pi^\circ_{n_0}\times\pi^\circ_{n_1})_!}
\rH^0_\fT(\ol\rS^\circ_{n_0}\times_{\ol\rT_\fp}\ol\rS^\circ_{n_1},L)=\rH^0_\fT(\ol\rS^\circ_{n_0},L)\otimes_L\rH^0_\fT(\ol\rS^\circ_{n_1},L) \\
\qquad\qquad\xrightarrow{\sim}L[\Sh(\rV^\circ_{n_0},\obj_{n_0}\rK^\circ_{n_0,p})]\otimes_L L[\Sh(\rV^\circ_{n_1},\obj_{n_1}\rK^\circ_{n_1,p})];
\end{dcases}
\end{align*}
\begin{align*}
\begin{dcases}
\inc^{\circ,\bullet}_!\colon L[\Sh(\rV^\circ_{n_0},\obj_{n_0}\rK^\circ_{n_0,p})]\otimes_L L[\Sh(\rV^\circ_{n_1},\obj_{n_1}\rK^\bullet_{n_1,p})] \\
\qquad\qquad\xrightarrow{\sim}
\rH^0_\fT(\ol\rS^\circ_{n_0},L)\otimes_L\rH^0_\fT(\ol\rS^\bullet_{n_1},L) =\rH^0_\fT(\ol\rS^\circ_{n_0}\times_{\ol\rT_\fp}\ol\rS^\bullet_{n_1},L) \\
\qquad\qquad\xrightarrow{(\pi^\circ_{n_0}\times\pi^\bullet_{n_1})^*}
\rH^0_\fT(\ol\rB^\circ_{n_0}\times_{\ol\rT_\fp}\ol\rB^\bullet_{n_1},L) \\
\qquad\qquad\xrightarrow{\cup\xi_0^{r_0-1}}
\rH^{2(r_0-1)}_\fT(\ol\rB^\circ_{n_0}\times_{\ol\rT_\fp}\ol\rB^\bullet_{n_1},L(r_0-1)) \\
\qquad\qquad\xrightarrow{(\iota^\circ_{n_0}\times\iota^\bullet_{n_1})_!}
\rH^{2(n-1)}_{\fT,c}(\ol\rM^\circ_{n_0}\times_{\ol\rT_\fp}\ol\rM^\bullet_{n_1},L(n-1))=\rH^{2(n-1)}_{\fT,c}(\ol\rP^{\circ,\bullet},L(n-1)), \\
\inc_{\circ,\bullet}^*\colon\rH^{2n}_\fT(\ol\rP^{\circ,\bullet},L(n))
=\rH^{2n}_\fT(\ol\rM^\circ_{n_0}\times_{\ol\rT_\fp}\ol\rM^\bullet_{n_1},L(n)) \\
\qquad\qquad\xrightarrow{(\iota^\circ_{n_0}\times\iota^\bullet_{n_1})^*}
\rH^{2n}_\fT(\ol\rB^\circ_{n_0}\times_{\ol\rT_\fp}\ol\rB^\bullet_{n_1},L(n)) \\
\qquad\qquad\xrightarrow{\cup\xi_0^{r_0-1}}\rH^{2(n_0-1+r_1)}_\fT
(\ol\rB^\circ_{n_0}\times_{\ol\rT_\fp}\ol\rB^\bullet_{n_1},L(n_0-1+r_1)) \\
\qquad\qquad\xrightarrow{(\pi^\circ_{n_0}\times\pi^\bullet_{n_1})_!}
\rH^0_\fT(\ol\rS^\circ_{n_0}\times_{\ol\rT_\fp}\ol\rS^\bullet_{n_1},L) =\rH^0_\fT(\ol\rS^\circ_{n_0},L)\otimes_L\rH^0_\fT(\ol\rS^\bullet_{n_1},L) \\
\qquad\qquad\xrightarrow{\sim}L[\Sh(\rV^\circ_{n_0},\obj_{n_0}\rK^\circ_{n_0,p})]\otimes_L L[\Sh(\rV^\circ_{n_1},\obj_{n_1}\rK^\bullet_{n_1,p})];
\end{dcases}
\end{align*}

\begin{align*}
&\begin{dcases}
\inc^{\bullet,\dag}_!\colon L[\Sh(\rV^\circ_{n_0},\obj_{n_0}\rK^\bullet_{n_0,p})]
\otimes_L L[\Sh(\rV^\circ_{n_1},\obj_{n_1}\rK^\circ_{n_1,p})] \\
\qquad\qquad\xrightarrow{\sim}\rH^0_\fT(\ol\rS^\bullet_{n_0},L)\otimes_L \rH^0_\fT(\ol\rS^\circ_{n_1},L) =\rH^0_\fT(\ol\rS^\bullet_{n_0}\times_{\ol\rT_\fp}\ol\rS^\circ_{n_1},L) \\
\qquad\qquad\xrightarrow{(\pi^\circ_{n_0}\times\pi^\circ_{n_1})^*}
\rH^0_\fT(\ol\rB^\bullet_{n_0}\times_{\ol\rT_\fp}\ol\rB^\circ_{n_1},L) \\
\qquad\qquad\xrightarrow{\cup\xi_1^{r_1-1}}
\rH^{2r_1-2}_\fT(\ol\rB^\bullet_{n_0}\times_{\ol\rT_\fp}\ol\rB^\circ_{n_1},L(r_1-1)) \\
\qquad\qquad\xrightarrow{(\id\times\iota^\circ_{n_1})_!}
\rH^{2r_1-2}_\fT(\ol\rB^\bullet_{n_0}\times_{\ol\rT_\fp}\ol\rM^\circ_{n_1},L(r_1-1)) \\
\qquad\qquad\xrightarrow{(\id\times\rm^{\dag\circ}_{n_1})^*}
\rH^{2r_1-2}_\fT(\ol\rB^\bullet_{n_0}\times_{\ol\rT_\fp}\ol\rM^\dag_{n_1},L(r_1-1)) \\
\qquad\qquad\xrightarrow{(\iota^\bullet_{n_0}\times\rm^{\dag\bullet}_{n_1})_!}
\rH^{2(n-1)}_{\fT,c}(\ol\rM^\bullet_{n_0}\times_{\ol\rT_\fp}\ol\rM^\bullet_{n_1},L(n-1))
=\rH^{2(n-1)}_{\fT,c}(\ol\rP^{\bullet,\bullet},L(n-1)), \\
\inc_{\bullet,\dag}^*\colon\rH^{2n}_\fT(\ol\rP^{\bullet,\bullet},L(n))
=\rH^{2n}_\fT(\ol\rM^\bullet_{n_0}\times_{\ol\rT_\fp}\ol\rM^\bullet_{n_1},L(n)) \\
\qquad\qquad\xrightarrow{(\iota^\bullet_{n_0}\times\rm^{\dag\bullet}_{n_1})^*}
\rH^{2n}_\fT(\ol\rB^\bullet_{n_0}\times_{\ol\rT_\fp}\ol\rM^\dag_{n_1},L(n)) \\
\qquad\qquad\xrightarrow{(\id\times\rm^{\dag\circ}_{n_1})_!}
\rH^{2n+2}_\fT(\ol\rB^\bullet_{n_0}\times_{\ol\rT_\fp}\ol\rM^\circ_{n_1},L(n+1)) \\
\qquad\qquad\xrightarrow{(\id\times\iota^\circ_{n_1})^*}
\rH^{2n+2}_\fT(\ol\rB^\bullet_{n_0}\times_{\ol\rT_\fp}\ol\rB^\circ_{n_1},L(n+1)) \\
\qquad\qquad\xrightarrow{\cup\xi_1^{r_1-1}}\rH^{2(r_0+n_1-1)}_\fT
(\ol\rB^\bullet_{n_0}\times_{\ol\rT_\fp}\ol\rB^\circ_{n_1},L(r_0+n_1-1)) \\
\qquad\qquad\xrightarrow{(\pi^\bullet_{n_0}\times\pi^\circ_{n_1})_!}
\rH^0_\fT(\ol\rS^\bullet_{n_0}\times_{\ol\rT_\fp}\ol\rS^\circ_{n_1},L) =\rH^0_\fT(\ol\rS^\bullet_{n_0},L)\otimes_L\rH^0_\fT(\ol\rS^\circ_{n_1},L) \\
\qquad\qquad\xrightarrow{\sim}L[\Sh(\rV^\circ_{n_0},\obj_{n_0}\rK^\bullet_{n_0,p})]\otimes_L L[\Sh(\rV^\circ_{n_1},\obj_{n_1}\rK^\circ_{n_1,p})];
\end{dcases}
\end{align*}
\begin{align*}
\begin{dcases}
\inc^{\bullet,\bullet}_!\colon L[\Sh(\rV^\circ_{n_0},\obj_{n_0}\rK^\bullet_{n_0,p})]
\otimes_L L[\Sh(\rV^\circ_{n_1},\obj_{n_1}\rK^\bullet_{n_1,p})] \\
\qquad\qquad\xrightarrow{\sim}\rH^0_\fT(\ol\rS^\bullet_{n_0},L)\otimes_L \rH^0_\fT(\ol\rS^\bullet_{n_1},L) =\rH^0_\fT(\ol\rS^\bullet_{n_0}\times_{\ol\rT_\fp}\ol\rS^\bullet_{n_1},L) \\
\qquad\qquad\xrightarrow{(\pi^\bullet_{n_0}\times\pi^\bullet_{n_1})^*}
\rH^0_\fT(\ol\rB^\bullet_{n_0}\times_{\ol\rT_\fp}\ol\rB^\bullet_{n_1},L) \\
\qquad\qquad\xrightarrow{(\iota^\bullet_{n_0}\times\iota^\bullet_{n_1})_!}
\rH^{2(n-1)}_{\fT,c}(\ol\rM^\bullet_{n_0}\times_{\ol\rT_\fp}\ol\rM^\bullet_{n_1},L(n-1))
=\rH^{2(n-1)}_{\fT,c}(\ol\rP^{\bullet,\bullet},L(n-1)), \\
\inc_{\bullet,\bullet}^*\colon\rH^{2n}_\fT(\ol\rP^{\bullet,\bullet},L(n))
=\rH^{2n}_\fT(\ol\rM^\bullet_{n_0}\times_{\ol\rT_\fp}\ol\rM^\bullet_{n_1},L(n)) \\
\qquad\qquad\xrightarrow{(\iota^\bullet_{n_0}\times\iota^\bullet_{n_1})^*}
\rH^{2n}_\fT(\ol\rB^\bullet_{n_0}\times_{\ol\rT_\fp}\ol\rB^\bullet_{n_1},L(n)) \\
\qquad\qquad\xrightarrow{(\pi^\bullet_{n_0}\times\pi^\bullet_{n_1})_!}
\rH^0_\fT(\ol\rS^\bullet_{n_0}\times_{\ol\rT_\fp}\ol\rS^\bullet_{n_1},L) =\rH^0_\fT(\ol\rS^\bullet_{n_0},L)\otimes_L\rH^0_\fT(\ol\rS^\bullet_{n_1},L) \\
\qquad\qquad\xrightarrow{\sim}L[\Sh(\rV^\circ_{n_0},\obj_{n_0}\rK^\bullet_{n_0,p})]\otimes_L L[\Sh(\rV^\circ_{n_1},\obj_{n_1}\rK^\bullet_{n_1,p})].
\end{dcases}
\end{align*}
In fact, the two maps in each pair are Poincar\'{e} dual to each other.
\end{construction}

\begin{theorem}[First geometric reciprocity law]\label{th:ns_reciprocity}
Take an object $\rK^{p\circ}\in\fK(\rV^\circ_n)_\sp^p$. For the class $\cl(\rP^\bullet_\graph)\in\rH^{2n}_\fT(\ol\rP^{\bullet,\bullet},L(n))$, we have
\begin{enumerate}
  \item For $f\in L[\Sh(\rV^\circ_{n_0},\rK^{p\circ}_{n_0}\rK^\bullet_{n_0,p})]\otimes_L L[\Sh(\rV^\circ_{n_1},\rK^{p\circ}_{n_1}\rK^\circ_{n_1,p})]$, the identity
      \[
      \int_{\ol\rP^{\bullet,\bullet}}^\fT\cl(\rP^\bullet_\graph)
      \cup\inc^{\bullet,\dag}_!(f)=\sum_{s\in\Sh(\rV^\circ_n,\rK^{p\circ}_n\rK^\bullet_{\sp,p})}
      (\tT^{\bullet\circ}_{n_1,\fp}f)(\sh^\bullet_\downarrow(s),\sh^\bullet_\uparrow(s))
      \]
      holds.

  \item For $f\in L[\Sh(\rV^\circ_{n_0},\rK^{p\circ}_{n_0}\rK^\bullet_{n_0,p})]\otimes_L L[\Sh(\rV^\circ_{n_1},\rK^{p\circ}_{n_1}\rK^\bullet_{n_1,p})]$, the identity
      \[
      \int_{\ol\rP^{\bullet,\bullet}}^\fT\cl(\rP^\bullet_\graph)
      \cup\inc^{\bullet,\bullet}_!(f)=
      \sum_{s\in\Sh(\rV^\circ_n,\rK^{p\circ}_n\rK^\bullet_{\sp,p})}(\tT^\bullet_{n_1,\fp}f)(\sh^\bullet_\downarrow(s),\sh^\bullet_\uparrow(s))
      \]
      holds.

  \item For $f\in L[\Sh(\rV^\circ_{n_0},\rK^{p\circ}_{n_0}\rK^\circ_{n_0,p})]\otimes_L L[\Sh(\rV^\circ_{n_1},\rK^{p\circ}_{n_1}\rK^\circ_{n_1,p})]$, the identity
      \begin{align*}
      &\quad\int_{\ol\rP^{\bullet,\bullet}}^\fT\cl(\rP^\bullet_\graph)\cup
      \(\inc^{\bullet,\dag}_!(\tT^{\bullet\circ}_{n_0,\fp}\otimes\tI^\circ_{n_1,\fp}f)+
      (p+1)^2\inc^{\bullet,\bullet}_!(\tT^{\bullet\circ}_{n_0,\fp}\otimes\tT^{\bullet\circ}_{n_1,\fp}f)\) \\
      &=\sum_{s\in\Sh(\rV^\circ_n,\rK^{p\circ}_n\rK^\circ_{n,p})}(\tI^\circ_{n_0,\fp}\otimes\tT^\circ_{n_1,\fp}f)(s,\sh^\circ_\uparrow(s))
      \end{align*}
      holds.
\end{enumerate}
Here, $\int_{\ol\rP^{\bullet,\bullet}}^\fT$ denotes the $\fT$-trace map in Definition \ref{de:trace}; and $\sh^\circ_\uparrow$, $\sh^\bullet_\uparrow$, and $\sh^\bullet_\downarrow$ are maps in Notation \ref{no:ns_uniformization}.
\end{theorem}

The intersection number in (3) is the actual one that is responsible for the first explicit reciprocity law which will be discussed in \S\ref{ss:first_reciprocity}.

\begin{proof}
We first show (3) assuming (1) and (2). By (1), (2), and Lemma \ref{le:enumeration_odd_0}, we have for $f\in L[\Sh(\rV^\circ_{n_0},\rK^{p\circ}_{n_0}\rK^\circ_{n_0,p})]\otimes_L L[\Sh(\rV^\circ_{n_1},\rK^{p\circ}_{n_1}\rK^\circ_{n_1,p})]$,
\begin{align*}
&\quad\int_{\ol\rP^{\bullet,\bullet}}^\fT\cl(\rP^\bullet_\graph)\cup
\(\inc^{\bullet,\dag}_!(\tT^{\bullet\circ}_{n_0,\fp}\otimes\tI^\circ_{n_1,\fp}f)+
(p+1)^2\inc^{\bullet,\bullet}_!(\tT^{\bullet\circ}_{n_0,\fp}\otimes\tT^{\bullet\circ}_{n_1,\fp}f)\) \\
&=\sum_{s\in\Sh(\rV^\circ_n,\rK^{p\circ}_n\rK^\bullet_{\sp,p})}
(\tT^{\bullet\circ}_{n_0,\fp}\otimes(\tT^{\bullet\circ}_{n_1,\fp}\circ\tI^\circ_{n_1,\fp})f)(\sh^\bullet_\downarrow(s),\sh^\bullet_\uparrow(s)) \\
&\quad +\sum_{s\in\Sh(\rV^\circ_n,\rK^{p\circ}_n\rK^\bullet_{\sp,p})}
(\tT^{\bullet\circ}_{n_0,\fp}\otimes((p+1)^2\tT^\bullet_{n_1,\fp}\circ\tT^{\bullet\circ}_{n_1,\fp})f)
(\sh^\bullet_\downarrow(s),\sh^\bullet_\uparrow(s)) \\
&=\sum_{s\in\Sh(\rV^\circ_n,\rK^{p\circ}_n\rK^\bullet_{\sp,p})}
(\tT^{\bullet\circ}_{n_0,\fp}\otimes(\tT^{\bullet\circ}_{n_1,\fp}\circ\tI^\circ_{n_1,\fp})f)(\sh^\bullet_\downarrow(s),\sh^\bullet_\uparrow(s)) \\
&\quad +\sum_{s\in\Sh(\rV^\circ_n,\rK^{p\circ}_n\rK^\bullet_{\sp,p})}
(\tT^{\bullet\circ}_{n_0,\fp}\otimes(\tT^{\bullet\circ}_{n_1,\fp}\circ\tT^\circ_{n_1,\fp}-\tT^{\bullet\circ}_{n_1,\fp}\circ\tI^\circ_{n_1,\fp})f)
(\sh^\bullet_\downarrow(s),\sh^\bullet_\uparrow(s)) \\
&=\sum_{s\in\Sh(\rV^\circ_n,\rK^{p\circ}_n\rK^\bullet_{\sp,p})}
(\tT^{\bullet\circ}_{n_0,\fp}\otimes(\tT^{\bullet\circ}_{n_1,\fp}\circ\tT^\circ_{n_1,\fp})f)(\sh^\bullet_\downarrow(s),\sh^\bullet_\uparrow(s))
\end{align*}
which, by Lemma \ref{le:ns_reciprocity} below, equals
\begin{align*}
\sum_{s\in\Sh(\rV^\circ_n,\rK^{p\circ}_n\rK^\circ_{n,p})}(\tI^\circ_{n_0,\fp}\otimes\tT^\circ_{n_1,\fp}f)(s,\sh^\circ_\uparrow(s)).
\end{align*}
Thus, (3) is proved.

Now we consider (1) and (2) simultaneously. Similar to the maps $\inc^\bullet_!$ and $\inc^\dag_!$ in Construction \ref{cs:ns_incidence}, we have maps
\begin{align*}
\inc^\bullet_\alpha &\colon L[\Sh(\rV^\circ_{n_\alpha},\rK^{p\circ}_{n_\alpha}\rK^\bullet_{n_\alpha,p})]
\to\rH^{2(r_\alpha+\alpha-1)}_{\fT,c}(\ol\rM^\bullet_{n_\alpha},L(r_\alpha+\alpha-1)),\\
\inc^\dag_\alpha & \colon L[\Sh(\rV^\circ_{n_\alpha},\rK^{p\circ}_{n_\alpha}\rK^\circ_{n_\alpha,p})]
\to\rH^{2(r_\alpha+\alpha-1)}_{\fT,c}(\ol\rM^\bullet_{n_\alpha},L(r_\alpha+\alpha-1)),
\end{align*}
for $\alpha=0,1$. Note that we now take $\rH_{\fT,c}$ for the target of the maps rather than $\rH_\fT$. Moreover, in the calculation below, we will frequently use the following formula for intersection number pairings: for a finite morphism $i\colon X\to Y$ of smooth schemes over an algebraically closed field, and proper smooth subschemes $X'$ of $X$ and $Y'$ of $Y$, we have
\[
\langle X_\graph, X'\times Y'\rangle_{X\times Y}=\langle X'_\graph,X'\times Y'\rangle_{X'\times Y}=\langle i_*X',Y'\rangle_Y
\]
where $X_\graph$ and $X'_\graph$ denote by the graphs of $i$ and $i\res X'$, respectively. The proof for (1) and (2) differs by the parity of $n$.

We first consider the case where $n=n_0$ is even. By Lemma \ref{le:ns_functoriality_source}(1) and Proposition \ref{pr:ns_uniformization_basic}, $\sh^\bullet_\downarrow$ is an isomorphism. Take a point $s^\bullet_n\in\Sh(\rV^\circ_n,\rK^{p\circ}_n\rK^\bullet_{n,p})$. Let $s^\bullet$ be the unique element in $\Sh(\rV^\circ_n,\rK^{p\circ}_n\rK^\bullet_{\sp,p})$ such that $s^\bullet_n=\sh^\bullet_\downarrow(s^\bullet)$, and put $s^\bullet_{n+1}\coloneqq\sh^\bullet_\uparrow(s^\bullet)$. By (the last assertion in) Proposition \ref{pr:ns_functoriality_dl}, we have
\[
\rm^\bullet_{\uparrow!}\inc^\bullet_0(1_{s^\bullet_n})=\inc^\bullet_1(1_{s^\bullet_{n+1}}).
\]
For (1), we have for every $s'_{n+1}\in\Sh(\rV^\circ_{n+1},\rK^{p\circ}_{n+1}\rK^\circ_{n+1,p})$ the identity
\begin{align*}
\int_{\ol\rP^{\bullet,\bullet}}^\fT\cl(\rP^\bullet_\graph)\cup\inc^{\bullet,\dag}_!(1_{(s^\bullet_n,s'_{n+1})})
&=\int_{\ol\rM^\bullet_{n+1}}^\fT\(\rm^\bullet_{\uparrow!}\inc^\bullet_0(1_{s^\bullet_n})\)\cup\inc^\dag_1(1_{s'_{n+1}}) \\
&=\int_{\ol\rM^\bullet_{n+1}}^\fT\inc^\bullet_0(1_{s^\bullet_{n+1}})\cup\inc^\dag_1(1_{s'_{n+1}}).
\end{align*}
Thus, (1) follows from Proposition \ref{pr:ns_incidence_odd}. For (2), we have for every $s'_{n+1}\in\Sh(\rV^\circ_{n+1},\rK^{p\circ}_{n+1}\rK^\bullet_{n+1,p})$ the identity
\begin{align*}
\int_{\ol\rP^{\bullet,\bullet}}^\fT\cl(\rP^\bullet_\graph)\cup\inc^{\bullet,\bullet}_!(1_{(s^\bullet_n,s'_{n+1})})
&=\int_{\ol\rM^\bullet_{n+1}}^\fT\(\rm^\bullet_{\uparrow!}\inc^\bullet_0(1_{s^\bullet_n})\)\cup\inc^\bullet_1(1_{s'_{n+1}}) \\
&=\int_{\ol\rM^\bullet_{n+1}}^\fT\inc^\bullet_0(1_{s^\bullet_{n+1}})\cup\inc^\bullet_1(1_{s'_{n+1}}).
\end{align*}
Thus, (2) follows from Proposition \ref{pr:ns_incidence_odd}.

We then consider the case where $n=n_1$ is odd. Take a point $s^\bullet_{n+1}\in\Sh(\rV^\circ_{n+1},\rK^{p\circ}_{n+1}\rK^\bullet_{n+1,p})$. By Proposition \ref{pr:ns_functoriality2}, Proposition \ref{pr:ns_functoriality_dl}, and Proposition \ref{pr:ns_uniformization_basic}, we have
\[
\rm^{\bullet*}_\uparrow\inc^\bullet_0(1_{s^\bullet_{n+1}})=
\inc^\bullet_1(\sh^\bullet_{\downarrow!}\sh^{\bullet*}_\uparrow1_{s^\bullet_{n+1}}).
\]
For (1), we have for every $s'_n\in\Sh(\rV^\circ_n,\rK^{p\circ}_n\rK^\circ_{n,p})$ the identity
\begin{align*}
\int_{\ol\rP^{\bullet,\bullet}}^\fT\cl(\rP^\bullet_\graph)\cup\inc^{\bullet,\dag}_!(1_{(s^\bullet_{n+1},s'_n)})
&=\int_{\ol\rM^\bullet_n}^\fT\(\rm^{\bullet*}_\uparrow\inc^\bullet_0(1_{s^\bullet_{n+1}})\)\cup\inc^\dag_1(1_{s'_n}) \\
&=\int_{\ol\rM^\bullet_n}^\fT\inc^\bullet_1(\sh^\bullet_{\downarrow!}\sh^{\bullet*}_\uparrow1_{s^\bullet_{n+1}})\cup\inc^\dag_1(1_{s'_n}).
\end{align*}
Thus, (1) follows from Proposition \ref{pr:ns_incidence_odd}. For (2), we have for every $s'_n\in\Sh(\rV^\circ_n,\rK^{p\circ}_n\rK^\bullet_{n,p})$ the identity
\begin{align*}
\int_{\ol\rP^{\bullet,\bullet}}^\fT\cl(\rP^\bullet_\graph)\cup\inc^{\bullet,\bullet}_!(1_{(s^\bullet_{n+1},s'_n)}) &=\int_{\ol\rM^\bullet_n}^\fT\(\rm^{\bullet*}_\uparrow\inc^\bullet_0(1_{s^\bullet_{n+1}})\)\cup\inc^\bullet_1(1_{s'_n})  \\ &=\int_{\ol\rM^\bullet_n}^\fT\inc^\bullet_1(\sh^\bullet_{\downarrow!}\sh^{\bullet*}_\uparrow1_{s^\bullet_{n+1}})\cup\inc^\bullet_1(1_{s'_n}).
\end{align*}
Thus, (1) follows from Proposition \ref{pr:ns_incidence_odd}.

The theorem is proved.
\end{proof}

\begin{lem}\label{le:ns_reciprocity}
For every $f\in L[\Sh(\rV^\circ_{n_0},\rK^{p\circ}_{n_0}\rK^\bullet_{n_0,p})]\otimes_LL[\Sh(\rV^\circ_{n_1},\rK^{p\circ}_{n_1}\rK^\circ_{n_1,p})]$, we have
\begin{align*}
\sum_{s\in\Sh(\rV^\circ_n,\rK^{p\circ}_n\rK^\bullet_{\sp,p})}(\tT^{\bullet\circ}_{n_1,\fp}f)(\sh^\bullet_\downarrow(s),\sh^\bullet_\uparrow(s))
=\sum_{s\in\Sh(\rV^\circ_n,\rK^{p\circ}_n\rK^\circ_{n,p})}(\tT^{\circ\bullet}_{n_0,\fp}f)(s,\sh^\circ_\uparrow(s)).
\end{align*}
\end{lem}

\begin{proof}
There are two cases.

When $n$ is even, by Lemma \ref{le:ns_functoriality_source_dagger}(1) and Proposition \ref{pr:ns_uniformization_basic}, we have
\begin{align*}
\sum_{s\in\Sh(\rV^\circ_n,\rK^{p\circ}_n\rK^\bullet_{\sp,p})}(\tT^{\bullet\circ}_{n_1,\fp}f)(\sh^\bullet_\downarrow(s),\sh^\bullet_\uparrow(s))
&=\sum_{s\in\Sh(\rV^\circ_n,\rK^{p\circ}_n\rK^\dag_{\sp,p})}
f(\sh^{\dag\bullet}_n(\sh^\dag_\downarrow(s)),\sh^{\dag\circ}_{n+1}(\sh^\dag_\uparrow(s))) \\
&=\sum_{s\in\Sh(\rV^\circ_n,\rK^{p\circ}_n\rK^\dag_{\sp,p})}
f(\sh^{\dag\bullet}_n(\sh^\dag_\downarrow(s)),\sh^\circ_\uparrow(\sh^{\dag\circ}_{n}(\sh^\dag_\downarrow(s)))),
\end{align*}
which, by Lemma \ref{le:ns_functoriality_source}(1), Definition \ref{de:ns_definite_special3}, and Proposition \ref{pr:ns_uniformization_basic}, equals
\[
\sum_{s\in\Sh(\rV^\circ_n,\rK^{p\circ}_n\rK^\dag_{n,p})}f(\sh^{\dag\bullet}_n(s),\sh^\circ_\uparrow(\sh^{\dag\circ}_n(s)))
=\sum_{s\in\Sh(\rV^\circ_n,\rK^{p\circ}_n\rK^\circ_{n,p})}(\tT^{\circ\bullet}_{n_0,\fp}f)(s,\sh^\circ_\uparrow(s)).
\]

When $n$ is odd, by Definition \ref{de:ns_definite_special3} and Proposition \ref{pr:ns_uniformization_basic}, we have
\begin{align*}
\sum_{s\in\Sh(\rV^\circ_n,\rK^{p\circ}_n\rK^\bullet_{\sp,p})}(\tT^{\bullet\circ}_{n_1,\fp}f)(\sh^\bullet_\downarrow(s),\sh^\bullet_\uparrow(s))
&=\sum_{s\in\Sh(\rV^\circ_n,\rK^{p\circ}_n\rK^\dag_{\sp,p})}
f(\sh^{\dag\circ}_n(\sh^\dag_\downarrow(s)),\sh^\bullet_\uparrow(\sh^{\dag\bullet}_\sp(s))) \\
&=\sum_{s\in\Sh(\rV^\circ_n,\rK^{p\circ}_n\rK^\dag_{\sp,p})}
f(\sh^{\dag\circ}_n(\sh^\dag_\downarrow(s)),\sh^{\dag\bullet}_{n+1}(\sh^\dag_\uparrow(s))),
\end{align*}
which, by Lemma \ref{le:ns_functoriality_source_dagger}(2) and Proposition \ref{pr:ns_uniformization_basic}, equals
\[
\sum_{s\in\Sh(\rV^\circ_n,\rK^{p\circ}_n\rK^\circ_{n,p})}(\tT^{\circ\bullet}_{n_0,\fp}f)(s,\sh^\circ_\uparrow(s)).
\]

The lemma is proved.
\end{proof}

\begin{construction}\label{cs:ns_nabla_product}
We constructs maps
\begin{align*}
\begin{dcases}
\Inc_{\circ,\dag}^*\colon\rH^{2n}_\fT(\ol\rQ^{(0)},L(n))\to\rH^{2n}_\fT(\ol\rQ^{\circ,\bullet},L(n))\xrightarrow{\sigma_!}
\rH^{2n}_\fT(\ol\rP^{\circ,\bullet},L(n)) \\
\qquad\qquad\xrightarrow{\inc_{\circ,\dag}^*}
L[\Sh(\rV^\circ_{n_0},\obj_{n_0}\rK^\circ_{n_0,p})]\otimes_L L[\Sh(\rV^\circ_{n_1},\obj_{n_1}\rK^\circ_{n_1,p})], \\
\Inc_{\circ,\bullet}^*\colon\rH^{2n}_\fT(\ol\rQ^{(0)},L(n))\to\rH^{2n}_\fT(\ol\rQ^{\circ,\bullet},L(n))\xrightarrow{\sigma_!}
\rH^{2n}_\fT(\ol\rP^{\circ,\bullet},L(n)) \\
\qquad\qquad\xrightarrow{\inc_{\circ,\bullet}^*}
L[\Sh(\rV^\circ_{n_0},\obj_{n_0}\rK^\circ_{n_0,p})]\otimes_L L[\Sh(\rV^\circ_{n_1},\obj_{n_1}\rK^\bullet_{n_1,p})], \\
\Inc_{\bullet,\dag}^*\colon\rH^{2n}_\fT(\ol\rQ^{(0)},L(n))\to\rH^{2n}_\fT(\ol\rQ^{\bullet,\bullet},L(n))\xrightarrow{\sigma_!}
\rH^{2n}_\fT(\ol\rP^{\bullet,\bullet},L(n)) \\
\qquad\qquad\xrightarrow{\inc_{\bullet,\dag}^*}
L[\Sh(\rV^\circ_{n_0},\obj_{n_0}\rK^\bullet_{n_0,p})]\otimes_L L[\Sh(\rV^\circ_{n_1},\obj_{n_1}\rK^\circ_{n_1,p})], \\
\Inc_{\bullet,\bullet}^*\colon\rH^{2n}_\fT(\ol\rQ^{(0)},L(n))\to\rH^{2n}_\fT(\ol\rQ^{\bullet,\bullet},L(n))\xrightarrow{\sigma_!}
\rH^{2n}_\fT(\ol\rP^{\bullet,\bullet},L(n)) \\
\qquad\qquad\xrightarrow{\inc_{\bullet,\bullet}^*}
L[\Sh(\rV^\circ_{n_0},\obj_{n_0}\rK^\bullet_{n_0,p})]\otimes_L L[\Sh(\rV^\circ_{n_1},\obj_{n_1}\rK^\bullet_{n_1,p})].
\end{dcases}
\end{align*}
Define the map
\begin{align*}
\nabla\colon\rH^{2n}_\fT(\ol\rQ^{(0)},L(n))\to
L[\Sh(\rV^\circ_{n_0},\obj_{n_0}\rK^\circ_{n_0,p})]\otimes_L L[\Sh(\rV^\circ_{n_1},\obj_{n_1}\rK^\circ_{n_1,p})]
\end{align*}
to be the sum of the following four maps
\begin{align*}
(\tI^\circ_{n_0,\fp}\otimes\tI^\circ_{n_1,\fp})\circ\Inc_{\circ,\dag}^*&,\quad
(p+1)^2(\tI^\circ_{n_0,\fp}\otimes\tT^{\circ\bullet}_{n_1,\fp})\circ\Inc_{\circ,\bullet}^*,\\
(p+1)(\tT^{\circ\bullet}_{n_0,\fp}\otimes\tI^\circ_{n_1,\fp})\circ\Inc_{\bullet,\dag}^*&,\quad
(p+1)^3(\tT^{\circ\bullet}_{n_0,\fp}\otimes\tT^{\circ\bullet}_{n_1,\fp})\circ\Inc_{\bullet,\bullet}^*.
\end{align*}
\end{construction}

At last, we recall the construction of potential map from \cite{Liu2}*{\S2.2}. For $r\in\dZ$, put
\[
B^r(\rQ,L)\coloneqq\Ker\(\delta_0^*\colon\rH^{2r}_\fT(\ol\rQ^{(0)},L(r))\to\rH^{2r}_\fT(\ol\rQ^{(1)},L(r))\)
\]
and
\[
B_r(\rQ,L)\coloneqq\coker\(\delta_{1!}\colon\rH^{2(2n-r-2)}_\fT(\ol\rQ^{(1)},L(2n-r-2))
\to\rH^{2(2n-r-1)}_\fT(\ol\rQ^{(0)},L(2n-r-1))\).
\]
Here, in our case,
\begin{align*}
\delta_0^*&=(\delta^{\circ,\dag}_{\circ,\bullet})^*-(\delta^{\circ,\dag}_{\circ,\circ})^*
+(\delta^{\dag,\bullet}_{\bullet,\bullet})^*-(\delta^{\dag,\bullet}_{\circ,\bullet})^*
+(\delta^{\bullet,\dag}_{\bullet,\circ})^*-(\delta^{\bullet,\dag}_{\bullet,\bullet})^*
+(\delta^{\dag,\circ}_{\bullet,\circ})^*-(\delta^{\dag,\circ}_{\circ,\circ})^*
+(\delta^{\dag,\dag}_{\bullet,\bullet})^*-(\delta^{\dag,\dag}_{\circ,\circ})^*, \\
\delta_{1!}&=(\delta^{\circ,\dag}_{\circ,\bullet})_!-(\delta^{\circ,\dag}_{\circ,\circ})_!
+(\delta^{\dag,\bullet}_{\bullet,\bullet})_!-(\delta^{\dag,\bullet}_{\circ,\bullet})_!
+(\delta^{\bullet,\dag}_{\bullet,\circ})_!-(\delta^{\bullet,\dag}_{\bullet,\bullet})_!
+(\delta^{\dag,\circ}_{\bullet,\circ})_!-(\delta^{\dag,\circ}_{\circ,\circ})_!
+(\delta^{\dag,\dag}_{\bullet,\bullet})_!-(\delta^{\dag,\dag}_{\circ,\circ})_!.
\end{align*}
We define $B^r(\rQ,L)^0$ and $B_{2n-r-1}(\rQ,L)_0$ to be the kernel and the cokernel of the tautological map
\[
B^r(\rQ,L)\to B_{2n-r-1}(\rQ,L),
\]
respectively. By \cite{Liu2}*{Lemma~2.4}, the composite map
\[
\rH^{2(r-1)}_\fT(\ol\rQ^{(0)},L(r-1))\xrightarrow{\delta_0^*}
\rH^{2(r-1)}_\fT(\ol\rQ^{(1)},L(r-1))\xrightarrow{\delta_{1!}}
\rH^{2r}_\fT(\ol\rQ^{(0)},L(r))
\]
factors through a unique map
\[
B_{2n-r}(\rQ,L)_0\to B^r(\rQ,L)^0
\]
in $\Fun(\fK(\rV^\circ_n)^p\times\fK(\rV^\circ_{n+1})^p,\Mod(L[\Gal(\ol\dF_p/\dF_p^\Phi)]))$. Put
\[
C_r(\rQ,L)\coloneqq B_r(\rQ,L)_0^{\Gal(\ol\dF_p/\dF_p^\Phi)}, \quad
C^r(\rQ,L)\coloneqq B^r(\rQ,L)^0_{\Gal(\ol\dF_p/\dF_p^\Phi)}.
\]
Then we obtain the \emph{potential map}
\begin{align}\label{eq:ns_potential}
\Delta^r\colon C_{2n-r}(\rQ,L)\to C^r(\rQ,L)
\end{align}
in $\Fun(\fK(\rV^\circ_n)^p\times\fK(\rV^\circ_{n+1})^p,\Mod(L))$.\footnote{In \cite{Liu2}, $C^r(\rQ,L)$ and $C_r(\rQ,L)$ are denoted by $A^r(\rQ,L)^0$ and $A_r(\rQ,L)_0$, respectively.} We will be most interested in the case where $r=n$.

\begin{remark}\label{re:ns_nabla_product}
By the descriptions of the Galois actions in Construction \ref{cs:ns_uniformization1} and Construction \ref{cs:ns_uniformization2}, the map $\nabla$ in Construction \ref{cs:ns_nabla_product} factors through the quotient map
\[
\rH^{2n}_\fT(\ol\rQ^{(0)},L(n))\to\rH^{2n}_\fT(\ol\rQ^{(0)},L(n))_{\Gal(\ol\dF_p/\dF_p^\Phi)},
\]
hence restricts to a map
\begin{align*}
\nabla\colon C^n(\rQ,L)\to
L[\Sh(\rV^\circ_{n_0},\obj_{n_0}\rK^\circ_{n_0,p})]\otimes_L L[\Sh(\rV^\circ_{n_1},\obj_{n_1}\rK^\circ_{n_1,p})]
\end{align*}
in $\Fun(\fK(\rV^\circ_n)^p\times\fK(\rV^\circ_{n+1})^p,\Mod(L))$, via the canonical map $C^n(\rQ,L)\to\rH^{2n}_\fT(\ol\rQ^{(0)},L(n))_{\Gal(\ol\dF_p/\dF_p^\Phi)}$.
\end{remark}

\section{Tate classes and arithmetic level-raising}
\label{ss:6}

In this section, we study two important arithmetic properties of semistable moduli schemes introduced in \S\ref{ss:ns}. The first is the existence of Tate cycles when the rank is odd, studied in \S\ref{ss:tate}. The second is the arithmetic level-raising when the rank is even, studied in \S\ref{ss:raising} and \S\ref{ss:raising2}. In \S\ref{ss:preliminaries}, we collect some preliminaries on automorphic representations and their motives.

Let $N\geq 2$ be an integer with $r\coloneqq\lfloor\tfrac{N}{2}\rfloor$.

\subsection{Preliminaries on automorphic representations}
\label{ss:preliminaries}

In this subsection, we consider
\begin{itemize}[label={\ding{109}}]
  \item a relevant representation $\Pi$ of $\GL_N(\dA_F)$ (Definition \ref{de:relevant}),

  \item a strong coefficient field $E\subseteq\dC$ of $\Pi$ (Definition \ref{de:weak_field}),

  \item a finite set $\Sigma^+_\mnm$ of nonarchimedean places of $F^+$ containing $\Sigma^+_\Pi$ (Notation \ref{no:satake}),

  \item a (possibly empty) finite set $\Sigma^+_\lr$ of nonarchimedean places of $F^+$ that are inert in $F$,\footnote{Here, the subscript ``lr'' standards for ``level-raising''.} strongly disjoint from $\Sigma^+_\mnm$ (Definition \ref{de:strongly_disjoint}),

  \item a finite set $\Sigma^+$ of nonarchimedean places of $F^+$ containing $\Sigma^+_\mnm\cup\Sigma^+_\lr$.
\end{itemize}
We then have, by Construction \ref{cs:satake_hecke}, the homomorphism
\[
\phi_\Pi\colon\dT^{\Sigma^+}_N\to O_E.
\]
For every prime $\lambda$ of $E$, we have a continuous homomorphism
\begin{align*}
\rho_{\Pi,\lambda}\colon\Gamma_F\to\GL_N(E_\lambda)
\end{align*}
from Proposition \ref{pr:galois}(2) and Definition \ref{de:weak_field}, such that $\rho_{\Pi,\lambda}^\tc$ and $\rho_{\Pi,\lambda}^\vee(1-N)$ are conjugate.

We choose
\begin{itemize}[label={\ding{109}}]
  \item a prime $\lambda$ of $E$, whose underlying rational prime $\ell$ satisfies $\Sigma^+_\mnm\cap\Sigma^+_\ell=\emptyset$ and $\ell\nmid\|v\|(\|v\|^2-1)$ for every $v\in\Sigma^+_\lr$,

  \item a positive integer $m$,

  \item a standard definite hermitian space $\rV^\circ_N$ of rank $N$ over $F$, together with a self-dual $\prod_{v\not\in\Sigma^+_\infty\cup\Sigma^+_\mnm\cup\Sigma^+_\lr}O_{F_v}$-lattice $\Lambda^\circ_N$ in $\rV^\circ_N\otimes_F\dA_F^{\Sigma^+_\infty\cup\Sigma^+_\mnm\cup\Sigma^+_\lr}$, satisfying that $(\rV^\circ_N)_v$ is not split for $v\in\Sigma^+_\lr$ when $N$ is even,

  \item an object $\rK^\circ_N\in\fK(\rV^\circ_N)$ of the form
      \[
      \rK^\circ_N=\prod_{v\in\Sigma^+_\mnm\cup\Sigma^+_\lr}(\rK^\circ_N)_v\times
      \prod_{v\not\in\Sigma^+_\infty\cup\Sigma^+_\mnm\cup\Sigma^+_\lr}\rU(\Lambda^\circ_N)(O_{F^+_v}),
      \]
      satisfying that when $N$ is even, $(\rK^\circ_N)_v$ is a transferable open compact subgroup of $\rU(\rV^\circ_N)(F^+_v)$ (Definition \ref{de:transferable})\footnote{By Lemma \ref{le:transferable}(3), every sufficiently small $(\rK^\circ_N)_v$ is transferable. So the readers may ignore this technical requirement.} for $v\in\Sigma^+_\mnm$ and is a special maximal subgroup of $\rU(\rV^\circ_N)(F^+_v)$ for $v\in\Sigma^+_\lr$,

  \item a special inert prime (Definition \ref{de:special_inert}) $\fp$ of $F^+$ (with the underlying rational prime $p$) satisfying
      \begin{description}
        \item[(P1)] $\Sigma^+$ does not contain $p$-adic places;

        \item[(P2)] $\ell$ does not divide $p(p^2-1)$;

        \item[(P3)] there exists a CM type $\Phi$ containing $\tau_\infty$ as in the initial setup of \S\ref{ss:ns} satisfying $\dQ_p^\Phi=\dQ_{p^2}$;

        \item[(P4)] if $N$ is even, then $P_{\balpha(\Pi_\fp)}\modulo\lambda^m$ is level-raising special at $\fp$ (Definition \ref{de:satake_condition});

            if $N$ is odd, then $P_{\balpha(\Pi_\fp)}\modulo\lambda$ is Tate generic at $\fp$ (Definition \ref{de:satake_condition});

        \item[(P5)] $P_{\balpha(\Pi_\fp)}\modulo\lambda$ is intertwining generic at $\fp$ (Definition \ref{de:satake_condition});

        \item[(P6)] if $N$ is even, the natural map
            \[
            \frac{(O_E/\lambda^m)[\Sh(\rV^\circ_N,\rK^\circ_N)]}{\dT^{\Sigma^+\cup\Sigma^+_p}_N\cap\Ker\phi_\Pi}
            \to\frac{(O_E/\lambda^m)[\Sh(\rV^\circ_N,\rK^\circ_N)]}{\Ker\phi_\Pi}
            \]
            is an isomorphism;
      \end{description}
      (So we can and will apply the setup in \S\ref{ss:ns} to the datum $(\rV^\circ_N,\{\Lambda^\circ_{N,\fq}\}\res_{\fq\mid p})$.)

  \item remaining data in \S\ref{ss:ns_initial} with $\dQ_p^\Phi=\dQ_{p^2}$;

  \item data as in Construction \ref{cs:ns_uniformization2}, which in particular give the open compact subgroup $\rK^\bullet_p$; and

  \item an indefinite uniformization datum $(\rV'_N,\tj_N,\{\Lambda'_{\fq,N}\}_{\fq\mid p})$ for $\rV^\circ_N$ as in Definition \ref{de:ns_uniformization_data}.
\end{itemize}

Put $\rK^{p\circ}_N\coloneqq(\rK^\circ_N)^p$ and $\rK^\bullet_N\coloneqq\rK^{p\circ}_N\times\rK^\bullet_p$. As in \S\ref{ss:ns_weight}, we put $\rX^?_N\coloneqq\rX^?_\fp(\rV^\circ_N,\rK^{p\circ}_N)$ for meaningful pairs $(\rX,?)\in\{\bM,\rM,\rB,\rS\}\times\{\;,\eta,\circ,\bullet,\dag\}$. Let $(\rE^{p,q}_s,\rd^{p,q}_s)$ be the weight spectral sequence abutting to the cohomology $\rH^{p+q}_\fT(\ol\rM_N,\rR\Psi O_\lambda(r))$ from \S\ref{ss:ns_weight}.

\begin{remark}
By Construction \ref{cs:satake_hecke} and (P2) (namely, $\ell\neq p$), we know that $P_{\balpha(\Pi_\fp)}$ is a polynomial with coefficients in $O_\lambda$.
\end{remark}

\begin{remark}
Note that when $N=2$, (P2) and (P4) together imply (P5).
\end{remark}

\begin{notation}\label{no:single_ideal}
We introduce the following ideals of $\dT^{\Sigma^+\cup\Sigma^+_p}_N$
\begin{align*}
\begin{dcases}
\fm\coloneqq\dT^{\Sigma^+\cup\Sigma^+_p}_N\cap\Ker\(\dT^{\Sigma^+}_N\xrightarrow{\phi_\Pi}O_E\to O_E/\lambda\),\\
\fn\coloneqq\dT^{\Sigma^+\cup\Sigma^+_p}_N\cap\Ker\(\dT^{\Sigma^+}_N\xrightarrow{\phi_\Pi}O_E\to O_E/\lambda^m\).
\end{dcases}
\end{align*}
\end{notation}

We then introduce the following assumptions.

\begin{assumption}\label{as:single_vanishing}
We have $\rH^i_\fT(\ol\rM_N,\rR\Psi O_\lambda)_\fm=0$ for $i\neq N-1$, and that $\rH^{N-1}_\fT(\ol\rM_N,\rR\Psi O_\lambda)_\fm$ is a finite free $O_\lambda$-module.
\end{assumption}

\begin{remark}\label{re:single_vanishing}
Assumption \ref{as:single_vanishing} holds, for example, when the composite homomorphism $\dT^{\Sigma^+}_N\xrightarrow{\phi_\Pi}O_E\to O_E/\lambda$ is cohomologically generic (Definition \ref{de:generic}). This follows from Lemma \ref{le:ns_pbc} and the universal coefficient theorem.
\end{remark}

\begin{assumption}\label{as:single_irreducible}
The Galois representation $\rho_{\Pi,\lambda}$ is residually absolutely irreducible.
\end{assumption}

\begin{remark}\label{re:single_irreducible}
Under Assumption \ref{as:single_irreducible}, we obtain a homomorphism
\[
\bar\rho_{\Pi,\lambda}\colon\Gamma_F\to\GL_N(O_\lambda/\lambda)
\]
from the residual homomorphism of $\rho_{\Pi,\lambda}$, which is unique to conjugation, absolutely irreducible, and $(1-N)$-polarizable (Definition \ref{de:polarization}). Applying Construction \ref{cs:extension}, we obtain an extension
\begin{align*}
\bar\rho_{\Pi,\lambda,+}\colon\Gamma_{F^+}\to\sG_N(O_\lambda/\lambda)
\end{align*}
of $\bar\rho_{\Pi,\lambda}$.
\end{remark}

We now fix an isomorphism $\iota_\ell\colon\dC\simeq\ol\dQ_\ell$ that induces the prime $\lambda$ of $E$, till the end of this section.

\begin{definition}\label{de:single_congruent}
We say that a standard pair $(\rV,\pi)$ (Definition \ref{de:standard_pair}) with $\dim_F\rV=N$ is \emph{$\Pi$-congruent} (outside $\Sigma^+\cup\Sigma^+_p$) if for every nonarchimedean place $v$ of $F^+$ not in $\Sigma^+\cup\Sigma^+_p\cup\Sigma^+_\ell$, $\pi_v$ is unramified; and the two homomorphisms $\iota_\ell\phi_{\balpha(\BC(\pi_v))}$ and $\iota_\ell\phi_{\balpha(\Pi_v)}$ from $\dT_{N,v}$ to $\ol\dQ_\ell$, which in fact take values in $\ol\dZ_\ell$, coincide in $\ol\dF_\ell$.
\end{definition}

\begin{lem}\label{le:single_intertwining}
The two maps
\begin{align*}
\tT^{\bullet\circ}_{N,\fp}&\colon O_E[\Sh(\rV^\circ_N,\rK^\circ_N)]_\fm\to O_E[\Sh(\rV^\circ_N,\rK^\bullet_N)]_\fm \\
\tT^{\circ\bullet}_{N,\fp}&\colon O_E[\Sh(\rV^\circ_N,\rK^\bullet_N)]_\fm\to O_E[\Sh(\rV^\circ_N,\rK^\circ_N)]_\fm
\end{align*}
are both isomorphisms, where $\tT^{\bullet\circ}_{N,\fp}$ and $\tT^{\circ\bullet}_{N,\fp}$ are introduced in Definition \ref{de:ns_hecke}.
\end{lem}

\begin{proof}
By Proposition \ref{pr:enumeration_odd}(1) (resp.\ Proposition \ref{pr:enumeration_even}(1)) when $N$ is odd (resp.\ even) and (P5), we know that the endomorphism $\tI^\circ_{N,\fp}=\tT^{\circ\bullet}_{N,\fp}\circ\tT^{\bullet\circ}_{N,\fp}$ of $O_E[\Sh(\rV^\circ_N,\rK^\circ_N)]_\fm$ is an isomorphism. Thus, it suffices to show that the free $O_\lambda$-modules $O_E[\Sh(\rV^\circ_N,\rK^\circ_N)]_\fm$ and $O_E[\Sh(\rV^\circ_N,\rK^\bullet_N)]_\fm$ have the same rank. We show that $O_E[\Sh(\rV^\circ_N,\rK^\circ_N)]_\fm\otimes_{O_\lambda}\ol\dQ_\ell$ and $O_E[\Sh(\rV^\circ_N,\rK^\bullet_N)]_\fm\otimes_{O_\lambda}\ol\dQ_\ell$ have the same dimension. We have
\begin{align*}
O_E[\Sh(\rV^\circ_N,\rK^\circ_N)]_\fm\otimes_{O_\lambda}\ol\dQ_\ell&\simeq\bigoplus_{\pi}m(\pi)\cdot\pi^{\rK^\circ_N},\\
O_E[\Sh(\rV^\circ_N,\rK^\bullet_N)]_\fm\otimes_{O_\lambda}\ol\dQ_\ell&\simeq\bigoplus_{\pi}m(\pi)\cdot\pi^{\rK^\bullet_N},
\end{align*}
where $\pi$ runs over all irreducible admissible representations of $\rU(\rV^\circ_N)(\dA_{F^+})$ with coefficients in $\ol\dQ_\ell$ such that $(\rV^\circ_N,\iota_\ell^{-1}\pi)$ is a $\Pi$-congruent standard pair (Definition \ref{de:single_congruent}); and $m(\pi)$ denotes the automorphic multiplicity of $\pi$.\footnote{Although we know that $m(\pi)=1$ by Proposition \ref{pr:arthur}(2), we do not need this fact here.} It suffices to show that if in the second direct sum $\pi_\fp^{\rK^\bullet_N}\neq\{0\}$, which has to be of dimension one since $\rK^\bullet_N$ is special maximal, then $\pi_\fp^{\rK^\circ_N}\neq\{0\}$ as well. Moreover, the Satake parameter $\balpha$ of $\pi_\fp$ does not contain the pair $\{-1,-1\}$ (resp.\ $\{-p,-p^{-1}\}$) when $N$ is even (resp.\ odd) by (P5). Let $\pi'_\fp$ be the unique constituent of the principal series of $\balpha$ such that $(\pi'_\fp)^{\rK^\circ_N}\neq\{0\}$, then by Proposition \ref{pr:enumeration_odd}(1) (resp.\ Proposition \ref{pr:enumeration_even}(1)) when $N$ is odd (resp.\ even) again, we see that $(\pi'_\fp)^{\rK^\bullet_N}\neq\{0\}$. Thus, we must have $\pi_\fp=\pi'_\fp$ as $\rK^\bullet_N$ is special maximal. The lemma follows.
\end{proof}

\begin{lem}\label{le:single_congruent}
Let $(\rV,\pi)$ be a $\Pi$-congruent standard pair. If Assumption \ref{as:single_irreducible} holds, then $\BC(\pi)$, which exists by Proposition \ref{pr:bc_global}, is a relevant representation of $\GL_N(\dA_F)$ (Definition \ref{de:relevant}); and moreover, $\rho_{\BC(\pi),\iota_\ell}$ is residually irreducible.
\end{lem}

\begin{proof}
Let $\rho_{\BC(\pi),\iota_\ell}\colon\Gamma_F\to\GL_N(\ol\dQ_\ell)$ be the associated Galois representation (Remark \ref{re:bc_global}). Since $\pi$ is $\Pi$-congruent, by the Chebotarev density theorem, $\rho_{\BC(\pi),\iota_\ell}$ admits a lattice whose residual representation is isomorphic to $\bar\rho_{\Pi,\lambda}\otimes_{O_\lambda/\lambda}\ol\dF_\ell$, which is irreducible. If $\BC(\pi)$ is not cuspidal, then $\rho_{\BC(\pi),\iota_\ell}$ is decomposable, which is a contradiction. Thus, $\BC(\pi)$ is cuspidal. Together with \cite{Shi}*{Theorem~1.1(iii,iv)}, we obtain that $\BC(\pi)$ is relevant. The lemma follows.
\end{proof}

\begin{lem}\label{le:single_compact}
Assume Assumption \ref{as:single_irreducible}. Then the natural maps
\begin{align*}
\rH^i_{\et,c}(\Sh(\rV'_N,\tj_N\rK^{p\circ}_N\rK'_{p,N})_{\ol{F}},O_\lambda)_{\fm}&\to
\rH^i_\et(\Sh(\rV'_N,\tj_N\rK^{p\circ}_N\rK'_{p,N})_{\ol{F}},O_\lambda)_{\fm},\\
\rH^i_{\fT,c}(\ol\rM^\bullet_N,O_\lambda)_{\fm}
&\to\rH^i_\fT(\ol\rM^\bullet_N,O_\lambda)_{\fm},
\end{align*}
are both isomorphisms for every $i\in\dZ$.
\end{lem}

\begin{proof}
By Lemma \ref{le:ns_pbc}, and the description of the weight spectral sequence $(\rE^{p,q}_s,\rd^{p,q}_s)$ in Lemma \ref{le:ns_weight_odd} (for $N$ odd) and Lemma \ref{le:ns_weight_even} (for $N$ even), it suffices to show that the natural map
\begin{align}\label{eq:single_compact}
\rH^i_{\et,c}(\Sh(\rV'_N,\tj_N\rK^{p\circ}_N\rK'_{p,N})_{\ol{F}},O_\lambda)_{\fm}\to
\rH^i_\et(\Sh(\rV'_N,\tj_N\rK^{p\circ}_N\rK'_{p,N})_{\ol{F}},O_\lambda)_{\fm}
\end{align}
is an isomorphism for every $i\in\dZ$. This is trivial when $\Sh(\rV'_N,\tj_N\rK^{p\circ}_N\rK'_{p,N})$ is proper.

If $\Sh(\rV'_N,\tj_N\rK^{p\circ}_N\rK'_{p,N})$ is not proper, then the Witt index of $\rV'_N$ is $1$. In this case, the Shimura variety $\Sh(\rV'_N,\tj_N\rK^{p\circ}_N\rK'_{p,N})$ has a \emph{unique} toroidal compactification \cite{AMRT}, which we denote by $\wt\Sh(\rV'_N,\tj_N\rK^{p\circ}_N\rK'_{p,N})$, since the choice of the relevant combinatorial data is unique (see also \cite{Lar92} for more details in the case where $N=3$); it is smooth over $F$. As $\tj_N\rK^{p\circ}_N\rK'_{p,N}$ is neat, the boundary $Z\coloneqq\wt\Sh(\rV'_N,\tj_N\rK^{p\circ}_N\rK'_{p,N})\setminus\Sh(\rV'_N,\tj_N\rK^{p\circ}_N\rK'_{p,N})$ is geometrically isomorphic to a disjoint union of abelian varieties (of dimension $N-2$). In particular, $\rH^i_\et(Z_{\ol{F}},O_\lambda)$ is a free $O_\lambda$-module (of finite rank). Let $\pi'^\infty$ be an irreducible admissible representation of $\rU(\rV'_N)(\dA_{F^+}^\infty)$ that appears in $\rH^i_\et(Z_{\ol{F}},O_\lambda)\otimes_{O_\lambda,\iota_\ell^{-1}}\dC$. Then $\pi'^\infty$ extends to an automorphic representation $\pi'$ of $\rU(\rV'_N)(\dA_{F^+})$ that is a subquotient of the parabolic induction of a cuspidal automorphic representation of $L(\dA_{F^+})$ where $L$ is the unique proper Levi subgroup of $\rU(\rV'_N)$ up to conjugation. In particular, $\BC(\pi')$ exists and is not cuspidal. Thus, by (the same argument of) Lemma \ref{le:single_congruent}, we have $\rH^i_\et(Z_{\ol{F}},O_\lambda)_\fm=0$ for every $i\in\dZ$. This implies that \eqref{eq:single_compact} is an isomorphism.
\end{proof}

\subsection{Tate classes in the odd rank case}
\label{ss:tate}

In this section, we assume that $N=2r+1$ is odd with $r\geq 1$. We study the properties of the localized spectral sequence $\rE^{p,q}_{s,\fm}$, after Lemma \ref{le:ns_weight_odd}.

\begin{lem}\label{le:tate_vanishing}
We have
\[
\rH^i_\fT(\ol\rM^\dag_N,O_\lambda)_\fm=0
\]
for every odd integer $i$.
\end{lem}

\begin{proof}
For $i\neq 2r-1$, it follows from Lemma \ref{le:ns_cohomology}(1). Now we assume $i=2r-1$.

Suppose that $\pi^{\infty,p}$ is an irreducible admissible representation of $\rU(\rV^\circ_N)(\dA^{\infty,p}_{F^+})$ that appears in the cohomology $\rH^{2r-1}_\fT(\ol\rM^\dag_N,O_\lambda)_\fm\otimes_{O_\lambda,\iota_\ell^{-1}}\dC$. By Proposition \ref{pr:ns_link_cohomology}, we may complete $\pi^{\infty,p}$ to an automorphic representation $\pi$ of $\rU(\rV^\circ_N)(\dA_{F^+})$ as in that proposition, such that $(\rV^\circ_N,\pi)$ is a $\Pi$-congruent standard pair, and that $\BC(\pi_\fp)$ is a constituent of an unramified principal series of $\GL_N(F_\fp)$, whose Satake parameter contains $\{-p,-p^{-1}\}$ which is then different from $\balpha(\Pi_\fp)$ in $\ol\dF_\ell$ by (P5).

On the other hand, by the Chebotarev density theorem, both $\rho_{\BC(\pi),\iota_\ell}$ and $\rho_{\Pi,\lambda}\otimes_{E_\lambda}\ol\dQ_\ell$ each admits a lattice such that their reductions are isomorphic. In particular, the residual representations of $\rho_{\BC(\pi),\iota_\ell}$ and $\rho_{\Pi,\lambda}\otimes_{E_\lambda}\ol\dQ_\ell$ have the same Frobenius eigenvalues at the unique place of $F$ above $\fp$. However, this is not possible by Proposition \ref{pr:arthur}(2) and Proposition \ref{pr:galois}(2). Therefore, we must have $\rH^{2r-1}_\fT(\ol\rM^\dag_N,O_\lambda)_\fm=0$. The lemma is proved.
\end{proof}

\begin{lem}\label{le:tate_weight}
Assume Assumption \ref{as:single_vanishing}. We have
\begin{enumerate}
  \item $\rE^{p,q}_{1,\fm}=0$ if $q$ is odd;

  \item $\rE^{p,q}_{1,\fm}$ is a free $O_\lambda$-module for every $(p,q)\in\dZ^2$;

  \item $\rE^{p,q}_{2,\fm}=0$ unless $(p,q)=(0,2r)$;

  \item $\rE^{0,2r}_{2,\fm}$ is canonically isomorphic to $\rH^{2r}_\fT(\ol\rM_N,\rR\Psi O_\lambda(r))_\fm$, which is a free $O_\lambda$-module;

  \item $\rE^{0,2r}_{s,\fm}$ degenerates at the second page.
\end{enumerate}
\end{lem}

\begin{proof}
Part (1) follows from Lemma \ref{le:tate_vanishing} and Assumption \ref{as:single_vanishing}. Part (3) follows since $\rd^{-1,2r}_1$ is injective and $\rd^{0,2r}_1$ is surjective. The remaining parts are immediate consequences of (1) and Assumption \ref{as:single_vanishing}.
\end{proof}

\begin{theorem}\label{th:tate}
The map
\[
\nabla^1_\fm\colon\rE^{0,2r}_{2,\fm}\to O_\lambda[\Sh(\rV^\circ_N,\rK^\circ_N)]_\fm
\]
(Construction \ref{cs:ns_nabla}) is surjective. Moreover, if we assume Assumptions \ref{as:single_vanishing}, \ref{as:single_irreducible}, and Hypothesis \ref{hy:unitary_cohomology} for $N$, then we have
\begin{enumerate}
  \item The generalized Frobenius eigenvalues of the $(O_\lambda/\lambda)[\Gal(\ol\dF_p/\dF_{p^2})]$-module $\rE^{0,2r}_{2,\fm}\otimes_{O_\lambda}O_\lambda/\lambda$ is contained in the set of roots of $P_{\balpha(\Pi_\fp)}\modulo\lambda$ in a finite extension of $O_\lambda/\lambda$.

  \item The $O_\lambda[\Gal(\ol\dF_p/\dF_{p^2})]$-module $\rE^{0,2r}_{2,\fm}$ is weakly semisimple (Definition \ref{de:weakly_semisimple}).

  \item The map $\nabla^1_\fm$ induces an isomorphism
     \[
     \nabla^1_\fm\colon(\rE^{0,2r}_{2,\fm})_{\Gal(\ol\dF_p/\dF_{p^2})}\xrightarrow{\sim} O_\lambda[\Sh(\rV^\circ_N,\rK^\circ_N)]_\fm.
     \]
\end{enumerate}
\end{theorem}

By Remark \ref{re:ns_nabla}, the map $\nabla^1_\fm$ always factors through the quotient map $\rE^{0,2r}_{2,\fm}\to(\rE^{0,2r}_{2,\fm})_{\Gal(\ol\dF_p/\dF_{p^2})}$.

\begin{proof}
We first show that $\nabla^1_\fm$ is surjective. From Construction \ref{cs:ns_weight}, we have a map
\[
(\Inc^\circ_!,\Inc^\dag_!,\Inc^\bullet_!\circ\tT^{\bullet\circ}_\fp)\coloneqq O_\lambda[\Sh(\rV^\circ_N,\rK^\circ_N)]^{\oplus 3}\to\rE^{0,2r}_1
\]
which induces a map
\[
\Ker\(\rd^{0,2r}_1\circ(\Inc^\circ_!,\Inc^\dag_!,\Inc^\bullet_!\circ\tT^{\bullet\circ}_\fp)\)\to\Ker\rd^{0,2r}_1.
\]
However, by Lemma \ref{le:ns_weight_pre}, the former kernel is simply the kernel of the map
\[
\begin{pmatrix}
p+1 & -1 & 0
\end{pmatrix}
\begin{pmatrix}
\Inc_\circ^* \\
\Inc_\dag^* \\
\Inc_\bullet^*
\end{pmatrix}
\begin{pmatrix}
\Inc^\circ_! & \Inc^\dag_! & \Inc^\bullet_!\circ\tT^{\bullet\circ}_\fp
\end{pmatrix}.
\]
Now since $(p+1,-1,0)$ and $(0,\tT^{\circ\bullet}_\fp\circ\tT^{\bullet\circ}_\fp,(p+1)^2\tT^{\circ\bullet}_\fp)\otimes O_\lambda$ are linearly independent, by Nakayama's lemma, $\nabla^1_\fm$ is surjective if the following matrix
\[
\begin{pmatrix}
\Inc_\circ^* \\
\Inc_\dag^* \\
\tT^{\circ\bullet}_\fp\circ\Inc_\bullet^*
\end{pmatrix}
\begin{pmatrix}
\Inc^\circ_! & \Inc^\dag_! & \Inc^\bullet_!\circ\tT^{\bullet\circ}_\fp
\end{pmatrix}
\]
in $\dT^\circ_{N,\fp}$ is nondegenerate modulo $\fm$. However, by Lemma \ref{le:ns_weight_odd}(2), the above matrix equals
\[
\begin{pmatrix}
1 & 0 & 0 \\
0 & -(p+1)^2 & \tI^\circ_{N,\fp} \\
0 &
\tI^\circ_{N,\fp} & \tT^{\circ\bullet}_{N,\fp}\circ\tT^\bullet_{N,\fp}\circ\tT^{\bullet\circ}_{N,\fp}
\end{pmatrix},
\]
whose non-degeneracy modulo $\fm$ follows from Lemma \ref{le:enumeration_odd_0}, Proposition \ref{pr:enumeration_odd}, and (P4,P5).

Now we consider the three remaining assertions. By Lemma \ref{le:ns_pbc} and Lemma \ref{le:tate_weight}, we have an isomorphism
\[
\rE^{0,2r}_{2,\fm}\simeq\rH^{2r}_\et(\Sh(\rV',\tj_N\rK^{p\circ}_N\rK'_{p,N})_{\ol{F}},O_\lambda(r))_\fm
\]
of $O_\lambda[\Gal(\ol\dQ_p/\dQ_{p^2})]$-modules. By Lemma \ref{le:single_congruent}, Lemma \ref{le:single_compact}, Proposition \ref{pr:arthur}(2), and Hypothesis \ref{hy:unitary_cohomology}, we have
\[
\rH^{2r}_\et(\Sh(\rV',\tj_N\rK^{p\circ}_N\rK'_{p,N})_{\ol{F}},O_\lambda(r))_\fm\otimes_{O_\lambda}\ol\dQ_\ell
\simeq\bigoplus_{\pi'}\rho_{\BC(\pi'),\iota_\ell}^\tc(r)^{\oplus d(\pi')}
\]
of representations of $\Gamma_F$ with coefficients in $\ol\dQ_\ell$, where $d(\pi')\coloneqq\dim(\pi'^{\infty,p})^{\tj_N\rK^{p\circ}_N}$; and the direct sum is taken over all automorphic representations $\pi'$ of $\rU(\rV')(\dA_{F^+})$ satisfying:
\begin{itemize}[label={\ding{109}}]
  \item $(\rV',\pi')$ is a $\Pi$-congruent standard pair;

  \item $\pi'_{\ul\tau_\infty}$ is a holomorphic discrete series representation of $\rU(\rV')(F^+_{\ul\tau_\infty})$ with the Harish-Chandra parameter $\{-r,1-r,\dots,r-1,r\}$; and

  \item $\pi'_{\ul\tau}$ is trivial for every archimedean place $\ul\tau\neq\ul\tau_\infty$.
\end{itemize}

For the proof of (1--3), we may replace $E_\lambda$ by a finite extension inside $\ol\dQ_\ell$ such that $\rho_{\BC(\pi'),\iota_\ell}$ is defined over $E_\lambda$ for every $\pi'$ appearing in the previous direct sum. Now we regard $\rho_{\BC(\pi'),\iota_\ell}$ as a representation over $E_\lambda$. Then $\rho_{\BC(\pi'),\iota_\ell}(r)$ admits a $\Gamma_F$-stable $O_\lambda$-lattice $\rR_{\BC(\pi')}$, unique up to homothety, whose reduction $\bar\rR_{\BC(\pi')}$ is isomorphic to $\bar\rho_{\Pi,\lambda}(r)$. Moreover, we have an inclusion
\[
\rE^{0,2r}_{2,\fm}\simeq\rH^{2r}_\et(\Sh(\rV',\tj_N\rK^{p\circ}_N\rK'_{p,N})_{\ol{F}},O_\lambda(r))_\fm
\subseteq\bigoplus_{\pi'}(\rR_{\BC(\pi')}^\tc)^{\oplus d(\pi')}
\]
of $O_\lambda[\Gal(\ol\dF_p/\dF_{p^2})]$-modules. This already implies (1).

By (P4), we know that $\bar\rho_{\Pi,\lambda}^\tc(r)$ is weakly semisimple and
\[
\dim_{O_\lambda/\lambda}\bar\rho_{\Pi,\lambda}^\tc(r)^{\Gal(\ol\dF_p/\dF_{p^2})}=1.
\]
On the other hand, we have
\[
\dim_{E_\lambda}\rho_{\BC(\pi'),\iota_\ell}^\tc(r)^{\Gal(\ol\dF_p/\dF_{p^2})}\geq 1.
\]
Thus by Lemma \ref{le:weakly_semisimple2}, for every $\pi'$ in the previous direct sum, $\rR_{\BC(\pi')}^\tc$ is weakly semisimple, and
\[
\dim_{E_\lambda}\rho_{\BC(\pi'),\iota_\ell}^\tc(r)^{\Gal(\ol\dF_p/\dF_{p^2})}=1.
\]
This implies (2) by Lemma \ref{le:weakly_semisimple1}(1).

The above discussion also implies that, for (3), it suffices to show
\[
\sum_{\pi'}d(\pi')\leq\dim_{E_\lambda}O_\lambda[\Sh(\rV^\circ_N,\rK^\circ_N)]_\fm\otimes_{O_\lambda}E_\lambda
\]
where $\pi'$ is taken over the same set as in the previous direct sum. However, this follows from Corollary \ref{pr:jacquet_langlands_1} and Lemma \ref{le:single_intertwining}. The theorem is proved.
\end{proof}

\subsection{Arithmetic level-raising in the even rank case}
\label{ss:raising}

In this subsection, we assume that $N=2r$ is even with $r\geq 1$. We study the properties of the localized spectral sequence $\rE^{p,q}_{s,\fm}$, after Lemma \ref{le:ns_weight_even}.

\begin{proposition}\label{pr:raising_weight}
Assume Assumptions \ref{as:single_vanishing}, \ref{as:single_irreducible}, and Hypothesis \ref{hy:unitary_cohomology} for $N$. Then we have
\begin{enumerate}
  \item The maps
     \begin{align*}
     (\Inc^\circ_!+\Inc^\dag_!+\Inc^\bullet_!)_\fm &\colon
     O_\lambda[\Sh(\rV^\circ_N,\rK^\circ_N)]_\fm^{\oplus 2}\bigoplus O_\lambda[\Sh(\rV^\circ_N,\rK^\bullet_N)]_\fm
     \to\rE^{0,2r-2}_{1,\fm}(-1) \\
     (\Inc^\circ_!+\Inc^\bullet_!)_\fm &\colon
     O_\lambda[\Sh(\rV^\circ_N,\rK^\circ_N)]_\fm\bigoplus O_\lambda[\Sh(\rV^\circ_N,\rK^\bullet_N)]_\fm
     \to\rE^{0,2r-2}_{1,\fm}(-1)
     \end{align*}
     from Construction \ref{cs:ns_weight} are isomorphisms when $N\geq 4$ and $N=2$, respectively.

  \item The maps
     \begin{align*}
     (\Inc_\circ^*,\Inc_\dag^*,\Inc_\bullet^*)_\fm &\colon\rE^{0,2r}_{1,\fm}\to
     O_\lambda[\Sh(\rV^\circ_N,\rK^\circ_N)]_\fm^{\oplus 2}\bigoplus O_\lambda[\Sh(\rV^\circ_N,\rK^\bullet_N)]_\fm \\
     (\Inc_\circ^*,\Inc_\bullet^*)_\fm &\colon\rE^{0,2r}_{1,\fm}\to
     O_\lambda[\Sh(\rV^\circ_N,\rK^\circ_N)]_\fm\bigoplus O_\lambda[\Sh(\rV^\circ_N,\rK^\bullet_N)]_\fm
     \end{align*}
     from Construction \ref{cs:ns_weight} are surjective with kernel the $O_\lambda$-torsion of $\rH^{2r}_\fT(\ol\rM^\bullet_N,O_\lambda(r))_\fm$ when $N\geq 4$ and $N=2$, respectively.

  \item The map $\nabla^0_\fm\colon\Ker\rd^{0,2r}_{1,\fm}\to O_\lambda[\Sh(\rV^\circ_N,\rK^\circ_N)]_\fm$ (Construction \ref{cs:ns_nabla}) is surjective.

  \item The map $\nabla^0_\fm\circ\rd^{-1,2r}_{1,\fm}\circ\rd^{0,2r-2}_{1,\fm}(-1)$ induces a map
      \[
      \rF_{-1}\rH^1(\rI_{\dQ_{p^2}},\rH^{2r-1}_\fT(\ol\rM_N,\rR\Psi O_\lambda(r))_\fm)\to
      O_\lambda[\Sh(\rV^\circ_N,\rK^\circ_N)]_\fm/((p+1)\tR^\circ_{N,\fp}-\tI^\circ_{N,\fp})
      \]
      which is surjective, whose kernel is canonically the $O_\lambda$-torsion of $\rH^{2r}_\fT(\ol\rM^\bullet_N,O_\lambda(r))_\fm$.
\end{enumerate}
\end{proposition}

\begin{proof}
We only prove the proposition when $N\geq 4$, and leave the much easier case where $N=2$ to the readers.

We first claim that the map
\[
(\inc^\dag_!+\inc^\bullet_!\circ\tT^{\bullet\circ}_{N,\fp})_\fm\colon
O_\lambda[\Sh(\rV^\circ_N,\rK^\circ_N)]_\fm^{\oplus 2}\to\rH^{2r-2}_\fT(\ol\rM^\bullet_N,O_\lambda(r-1))_\fm
\]
is an isomorphism. In fact, by Lemma \ref{le:raising_bullet} below, it suffices to find a line bundle $\cL$ as in Definition \ref{de:ns_incidence_even} such that $(\inc_\cL)_\fm$ is surjective, where
\[
\inc_\cL\coloneqq(\inc_\dag^*,\tT^{\circ\bullet}_{N,\fp}\circ\inc_\bullet^*)
\circ\Theta_\cL\circ(\inc^\dag_!+\inc^\bullet_!\circ\tT^{\bullet\circ}_{N,\fp})
\]
in which $\Theta_\cL$ is defined in Definition \ref{de:ns_incidence_even}. We take $\cL$ to be $\cO(\rM^\dag_N)^{\otimes 2}\otimes(\Lie_{\cA,\tau_\infty^\tc})^{\otimes p+1}$. Then by Proposition \ref{pr:ns_incidence_even1} and Proposition \ref{pr:ns_incidence_even2}, the endomorphism $\inc_\cL$ is given the matrix
\[
\begin{pmatrix}
(p+1)^3 & -(p+1)\tI^\circ_{N,\fp} \\
-(p+1)\tI^\circ_{N,\fp} &
\tT^{\circ\bullet}_{N,\fp}\circ(\tR^\bullet_{N,\fp}+(\tR^\bullet_{N,\fp}+(p+1)\tT^\bullet_{N,\fp}))\circ\tT^{\bullet\circ}_{N,\fp}
\end{pmatrix}
\]
in $\dT^\circ_{N,\fp}$. Now, by Lemma \ref{le:enumeration_even_0} and Proposition \ref{pr:enumeration_even}, the determinant of the above matrix mod $\fm$ is equal to
\[
\resizebox{\hsize}{!}{
\xymatrix{
\displaystyle
-p^{r^2}\prod_{i=1}^r\(\alpha_i+\frac{1}{\alpha_i}+2\)\cdot
\((p+1)^2p^{r^2}\prod_{i=1}^r\(\alpha_i+\frac{1}{\alpha_i}-p-\frac{1}{p}\)
+(p+1)^3\(p^{r^2+1}-p^{r^2-1}\)\sum_{j=1}^r\prod_{\substack{i=1\\i\neq j}}^r\(\alpha_i+\frac{1}{\alpha_i}-p-\frac{1}{p}\)\)
}
}
\]
where $\{\alpha_r,\dots,\alpha_1,\alpha_1^{-1},\dots,\alpha_r^{-1}\}$ are the roots of $P_{\balpha(\Pi_{N,\fp})}\modulo\lambda$ in a finite extension of $O_\lambda/\lambda$. By (P2), we have
\[
p^{r^2}(p+1)^3\(p^{r^2+1}-p^{r^2-1}\)\not\equiv0 \mod\lambda;
\]
by (P4), we have
\[
\prod_{i=1}^r\(\alpha_i+\frac{1}{\alpha_i}-p-\frac{1}{p}\)\equiv0 \mod\lambda,\quad
\sum_{j=1}^r\prod_{\substack{i=1\\i\neq j}}^r\(\alpha_i+\frac{1}{\alpha_i}-p-\frac{1}{p}\)\not\equiv0 \mod\lambda;
\]
and by (P5), we have
\[
\prod_{i=1}^r\(\alpha_i+\frac{1}{\alpha_i}+2\)\not\equiv0 \mod\lambda.
\]
In particular, the matrix representing $\inc_\cL$ is nondegenerate modulo $\fm$, hence the claim follows from Nakayama's lemma.

Part (1) follows immediately from the above claim and Lemma \ref{le:single_intertwining}. Part (2) follows from (1) by the Poincar\'{e} duality theorem, together with Lemma \ref{le:single_compact}.

For (3), by definition, $\nabla^0_\fm$ is the restriction to $\Ker\rd^{0,2r}_{1,\fm}$ of the composition of
\[
(\tT^{\circ\bullet}_{N,\fp}\circ\tT^{\bullet\circ}_{N,\fp}\circ\Inc_\circ^*,\Inc_\dag^*,\tT^{\circ\bullet}_{N,\fp}\circ\Inc_\bullet^*)_\fm\colon
\rE^{0,2r}_{1,\fm}\to O_\lambda[\Sh(\rV^\circ_N,\rK^\circ_N)]_\fm^{\oplus 3}
\]
and the obviously surjective map
\[
(1,0,p+1)\colon O_\lambda[\Sh(\rV^\circ_N,\rK^\circ_N)]_\fm^{\oplus 3}\to O_\lambda[\Sh(\rV^\circ_N,\rK^\circ_N)]_\fm.
\]
By (2) and Lemma \ref{le:single_intertwining}, the map $(\tT^{\circ\bullet}_{N,\fp}\circ\tT^{\bullet\circ}_{N,\fp}\circ\Inc_\circ^*,\Inc_\dag^*,\tT^{\circ\bullet}_{N,\fp}\circ\Inc_\bullet^*)_\fm$ is surjective. On the other hand, the restriction of $\rd^{0,2r}_1$ to $\rH^{2r}_\fT(\ol\rM^\bullet_N,O_\lambda(r))$ coincides with $\inc_\dag^*$ (Construction \ref{cs:ns_incidence}), after composing with the isomorphism $\rH^{2r}_\fT(\ol\rM^\dag_N,O_\lambda(r))\xrightarrow{\sim}O_\lambda[\Sh(\rV^\circ_N,\rK^\circ_N)]$ as in the construction of $\inc_\dag^*$. Thus, by (2), the restriction of $\rd^{0,2r}_{1,\fm}$ to $\rH^{2r}_\fT(\ol\rM^\bullet_N,O_\lambda(r))_\fm$ is surjective, hence $\Delta^0_\fm$ is surjective.

Now we consider (4). Let $(\rE^{0,2r}_{1,\fm})_\free$ be the free $O_\lambda$-quotient of $\rE^{0,2r}_{1,\fm}$, which is simply the quotient by the $O_\lambda$-torsion $(\rH^{2r}_\fT(\ol\rM^\bullet_N,O_\lambda(r))_\fm)_\tor$ of $\rH^{2r}_\fT(\ol\rM^\bullet_N,O_\lambda(r))_\fm$. Thus by (2), we obtain an isomorphism
\[
(\Inc_\circ^*,\Inc_\dag^*,\Inc_\bullet^*)_\fm\colon(\rE^{0,2r}_{1,\fm})_\free\xrightarrow{\sim}
O_\lambda[\Sh(\rV^\circ_N,\rK^\circ_N)]_\fm^{\oplus 2}\bigoplus O_\lambda[\Sh(\rV^\circ_N,\rK^\bullet_N)]_\fm
\]
through which we identify the two sides. If we let $(\Ker\rd^{0,2r}_{1,\fm})_\free$ be the free $O_\lambda$-quotient of $\Ker\rd^{0,2r}_{1,\fm}$, then by Lemma \ref{le:ns_weight_pre}, the above isomorphism maps the submodule $(\Ker\rd^{0,2r}_{1,\fm})_\free$ to the kernel of the map
\[
(p+1,-1,0)\colon O_\lambda[\Sh(\rV^\circ_N,\rK^\circ_N)]_\fm^{\oplus 2}\bigoplus O_\lambda[\Sh(\rV^\circ_N,\rK^\bullet_N)]_\fm
\to O_\lambda[\Sh(\rV^\circ_N,\rK^\circ_N)]_\fm.
\]
By Assumption \ref{as:single_vanishing}, we have $\IM\rd^{-1,2r}_{1,\fm}=\Ker\rd^{0,2r}_{1,\fm}$. Combining Lemma \ref{le:ns_weight_even}(5), we see that the map $\rd^{-1,2r}_{1,\fm}$ induces a canonical isomorphism
\[
\rF_{-1}\rH^1(\rI_{\dQ_{p^2}},\rH^{2r-1}_\fT(\ol\rM_N,\rR\Psi O_\lambda(r))_\fm)
\simeq\frac{\IM\rd^{-1,2r}_{1,\fm}}{\IM(\rd^{-1,2r}_{1,\fm}\circ\rd^{0,2r-2}_{1,\fm}(-1))}
=\frac{\Ker\rd^{0,2r}_{1,\fm}}{\IM(\rd^{-1,2r}_{1,\fm}\circ\rd^{0,2r-2}_{1,\fm}(-1))}
\]
induced by $\rd^{-1,2r}_{1,\fm}$. Thus, we have a canonical surjective map
\[
\rF_{-1}\rH^1(\rI_{\dQ_{p^2}},\rH^{2r-1}_\fT(\ol\rM_N,\rR\Psi O_\lambda(r))_\fm)
\to\frac{(\Ker\rd^{0,2r}_{1,\fm})_\free}{\IM(\rd^{-1,2r}_{1,\fm}\circ\rd^{0,2r-2}_{1,\fm}(-1))}
\]
whose kernel is
\[
\frac{(\rH^{2r}_\fT(\ol\rM^\bullet_N,O_\lambda(r))_\fm)_\tor}
{(\rH^{2r}_\fT(\ol\rM^\bullet_N,O_\lambda(r))_\fm)_\tor\cap\IM(\rd^{-1,2r}_{1,\fm}\circ\rd^{0,2r-2}_{1,\fm}(-1))}.
\]
By Lemma \ref{le:single_intertwining} and Lemma \ref{le:ns_weight_even}(7), we see that $(\Ker\rd^{0,2r}_{1,\fm})_\free\cap\Ker\nabla^0_\fm$ is contained in the image $\rd^{-1,2r}_{1,\fm}\circ\rd^{0,2r-2}_{1,\fm}(-1)$, as modules of $(\rE^{0,2r}_{1,\fm})_\free$. Thus, by (3), the map $\nabla^0_\fm$ induces an isomorphism
\[
\frac{(\Ker\rd^{0,2r}_{1,\fm})_\free}{\IM(\rd^{-1,2r}_{1,\fm}\circ\rd^{0,2r-2}_{1,\fm}(-1))}\xrightarrow{\sim}
\frac{O_\lambda[\Sh(\rV^\circ_N,\rK^\circ_N)]_\fm}{\IM(\nabla^0_\fm\circ\rd^{-1,2r}_{1,\fm}\circ\rd^{0,2r-2}_{1,\fm}(-1))}.
\]
By Lemma \ref{le:ns_weight_even}(8), $\IM(\nabla^0_\fm\circ\rd^{-1,2r}_{1,\fm}\circ\rd^{0,2r-2}_{1,\fm}(-1))$ coincides with the submodule
\[
\(\tT^{\circ\bullet}_{N,\fp}\circ((p+1)\tR^\bullet_{N,\fp}-\tT^{\bullet\circ}_{N,\fp}\circ\tT^{\circ\bullet}_{N,\fp})
\circ\tT^{\bullet\circ}_{N,\fp}\).O_\lambda[\Sh(\rV^\circ_N,\rK^\circ_N)]_\fm.
\]
Note that, by Lemma \ref{le:enumeration_even_0}, we have
\[
\tT^{\circ\bullet}_{N,\fp}\circ((p+1)\tR^\bullet_{N,\fp}-\tT^{\bullet\circ}_{N,\fp}\circ\tT^{\circ\bullet}_{N,\fp})
\circ\tT^{\bullet\circ}_{N,\fp}=\tI^\circ_{N,\fp}\((p+1)\tR^\circ_{N,\fp}-\tI^\circ_{N,\fp}\).
\]
Thus, to conclude (4), it remains to show that
\begin{align}\label{eq:raising_weight}
(\rH^{2r}_\fT(\ol\rM^\bullet_N,O_\lambda(r))_\fm)_\tor\cap\IM(\rd^{-1,2r}_{1,\fm}\circ\rd^{0,2r-2}_{1,\fm}(-1))=0.
\end{align}
By Lemma \ref{le:ns_pbc}, Hypothesis \ref{hy:unitary_cohomology}, Lemma \ref{le:single_congruent}, Lemma \ref{le:single_compact}, and Proposition \ref{pr:arthur}(2), we know that the $\ol\dQ_\ell[\Gamma_F]$-module $\rH^{2r-1}_\fT(\ol\rM_N,\rR\Psi O_\lambda(r))_\fm\otimes_{O_\lambda}\ol\dQ_\lambda$ is isomorphic to a direct sum of $\rho_{\Pi',\iota_\ell}(r)$ for some relevant representations $\Pi'$ of $\GL_N(\dA_F)$. By Proposition \ref{pr:galois} and \cite{TY07}*{Lemma~1.4(3)}, we know that $\rho_{\Pi',\iota_\ell}(r)$ is pure of weight $-1$ at $\fp$ (Definition \ref{de:purity}). In particular, we have $\rH^1(\dQ_{p^2},\rho_{\Pi',\iota_\ell}(r))=0$ by \cite{Nek07}*{Proposition~4.2.2(1)}, hence that both sides of the inclusion
\[
\rF_{-1}\rH^1(\rI_{\dQ_{p^2}},\rH^{2r-1}_\fT(\ol\rM_N,\rR\Psi O_\lambda(r))_\fm)\subseteq
\rH^1_\sing(\dQ_{p^2},\rH^{2r-1}_\fT(\ol\rM_N,\rR\Psi O_\lambda(r))_\fm)
\]
are torsion $O_\lambda$-modules. Thus, the $O_\lambda$-rank of $\IM(\rd^{-1,2r}_{1,\fm}\circ\rd^{0,2r-2}_{1,\fm}(-1))$ is equal to the $O_\lambda$-rank of $\Ker\rd^{0,2r}_{1,\fm}$, which in turn is equal to the sum of $O_\lambda$-ranks of $O_\lambda[\Sh(\rV^\circ_N,\rK^\circ_N)]_\fm$ and $O_\lambda[\Sh(\rV^\circ_N,\rK^\bullet_N)]_\fm$. However, the source of the map $\rd^{-1,2r}_{1,\fm}\circ\rd^{0,2r-2}_{1,\fm}(-1)$, which is $\rE^{0,2r-2}_{1,\fm}/\IM\rd^{-1,2r-2}_{1,\fm}$, is also a free $O_\lambda$-module of the same rank. Therefore, we must have \eqref{eq:raising_weight}. Part (4) is proved.
\end{proof}

\begin{lem}\label{le:raising_bullet}
Suppose that $N\geq 4$. Assume Assumptions \ref{as:single_vanishing}, \ref{as:single_irreducible}, and Hypothesis \ref{hy:unitary_cohomology} for $N$. Then $\rH^{2r-2}_\fT(\ol\rM^\bullet_N,O_\lambda)_\fm$ is a free $O_\lambda$-module; and its rank over $O_\lambda$ is at most twice the rank of the (free) $O_\lambda$-module $O_\lambda[\Sh(\rV^\circ_N,\rK^\circ_N)]_\fm$.
\end{lem}

\begin{proof}
By Assumption \ref{as:single_vanishing}, Lemma \ref{le:ns_weight_even}(2), and Lemma \ref{le:ns_cohomology}(2), we have an injective map
\[
\rH^{2r-2}_\fT(\ol\rM^\bullet_N,O_\lambda)_\fm\hookrightarrow\rH^{2r-2}_\fT(\ol\rM^\dag_N,O_\lambda)_\fm
\]
induced by $\rd^{0,2r-2}_1$. For the target, we have an isomorphism
\[
\rH^{2r-2}_\fT(\ol\rM^\dag_N,O_\lambda)_\fm\simeq O_\lambda[\Sh(\rV^\circ_N,\rK^\circ_N)]_\fm\oplus
\rH^\prim(\ol\rM^\dag_N,O_\lambda)_\fm.
\]
In particular, $\rH^{2r-2}_\fT(\ol\rM^\dag_N,O_\lambda)_\fm$, hence $\rH^{2r-2}_\fT(\ol\rM^\bullet_N,O_\lambda)_\fm$ are free $O_\lambda$-modules.

Suppose that $\pi^{\infty,p}$ is an irreducible admissible representation of $\rU(\rV^\circ_N)(\dA^{\infty,p}_{F^+})$ that appears in $\rH^{2r-2}_\fT(\ol\rM^\bullet_N,O_\lambda)_\fm\otimes_{O_\lambda,\iota_\ell^{-1}}\dC$. Then, by Proposition \ref{pr:ns_link_cohomology}, one can complete $\pi^{\infty,p}$ to an automorphic representation $\pi=\pi^{\infty,p}\otimes\pi_\infty\otimes\prod_{\fq\mid p}\pi_\fq$ such that $\pi_\infty$ is trivial; $\pi_\fq$ is unramified for $\fq\neq\fp$; and $\pi_\fp$ is a constituent of an unramified principal series. Moreover, $(\rV^\circ_N,\pi)$ is a $\Pi$-congruent standard pair. By Assumption \ref{as:single_irreducible} and Lemma \ref{le:single_congruent}, we know that $\BC(\pi)$ is relevant.

To prove the lemma, it suffices to show that for such $\pi$ as above, we have
\begin{align}\label{eq:raising_bullet_0}
\dim_{\ol\dQ_\ell}\rH^{2r-2}_\fT(\ol\rM^\bullet_N,\ol\dQ_\ell)[\iota_\ell\pi^\infty]
\leq 2\dim_{\ol\dQ_\ell}\ol\dQ_\ell[\Sh(\rV^\circ_N,\rK^\circ_N)][\iota_\ell\pi^\infty].
\end{align}
Recall from Proposition \ref{pr:ns_link_cohomology} that we have an isomorphism
\begin{align}\label{eq:raising_bullet_1}
\iota_\ell^{-1}\rH^\prim(\ol\rM^\dag_N,\ol\dQ_\ell)\simeq
\Map_{\rK^\circ_{N,\fp}}\(\rU(\rV^\circ_N)(F^+)\backslash\rU(\rV^\circ_N)(\dA^\infty_{F^+})/\rK^{p\circ}_N\prod_{\fq\mid p,\fq\neq\fp}\rK^\circ_{N,\fq},\Omega_N\).
\end{align}

By Proposition \ref{pr:arthur}(2), we have $\BC(\pi_\fp)\simeq\BC(\pi)_\fp$. Let $\rho_{\BC(\pi),\iota_\ell}\colon\Gamma_F\to\GL_N(\ol\dQ_\ell)$ be the associated Galois representation. Since $(\rV^\circ_N,\pi)$ is $\Pi$-congruent, by the Chebotarev density theorem, $\rho_{\BC(\pi),\iota_\ell}$ admits a lattice whose residual representation is isomorphic to $\bar\rho_{\Pi,\lambda}\otimes_{O_\lambda/\lambda}\ol\dF_\ell$, which is irreducible by Assumption \ref{as:single_irreducible}. Thus, by Proposition \ref{pr:galois}(2), $\balpha(\BC(\pi_\fp))$ does not contain $\{-1,-1\}$ due to (P5) and contains $\{p,p^{-1}\}$ with multiplicity at most one by (P4). We now have three cases.

Case 1: $\pi_\fp$ is unramified. Then \eqref{eq:raising_bullet_0} follows by \eqref{eq:raising_bullet_1} and the fact that the multiplicity of $\Omega_N$ in $\pi_\fp\res_{\rK^\circ_{N,\fp}}$ is at most $1$ by Proposition \ref{pr:omega}(2).

Case 2: $\pi_\fp$ is not unramified and $\pi_\fp\not\in\cS$, where $\cS$ is introduced in Proposition \ref{pr:packet}. Then by Lemma \ref{le:base_change_0}(1), $\pi_\fp\res_{\rK^\circ_{N,\fp}}$ does not contain $\Omega_N$. Thus, both sides of \eqref{eq:raising_bullet_0} are zero by \eqref{eq:raising_bullet_1}.

Case 3: $\pi_\fp$ belongs to $\cS$. Then we have $\ol\dQ_\ell[\Sh(\rV^\circ_N,\rK^\circ_N)][\iota_\ell\pi^\infty]=0$, hence an inclusion
\begin{align}\label{eq:raising_bullet}
\iota_\ell^{-1}\rH^{2r-2}_\fT(\ol\rM^\bullet_N,\ol\dQ_\ell)[\pi^\infty]\hookrightarrow
\Map_{\rK^\circ_{N,\fp}}\(\rU(\rV^\circ_N)(F^+)\backslash\rU(\rV^\circ_N)(\dA^\infty_{F^+})/\rK^{p\circ}_N\prod_{\fq\mid p,\fq\neq\fp}\rK^\circ_{N,\fq},\Omega_N\)[\pi^\infty]
\end{align}
by \eqref{eq:raising_bullet_1}. Note that, by Proposition \ref{pr:omega}(2), the multiplicity of $\Omega_N$ in $\pi_\fp\res_{\rK^\circ_{N,\fp}}$ is one, hence we have
\[
\Map_{\rK^\circ_{N,\fp}}\(\rU(\rV^\circ_N)(F^+)\backslash\rU(\rV^\circ_N)(\dA^\infty_{F^+})/\rK^{p\circ}_N\prod_{\fq\mid p,\fq\neq\fp}\rK^\circ_{N,\fq},\Omega_N\)[\pi^\infty]
\simeq (\pi^{\infty,p})^{\rK^{p\circ}_N}
\]
by Proposition \ref{pr:arthur}(2).

On the other hand, by Lemma \ref{le:single_compact}, Proposition \ref{pr:arthur}(2), Corollary \ref{pr:jacquet_langlands_0}, and Hypothesis \ref{hy:unitary_cohomology}, we know that the $\ol\dQ_\ell[\Gamma_F]$-module
\[
\rH^{2r-1}_\et(\Sh(\rV'_N,\tj_N\rK^{p\circ}_N\rK'_{p,N})_{\ol{F}},\ol\dQ_\ell)[\iota_\ell\pi^{\infty,p}]
\]
is isomorphic to $\dim(\pi^{\infty,p})^{\rK^{p\circ}_N}$ copies of $\rho_{\BC(\pi),\iota_\ell}^\tc$. By Proposition \ref{pr:galois}(2), $\rho_{\BC(\pi),\iota_\ell}^\tc\res_{\Gal(\ol\dQ_p/\dQ_{p^2})}$ has nontrivial monodromy action. Thus, by Lemma \ref{le:ns_pbc} and the spectral sequence $\rE^{p,q}_s$, the cokernel of \eqref{eq:raising_bullet} has dimension $\dim(\pi^{\infty,p})^{\rK^{p\circ}_N}$, which forces the source of \eqref{eq:raising_bullet} to vanish. In particular, \eqref{eq:raising_bullet_0} holds.

The lemma is proved.
\end{proof}

\begin{remark}
Proposition \ref{pr:raising_weight} has an amazing corollary which we now explain.

Suppose that $\ell\nmid p\prod_{i=1}^N(p^i-(-1)^i)$. Then the Tate--Thompson representation of $\Omega_N$ from \S\ref{ss:omega} of $\rK^\circ_{N,\fp}$ has a model $\Omega_{N,\ol\dF_\ell}$ over $\ol\dF_\ell$, which is again an irreducible summand of $\Ind_{\rK^\circ_{N,\fp}\cap\rK^\bullet_{N,\fp}}^{\rK^\circ_{N,\fp}}\ol\dF_\ell$. Thus, we obtain a natural map
\[
\bi\colon\ol\dF_\ell[\Sh(\rV^\circ_N,\rK^\bullet_N)]\to
\Map_{\rK^\circ_{N,\fp}}\(\rU(\rV^\circ_N)(F^+)\backslash\rU(\rV^\circ_N)(\dA^\infty_{F^+})/\rK_N^{p\circ}\prod_{\fq\mid p,\fq\neq\fp}\rK^\circ_{N,\fq},\Omega_{N,\ol\dF_\ell}\)
\]
of $\ol\dF_\ell[\dT^{\Sigma^+\cup\Sigma^+_p}_N]$-modules. Let the setup be as in \S\ref{ss:preliminaries} but replacing (P4) with a weaker condition that $\balpha(\Pi_\fp)\mod\lambda$ contains the pair $\{p,p^{-1}\}$ \emph{at most} once. Assume Assumptions \ref{as:single_vanishing}, \ref{as:single_irreducible}, and Hypothesis \ref{hy:unitary_cohomology} for $N$. Then $\bi_\fm$ is \emph{injective}.

Note that this result can be regarded as an Ihara type lemma for the definite unitary Shimura sets. Now we explain how to deduce it. For simplicity, we only consider the case where $N\geq 4$, and leave the much easier case where $N=2$ to the readers. First, we point out that since $\ell\nmid p\prod_{i=1}^N(p^i-(-1)^i)$, \eqref{eq:tate} holds with $\ol\dQ_\ell$ replaced by $\ol\dF_\ell$, under which the map $\bi$ coincides with the composite map
\[
\ol\dF_\ell[\Sh(\rV^\circ_N,\rK^\bullet_N)]\xrightarrow{\inc^\bullet_!}\rH^{2r-2}_\fT(\ol\rM^\bullet_N,\ol\dF_\ell(r-1))
\xrightarrow{(\rm^{\dag\bullet})^*}\rH^{2r-2}_\fT(\ol\rM^\dag_N,\ol\dF_\ell(r-1))\to\rH^\prim_\fT(\ol\rM^\dag_N,\ol\dF_\ell).
\]
As pointed out in the proof of Lemma \ref{le:raising_bullet}, the map
$(\rm^{\dag\bullet})^*_\fm\colon\rH^{2r-2}_\fT(\ol\rM^\bullet_N,O_\lambda)_\fm\to\rH^{2r-2}_\fT(\ol\rM^\dag_N,O_\lambda)_\fm$ is injective. Thus, it suffices to show that the map
\begin{align*}
(\inc^\dag_!+\inc^\bullet_!)_\fm\colon
\ol\dF_\ell[\Sh(\rV^\circ_N,\rK^\circ_N)]_\fm\bigoplus\ol\dF_\ell[\Sh(\rV^\circ_N,\rK^\bullet_N)]_\fm
\to\rH^{2r-2}_\fT(\ol\rM^\bullet_N,\ol\dF_\ell(r-1))_\fm
\end{align*}
is injective. When $\balpha(\Pi_\fp)\mod\lambda$ contains the pair $\{p,p^{-1}\}$ (exactly once), this follows from Proposition \ref{pr:raising_weight} (1). When $\balpha(\Pi_\fp)\mod\lambda$ does not contain the pair $\{p,p^{-1}\}$, it suffices to show that $(\inc_\cL)_\fm$ (Definition \ref{de:ns_incidence_even}) is injective with $\cL=\cO(\rM^\dag_N)$ and the coefficients $\ol\dF_\ell$. It is straightforward to see that such injectivity follows from Proposition \ref{pr:ns_incidence_even1}, Lemma \ref{le:single_intertwining}, Proposition \ref{pr:enumeration_even}(2), and Lemma \ref{le:enumeration_even_0}.
\end{remark}

\if false

\begin{remark}\label{re:raising_weight}
Following the well-known computation of level-raising of Shimura curves (see, for example, Step 5 of the proof of \cite{Liu2}*{Proposition~3.32}), we know that Proposition \ref{pr:raising_weight}(4) also holds when $N=2$. Moreover, as $\ol\rM^\bullet_2$ is a disjoint union of projective lines, the kernel of the map is trivial, hence the map is an isomorphism.
\end{remark}

\fi

Before stating the main theorem on the arithmetic level raising, we recall the following definition from \cite{LTXZZ}*{\S3.6}.

\begin{definition}\label{de:rigid}
Let $\bar{r}\colon\Gamma_{F^+}\to\sG_N(O_\lambda/\lambda)$ be a continuous homomorphism subject to the relation $\bar{r}^{-1}(\GL_N(O_\lambda/\lambda)\times(O_\lambda/\lambda)^\times)=\Gamma_F$ and $\nu\circ\bar{r}=\eta_{F/F^+}^N\epsilon_\ell^{1-N}$. We say that $\bar{r}$ is \emph{rigid for $(\Sigma^+_\mnm,\Sigma^+_\lr)$} if the following are satisfied:
\begin{enumerate}
  \item For $v$ in $\Sigma^+_\mnm$, every lifting of $\bar{r}_v$ is minimally ramified \cite{LTXZZ}*{Definition~3.4.8}.

  \item For $v$ in $\Sigma^+_\lr$, the generalized eigenvalues of $\bar{r}^\natural_v(\phi_w)$ in $\ol\dF_\ell$ contain the pair $\{\|v\|^{-N},\|v\|^{-N+2}\}$ exactly once, where $w$ is the unique place of $F$ above $v$.

  \item For $v$ in $\Sigma^+_\ell$, $\bar{r}_v^\natural$ is regular Fontaine--Laffaille crystalline \cite{LTXZZ}*{Definition~3.2.4}.

  \item For a nonarchimedean place $v$ of $F^+$ not in $\Sigma^+_\mnm\cup\Sigma^+_\lr\cup\Sigma^+_\ell$, the homomorphism $\ol{r}_v$ is unramified.
\end{enumerate}
Here, all liftings are with respect to the similitude character $\eta_{F/F^+}^N\epsilon_\ell^{1-N}$.
\end{definition}

Recall that we have fixed a positive integer $m$ at the beginning of \S\ref{ss:preliminaries}.

\begin{theorem}\label{th:raising}
Assume Assumptions \ref{as:single_vanishing}, \ref{as:single_irreducible}, and Hypothesis \ref{hy:unitary_cohomology} for $N$. We further assume that
\begin{enumerate}[label=(\alph*)]
  \item $\ell\geq 2(N+1)$ and $\ell$ is unramified in $F$;

  \item $\bar\rho_{\Pi,\lambda,+}$ (Remark \ref{re:single_irreducible}) is rigid for $(\Sigma^+_\mnm,\Sigma^+_\lr)$ (Definition \ref{de:rigid}), and $\bar\rho_{\Pi,\lambda}\res_{\Gal(\ol{F}/F(\zeta_\ell))}$ is absolutely irreducible;

  \item the composite homomorphism $\dT^{\Sigma^+}_N\xrightarrow{\phi_\Pi}O_E\to O_E/\lambda$ is cohomologically generic (Definition \ref{de:generic}); and

  \item $O_\lambda[\Sh(\rV^\circ_N,\rK^\circ_N)]_\fm$ is nontrivial.
\end{enumerate}
Then we have
\begin{enumerate}
  \item $\rH^i_\fT(\ol\rM^\bullet_N,O_\lambda)_\fm$ is a free $O_\lambda$-module for every $i\in\dZ$.

  \item $\rE^{p,q}_{2,\fm}$ is a free $O_\lambda$-module, and vanishes if $(p,q)\not\in\{(-1,2r),(0,2r-1),(1,2r-2)\}$.

  \item If we denote by $\{\alpha_1^{\pm1},\dots,\alpha_r^{\pm1}\}$ the roots of $P_{\balpha(\Pi_\fp)}\modulo\lambda$ in a finite extension of $O_\lambda/\lambda$, then the generalized Frobenius eigenvalues of the $(O_\lambda/\lambda)[\Gal(\ol\dF_p/\dF_{p^2})]$-module $\rH^{2r-1}_\fT(\ol\rM^\bullet_N,O_\lambda(r))_\fm\otimes_{O_\lambda}O_\lambda/\lambda$ is contained in $\{p\alpha_1^{\pm1},\dots,p\alpha_r^{\pm1}\}\setminus\{1,p^2\}$.

  \item The map in Proposition \ref{pr:raising_weight}(4) factors through a map
     \[
     \nabla^0_{/\fn}\colon\rF_{-1}\rH^1(\rI_{\dQ_{p^2}},\rH^{2r-1}_\fT(\ol\rM_N,\rR\Psi O_\lambda(r))/\fn)
     \to O_\lambda[\Sh(\rV^\circ_N,\rK^\circ_N)]/\fn
     \]
     which is an isomorphism, where $\fn$ is the ideal in Notation \ref{no:single_ideal}. The map from Lemma \ref{le:ns_weight_even}(6) induces a canonical isomorphism
      \[
      \rF_{-1}\rH^1(\rI_{\dQ_{p^2}},\rH^{2r-1}_\fT(\ol\rM_N,\rR\Psi O_\lambda(r))/\fn)
      \xrightarrow{\sim}\rH^1_\sing(\dQ_{p^2},\rH^{2r-1}_\fT(\ol\rM_N,\rR\Psi O_\lambda(r))/\fn).
      \]

  \item There exists a positive integer $\mu$ such that
     \[
     \rH^{2r-1}_\et(\Sh(\rV'_N,\tj_N\rK^{p\circ}_N\rK'_{p,N})_{\ol{F}},O_\lambda(r))/\fn
     \simeq\(\bar\rR^{(m)\tc}\)^{\oplus\mu}
     \]
     of $O_\lambda[\Gamma_F]$-modules, where $\rR$ is the $\Gamma_F$-stable $O_\lambda$-lattice in $\rho_{\Pi,\lambda}(r)$, unique up to homothety.
\end{enumerate}
\end{theorem}

\subsection{Proof of Theorem \ref{th:raising}}
\label{ss:raising2}

We apply the discussion of \cite{LTXZZ}*{\S3} to the pair $(\bar{r},\chi)$, where
\[
\bar{r}\coloneqq\bar\rho_{\Pi,\lambda,+}\colon\Gamma_{F^+}\to\sG_N(O_\lambda/\lambda)
\]
and $\chi\coloneqq\epsilon_\ell^{1-N}$ for the similitude character. Then $\bar{r}$ is rigid for $(\Sigma^+_\mnm,\Sigma^+_\lr)$, and also for $(\Sigma^+_\mnm,\Sigma^+_\lr\cup\{\fp\})$ by (P4).

For $?=\mix,\unr,\ram$, consider a global deformation problem \cite{LTXZZ}*{Definition~3.1.6}
\[
\sS^?\coloneqq(\bar{r},\eta_{F/F^+}^\mu\epsilon_\ell^{1-N},\Sigma^+_\mnm\cup\Sigma^+_\lr\cup\{\fp\}\cup\Sigma^+_\ell,
\{\sD_v\}_{v\in\Sigma^+_\mnm\cup\Sigma^+_\lr\cup\{\fp\}\cup\Sigma^+_\ell})
\]
where
\begin{itemize}[label={\ding{109}}]
  \item for $v\in\Sigma^+_\mnm$, $\sD_v$ is the local deformation problem classifying all liftings of $\bar{r}_v$;

  \item for $v\in\Sigma^+_\lr$, $\sD_v$ is the local deformation problem $\sD^\ram$ of $\bar{r}_v$ from \cite{LTXZZ}*{Definition~3.5.1};

  \item for $v=\fp$, $\sD_v$ is the local deformation problem $\sD^?$ of $\bar{r}_v$ from \cite{LTXZZ}*{Definition~3.5.1};

  \item for $v\in\Sigma^+_\ell$, $\sD_v$ is the local deformation problem $\sD^\FL$ of $\bar{r}_v$ from \cite{LTXZZ}*{Definition~3.2.5}.
\end{itemize}
Then we have the \emph{global universal deformation ring} $\sfR^\univ_{\sS^?}$ from \cite{LTXZZ}*{Proposition~3.1.7}. Put $\sfR^?\coloneqq\sfR^\univ_{\sS^?}$ for short. Then we have canonical surjective homomorphisms $\sfR^\mix\to\sfR^\unr$ and $\sfR^\mix\to\sfR^\ram$ of $O_\lambda$-rings. Finally, put
\[
\sfR^{\r{cong}}\coloneqq\sfR^\unr\otimes_{\sfR^\mix}\sfR^\ram.
\]
We fix a universal lifting
\[
r_\mix\colon\Gamma_{F^+}\to\sG_N(\sfR^\mix)
\]
of $\bar{r}$, which induce a continuous homomorphism
\[
r_\mix^\natural\colon\Gamma_F\to\GL_N(\sfR^\mix)
\]
by restriction (Notation \ref{no:sg_extension}). By pushforward, $\sfR^{\r{cong}}$ also induces homomorphisms
\[
r_\unr\colon\Gamma_{F^+}\to\sG_N(\sfR^\unr),\quad
r_\ram\colon\Gamma_{F^+}\to\sG_N(\sfR^\ram).
\]
Denote by $\rP_{F^+_\fp}$ the maximal closed subgroup of the inertia subgroup $\rI_{F^+_\fp}\subseteq\Gamma_{F^+_\fp}$ of pro-order coprime to $\ell$. Then $\Gamma_{F^+_\fp}/\rP_{F^+_\fp}\simeq t^{\dZ_\ell}\rtimes\phi_p^{\widehat\dZ}$ is a $p$-tame group \cite{LTXZZ}*{Definition~3.3.1}. By definition, the homomorphism $r_\mix^\natural$ is trivial on $\rP_{F^+_\fp}$. Let $\bar\sfv$ and $\bar\sfv'$ be eigenvectors in $(O_\lambda/\lambda)^{\oplus N}$ for $\bar{r}^\natural(\phi_p^2)$ with eigenvalues $p^{-2r}$ and $p^{-2r+2}$, respectively. By Hensel's lemma, $\bar\sfv$ and $\bar\sfv'$ lift to eigenvectors $\sfv$ and $\sfv'$ in $(\sfR^\mix)^{\oplus N}$ for $r_\mix^\natural(\phi_p^2)$, with eigenvalues $\sfs$ and $\sfs'$ in $\sfR^\mix$ lifting $p^{-2r}$ and $p^{-2r+2}$, respectively. Let $\sfx\in\sfR^\mix$ be the unique element such that $r_\mix^\natural(t)\sfv'=\sfx\sfv+\sfv'$. Then we must have $\sfx(\sfs-p^{-2r})=0$. By \cite{LTXZZ}*{Definition~3.5.1}, we have
\[
\sfR^\unr=\sfR^\mix/(\sfx),\quad
\sfR^\ram=\sfR^\mix/(\sfs-p^{-2r}),\quad
\sfR^{\r{cong}}=\sfR^\mix/(\sfs-p^{-2r},\sfx).
\]

Let $\sfT^\unr$ be the image of $\dT^{\Sigma^+}_N$ in $\End_{O_\lambda}(O_\lambda[\Sh(\rV^\circ_N,\rK^\circ_N)])$. By (d) in Theorem \ref{th:raising}, we know that $\sfT^\unr_\fm\neq 0$. Thus by \cite{LTXZZ}*{Theorem~3.6.3}, we have a canonical isomorphism $\sfR^\unr\xrightarrow{\sim}\sfT^\unr_\fm$ such that $O_\lambda[\Sh(\rV^\circ_N,\rK^\circ_N)]_\fm$ is canonically a free $\sfR^\unr$-module of rank $d_\unr>0$.\footnote{Here, we also need the easy fact that $\sfT^\unr_\fm$ and $O_\lambda[\Sh(\rV^\circ_N,\rK^\circ_N)]_\fm$ do not change if we replace $\fm$ by $\fm\cap\dT^{\Sigma^+\cup\Sigma^+_p\cup\Sigma^+_\ell}_N$} We may write the characteristic polynomial of $r_\unr^\natural(\phi_p^2)$ as $(T-\sfs)(T-p^{-4r+2}\sfs^{-1})Q(T)$, with $Q(T)\in\sfR^\unr[T]$ whose reduction in $(O_\lambda/\lambda)[T]$ does not have $p^{-2r}$ or $p^{-2r+2}$ as roots. By Proposition \ref{pr:enumeration_even}(2), we have
\[
((p+1)\tR^\circ_{N,\fp}-\tI^\circ_{N,\fp}).O_\lambda[\Sh(\rV^\circ_N,\rK^\circ_N)]_\fm
=(\sfs-p^{-2r}).O_\lambda[\Sh(\rV^\circ_N,\rK^\circ_N)]_\fm.
\]
In particular, we have
\[
O_\lambda[\Sh(\rV^\circ_N,\rK^\circ_N)]_\fm/((p+1)\tR^\circ_{N,\fp}-\tI^\circ_{N,\fp})
=O_\lambda[\Sh(\rV^\circ_N,\rK^\circ_N)]_\fm\otimes_{\sfR^\unr}\sfR^{\r{cong}},
\]
which is a free $\sfR^{\r{cong}}$-module of rank $d_\unr$.

On the other hand, let $\sfT^\ram$ be the image of $\dT^{\Sigma^+\cap\Sigma^+_p}_N$ in $\End_{O_\lambda}(\rH^{2r-1}_\fT(\ol\rM_N,\rR\Psi O_\lambda))$. By Proposition \ref{pr:raising_weight}(4) and Lemma \ref{le:ns_weight_even}(6), we know that $\sfT^\ram_\fm\neq0$. Thus by Lemma \ref{le:ns_pbc} and \cite{LTXZZ}*{Theorem~3.6.3} (with $(\Sigma^+_\mnm,\Sigma^+_\lr)$ replaced by $(\Sigma^+_\mnm,\Sigma^+_\lr\cup\{\fp\})$), we have a canonical isomorphism $\sfR^\ram\xrightarrow{\sim}\sfT^\ram_\fm$ such that $\rH^{2r-1}_\fT(\ol\rM_N,\rR\Psi O_\lambda)_\fm$ is canonically a free $\sfR^\ram$-module.\footnote{Here, we also need the fact that $\sfT^\ram_\fm$ and $\rH^{2r-1}_\fT(\ol\rM_N,\rR\Psi O_\lambda)_\fm$ do not change if we replace $\fm$ by $\fm\cap\dT^{\Sigma^+\cup\Sigma^+_p\cup\Sigma^+_\ell}_N$, which is a consequence of Theorem \ref{th:raising}(c).} Define the $\sfR^\ram$-module
\[
\sfH\coloneqq\Hom_{\Gamma_F}\((\sfR^\ram)^{\oplus N},\rH^{2r-1}_\fT(\ol\rM_N,\rR\Psi O_\lambda)_\fm\)
\]
where $\Gamma_F$ acts on $(\sfR^\ram)^{\oplus N}$ via the homomorphism $r_\ram^{\natural,\tc}$. By the same argument for \cite{Sch18}*{Theorem~5.6} (using Proposition \ref{pr:arthur} and Hypothesis \ref{hy:unitary_cohomology} here), we have a canonical isomorphism
\[
\rH^{2r-1}_\fT(\ol\rM_N,\rR\Psi O_\lambda)_\fm\simeq\sfH\otimes_{\sfR^\ram}(\sfR^\ram)^{\oplus N}
\]
of $\sfR^\ram[\Gamma_F]$-modules. Since $\sfR^\ram$ is a local ring, $\sfH$ is a free $\sfR^\ram$-module, say of rank $d_\ram$. If we still denote by $\sfv$ and $\sfv'$ for their projection in $(\sfR^\ram)^{\oplus N}$, then it is easy to see that
\[
\rH^1_\sing(\dQ_{p^2},(\sfR^\ram)^{\oplus N}(r))=\sfR^\ram\sfv/\sfx\sfv\simeq\sfR^\ram/(\sfx)=\sfR^{\r{cong}}.
\]
Thus, we obtain
\[
\rH^1_\sing(\dQ_{p^2},\rH^{2r-1}_\fT(\ol\rM_N,\rR\Psi O_\lambda(r))_\fm)\simeq \sfH\otimes_{\sfR^\ram}\rH^1_\sing(\dQ_{p^2},(\sfR^\ram)^{\oplus N}(r))
\simeq\sfH\otimes_{\sfR^\ram}\sfR^{\r{cong}},
\]
which is a free $\sfR^{\r{cong}}$-module of rank $d_\ram>0$.

\begin{proposition}\label{pr:old_same}
Under the assumptions of Theorem \ref{th:raising}, we have $d_\unr=d_\ram$. In particular, the two canonical maps
\begin{align*}
&\rF_{-1}\rH^1(\rI_{\dQ_{p^2}},\rH^{2r-1}_\fT(\ol\rM_N,\rR\Psi O_\lambda(r))_\fm)\to
O_\lambda[\Sh(\rV^\circ_N,\rK^\circ_N)]_\fm/((p+1)\tR^\circ_{N,\fp}-\tI^\circ_{N,\fp}), \\
&\rF_{-1}\rH^1(\rI_{\dQ_{p^2}},\rH^{2r-1}_\fT(\ol\rM_N,\rR\Psi O_\lambda(r))_\fm)\to
\rH^1_\sing(\dQ_{p^2},\rH^{2r-1}_\fT(\ol\rM_N,\rR\Psi O_\lambda(r))_\fm),
\end{align*}
from Proposition \ref{pr:raising_weight}(4) and Lemma \ref{le:ns_weight_even}(6), respectively, are both isomorphisms.
\end{proposition}

\begin{proof}
By Proposition \ref{pr:raising_weight}(4), the first map is surjective. By Lemma \ref{le:ns_weight_even}(6), the second map is injective. Thus, we must have $d_\ram\geq d_\unr>0$ by the previous discussion.

Take a geometric point $\eta_1\in(\Spec\sfR^\unr)(\ol\dQ_\ell)$ in the support of $O_\lambda[\Sh(\rV^\circ_N,\rK^\circ_N)]_\fm$, which corresponds to a relevant representation $\Pi_1$ of $\GL_N(\dA_F)$ by Lemma \ref{le:single_congruent}, such that $\rho_{\Pi_1,\iota_\ell}$ is residually isomorphic to $\bar\rho_{\Pi,\lambda}\otimes_{O_\lambda/\lambda}\ol\dF_\ell$. Then we have
\[
d_\unr=\dim\ol\dQ_\ell[\Sh(\rV^\circ_N,\rK^\circ_N)][\iota_\ell\phi_{\Pi_1}].
\]

Take a geometric point $\eta_2\in(\Spec\sfR^\ram)(\ol\dQ_\ell)$ in the support of $\rH^{2r-1}_\fT(\ol\rM_N,\rR\Psi O_\lambda)_\fm$, which corresponds to a relevant representation $\Pi_2$ of $\GL_N(\dA_F)$ by Lemma \ref{le:single_congruent}, such that $\rho_{\Pi_2,\iota_\ell}$ is residually isomorphic to $\bar\rho_{\Pi,\lambda}\otimes_{O_\lambda/\lambda}\ol\dF_\ell$. Then we have
\[
Nd_\ram=\dim\rH^{2r-1}_\fT(\ol\rM_N,\rR\Psi\ol\dQ_\ell)[\iota_\ell\phi_{\Pi_2}]
=\dim\rH^{2r-1}_\et(\Sh(\rV'_N,\tj_N\rK^{p\circ}_N\rK'_{p,N})_{\ol{F}},\ol\dQ_\ell)[\iota_\ell\phi_{\Pi_2}]
\]
by Lemma \ref{le:ns_pbc}. By Proposition \ref{pr:old} and Lemma \ref{le:old} below, we have $d_\unr=d_\ram$. The proposition follows.
\end{proof}

\begin{lem}\label{le:old}
Let $\Pi_1$ and $\Pi_2$ be two relevant representations of $\GL_N(\dA_F)$ such that the associated Galois representations $\rho_{\Pi_1,\iota_\ell}$ and $\rho_{\Pi_2,\iota_\ell}$ are both residually isomorphic to $\bar\rho_{\Pi,\lambda}\otimes_{O_\lambda/\lambda}\ol\dF_\ell$. For every $v\in\Sigma^+_\mnm$ (so that every lifting of $\bar\rho_{\Pi,\lambda,+,v}$ is minimally ramified), if we realize $\Pi_{1,v}$ and $\Pi_{2,v}$ on vector spaces $V_1$ and $V_2$, respectively, then there exist normalized intertwining operators $A_{\Pi_{1,v}}$ and $A_{\Pi_{2,v}}$ for $\Pi_{1,v}$ and $\Pi_{2,v}$ \cite{Shi11}*{\S4.1}, respectively, such that we have an $\GL_N(O_{F_v})$-equivariant isomorphism $i\colon V_1\xrightarrow{\sim}V_2$ satisfying $i\circ A_{\Pi_{1,v}}=A_{\Pi_{2,v}}\circ i$.
\end{lem}

\begin{proof}
We will give the proof when $v$ is nonsplit in $F$, and leave the other similar case to the readers. Let $w$ be the unique place of $F$ above $v$.

By Proposition \ref{pr:galois}(1), both $\Pi_{1,w}$ and $\Pi_{2,w}$ are tempered. Thus by the Bernstein--Zelevinsky classification, for $\alpha=1,2$, we can write
\[
\Pi_{\alpha,w}=\rI^{\GL_N(F_w)}_{P_\alpha}\(\sigma_{\alpha,-t_\alpha}\boxtimes\cdots\boxtimes\sigma_{\alpha,-1}
\boxtimes\sigma_{\alpha,0}\boxtimes\sigma_{\alpha,1}\boxtimes\cdots\boxtimes\sigma_{\alpha,t_\alpha}\)
\]
for some integer $t_\alpha\geq 0$, some standard parabolic subgroup $P_\alpha\subseteq\GL_N(F_w)$, and some (unitary) discrete series representations $\{\sigma_{\alpha,-t_\alpha},\dots,\sigma_{\alpha,t_\alpha}\}$ satisfying $\sigma_{\alpha,-i}\simeq\sigma_{\alpha,i}^{\vee\tc}$. See \S\ref{ss:local_bc} for the notation on parabolic induction.

By \cite{LTXZZ}*{Proposition~3.4.12(3)} and \cite{BLGGT}*{Lemma~1.3.4(2)}, we know that $\rho_{\Pi_1,\iota_\ell}\res_{\rI_{F_w}}$ and $\rho_{\Pi_2,\iota_\ell}\res_{\rI_{F_w}}$ are conjugate. Thus, by \cite{Yao}*{Lemma~3.6}, we have $P_1=P_2$ (say $P$) and $t_1=t_2$ (say $t$), and we assume that there are unramified (unitary) characters $\{\chi_{-t},\dots,\chi_t\}$ of $F_w^\times$ satisfying $\chi_{-i}\simeq\chi_i^{-1}$ such that $\sigma_{2,i}=\sigma_{1,i}\otimes\chi_i$. For every $i$, we choose a vector space $W_i$ on which $\sigma_{1,i}$ realizes (and also realize $\sigma_{1,i}^{\vee\tc}$ on $W_i$ via $g\mapsto\tp{g}^{-1,\tc}$), and fix a linear map $A_i\colon W_i\to W_{-i}$ intertwining $\sigma_i$ and $\sigma_{-i}^{\vee\tc}$ satisfying $A_{-i}\circ A_i=\id_{W_i}$. Put $\sigma\coloneqq\boxtimes_{i=-t}^t\sigma_{1,i}$ regarded as a representation of $P$ by inflation, which realizes on the space $W\coloneqq\bigotimes_{i=-t}^tW_i$; and put $A_\sigma\coloneqq\otimes_{i=-t}^tA_i\in\End(W)$. Choose an element $w\in\GL_N(F_w)$ satisfying $w=\tp{w}^\tc$, that $wPw^{-1}\cap P$ is the standard Levi subgroup of $P$, and that for $(a_{-t},\dots,a_t)\in wPw^{-1}\cap P$, we have $w(a_{-t},\dots,a_t)w^{-1}=(a_t,\dots,a_{-t})$.

We realize $\Pi_{1,w}$ on the space
\[
V_1\coloneqq\{f\colon\GL_N(F_w)\to W\res f(pg)=\delta_P^{1/2}(p)\sigma(p)f(g),p\in P,g\in\GL_N(F_w)\}.
\]
Define a linear map $A_{\Pi_{1,w}}\colon V_1\to V_1$ by the formula
\[
\(A_{\Pi_{1,w}}(f)\)(g)=A_\sigma\(f(w\tp{g}^{-1,\tc})\).
\]
Then it is clear that $A_{\Pi_{1,w}}$ is a intertwining operator for $\Pi_{1,w}$ satisfying $A_{\Pi_{1,w}}^2=1$. Similarly, we realize $\Pi_{2,w}$ on the space
\[
V_2\coloneqq\{f\colon\GL_N(F_w)\to W\res f(pg)=\delta_P^{1/2}(p)\chi(p)\sigma(p)f(g),p\in P,g\in\GL_N(F_w)\},
\]
where we put $\chi\coloneqq\boxtimes_{i=-t}^t\chi_i$ regarded as a character of $P$. We define $A_{\Pi_{2,w}}\colon V_2\to V_2$ by the same formula, which is a normalized intertwining operator for $\Pi_{2,w}$. The desired isomorphism $i$ is the map sending $f\in V_1$ to the unique function $i(f)$ such that $i(f)(g)=f(g)$ for $g\in\GL_N(O_{F_w})$. The lemma is proved.
\end{proof}

Now we can prove Theorem \ref{th:raising}.

\begin{proof}[Proof of Theorem \ref{th:raising}]
For (1), Assumption \ref{as:single_vanishing}, Lemma \ref{le:ns_cohomology}, and the spectral sequence in Lemma \ref{le:ns_weight_even} imply that $\rH^i_\fT(\ol\rM^\bullet_N,O_\lambda)_\fm$ is $O_\lambda$-torsion free for $i\neq 2r-1,2r$. By Proposition \ref{pr:raising_weight}(4) and Proposition \ref{pr:old_same}, we know that $\rH^{2r}_\fT(\ol\rM^\bullet_N,O_\lambda)_\fm$ is $O_\lambda$-torsion free. By the Poincar\'{e} duality theorem and Lemma \ref{le:single_compact}, we have
\begin{align*}
\rank_{O_\lambda}\rH^{2r}_\fT(\ol\rM^\bullet_N,O_\lambda)_\fm&=\rank_{O_\lambda}\rH^{2r-2}_\fT(\ol\rM^\bullet_N,O_\lambda)_\fm, \\
\dim_{O_\lambda/\lambda}\rH^{2r}_\fT(\ol\rM^\bullet_N,O_\lambda/\lambda)_\fm
&=\dim_{O_\lambda/\lambda}\rH^{2r-2}_\fT(\ol\rM^\bullet_N,O_\lambda/\lambda)_\fm,
\end{align*}
which imply that $\rH^{2r-1}_\fT(\ol\rM^\bullet_N,O_\lambda)_\fm$ is $O_\lambda$-torsion free as well by the universal coefficient theorem.

\if false

By definition, the universal homomorphism $r_\ram^\natural\colon\Gamma_F\to\GL_N(\sfR^\ram)$ has a direct sum decomposition $(\sfR^\ram)^{\oplus N}=\sfR_1\oplus\sfR_2$ in which $\sfR_1$ is a free $\sfR^\ram$-submodule of rank $N-2$, and $\sfR_2$ is the free $\sfR^\ram$-submodule of rank $2$ generated by (the image of) $\sfv$ and $\sfv'$. We have
\[
\rF_{-1}\rH^{2r-1}_\fT(\ol\rM_N,\rR\Psi O_\lambda)_\fm\subseteq\sfH\otimes_{\sfR^\ram}\sfR^\ram\sfv.
\]
On the other hand, we have
\[
\rH^{2r-1}_\fT(\ol\rM^\bullet_N,O_\lambda)_\fm=
\frac{\rF_0\rH^{2r-1}_\fT(\ol\rM_N,\rR\Psi O_\lambda)_\fm}{\rF_{-1}\rH^{2r-1}_\fT(\ol\rM_N,\rR\Psi O_\lambda)_\fm},
\]
and that the quotient $\tfrac{\rF_1\rH^{2r-1}_\fT(\ol\rM_N,\rR\Psi O_\lambda)_\fm}{\rF_0\rH^{2r-1}_\fT(\ol\rM_N,\rR\Psi O_\lambda)_\fm}$ is torsion free by Lemma \ref{le:ns_cohomology}. Thus, the $O_\lambda$-torsion of $\rH^{2r-1}_\fT(\ol\rM^\bullet_N,O_\lambda)_\fm$ coincides with
\[
\sfH\otimes_{\sfR^\ram}\sfR^\ram\sfv/\rF_{-1}\rH^{2r-1}_\fT(\ol\rM_N,\rR\Psi O_\lambda)_\fm,
\]
which is fixed by $\Gal(\ol\dF_p/\dF_{p^2})$. However, by Lemma \ref{le:ns_weight_even}(6) and Proposition \ref{pr:old_same}, the $O_\lambda$-torsion of $\rH^{2r-1}_\fT(\ol\rM^\bullet_N,O_\lambda)_\fm$ has to vanish. Part (1) is proved.

\fi

Part (2) is an immediate consequence of (1), Assumption \ref{as:single_vanishing}, Lemma \ref{le:ns_cohomology}, and the spectral sequence in Lemma \ref{le:ns_weight_even}.

Part (3) is a consequence of (1) and (P4) that $P_{\balpha(\Pi_\fp)}\modulo\lambda^m$ is level-raising special at $\fp$. In fact, we have an isomorphism
\[
\rH^{2r-1}_\fT(\ol\rM^\bullet_N,O_\lambda(r))\simeq\sfH\otimes_{\sfR^\ram}\sfR_1(r)
\]
of $O_\lambda[\Gal(\ol\dF_p/\dF_{p^2})]$-modules.

For (4), by Proposition \ref{pr:old_same} and (P6), it suffices to show that the two natural maps
\begin{align*}
\rF_{-1}\rH^1(\rI_{\dQ_{p^2}},\rH^{2r-1}_\fT(\ol\rM_N,\rR\Psi O_\lambda(r))_\fm)/\fn
&\to\rF_{-1}\rH^1(\rI_{\dQ_{p^2}},\rH^{2r-1}_\fT(\ol\rM_N,\rR\Psi O_\lambda(r))/\fn), \\
\rH^1_\sing(\dQ_{p^2},\rH^{2r-1}_\fT(\ol\rM_N,\rR\Psi O_\lambda(r))_\fm)/\fn
&\to\rH^1_\sing(\dQ_{p^2},\rH^{2r-1}_\fT(\ol\rM_N,\rR\Psi O_\lambda(r))/\fn),
\end{align*}
are both isomorphisms. Note that we have a short exact sequence
\[
\resizebox{\hsize}{!}{
\xymatrix{
0\to\rF_{-1}\rH^1(\rI_{\dQ_{p^2}},\rH^{2r-1}_\fT(\ol\rM_N,\rR\Psi O_\lambda(r))_\fm)
\to\rH^1(\rI_{\dQ_{p^2}},\rH^{2r-1}_\fT(\ol\rM_N,\rR\Psi O_\lambda(r))_\fm)
\to\frac{\rH^{2r-1}_\fT(\ol\rM_N,\rR\Psi O_\lambda(r))_\fm}{\rF_{-1}\rH^{2r-1}_\fT(\ol\rM_N,\rR\Psi O_\lambda(r))_\fm}
\to 0
}
}
\]
of $\dT^{\Sigma^+\cup\Sigma^+_p}_N$-modules, which is split by considering $\Gal(\ol\dF_p/\dF_{p^2})$ actions and (3). Thus, the first isomorphism is confirmed. The second one is also confirmed as, by (3), one can replace $\Gal(\ol\dF_p/\dF_{p^2})$-invariants by $\Gal(\ol\dF_p/\dF_{p^2})$-coinvariants. Part (4) is proved.

For (5), we have
\[
\rH^{2r-1}_\et(\Sh(\rV'_N,\tj_N\rK^{p\circ}_N\rK'_{p,N})_{\ol{F}},O_\lambda(r))/\fn
\simeq\sfH\otimes_{\sfR^\ram/\fn}(\sfR^\ram/\fn)^{\oplus N}(r)
\]
by Lemma \ref{le:ns_pbc}. Here, we regard $\fn$ as its image in $\sfT^\ram_\fm$, where the latter is canonically isomorphic to $\sfR^\ram$. We claim that $O_\lambda/\lambda^m=\sfR^\ram/\fn$ and $(\sfR^\ram/\fn)^{\oplus N}(r)\simeq\bar\rR^{(m)\tc}$ as $(O_\lambda/\lambda^m)[\Gamma_F]$-modules, where we recall that $\Gamma_F$ acts on $(\sfR^\ram/\fn)^{\oplus N}$ via $r_\ram^{\natural,\tc}$. Since $\fn$ satisfies $\fn\cap O_\lambda=\lambda^mO_\lambda$, the structure homomorphism $O_\lambda\to\sfR^\ram$ induces an equality $O_\lambda/\lambda^m=\sfR^\ram/\fn$. Now by the Chebotarev density theorem and \cite{Car94}*{Th\'{e}or\`{e}me~1}, we know that the two liftings $(\sfR^\ram/\fn)^{\oplus N}(r)$ and $\bar\rR^{(m)\tc}$ of $\bar\rho_{\Pi,\lambda}^\tc(r)$ to $O_\lambda/\lambda^m$ have to be isomorphic.

Theorem \ref{th:raising} is all proved.
\end{proof}

\if false

\begin{theorem}\label{th:raising}
Assume Assumptions \ref{as:single_vanishing}, \ref{as:single_irreducible}, and Hypothesis \ref{hy:unitary_cohomology} for $N$. Moreover, if $N\geq 4$, we further assume that
\begin{enumerate}[label=(\alph*)]
  \item $\ell\geq 2(N+1)$ and $\ell$ is unramified in $F$;

  \item $\bar\rho_{\Pi,\lambda,+}$ (Remark \ref{re:single_irreducible}) is rigid for $(\Sigma^+_\mnm,\Sigma^+_\lr)$ (Definition \ref{de:rigid}), and $\bar\rho_{\Pi,\lambda}\res_{\Gal(\ol{F}/F(\zeta_\ell))}$ is absolutely irreducible;

  \item the composite homomorphism $\dT^{\Sigma^+}_N\xrightarrow{\phi_\Pi}O_E\to O_E/\lambda$ is cohomologically generic (Definition \ref{de:generic}); and

  \item $O_\lambda[\Sh(\rV^\circ_N,\rK^\circ_N)]_\fm$ is nontrivial.
\end{enumerate}
Then we have
\begin{enumerate}
  \item $\rH^i_\fT(\ol\rM^\bullet_N,O_\lambda)_\fm$ is a free $O_\lambda$-module for every $i\in\dZ$.

  \item $\rE^{p,q}_{2,\fm}$ is a free $O_\lambda$-module, and vanishes if $(p,q)\not\in\{(-1,2r),(0,2r-1),(1,2r-2)\}$.

  \item If we denote by $\{\alpha_1^{\pm1},\dots,\alpha_r^{\pm1}\}$ the roots of $P_{\balpha(\Pi_\fp)}\modulo\lambda$ in a finite extension of $O_\lambda/\lambda$, then the generalized Frobenius eigenvalues of the $(O_\lambda/\lambda)[\Gal(\ol\dF_p/\dF_{p^2})]$-module $\rH^{2r-1}_\fT(\ol\rM^\bullet_N,O_\lambda(r))_\fm\otimes_{O_\lambda}O_\lambda/\lambda$ is contained in $\{p\alpha_1^{\pm1},\dots,p\alpha_r^{\pm1}\}\setminus\{1,p^2\}$.

  \item The map in Proposition \ref{pr:raising_weight}(4) (see Remark \ref{re:raising_weight} for $N=2$) factors through a map
     \[
     \nabla^0_{/\fn}\colon\rF_{-1}\rH^1(\rI_{\dQ_{p^2}},\rH^{2r-1}_\fT(\ol\rM_N,\rR\Psi O_\lambda(r))/\fn)
     \to O_\lambda[\Sh(\rV^\circ_N,\rK^\circ_N)]/\fn
     \]
     which is an isomorphism. The map from Lemma \ref{le:ns_weight_even}(6) induces a canonical isomorphism
      \[
      \rF_{-1}\rH^1(\rI_{\dQ_{p^2}},\rH^{2r-1}_\fT(\ol\rM_N,\rR\Psi O_\lambda(r))/\fn)
      \xrightarrow{\sim}\rH^1_\sing(\dQ_{p^2},\rH^{2r-1}_\fT(\ol\rM_N,\rR\Psi O_\lambda(r))/\fn).
      \]

  \item There exist finitely many positive integers $m_1,\dots,m_\mu$ at most $m$ such that we have an isomorphism
     \[
     \rH^{2r-1}_\et(\Sh(\rV'_N,\tj_N\rK^{p\circ}_N\rK'_{p,N})_{\ol{F}},O_\lambda(r))/\fn
     \simeq\bigoplus_{i=1}^\mu\bar\rR^{(m_i)\tc}
     \]
     of $O_\lambda[\Gamma_F]$-modules, where $\rR$ is the $\Gamma_F$-stable $O_\lambda$-lattice in $\rho_{\Pi,\lambda}(r)$, unique up to homothety.
\end{enumerate}
\end{theorem}

\begin{remark}
In fact, from the proof, one sees that when $N\geq 4$, we can take $m_1=\cdots=m_\mu=m$ due to our strong extra assumptions.
\end{remark}

\subsection{Proof of Theorem \ref{th:raising}}
\label{ss:raising2}

The proof differs according to $N=2$ or $N\geq 4$ as we can see from the extra assumptions. We start from the much more difficult case where $N\geq 4$.

We apply the discussion of \cite{LTXZZ}*{\S3} to the pair $(\bar{r},\chi)$, where
\[
\bar{r}\coloneqq\bar\rho_{\Pi,\lambda,+}\colon\Gamma_{F^+}\to\sG_N(O_\lambda/\lambda)
\]
and $\chi\coloneqq\epsilon_\ell^{1-N}$ for the similitude character. Then $\bar{r}$ is rigid for $(\Sigma^+_\mnm,\Sigma^+_\lr)$, and also for $(\Sigma^+_\mnm,\Sigma^+_\lr\cup\{\fp\})$ by (P4).

For $?=\mix,\unr,\ram$, consider a global deformation problem \cite{LTXZZ}*{Definition~3.1.6}
\[
\sS^?\coloneqq(\bar{r},\eta_{F/F^+}^\mu\epsilon_\ell^{1-N},\Sigma^+_\mnm\cup\Sigma^+_\lr\cup\{\fp\}\cup\Sigma^+_\ell,
\{\sD_v\}_{v\in\Sigma^+_\mnm\cup\Sigma^+_\lr\cup\{\fp\}\cup\Sigma^+_\ell})
\]
where
\begin{itemize}[label={\ding{109}}]
  \item for $v\in\Sigma^+_\mnm$, $\sD_v$ is the local deformation problem classifying all liftings of $\bar{r}_v$;

  \item for $v\in\Sigma^+_\lr$, $\sD_v$ is the local deformation problem $\sD^\ram$ of $\bar{r}_v$ from \cite{LTXZZ}*{Definition~3.5.1};

  \item for $v=\fp$, $\sD_v$ is the local deformation problem $\sD^?$ of $\bar{r}_v$ from \cite{LTXZZ}*{Definition~3.5.1};

  \item for $v\in\Sigma^+_\ell$, $\sD_v$ is the local deformation problem $\sD^\FL$ of $\bar{r}_v$ from \cite{LTXZZ}*{Definition~3.2.6}.
\end{itemize}
Then we have the \emph{global universal deformation ring} $\sfR^\univ_{\sS^?}$ from \cite{LTXZZ}*{Proposition~3.1.7}. Put $\sfR^?\coloneqq\sfR^\univ_{\sS^?}$ for short. Then we have canonical surjective homomorphisms $\sfR^\mix\to\sfR^\unr$ and $\sfR^\mix\to\sfR^\ram$ of $O_\lambda$-rings. Finally, put
\[
\sfR^{\r{cong}}\coloneqq\sfR^\unr\otimes_{\sfR^\mix}\sfR^\ram.
\]
We fix a universal lifting
\[
r_\mix\colon\Gamma_{F^+}\to\sG_N(\sfR^\mix)
\]
of $\bar{r}$, which induce a continuous homomorphism
\[
r_\mix^\natural\colon\Gamma_F\to\GL_N(\sfR^\mix)
\]
by restriction (Notation \ref{no:sg_extension}). By pushforward, $\sfR^{\r{cong}}$ also induces homomorphisms
\[
r_\unr\colon\Gamma_{F^+}\to\sG_N(\sfR^\unr),\quad
r_\ram\colon\Gamma_{F^+}\to\sG_N(\sfR^\ram).
\]
Denote by $\rP_{F^+_\fp}$ the maximal closed subgroup of the inertia subgroup $\rI_{F^+_\fp}\subseteq\Gamma_{F^+_\fp}$ of pro-order coprime to $\ell$. Then $\Gamma_{F^+_\fp}/\rP_{F^+_\fp}\simeq t^{\dZ_\ell}\rtimes\phi_p^{\widehat\dZ}$ is a $p$-tame group \cite{LTXZZ}*{Definition~3.3.1}. By definition, the homomorphism $r_\mix^\natural$ is trivial on $\rP_{F^+_\fp}$. Let $\bar\sfv$ and $\bar\sfv'$ be eigenvectors in $(O_\lambda/\lambda)^{\oplus N}$ for $\bar{r}^\natural(\phi_p^2)$ with eigenvalues $p^{-2r}$ and $p^{-2r+2}$, respectively. By Hensel's lemma, $\bar\sfv$ and $\bar\sfv'$ lift to eigenvectors $\sfv$ and $\sfv'$ in $(\sfR^\mix)^{\oplus N}$ for $r_\mix^\natural(\phi_p^2)$, with eigenvalues $\sfs$ and $\sfs'$ in $\sfR^\mix$ lifting $p^{-2r}$ and $p^{-2r+2}$, respectively. Let $\sfx\in\sfR^\mix$ be the unique element such that $r_\mix^\natural(t)\sfv'=\sfx\sfv+\sfv'$. Then we must have $\sfx(\sfs-p^{-2r})=0$. By \cite{LTXZZ}*{Definition~3.5.1}, we have
\[
\sfR^\unr=\sfR^\mix/(\sfx),\quad
\sfR^\ram=\sfR^\mix/(\sfs-p^{-2r}),\quad
\sfR^{\r{cong}}=\sfR^\mix/(\sfs-p^{-2r},\sfx).
\]

Let $\sfT^\unr$ be the image of $\dT^{\Sigma^+}_N$ in $\End_{O_\lambda}(O_\lambda[\Sh(\rV^\circ_N,\rK^\circ_N)])$. By (d) in Theorem \ref{th:raising}, we know that $\sfT^\unr_\fm\neq 0$. Thus by \cite{LTXZZ}*{Theorem~3.6.3}, we have a canonical isomorphism $\sfR^\unr\xrightarrow{\sim}\sfT^\unr_\fm$ so that $O_\lambda[\Sh(\rV^\circ_N,\rK^\circ_N)]_\fm$ is canonically a free $\sfR^\unr$-module of rank $d_\unr>0$.\footnote{Here, we also need the easy fact that $\sfT^\unr_\fm$ and $O_\lambda[\Sh(\rV^\circ_N,\rK^\circ_N)]_\fm$ do not change if we replace $\fm$ by $\fm\cap\dT^{\Sigma^+\cup\Sigma^+_p\cup\Sigma^+_\ell}_N$} We may write the eigenvalues of $r_\unr^\natural(\phi_p^2)$ by $\{p^{-2r+1}\alpha_1^{\pm 1},\dots,p^{-2r+1}\alpha_{r-1}^{\pm 1},\sfs,\sfs'=p^{-4r+2}\sfs^{-1}\}$ with $\alpha_1,\dots,\alpha_{r-1}$ in a certain finite flat extension of $\sfR^\unr$ that are not congruent to $p$ or $p^{-1}$ in $O_\lambda/\lambda$. By Proposition \ref{pr:enumeration_even}(2), we have
\[
((p+1)\tR^\circ_{N,\fp}-\tI^\circ_{N,\fp}).O_\lambda[\Sh(\rV^\circ_N,\rK^\circ_N)]_\fm
=(\sfs-p^{-2r}).O_\lambda[\Sh(\rV^\circ_N,\rK^\circ_N)]_\fm.
\]
In particular, we have
\[
O_\lambda[\Sh(\rV^\circ_N,\rK^\circ_N)]_\fm/((p+1)\tR^\circ_{N,\fp}-\tI^\circ_{N,\fp})
=O_\lambda[\Sh(\rV^\circ_N,\rK^\circ_N)]_\fm\otimes_{\sfR^\unr}\sfR^{\r{cong}},
\]
which is a free $\sfR^{\r{cong}}$-module of rank $d_\unr$.

On the other hand, let $\sfT^\ram$ be the image of $\dT^{\Sigma^+\cap\Sigma^+_p}_N$ in $\End_{O_\lambda}(\rH^{2r-1}_\fT(\ol\rM_N,\rR\Psi O_\lambda))$. By Proposition \ref{pr:raising_weight}(4) and Lemma \ref{le:ns_weight_even}(6), we know that $\sfT^\ram_\fm\neq0$. Thus by Lemma \ref{le:ns_pbc} and \cite{LTXZZ}*{Theorem~3.6.3} (with $(\Sigma^+_\mnm,\Sigma^+_\lr)$ replaced by $(\Sigma^+_\mnm,\Sigma^+_\lr\cup\{\fp\})$), we have a canonical isomorphism $\sfR^\ram\xrightarrow{\sim}\sfT^\ram_\fm$ so that $\rH^{2r-1}_\fT(\ol\rM_N,\rR\Psi O_\lambda)_\fm$ is canonically a free $\sfR^\ram$-module.\footnote{Here, we also need the fact that $\sfT^\ram_\fm$ and $\rH^{2r-1}_\fT(\ol\rM_N,\rR\Psi O_\lambda)_\fm$ do not change if we replace $\fm$ by $\fm\cap\dT^{\Sigma^+\cup\Sigma^+_p\cup\Sigma^+_\ell}_N$, which is a consequence of Theorem \ref{th:raising}(c).} Define the $\sfR^\ram$-module
\[
\sfH\coloneqq\Hom_{\Gamma_F}\((\sfR^\ram)^{\oplus N},\rH^{2r-1}_\fT(\ol\rM_N,\rR\Psi O_\lambda)_\fm\)
\]
where $\Gamma_F$ acts on $(\sfR^\ram)^{\oplus N}$ via the homomorphism $r_\ram^{\natural,\tc}$. By the same argument for \cite{Sch18}*{Theorem~5.6} (using Proposition \ref{pr:arthur} and Hypothesis \ref{hy:unitary_cohomology}, here), we have a canonical isomorphism
\[
\rH^{2r-1}_\fT(\ol\rM_N,\rR\Psi O_\lambda)_\fm\simeq\sfH\otimes_{\sfR^\ram}(\sfR^\ram)^{\oplus N}
\]
of $\sfR^\ram[\Gamma_F]$-modules. Since $\sfR^\ram$ is a local ring, $\sfH$ is a free $\sfR^\ram$-module, say of rank $d_\ram$. If we still denote by $\sfv$ and $\sfv'$ for their projection in $(\sfR^\ram)^{\oplus N}$, then it is easy to see that
\[
\rH^1_\sing(\dQ_{p^2},(\sfR^\ram)^{\oplus N}(r))=\sfR^\ram\sfv/\sfx\sfv\simeq\sfR^\ram/(\sfx)=\sfR^{\r{cong}}.
\]
Thus, we obtain
\[
\rH^1_\sing(\dQ_{p^2},\rH^{2r-1}_\fT(\ol\rM_N,\rR\Psi O_\lambda(r))_\fm)\simeq \sfH\otimes_{\sfR^\ram}\rH^1_\sing(\dQ_{p^2},(\sfR^\ram)^{\oplus N}(r))
\simeq\sfH\otimes_{\sfR^\ram}\sfR^{\r{cong}},
\]
which is a free $\sfR^{\r{cong}}$-module of rank $d_\ram>0$.

\begin{proposition}\label{pr:old_same}
Under the assumptions of Theorem \ref{th:raising}, we have $d_\unr=d_\ram$. In particular, the two canonical maps
\begin{align*}
&\rF_{-1}\rH^1(\rI_{\dQ_{p^2}},\rH^{2r-1}_\fT(\ol\rM_N,\rR\Psi O_\lambda(r))_\fm)\to
O_\lambda[\Sh(\rV^\circ_N,\rK^\circ_N)]_\fm/((p+1)\tR^\circ_{N,\fp}-\tI^\circ_{N,\fp}), \\
&\rF_{-1}\rH^1(\rI_{\dQ_{p^2}},\rH^{2r-1}_\fT(\ol\rM_N,\rR\Psi O_\lambda(r))_\fm)\to
\rH^1_\sing(\dQ_{p^2},\rH^{2r-1}_\fT(\ol\rM_N,\rR\Psi O_\lambda(r))_\fm),
\end{align*}
from Proposition \ref{pr:raising_weight}(4) and Lemma \ref{le:ns_weight_even}(6), respectively, are both isomorphisms.
\end{proposition}

\begin{proof}
By Proposition \ref{pr:raising_weight}(4), the first map is surjective. By Lemma \ref{le:ns_weight_even}(6), the second map is injective. Thus, we must have $d_\ram\geq d_\unr>0$ by the previous discussion.

Take a geometric point $\eta_1\in(\Spec\sfR^\unr)(\ol\dQ_\ell)$ in the support of $O_\lambda[\Sh(\rV^\circ_N,\rK^\circ_N)]_\fm$, which corresponds to a relevant representation $\Pi_1$ of $\GL_N(\dA_F)$ by Lemma \ref{le:single_congruent}, such that $\rho_{\Pi_1,\iota_\ell}$ is residually isomorphic to $\bar\rho_{\Pi,\lambda}\otimes_{O_\lambda/\lambda}\ol\dF_\ell$. Then we have
\[
d_\unr=\dim\ol\dQ_\ell[\Sh(\rV^\circ_N,\rK^\circ_N)][\iota_\ell\phi_{\Pi_1}].
\]

Take a geometric point $\eta_2\in(\Spec\sfR^\ram)(\ol\dQ_\ell)$ in the support of $\rH^{2r-1}_\fT(\ol\rM_N,\rR\Psi O_\lambda)_\fm$, which corresponds to a relevant representation $\Pi_2$ of $\GL_N(\dA_F)$ by Lemma \ref{le:single_congruent}, such that $\rho_{\Pi_2,\iota_\ell}$ is residually isomorphic to $\bar\rho_{\Pi,\lambda}\otimes_{O_\lambda/\lambda}\ol\dF_\ell$. Then we have
\[
d_\ram=\dim\rH^{2r-1}_\fT(\ol\rM_N,\rR\Psi\ol\dQ_\ell)[\iota_\ell\phi_{\Pi_2}]
=\dim\rH^{2r-1}_\et(\Sh(\rV'_N,\tj_N\rK^{p\circ}_N\rK'_{p,N})_{\ol{F}},\ol\dQ_\ell)[\iota_\ell\phi_{\Pi_2}]
\]
by Lemma \ref{le:ns_pbc}. By Proposition \ref{pr:old} and Lemma \ref{le:old} below, we have $d_\unr=d_\ram$. The proposition follows.
\end{proof}

\begin{lem}\label{le:old}
Let $\Pi_1$ and $\Pi_2$ be two relevant representations of $\GL_N(\dA_F)$ such that the associated Galois representations $\rho_{\Pi_1,\iota_\ell}$ and $\rho_{\Pi_2,\iota_\ell}$ are both residually isomorphic to $\bar\rho_{\Pi,\lambda}\otimes_{O_\lambda/\lambda}\ol\dF_\ell$. For every $v\in\Sigma^+_\mnm$ (so that every lifting of $\bar\rho_{\Pi,\lambda,+,v}$ is minimally ramified), if we realize $\Pi_{1,v}$ and $\Pi_{2,v}$ on vector spaces $V_1$ and $V_2$, respectively, then there exist normalized intertwining operators $A_{\Pi_{1,v}}$ and $A_{\Pi_{2,v}}$ for $\Pi_{1,v}$ and $\Pi_{2,v}$ \cite{Shi11}*{\S4.1}, respectively, such that we have an $\GL_N(O_{F_v})$-equivariant isomorphism $i\colon V_1\xrightarrow{\sim}V_2$ satisfying $i\circ A_{\Pi_{1,v}}=A_{\Pi_{2,v}}\circ i$.
\end{lem}

\begin{proof}
We will give the proof when $v$ is nonsplit in $F$, and leave the other similar case to the readers. Let $w$ be the unique place of $F$ above $v$.

By Proposition \ref{pr:galois}(1), both $\Pi_{1,w}$ and $\Pi_{2,w}$ are tempered. Thus by the Bernstein--Zelevinsky classification, for $\alpha=1,2$, we can write
\[
\Pi_{\alpha,w}=\rI^{\GL_N(F_w)}_{P_\alpha}\(\sigma_{\alpha,-t_\alpha}\boxtimes\cdots\boxtimes\sigma_{\alpha,-1}
\boxtimes\sigma_{\alpha,0}\boxtimes\sigma_{\alpha,1}\boxtimes\cdots\boxtimes\sigma_{\alpha,t_\alpha}\)
\]
for some integer $t_\alpha\geq 0$, some standard parabolic subgroup $P_\alpha\subseteq\GL_N(F_w)$, and some (unitary) discrete series representations $\{\sigma_{\alpha,-t_\alpha},\dots,\sigma_{\alpha,t_\alpha}\}$ satisfying $\sigma_{\alpha,-i}\simeq\sigma_{\alpha,i}^{\vee\tc}$. See \S\ref{ss:local_bc} for the notation on parabolic induction.

By \cite{LTXZZ}*{Proposition~3.4.12(3)} and \cite{BLGGT}*{Lemma~1.3.4(2)}, we know that $\rho_{\Pi_1,\iota_\ell}\res_{\rI_{F_w}}$ and $\rho_{\Pi_2,\iota_\ell}\res_{\rI_{F_w}}$ are conjugate. Thus, by \cite{Yao}*{Lemma~3.6}, we have $P_1=P_2$ (say $P$) and $t_1=t_2$ (say $t$), and we assume that there are unramified (unitary) characters $\{\chi_{-t},\dots,\chi_t\}$ of $F_w^\times$ satisfying $\chi_{-i}\simeq\chi_i^{-1}$ such that $\sigma_{2,i}=\sigma_{1,i}\otimes\chi_i$. For every $i$, we choose a vector space $W_i$ on which $\sigma_{1,i}$ realizes (and also realize $\sigma_{1,i}^{\vee\tc}$ on $W_i$ via $g\mapsto\tp{g}^{-1,\tc}$), and fix a linear map $A_i\colon W_i\to W_{-i}$ intertwining $\sigma_i$ and $\sigma_{-i}^{\vee\tc}$ satisfying $A_{-i}\circ A_i=\id_{W_i}$. Put $\sigma\coloneqq\boxtimes_{i=-t}^t\sigma_{1,i}$ regarded as a representation of $P$ by inflation, which realizes on the space $W\coloneqq\bigotimes_{i=-t}^tW_i$; and put $A_\sigma\coloneqq\otimes_{i=-t}^tA_i\in\End(W)$. Choose an element $w\in\GL_N(F_w)$ satisfying $w=\tp{w}^\tc$, that $wPw^{-1}\cap P$ is the standard Levi subgroup of $P$, and that for $(a_{-t},\dots,a_t)\in wPw^{-1}\cap P$, we have $w(a_{-t},\dots,a_t)w^{-1}=(a_t,\dots,a_{-t})$.

We realize $\Pi_{1,w}$ on the space
\[
V_1\coloneqq\{f\colon\GL_N(F_w)\to W\res f(pg)=\delta_P^{1/2}(p)\sigma(p)f(g),p\in P,g\in\GL_N(F_w)\}.
\]
Define a linear map $A_{\Pi_{1,w}}\colon V_1\to V_1$ by the formula
\[
\(A_{\Pi_{1,w}}(f)\)(g)=A_\sigma\(f(w\tp{g}^{-1,\tc})\).
\]
Then it is clear that $A_{\Pi_{1,w}}$ is a intertwining operator for $\Pi_{1,w}$ satisfying $A_{\Pi_{1,w}}^2=1$. Similarly, we realize $\Pi_{2,w}$ on the space
\[
V_2\coloneqq\{f\colon\GL_N(F_w)\to W\res f(pg)=\delta_P^{1/2}(p)\chi(p)\sigma(p)f(g),p\in P,g\in\GL_N(F_w)\},
\]
where we put $\chi\coloneqq\boxtimes_{i=-t}^t\chi_i$ regarded as a character of $P$. We define $A_{\Pi_{2,w}}\colon V_2\to V_2$ by the same formula, which is a normalized intertwining operator for $\Pi_{2,w}$. The desired isomorphism $i$ is the map sending $f\in V_1$ to the unique function $i(f)$ such that $i(f)(g)=f(g)$ for $g\in\GL_N(O_{F_w})$. The lemma is proved.
\end{proof}

Now we can prove Theorem \ref{th:raising} when $N\geq 4$.

\begin{proof}[Proof of Theorem \ref{th:raising} when $N\geq 4$]
For (1), Assumption \ref{as:single_vanishing}, Lemma \ref{le:ns_cohomology}, and the spectral sequence in Lemma \ref{le:ns_weight_even} imply that $\rH^i_\fT(\ol\rM^\bullet_N,O_\lambda)_\fm$ is $O_\lambda$-torsion free for $i\neq 2r-1,2r$. By Proposition \ref{pr:raising_weight}(4) and Proposition \ref{pr:old_same}, we know that $\rH^{2r}_\fT(\ol\rM^\bullet_N,O_\lambda)_\fm$ is $O_\lambda$-torsion free. It remains to show that $\rH^{2r-1}_\fT(\ol\rM^\bullet_N,O_\lambda)_\fm$ is $O_\lambda$-torsion free as well.

By definition, the universal homomorphism $r_\ram^\natural\colon\Gamma_F\to\GL_N(\sfR^\ram)$ has a direct sum decomposition $(\sfR^\ram)^{\oplus N}=\sfR_1\oplus\sfR_2$ in which $\sfR_1$ is a free $\sfR^\ram$-submodule of rank $N-2$, and $\sfR_2$ is the free $\sfR^\ram$-submodule of rank $2$ generated by (the image of) $\sfv$ and $\sfv'$. We have
\[
\rF_{-1}\rH^{2r-1}_\fT(\ol\rM_N,\rR\Psi O_\lambda)_\fm\subseteq\sfH\otimes_{\sfR^\ram}\sfR^\ram\sfv.
\]
On the other hand, we have
\[
\rH^{2r-1}_\fT(\ol\rM^\bullet_N,O_\lambda)_\fm=
\frac{\rF_0\rH^{2r-1}_\fT(\ol\rM_N,\rR\Psi O_\lambda)_\fm}{\rF_{-1}\rH^{2r-1}_\fT(\ol\rM_N,\rR\Psi O_\lambda)_\fm},
\]
and that the quotient $\tfrac{\rF_1\rH^{2r-1}_\fT(\ol\rM_N,\rR\Psi O_\lambda)_\fm}{\rF_0\rH^{2r-1}_\fT(\ol\rM_N,\rR\Psi O_\lambda)_\fm}$ is torsion free by Lemma \ref{le:ns_cohomology}. Thus, the $O_\lambda$-torsion of $\rH^{2r-1}_\fT(\ol\rM^\bullet_N,O_\lambda)_\fm$ coincides with
\[
\sfH\otimes_{\sfR^\ram}\sfR^\ram\sfv/\rF_{-1}\rH^{2r-1}_\fT(\ol\rM_N,\rR\Psi O_\lambda)_\fm,
\]
which is fixed by $\Gal(\ol\dF_p/\dF_{p^2})$. However, by Lemma \ref{le:ns_weight_even}(6) and Proposition \ref{pr:old_same}, the $O_\lambda$-torsion of $\rH^{2r-1}_\fT(\ol\rM^\bullet_N,O_\lambda)_\fm$ has to vanish. Part (1) is proved.

Part (2) is an immediate consequence of (1), Assumption \ref{as:single_vanishing}, Lemma \ref{le:ns_cohomology}, and the spectral sequence in Lemma \ref{le:ns_weight_even}.

Part (3) is a consequence of (1) and (P4) that $P_{\balpha(\Pi_\fp)}\modulo\lambda^m$ is level-raising special at $\fp$. In fact, we have an isomorphism
\[
\rH^{2r-1}_\fT(\ol\rM^\bullet_N,O_\lambda(r))\simeq\sfH\otimes_{\sfR^\ram}\sfR_1(r)
\]
of $O_\lambda[\Gal(\ol\dF_p/\dF_{p^2})]$-modules.

For (4), by Proposition \ref{pr:old_same} and (P6), it suffices to show that the two natural maps
\begin{align*}
\rF_{-1}\rH^1(\rI_{\dQ_{p^2}},\rH^{2r-1}_\fT(\ol\rM_N,\rR\Psi O_\lambda(r))_\fm)/\fn
&\to\rF_{-1}\rH^1(\rI_{\dQ_{p^2}},\rH^{2r-1}_\fT(\ol\rM_N,\rR\Psi O_\lambda(r))/\fn), \\
\rH^1_\sing(\dQ_{p^2},\rH^{2r-1}_\fT(\ol\rM_N,\rR\Psi O_\lambda(r))_\fm)/\fn
&\to\rH^1_\sing(\dQ_{p^2},\rH^{2r-1}_\fT(\ol\rM_N,\rR\Psi O_\lambda(r))/\fn),
\end{align*}
are both isomorphisms. Note that we have a short exact sequence
\[
\resizebox{\hsize}{!}{
\xymatrix{
0\to\rF_{-1}\rH^1(\rI_{\dQ_{p^2}},\rH^{2r-1}_\fT(\ol\rM_N,\rR\Psi O_\lambda(r))_\fm)
\to\rH^1(\rI_{\dQ_{p^2}},\rH^{2r-1}_\fT(\ol\rM_N,\rR\Psi O_\lambda(r))_\fm)
\to\frac{\rH^{2r-1}_\fT(\ol\rM_N,\rR\Psi O_\lambda(r))_\fm}{\rF_{-1}\rH^{2r-1}_\fT(\ol\rM_N,\rR\Psi O_\lambda(r))_\fm}
\to 0
}
}
\]
of $\dT^{\Sigma^+\cup\Sigma^+_p}_N$-modules, which is split by considering $\Gal(\ol\dF_p/\dF_{p^2})$ actions and (3). Thus, the first isomorphism is confirmed. The second one is also confirmed as, by (3), one can replace $\Gal(\ol\dF_p/\dF_{p^2})$-invariants by $\Gal(\ol\dF_p/\dF_{p^2})$-coinvariants. Part (4) is proved.

For (5), we have
\[
\rH^{2r-1}_\et(\Sh(\rV'_N,\tj_N\rK^{p\circ}_N\rK'_{p,N})_{\ol{F}},O_\lambda(r))/\fn
\simeq\sfH\otimes_{\sfR^\ram/\fn}(\sfR^\ram/\fn)^{\oplus N}(r)
\]
by Lemma \ref{le:ns_pbc}. Here, we regard $\fn$ as its image in $\sfT^\ram_\fm$, where the latter is canonically isomorphic to $\sfR^\ram$. We claim that $O_\lambda/\lambda^m=\sfR^\ram/\fn$ and $(\sfR^\ram/\fn)^{\oplus N}(r)\simeq\bar\rR^{(m)\tc}$ as $(O_\lambda/\lambda^m)[\Gamma_F]$-modules, where we recall that $\Gamma_F$ acts on $(\sfR^\ram/\fn)^{\oplus N}$ via $r_\ram^{\natural,\tc}$. By definition, $\fn$ is the kernel of homomorphism
\[
\dT^{\Sigma^+\cup\Sigma^+_p}_N\xrightarrow{\phi_\Pi}O_E\to O_E/\lambda^m,
\]
which satisfies $\fn\cap O_\lambda=\lambda^mO_\lambda$. Thus, the structure homomorphism $O_\lambda\to\sfR^\ram$ induces an equality $O_\lambda/\lambda^m=\sfR^\ram/\fn$. Now by the Chebotarev density theorem and \cite{Car94}*{Th\'{e}or\`{e}me~1}, we know that the two liftings $(\sfR^\ram/\fn)^{\oplus N}(r)$ and $\bar\rR^{(m)\tc}$ of $\bar\rho_{\Pi,\lambda}^\tc(r)$ to $O_\lambda/\lambda^m$ have to be isomorphic.

Theorem \ref{th:raising} is all proved when $N\geq 4$.
\end{proof}

\begin{proof}[Proof of Theorem \ref{th:raising} when $N=2$]
Part (1) is trivial since $\ol\rM^\bullet_2$ is a disjoint union of projective lines.

Part (2) follows from (1) by the same reason as for $N\geq 4$.

Part (3) is trivial as $\rH^1_\fT(\ol\rM^\bullet_2,O_\lambda(1))=0$.

For (4), from Remark \ref{re:raising_weight}, we know that the natural map
\[
\rF_{-1}\rH^1(\rI_{\dQ_{p^2}},\rH^1_\fT(\ol\rM_2,\rR\Psi O_\lambda(1))_\fm)\to
O_\lambda[\Sh(\rV^\circ_2,\rK^\circ_2)]_\fm/((p+1)\tR^\circ_{2,\fp}-\tI^\circ_{2,\fp})
\]
is an isomorphism. By (3) and Lemma \ref{le:ns_weight_even}(6), the natural map
\[
\rF_{-1}\rH^1(\rI_{\dQ_{p^2}},\rH^1_\fT(\ol\rM_2,\rR\Psi O_\lambda(1))_\fm)\to
\rH^1_\sing(\dQ_{p^2},\rH^1_\fT(\ol\rM_2,\rR\Psi O_\lambda(1))_\fm)
\]
is an isomorphism as well. Passing to the quotient by $\fn$ follows from the same argument as for $N\geq 4$.

For (5), let $\sfT^\ram$ be the image of $\dT^{\Sigma^+\cap\Sigma^+_p}_2$ in $\End_{O_\lambda}(\rH^1_\fT(\ol\rM_2,\rR\Psi O_\lambda(1)))$. Then by the same argument for \cite{Sch18}*{Theorem~5.6}, one have an isomorphism
\[
\rH^1_\fT(\ol\rM_2,\rR\Psi O_\lambda(1))_\fm\simeq\sfH\otimes_{\sfT^\ram_\fm}(\sfT^\ram_\fm)^{\oplus 2}
\]
of $\sfT^\ram_\fm[\Gamma_F]$-modules, where $\sfH$ is a (finitely generated) $\sfT^\ram_\fm$-module, and $\Gamma_F$ acts on the factor $(\sfT^\ram_\fm)^{\oplus 2}$ by some continuous homomorphism which lifts $\bar\rho^\tc_{\Pi,\lambda}(r)$ (from $O_\lambda/\lambda=\sfT^\ram_\fm/\fm$ to $\sfT^\ram_\fm$). Clearly, the natural homomorphism $O_\lambda/\lambda^m\to\sfT^\ram_\fm/\fn$ is an isomorphism. Then as an $O_\lambda/\lambda^m$-module, $\sfH\otimes_{\sfT^\ram_\fm}\sfT^\ram_\fm/\fn$ is isomorphic to $\bigoplus_{i=1}^\mu O_\lambda/\lambda^{m_i}$ for some positive integers $m_1,\dots,m_\mu$ at most $m$. Thus, it remains to show that $(\sfT^\ram_\fm/\fn)^{\oplus 2}$ and $\bar\rR^{(m)\tc}$ are isomorphic as deformations of $\bar\rho^\tc_{\Pi,\lambda}(r)$. But this is a consequence of the Eichler--Shimura relation for the unitary Shimura curve $\Sh(\rV'_2,\tj_2\rK^{p\circ}_2\rK'_{p,2})$ \cite{Liu3}*{Corollary~D.9}, the Chebotarev density theorem, and \cite{Car94}*{Th\'{e}or\`{e}me~1}.

Theorem \ref{th:raising} is all proved when $N=2$.
\end{proof}

\fi

\section{Explicit reciprocity laws for Rankin--Selberg motives}
\label{ss:7}

In this section, we state and prove the two explicit reciprocity laws for automorphic Rankin--Selberg motives. In \S\ref{ss:setup}, we setup the stage for automorphic Rankin--Selberg motives. In \S\ref{ss:first_reciprocity} and \S\ref{ss:second_reciprocity}, we state and prove our first and second explicit reciprocity law, respectively.

\subsection{Setup for automorphic Rankin--Selberg motives}
\label{ss:setup}

Let $n\geq 2$ be an integer. We denote by $n_0$ and $n_1$ the unique even and odd numbers in $\{n,n+1\}$, respectively. Write $n_0=2r_0$ and $n_1=2r_1+1$ for unique integers $r_0,r_1\geq 1$. In particular, we have $n=r_0+r_1$.

In this and the next sections, we consider
\begin{itemize}[label={\ding{109}}]
  \item for $\alpha=0,1$, a relevant representation $\Pi_\alpha$ of $\GL_{n_\alpha}(\dA_F)$ (Definition \ref{de:relevant}),

  \item a strong coefficient field $E\subseteq\dC$ of both $\Pi_0$ and $\Pi_1$ (Definition \ref{de:weak_field}).
\end{itemize}
Put $\Sigma^+_\mnm\coloneqq\Sigma^+_{\Pi_0}\cup\Sigma^+_{\Pi_1}$ (Notation \ref{no:satake}). We then have the homomorphism
\[
\phi_{\Pi_\alpha}\colon\dT^{\Sigma^+_\mnm}_{n_\alpha}\to O_E
\]
for $\alpha=0,1$. For $\alpha=0,1$ and every prime $\lambda$ of $E$, we have a continuous homomorphism
\[
\rho_{\Pi_\alpha,\lambda}\colon\Gamma_F\to\GL_{n_\alpha}(E_\lambda)
\]
from Proposition \ref{pr:galois}(2) and Definition \ref{de:weak_field}, such that $\rho_{\Pi_\alpha,\lambda}^\tc$ and $\rho_{\Pi_\alpha,\lambda}^\vee(1-n_\alpha)$ are conjugate.

\begin{assumption}\label{as:first_irreducible}
For $\alpha=0,1$, the Galois representation $\rho_{\Pi_\alpha,\lambda}$ is residually absolutely irreducible.
\end{assumption}

\subsection{First explicit reciprocity law}
\label{ss:first_reciprocity}

We start by choosing
\begin{itemize}[label={\ding{109}}]
  \item a prime $\lambda$ of $E$, whose underlying rational prime $\ell$ satisfies $\Sigma^+_\mnm\cap\Sigma^+_\ell=\emptyset$, $\ell\geq 2(n_0+1)$, and that $\ell$ is unramified in $F$,

  \item a positive integer $m$,

  \item a (possibly empty) finite set $\Sigma^+_{\lr,\rI}$ of nonarchimedean places of $F^+$ that are inert in $F$,\footnote{Here, the subscript ``lr'' stands for ``level-raising'', while the subscript ``I'' (Roman number one) stands for the ``first''. In the next subsection, we will have $\Sigma^+_{\lr,\r{II}}$ for the second reciprocity law.} strongly disjoint from $\Sigma^+_\mnm$ (Definition \ref{de:strongly_disjoint}), satisfying $\ell\nmid\|v\|(\|v\|^2-1)$ for $v\in\Sigma^+_{\lr,\rI}$,

  \item a finite set $\Sigma^+_\rI$ of nonarchimedean places of $F^+$ containing $\Sigma^+_\mnm\cup\Sigma^+_{\lr,\rI}$,

  \item a standard definite hermitian space $\rV^\circ_n$ of rank $n$ over $F$, together with a self-dual $\prod_{v\not\in\Sigma^+_\infty\cup\Sigma^+_\mnm\cup\Sigma^+_{\lr,\rI}}O_{F_v}$-lattice $\Lambda^\circ_n$ in $\rV^\circ_n\otimes_F\dA_F^{\Sigma^+_\infty\cup\Sigma^+_\mnm\cup\Sigma^+_{\lr,\rI}}$ (and put $\rV^\circ_{n+1}\coloneqq(\rV^\circ_n)_\sharp$ and $\Lambda^\circ_{n+1}\coloneqq(\Lambda^\circ_n)_\sharp$), satisfying that the hermitian space $(\rV^\circ_{n_0})_v$ is not split for $v\in\Sigma^+_{\lr,\rI}$,

  \item objects $\rK^\circ_n\in\fK(\rV^\circ_n)$ and $(\rK^\circ_\sp,\rK^\circ_{n+1})\in\fK(\rV^\circ_n)_\sp$ of the forms
      \begin{align*}
      \rK^\circ_n&=\prod_{v\in\Sigma^+_\mnm\cup\Sigma^+_{\lr,\rI}}(\rK^\circ_n)_v
      \times\prod_{v\not\in\Sigma^+_\infty\cup\Sigma^+_\mnm\cup\Sigma^+_{\lr,\rI}}\rU(\Lambda^\circ_n)(O_{F^+_v}), \\
      \rK^\circ_\sp&=\prod_{v\in\Sigma^+_\mnm\cup\Sigma^+_{\lr,\rI}}(\rK^\circ_\sp)_v
      \times\prod_{v\not\in\Sigma^+_\infty\cup\Sigma^+_\mnm\cup\Sigma^+_{\lr,\rI}}\rU(\Lambda^\circ_n)(O_{F^+_v}), \\
      \rK^\circ_{n+1}&=\prod_{v\in\Sigma^+_\mnm\cup\Sigma^+_{\lr,\rI}}(\rK^\circ_{n+1})_v
      \times\prod_{v\not\in\Sigma^+_\infty\cup\Sigma^+_\mnm\cup\Sigma^+_{\lr,\rI}}\rU(\Lambda^\circ_{n+1})(O_{F^+_v}),
      \end{align*}
      satisfying
      \begin{itemize}
        \item $(\rK^\circ_\sp)_v=(\rK^\circ_n)_v$ for $v\in\Sigma^+_\mnm$,

        \item $(\rK^\circ_\sp)_v\subseteq(\rK^\circ_n)_v$ for $v\in\Sigma^+_{\lr,\rI}$, and

        \item $(\rK^\circ_{n_0})_v$ is a transferable open compact subgroup (Definition \ref{de:transferable}) of $\rU(\rV^\circ_{n_0})(F^+_v)$ for $v\in\Sigma^+_\mnm$ and is a special maximal subgroup of $\rU(\rV^\circ_{n_0})(F^+_v)$ for $v\in\Sigma^+_{\lr,\rI}$,
      \end{itemize}

  \item a special inert prime (Definition \ref{de:special_inert}) $\fp$ of $F^+$ (with the underlying rational prime $p$) satisfying
      \begin{description}
        \item[(PI1)] $\Sigma^+_\rI$ does not contain $p$-adic places;

        \item[(PI2)] $\ell$ does not divide $p(p^2-1)$;

        \item[(PI3)] there exists a CM type $\Phi$ containing $\tau_\infty$ as in the initial setup of \S\ref{ss:ns} satisfying $\dQ_p^\Phi=\dQ_{p^2}$;

        \item[(PI4)] $P_{\balpha(\Pi_{0,\fp})}\modulo\lambda^m$ is level-raising special at $\fp$ (Definition \ref{de:satake_condition});

            $P_{\balpha(\Pi_{1,\fp})}\modulo\lambda$ is Tate generic at $\fp$ (Definition \ref{de:satake_condition});

        \item[(PI5)] $P_{\balpha(\Pi_{\alpha,\fp})}\modulo\lambda$ is intertwining generic at $\fp$ (Definition \ref{de:satake_condition}) for $\alpha=0,1$;

        \item[(PI6)] the natural map
            \[
            \frac{(O_E/\lambda^m)[\Sh(\rV^\circ_{n_\alpha},\rK^\circ_{n_\alpha})]}
            {\dT^{\Sigma^+_\rI\cup\Sigma^+_p}_{n_\alpha}\cap\Ker\phi_{\Pi_\alpha}}
            \to\frac{(O_E/\lambda^m)[\Sh(\rV^\circ_{n_\alpha},\rK^\circ_{n_\alpha})]}
            {\dT^{\Sigma^+_\rI}_{n_\alpha}\cap\Ker\phi_{\Pi_\alpha}}
            \]
            is an isomorphism of \emph{nontrivial} $O_E/\lambda^m$-modules for $\alpha=0,1$;

        \item[(PI7)] $P_{\balpha(\Pi_{0,\fp})\otimes\balpha(\Pi_{1,\fp})}\modulo\lambda^m$ is level-raising special at $\fp$ (Definition \ref{de:satake_condition});
      \end{description}
      (So we can and will apply the setup in \S\ref{ss:ns_functoriality} to the datum $(\rV^\circ_n,\{\Lambda^\circ_{n,\fq}\}\res_{\fq\mid p})$.)

  \item remaining data in \S\ref{ss:ns_initial} with $\dQ_p^\Phi=\dQ_{p^2}$; and

  \item data as in Notation \ref{no:ns_uniformization}, which in particular give open compact subgroups $\rK^\bullet_{n,p}$ and $\rK^\bullet_{n+1,p}$.
\end{itemize}

Put $\rK^{p\circ}_\sp\coloneqq(\rK_\sp^\circ)^p$ and $\rK^\bullet_\sp\coloneqq\rK^{p\circ}_\sp\times\rK^\bullet_{n_0,p}$; put $\rK^{p\circ}_{n_\alpha}\coloneqq(\rK_{n_\alpha}^\circ)^p$ and $\rK^\bullet_{n_\alpha}\coloneqq\rK^{p\circ}_{n_\alpha}\times\rK^\bullet_{n_\alpha,p}$ for $\alpha=0,1$. As in \S\ref{ss:ns_reciprocity}, we put $\rX^?_{n_\alpha}\coloneqq\rX^?_\fp(\rV^\circ_{n_\alpha},\rK^{p\circ}_{n_\alpha})$ for meaningful triples $(\rX,?,\alpha)\in\{\bM,\rM,\rB,\rS\}\times\{\;,\eta,\circ,\bullet,\dag\}\times\{0,1\}$. For $\alpha=0,1$, let $(\pres{\alpha}\rE^{p,q}_s,\pres{\alpha}\rd^{p,q}_s)$ be the weight spectral sequence abutting to the cohomology $\rH^{p+q}_\fT(\ol\rM_{n_\alpha},\rR\Psi O_\lambda(r_\alpha))$ from \S\ref{ss:ns_weight}.

\begin{notation}\label{no:ideal}
We introduce the following ideals of $\dT^{\Sigma^+_\rI\cup\Sigma^+_p}_{n_\alpha}$, for $\alpha=0,1$
\begin{align*}
\begin{dcases}
\fm_\alpha\coloneqq\dT^{\Sigma^+_\rI\cup\Sigma^+_p}_{n_\alpha}\cap\Ker\(\dT^{\Sigma^+}_{n_\alpha}\xrightarrow{\phi_{\Pi_\alpha}}O_E\to O_E/\lambda\),\\
\fn_\alpha\coloneqq\dT^{\Sigma^+_\rI\cup\Sigma^+_p}_{n_\alpha}\cap\Ker\(\dT^{\Sigma^+}_{n_\alpha}\xrightarrow{\phi_{\Pi_\alpha}}O_E\to O_E/\lambda^m\).
\end{dcases}
\end{align*}
\end{notation}

We then introduce the following assumptions.

\begin{assumption}\label{as:first_minimal}
Under Assumption \ref{as:first_irreducible}, $\bar\rho_{\Pi_0,\lambda,+}$ (Remark \ref{re:single_irreducible}) is rigid for $(\Sigma^+_\mnm,\Sigma^+_{\lr,\rI})$ (Definition \ref{de:rigid}); and $\bar\rho_{\Pi_0,\lambda}\res_{\Gal(\ol{F}/F(\zeta_\ell))}$ is absolutely irreducible.
\end{assumption}

\begin{assumption}\label{as:first_vanishing}
For $\alpha=0,1$, we have $\rH^i_\fT(\ol\rM_{n_\alpha},\rR\Psi O_\lambda)_{\fm_\alpha}=0$ for $i\neq n_\alpha-1$, and that $\rH^{n_\alpha-1}_\fT(\ol\rM_{n_\alpha},\rR\Psi O_\lambda)_{\fm_\alpha}$ is a finite free $O_\lambda$-module.
\end{assumption}

\begin{assumption}\label{as:first_generic}
The composite homomorphism $\dT^{\Sigma^+_\mnm}_{n_0}\xrightarrow{\phi_{\Pi_0}}O_E\to O_E/\lambda$ is cohomologically generic (Definition \ref{de:generic}).
\end{assumption}

Now we apply constructions in \S\ref{ss:ns_reciprocity}, evaluating on the object $(\rK^{p\circ}_n,\rK^{p\circ}_{n+1})$ of $\fK(\rV^\circ_n)^p\times\fK(\rV^\circ_{n+1})^p$. In particular, we have the blow-up morphism $\sigma\colon\bQ\to\bP$ from Notation \ref{no:ns_product}, and the localized spectral sequence $(\dE^{p,q}_{s,(\fm_0,\fm_1)},\rd^{p,q}_{s,(\fm_0,\fm_1)})$ from \eqref{eq:ns_product_weight}.

\begin{lem}\label{le:first_weight}
Assume Assumptions \ref{as:first_irreducible}, \ref{as:first_minimal}, \ref{as:first_vanishing}, \ref{as:first_generic} and Hypothesis \ref{hy:unitary_cohomology} for both $n$ and $n+1$. Then
\begin{enumerate}
  \item For $(?_0,?_1)\in\{\circ,\bullet,\dag\}^2$ and $i\in\dZ$, we have a canonical isomorphism
     \[
     \rH^i_\fT(\ol\rP^{?_0,?_1},O_\lambda)_{(\fm_0,\fm_1)}
     \simeq\bigoplus_{i_0+i_1=i}\rH^{i_0}_\fT(\ol\rM^{?_0}_{n_0},O_\lambda)_{\fm_0}
     \otimes_{O_\lambda}\rH^{i_1}_\fT(\ol\rM^{?_1}_{n_1},O_\lambda)_{\fm_1}
     \]
     in $\Mod(\Gal(\ol\dF_p/\dF_{p^2}),O_\lambda)_\free$.

  \item We have $\dE^{p,q}_{2,(\fm_0,\fm_1)}=0$ if $(p,q)\not\in\{(-1,2n),(0,2n-1),(1,2n-2)\}$, and canonical isomorphisms
     \[
     \begin{dcases}
     \dE^{-1,2n}_{2,(\fm_0,\fm_1)}\simeq\pres{0}\rE^{-1,2r_0}_{2,\fm_0}\otimes_{O_\lambda}\pres{1}\rE^{0,2r_1}_{2,\fm_1}, \\
     \dE^{0,2n-1}_{2,(\fm_0,\fm_1)}\simeq\pres{0}\rE^{0,2r_0-1}_{2,\fm_0}\otimes_{O_\lambda}\pres{1}\rE^{0,2r_1}_{2,\fm_1}, \\
     \dE^{1,2n-2}_{2,(\fm_0,\fm_1)}\simeq\pres{0}\rE^{1,2r_0-2}_{2,\fm_0}\otimes_{O_\lambda}\pres{1}\rE^{0,2r_1}_{2,\fm_1},
     \end{dcases}
     \]
     in $\Mod(\Gal(\ol\dF_p/\dF_{p^2}),O_\lambda)_\free$.

  \item If $\dE^{i,2n-1-i}_{2,(\fm_0,\fm_1)}(-1)$ has a nontrivial subquotient on which $\Gal(\ol\dF_p/\dF_{p^2})$ acts trivially, then $i=1$.

  \item For $(?_0,?_1)\in\{\circ,\bullet,\dag\}^2$ and $i\in\dZ$, both $\rH^{2i}_\fT(\ol\rP^{?_0,?_1},O_\lambda(i))_{(\fm_0,\fm_1)}$ and $\rH^{2i}_\fT(\ol\rQ^{?_0,?_1},O_\lambda(i))_{(\fm_0,\fm_1)}$ are weakly semisimple.

  \item We have $\rH^i_\fT(\ol\rQ,\rR\Psi O_\lambda)_{(\fm_0,\fm_1)}=0$ for $i\neq 2n-1$.

  \item The canonical map $\rH^i_{\fT,c}(\ol\rQ^{(c)},O_\lambda)_{(\fm_0,\fm_1)}\to\rH^i_\fT(\ol\rQ^{(c)},O_\lambda)_{(\fm_0,\fm_1)}$ is an isomorphism for every integers $c$ and $i$.
\end{enumerate}
\end{lem}

\begin{proof}
For (1), by Lemma \ref{le:ns_cohomology}, Lemma \ref{le:tate_weight}(2), Theorem \ref{th:raising}(1), we know that $\rH^{i_\alpha}_\fT(\ol\rM^{?_\alpha}_{n_\alpha},O_\lambda)_{\fm_\alpha}$ is a finitely generated free $O_\lambda$-module for $\alpha=0,1$ and every $i_\alpha\in\dZ$. Thus, (1) follows from Lemma \ref{le:single_compact} and the K\"{u}nneth formula.

For (2), we first show that $\dE^{p,q}_{s,(\fm_0,\fm_1)}$ degenerates at the second page. By (1), Lemma \ref{le:ns_product_weight}(2), Lemma \ref{le:ns_cohomology}, and Lemma \ref{le:tate_vanishing}, the composition of $\rd^{-2,q}_{1,(\fm_0,\fm_1)}$ and the natural projection
\[
\dE^{-1,q}_{1,(\fm_0,\fm_1)}\to\rH^{q-2}_\fT(\ol\rQ^{\dag,\dag},O_\lambda(n-1))\bigoplus\rH^{q-2}_\fT(\ol\rQ^{\dag,\circ},O_\lambda(n-1))
\]
is injective for every $q\in\dZ$. Thus, $\rd^{-2,q}_{1,(\fm_0,\fm_1)}$ is injective, which implies $\dE^{-2,q}_{2,(\fm_0,\fm_1)}=0$ for every $q\in\dZ$. By a dual argument, we have $\dE^{2,q}_{2,(\fm_0,\fm_1)}=0$ for every $q\in\dZ$ as well. For the degeneration, it suffices to show that $\rd^{-1,q}_{1,(\fm_0,\fm_1)}$ is injective and $\rd^{0,q}_{1,(\fm_0,\fm_1)}$ is surjective for $q$ odd. By Lemma \ref{le:ns_product_weight}(2), Lemma \ref{le:ns_cohomology}, and Lemma \ref{le:tate_weight}(1), we have $\rH^{q-2}_\fT(\ol\rQ^{(1)},O_\lambda(n-1))=\rH^{q-2}_\fT(\ol\rQ^{\bullet,\dag},O_\lambda(n-1))$ for $q$ odd, which easily implies the injectivity of $\rd^{-1,q}_{1,(\fm_0,\fm_1)}$. By a dual argument, $\rd^{0,q}_{1,(\fm_0,\fm_1)}$ is surjective for $q$ odd.

Now for every $q\in\dZ$, the morphism $\sigma$ induces a map
\[
\sigma^*_1\colon\bigoplus_{q_0+q_1=q}\pres{0}\rE^{*,q_0}_{1,\fm_0}\otimes_{O_\lambda}\pres{1}\rE^{*,q_1}_{1,\fm_1}
\to\dE^{*,q}_{1,(\fm_0,\fm_1)}
\]
of complexes of $O_\lambda[\Gal(\ol\dF/\dF_{p^2})]$-modules, hence a map
\[
\sigma^*_2\colon\bigoplus_{p_0+p_1=p}\bigoplus_{q_0+q_1=q}\pres{0}\rE^{p_0,q_0}_{2,\fm_0}\otimes_{O_\lambda}\pres{1}\rE^{p_1,q_1}_{2,\fm_1}
\to\dE^{p,q}_{2,(\fm_0,\fm_1)}
\]
of $O_\lambda[\Gal(\ol\dF/\dF_{p^2})]$-modules for $(p,q)\in\dZ^2$. By Lemma \ref{le:tate_weight} and Theorem \ref{th:raising}(2), to show (2), it suffices to show that $\sigma^*_2$ is an isomorphism, or the natural map
\[
\bigoplus_{i_0+i_1=i}\rH^{i_0}_\fT(\ol\rM_{n_0},\rR\Psi O_\lambda(r_0))_{\fm_0}
\otimes_{O_\lambda}\rH^{i_1}_\fT(\ol\rM_{n_1},\rR\Psi O_\lambda(r_1))_{\fm_1}
\to\rH^{i}_\fT(\ol\rQ,\rR\Psi O_\lambda(n))_{(\fm_0,\fm_1)}
\]
induced by $\sigma$ is an isomorphism for every $i\in\dZ$. By Lemma \ref{le:ns_pbc} and Lemma \ref{le:ns_product_pbc}, the above map is identified with
\[
\bigoplus_{i_0+i_1=i}\rH^{i_0}_\fT(\bM^\eta_{n_0}\otimes_{\dQ_{p^2}}\ol\dQ_p,O_\lambda(r_0))_{\fm_0}
\otimes_{O_\lambda}\rH^{i_1}_\fT(\bM^\eta_{n_1}\otimes_{\dQ_{p^2}}\ol\dQ_p,O_\lambda(r_1))_{\fm_1}
\to\rH^{i}_\fT(\bQ^\eta\otimes_{\dQ_{p^2}}\ol\dQ_p,O_\lambda(n))_{(\fm_0,\fm_1)},
\]
which is an isomorphism by Lemma \ref{le:single_compact}, and the K\"{u}nneth formula. Thus, (2) follows.

For (3), let $\{\alpha_{0,1}^{\pm1},\dots,\alpha_{0,r_0}^{\pm1}\}$ and $\{\alpha_{1,1}^{\pm1},\dots,\alpha_{1,r_1}^{\pm1},1\}$ be the roots of $P_{\balpha(\Pi_{0,\fp})}\modulo\lambda$ and $P_{\balpha(\Pi_{1,\fp})}\modulo\lambda$ in a finite extension of $O_\lambda/\lambda$, respectively. By (PI4), we may assume $\alpha_{0,r_0}=p$. By (2), Theorem \ref{th:tate}(1), and Theorem \ref{th:raising}(3), the generalized Frobenius eigenvalues of the $(O_\lambda/\lambda)[\Gal(\ol\dF_p/\dF_{p^2})]$-modules $\dE^{-1,2n}_{2,(\fm_0,\fm_1)}(-1)\otimes_{O_\lambda}O_\lambda/\lambda$ and $\dE^{0,2n-1}_{2,(\fm_0,\fm_1)}(-1)\otimes_{O_\lambda}O_\lambda/\lambda$ are contained in $\{p^{-2}\alpha_{1,1}^{\pm1},\dots,p^{-2}\alpha_{1,r_1}^{\pm1},p^{-2}\}$ and $\{p^{-1}\alpha_{0,1}^{\pm1}\alpha_{1,1}^{\pm1},\dots,p^{-1}\alpha_{0,r_0-1}^{\pm1}\alpha_{1,r_1}^{\pm1}\}\cup
\{p^{-1}\alpha_{0,1}^{\pm1},\dots,p^{-1}\alpha_{0,r_0-1}^{\pm1}\}$, respectively. By (PI2), we have $p^2\neq 1$ in $O_\lambda/\lambda$. By (PI7), we have $\alpha_{1,i_1}\not\in\{p^2,p^{-2}\}$ for $1\leq i_1\leq r_1$, which implies $1\not\in\{p^{-2}\alpha_{1,1}^{\pm1},\dots,p^{-2}\alpha_{1,r_1}^{\pm1},p^{-2}\}$. Again by (PI7), we have $\alpha_{0,i_0}\alpha_{1,i_1}\not\in\{p,p^{-1}\}$ for $1\leq i_0< r_0$ and $1\leq i_1\leq r_1$, which implies $1\not\in\{p^{-1}\alpha_{0,1}^{\pm1}\alpha_{1,1}^{\pm1},\dots,p^{-1}\alpha_{0,r_0-1}^{\pm1}\alpha_{1,r_1}^{\pm1}\}$. By (PI4), we have $\alpha_{0,i_0}\not\in\{p,p^{-1}\}$ for $1\leq i_0< r_0$, which implies $1\not\in\{p^{-1}\alpha_{0,1}^{\pm1},\dots,p^{-1}\alpha_{0,r_0-1}^{\pm1}\}$. Thus, (3) follows.

For (4), by Lemma \ref{le:ns_product_weight} (3--5) and Lemma \ref{le:weakly_semisimple1}(1), it suffices to show that $\rH^{2i}_\fT(\ol\rP^{?_0,?_1},O_\lambda(i))_{(\fm_0,\fm_1)}$ is weakly semisimple. By (1) and Lemma \ref{le:tate_weight}(1), it suffices to show that $\rH^{2i_0}_\fT(\ol\rM^{?_0}_{n_0},O_\lambda(i_0))_{\fm_0}\otimes_{O_\lambda}\rH^{2i_1}_\fT(\ol\rM^{?_1}_{n_1},O_\lambda(i_1))_{\fm_1}$ is weakly semisimple for $i_0,i_1\in\dZ$. By Lemma \ref{le:ns_cohomology}, the action of $\Gal(\ol\dF_p/\dF_{p^2})$ on $\rH^{2i_\alpha}_\fT(\ol\rM^?_{n_\alpha},O_\lambda(i_\alpha))_{\fm_\alpha}$ is trivial for $\alpha=0,1$, $?=\circ,\dag$, and every $i_\alpha\in\dZ$. On the other hand, it is a consequence of Theorem \ref{th:raising}(2) (for $i_0$) and Lemma \ref{le:tate_weight}(3) (for $i_1$) that the action of $\Gal(\ol\dF_p/\dF_{p^2})$ on $\rH^{2i_\alpha}_\fT(\ol\rM^\bullet_{n_\alpha},O_\lambda(i_\alpha))_{\fm_\alpha}$ is trivial if $i_0\not\in\{r_0-1,r_0\}$ or $i_1\neq r_1$. By Proposition \ref{pr:raising_weight}(1,2) and Theorem \ref{th:raising}(1), the actions of $\Gal(\ol\dF_p/\dF_{p^2})$ on both $\rH^{2r_0-2}_\fT(\ol\rM^\bullet_{n_0},O_\lambda(r_0-1))_{\fm_0}$ and $\rH^{2r_0}_\fT(\ol\rM^\bullet_{n_0},O_\lambda(r_0))_{\fm_0}$ are also trivial. Thus, by Lemma \ref{le:weakly_semisimple1}(1), it remains to show that $\rH^{2r_1}_\fT(\ol\rM^\bullet_{n_1},O_\lambda(r_1))_{\fm_1}$ is weakly semisimple, which follows from Theorem \ref{th:tate}(2) as it is isomorphic to the direct sum of $\pres{1}\rE^{0,2r_1}_{2,\fm_1}$ and $\rH^{2r_1}_\fT(\ol\rM^\dag_{n_1},O_\lambda(r_1))_{\fm_1}$.

Part (5) is a direct consequence of (2).

Part (6) follows from (1), Lemma \ref{le:single_compact}, and Lemma \ref{le:ns_product_weight}(3--5).
\end{proof}

\begin{remark}
In fact, Lemma \ref{le:first_weight}(5) holds under only Assumption \ref{as:first_vanishing}; and Lemma \ref{le:first_weight}(6) holds under only Assumption \ref{as:first_irreducible}.
\end{remark}

Lemma \ref{le:first_weight}(5) induces a coboundary map
\begin{align*}
\AJ_\bQ\colon\rZ^n_\fT(\bQ^\eta)\to\rH^1(\dQ_{p^2},\rH^{2n-1}_\fT(\ol\rQ,\rR\Psi O_\lambda(n))_{(\fm_0,\fm_1)}).
\end{align*}
We also recall the singular quotient map
\begin{align}\label{eq:partial}
\partial\colon\rH^1(\dQ_{p^2},\rH^{2n-1}_\fT(\ol\rQ,\rR\Psi O_\lambda(n))_{(\fm_0,\fm_1)})\to
\rH^1_\sing(\dQ_{p^2},\rH^{2n-1}_\fT(\ol\rQ,\rR\Psi O_\lambda(n))_{(\fm_0,\fm_1)})
\end{align}
from Definition \ref{de:local_finite}.

By our choice of $\rK^\circ_n$ and $(\rK^\circ_\sp,\rK^\circ_{n+1})$, we obtain a morphism
\[
\bM_\fp(\rV^\circ_n,\rK^\circ_\sp)\to\bP
\]
which is finite. Denote by $\bP_\sp$ the corresponding cycle; and let $\bQ_\sp$ be the strict transform of $\bP_\sp$ under $\sigma$, which is a $\bT_\fp$-invariant cycle of $\bQ$. Our main goal is to compute $\partial\AJ_\bQ(\bQ^\eta_\sp)$ in $\rH^1_\sing(\dQ_{p^2},\rH^{2n-1}_\fT(\ol\rQ,\rR\Psi O_\lambda(n))/(\fn_0,\fn_1))$. The cycle $\rQ_\sp$ gives rise to a class $\cl(\rQ_\sp)\in C^n(\rQ,L)$, where $C^n(\rQ,L)$ is the target of the map $\Delta^n$ \eqref{eq:ns_potential}.

\begin{proposition}\label{pr:first_potential}
Assume Assumptions \ref{as:first_irreducible}, \ref{as:first_minimal}, \ref{as:first_vanishing}, \ref{as:first_generic}, and Hypothesis \ref{hy:unitary_cohomology} for both $n$ and $n+1$. There is a canonical isomorphism
\[
\rH^1_\sing(\dQ_{p^2},\rH^{2n-1}_\fT(\ol\rQ,\rR\Psi O_\lambda(n))_{(\fm_0,\fm_1)})
\simeq\coker\Delta^n_{(\fm_0,\fm_1)}
\]
under which $\partial\AJ_\bQ(\bQ^\eta_\sp)$ coincides with the image of $\cl(\rQ_\sp)$ in $\coker\Delta^n_{(\fm_0,\fm_1)}$.
\end{proposition}

\begin{proof}
By \cite{Liu2}*{Theorem~2.16 and Theorem~2.18},\footnote{Although it is assumed that the underlying strictly semistable scheme $X$ is proper over the base in \cite{Liu2}, the proof of relevant results works without change in our case even when $\bQ$ is not proper in view of Lemma \ref{le:first_weight}(6).} it suffices to show that $O_\lambda$ is a very nice coefficient ring for $\dE^{p,q}_{s,(\fm_0,\fm_1)}$ in the sense of \cite{Liu2}*{Definition~2.15}. In fact, in \cite{Liu2}*{Definition~2.15}, (N1) is satisfied due to Lemma \ref{le:first_weight}(2); (N2) is satisfied due to Lemma \ref{le:first_weight}(3); and (N3) is satisfied due to Lemma \ref{le:first_weight}(4) and Lemma \ref{le:weakly_semisimple1}(2).

The proposition is proved.
\end{proof}

By Construction \ref{cs:ns_nabla_product} and Remark \ref{re:ns_nabla_product}, we have a map
\begin{align*}
\nabla\colon C^n(\rQ,O_\lambda)\to
O_\lambda[\Sh(\rV^\circ_{n_0},\rK^\circ_{n_0})]\otimes_{O_\lambda}O_\lambda[\Sh(\rV^\circ_{n_1},\rK^\circ_{n_1})].
\end{align*}

\begin{theorem}[First explicit reciprocity law]\label{th:first}
Assume Assumptions \ref{as:first_irreducible}, \ref{as:first_minimal}, \ref{as:first_vanishing}, \ref{as:first_generic}, and Hypothesis \ref{hy:unitary_cohomology} for both $n$ and $n+1$.
\begin{enumerate}
  \item The image of the composite map $\nabla_{(\fm_0,\fm_1)}\circ\Delta^n_{(\fm_0,\fm_1)}$ is contained in $\fn_0.O_\lambda[\Sh(\rV^\circ_{n_0},\rK^\circ_{n_0})]_{\fm_0}\otimes_{O_\lambda}O_\lambda[\Sh(\rV^\circ_{n_1},\rK^\circ_{n_1})]_{\fm_1}$.

  \item In view of (1), the induced map
      \[
      \nabla_{\fm_1/\fn_0}\colon\coker\Delta^n_{\fm_1}/\fn_0\to
      O_\lambda[\Sh(\rV^\circ_{n_0},\rK^\circ_{n_0})]/\fn_0\otimes_{O_\lambda}O_\lambda[\Sh(\rV^\circ_{n_1},\rK^\circ_{n_1})]_{\fm_1}
      \]
      is an isomorphism.

  \item Under the natural pairing
      \[
      \resizebox{\hsize}{!}{
      \xymatrix{
      O_\lambda[\Sh(\rV^\circ_{n_0},\rK^\circ_{n_0})]/\fn_0\otimes_{O_\lambda}O_\lambda[\Sh(\rV^\circ_{n_1},\rK^\circ_{n_1})]_{\fm_1}\times
      (O_\lambda/\lambda^m)[\Sh(\rV^\circ_{n_0},\rK^\circ_{n_0})][\fn_0]
      \otimes_{O_\lambda}O_\lambda[\Sh(\rV^\circ_{n_1},\rK^\circ_{n_1})]_{\fm_1}
      \to O_\lambda/\lambda^m
      }
      }
      \]
      obtained by taking inner product, the pairing of $\nabla_{/(\fn_0,\fn_1)}(\partial\AJ_\bQ(\bQ^\eta_\graph))$ and every function $f\in (O_\lambda/\lambda^m)[\Sh(\rV^\circ_{n_0},\rK^\circ_{n_0})][\fn_0]\otimes_{O_\lambda}
      (O_\lambda/\lambda^m)[\Sh(\rV^\circ_{n_1},\rK^\circ_{n_1})][\fn_1]$ is equal to
      \[
      (p+1)\cdot\phi_{\Pi_0}(\tI^\circ_{n_0,\fp})\cdot\phi_{\Pi_1}(\tT^\circ_{n_1,\fp})\cdot
      \sum_{s\in\Sh(\rV^\circ_n,\rK^\circ_\sp)}f(s,\sh^\circ_\uparrow(s)).
      \]
      Here, we regard $\partial\AJ_\bQ(\bQ^\eta_\sp)$ as an element in $\coker\Delta^n_{(\fm_0,\fm_1)}$ (hence in $\coker\Delta^n_{\fm_1}/\fn_0$) via the canonical isomorphism in Proposition \ref{pr:first_potential}.
\end{enumerate}
\end{theorem}

\begin{proof}
We first consider (1). By Lemma \ref{le:ns_product_weight}(3,4), we have
\begin{align*}
&\rH^{2(n-1)}_\fT(\ol\rQ^{(0)},O_\lambda(n-1))_{(\fm_0,\fm_1)}
=\bigoplus_{(?_0,?_1)\in\{\circ,\bullet\}^2}\sigma^*\rH^{2(n-1)}_\fT(\ol\rP^{?_0,?_1},O_\lambda(n-1))_{(\fm_0,\fm_1)} \\
&\qquad\qquad\bigoplus(\delta^{\dag,\dag}_{\circ,\circ})_!\sigma^*\rH^{2(n-2)}_\fT(\ol\rP^{\dag,\dag},O_\lambda(n-2))_{(\fm_0,\fm_1)}
\bigoplus(\delta^{\dag,\dag}_{\bullet,\bullet})_!\sigma^*\rH^{2(n-2)}_\fT(\ol\rP^{\dag,\dag},O_\lambda(n-2))_{(\fm_0,\fm_1)}.
\end{align*}
Thus, it suffices to show that
\begin{enumerate}
  \item[(1a)] The image of $\sigma^*\rH^{2(n-1)}_\fT(\ol\rP^{\circ,\bullet},O_\lambda(n-1))_{(\fm_0,\fm_1)}
      \bigoplus\sigma^*\rH^{2(n-1)}_\fT(\ol\rP^{\bullet,\bullet},O_\lambda(n-1))_{(\fm_0,\fm_1)}$ under the map $(\nabla\circ\delta_{1!}\circ\delta_0^*)_{(\fm_0,\fm_1)}$ is contained in $\fn_0.O_\lambda[\Sh(\rV^\circ_{n_0},\rK^\circ_{n_0})]_{\fm_0}
      \otimes_{O_\lambda}O_\lambda[\Sh(\rV^\circ_{n_1},\rK^\circ_{n_1})]_{\fm_1}$.

  \item[(1b)] The image of $\sigma^*\rH^{2(n-1)}_\fT(\ol\rP^{\circ,\circ},O_\lambda(n-1))_{(\fm_0,\fm_1)}
      \bigoplus\sigma^*\rH^{2(n-1)}_\fT(\ol\rP^{\bullet,\circ},O_\lambda(n-1))_{(\fm_0,\fm_1)}$ under the map $(\nabla\circ\delta_{1!}\circ\delta_0^*)_{(\fm_0,\fm_1)}$ is zero.

  \item[(1c)] The image of $(\delta^{\dag,\dag}_{\circ,\circ})_!\sigma^*\rH^{2(n-2)}_\fT(\ol\rP^{\dag,\dag},O_\lambda(n-2))_{(\fm_0,\fm_1)}$ under the map $(\nabla\circ\delta_{1!}\circ\delta_0^*)_{(\fm_0,\fm_1)}$ is zero.

  \item[(1d)] The image of  $(\delta^{\dag,\dag}_{\bullet,\bullet})_!\sigma^*\rH^{2(n-2)}_\fT(\ol\rP^{\dag,\dag},O_\lambda(n-2))_{(\fm_0,\fm_1)}$ under the map $(\nabla\circ\delta_{1!}\circ\delta_0^*)_{(\fm_0,\fm_1)}$ is zero.
\end{enumerate}

For (1a), we have a commutative diagram
\begin{align*}
\resizebox{\hsize}{!}{
\xymatrix{
\rH^{2(n-1)}_\fT(\ol\rP^{\circ,\bullet},O_\lambda(n-1))_{(\fm_0,\fm_1)}\bigoplus \rH^{2(n-1)}_\fT(\ol\rP^{\bullet,\bullet},O_\lambda(n-1))_{(\fm_0,\fm_1)} \ar[r]\ar[d]_-{\sigma^*}
& \pres{0}\rE^{0,2r_0-2}_{1,\fm_0}\otimes_{O_\lambda}\rH^{2r_1}_\fT(\ol\rM^\bullet_{n_1},O_\lambda(r_1))_{\fm_1} \ar[d] \\
\rH^{2(n-1)}_\fT(\ol\rQ^{\circ,\bullet},O_\lambda(n-1))_{(\fm_0,\fm_1)}\bigoplus
\rH^{2(n-1)}_\fT(\ol\rQ^{\bullet,\bullet},O_\lambda(n-1))_{(\fm_0,\fm_1)} \ar[r]&
O_\lambda[\Sh(\rV^\circ_{n_0},\rK^\circ_{n_0})]_{\fm_0}
\otimes_{O_\lambda}O_\lambda[\Sh(\rV^\circ_{n_1},\rK^\circ_{n_1})]_{\fm_1}
}
}
\end{align*}
in which
\begin{itemize}[label={\ding{109}}]
  \item the upper horizontal arrow is the map
    \begin{align*}
    &\quad\rH^{2(n-1)}_\fT(\ol\rP^{\circ,\bullet},O_\lambda(n-1))_{(\fm_0,\fm_1)}\bigoplus \rH^{2(n-1)}_\fT(\ol\rP^{\bullet,\bullet},O_\lambda(n-1))_{(\fm_0,\fm_1)} \\
    &\to
    \rH^{2(r_0-1)}_\fT(\ol\rM^\circ_{n_0},O_\lambda(r_0-1))_{\fm_0}
    \otimes_{O_\lambda}\rH^{2r_1}_\fT(\ol\rM^\bullet_{n_1},O_\lambda(r_1))_{\fm_1} \\
    &\quad\quad\bigoplus
    \rH^{2(r_0-1)}_\fT(\ol\rM^\bullet_{n_0},O_\lambda(r_0-1))_{\fm_0}
    \otimes_{O_\lambda}\rH^{2r_1}_\fT(\ol\rM^\bullet_{n_1},O_\lambda(r_1))_{\fm_1} \\
    &=\pres{0}\rE^{0,2r_0-2}_{1,\fm_0}\otimes_{O_\lambda}\rH^{2r_1}_\fT(\ol\rM^\bullet_{n_1},O_\lambda(r_1))_{\fm_1}
    \end{align*}
    given by Lemma \ref{le:first_weight}(1) and the K\"{u}nneth formula;

  \item the right vertical arrow is
      \[
      (\nabla^0\circ\pres{0}\rd^{-1,2r_0}_1\circ\pres{0}\rd^{0,2r_0-2}_1(-1))_{\fm_0}\otimes (\tI^\circ_{n_1,\fp}\circ\inc_\dag^*+(p+1)^2\tT^{\circ\bullet}_{n_1,\fp}\circ\inc_\bullet^*)_{\fm_1};
      \]
      and

  \item the lower horizontal arrow is $(\nabla\circ\delta_{1!}\circ\delta_0^*)_{(\fm_0,\fm_1)}$.
\end{itemize}
For (1a), by Proposition \ref{pr:enumeration_even}(2) and (PI4), we have
\[
((p+1)\tR^\circ_{n_0,\fp}-\tI^\circ_{n_0,\fp}).O_\lambda[\Sh(\rV^\circ_{n_0},\rK^\circ_{n_0})]_{\fm_0}\subseteq
\fn_0.O_\lambda[\Sh(\rV^\circ_{n_0},\rK^\circ_{n_0})]_{\fm_0}.
\]
Thus, (1a) follows from Proposition \ref{pr:raising_weight}(4) and Lemma \ref{le:ns_product_weight}(3).

For (1b) and (1c), both images are actually contained in the sum of
\[
(\tI^\circ_{n_1,\fp}\circ\inc_{\circ,\dag}^*+(p+1)^2\tT^{\circ\bullet}_{n_1,\fp}
\circ\inc_{\circ,\bullet}^*)(\gamma^{\circ,\dag}_{\circ,\bullet})_!\rH^{2(n-1)}_\fT(\ol\rP^{\circ,\dag},O_\lambda(n-1))_{(\fm_0,\fm_1)}
\]
and
\[
(\tI^\circ_{n_1,\fp}\circ\inc_{\circ,\dag}^*+(p+1)^2\tT^{\circ\bullet}_{n_1,\fp}
\circ\inc_{\bullet,\bullet}^*)(\gamma^{\bullet,\dag}_{\bullet,\bullet})_!\rH^{2(n-1)}_\fT(\ol\rP^{\bullet,\dag},O_\lambda(n-1))_{(\fm_0,\fm_1)},
\]
which by Lemma \ref{le:first_weight}(1) coincide with
\[
\rH^{2r_0}_\fT(\ol\rM^\circ_{n_0},O_\lambda(r_0))_{\fm_0}\otimes_{O_\lambda}
\((\tI^\circ_{n_1,\fp}\circ\Inc_\dag^*+(p+1)^2\tT^{\circ\bullet}_{n_1,\fp}\circ\Inc_\bullet^*)\pres{1}\rd^{-1,2r_1}_1
\rH^{2(r_1-1)}_\fT(\ol\rM^\dag_{n_1},O_\lambda(r_1-1))_{\fm_1}\)
\]
and
\[
\rH^{2r_0}_\fT(\ol\rM^\bullet_{n_0},O_\lambda(r_0))_{\fm_0}\otimes_{O_\lambda}
\((\tI^\circ_{n_1,\fp}\circ\Inc_\dag^*+(p+1)^2\tT^{\circ\bullet}_{n_1,\fp}\circ\Inc_\bullet^*)\pres{1}\rd^{-1,2r_1}_1
\rH^{2(r_1-1)}_\fT(\ol\rM^\dag_{n_1},O_\lambda(r_1-1))_{\fm_1}\),
\]
respectively. However, they vanish by Lemma \ref{le:ns_weight_odd}(3). Thus, (1b) and (1c) follow.

For (1d), by \cite{Liu2}*{Lemma~2.4}, it follows from (1c). Thus, (1) is proved.

Now we consider (2). We claim that the map $\nabla_{(\fm_0,\fm_1)}$ (with domain $C^n(\rQ,O_\lambda)_{(\fm_0,\fm_1)}$) is surjective. In fact, consider the submodule
\[
\Ker\pres{0}\rd^{0,2r_0}_{1,\fm_0}\otimes_{O_\lambda}\Ker\pres{1}\rd^{0,2r_1}_{1,\fm_1}\subseteq \bigoplus_{(?_0,?_1)\in\{\circ,\bullet\}^2}\rH^{2(n-1)}_\fT(\ol\rP^{?_0,?_1},O_\lambda(n-1))_{(\fm_0,\fm_1)}
\]
in view of Lemma \ref{le:first_weight}(1). Then $\sigma^*\(\Ker\pres{0}\rd^{0,2r_0}_{1,\fm_0}\otimes_{O_\lambda}\Ker\pres{1}\rd^{0,2r_1}_{1,\fm_1}\)$ is contained in $C^n(\rQ,O_\lambda)_{(\fm_0,\fm_1)}$. On the other hand, the map $\nabla_{(\fm_0,\fm_1)}\circ\sigma^*$ (with domain $\Ker\pres{0}\rd^{0,2r_0}_{1,\fm_0}\otimes_{O_\lambda}\Ker\pres{1}\rd^{0,2r_1}_{1,\fm_1}$) coincides with $\nabla^0_{\fm_0}\otimes\nabla^1_{\fm_1}$, which is surjective by Proposition \ref{pr:raising_weight}(3) and Theorem \ref{th:tate}. The claim follows.

Thus, it remains to show that the domain and the target of $\nabla_{\fm_1/\fn_0}$ have the same cardinality. By Proposition \ref{pr:first_potential}, we have an isomorphism
\begin{align}\label{eq:first_1}
\coker\Delta^n_{\fm_1}/\fn_0=\coker\Delta^n_{(\fm_0,\fm_1)}/\fn_0\simeq
\rH^1_\sing(\dQ_{p^2},\rH^{2n-1}_\fT(\ol\rQ,\rR\Psi O_\lambda(n))_{(\fm_0,\fm_1)})/\fn_0
\end{align}
of $O_\lambda/\lambda^m$-modules. By Lemma \ref{le:first_weight}(2,3) and Theorem \ref{th:tate}(2), we have
\[
\resizebox{\hsize}{!}{
\xymatrix{
\rH^1_\sing(\dQ_{p^2},\rH^{2n-1}_\fT(\ol\rQ,\rR\Psi O_\lambda(n))_{(\fm_0,\fm_1)})
\simeq
\rH^1_\sing(\dQ_{p^2},\rH^{2r_0-1}_\fT(\ol\rM_{n_0},\rR\Psi O_\lambda(r_0))_{\fm_0})\otimes_{O_\lambda}
(\pres{1}\rE^{0,2r_1}_{2,\fm_1})^{\Gal(\ol\dF_p/\dF_{p^2})}.
}
}
\]
Then by Theorem \ref{th:tate}(3) and Theorem \ref{th:raising}(4), we have
\[
\eqref{eq:first_1}\simeq
O_\lambda[\Sh(\rV^\circ_{n_0},\rK^\circ_{n_0})]/\fn_0\otimes_{O_\lambda}O_\lambda[\Sh(\rV^\circ_{n_1},\rK^\circ_{n_1})]_{\fm_1}.
\]
Thus, (2) is proved.

Finally we consider (3). As $\rQ_\sp$ does not intersect with $\rQ^{\circ,\bullet}$, we have
\[
\nabla(\cl(\rQ_\sp))=\nabla(\cl(\rQ^\bullet_\graph))
\]
where $\cl(\rQ^\bullet_\graph)\in\rH^{2n}_\fT(\ol\rQ^{\bullet,\bullet},O_\lambda(n))$. Then by Construction \ref{cs:ns_nabla_product}, we have
\[
\nabla(\cl(\rQ_\sp))=\((p+1)(\tT^{\circ\bullet}_{n_0,\fp}\otimes\tI^\circ_{n_1,\fp})\circ\inc_{\bullet,\dag}^*+
(p+1)^3(\tT^{\circ\bullet}_{n_0,\fp}\otimes\tT^{\circ\bullet}_{n_1,\fp})\circ\inc_{\bullet,\bullet}^*\)(\cl(\rP^\bullet_\sp)).
\]
Applying Theorem \ref{th:ns_reciprocity}(3) to the object $(\rK^\circ_\sp,\rK^\circ_{n+1})\in\fK(\rV^\circ_n)_\sp$ followed by pushforward, we know that the pairing between $\nabla_{\fm_1/\fn_0}(\cl(\rQ_\sp))$ and any function
\[
f\in(O_\lambda/\lambda^m)[\Sh(\rV^\circ_{n_0},\rK^\circ_{n_0})][\fn_0]
\otimes_{O_\lambda}(O_\lambda/\lambda^m)[\Sh(\rV^\circ_{n_1},\rK^\circ_{n_1})][\fn_1]
\]
is given by the formula
\[
(p+1)\cdot\phi_{\Pi_0}(\tI^\circ_{n_0,\fp})\cdot\phi_{\Pi_1}(\tT^\circ_{n_1,\fp})\cdot
\sum_{s\in\Sh(\rV^\circ_n,\rK^\circ_\sp)}f(s,\sh^\circ_\uparrow(s))
\]
in view of (PI6). We then obtain (3) by Proposition \ref{pr:first_potential}.

The theorem is proved.
\end{proof}

We state a corollary for later application. We choose an indefinite uniformization datum as in Notation \ref{no:ns_uniformization_indefinite}, and put $\Sh'_{n_\alpha}\coloneqq\Sh(\rV'_{n_\alpha},\tj_{n_\alpha}\rK^{p\circ}_{n_\alpha}\rK'_{n_\alpha,p})$ for $\alpha=0,1$.

Assume Assumption \ref{as:first_irreducible} and Assumption \ref{as:first_vanishing}. By Lemma \ref{le:single_compact}, Lemma \ref{le:ns_pbc}, and the K\"{u}nneth formula, we have $\rH^i_\et((\Sh'_{n_0}\times_{\Spec{F}}\Sh'_{n_1})_{\ol{F}},O_\lambda)_{(\fm_0,\fm_1)}=0$ if $i\neq 2n-1$. In particular, we obtain the Abel--Jacobi map
\[
\AJ\colon\rZ^n(\Sh'_{n_0}\times_{\Spec{F}}\Sh'_{n_1})\to
\rH^1(F,\rH^{2n-1}_\et((\Sh'_{n_0}\times_{\Spec{F}}\Sh'_{n_1})_{\ol{F}},O_\lambda(n))/(\fn_0,\fn_1)).
\]
Let $\Sh'_\sp$ be the cycle given by the finite morphism $\Sh(\rV'_n,\tj_n\rK^{p\circ}_\sp\rK'_{n,p})\to\Sh'_n\times_{\Spec{F}}\Sh'_{n+1}$, which is an element in $\rZ^n(\Sh'_{n_0}\times_{\Spec{F}}\Sh'_{n_1})$.

\begin{corollary}\label{co:first}
Assume Assumptions \ref{as:first_irreducible}, \ref{as:first_minimal}, \ref{as:first_vanishing}, \ref{as:first_generic}, and Hypothesis \ref{hy:unitary_cohomology} for both $n$ and $n+1$. Then we have
\begin{align*}
&\quad\exp_\lambda\(\partial_\fp\loc_\fp\AJ(\Sh'_\sp),
\rH^1_\sing(F_\fp,\rH^{2n-1}_\et((\Sh'_{n_0}\times_{\Spec{F}}\Sh'_{n_1})_{\ol{F}},O_\lambda(n))/(\fn_0,\fn_1))\)\\
&=\exp_\lambda\(\CF_{\Sh(\rV^\circ_n,\rK^\circ_\sp)},
O_\lambda[\Sh(\rV^\circ_{n_0},\rK^\circ_{n_0})\times\Sh(\rV^\circ_{n_1},\rK^\circ_{n_1})]/(\fn_0,\fn_1)\)
\end{align*}
where $\exp_\lambda$ is introduced in Definition \ref{de:divisibility}. Here, we regard $\CF_{\Sh(\rV^\circ_n,\rK^\circ_\sp)}$ as the pushforward of the characteristic function along the map $\Sh(\rV^\circ_n,\rK^\circ_\sp)\to\Sh(\rV^\circ_n,\rK^\circ_n)\times\Sh(\rV^\circ_{n+1},\rK^\circ_{n+1})$.
\end{corollary}

\begin{proof}
Note that the isomorphism \eqref{eq:ns_moduli_scheme_shimura} induces a map
\[
\rH^{2n-1}_\et((\Sh'_{n_0}\times_{\Spec{F}}\Sh'_{n_1})_{\ol{F}},O_\lambda(n))_{(\fm_0,\fm_1)}
\to\rH^{2n-1}_\fT(\ol\rQ,\rR\Psi O_\lambda(n))_{(\fm_0,\fm_1)}
\]
of $O_\lambda[\Gal(\ol\dQ_p/\dQ_{p^2})]$-modules, which is an isomorphism by Lemma \ref{le:ns_product_pbc}. Combining with the diagram \eqref{eq:ns_functoriality_shimura1}, we have
\begin{align*}
&\quad\exp_\lambda\(\partial_\fp\loc_\fp\AJ(\Sh'_\sp),
\rH^1_\sing(F_\fp,\rH^{2n-1}_\et((\Sh'_{n_0}\times_{\Spec{F}}\Sh'_{n_1})_{\ol{F}},O_\lambda(n))/(\fn_0,\fn_1))\)\\
&=\exp_\lambda\(\partial\AJ_\bQ(\bQ^\eta_\sp),\rH^1_\sing(\dQ_{p^2},\rH^{2n-1}_\fT(\ol\rQ,\rR\Psi O_\lambda(n))/(\fn_0,\fn_1))\),
\end{align*}
where $\partial$ is the map \eqref{eq:partial}. Now Theorem \ref{th:first} implies
\begin{align*}
&\quad\exp_\lambda\(\partial\AJ_\bQ(\bQ^\eta_\sp),\rH^1_\sing(\dQ_{p^2},\rH^{2n-1}_\fT(\ol\rQ,\rR\Psi O_\lambda(n))/(\fn_0,\fn_1))\) \\
&=\exp_\lambda\((p+1)\phi_{\Pi_0}(\tI^\circ_{n_0,\fp})\phi_{\Pi_1}(\tT^\circ_{n_1,\fp})
\CF_{\Sh(\rV^\circ_n,\rK^\circ_\sp)},
O_\lambda[\Sh(\rV^\circ_{n_0},\rK^\circ_{n_0})]/\fn_0\otimes_{O_\lambda}O_\lambda[\Sh(\rV^\circ_{n_1},\rK^\circ_{n_1})]/\fn_1\).
\end{align*}
Note that $(p+1)$ is invertible in $O_\lambda$ by (PI2); $\phi_{\Pi_0}(\tI^\circ_{n_0,\fp})$ is invertible in $O_\lambda$ by (PI5) and Proposition \ref{pr:enumeration_even}(1); and $\phi_{\Pi_1}(\tT^\circ_{n_1,\fp})$ is invertible in $O_\lambda$ by (PI4) and Proposition \ref{pr:enumeration_odd}(2). Thus, the corollary follows.
\end{proof}

\if false

\subsection{First explicit reciprocity law}
\label{ss:first_reciprocity}

We start by choosing
\begin{itemize}[label={\ding{109}}]
  \item a prime $\lambda$ of $E$, whose underlying rational prime $\ell$ satisfies $\Sigma^+_\mnm\cap\Sigma^+_\ell=\emptyset$,

  \item a positive integer $m$,

  \item a (possibly empty) finite set $\Sigma^+_{\lr,\rI}$ of nonarchimedean places of $F^+$ that are inert in $F$,\footnote{Here, the subscript ``lr'' stands for ``level-raising'', while the subscript ``I'' (Roman number one) stands for the ``first''. In the next subsection, we will have $\Sigma^+_{\lr,\r{II}}$ for the second reciprocity law.} strongly disjoint from $\Sigma^+_\mnm$ (Definition \ref{de:strongly_disjoint}), satisfying $\ell\nmid\|v\|(\|v\|^2-1)$ for $v\in\Sigma^+_{\lr,\rI}$,

  \item a finite set $\Sigma^+_\rI$ of nonarchimedean places of $F^+$ containing $\Sigma^+_\mnm\cup\Sigma^+_{\lr,\rI}$,

  \item a standard definite hermitian space $\rV^\circ_n$ of rank $n$ over $F$, together with a self-dual $\prod_{v\not\in\Sigma^+_\infty\cup\Sigma^+_\mnm\cup\Sigma^+_{\lr,\rI}}O_{F_v}$-lattice $\Lambda^\circ_n$ in $\rV^\circ_n\otimes_F\dA_F^{\Sigma^+_\infty\cup\Sigma^+_\mnm\cup\Sigma^+_{\lr,\rI}}$ (and put $\rV^\circ_{n+1}\coloneqq(\rV^\circ_n)_\sharp$ and $\Lambda^\circ_{n+1}\coloneqq(\Lambda^\circ_n)_\sharp$), satisfying that $(\rV^\circ_{n_0})_v$ is not split for $v\in\Sigma^+_{\lr,\rI}$,

  \item objects $\rK^\circ_n\in\fK(\rV^\circ_n)$ and $(\rK^\circ_\sp,\rK^\circ_{n+1})\in\fK(\rV^\circ_n)_\sp$ of the forms
      \begin{align*}
      \rK^\circ_n&=\prod_{v\in\Sigma^+_\mnm\cup\Sigma^+_{\lr,\rI}}(\rK^\circ_n)_v
      \times\prod_{v\not\in\Sigma^+_\infty\cup\Sigma^+_\mnm\cup\Sigma^+_{\lr,\rI}}\rU(\Lambda^\circ_n)(O_{F^+_v}), \\
      \rK^\circ_\sp&=\prod_{v\in\Sigma^+_\mnm\cup\Sigma^+_{\lr,\rI}}(\rK^\circ_\sp)_v
      \times\prod_{v\not\in\Sigma^+_\infty\cup\Sigma^+_\mnm\cup\Sigma^+_{\lr,\rI}}\rU(\Lambda^\circ_n)(O_{F^+_v}), \\
      \rK^\circ_{n+1}&=\prod_{v\in\Sigma^+_\mnm\cup\Sigma^+_{\lr,\rI}}(\rK^\circ_{n+1})_v
      \times\prod_{v\not\in\Sigma^+_\infty\cup\Sigma^+_\mnm\cup\Sigma^+_{\lr,\rI}}\rU(\Lambda^\circ_{n+1})(O_{F^+_v}),
      \end{align*}
      satisfying
      \begin{itemize}
        \item $(\rK^\circ_\sp)_v=(\rK^\circ_n)_v$ for $v\in\Sigma^+_\mnm$,

        \item $(\rK^\circ_\sp)_v\subseteq(\rK^\circ_n)_v$ for $v\in\Sigma^+_{\lr,\rI}$, and

        \item $(\rK^\circ_{n_0})_v$ is a transferable open compact subgroup (Definition \ref{de:transferable}) of $\rU(\rV^\circ_{n_0})(F^+_v)$ for $v\in\Sigma^+_\mnm$ and is a special maximal subgroup of $\rU(\rV^\circ_{n_0})(F^+_v)$ for $v\in\Sigma^+_{\lr,\rI}$,
      \end{itemize}

  \item a special inert prime (Definition \ref{de:special_inert}) $\fp$ of $F^+$ (with the underlying rational prime $p$) satisfying
      \begin{description}
        \item[(PI1)] $\Sigma^+_\rI$ does not contain $p$-adic places;

        \item[(PI2)] $\ell$ does not divide $p(p^2-1)$;

        \item[(PI3)] there exists a CM type $\Phi$ containing $\tau_\infty$ as in the initial setup of \S\ref{ss:ns} satisfying $\dQ_p^\Phi=\dQ_{p^2}$;

        \item[(PI4)] $P_{\balpha(\Pi_{0,\fp})}\modulo\lambda^m$ is level-raising special at $\fp$ (Definition \ref{de:satake_condition});

            $P_{\balpha(\Pi_{1,\fp})}\modulo\lambda$ is Tate generic at $\fp$ (Definition \ref{de:satake_condition});

        \item[(PI5)] $P_{\balpha(\Pi_{\alpha,\fp})}\modulo\lambda$ is intertwining generic at $\fp$ (Definition \ref{de:satake_condition}) for $\alpha=0,1$;

        \item[(PI6)] the natural map
            \[
            \frac{(O_E/\lambda^m)[\Sh(\rV^\circ_{n_\alpha},\rK^\circ_{n_\alpha})]}
            {\dT^{\Sigma^+_\rI\cup\Sigma^+_p}_{n_\alpha}\cap\Ker\phi_{\Pi_\alpha}}
            \to\frac{(O_E/\lambda^m)[\Sh(\rV^\circ_{n_\alpha},\rK^\circ_{n_\alpha})]}
            {\dT^{\Sigma^+_\rI}_{n_\alpha}\cap\Ker\phi_{\Pi_\alpha}}
            \]
            is an isomorphism of \emph{nontrivial} $O_E/\lambda^m$-modules for $\alpha=0,1$;

        \item[(PI7)] $P_{\balpha(\Pi_{0,\fp})\otimes\balpha(\Pi_{1,\fp})}\modulo\lambda^m$ is level-raising special at $\fp$ (Definition \ref{de:satake_condition});
      \end{description}
      (So we can and will apply the setup in \S\ref{ss:ns_functoriality} to the datum $(\rV^\circ_n,\{\Lambda^\circ_{n,\fq}\}\res_{\fq\mid p})$.)

  \item remaining data in \S\ref{ss:ns_initial} with $\dQ_p^\Phi=\dQ_{p^2}$; and

  \item data as in Notation \ref{no:ns_uniformization}, which in particular give open compact subgroups $\rK^\bullet_{n,p}$ and $\rK^\bullet_{n+1,p}$.
\end{itemize}

Put $\rK^{p\circ}_\sp\coloneqq(\rK_\sp^\circ)^p$ and $\rK^\bullet_\sp\coloneqq\rK^{p\circ}_\sp\times\rK^\bullet_{n_0,p}$; put $\rK^{p\circ}_{n_\alpha}\coloneqq(\rK_{n_\alpha}^\circ)^p$ and $\rK^\bullet_{n_\alpha}\coloneqq\rK^{p\circ}_{n_\alpha}\times\rK^\bullet_{n_\alpha,p}$ for $\alpha=0,1$. As in \S\ref{ss:ns_reciprocity}, we put $\rX^?_{n_\alpha}\coloneqq\rX^?_\fp(\rV^\circ_{n_\alpha},\rK^{p\circ}_{n_\alpha})$ for meaningful triples $(\rX,?,\alpha)\in\{\bM,\rM,\rB,\rS\}\times\{\;,\eta,\circ,\bullet,\dag\}\times\{0,1\}$. For $\alpha=0,1$, let $(\pres{\alpha}\rE^{p,q}_s,\pres{\alpha}\rd^{p,q}_s)$ be the weight spectral sequence abutting to the cohomology $\rH^{p+q}_\fT(\ol\rM_{n_\alpha},\rR\Psi O_\lambda(r_\alpha))$ from \S\ref{ss:ns_weight}.

\begin{notation}\label{no:ideal}
We introduce the following ideals of $\dT^{\Sigma^+_\rI\cup\Sigma^+_p}_{n_\alpha}$, for $\alpha=0,1$
\begin{align*}
\begin{dcases}
\fm_\alpha\coloneqq\dT^{\Sigma^+_\rI\cup\Sigma^+_p}_{n_\alpha}\cap\Ker\(\dT^{\Sigma^+}_{n_\alpha}\xrightarrow{\phi_{\Pi_\alpha}}O_E\to O_E/\lambda\),\\
\fn_\alpha\coloneqq\dT^{\Sigma^+_\rI\cup\Sigma^+_p}_{n_\alpha}\cap\Ker\(\dT^{\Sigma^+}_{n_\alpha}\xrightarrow{\phi_{\Pi_\alpha}}O_E\to O_E/\lambda^m\).
\end{dcases}
\end{align*}
\end{notation}

We then introduce the following assumptions.

\begin{assumption}\label{as:first_vanishing}
For $\alpha=0,1$, we have $\rH^i_\fT(\ol\rM_{n_\alpha},\rR\Psi O_\lambda)_{\fm_\alpha}=0$ for $i\neq n_\alpha-1$, and that $\rH^{n_\alpha-1}_\fT(\ol\rM_{n_\alpha},\rR\Psi O_\lambda)_{\fm_\alpha}$ is a finite free $O_\lambda$-module.
\end{assumption}

\begin{assumption}\label{as:first_minimal}
Under Assumption \ref{as:first_irreducible}, if $n_0\geq 4$, then
\begin{enumerate}[label=(\alph*)]
  \item $\ell\geq 2(n_0+1)$ and $\ell$ is unramified in $F$;

  \item $\bar\rho_{\Pi_0,\lambda,+}$ (Remark \ref{re:single_irreducible}) is rigid for $(\Sigma^+_\mnm,\Sigma^+_{\lr,\rI})$ (Definition \ref{de:rigid}), and $\bar\rho_{\Pi_0,\lambda}\res_{\Gal(\ol{F}/F(\zeta_\ell))}$ is absolutely irreducible; and

  \item the composite homomorphism $\dT^{\Sigma^+_\mnm}_{n_0}\xrightarrow{\phi_{\Pi_0}}O_E\to O_E/\lambda$ is cohomologically generic (Definition \ref{de:generic}).
\end{enumerate}
\end{assumption}

Now we apply constructions in \S\ref{ss:ns_reciprocity}, evaluating on the object $(\rK^{p\circ}_n,\rK^{p\circ}_{n+1})$ of $\fK(\rV^\circ_n)^p\times\fK(\rV^\circ_{n+1})^p$. In particular, we have the blow-up morphism $\sigma\colon\bQ\to\bP$ from Notation \ref{no:ns_product}, and the localized spectral sequence $(\dE^{p,q}_{s,(\fm_0,\fm_1)},\rd^{p,q}_{s,(\fm_0,\fm_1)})$ from \eqref{eq:ns_product_weight}.

\begin{lem}\label{le:first_weight}
Assume Assumptions \ref{as:first_irreducible}, \ref{as:first_minimal}, \ref{as:first_vanishing}, and Hypothesis \ref{hy:unitary_cohomology} for both $n$ and $n+1$. Then
\begin{enumerate}
  \item For $(?_0,?_1)\in\{\circ,\bullet,\dag\}^2$ and $i\in\dZ$, we have a canonical isomorphism
     \[
     \rH^i_\fT(\ol\rP^{?_0,?_1},O_\lambda)_{(\fm_0,\fm_1)}
     \simeq\bigoplus_{i_0+i_1=i}\rH^{i_0}_\fT(\ol\rM^{?_0}_{n_0},O_\lambda)_{\fm_0}
     \otimes_{O_\lambda}\rH^{i_1}_\fT(\ol\rM^{?_1}_{n_1},O_\lambda)_{\fm_1}
     \]
     in $\Mod(\Gal(\ol\dF_p/\dF_{p^2}),O_\lambda)_\free$.

  \item We have $\dE^{p,q}_{2,(\fm_0,\fm_1)}=0$ if $(p,q)\not\in\{(-1,2n),(0,2n-1),(1,2n-2)\}$, and canonical isomorphisms
     \[
     \begin{dcases}
     \dE^{-1,2n}_{2,(\fm_0,\fm_1)}\simeq\pres{0}\rE^{-1,2r_0}_{2,\fm_0}\otimes_{O_\lambda}\pres{1}\rE^{0,2r_1}_{2,\fm_1}, \\
     \dE^{0,2n-1}_{2,(\fm_0,\fm_1)}\simeq\pres{0}\rE^{0,2r_0-1}_{2,\fm_0}\otimes_{O_\lambda}\pres{1}\rE^{0,2r_1}_{2,\fm_1}, \\
     \dE^{1,2n-2}_{2,(\fm_0,\fm_1)}\simeq\pres{0}\rE^{1,2r_0-2}_{2,\fm_0}\otimes_{O_\lambda}\pres{1}\rE^{0,2r_1}_{2,\fm_1},
     \end{dcases}
     \]
     in $\Mod(\Gal(\ol\dF_p/\dF_{p^2}),O_\lambda)_\free$.

  \item If $\dE^{i,2n-1-i}_{2,(\fm_0,\fm_1)}(-1)$ has a nontrivial subquotient on which $\Gal(\ol\dF_p/\dF_{p^2})$ acts trivially, then $i=1$.

  \item For $(?_0,?_1)\in\{\circ,\bullet,\dag\}^2$ and $i\in\dZ$, both $\rH^{2i}_\fT(\ol\rP^{?_0,?_1},O_\lambda(i))_{(\fm_0,\fm_1)}$ and $\rH^{2i}_\fT(\ol\rQ^{?_0,?_1},O_\lambda(i))_{(\fm_0,\fm_1)}$ are weakly semisimple.

  \item We have $\rH^i_\fT(\ol\rQ,\rR\Psi O_\lambda)_{(\fm_0,\fm_1)}=0$ for $i\neq 2n-1$.

  \item The canonical map $\rH^i_{\fT,c}(\ol\rQ^{(c)},O_\lambda)_{(\fm_0,\fm_1)}\to\rH^i_\fT(\ol\rQ^{(c)},O_\lambda)_{(\fm_0,\fm_1)}$ is an isomorphism for every integers $c$ and $i$.
\end{enumerate}
\end{lem}

\begin{proof}
For (1), by Lemma \ref{le:ns_cohomology}, Lemma \ref{le:tate_weight}(2), Theorem \ref{th:raising}(1), we know that $\rH^{i_\alpha}_\fT(\ol\rM^{?_\alpha}_{n_\alpha},O_\lambda)_{\fm_\alpha}$ is a finitely generated free $O_\lambda$-module for $\alpha=0,1$ and every $i_\alpha\in\dZ$. Thus, (1) follows from Lemma \ref{le:single_compact} and the K\"{u}nneth formula.

For (2), we first show that $\dE^{p,q}_{s,(\fm_0,\fm_1)}$ degenerates at the second page. By (1), Lemma \ref{le:ns_product_weight}(2), Lemma \ref{le:ns_cohomology}, and Lemma \ref{le:tate_vanishing}, the composition of $\rd^{-2,q}_{1,(\fm_0,\fm_1)}$ and the natural projection
\[
\dE^{-1,q}_{1,(\fm_0,\fm_1)}\to\rH^{q-2}_\fT(\ol\rQ^{\dag,\dag},O_\lambda(n-1))\bigoplus\rH^{q-2}_\fT(\ol\rQ^{\dag,\circ},O_\lambda(n-1))
\]
is injective for every $q\in\dZ$. Thus, $\rd^{-2,q}_{1,(\fm_0,\fm_1)}$ is injective, which implies $\dE^{-2,q}_{2,(\fm_0,\fm_1)}=0$ for every $q\in\dZ$. By a dual argument, we have $\dE^{2,q}_{2,(\fm_0,\fm_1)}=0$ for every $q\in\dZ$ as well. For the degeneration, it suffices to show that $\rd^{-1,q}_{1,(\fm_0,\fm_1)}$ is injective and $\rd^{0,q}_{1,(\fm_0,\fm_1)}$ is surjective for $q$ odd. By Lemma \ref{le:ns_product_weight}(2), Lemma \ref{le:ns_cohomology}, and Lemma \ref{le:tate_weight}(1), we have $\rH^{q-2}_\fT(\ol\rQ^{(1)},O_\lambda(n-1))=\rH^{q-2}_\fT(\ol\rQ^{\bullet,\dag},O_\lambda(n-1))$ for $q$ odd, which easily implies the injectivity of $\rd^{-1,q}_{1,(\fm_0,\fm_1)}$. By a dual argument, $\rd^{0,q}_{1,(\fm_0,\fm_1)}$ is surjective for $q$ odd.

Now for every $q\in\dZ$, the morphism $\sigma$ induces a map
\[
\sigma^*_1\colon\bigoplus_{q_0+q_1=q}\pres{0}\rE^{*,q_0}_{1,\fm_0}\otimes_{O_\lambda}\pres{1}\rE^{*,q_1}_{1,\fm_1}
\to\dE^{*,q}_{1,(\fm_0,\fm_1)}
\]
of complexes of $O_\lambda[\Gal(\ol\dF/\dF_{p^2})]$-modules, hence a map
\[
\sigma^*_2\colon\bigoplus_{p_0+p_1=p}\bigoplus_{q_0+q_1=q}\pres{0}\rE^{p_0,q_0}_{2,\fm_0}\otimes_{O_\lambda}\pres{1}\rE^{p_1,q_1}_{2,\fm_1}
\to\dE^{p,q}_{2,(\fm_0,\fm_1)}
\]
of $O_\lambda[\Gal(\ol\dF/\dF_{p^2})]$-modules for $(p,q)\in\dZ^2$. By Lemma \ref{le:tate_weight} and Theorem \ref{th:raising}(2), to show (2), it suffices to show that $\sigma^*_2$ is an isomorphism, or the natural map
\[
\bigoplus_{i_0+i_1=i}\rH^{i_0}_\fT(\ol\rM_{n_0},\rR\Psi O_\lambda(r_0))_{\fm_0}
\otimes_{O_\lambda}\rH^{i_1}_\fT(\ol\rM_{n_1},\rR\Psi O_\lambda(r_1))_{\fm_1}
\to\rH^{i}_\fT(\ol\rQ,\rR\Psi O_\lambda(n))_{(\fm_0,\fm_1)}
\]
induced by $\sigma$ is an isomorphism for every $i\in\dZ$. By Lemma \ref{le:ns_pbc} and Lemma \ref{le:ns_product_pbc}, the above map is identified with
\[
\bigoplus_{i_0+i_1=i}\rH^{i_0}_\fT(\bM^\eta_{n_0}\otimes_{\dQ_{p^2}}\ol\dQ_p,O_\lambda(r_0))_{\fm_0}
\otimes_{O_\lambda}\rH^{i_1}_\fT(\bM^\eta_{n_1}\otimes_{\dQ_{p^2}}\ol\dQ_p,O_\lambda(r_1))_{\fm_1}
\to\rH^{i}_\fT(\bQ^\eta\otimes_{\dQ_{p^2}}\ol\dQ_p,O_\lambda(n))_{(\fm_0,\fm_1)},
\]
which is an isomorphism by Lemma \ref{le:single_compact}, and the K\"{u}nneth formula. Thus, (2) follows.

For (3), let $\{\alpha_{0,1}^{\pm1},\dots,\alpha_{0,r_0}^{\pm1}\}$ and $\{\alpha_{1,1}^{\pm1},\dots,\alpha_{1,r_1}^{\pm1},1\}$ be the roots of $P_{\balpha(\Pi_{0,\fp})}\modulo\lambda$ and $P_{\balpha(\Pi_{1,\fp})}\modulo\lambda$ in a finite extension of $O_\lambda/\lambda$, respectively. By (PI4), we may assume $\alpha_{0,r_0}=p$. By (2), Theorem \ref{th:tate}(1), and Theorem \ref{th:raising}(3), the generalized Frobenius eigenvalues of the $(O_\lambda/\lambda)[\Gal(\ol\dF_p/\dF_{p^2})]$-modules $\dE^{-1,2n}_{2,(\fm_0,\fm_1)}(-1)\otimes_{O_\lambda}O_\lambda/\lambda$ and $\dE^{0,2n-1}_{2,(\fm_0,\fm_1)}(-1)\otimes_{O_\lambda}O_\lambda/\lambda$ are contained in $\{p^{-2}\alpha_{1,1}^{\pm1},\dots,p^{-2}\alpha_{1,r_1}^{\pm1},p^{-2}\}$ and $\{p^{-1}\alpha_{0,1}^{\pm1}\alpha_{1,1}^{\pm1},\dots,p^{-1}\alpha_{0,r_0-1}^{\pm1}\alpha_{1,r_1}^{\pm1}\}\cup
\{p^{-1}\alpha_{0,1}^{\pm1},\dots,p^{-1}\alpha_{0,r_0-1}^{\pm1}\}$, respectively. By (PI2), we have $p^2\neq 1$ in $O_\lambda/\lambda$. By (PI7), we have $\alpha_{1,i_1}\not\in\{p^2,p^{-2}\}$ for $1\leq i_1\leq r_1$, which implies $1\not\in\{p^{-2}\alpha_{1,1}^{\pm1},\dots,p^{-2}\alpha_{1,r_1}^{\pm1},p^{-2}\}$. Again by (PI7), we have $\alpha_{0,i_0}\alpha_{1,i_1}\not\in\{p,p^{-1}\}$ for $1\leq i_0< r_0$ and $1\leq i_1\leq r_1$, which implies $1\not\in\{p^{-1}\alpha_{0,1}^{\pm1}\alpha_{1,1}^{\pm1},\dots,p^{-1}\alpha_{0,r_0-1}^{\pm1}\alpha_{1,r_1}^{\pm1}\}$. By (PI4), we have $\alpha_{0,i_0}\not\in\{p,p^{-1}\}$ for $1\leq i_0< r_0$, which implies $1\not\in\{p^{-1}\alpha_{0,1}^{\pm1},\dots,p^{-1}\alpha_{0,r_0-1}^{\pm1}\}$. Thus, (3) follows.

For (4), by Lemma \ref{le:ns_product_weight} (3--5) and Lemma \ref{le:weakly_semisimple1}(1), it suffices to show that $\rH^{2i}_\fT(\ol\rP^{?_0,?_1},O_\lambda(i))_{(\fm_0,\fm_1)}$ is weakly semisimple. By (1) and Lemma \ref{le:tate_weight}(1), it suffices to show that $\rH^{2i_0}_\fT(\ol\rM^{?_0}_{n_0},O_\lambda(i_0))_{\fm_0}\otimes_{O_\lambda}\rH^{2i_1}_\fT(\ol\rM^{?_1}_{n_1},O_\lambda(i_1))_{\fm_1}$ is weakly semisimple for $i_0,i_1\in\dZ$. By Lemma \ref{le:ns_cohomology}, the action of $\Gal(\ol\dF_p/\dF_{p^2})$ on $\rH^{2i_\alpha}_\fT(\ol\rM^?_{n_\alpha},O_\lambda(i_\alpha))_{\fm_\alpha}$ is trivial for $\alpha=0,1$, $?=\circ,\dag$, and every $i_\alpha\in\dZ$. On the other hand, it is a consequence of Theorem \ref{th:raising}(2) (for $i_0$) and Lemma \ref{le:tate_weight}(3) (for $i_1$) that the action of $\Gal(\ol\dF_p/\dF_{p^2})$ on $\rH^{2i_\alpha}_\fT(\ol\rM^\bullet_{n_\alpha},O_\lambda(i_\alpha))_{\fm_\alpha}$ is trivial if $i_0\not\in\{r_0-1,r_0\}$ or $i_1\neq r_1$. By Proposition \ref{pr:raising_weight}(1,2) and Theorem \ref{th:raising}(1), the actions of $\Gal(\ol\dF_p/\dF_{p^2})$ on both $\rH^{2r_0-2}_\fT(\ol\rM^\bullet_{n_0},O_\lambda(r_0-1))_{\fm_0}$ and $\rH^{2r_0}_\fT(\ol\rM^\bullet_{n_0},O_\lambda(r_0))_{\fm_0}$ are also trivial. Thus, by Lemma \ref{le:weakly_semisimple1}(1), it remains to show that $\rH^{2r_1}_\fT(\ol\rM^\bullet_{n_1},O_\lambda(r_1))_{\fm_1}$ is weakly semisimple, which follows from Theorem \ref{th:tate}(2) as it is isomorphic to the direct sum of $\pres{1}\rE^{0,2r_1}_{2,\fm_1}$ and $\rH^{2r_1}_\fT(\ol\rM^\dag_{n_1},O_\lambda(r_1))_{\fm_1}$.

Part (5) is a direct consequence of (2).

Part (6) follows from (1), Lemma \ref{le:single_compact}, and Lemma \ref{le:ns_product_weight}(3--5).
\end{proof}

\begin{remark}
In fact, Lemma \ref{le:first_weight}(5) holds under only Assumption \ref{as:first_vanishing}; and Lemma \ref{le:first_weight}(6) holds under only Assumption \ref{as:first_irreducible}.
\end{remark}

Lemma \ref{le:first_weight}(5) induces a coboundary map
\begin{align*}
\AJ_\bQ\colon\rZ^n_\fT(\bQ^\eta)\to\rH^1(\dQ_{p^2},\rH^{2n-1}_\fT(\ol\rQ,\rR\Psi O_\lambda(n))_{(\fm_0,\fm_1)}).
\end{align*}
We also recall the singular quotient map
\[
\partial\colon\rH^1(\dQ_{p^2},\rH^{2n-1}_\fT(\ol\rQ,\rR\Psi O_\lambda(n))_{(\fm_0,\fm_1)})\to
\rH^1_\sing(\dQ_{p^2},\rH^{2n-1}_\fT(\ol\rQ,\rR\Psi O_\lambda(n))_{(\fm_0,\fm_1)}).
\]

By our choice of $\rK^\circ_n$ and $(\rK^\circ_\sp,\rK^\circ_{n+1})$, we obtain a morphism
\[
\bM_\fp(\rV^\circ_n,\rK^\circ_\sp)\to\bP
\]
which is finite. Denote by $\bP_\sp$ the corresponding cycle; and let $\bQ_\sp$ be the strict transform of $\bP_\sp$ under $\sigma$, which is a $\bT_\fp$-invariant cycle of $\bQ$. Our main goal is to compute $\partial\AJ_\bQ(\bQ^\eta_\sp)$ in $\rH^1_\sing(\dQ_{p^2},\rH^{2n-1}_\fT(\ol\rQ,\rR\Psi O_\lambda(n))/(\fn_0,\fn_1))$. Recall the map $\Delta^n$ \eqref{eq:ns_potential}; the cycle $\rQ_\sp$ gives rise to a class $\cl(\rQ_\sp)\in C^n(\rQ,L)$ (see \S\ref{ss:ns_reciprocity} for the target).

\begin{proposition}\label{pr:first_potential}
Assume Assumptions \ref{as:first_irreducible}, \ref{as:first_minimal}, \ref{as:first_vanishing}, and Hypothesis \ref{hy:unitary_cohomology} for both $n$ and $n+1$. There is a canonical isomorphism
\[
\rH^1_\sing(\dQ_{p^2},\rH^{2n-1}_\fT(\ol\rQ,\rR\Psi O_\lambda(n))_{(\fm_0,\fm_1)})
\simeq\coker\Delta^n_{(\fm_0,\fm_1)}
\]
under which $\partial\AJ_\bQ(\bQ^\eta_\sp)$ coincides with the image of $\cl(\rQ_\sp)$ in $\coker\Delta^n_{(\fm_0,\fm_1)}$.
\end{proposition}

\begin{proof}
By \cite{Liu2}*{Theorem~2.16 and Theorem~2.18},\footnote{Although it is assumed that the underlying strictly semistable scheme $X$ is proper over the base in \cite{Liu2}, the proof of relevant results works without change in our case even when $\bQ$ is not proper in view of Lemma \ref{le:first_weight}(6).} it suffices to show that $O_\lambda$ is a very nice coefficient ring for $\dE^{p,q}_{s,(\fm_0,\fm_1)}$ in the sense of \cite{Liu2}*{Definition~2.15}. In fact, in \cite{Liu2}*{Definition~2.15}, (N1) is satisfied due to Lemma \ref{le:first_weight}(2); (N2) is satisfied due to Lemma \ref{le:first_weight}(3); and (N3) is satisfied due to Lemma \ref{le:first_weight}(4) and Lemma \ref{le:weakly_semisimple1}(2).

The proposition is proved.
\end{proof}

By Construction \ref{cs:ns_nabla_product} and Remark \ref{re:ns_nabla_product}, we have a map
\begin{align*}
\nabla\colon C^n(\rQ,O_\lambda)\to
O_\lambda[\Sh(\rV^\circ_{n_0},\rK^\circ_{n_0})]\otimes_{O_\lambda}O_\lambda[\Sh(\rV^\circ_{n_1},\rK^\circ_{n_1})].
\end{align*}

\begin{theorem}[First explicit reciprocity law]\label{th:first}
Assume Assumptions \ref{as:first_irreducible}, \ref{as:first_minimal}, \ref{as:first_vanishing}, and Hypothesis \ref{hy:unitary_cohomology} for both $n$ and $n+1$.
\begin{enumerate}
  \item The image of the composite map $\nabla_{(\fm_0,\fm_1)}\circ\Delta^n_{(\fm_0,\fm_1)}$ is contained in $\fn_0.O_\lambda[\Sh(\rV^\circ_{n_0},\rK^\circ_{n_0})]_{\fm_0}\otimes_{O_\lambda}O_\lambda[\Sh(\rV^\circ_{n_1},\rK^\circ_{n_1})]_{\fm_1}$.

  \item In view of (1), the induced map
      \[
      \nabla_{\fm_1/\fn_0}\colon\coker\Delta^n_{\fm_1}/\fn_0\to
      O_\lambda[\Sh(\rV^\circ_{n_0},\rK^\circ_{n_0})]/\fn_0\otimes_{O_\lambda}O_\lambda[\Sh(\rV^\circ_{n_1},\rK^\circ_{n_1})]_{\fm_1}
      \]
      is an isomorphism.

  \item Under the natural pairing
      \[
      \resizebox{\hsize}{!}{
      \xymatrix{
      O_\lambda[\Sh(\rV^\circ_{n_0},\rK^\circ_{n_0})]/\fn_0\otimes_{O_\lambda}O_\lambda[\Sh(\rV^\circ_{n_1},\rK^\circ_{n_1})]_{\fm_1}\times
      (O_\lambda/\lambda^m)[\Sh(\rV^\circ_{n_0},\rK^\circ_{n_0})][\fn_0]
      \otimes_{O_\lambda}O_\lambda[\Sh(\rV^\circ_{n_1},\rK^\circ_{n_1})]_{\fm_1}
      \to O_\lambda/\lambda^m
      }
      }
      \]
      obtained by taking inner product, the pairing of $\nabla_{/(\fn_0,\fn_1)}(\partial\AJ_\bQ(\bQ^\eta_\graph))$ and every function $f\in (O_\lambda/\lambda^m)[\Sh(\rV^\circ_{n_0},\rK^\circ_{n_0})][\fn_0]\otimes_{O_\lambda}
      (O_\lambda/\lambda^m)[\Sh(\rV^\circ_{n_1},\rK^\circ_{n_1})][\fn_1]$ is equal to
      \[
      (p+1)\cdot\phi_{\Pi_0}(\tI^\circ_{n_0,\fp})\cdot\phi_{\Pi_1}(\tT^\circ_{n_1,\fp})\cdot
      \sum_{s\in\Sh(\rV^\circ_n,\rK^\circ_\sp)}f(s,\sh^\circ_\uparrow(s)).
      \]
      Here, we regard $\partial\AJ_\bQ(\bQ^\eta_\sp)$ as an element in $\coker\Delta^n_{(\fm_0,\fm_1)}$ (hence in $\coker\Delta^n_{\fm_1}/\fn_0$) via the canonical isomorphism in Proposition \ref{pr:first_potential}.
\end{enumerate}
\end{theorem}

\begin{proof}
We first consider (1). By Lemma \ref{le:ns_product_weight}(3,4), we have
\begin{align*}
&\rH^{2(n-1)}_\fT(\ol\rQ^{(0)},O_\lambda(n-1))_{(\fm_0,\fm_1)}
=\bigoplus_{(?_0,?_1)\in\{\circ,\bullet\}^2}\sigma^*\rH^{2(n-1)}_\fT(\ol\rP^{?_0,?_1},O_\lambda(n-1))_{(\fm_0,\fm_1)} \\
&\qquad\qquad\bigoplus(\delta^{\dag,\dag}_{\circ,\circ})_!\sigma^*\rH^{2(n-2)}_\fT(\ol\rP^{\dag,\dag},O_\lambda(n-2))_{(\fm_0,\fm_1)}
\bigoplus(\delta^{\dag,\dag}_{\bullet,\bullet})_!\sigma^*\rH^{2(n-2)}_\fT(\ol\rP^{\dag,\dag},O_\lambda(n-2))_{(\fm_0,\fm_1)}.
\end{align*}
Thus, it suffices to show that
\begin{enumerate}
  \item[(1a)] The image of $\sigma^*\rH^{2(n-1)}_\fT(\ol\rP^{\circ,\bullet},O_\lambda(n-1))_{(\fm_0,\fm_1)}
      \bigoplus\sigma^*\rH^{2(n-1)}_\fT(\ol\rP^{\bullet,\bullet},O_\lambda(n-1))_{(\fm_0,\fm_1)}$ under the map $(\nabla\circ\delta_{1!}\circ\delta_0^*)_{(\fm_0,\fm_1)}$ is contained in $\fn_0.O_\lambda[\Sh(\rV^\circ_{n_0},\rK^\circ_{n_0})]_{\fm_0}
      \otimes_{O_\lambda}O_\lambda[\Sh(\rV^\circ_{n_1},\rK^\circ_{n_1})]_{\fm_1}$.

  \item[(1b)] The image of $\sigma^*\rH^{2(n-1)}_\fT(\ol\rP^{\circ,\circ},O_\lambda(n-1))_{(\fm_0,\fm_1)}
      \bigoplus\sigma^*\rH^{2(n-1)}_\fT(\ol\rP^{\bullet,\circ},O_\lambda(n-1))_{(\fm_0,\fm_1)}$ under the map $(\nabla\circ\delta_{1!}\circ\delta_0^*)_{(\fm_0,\fm_1)}$ is zero.

  \item[(1c)] The image of $(\delta^{\dag,\dag}_{\circ,\circ})_!\sigma^*\rH^{2(n-2)}_\fT(\ol\rP^{\dag,\dag},O_\lambda(n-2))_{(\fm_0,\fm_1)}$ under the map $(\nabla\circ\delta_{1!}\circ\delta_0^*)_{(\fm_0,\fm_1)}$ is zero.

  \item[(1d)] The image of  $(\delta^{\dag,\dag}_{\bullet,\bullet})_!\sigma^*\rH^{2(n-2)}_\fT(\ol\rP^{\dag,\dag},O_\lambda(n-2))_{(\fm_0,\fm_1)}$ under the map $(\nabla\circ\delta_{1!}\circ\delta_0^*)_{(\fm_0,\fm_1)}$ is zero.
\end{enumerate}

For (1a), we have a commutative diagram
\begin{align*}
\resizebox{\hsize}{!}{
\xymatrix{
\rH^{2(n-1)}_\fT(\ol\rP^{\circ,\bullet},O_\lambda(n-1))_{(\fm_0,\fm_1)}\bigoplus \rH^{2(n-1)}_\fT(\ol\rP^{\bullet,\bullet},O_\lambda(n-1))_{(\fm_0,\fm_1)} \ar[r]\ar[d]_-{\sigma^*}
& \pres{0}\rE^{0,2r_0-2}_{1,\fm_0}\otimes_{O_\lambda}\rH^{2r_1}_\fT(\ol\rM^\bullet_{n_1},O_\lambda(r_1))_{\fm_1} \ar[d] \\
\rH^{2(n-1)}_\fT(\ol\rQ^{\circ,\bullet},O_\lambda(n-1))_{(\fm_0,\fm_1)}\bigoplus
\rH^{2(n-1)}_\fT(\ol\rQ^{\bullet,\bullet},O_\lambda(n-1))_{(\fm_0,\fm_1)} \ar[r]&
O_\lambda[\Sh(\rV^\circ_{n_0},\rK^\circ_{n_0})]_{\fm_0}
\otimes_{O_\lambda}O_\lambda[\Sh(\rV^\circ_{n_1},\rK^\circ_{n_1})]_{\fm_1}
}
}
\end{align*}
in which
\begin{itemize}[label={\ding{109}}]
  \item the upper horizontal arrow is the map
    \begin{align*}
    &\quad\rH^{2(n-1)}_\fT(\ol\rP^{\circ,\bullet},O_\lambda(n-1))_{(\fm_0,\fm_1)}\bigoplus \rH^{2(n-1)}_\fT(\ol\rP^{\bullet,\bullet},O_\lambda(n-1))_{(\fm_0,\fm_1)} \\
    &\to
    \rH^{2(r_0-1)}_\fT(\ol\rM^\circ_{n_0},O_\lambda(r_0-1))_{\fm_0}
    \otimes_{O_\lambda}\rH^{2r_1}_\fT(\ol\rM^\bullet_{n_1},O_\lambda(r_1))_{\fm_1} \\
    &\quad\quad\bigoplus
    \rH^{2(r_0-1)}_\fT(\ol\rM^\bullet_{n_0},O_\lambda(r_0-1))_{\fm_0}
    \otimes_{O_\lambda}\rH^{2r_1}_\fT(\ol\rM^\bullet_{n_1},O_\lambda(r_1))_{\fm_1} \\
    &=\pres{0}\rE^{0,2r_0-2}_{1,\fm_0}\otimes_{O_\lambda}\rH^{2r_1}_\fT(\ol\rM^\bullet_{n_1},O_\lambda(r_1))_{\fm_1}
    \end{align*}
    given by Lemma \ref{le:first_weight}(1) and the K\"{u}nneth formula;

  \item the right vertical arrow is
      \[
      (\nabla^0\circ\pres{0}\rd^{-1,2r_0}_1\circ\pres{0}\rd^{0,2r_0-2}_1(-1))_{\fm_0}\otimes (\tI^\circ_{n_1,\fp}\circ\inc_\dag^*+(p+1)^2\tT^{\circ\bullet}_{n_1,\fp}\circ\inc_\bullet^*)_{\fm_1};
      \]
      and

  \item the lower horizontal arrow is $(\nabla\circ\delta_{1!}\circ\delta_0^*)_{(\fm_0,\fm_1)}$.
\end{itemize}
For (1a), by Proposition \ref{pr:enumeration_even}(2) and (PI4), we have
\[
((p+1)\tR^\circ_{n_0,\fp}-\tI^\circ_{n_0,\fp}).O_\lambda[\Sh(\rV^\circ_{n_0},\rK^\circ_{n_0})]_{\fm_0}\subseteq
\fn_0.O_\lambda[\Sh(\rV^\circ_{n_0},\rK^\circ_{n_0})]_{\fm_0}.
\]
Thus, (1a) follows from Proposition \ref{pr:raising_weight}(4) and Lemma \ref{le:ns_product_weight}(3).

For (1b) and (1c), both images are actually contained in the sum of
\[
(\tI^\circ_{n_1,\fp}\circ\inc_{\circ,\dag}^*+(p+1)^2\tT^{\circ\bullet}_{n_1,\fp}
\circ\inc_{\circ,\bullet}^*)(\gamma^{\circ,\dag}_{\circ,\bullet})_!\rH^{2(n-1)}_\fT(\ol\rP^{\circ,\dag},O_\lambda(n-1))_{(\fm_0,\fm_1)}
\]
and
\[
(\tI^\circ_{n_1,\fp}\circ\inc_{\circ,\dag}^*+(p+1)^2\tT^{\circ\bullet}_{n_1,\fp}
\circ\inc_{\bullet,\bullet}^*)(\gamma^{\bullet,\dag}_{\bullet,\bullet})_!\rH^{2(n-1)}_\fT(\ol\rP^{\bullet,\dag},O_\lambda(n-1))_{(\fm_0,\fm_1)},
\]
which by Lemma \ref{le:first_weight}(1) coincide with
\[
\rH^{2r_0}_\fT(\ol\rM^\circ_{n_0},O_\lambda(r_0))_{\fm_0}\otimes_{O_\lambda}
\((\tI^\circ_{n_1,\fp}\circ\Inc_\dag^*+(p+1)^2\tT^{\circ\bullet}_{n_1,\fp}\circ\Inc_\bullet^*)\pres{1}\rd^{-1,2r_1}_1
\rH^{2(r_1-1)}_\fT(\ol\rM^\dag_{n_1},O_\lambda(r_1-1))_{\fm_1}\)
\]
and
\[
\rH^{2r_0}_\fT(\ol\rM^\bullet_{n_0},O_\lambda(r_0))_{\fm_0}\otimes_{O_\lambda}
\((\tI^\circ_{n_1,\fp}\circ\Inc_\dag^*+(p+1)^2\tT^{\circ\bullet}_{n_1,\fp}\circ\Inc_\bullet^*)\pres{1}\rd^{-1,2r_1}_1
\rH^{2(r_1-1)}_\fT(\ol\rM^\dag_{n_1},O_\lambda(r_1-1))_{\fm_1}\),
\]
respectively. However, they vanish by Lemma \ref{le:ns_weight_odd}(3). Thus, (1b) and (1c) follow.

For (1d), by \cite{Liu2}*{Lemma~2.4}, it follows from (1c). Thus, (1) is proved.

Now we consider (2). We claim that the map $\nabla_{(\fm_0,\fm_1)}$ (with domain $C^n(\rQ,O_\lambda)_{(\fm_0,\fm_1)}$) is surjective. In fact, consider the submodule
\[
\Ker\pres{0}\rd^{0,2r_0}_{1,\fm_0}\otimes_{O_\lambda}\Ker\pres{1}\rd^{0,2r_1}_{1,\fm_1}\subseteq \bigoplus_{(?_0,?_1)\in\{\circ,\bullet\}^2}\rH^{2(n-1)}_\fT(\ol\rP^{?_0,?_1},O_\lambda(n-1))_{(\fm_0,\fm_1)}
\]
in view of Lemma \ref{le:first_weight}(1). Then $\sigma^*\(\Ker\pres{0}\rd^{0,2r_0}_{1,\fm_0}\otimes_{O_\lambda}\Ker\pres{1}\rd^{0,2r_1}_{1,\fm_1}\)$ is contained in $C^n(\rQ,O_\lambda)_{(\fm_0,\fm_1)}$. On the other hand, the map $\nabla_{(\fm_0,\fm_1)}\circ\sigma^*$ (with domain $\Ker\pres{0}\rd^{0,2r_0}_{1,\fm_0}\otimes_{O_\lambda}\Ker\pres{1}\rd^{0,2r_1}_{1,\fm_1}$) coincides with $\nabla^0_{\fm_0}\otimes\nabla^1_{\fm_1}$, which is surjective by Proposition \ref{pr:raising_weight}(3) and Theorem \ref{th:tate}. The claim follows.

Thus, it remains to show that the domain and the target of $\nabla_{\fm_1/\fn_0}$ have the same cardinality. By Proposition \ref{pr:first_potential}, we have an isomorphism
\begin{align}\label{eq:first_1}
\coker\Delta^n_{\fm_1}/\fn_0=\coker\Delta^n_{(\fm_0,\fm_1)}/\fn_0\simeq
\rH^1_\sing(\dQ_{p^2},\rH^{2n-1}_\fT(\ol\rQ,\rR\Psi O_\lambda(n))_{(\fm_0,\fm_1)})/\fn_0
\end{align}
of $O_\lambda/\lambda^m$-modules. By Lemma \ref{le:first_weight}(2,3) and Theorem \ref{th:tate}(2), we have
\[
\resizebox{\hsize}{!}{
\xymatrix{
\rH^1_\sing(\dQ_{p^2},\rH^{2n-1}_\fT(\ol\rQ,\rR\Psi O_\lambda(n))_{(\fm_0,\fm_1)})
\simeq
\rH^1_\sing(\dQ_{p^2},\rH^{2r_0-1}_\fT(\ol\rM_{n_0},\rR\Psi O_\lambda(r_0))_{\fm_0})\otimes_{O_\lambda}
(\pres{1}\rE^{0,2r_1}_{2,\fm_1})^{\Gal(\ol\dF_p/\dF_{p^2})}.
}
}
\]
Then by Theorem \ref{th:tate}(3) and Theorem \ref{th:raising}(4), we have
\[
\eqref{eq:first_1}\simeq
O_\lambda[\Sh(\rV^\circ_{n_0},\rK^\circ_{n_0})]/\fn_0\otimes_{O_\lambda}O_\lambda[\Sh(\rV^\circ_{n_1},\rK^\circ_{n_1})]_{\fm_1}.
\]
Thus, (2) is proved.

Finally we consider (3). As $\rQ_\sp$ does not intersect with $\rQ^{\circ\bullet}$, we have
\[
\nabla(\cl(\rQ_\sp))=\nabla(\cl(\rQ^\bullet_\graph))
\]
where $\cl(\rQ^\bullet_\graph)\in\rH^{2n}_\fT(\ol\rQ^{\bullet,\bullet},O_\lambda(n))$. Then by Construction \ref{cs:ns_nabla_product}, we have
\[
\nabla(\cl(\rQ_\sp))=\((p+1)(\tT^{\circ\bullet}_{n_0,\fp}\otimes\tI^\circ_{n_1,\fp})\circ\inc_{\bullet,\dag}^*+
(p+1)^3(\tT^{\circ\bullet}_{n_0,\fp}\otimes\tT^{\circ\bullet}_{n_1,\fp})\circ\inc_{\bullet,\bullet}^*\)(\cl(\rP^\bullet_\sp)).
\]
Applying Theorem \ref{th:ns_reciprocity}(3) to the object $(\rK^\circ_\sp,\rK^\circ_{n+1})\in\fK(\rV^\circ_n)_\sp$ followed by pushforward, we know that the pairing between $\nabla_{\fm_1/\fn_0}(\cl(\rQ_\sp))$ and any function
\[
f\in(O_\lambda/\lambda^m)[\Sh(\rV^\circ_{n_0},\rK^\circ_{n_0})][\fn_0]
\otimes_{O_\lambda}(O_\lambda/\lambda^m)[\Sh(\rV^\circ_{n_1},\rK^\circ_{n_1})][\fn_1]
\]
is given by the formula
\[
(p+1)\cdot\phi_{\Pi_0}(\tI^\circ_{n_0,\fp})\cdot\phi_{\Pi_1}(\tT^\circ_{n_1,\fp})\cdot
\sum_{s\in\Sh(\rV^\circ_n,\rK^\circ_\sp)}f(s,\sh^\circ_\uparrow(s))
\]
in view of (PI6). We then obtain (3) by Proposition \ref{pr:first_potential}.

The theorem is proved.
\end{proof}

We state a corollary for later application. We choose an indefinite uniformization datum as in Notation \ref{no:ns_uniformization_indefinite}, and put $\Sh'_{n_\alpha}\coloneqq\Sh(\rV'_{n_\alpha},\tj_{n_\alpha}\rK^{p\circ}_{n_\alpha}\rK'_{n_\alpha,p})$ for $\alpha=0,1$.

Assume Assumption \ref{as:first_irreducible} and Assumption \ref{as:first_vanishing}. By Lemma \ref{le:single_compact}, Lemma \ref{le:ns_pbc}, and the K\"{u}nneth formula, we have $\rH^i_\et((\Sh'_{n_0}\times_{\Spec{F}}\Sh'_{n_1})_{\ol{F}},O_\lambda)_{(\fm_0,\fm_1)}=0$ if $i\neq 2n-1$. In particular, we obtain the Abel--Jacobi map
\[
\AJ\colon\rZ^n(\Sh'_{n_0}\times_{\Spec{F}}\Sh'_{n_1})\to
\rH^1(F,\rH^{2n-1}_\et((\Sh'_{n_0}\times_{\Spec{F}}\Sh'_{n_1})_{\ol{F}},O_\lambda(n))/(\fn_0,\fn_1)).
\]
Let $\Sh'_\sp$ be the cycle given by the finite morphism $\Sh(\rV'_n,\tj_n\rK^{p\circ}_\sp\rK'_{n,p})\to\Sh'_n\times_{\Spec{F}}\Sh'_{n+1}$, which is an element in $\rZ^n(\Sh'_{n_0}\times_{\Spec{F}}\Sh'_{n_1})$.

\begin{corollary}\label{co:first}
Assume Assumptions \ref{as:first_irreducible}, \ref{as:first_minimal}, \ref{as:first_vanishing}, and Hypothesis \ref{hy:unitary_cohomology} for both $n$ and $n+1$. Then we have
\begin{align*}
&\quad\exp_\lambda\(\partial_\fp\loc_\fp\AJ(\Sh'_\sp),
\rH^1_\sing(F_\fp,\rH^{2n-1}_\et((\Sh'_{n_0}\times_{\Spec{F}}\Sh'_{n_1})_{\ol{F}},O_\lambda(n))/(\fn_0,\fn_1))\)\\
&=\exp_\lambda\(\CF_{\Sh(\rV^\circ_n,\rK^\circ_\sp)},
O_\lambda[\Sh(\rV^\circ_{n_0},\rK^\circ_{n_0})\times\Sh(\rV^\circ_{n_1},\rK^\circ_{n_1})]/(\fn_0,\fn_1)\)
\end{align*}
where $\exp_\lambda$ is introduced in Definition \ref{de:divisibility}. Here, we regard $\CF_{\Sh(\rV^\circ_n,\rK^\circ_\sp)}$ as the pushforward of the characteristic function along the map $\Sh(\rV^\circ_n,\rK^\circ_\sp)\to\Sh(\rV^\circ_n,\rK^\circ_n)\times\Sh(\rV^\circ_{n+1},\rK^\circ_{n+1})$.
\end{corollary}

\begin{proof}
Note that the isomorphism \eqref{eq:ns_moduli_scheme_shimura} induces a map
\[
\rH^{2n-1}_\et((\Sh'_{n_0}\times_{\Spec{F}}\Sh'_{n_1})_{\ol{F}},O_\lambda(n))_{(\fm_0,\fm_1)}
\to\rH^{2n-1}_\fT(\ol\rQ,\rR\Psi O_\lambda(n))_{(\fm_0,\fm_1)}
\]
of $O_\lambda[\Gal(\ol\dQ_p/\dQ_{p^2})]$-modules, which is an isomorphism by Lemma \ref{le:ns_product_pbc}. Combining with the diagram \eqref{eq:ns_functoriality_shimura1}, we have
\begin{align*}
&\quad\exp_\lambda\(\partial_\fp\loc_\fp\AJ(\Sh'_\sp),
\rH^1_\sing(F_\fp,\rH^{2n-1}_\et((\Sh'_{n_0}\times_{\Spec{F}}\Sh'_{n_1})_{\ol{F}},O_\lambda(n))/(\fn_0,\fn_1))\)\\
&=\exp_\lambda\(\partial\AJ_\bQ(\bQ^\eta_\sp),\rH^1_\sing(\dQ_{p^2},\rH^{2n-1}_\fT(\ol\rQ,\rR\Psi O_\lambda(n))/(\fn_0,\fn_1))\).
\end{align*}
Now Theorem \ref{th:first} implies
\begin{align*}
&\quad\exp_\lambda\(\partial\AJ_\bQ(\bQ^\eta_\sp),\rH^1_\sing(\dQ_{p^2},\rH^{2n-1}_\fT(\ol\rQ,\rR\Psi O_\lambda(n))/(\fn_0,\fn_1))\) \\
&=\exp_\lambda\((p+1)\phi_{\Pi_0}(\tI^\circ_{n_0,\fp})\phi_{\Pi_1}(\tT^\circ_{n_1,\fp})
\CF_{\Sh(\rV^\circ_n,\rK^\circ_\sp)},
O_\lambda[\Sh(\rV^\circ_{n_0},\rK^\circ_{n_0})]/\fn_0\otimes_{O_\lambda}O_\lambda[\Sh(\rV^\circ_{n_1},\rK^\circ_{n_1})]/\fn_1\).
\end{align*}
Note that $(p+1)$ is invertible in $O_\lambda$ by (PI2); $\phi_{\Pi_0}(\tI^\circ_{n_0,\fp})$ is invertible in $O_\lambda$ by (PI5) and Proposition \ref{pr:enumeration_even}(1); and $\phi_{\Pi_1}(\tT^\circ_{n_1,\fp})$ is invertible in $O_\lambda$ by (PI4) and Proposition \ref{pr:enumeration_odd}(2). Thus, the corollary follows.
\end{proof}

\fi

\subsection{Second explicit reciprocity law}
\label{ss:second_reciprocity}

We start by choosing
\begin{itemize}[label={\ding{109}}]
  \item a prime $\lambda$ of $E$, whose underlying rational prime $\ell$ satisfies $\Sigma^+_\mnm\cap\Sigma^+_\ell=\emptyset$,

  \item a positive integer $m$,

  \item a (possibly empty) finite set $\Sigma^+_{\lr,\r{II}}$ of nonarchimedean places of $F^+$ that are inert in $F$, strongly disjoint from $\Sigma^+_\mnm$  (Definition \ref{de:strongly_disjoint}), satisfying $\ell\nmid\|v\|(\|v\|^2-1)$ for $v\in\Sigma^+_{\lr,\r{II}}$,

  \item a finite set $\Sigma^+_{\r{II}}$ of nonarchimedean places of $F^+$ containing $\Sigma^+_\mnm\cup\Sigma^+_{\lr,\r{II}}$,

  \item a standard indefinite hermitian space $\rV_n$ of rank $n$ over $F$, together with a self-dual $\prod_{v\not\in\Sigma^+_\infty\cup\Sigma^+_\mnm\cup\Sigma^+_{\lr,\r{II}}}O_{F_v}$-lattice $\Lambda_n$ in $\rV_n\otimes_F\dA_F^{\Sigma^+_\infty\cup\Sigma^+_\mnm\cup\Sigma^+_{\lr,\r{II}}}$ (and put $\rV_{n+1}\coloneqq(\rV_n)_\sharp$ and $\Lambda_{n+1}\coloneqq(\Lambda_n)_\sharp$), satisfying that the hermitian space $(\rV_{n_0})_v$ is not split for $v\in\Sigma^+_{\lr,\r{II}}$,

  \item objects $\rK_n\in\fK(\rV_n)$ and $(\rK_\sp,\rK_{n+1})\in\fK(\rV_n)_\sp$ of the forms
      \begin{align*}
      \rK_n&=\prod_{v\in\Sigma^+_\mnm\cup\Sigma^+_{\lr,\r{II}}}(\rK_n)_v
      \times\prod_{v\not\in\Sigma^+_\infty\cup\Sigma^+_\mnm\cup\Sigma^+_{\lr,\r{II}}}\rU(\Lambda_n)(O_{F^+_v}), \\
      \rK_\sp&=\prod_{v\in\Sigma^+_\mnm\cup\Sigma^+_{\lr,\r{II}}}(\rK_\sp)_v
      \times\prod_{v\not\in\Sigma^+_\infty\cup\Sigma^+_\mnm\cup\Sigma^+_{\lr,\r{II}}}\rU(\Lambda_n)(O_{F^+_v}), \\
      \rK_{n+1}&=\prod_{v\in\Sigma^+_\mnm\cup\Sigma^+_{\lr,\r{II}}}(\rK_{n+1})_v
      \times\prod_{v\not\in\Sigma^+_\infty\cup\Sigma^+_\mnm\cup\Sigma^+_{\lr,\r{II}}}\rU(\Lambda_{n+1})(O_{F^+_v}),
      \end{align*}
      satisfying
      \begin{itemize}
        \item $(\rK_\sp)_v=(\rK_n)_v$ for $v\in\Sigma^+_\mnm$,

        \item $(\rK_\sp)_v\subseteq(\rK_n)_v$ for $v\in\Sigma^+_{\lr,\r{II}}$, and

        \item $(\rK_{n_0})_v$ is a transferable open compact subgroup (Definition \ref{de:transferable}) of $\rU(\rV_{n_0})(F^+_v)$ for $v\in\Sigma^+_\mnm$ and is a special maximal subgroup of $\rU(\rV_{n_0})(F^+_v)$ for $v\in\Sigma^+_{\lr,\r{II}}$,
      \end{itemize}

  \item a special inert prime (Definition \ref{de:special_inert}) $\fp$ of $F^+$ (with the underlying rational prime $p$) satisfying\footnote{In what follows, we will also regard $\fp$ as the unique place of $F$ above $\fp$, according to the context.}
      \begin{description}
        \item[(PII1)] $\Sigma^+_{\r{II}}$ does not contain $p$-adic places;

        \item[(PII2)] $\ell$ does not divide $p(p^2-1)$;

        \item[(PII3)] there exists a CM type $\Phi$ containing $\tau_\infty$ as in the initial setup of \S\ref{ss:ns} satisfying $\dQ_p^\Phi=\dQ_{p^2}$;

        \item[(PII4)] $P_{\balpha(\Pi_{0,\fp})}\modulo\lambda^m$ is level-raising special at $\fp$ (Definition \ref{de:satake_condition});

            $P_{\balpha(\Pi_{1,\fp})}\modulo\lambda$ is Tate generic at $\fp$ (Definition \ref{de:satake_condition});

        \item[(PII7)] $P_{\balpha(\Pi_{0,\fp})\otimes\balpha(\Pi_{1,\fp})}\modulo\lambda^m$ is level-raising special at $\fp$ (Definition \ref{de:satake_condition});
      \end{description}
      (So we can and will apply the setup in \S\ref{ss:qs_functoriality} to the datum $(\rV_n,\{\Lambda_{n,\fq}\}\res_{\fq\mid p})$.)

  \item remaining data in \S\ref{ss:qs_initial} with $\dQ_p^\Phi=\dQ_{p^2}$; and

  \item a definite uniformization datum as in Notation \ref{no:qs_uniformization}, which in particular gives open compact subgroups $\rK^\star_{n,p}$, $\rK^\star_{n+1,p}$, and $\rK^\star_{\sp,p}$.
\end{itemize}

Put $\rK^\star_\sp\coloneqq(\ti_n\rK^p_\sp)\times\rK^\star_{n,p}$, and $\rK^\star_{n_\alpha}\coloneqq(\ti_{n_\alpha}\rK^p_{n_\alpha})\times\rK^\star_{n_\alpha,p}$ for $\alpha=0,1$. Put $\rK^\star_{\sp,\sp}\coloneqq(\ti_n\rK^p_\sp)\times\rK^\star_{\sp,p}$ and $\rK^\star_{n,\sp}\coloneqq(\ti_n\rK^p_n)\times\rK^\star_{\sp,p}$. As in \S\ref{ss:qs_reciprocity}, we put $\rX^?_{n_\alpha}\coloneqq\rX^?_\fp(\rV_{n_\alpha},\rK^p_{n_\alpha})$ for meaningful triples $(\rX,?,\alpha)\in\{\bM,\rM,\rB,\rS\}\times\{\;,\eta\}\times\{0,1\}$.

\begin{notation}
We introduce the following ideals $\fm_\alpha$ and $\fn_\alpha$ of $\dT^{\Sigma^+_\r{II}\cup\Sigma^+_p}_{n_\alpha}$ for $\alpha=0,1$ in the same way as in Notation \ref{no:ideal} (but replacing $\Sigma^+_{\rI}$ with $\Sigma^+_{\r{II}}$).
\end{notation}

We then introduce the following assumption.

\begin{assumption}\label{as:second_vanishing}
For $\alpha=0,1$, we have $\rH^i_\fT(\ol\rM_{n_\alpha},O_\lambda)_{\fm_\alpha}=0$ for $i\neq n_\alpha-1$, and that $\rH^{n_\alpha-1}_\fT(\ol\rM_{n_\alpha},O_\lambda)_{\fm_\alpha}$ is a finite free $O_\lambda$-module.
\end{assumption}

\begin{lem}\label{le:second}
Assume Assumptions \ref{as:first_irreducible}, \ref{as:second_vanishing}, and Hypothesis \ref{hy:unitary_cohomology} for $n_1$.
\begin{enumerate}
  \item The $O_\lambda[\Gal(\ol\dF_p/\dF_{p^2})]$-module $\rH^{2r_1}_\fT(\ol\rM_{n_1},O_\lambda(r_1))_{\fm_1}$ is weakly semisimple (Definition \ref{de:weakly_semisimple}).

  \item The map
     \[
     \pi_{n_1!}\circ\iota_{n_1}^*\colon(\rH^{2r_1}_\fT(\ol\rM_{n_1},O_\lambda(r_1))_{\fm_1})_{\Gal(\ol\dF_p/\dF_{p^2})}
     \to\rH^0_\fT(\ol\rS_{n_1},O_\lambda)_{\fm_1}
     \]
     is an isomorphism.
\end{enumerate}
\end{lem}

\begin{proof}
The proof of the lemma is similar to Theorem \ref{th:tate}. For the readers' convenience, we reproduce the details under the current setup.

For (1), by Lemma \ref{le:qs_pbc}, we have an isomorphism
\[
\rH^{2r_1}_\fT(\ol\rM_{n_1},O_\lambda(r_1))_{\fm_1}\simeq
\rH^{2r_1}_\et(\Sh(\rV_{n_1},\rK_{n_1})_{\ol{F}},O_\lambda(r_1))_{\fm_1}
\]
of $O_\lambda[\Gal(\ol\dQ_p/\dQ_{p^2})]$-modules. By Lemma \ref{le:single_congruent}, Lemma \ref{le:single_compact}, Proposition \ref{pr:arthur}(2), and Hypothesis \ref{hy:unitary_cohomology}, we have an isomorphism
\[
\rH^{2r_1}_\et(\Sh(\rV_{n_1},\rK_{n_1})_{\ol{F}},O_\lambda(r_1))_{\fm_1}\otimes_{O_\lambda}\ol\dQ_\ell
\simeq\bigoplus_{\pi_1}\rho_{\BC(\pi_1),\iota_\ell}^\tc(r_1)^{\oplus d(\pi_1)}
\]
of representations of $\Gamma_F$ with coefficients in $\ol\dQ_\ell$, where $d(\pi_1)\coloneqq\dim(\pi^{\infty,p}_1)^{\rK^p_{n_1}}$. Here, the direct sum is taken over all automorphic representations $\pi_1$ of $\rU(\rV_{n_1})(\dA_{F^+})$ satisfying:
\begin{itemize}[label={\ding{109}}]
  \item $(\rV_{n_1},\pi_1)$ is a $\Pi_1$-congruent standard pair (Definition \ref{de:single_congruent} with $\Sigma^+=\Sigma^+_{\r{II}}$);

  \item $\pi_{1\ul\tau_\infty}$ is a holomorphic discrete series representation of $\rU(\rV_{n_1})(F^+_{\ul\tau_\infty})$ with the Harish-Chandra parameter $\{-r_1,1-r_1,\dots,r_1-1,r_1\}$; and

  \item $\pi_{1\ul\tau}$ is trivial for every archimedean place $\ul\tau\neq\ul\tau_\infty$.
\end{itemize}
We may replace $E_\lambda$ by a finite extension inside $\ol\dQ_\ell$ such that $\rho_{\BC(\pi_1),\iota_\ell}$ is defined over $E_\lambda$ for every $\pi_1$ appearing in the previous direct sum. Now we regard $\rho_{\BC(\pi_1),\iota_\ell}$ as a representation over $E_\lambda$. Then $\rho_{\BC(\pi_1),\iota_\ell}(r_1)$ admits a $\Gamma_F$-stable $O_\lambda$-lattice $\rR_{\BC(\pi_1)}$, unique up to homothety, whose reduction $\bar\rR_{\BC(\pi_1)}$ is isomorphic to $\bar\rho_{\Pi_1,\lambda}(r_1)$. Moreover, we have an inclusion
\[
\rH^{2r_1}_\et(\Sh(\rV_{n_1},\rK_{n_1})_{\ol{F}},O_\lambda(r_1))_{\fm_1}
\subseteq\bigoplus_{\pi_1}(\rR_{\BC(\pi_1)}^\tc)^{\oplus d(\pi_1)}
\]
of $O_\lambda[\Gal(\ol\dF_p/\dF_{p^2})]$-modules. By (PII4), we know that $\bar\rho_{\Pi_1,\lambda}^\tc(r_1)$ is weakly semisimple and
\[
\dim_{O_\lambda/\lambda}\bar\rho_{\Pi_1,\lambda}^\tc(r_1)^{\Gal(\ol\dF_p/\dF_{p^2})}=1.
\]
On the other hand, we have
\[
\dim_{E_\lambda}\rho_{\BC(\pi_1),\iota_\ell}^\tc(r_1)^{\Gal(\ol\dF_p/\dF_{p^2})}\geq 1.
\]
Thus by Lemma \ref{le:weakly_semisimple2}, for every $\pi_1$ in the previous direct sum, $\rR_{\BC(\pi_1)}^\tc$ is weakly semisimple. Thus,  $\rH^{2r_1}_\fT(\ol\rM_{n_1},O_\lambda(r_1))_{\fm_1}$ is weakly semisimple by Lemma \ref{le:weakly_semisimple1}(1). Thus, (1) follows.

For (2), we note that in (1) we have also proved that $(\rH^{2r_1}_\fT(\ol\rM_{n_1},O_\lambda(r_1))_{\fm_1})_{\Gal(\ol\dF_p/\dF_{p^2})}$ is a free $O_\lambda$-module of rank $\sum_{\pi_1}d(\pi_1)$. By Theorem \ref{th:qs_tate}, Proposition \ref{pr:enumeration_odd}(2), and (PII4), we know that $\pi_{n_1!}\circ\iota_{n_1}^*$ is surjective. Thus, it remains to show that
\[
\sum_{\pi_1}d(\pi_1)\leq\dim_{E_\lambda}\rH^0_\fT(\ol\rS_{n_1},O_\lambda)_{\fm_1}\otimes_{O_\lambda}E_\lambda.
\]
However, the above inequality is a consequence of Proposition \ref{pr:qs_uniformization} and Corollary \ref{pr:jacquet_langlands_1}.

The lemma is proved.
\end{proof}

We have a finite morphism $\Sh(\rV_n,\rK_\sp)\to\Sh(\rV_n,\rK_n)\times_{\Spec{F}}\Sh(\rV_{n+1},\rK_{n+1})$, which gives rise to a class
\begin{align*}
[\Sh(\rV_n,\rK_\sp)]\in\rH^{2n}_\et(\Sh(\rV_{n_0},\rK_{n_0})\times_{\Spec{F}}\Sh(\rV_{n_1},\rK_{n_1}),O_\lambda(n))
\end{align*}
by the absolute cycle class map.

\begin{theorem}[Second explicit reciprocity law]\label{th:second}
Assume Assumptions \ref{as:first_irreducible}, \ref{as:second_vanishing}, and Hypothesis \ref{hy:unitary_cohomology} for both $n$ and $n+1$. Then we have
\begin{align*}
&\quad\exp_\lambda\(\loc_\fp([\Sh(\rV_n,\rK_\sp)]),
\rH^{2n}_\et((\Sh(\rV_{n_0},\rK_{n_0})\times_{\Spec{F}}\Sh(\rV_{n_1},\rK_{n_1}))_{F_\fp},O_\lambda(n))/(\fn_0,\fn_1)\)\\
&\leq\exp_\lambda\(\CF_{\Sh(\rV^\star_n,\rK^\star_{\sp,\sp})},
O_\lambda[\Sh(\rV^\star_{n_0},\rK^\star_{n_0})\times\Sh(\rV^\star_{n_1},\rK^\star_{n_1})]/(\fn_0,\fn_1)\)
\end{align*}
where $\loc_\fp$ is introduced in Construction \ref{cs:qs_reciprocity}; $\exp_\lambda$ is introduced in Definition \ref{de:divisibility}; and the element $\CF_{\Sh(\rV^\star_n,\rK^\star_{\sp,\sp})}$ is regarded as the pushforward of the characteristic function along the map $\Sh(\rV^\star_n,\rK^\star_{\sp,\sp})\to\Sh(\rV^\star_n,\rK^\star_n)\times\Sh(\rV^\star_{n+1},\rK^\star_{n+1})$.
\end{theorem}

\begin{proof}
We claim that
\begin{enumerate}
  \item the action of $\tT^\star_{n_1,\fp}$ on $\rH^{2r_0}_\fT(\rM_{n_0}\times_{\rT_\fp}\rS_{n_1},O_\lambda(r_0))_{(\fm_0,\fm_1)}$ is invertible; and

  \item the composite map
       \[
       (\id\times\pi_{n_1})_!\circ(\id\times\iota_{n_1})^*\colon\rH^{2n}_\fT(\rM_{n_0}\times_{\rT_\fp}\rM_{n_1},O_\lambda(n))_{(\fm_0,\fm_1)}
       \to\rH^{2r_0}_\fT(\rM_{n_0}\times_{\rT_\fp}\rS_{n_1},O_\lambda(r_0))_{(\fm_0,\fm_1)}
       \]
       is an isomorphism.
\end{enumerate}

We prove the theorem assuming these two claims. Take a uniformizer $\lambda_0$ of $E_\lambda$. Suppose that $\lambda_0^e\CF_{\Sh(\rV^\star_n,\rK^\star_{\sp,\sp})}=0$ in $O_\lambda[\Sh(\rV^\star_{n_0},\rK^\star_{n_0})\times\Sh(\rV^\star_{n_1},\rK^\star_{n_1})]/(\fn_0,\fn_1)$ for some integer $e\geq 0$. Applying Theorem \ref{th:qs_reciprocity} to the object $(\rK_\sp,\rK_{n+1})\in\fK(\rV_n)_\sp$ followed by pushforward, we have
\[
\lambda_0^e\tT^\star_{n_1,\fp}.(\id\times\pi_{n_1})_!(\id\times\iota_{n_1})^*\loc'_\fp([\Sh(\rV_n,\rK_\sp)])=0
\]
in $\rH^{2n}_\fT(\rM_{n_0}\times_{\rT_\fp}\rS_{n_1},O_\lambda(n))/(\fn_0,\fn_1)$. By the above two claims, we must have
\[
\lambda_0^e\loc'_\fp([\Sh(\rV_n,\rK_\sp)])=0
\]
in $\rH^{2n}_\fT(\rM_{n_0}\times_{\rT_\fp}\rM_{n_1},O_\lambda(n))/(\fn_0,\fn_1)$. Thus, we have
\[
\lambda_0^e\loc_\fp([\Sh(\rV_n,\rK_\sp)])=0
\]
as the map $\rH^{2n}_\et((\Sh(\rV_{n_0},\rK_{n_0})\times_{\Spec{F}}\Sh(\rV_{n_1},\rK_{n_1}))_{F_\fp},O_\lambda(n))
\to\rH^{2n}_\fT(\rM_{n_0}\times_{\rT_\fp}\rM_{n_1},O_\lambda(n))$ is an isomorphism. The theorem follows.

Now we consider the two claims. By the Hochschild--Serre spectral sequence, we have a short exact sequence
\[
\resizebox{\hsize}{!}{
\xymatrix{
0 \ar[r]& \rH^1(\dF_{p^2},\rH^{2n-1}_\fT(\ol\rM_{n_0}\times_{\ol\rT_\fp}\ol\rM_{n_1},O_\lambda(n))_{(\fm_0,\fm_1)})
\ar[r]& \rH^{2n}_\fT(\rM_{n_0}\times_{\rT_\fp}\rM_{n_1},O_\lambda(n))_{(\fm_0,\fm_1)}
\ar[r]& \rH^{2n}_\fT(\ol\rM_{n_0}\times_{\ol\rT_\fp}\ol\rM_{n_1},O_\lambda(n))_{(\fm_0,\fm_1)}^{\Gal(\ol\dF_p/\dF_{p^2})}
\ar[r]& 0
}
}
\]
of $O_\lambda$-modules. By the K\"{u}nneth formula and (an analog of) Lemma \ref{le:single_compact}, we have
\[
\rH^i_\fT(\ol\rM_{n_0}\times_{\ol\rT_\fp}\ol\rM_{n_1},O_\lambda)_{(\fm_0,\fm_1)}\simeq
\bigoplus_{i_0+i_1=i}\rH^{i_0}_\fT(\ol\rM_{n_0},O_\lambda)\otimes_{O_\lambda}\rH^{i_1}_\fT(\ol\rM_{n_1},O_\lambda)
\]
for every $i\in\dZ$. This implies $\rH^{2n}_\fT(\ol\rM_{n_0}\times_{\ol\rT_\fp}\ol\rM_{n_1},O_\lambda(n))_{(\fm_0,\fm_1)}=0$ and
\[
\rH^{2n-1}_\fT(\ol\rM_{n_0}\times_{\ol\rT_\fp}\ol\rM_{n_1},O_\lambda(n))_{(\fm_0,\fm_1)}
\simeq
\rH^{2r_0-1}_\fT(\ol\rM_{n_0},O_\lambda(r_0))_{\fm_0}\otimes_{O_\lambda}\rH^{2r_1}_\fT(\ol\rM_{n_1},O_\lambda(r_1))_{\fm_1}.
\]
In particular, we have a canonical isomorphism
\begin{align}\label{eq:second_1}
\rH^{2n}_\fT(\rM_{n_0}\times_{\rT_\fp}\rM_{n_1},O_\lambda(n))_{(\fm_0,\fm_1)}\simeq
\rH^1(\dF_{p^2},\rH^{2r_0-1}_\fT(\ol\rM_{n_0},O_\lambda(r_0))_{\fm_0}\otimes_{O_\lambda}\rH^{2r_1}_\fT(\ol\rM_{n_1},O_\lambda(r_1))_{\fm_1}).
\end{align}
Similarly, we have
\begin{align}\label{eq:second_2}
\rH^{2r_0}_\fT(\rM_{n_0}\times_{\rT_\fp}\rS_{n_1},O_\lambda(r_0))_{(\fm_0,\fm_1)}&\simeq
\rH^1(\dF_{p^2},\rH^{2r_0-1}_\fT(\ol\rM_{n_0},O_\lambda(r_0))_{\fm_0}\otimes_{O_\lambda}\rH^0_\fT(\ol\rS_{n_1},O_\lambda)_{\fm_1}) \\
&=\rH^1(\dF_{p^2},\rH^{2r_0-1}_\fT(\ol\rM_{n_0},O_\lambda(r_0))_{\fm_0})\otimes_{O_\lambda}\rH^0_\fT(\ol\rS_{n_1},O_\lambda)_{\fm_1}. \notag
\end{align}

For claim (1), note that the action of $\dT_{n_1,\fp}$ on $\rH^{2r_0}_\fT(\rM_{n_0}\times_{\rT_\fp}\rS_{n_1},O_\lambda(r_0))_{(\fm_0,\fm_1)}$ factors through the second factor under the isomorphism \eqref{eq:second_2}. By Proposition \ref{pr:enumeration_odd}(2) and (PII4), we know that the action of $\tT^\star_{n_1,\fp}$ on $\rH^0_\fT(\ol\rS_{n_1},O_\lambda)_{\fm_1}$ is invertible. Thus, (1) follows.

For claim (2), by (PII7) and a similar argument for the proof of Lemma \ref{le:first_weight}(3), we know that the $O_\lambda[\Gal(\ol\dF_p/\dF_{p^2})]$-module
\[
\rH^{2r_0-1}_\fT(\ol\rM_{n_0},O_\lambda(r_0))_{\fm_0}\otimes_{O_\lambda}
\Ker\((\rH^{2r_1}_\fT(\ol\rM_{n_1},O_\lambda(r_1))_{\fm_1})
\to(\rH^{2r_1}_\fT(\ol\rM_{n_1},O_\lambda(r_1))_{\fm_1})_{\Gal(\ol\dF_p/\dF_{p^2})}\)
\]
has zero $\Gal(\ol\dF_p/\dF_{p^2})$-coinvariants. Combining with Lemma \ref{le:second}, we obtain an isomorphism
\[
\resizebox{\hsize}{!}{
\xymatrix{
\rH^{2n}_\fT(\rM_{n_0}\times_{\rT_\fp}\rM_{n_1},O_\lambda(n))_{(\fm_0,\fm_1)}\simeq
\rH^1(\dF_{p^2},\rH^{2r_0-1}_\fT(\ol\rM_{n_0},O_\lambda(r_0))_{\fm_0})\otimes_{O_\lambda}
(\rH^{2r_1}_\fT(\ol\rM_{n_1},O_\lambda(r_1))_{\fm_1})_{\Gal(\ol\dF_p/\dF_{p^2})}
}
}
\]
from \eqref{eq:second_1}, under which the map $(\id\times\pi_{n_1})_!\circ(\id\times\iota_{n_1})^*$ coincides with $\id\otimes(\pi_{n_1!}\circ\iota_{n_1}^*)$. Thus, (2) follows.

The theorem is proved.
\end{proof}

\section{Proof of main theorems}
\label{ss:8}

In the section, we prove our main theorems on bounding Selmer groups. In \S\ref{ss:admissible_prime}, we introduce the notation of admissible primes for the coefficient field, and make some additional preparation for the main theorems. In \S\ref{ss:main_0} and \S\ref{ss:main_1}, we prove our main theorems in the (Selmer) rank $0$ and $1$ cases, respectively.

\subsection{Admissible primes for coefficient fields}
\label{ss:admissible_prime}

We keep the setup in \S\ref{ss:setup}.

\begin{definition}\label{de:lambda}
We say that a prime $\lambda$ of $E$, with the underlying rational prime $\ell$ (and the ring of integers $O_\lambda$ of $E_\lambda$), is \emph{admissible} (with respect to $(\Pi_0,\Pi_1)$) if
\begin{description}
  \item[(L1)] $\ell>4n$ and $\ell$ is unramified in $F$;

  \item[(L2)] $\Sigma^+_\mnm$ does not contain $\ell$-adic places;

  \item[(L3)] the Galois representation $\rho_{\Pi_0,\lambda}\otimes_{E_\lambda}\rho_{\Pi_1,\lambda}$ is absolutely irreducible;

  \item[(L4)] Assumption \ref{as:first_irreducible} is satisfied, that is, both $\rho_{\Pi_0,\lambda}$ and $\rho_{\Pi_1,\lambda}$ are residually absolutely irreducible;

  \item[(L5)] under (L4), for $\alpha=0,1$, we have a $\Gamma_F$-stable $O_\lambda$-lattice $\rR_\alpha$ in $\rho_{\Pi_\alpha,\lambda}(r_\alpha)$, unique up to homothety, that is $(1-\alpha)$-polarizable, for which we choose a $(1-\alpha)$-polarization $\Xi_\alpha\colon\rR_\alpha^\tc\xrightarrow{\sim}\rR_\alpha^\vee(1-\alpha)$ and an isomorphism $\rR_\alpha\simeq O_\lambda^{\oplus n_\alpha}$ of $O_\lambda$-modules.\footnote{In fact, (L5) does not depend on the choice of $\Xi_\alpha$ and the basis, since $\Xi_\alpha$ is unique up to units in $O_\lambda$ and the basis is unique up to conjugation in $\GL_{n_\alpha}(O_\lambda)$.} After adopting the notation in \S\ref{ss:rankin_selberg}, we have
      \begin{description}
        \item[(L5-1)] either one of the two assumptions in Lemma \ref{le:restriction} is satisfied;

        \item[(L5-2)] $(\r{GI}_{F',\sP}^1)$ from Lemma \ref{le:galois_element} holds with $F'=F^+_{\r{rflx}}$ (Definition \ref{de:reflexive}) and $\sP(T)=T^2-1$ (see Remark \ref{re:lambda} below for a more explicit description);
      \end{description}

  \item[(L6)] under (L4), the homomorphism $\bar\rho_{\Pi_0,\lambda,+}$ (Remark \ref{re:single_irreducible}) is rigid for $(\Sigma^+_\mnm,\emptyset)$ (Definition \ref{de:rigid}), and $\bar\rho_{\Pi_0,\lambda}\res_{\Gal(\ol{F}/F(\zeta_\ell))}$ is absolutely irreducible;

  \item[(L7)] for $\alpha=0,1$, the composite homomorphism $\dT^{\Sigma^+_\mnm}_{n_\alpha}\xrightarrow{\phi_{\Pi_\alpha}}O_E\to O_E/\lambda$ is cohomologically generic (Definition \ref{de:generic}).
\end{description}
\end{definition}

\begin{remark}\label{re:lambda}
In Definition \ref{de:lambda}, (L5-2) is equivalent to the following assertion: the image of the restriction of the homomorphism
\[
(\bar\rho_{0+},\bar\rho_{1+},\bar\epsilon_\ell)\colon\Gamma_{F^+}\to
\sG_{n_0}(O_\lambda/\lambda)\times\sG_{n_1}(O_\lambda/\lambda)\times(O_\lambda/\lambda)^\times
\]
(see Notation \ref{no:reduction} for the notation) to $\Gal(\ol{F}/F^+_{\r{rflx}})$ contains an element $(\gamma_0,\gamma_1,\xi)$ satisfying
\begin{enumerate}[label=(\alph*)]
  \item $\xi^2-1\neq 0$;

  \item for $\alpha=0,1$, $\gamma_\alpha$ belongs to $(\GL_{n_\alpha}(O_\lambda/\lambda)\times(O_\lambda/\lambda)^\times)\fc$ with order coprime to $\ell$;

  \item $1$ appears in the eigenvalues of each of $h_{\gamma_0}$, $h_{\gamma_1}$, and $h_{\gamma_0}\otimes h_{\gamma_1}$ (Notation \ref{no:h_gamma}) with multiplicity one;

  \item $h_{\gamma_0}$ does not have an eigenvalue that is equal to $-1$ in $O_\lambda/\lambda$;

  \item $h_{\gamma_1}$ does not have an eigenvalue that is equal to $-\xi$ in $O_\lambda/\lambda$.
\end{enumerate}
\end{remark}

\begin{lem}\label{pr:elliptic_admissible}
Suppose that $F^+\neq\dQ$, that $E=\dQ$, and that there are two elliptic curves $A_0$ and $A_1$ over $F^+$ such that for every rational prime $\ell$ of $E$ and $\alpha=0,1$, we have $\rho_{\Pi_\alpha,\ell}\simeq\Sym^{n_\alpha-1}\rH^1_\et({A_\alpha}_{\ol{F}},\dQ_\ell)\res_{\Gamma_F}$. If ${A_0}_{\ol{F}}$ and ${A_1}_{\ol{F}}$ are not isogenous to each other and $\End({A_0}_{\ol{F}})=\End({A_1}_{\ol{F}})=\dZ$, then all but finitely many rational primes $\ell$ are admissible.
\end{lem}

\begin{proof}
We need to show that every condition in Definition \ref{de:lambda} excludes only finitely many $\ell$. By \cite{Ser72}*{Th\'{e}or\`{e}me~6}, for sufficiently large $\ell$, the homomorphisms
\[
\Gamma_{F^+}\to\GL(\rH^1_\et({A_\alpha}_{\ol{F}},\dF_\ell))\simeq\GL_2(\dF_\ell)
\]
are both surjective for $\alpha=0,1$. Thus, we may assume that this is the case.

For (L1) and (L2), this is trivial.

For (L3), (L4), and (L5), this has been proved in Proposition \ref{pr:elliptic_curve}.

For (L6), by \cite{LTXZZ}*{Corollary~4.1.2}, the condition that $\bar\rho_{\Pi_0,\lambda,+}$ is rigid for $(\Sigma^+_\mnm,\emptyset)$ excludes only finitely many $\ell$. It is clear that the remaining two conditions also exclude only finitely many $\ell$.

For (L7), this follows from Corollary \ref{co:generic}.
\end{proof}

\begin{lem}\label{le:abstract}
Keep the setup in \S\ref{ss:setup}. Suppose that
\begin{enumerate}[label=(\alph*)]
  \item there exists a very special inert prime $\fp$ of $F^+$ (Definition \ref{de:special_inert}) such that $\Pi_{0,\fp}$ is Steinberg, and $\Pi_{1,\fp}$ is unramified whose Satake parameter contains $1$ exactly once;

  \item for $\alpha=0,1$, there exists a nonarchimedean place $w_\alpha$ of $F$ such that $\Pi_{\alpha,w_\alpha}$ is supercuspidal; and

  \item $F^+\neq\dQ$.
\end{enumerate}
Then all but finitely many primes $\lambda$ of $E$ are admissible.
\end{lem}

\begin{proof}
We need to show that every condition in Definition \ref{de:lambda} excludes only finitely many $\lambda$.

For (L1) and (L2), this is trivial.

For (L4), this follows from \cite{LTXZZ}*{Proposition~4.2.3(1)} by (b).

For (L3), this follows from Lemma \ref{le:jordan} below by (L4) and (a).

For (L6), this follows from \cite{LTXZZ}*{Theorem~4.2.6} by (b).

For (L7), this follows from Corollary \ref{co:generic} by (c).

For (L5-1), let $\lambda$ be a prime of $E$ satisfying (L4) and (L6), whose underlying rational prime is at least $2n(n+1)-1$. Then by (a), $\bar\rho_{\Pi_0,\lambda}$ and $\bar\rho_{\Pi_1,\lambda}$ satisfy the assumptions in Lemma \ref{le:jordan} below, with $k=O_\lambda/\lambda$ and $\Gamma=\Gamma_F$. Thus, by Lemma \ref{le:jordan}(2), assumption (b) of Lemma \ref{le:restriction}, hence (L5-1) hold.

For (L5-2), take an arithmetic Frobenius element $\phi_\fp\in\Gamma_{F^+_\fp}$. By Definition \ref{de:special_inert}, $\phi_\fp$ belongs to $\Gal(\ol{F}/F^+_{\r{rflx}})$. For $\alpha=0,1$, put $r_\alpha\coloneqq\lfloor\tfrac{n_\alpha}{2}\rfloor$ as always. By (a), the Satake parameter of $\Pi_{0,\fp}$ is $\{p^{\pm{1}},\dots,p^{\pm(2r_0-1)}\}$; and we may write the Satake parameter of $\Pi_{1,\fp}$ as $\{1,\alpha_1^{\pm{1}},\dots,\alpha_{r_1}^{\pm{1}}\}$ in which $\alpha_i$ is an algebraic number other than $1$ for $1\leq i\leq r_1$. For our purpose, we may replace $E$ by a finite extension in $\dC$ such that $\alpha_i\in E$ for $1\leq i\leq r_1$. By Proposition \ref{pr:galois}(1), we have $|\alpha_i|=1$ for $1\leq i\leq r_1$. Therefore, for all but finitely many prime $\lambda$ of $E$, we have
\begin{itemize}[label={\ding{109}}]
  \item $\{p,\alpha_1,\dots,\alpha_{r_1}\}$ is contained in $O_\lambda^\times$;

  \item $\{p^{\pm{1}}\modulo\lambda,\dots,p^{\pm(2r_0-1)}\modulo\lambda\}$ consists of distinct elements and does not contain $-1$;

  \item $\{\alpha_i\modulo\lambda\res 1\leq i\leq r_1\}$ is disjoint from $\{1,-p,-p^{-1}\}$;

  \item $\{p^{\pm{1}}\alpha_i\modulo\lambda,\dots,p^{\pm(2r_0-1)}\alpha_i\modulo\lambda\res 1\leq i\leq r_1\}$ is disjoint from $\{p,p^{-1}\}$.
\end{itemize}
Then for every prime $\lambda$ satisfying (L4) and the above properties, (L5-2) (that is, $(\r{GI}_{F',\sP}^1)$ from Lemma \ref{le:galois_element}) is satisfied by taking the element $(\bar\rho_{0+},\bar\rho_{1+},\bar\epsilon_\ell)(\phi_\fp)$.

The lemma is proved.
\end{proof}

For every integer $m\geq 1$, we denote by $J_m$ the standard upper triangular nilpotent Jordan block
\[
\begin{pmatrix}
0 & 1 & 0 & \cdots & 0 \\
& 0 & 1 & \cdots & 0 \\
&  & \ddots & \ddots & \vdots \\
&&& 0 & 1\\
&&&& 0
\end{pmatrix}
\]
or size $m$.

\begin{lem}\label{le:jordan}
Let $\Gamma$ be a group, and $k$ a field of characteristic either zero or at least $2n(n+1)-1$. Let $\rho_0\colon\Gamma\to\GL_{n_0}(k)$ and $\rho_1\colon\Gamma\to\GL_{n_1}(k)$ be two homomorphisms that are absolutely irreducible. Suppose that there exists an element $t\in\Gamma$ such that $\rho_0(t)=1+J_{n_0}$ and $\rho_1(t)=1$. Then we have
\begin{enumerate}
  \item $\rho_0\otimes\rho_1$ is absolutely irreducible;

  \item $\rho_0\otimes\rho_1$ is not a subquotient of $\ad(\rho_0\otimes\rho_1)$.
\end{enumerate}
\end{lem}

\begin{proof}
We may assume that $k$ is algebraically closed. For $\alpha=0,1$, let $V_i=k^{\oplus n_i}$ be the space which $\Gamma$ acts on through $\rho_\alpha$. By \cite{Ser94}*{Corollaire~1}, we know that both $\rho_0\otimes\rho_1$ and $\ad(\rho_0\otimes\rho_1)$ are semisimple.

For (1), we fix an element $e\in V_0$ such that the $t$-invariant subspace of $V_0$ is spanned by $e$. Then it is clear that the $t$-invariant subspace of $V_0\otimes_kV_1$ is $k.e\otimes_k V_1$. Now suppose that $W$ is a nonzero direct summand of the $k[\Gamma]$-module $V_0\otimes_kV_1$. Let $V'_1\subseteq V_1$ be the subspace such that $k.e\otimes_k V'_1$ is the $t$-invariant subspace of $W$. Then it is easy to see that $V'_1$ is closed under the action of $\Gamma$, which forces $V'_1=V_1$ since $\rho_1$ is irreducible. This further implies that $W=V_0\otimes_kV_1$ by looking at the Jordan decomposition of $t$ on $W$, hence $\rho_0\otimes\rho_1$ is irreducible.

For (2), note that $(\rho_0\otimes\rho_1)(t)$ is conjugate to $(1+J_{n_0})^{\oplus n_1}$. On the other hand, $\ad(\rho_0\otimes\rho_1)(t)$ is conjugate to
\[
\bigoplus_{i=1}^{n_0}(1+J_{2i-1})^{\oplus n_1^2}.
\]
Since $n_0$ is even and $1,3,\dots,2n_0-1$ are odd, $\rho_0\otimes\rho_1$ is not a subquotient of $\ad(\rho_0\otimes\rho_1)$ as $\ad(\rho_0\otimes\rho_1)$ is semisimple.

The lemma is proved.
\end{proof}

The following two lemmas will be used in later subsections.

\begin{lem}\label{le:purity}
The representation $\rho_{\Pi_0,\lambda}\otimes_{E_\lambda}\rho_{\Pi_1,\lambda}(n)$ is pure of weight $-1$ at every nonarchimedean place $w$ of $F$ not above $\ell$ (Definition \ref{de:purity}).
\end{lem}

\begin{proof}
It suffices to show that for $\alpha=0,1$, $\rho_{\Pi_\alpha,\lambda}\res_{\Gamma_{F_w}}$ is pure of some weight. By \cite{TY07}*{Lemma~1.4(3)} and Proposition \ref{pr:galois}(2), it follows from the fact that $\Pi_{\alpha,w}$ is tempered, which is ensured by Proposition \ref{pr:galois}(1).
\end{proof}

\begin{lem}\label{le:intertwining}
Assume Hypothesis \ref{hy:unitary_cohomology} for $n_1$. Let $\rV_{n_1}$ be a standard indefinite hermitian space of rank $n_1$ over $F$, $\Lambda_{n_1}$ a self-dual $\prod_{v\not\in\Sigma^+_\infty\cup\Sigma^+_\mnm}O_{F_v}$-lattice in $\rV_{n_1}\otimes_F\dA_F^{\Sigma^+_\infty\cup\Sigma^+_\mnm}$, and $\lambda$ a prime of $E$. Consider a finite set $\fP$ of special inert primes of $F^+$ whose underlying rational primes are distinct and coprime to $\Sigma^+_\mnm$, and an object $\rK_{n_1}\in\fK(\rV_{n_1})$ of the form $(\rK_{n_1})_{\Sigma^+_\mnm}\times\prod_{v\not\in\Sigma^+_\infty\cup\Sigma^+_\mnm}\rU(\Lambda_{n_1})(O_{F^+_v})$. Put
\[
\fm_1\coloneqq\dT^{\Sigma^+_\mnm\cup\Sigma^+_\fP}_{n_1}\cap\Ker\(\dT^{\Sigma^+_\mnm}_{n_1}\xrightarrow{\phi_{\Pi_1}}O_E\to O_E/\lambda\)
\]
where $\Sigma^+_\fP$ is the union of $\Sigma^+_p$ for all underlying rational primes $p$ of $\fP$. Suppose that $P_{\balpha(\Pi_{1,\fp})}\modulo\lambda$ is intertwining generic (Definition \ref{de:satake_condition}) for every $\fp\in\fP$, and that the composite homomorphism $\dT^{\Sigma^+_\mnm}_{n_1}\xrightarrow{\phi_{\Pi_1}}O_E\to O_E/\lambda$ is cohomologically generic. Then for every special maximal subgroup $\rK'_{n_1,\fP}$ of $\prod_{\fp\in\fP}\rU(\rV_{n_1})(F^+_{\fp})$ and every $i\in\dZ$, we have an isomorphism
\begin{align*}
\rH^i_\et(\Sh(\rV_{n_1},\rK_{n_1})_{\ol{F}},O_\lambda)_{\fm_1}
\simeq\rH^i_\et(\Sh(\rV_{n_1},\rK^\fP_{n_1}\rK'_{n_1,\fP})_{\ol{F}},O_\lambda)_{\fm_1}
\end{align*}
of $O_\lambda[\Gamma_F]$-modules.
\end{lem}

\begin{proof}
We first note that for every $\fp\in\fP$, $\rU(\rV_{n_1})(F^+_\fp)$ has two special maximal subgroups up to conjugation, exact one of which is hyperspecial maximal.

For the lemma, it suffices to show the following: For every $\fp\in\fP$, every special maximal subgroup $\rK'^\fp_{n_1,\fP}$ of $\prod_{\fp'\in\fP\setminus\{\fp\}}\rU(\rV_{n_1})(F^+_{\fp'})$, every hyperspecial maximal subgroup $\rK^\circ_{n_1,\fp}$ of $\rU(\rV_{n_1})(F^+_\fp)$, and every non-hyperspecial special maximal subgroup $\rK^\bullet_{n_1,\fp}$ of $\rU(\rV_{n_1})(F^+_\fp)$, there is an isomorphism
\begin{align*}
\rH^i_\et(\Sh(\rV_{n_1},\rK_{n_1}^\fP\rK'^\fp_{n_1,\fP}\rK^\circ_{n_1,\fp})_{\ol{F}},O_\lambda)_{\fm_1}
\simeq
\rH^i_\et(\Sh(\rV_{n_1},\rK_{n_1}^\fP\rK'^\fp_{n_1,\fP}\rK^\bullet_{n_1,\fp})_{\ol{F}},O_\lambda)_{\fm_1}
\end{align*}
of $O_\lambda[\Gamma_F]$-modules for every $i\in\dZ$.

Fix an isomorphism $\iota_\ell\colon\dC\simeq\ol\dQ_\ell$ that induces the prime $\lambda$ of $E$. Since the composite homomorphism $\dT^{\Sigma^+_\mnm}_{n_1}\xrightarrow{\phi_{\Pi_1}}O_E\to O_E/\lambda$ is cohomologically generic, we have for $?\in\{\circ,\bullet\}$, $\rH^i_\et(\Sh(\rV_{n_1},\rK_{n_1}^\fP\rK'^\fp_{n_1,\fP}\rK^?_{n_1,\fp})_{\ol{F}},O_E/\lambda)_{\fm_1}=0$ for $i\neq 2r_1$, hence $\rH^i_\et(\Sh(\rV_{n_1},\rK_{n_1}^\fP\rK'^\fp_{n_1,\fP}\rK^?_{n_1,\fp})_{\ol{F}},O_\lambda)_{\fm_1}$ is $O_\lambda$-torsion free for every $i\in\dZ$. Thus, it suffices to show that there is an isomorphism
\begin{align}\label{eq:intertwining}
\rH^i_\et(\Sh(\rV_{n_1},\rK_{n_1}^\fP\rK'^\fp_{n_1,\fP}\rK^\circ_{n_1,\fp})_{\ol{F}},O_\lambda)_{\fm_1}\otimes_{O_\lambda}\ol\dQ_\ell
\simeq
\rH^i_\et(\Sh(\rV_{n_1},\rK_{n_1}^\fP\rK'^\fp_{n_1,\fP}\rK^\bullet_{n_1,\fp})_{\ol{F}},O_\lambda)_{\fm_1}\otimes_{O_\lambda}\ol\dQ_\ell
\end{align}
of $\ol\dQ_\ell[\Gamma_F]$-modules for every $i\in\dZ$. Let $\Lambda^\circ_{n_1,\fp}$ be the self-dual $O_{F_\fp}$-lattice in $\rV_{n_1}\otimes_FF_\fp$ whose stabilizer is $\rK^\circ_{n_1,\fp}$. Without loss of generality, we may assume that $\rK^\bullet_{n_1,\fp}$ is the stabilizer of a lattice $\Lambda^\bullet_{n_1,\fp}$ satisfying $\Lambda^\circ_{n_1,\fp}\subseteq\Lambda^\bullet_{n_1,\fp}$ and $(\Lambda^\bullet_{n_1,\fp})^\vee/\fp\Lambda^\bullet_{n_1,\fp}\simeq\dF_{p^2}$. To show \eqref{eq:intertwining}, it suffices to show that for every (necessarily cuspidal) automorphic representation $\pi_1$ of $\rU(\rV_{n_1})(\dA_{F^+})$ that appears in either side of \eqref{eq:intertwining}, the maps
\begin{align}\label{eq:intertwining_1}
\tT^{\bullet\circ}_{n_1,\fp}\colon\pi_{1,\fp}^{\rK^\circ_{n_1,\fp}}\to\pi_{1,\fp}^{\rK^\bullet_{n_1,\fp}},\quad
\tT^{\circ\bullet}_{n_1,\fp}\colon\pi_{1,\fp}^{\rK^\bullet_{n_1,\fp}}\to\pi_{1,\fp}^{\rK^\circ_{n_1,\fp}}
\end{align}
are both isomorphisms. Here, $\tT^{\bullet\circ}_{n_1,\fp}$ and $\tT^{\circ\bullet}_{n_1,\fp}$ are introduced in Definition \ref{de:ns_hecke}. By the Chebotarev density theorem, $\rho_{\BC(\pi_1),\iota_\ell}$ and $\rho_{\Pi_1,\lambda}\otimes_{E_\lambda}\ol\dQ_\ell$ have the isomorphic (irreducible) residual representations. In particular, the Satake parameter of $\BC(\pi_1)_\fp$ does not contain $\{-p,-p^{-1}\}$ by Proposition \ref{pr:galois}(2) and the assumption that $P_{\balpha(\Pi_{1,\fp})}\modulo\lambda$ is intertwining generic. Let $\tilde\pi$ be an (unramified) principal series representation of $\rU(\rV_{n_1})(F^+_\fp)$ that has $\pi_{1,\fp}$ as a constituent. By Proposition \ref{pr:enumeration_odd}(1) and the definition of the intertwining Hecke operator $\rI^\circ_{n_1,\fp}\coloneqq\tT^{\circ\bullet}_{n_1,\fp}\circ\tT^{\bullet\circ}_{n_1,\fp}$ from Definition \ref{de:ns_hecke} or Definition \ref{de:lattice}, the composite map $\tT^{\circ\bullet}_{n_1,\fp}\circ\tT^{\bullet\circ}_{n_1,\fp}\colon\tilde\pi^{\rK^\circ_{n_1,\fp}}\to\tilde\pi^{\rK^\circ_{n_1,\fp}}$ is an isomorphism. Since both $\rK^\circ_{n_1,\fp}$ and $\rK^\bullet_{n_1,\fp}$ are special maximal subgroups of $\rU(\rV_{n_1})(F^+_\fp)$, both $\tilde\pi^{\rK^\circ_{n_1,\fp}}$ and $\tilde\pi^{\rK^\bullet_{n_1,\fp}}$ are one-dimensional. It follows that the constituent of $\tilde\pi$ that has nonzero $\rK^\circ_{n_1,\fp}$-invariants is the same as the constituent that has nonzero $\rK^\bullet_{n_1,\fp}$-invariants, which further implies that the two maps in \eqref{eq:intertwining_1} are both isomorphisms. Thus, we obtain the isomorphism \eqref{eq:intertwining}.

The lemma is proved.
\end{proof}

\subsection{Main theorem in the Selmer rank 0 case}
\label{ss:main_0}

The following lemma is a key ingredient in the proof of Theorem \ref{th:selmer0}, which is essentially the solution of the Gan--Gross--Prasad conjecture for $\Pi_0\times\Pi_1$.

\begin{lem}\label{le:ggp}
Keep the setup in \S\ref{ss:setup}. If $L(\frac{1}{2},\Pi_0\times\Pi_1)\neq 0$, then there exist
\begin{itemize}[label={\ding{109}}]
  \item a standard definite hermitian space $\rV^\circ_n$ of rank $n$ over $F$, together with a self-dual $\prod_{v\not\in\Sigma^+_\infty\cup\Sigma^+_\mnm}O_{F_v}$-lattice $\Lambda^\circ_n$ in $\rV^\circ_n\otimes_F\dA_F^{\Sigma^+_\infty\cup\Sigma^+_\mnm}$ (and put $\rV^\circ_{n+1}\coloneqq(\rV^\circ_n)_\sharp$ and $\Lambda^\circ_{n+1}\coloneqq(\Lambda^\circ_n)_\sharp$),

  \item an object $(\rK^\circ_n,\rK^\circ_{n+1})\in\fK(\rV^\circ_n)_\sp$ in which $\rK^\circ_{n_\alpha}$ is of the form
      \[
      \rK^\circ_{n_\alpha}=\prod_{v\in\Sigma^+_\mnm}(\rK^\circ_{n_\alpha})_v\times
      \prod_{v\not\in\Sigma^+_\infty\cup\Sigma^+_\mnm}\rU(\Lambda^\circ_{n_\alpha})(O_{F^+_v})
      \]
      for $\alpha=0,1$,
\end{itemize}
such that
\[
\sum_{s\in\Sh(\rV^\circ_n,\rK^\circ_n)}f(s,\sh_\uparrow(s))\neq 0
\]
for some element $f\in O_E[\Sh(\rV^\circ_{n_0},\rK^\circ_{n_0})][\Ker\phi_{\Pi_0}]\otimes_{O_E}O_E[\Sh(\rV^\circ_{n_1},\rK^\circ_{n_1})][\Ker\phi_{\Pi_1}]$.
\end{lem}

\begin{proof}
In view of Remark \ref{re:relevant}, this follows from the direction (1)$\Rightarrow$(2) of \cite{BPLZZ}*{Theorem~1.8}, together with \cite{BPLZZ}*{Remark~4.17}. Note that since our $\Pi_0$ and $\Pi_1$ are relevant representations of $\GL_{n_0}(\dA_F)$ and $\GL_{n_1}(\dA_F)$, respectively, both members in the pair of hermitian spaces in (2) of \cite{BPLZZ}*{Theorem~1.8} have to be standard definite.
\end{proof}

\begin{theorem}\label{th:selmer0}
Keep the setup in \S\ref{ss:setup}. Assume Hypothesis \ref{hy:unitary_cohomology} for both $n$ and $n+1$. If $L(\frac{1}{2},\Pi_0\times\Pi_1)\neq 0$, then for all admissible primes $\lambda$ of $E$, we have
\[
\rH^1_f(F,\rho_{\Pi_0,\lambda}\otimes_{E_\lambda}\rho_{\Pi_1,\lambda}(n))=0.
\]
\end{theorem}

\begin{proof}
By Lemma \ref{le:ggp}, we may fix the choices of $\rV^\circ_n$, $\Lambda^\circ_n$, $(\rK^\circ_n,\rK^\circ_{n+1})$ in that lemma such that
\[
\sum_{s\in\Sh(\rV^\circ_n,\rK^\circ_n)}f(s,\sh_\uparrow(s))\neq 0
\]
for some $f\in O_E[\Sh(\rV^\circ_{n_0},\rK^\circ_{n_0})][\Ker\phi_{\Pi_0}]\otimes_{O_E}O_E[\Sh(\rV^\circ_{n_1},\rK^\circ_{n_1})][\Ker\phi_{\Pi_1}]$. Moreover, by Lemma \ref{le:transferable}(3), we may assume that $(\rK^\circ_{n_0})_v$ is transferable (Definition \ref{de:transferable}) for $v\in\Sigma^+_\mnm$.

We take a prime $\lambda$ of $E$ with the underlying rational prime $\ell$. We adopt notation in \S\ref{ss:rankin_selberg} with the initial data in Definition \ref{de:lambda}. Define two nonnegative integers $m_{\r{per}}$ and $m_{\r{lat}}$ as follows.
\begin{enumerate}
  \item Let $m_{\r{per}}$ be the largest (nonnegative) integer such that
     \[
     \sum_{s\in\Sh(\rV^\circ_n,\rK^\circ_n)}f(s,\sh_\uparrow(s))\in \lambda^{m_{\r{per}}}O_E
     \]
     for every $f\in O_E[\Sh(\rV^\circ_{n_0},\rK^\circ_{n_0})][\Ker\phi_{\Pi_0}]\otimes_{O_E}O_E[\Sh(\rV^\circ_{n_1},\rK^\circ_{n_1})][\Ker\phi_{\Pi_1}]$.

  \item We choose a standard \emph{indefinite} hermitian space $\rV_{n_1}$ over $F$ of rank $n_1$, together with an identification $\rU((\rV^\circ_{n_1})^\infty)\simeq\rU(\rV_{n_1}^\infty)$ of reductive groups over $\dA^\infty_{F^+}$.\footnote{There are many choices of such $\rV_{n_1}$ and the isomorphism. We choose one only to get some control on the discrepancy of the integral cohomology of Shimura varieties and the lattice coming from Galois representations.} In particular, we have the Shimura variety $\Sh(\rV_{n_1},\rK^\circ_{n_1})$. By Hypothesis \ref{hy:unitary_cohomology}, we have an isomorphism
      \[
      \rH^{2r_1}_\et(\Sh(\rV_{n_1},\rK^\circ_{n_1})_{\ol{F}},E_\lambda(r_1))/\Ker\phi_{\Pi_1}
      \simeq(\rR_1^\tc\otimes_{O_\lambda}E_\lambda)^{\oplus\mu_1}
      \]
      of $E_\lambda[\Gamma_F]$-modules for some integer $\mu_1>0$. We fix a map
      \[
      \rH^{2r_1}_\et(\Sh(\rV_{n_1},\rK^\circ_{n_1})_{\ol{F}},O_\lambda(r_1))/\Ker\phi_{\Pi_1}
      \to(\rR_1^\tc)^{\oplus\mu_1}
      \]
      of $O_\lambda[\Gamma_F]$-modules whose kernel and cokernel are both $O_\lambda$-torsion. Then we let $m_{\r{lat}}$ be the smallest nonnegative integer such that both the kernel and the cokernel are annihilated by $\lambda^{m_{\r{lat}}}$.
\end{enumerate}
Now we assume that $\lambda$ is admissible.

We start to prove the theorem by contradiction, hence assume
\[
\dim_{E_\lambda}\rH^1_f(F,\rho_{\Pi_0,\lambda}\otimes_{E_\lambda}\rho_{\Pi_1,\lambda}(n))\geq 1.
\]
Take a sufficiently large positive integer $m$ which will be determined later. By Lemma \ref{le:purity}, we may apply Proposition \ref{pr:selmer_reduction} by taking $\Sigma$ to be the set of places of $F$ above $\Sigma^+_\mnm\cup\Sigma^+_\ell$. Then we obtain a submodule $S$ of $\rH^1_{f,\rR}(F,\bar\rR^{(m)})$ that is free of rank $1$ over $O_\lambda/\lambda^{m-m_\Sigma}$ such that $\loc_w\res_S=0$ for every nonarchimedean place $w\in\Sigma$ not above $\ell$. Now we apply the discussion in \S\ref{ss:galois_lemma} to the submodule $S\subseteq\rH^1(F,\bar\rR^{(m)})$. By (L5-1) and Lemma \ref{le:image}, we obtain an injective map
\[
\theta_S\colon\Gal(F_S/F_{\bar\rho^{(m)}})\to\Hom_{O_\lambda}(S,\bar\rR^{(m)})
\]
whose image generates an $O_\lambda$-submodule containing $\lambda^{\fr_{\bar\rR^{(m)}}}\Hom_{O_\lambda}(S,\bar\rR^{(m)})$, which further contains $\lambda^{\fr_\rR}\Hom_{O_\lambda}(S,\bar\rR^{(m)})$ by Lemma \ref{le:reducibility} and (L3). By (L5-2) and Lemma \ref{le:galois_element}, we may choose an element $(\gamma_1,\gamma_2,\xi)$ in the image of $(\bar\rho^{(m)}_{1+},\bar\rho^{(m)}_{2+},\bar\epsilon_\ell^{(m)})\res_{\Gal(\ol{F}/F^+_{\r{rflx}})}$ satisfying (a--e) in Lemma \ref{le:galois_element}. It then gives rise to an element $\gamma\in(\GL_{n_0n_1}(O_\lambda/\lambda^m)\times(O_\lambda/\lambda^m)^\times)\fc$ as in Notation \ref{no:h_gamma} such that $(\bar\rR^{(m)})^{h_\gamma}$ is a free $O_\lambda/\lambda^m$-module of rank $1$. Now we apply the discussion in \S\ref{ss:localization}. By Proposition \ref{pr:selmer_localization} (with $m_0=m_\Sigma$ and $r_S=1$), we may fix an $(S,\gamma)$-abundant element $\Psi\in G_{S,\gamma}$ (Definition \ref{de:abundant}).

We apply the discussion and notation in \S\ref{ss:first_reciprocity} to our situation with $\lambda$, $m$, $\Sigma^+_{\lr,\rI}=\emptyset$, $\Sigma^+_\rI=\Sigma^+_\mnm$, $(\rV^\circ_n,\Lambda^\circ_n)$, $\rK^\circ_n$ and $(\rK^\circ_n,\rK^\circ_{n+1})$. By the Chebotarev density theorem, we can choose a $\gamma$-associated place (Definition \ref{de:associated}) $w^{(m)}_+$ of $F^{(m)}_+$ satisfying $\Psi_{w^{(m)}}=\Psi$ and whose underlying prime $\fp$ of $F^+$ (and the underlying rational prime $p$) is a special inert prime satisfying (PI1)--(PI7) and
\begin{description}
  \item[(PI8)] the natural map
     \[
     \rH^i_\et(\Sh(\rV_{n_1},\rK^\circ_{n_1})_{\ol{F}},O_\lambda(r_1))/(\dT^{\Sigma^+_\rI\cup\Sigma^+_p}_{n_1}\cap\Ker\phi_{\Pi_1})
     \to\rH^i_\et(\Sh(\rV_{n_1},\rK^\circ_{n_1})_{\ol{F}},O_\lambda(r_1))/\Ker\phi_{\Pi_1}
     \]
     is an isomorphism for every integer $i$.
\end{description}
We also choose remaining data in \S\ref{ss:ns_initial} with $\dQ_p^\Phi=\dQ_{p^2}$, data as in Notation \ref{no:ns_uniformization}, and an indefinite uniformization datum as in Notation \ref{no:ns_uniformization_indefinite}. By the definition of $m_{\r{per}}$, we have
\begin{align}\label{eq:selmer0_5}
\exp_\lambda\(\CF_{\Sh(\rV^\circ_n,\rK^\circ_\sp)},
O_E[\Sh(\rV^\circ_{n_0},\rK^\circ_{n_0})\times\Sh(\rV^\circ_{n_1},\rK^\circ_{n_1})]/(\fn_0,\fn_1)\)\geq m-m_{\r{per}},
\end{align}
where we recall that
\[
\fn_\alpha=\dT^{\Sigma^+_\rI\cup\Sigma^+_p}_{n_\alpha}\cap\Ker\(\dT^{\Sigma^+_\mnm}_{n_\alpha}\xrightarrow{\phi_{\Pi_\alpha}}O_E\to O_E/\lambda^m\)
\]
for $\alpha=0,1$. Here, $\CF_{\Sh(\rV^\circ_n,\rK^\circ_\sp)}$ is nothing but the characteristic function of the graph $\graph\Sh(\rV^\circ_n,\rK^\circ_n)$ of the map $\Sh(\rV^\circ_n,\rK^\circ_n)\to\Sh(\rV^\circ_{n+1},\rK^\circ_{n+1})$.

We claim that there exists an element $c_1\in\rH^1(F,\bar\rR^{(m)\tc})$ satisfying \begin{align}\label{eq:selmer0_4}
\exp_\lambda\(\partial_\fp\loc_\fp(c_1),\rH^1_\sing(F_\fp,\bar\rR^{(m)\tc})\)\geq m-m_{\r{per}}-m_{\r{lat}};
\end{align}
and such that for every nonarchimedean place $w$ of $F$ not above $\Sigma^+\cup\{\fp\}$,
\begin{align}\label{eq:selmer0_6}
\loc_w(c_1)\in\rH^1_\ns(F_w,\bar\rR^{(m)\tc})
\end{align}
holds.

We first prove the theorem assuming the existence of such $c_1$. Fix a generator of the submodule $S\subseteq\rH^1_{f,\rR}(F,\bar\rR^{(m)})$ and denote by its image in $\rH^1(F,\bar\rR^{(m)})$ by $s_1$. We also identify $\bar\rR^{(m)\tc}$ with $(\bar\rR^{(m)})^*$ via the polarization $\Xi$. Now we compute the local Tate pairing $\langle s_1,c_1\rangle_w$ \eqref{eq:local_tate} for every nonarchimedean place $w$ of $F$.
\begin{itemize}[label={\ding{109}}]
  \item Suppose that $w$ is above $\Sigma^+_\mnm$. Then we have $\loc_w(s_1)=0$ by our choice of $S$. Thus, $\langle s_1,c_1\rangle_w=0$.

  \item Suppose that $w$ is above $\Sigma^+_\ell$. Then by (L2), $\rR_\dQ$ is crystalline with Hodge--Tate weights in $[-n,n-1]$. Thus, we have $\loc_w(s_1)\in\rH^1_\ns(F_w,\bar\rR^{(m)})$ by Lemma \ref{le:global_local}(2) and (L1). By \eqref{eq:selmer0_6}, Lemma \ref{le:crystalline_tate} and (L1), we have $\lambda^{m_{\r{dif}}}\langle s_1,c_1\rangle_w=0$ where $\fd_\lambda=\lambda^{m_{\r{dif}}}\subseteq O_\lambda$ is the different ideal of $E_\lambda$ over $\dQ_\ell$.

  \item Suppose that $w$ is not above $\Sigma^+_\mnm\cup\Sigma^+_\ell\cup\{\fp\}$. Then by (L2), $\rR$ is unramified. Thus, we have $\loc_w(s_1)\in\rH^1_\ns(F_w,\bar\rR^{(m)})$ by Lemma \ref{le:global_local}(1). By \eqref{eq:selmer0_6} and Lemma \ref{le:tate}, we have $\langle s_1,c_1\rangle_w=0$.

  \item Suppose that $w$ is the unique place above $\fp$. By Proposition \ref{co:selmer_localization}, we have
      \[
      \exp_\lambda\(\loc_w(s_1),\rH^1_\ns(F_w,\bar\rR^{(m)})\)\geq m-m_\Sigma-\fr_\rR.
      \]
      By \eqref{eq:selmer0_4} and Lemma \ref{le:tate} again, we have
      \[
      \exp_\lambda\(\langle s_1,c_1\rangle_w,O_\lambda/\lambda^{m}\)\geq m-m_{\r{per}}-m_{\r{lat}}-m_\Sigma-\fr_\rR.
      \]
\end{itemize}
Therefore, as long as we take $m$ such that $m>m_{\r{per}}+m_{\r{lat}}+m_\Sigma+\fr_\rR+m_{\r{dif}}$, we will have a contradiction to the relation
\[
\sum_{w}\langle s_1,c_1\rangle_w=0,
\]
where the sum is taken over all nonarchimedean places $w$ of $F$. The theorem is proved.

Now we consider the claim on the existence of $c_1$. First note that by Remark \ref{re:single_vanishing}, Assumption \ref{as:first_vanishing} is satisfied by Lemma \ref{le:ns_pbc} and (L7).

By (L4), (L6), and Theorem \ref{th:raising}(5), we have an isomorphism
\begin{align}\label{eq:selmer0_0}
\rH^{2r_0-1}_\et((\Sh(\rV'_{n_0},\tj_{n_0}\rK^{p\circ}_{n_0}\rK'_{n_0,p})_{\ol{F}},O_\lambda(r_0))/\fn_0
\xrightarrow{\sim}\(\bar\rR_0^{(m)\tc}\)^{\oplus\mu_0}
\end{align}
of $O_\lambda[\Gamma_F]$-modules, for some positive integer $\mu_0$.

By Lemma \ref{le:intertwining}, we have an isomorphism
\[
\rH^i_\et(\Sh(\rV_{n_1},\rK^\circ_{n_1})_{\ol{F}},O_\lambda)_{\fm_1}\simeq
\rH^i_\et(\Sh(\rV'_{n_1},\tj_{n_1}\rK^{p\circ}_{n_1}\rK'_{n_1,p})_{\ol{F}},O_\lambda)_{\fm_1}
\]
of $O_\lambda[\Gamma_F]$-modules. Moreover, by (PI8), we may fix a map
\[
\rH^{2r_1}_\et(\Sh(\rV'_{n_1},\tj_{n_1}\rK^{p\circ}_{n_1}\rK'_{n_1,p})_{\ol{F}},O_\lambda(r_1))/
(\dT^{\Sigma^+_\rI\cup\Sigma^+_p}_{n_1}\cap\Ker\phi_{\Pi_1})
\to\(\rR_1^\tc\)^{\oplus\mu_1}
\]
of $O_\lambda[\Gamma_F]$-modules whose kernel and cokernel are both annihilated by $\lambda^{m_{\r{lat}}}$. Taking quotient by $\lambda^m$, we obtain a map
\begin{align}\label{eq:selmer0_1}
\rH^{2r_1}_\et(\Sh(\rV'_{n_1},\tj_{n_1}\rK^{p\circ}_{n_1}\rK'_{n_1,p})_{\ol{F}},O_\lambda(r_1))/\fn_1
\to\(\bar\rR_1^{(m)\tc}\)^{\oplus\mu_1}
\end{align}
of $O_\lambda[\Gamma_F]$-modules whose kernel and cokernel are both annihilated by $\lambda^{m_{\r{lat}}}$.

To continue, we adopt the notational abbreviation prior to Corollary \ref{co:first}. By Lemma \ref{le:single_compact} and the K\"{u}nneth formula, we obtain a map
\begin{align}\label{eq:selmer0_2}
\Upsilon\colon\rH^{2n-1}_\et((\Sh'_{n_0}\times_{\Spec{F}}\Sh'_{n_1})_{\ol{F}},O_\lambda(n))/(\fn_0,\fn_1)\to
\(\bar\rR^{(m)\tc}\)^{\oplus\mu_0\mu_1}
\end{align}
of $O_\lambda[\Gamma_F]$-modules whose kernel and cokernel are both annihilated by $\lambda^{m_{\r{lat}}}$, from \eqref{eq:selmer0_0} and \eqref{eq:selmer0_1}. Recall that we have a class
\[
\AJ(\Sh'_\sp)\in\rH^1(F,\rH^{2n-1}_\et((\Sh'_{n_0}\times_{\Spec{F}}\Sh'_{n_1})_{\ol{F}},O_\lambda(n))/(\fn_0,\fn_1)),
\]
where $\Sh'_\sp$ is nothing but the graph of the morphism $\Sh'_n\to\Sh'_{n+1}$. By Corollary \ref{co:first} and \eqref{eq:selmer0_5}, we have
\begin{align}\label{eq:selmer0_3}
\exp_\lambda\(\partial_\fp\loc_\fp\AJ(\Sh'_\sp),
\rH^1_\sing(F_\fp,\rH^{2n-1}_\et((\Sh'_{n_0}\times_{\Spec{F}}\Sh'_{n_1})_{\ol{F}},O_\lambda(n))/(\fn_0,\fn_1))\)\geq m-m_{\r{per}}.
\end{align}
For every $1\leq i\leq\mu_0\mu_1$, let
\[
\Upsilon_i\colon\rH^{2n-1}_\et((\Sh'_{n_0}\times_{\Spec{F}}\Sh'_{n_1})_{\ol{F}},O_\lambda(n))/(\fn_0,\fn_1)\to\bar\rR^{(m)\tc}
\]
be the composition of $\Upsilon$ \eqref{eq:selmer0_2} with the projection to the $i$-th factor; and put
\[
c_i\coloneqq\rH^1(F,\Upsilon_i)(\AJ(\Sh'_\sp))\in\rH^1(F,\bar\rR^{(m)\tc}).
\]
Then \eqref{eq:selmer0_3} implies
\[
\max_{1\leq i\leq\mu_0\mu_1}\exp_\lambda\(\partial_\fp\loc_\fp(c_i),\rH^1_\sing(F_\fp,\bar\rR^{(m)\tc})\)\geq m-m_{\r{per}}-m_{\r{lat}}.
\]
Without loss of generality, we obtain \eqref{eq:selmer0_4}. On the other hand, as both $\Sh'_n$ and $\Sh'_{n+1}$ have smooth models over $O_{F_w}$ for which (an analogue of) Lemma \ref{le:qs_pbc} holds, we obtain \eqref{eq:selmer0_6}.
\end{proof}

Now we deduce two concrete consequences from Theorem \ref{th:selmer0}.

\begin{corollary}\label{co:sym}
Let $n\geq 2$ be an integer and denote by $n_0$ and $n_1$ the unique even and odd numbers in $\{n,n+1\}$, respectively. Let $A_0$ and $A_1$ be two modular elliptic curves over $F^+$ such that $\End({A_0}_{\ol{F}})=\End({A_1}_{\ol{F}})=\dZ$. Suppose that
\begin{enumerate}[label=(\alph*)]
  \item ${A_0}_{\ol{F}}$ and ${A_1}_{\ol{F}}$ are not isogenous to each other;

  \item both $\Sym^{n_0-1}A_0$ and $\Sym^{n_1-1}A_1$ are modular; and

  \item $F^+\neq\dQ$ if $n\geq 3$.
\end{enumerate}
If the (central critical) $L$-value $L(n,\Sym^{n_0-1}{A_0}_F\times\Sym^{n_1-1}{A_1}_F)$ does not vanish, then we have
\[
\rH^1_f(F,\Sym^{n_0-1}\rH^1_\et({A_0}_{\ol{F}},\dQ_\ell)\otimes_{\dQ_\ell}\Sym^{n_1-1}\rH^1_\et({A_1}_{\ol{F}},\dQ_\ell)(n))=0
\]
for all but finitely many rational primes $\ell$.
\end{corollary}

\begin{proof}
By (b) and \cite{AC89}, both $\Sym^{n_0-1}{A_0}_F$ and $\Sym^{n_1-1}{A_1}_F$ are modular. Thus, we may let $\Pi_\alpha$ be the (cuspidal) automorphic representation of $\GL_{n_\alpha}(\dA_F)$ associated to $\Sym^{n_\alpha-1}{A_\alpha}_F$ for $\alpha=0,1$, which is a relevant representation (Definition \ref{de:relevant}). We also have the identity
\[
L(n+s,\Sym^{n_0-1}{A_0}_F\times\Sym^{n_1-1}{A_1}_F)=L(\tfrac{1}{2}+s,\Pi_0\times\Pi_1)
\]
of $L$-functions, and that the representation of $\Gamma_F$ on $\Sym^{n_\alpha-1}\rH^1_\et({A_\alpha}_{\ol{F}},\dQ_\ell)$ is isomorphic to $\rho_{\Pi_\alpha,\ell}$ for $\alpha=0,1$. By Proposition \ref{pr:unitary_cohomology} and (c), Hypothesis \ref{hy:unitary_cohomology} is known in this case. Then the corollary follows immediately from Theorem \ref{th:selmer0} and Lemma \ref{pr:elliptic_admissible} (where we use (a) and (c)) with $E=\dQ$.
\end{proof}

\begin{remark}\label{re:sym}
In this remark, we summarize the current knowledge on the modularity of symmetric powers of elliptic curves, namely, condition (a) in Corollary \ref{co:sym}. Let $A$ be a modular elliptic curve over $F^+$ such that $\End(A_{\ol{F}})=\dZ$. We have
\begin{itemize}[label={\ding{109}}]
  \item $\Sym^2A$ is modular by \cite{GJ76};

  \item $\Sym^3A$ is modular by \cite{KS02};

  \item $\Sym^4A$ is modular by \cite{Kim03};

  \item $\Sym^5A$ and $\Sym^6A$ are modular if $F^+$ is linearly disjoint from $\dQ(\zeta_5)$ over $\dQ$;

  \item $\Sym^7A$ is modular if $F^+$ is linearly disjoint from $\dQ(\zeta_{35})$ over $\dQ$;

  \item $\Sym^8A$ is modular if $F^+$ is linearly disjoint from $\dQ(\zeta_7)$ over $\dQ$;
\end{itemize}
in which the last three cases are obtained in a series of recent work \cites{CT1,CT2,CT3} of Clozel and Thorne.

After we completed this article, we have learnt the groundbreaking result of Newton--Thorne \cites{NT1,NT2} where they prove the modularity of all symmetric powers of elliptic curves over $\dQ$ without complex multiplication. In particular, it follows that $\Sym^n A$ is modular if $F^+/\dQ$ is solvable and $A$ is the base change of an elliptic curve over $\dQ$.
\end{remark}

\begin{corollary}\label{co:abstract}
Keep the setup in \S\ref{ss:setup}. Suppose that
\begin{enumerate}[label=(\alph*)]
  \item there exists a very special inert prime $\fp$ of $F^+$ (Definition \ref{de:special_inert}) such that $\Pi_{0,\fp}$ is Steinberg, and $\Pi_{1,\fp}$ is unramified whose Satake parameter contains $1$ exactly once;

  \item for $\alpha=0,1$, there exists a nonarchimedean place $w_\alpha$ of $F$ such that $\Pi_{\alpha,w_\alpha}$ is supercuspidal; and

  \item $F^+\neq\dQ$ if $n\geq 3$.
\end{enumerate}
If $L(\frac{1}{2},\Pi_0\times\Pi_1)\neq 0$, then for all but finitely many primes $\lambda$ of $E$, we have
\[
\rH^1_f(F,\rho_{\Pi_0,\lambda}\otimes_{E_\lambda}\rho_{\Pi_1,\lambda}(n))=0.
\]
\end{corollary}

\begin{proof}
This follows from Theorem \ref{th:selmer0} and Lemma \ref{le:abstract}.
\end{proof}

\subsection{Main theorem in the Selmer rank 1 case}
\label{ss:main_1}

We state the following weak version of the arithmetic Gan--Gross--Prasad conjecture.

\begin{conjecture}\label{co:aggp}
Suppose that $L(\frac{1}{2},\Pi_0\times\Pi_1)=0$ but $L'(\frac{1}{2},\Pi_0\times\Pi_1)\neq 0$. Then there exist
\begin{itemize}[label={\ding{109}}]
  \item a standard indefinite hermitian space $\rV_n$ of rank $n$ over $F$, together with a self-dual $\prod_{v\not\in\Sigma^+_\infty\cup\Sigma^+_\mnm}O_{F_v}$-lattice $\Lambda_n$ in $\rV_n\otimes_F\dA_F^{\Sigma^+_\infty\cup\Sigma^+_\mnm}$ (and put $\rV_{n+1}\coloneqq(\rV_n)_\sharp$ and $\Lambda_{n+1}\coloneqq(\Lambda_n)_\sharp$),

  \item an object $(\rK_n,\rK_{n+1})\in\fK(\rV_n)_\sp$ in which $\rK_{n_\alpha}$ is of the form
      \[
      \rK_{n_\alpha}=\prod_{v\in\Sigma^+_\mnm}(\rK_{n_\alpha})_v\times
      \prod_{v\not\in\Sigma^+_\infty\cup\Sigma^+_\mnm}\rU(\Lambda_{n_\alpha})(O_{F^+_v})
      \]
      for $\alpha=0,1$,
\end{itemize}
such that for every prime $\lambda$ of $E$, the graph $\graph\Sh(\rV_n,\rK_n)$ of the morphism $\sh_\uparrow\colon\Sh(\rV_n,\rK_n)\to\Sh(\rV_{n+1},\rK_{n+1})$ \eqref{eq:qs_functoriality_shimura} is nonvanishing in the quotient Chow group
\[
\CH^n(\Sh(\rV_{n_0},\rK_{n_0})\times_{\Spec{F}}\Sh(\rV_{n_1},\rK_{n_1}))_E/(\Ker\phi_{\Pi_0},\Ker\phi_{\Pi_1}).
\]
\end{conjecture}

In the situation of the above conjecture, since both $\Pi_0$ and $\Pi_1$ are cuspidal, we have
\[
\rH^i_\et((\Sh(\rV_{n_0},\rK_{n_0})\times_{\Spec{F}}\Sh(\rV_{n_1},\rK_{n_1}))_{\ol{F}},E_\lambda)/(\Ker\phi_{\Pi_0},\Ker\phi_{\Pi_1})=0
\]
if $i\neq 2n-1$. In particular, the Hochschild--Serre spectral sequence gives rise to a coboundary map
\begin{align*}
\AJ^{\Pi_0,\Pi_1}_\lambda&\colon\rZ^n(\Sh(\rV_{n_0},\rK_{n_0})\times_{\Spec{F}}\Sh(\rV_{n_1},\rK_{n_1})) \\
&\to\rH^1(F,\rH^{2n-1}_\et((\Sh(\rV_{n_0},\rK_{n_0})\times_{\Spec{F}}\Sh(\rV_{n_1},\rK_{n_1}))_{\ol{F}},E_\lambda(n))
/(\Ker\phi_{\Pi_0},\Ker\phi_{\Pi_1})).
\end{align*}

\begin{theorem}\label{th:selmer1}
Keep the setup in \S\ref{ss:setup}. Assume Hypothesis \ref{hy:unitary_cohomology} for both $n$ and $n+1$. Let $\lambda$ be a prime of $E$ for which there exist
\begin{itemize}[label={\ding{109}}]
  \item a standard indefinite hermitian space $\rV_n$ of rank $n$ over $F$, together with a self-dual $\prod_{v\not\in\Sigma^+_\infty\cup\Sigma^+_\mnm}O_{F_v}$-lattice $\Lambda_n$ in $\rV_n\otimes_F\dA_F^{\Sigma^+_\infty\cup\Sigma^+_\mnm}$ (and put $\rV_{n+1}\coloneqq(\rV_n)_\sharp$ and $\Lambda_{n+1}\coloneqq(\Lambda_n)_\sharp$),

  \item an object $(\rK_n,\rK_{n+1})\in\fK(\rV_n)_\sp$ in which $\rK_{n_\alpha}$ is of the form
      \[
      \rK_{n_\alpha}=\prod_{v\in\Sigma^+_\mnm}(\rK_{n_\alpha})_v\times
      \prod_{v\not\in\Sigma^+_\infty\cup\Sigma^+_\mnm}\rU(\Lambda_{n_\alpha})(O_{F^+_v})
      \]
      for $\alpha=0,1$, satisfying that $(\rK_{n_0})_v$ is a transferable open compact subgroup (Definition \ref{de:transferable}) of $\rU(\rV^\circ_{n_0})(F^+_v)$ for $v\in\Sigma^+_\mnm$,
\end{itemize}
such that
\begin{align}\label{eq:selmer1}
\AJ^{\Pi_0,\Pi_1}_\lambda(\graph\Sh(\rV_n,\rK_n))\neq 0.
\end{align}
If $\lambda$ is admissible, then we have
\[
\dim_{E_\lambda}\rH^1_f(F,\rho_{\Pi_0,\lambda}\otimes_{E_\lambda}\rho_{\Pi_1,\lambda}(n))=1.
\]
\end{theorem}

\begin{remark}
In fact, \eqref{eq:selmer1} already implies that the global epsilon factor of $\Pi_0\times\Pi_1$ is $-1$.
\end{remark}

\begin{proof}[Proof of Theorem \ref{th:selmer1}]
We take an admissible prime $\lambda$ of $E$ for which we may choose data $\rV_n$, $\Lambda_n$, $(\rK_n,\rK_{n+1})$ as in the statement of the theorem such that $\AJ^{\Pi_0,\Pi_1}_\lambda(\graph\Sh(\rV_n,\rK_n))\neq 0$. Lemma \ref{le:purity} and (L2) imply that $\AJ^{\Pi_0,\Pi_1}_\lambda(\graph\Sh(\rV_n,\rK_n))$ belongs to the subspace
\[
\rH^1_f(F,\rH^{2n-1}_\et((\Sh(\rV_{n_0},\rK_{n_0})
\times_{\Spec{F}}\Sh(\rV_{n_1},\rK_{n_1}))_{\ol{F}},E_\lambda(n))/(\Ker\phi_{\Pi_0},\Ker\phi_{\Pi_1}))
\]
and hence to the submodule
\[
\rH^1_f(F,\rH^{2n-1}_\et((\Sh(\rV_{n_0},\rK_{n_0})
\times_{\Spec{F}}\Sh(\rV_{n_1},\rK_{n_1}))_{\ol{F}},O_\lambda(n))/(\Ker\phi_{\Pi_0},\Ker\phi_{\Pi_1}))
\]
by Definition \ref{de:bk_integral}.

We adopt notation in \S\ref{ss:rankin_selberg} with the initial data in Definition \ref{de:lambda}. Define two nonnegative integers $m_{\r{per}}$ and $m_{\r{lat}}$ as follows.
\begin{enumerate}
  \item By Hypothesis \ref{hy:unitary_cohomology}, we may choose a map
     \[
     \rH^{2n-1}_\et((\Sh(\rV_{n_0},\rK_{n_0})
     \times_{\Spec{F}}\Sh(\rV_{n_1},\rK_{n_1}))_{\ol{F}},O_\lambda(n))/(\Ker\phi_{\Pi_0},\Ker\phi_{\Pi_1})\to\rR^\tc
     \]
     of $O_\lambda[\Gamma_F]$-modules such that the induced image of $\AJ^{\Pi_0,\Pi_1}_\lambda(\graph\Sh(\rV_n,\rK_n))$ in $\rH^1_f(F,\rR^\tc)$, denoted by $s^\tc$, is non-torsion. Let $s\in\rH^1_f(F,\rR)$ be the element corresponding to $s^\tc$ under the isomorphism in Lemma \ref{le:conjugate}. We put
     \[
     m_{\r{per}}\coloneqq\ord_\lambda\(s,\rH^1_f(F,\rR)/\rH^1_f(F,\rR)_\tor\)
     \]
     (Definition \ref{de:divisibility}), which is a nonnegative integer.

  \item By Hypothesis \ref{hy:unitary_cohomology}, we have an isomorphism
      \[
      \rH^{2r_1}_\et(\Sh(\rV_{n_1},\rK_{n_1})_{\ol{F}},E_\lambda(r_1))/\Ker\phi_{\Pi_1}
      \simeq(\rR_1^\tc\otimes_{O_\lambda}E_\lambda)^{\oplus\mu_1}
      \]
      of $E_\lambda[\Gamma_F]$-modules for some integer $\mu_1>0$. We fix a map
      \[
      \rH^{2r_1}_\et(\Sh(\rV_{n_1},\rK_{n_1})_{\ol{F}},O_\lambda(r_1))/\Ker\phi_{\Pi_1}\to(\rR_1^\tc)^{\oplus\mu_1}
      \]
      of $O_\lambda[\Gamma_F]$-modules whose kernel and cokernel are both $O_\lambda$-torsion. Then we let $m_{\r{lat}}$ be the smallest nonnegative integer such that both the kernel and the cokernel are annihilated by $\lambda^{m_{\r{lat}}}$.
\end{enumerate}
Note that in (1), we obtain an element $s\in\rH^1_f(F,\rR)_\dQ=\rH^1_f(F,\rR_\dQ)=\rH^1_f(F,\rho_{\Pi_0,\lambda}\otimes_{E_\lambda}\rho_{\Pi_1,\lambda}(n))$ that is nonzero. In particular, we have $\dim_{E_\lambda}\rH^1_f(F,\rho_{\Pi_0,\lambda}\otimes_{E_\lambda}\rho_{\Pi_1,\lambda}(n))\geq 1$.

We start to prove the theorem by contradiction, hence assume
\[
\dim_{E_\lambda}\rH^1_f(F,\rho_{\Pi_0,\lambda}\otimes_{E_\lambda}\rho_{\Pi_1,\lambda}(n))\geq 2.
\]
Take a sufficiently large positive integer $m$ which will be determined later. We fix a uniformizer $\lambda_0$ of $E_\lambda$. By Lemma \ref{le:purity}, we may apply Proposition \ref{pr:selmer_reduction} by taking $\Sigma$ to be the set of places of $F$ above $\Sigma^+_\mnm\cup\Sigma^+_\ell$. Then we obtain a submodule $S$ of $\rH^1_{f,\rR}(F,\bar\rR^{(m)})$ containing (the image of) $\lambda_0^{m_\Sigma-m_{\r{per}}}s$ of order $0$,\footnote{Here, $\lambda_0^{-m_{\r{per}}}s$ is any element in $\rH^1_f(F,\rR)$ satisfying $\lambda_0^{m_{\r{per}}}(\lambda_0^{-m_{\r{per}}}s)=s$.} that is free of rank $2$ over $O_\lambda/\lambda^{m-m_\Sigma}$, and such that $\loc_w\res_S=0$ for every nonarchimedean place $w\in\Sigma$ not above $\ell$. Now we apply the discussion in \S\ref{ss:galois_lemma} to the submodule $S\subseteq\rH^1(F,\bar\rR^{(m)})$. By (L5-1) and Lemma \ref{le:image}, we obtain an injective map
\[
\theta_S\colon\Gal(F_S/F_{\bar\rho^{(m)}})\to\Hom_{O_\lambda}(S,\bar\rR^{(m)})
\]
whose image generates an $O_\lambda$-submodule containing $\lambda^{4\fr_{\bar\rR^{(m)}}}\Hom_{O_\lambda}(S,\bar\rR^{(m)})$, which further contains $\lambda^{4\fr_\rR}\Hom_{O_\lambda}(S,\bar\rR^{(m)})$ by Lemma \ref{le:reducibility} and (L3). By (L5-2) and Lemma \ref{le:galois_element}, we may choose an element $(\gamma_1,\gamma_2,\xi)$ in the image of $(\bar\rho^{(m)}_{1+},\bar\rho^{(m)}_{2+},\bar\epsilon_\ell^{(m)})\res_{\Gal(\ol{F}/F^+_{\r{rflx}})}$ satisfying (a--e) in Lemma \ref{le:galois_element}. It then gives rise to an element $\gamma\in(\GL_{n_0n_1}(O_\lambda/\lambda^m)\times(O_\lambda/\lambda^m)^\times)\fc$ as in Notation \ref{no:h_gamma} such that $(\bar\rR^{(m)})^{h_\gamma}$ is a free $O_\lambda/\lambda^m$-module of rank $1$. Now we apply the discussion in \S\ref{ss:localization}. By Proposition \ref{pr:selmer_localization} (with $m_0=m_\Sigma$ and $r_S=2$), we may fix an $(S,\gamma)$-abundant pair $(\Psi_1,\Psi_2)\in G_{S,\gamma}^2$ (Definition \ref{de:abundant}). By Proposition \ref{co:selmer_localization}, we may choose a basis $\{s_1,s_2\}$ of $S$ such that $\theta_S(\Psi_1)(s_2)=\theta_S(\Psi_2)(s_1)=0$, and
\begin{align}\label{eq:selmer1_1}
\exp_\lambda\(\theta_S(\Psi_j)(s_j),(\bar\rR^{(m)})^{h_\gamma}\)\geq m-m_\Sigma-4\fr_\rR
\end{align}
for $j=1,2$. Moreover, without loss of generality, we may assume $\lambda_0^{m_\Sigma-m_{\r{per}}}s=a_1s_1+a_2s_2$ in which $a_1\in O_\lambda^\times$.

First, we apply the discussion and notation in \S\ref{ss:second_reciprocity} to our situation with $\lambda$, $m$, $\Sigma^+_{\lr,\r{II}}=\emptyset$, $\Sigma^+_{\r{II}}=\Sigma^+_\r{min}$, $(\rV_n,\Lambda_n)$, $\rK_n$ and $(\rK_n,\rK_{n+1})$. By the Chebotarev density theorem, we can choose a $\gamma$-associated place (Definition \ref{de:associated}) $w^{(m)}_{1+}$ of $F^{(m)}_+$ satisfying $\Psi_{w^{(m)}_1}=\Psi_1$ and whose underlying prime $\fp_1$ of $F^+$ (and the underlying rational prime $p_1$) is a special inert prime satisfying (PII1)--(PII7) and
\begin{description}
  \item[(PII8)] the natural map
     \[
     \rH^i_\et(\Sh(\rV_{n_1},\rK_{n_1})_{\ol{F}},O_\lambda(r_1))/(\dT^{\Sigma^+_{\r{II}}\cup\Sigma^+_{p_1}}_{n_1}\cap\Ker\phi_{\Pi_1})\to
     \rH^i_\et(\Sh(\rV_{n_1},\rK_{n_1})_{\ol{F}},O_\lambda(r_1))/\Ker\phi_{\Pi_1}
     \]
     is an isomorphism for every integer $i$.
\end{description}
We also choose remaining data in \S\ref{ss:qs_initial} with $\dQ_{p_1}^\Phi=\dQ_{p_1^2}$, a definite uniformization datum $(\rV^\star_{n_\alpha},\ti_{n_\alpha},\{\Lambda^\star_{n_\alpha,\fq}\}_{\fq\mid p_1})$ for $\alpha=0,1$ as in Notation \ref{no:qs_uniformization}. By \eqref{eq:selmer1_1} and our choice of $S$, we have
\[
\exp_\lambda\(s,\rH^1_\ns(F_{w_1},\bar\rR^{(m)})\)\geq m-m_{\r{per}}-4\fr_\rR,
\]
which implies that
\[
\resizebox{\hsize}{!}{
\xymatrix{
\exp_\lambda\(\loc_{\fp_1}([\graph\Sh(\rV_n,\rK_n)]),
\rH^{2n}_\et((\Sh(\rV_{n_0},\rK_{n_0})\times_{\Spec{F}}\Sh(\rV_{n_1},\rK_{n_1}))_{F_{\fp_1}},L(n))/(\fn_0,\fn_1)\)\geq m-m_{\r{per}}-4\fr_\rR.
}
}
\]
Here, we recall that
\[
\fn_\alpha=\dT^{\Sigma^+_{\r{II}}\cup\Sigma^+_{p_1}}_{n_\alpha}\cap\Ker\(\dT^{\Sigma^+_\mnm}_{n_\alpha}\xrightarrow{\phi_{\Pi_\alpha}}O_E\to O_E/\lambda^m\)
\]
for $\alpha=0,1$. Note that, similar to Remark \ref{re:single_vanishing}, Assumption \ref{as:second_vanishing} is satisfied by Lemma \ref{le:qs_pbc} and (L7). Thus, we may apply Theorem \ref{th:second}, hence obtain
\begin{align}\label{eq:selmer1_0}
\exp_\lambda\(\CF_{\Sh(\rV^\star_n,\rK^\star_\sp)},
O_E[\Sh(\rV^\star_{n_0},\rK^\star_{n_0})\times\Sh(\rV^\star_{n_1},\rK^\star_{n_1})]/(\fn_0,\fn_1)\)\geq m-m_{\r{per}}-4\fr_\rR.
\end{align}

Second, we apply the discussion and notation in \S\ref{ss:first_reciprocity} to our situation with $\lambda$, $m$, $\Sigma^+_{\lr,\rI}=\{\fp_1\}$, $\Sigma^+_\rI=\Sigma^+_\mnm\cup\Sigma^+_{p_1}$, $\rV^\circ_n=\rV^\star_n$, $\rK^\circ_n=\rK^\star_n$ and $(\rK^\circ_\sp,\rK^\circ_{n+1})=(\rK^\star_\sp,\rK^\star_{n+1})$. By the Chebotarev density theorem, we can choose a $\gamma$-associated place $w^{(m)}_{2+}$ of $F^{(m)}_+$ satisfying $\Psi_{w^{(m)}_2}=\Psi_2$ and whose underlying prime $\fp_2$ of $F^+$ (and the underlying rational prime $p_2$) is a special inert prime satisfying (PI1)--(PI7), $p_2\neq p_1$, and
\begin{description}
  \item[(PI8)] the natural map
     \[
     \rH^{2r_1}_\et(\Sh(\rV_{n_1},\rK_{n_1})_{\ol{F}},O_\lambda(r_1))
     /(\dT^{\Sigma^+_\rI\cup\Sigma^+_{p_2}}_{n_1}\cap\Ker\phi_{\Pi_1})\to
     \rH^{2r_1}_\et(\Sh(\rV_{n_1},\rK_{n_1})_{\ol{F}},O_\lambda(r_1))/\Ker\phi_{\Pi_1}
     \]
     is an isomorphism.
\end{description}

We claim that there exists an element $c_2\in\rH^1(F,\bar\rR^{(m)\tc})$ satisfying \begin{align}\label{eq:selmer1_4}
\exp_\lambda\(\partial_{\fp_2}\loc_{\fp_2}(c_2),\rH^1_\sing(F_{\fp_2},\bar\rR^{(m)\tc})\)\geq m-m_{\r{per}}-4\fr_\rR-m_{\r{lat}};
\end{align}
and such that for every nonarchimedean place $w$ of $F$ not above $\Sigma^+\cup\{\fp_1,\fp_2\}$,
\begin{align}\label{eq:selmer1_6}
\loc_w(c_2)\in\rH^1_\ns(F_w,\bar\rR^{(m)\tc})
\end{align}
holds.

By Remark \ref{re:qs_tate_hecke} and Remark \ref{re:qs_uniformization}, we know that there exists an isomorphism $\rU((\rV^\circ_{n_1})^\infty)\simeq\rU(\rV_{n_1}^\infty)$ sending $\rK^\circ_{n_1}$ to $\rK_{n_1}$. Then the claim can be proved by the exactly same argument for the parallel claim in the proof of Theorem \ref{th:selmer0}, using \eqref{eq:selmer1_0} and the fact that $\bar\rho_{\Pi_0,\lambda,+}$ is rigid for $(\Sigma^+_\mnm,\Sigma^+_{\lr,\rI})$.\footnote{In fact, one needs to use the additional fact that when $F^+\neq\dQ$, both Shimura varieties $\Sh'_{n_0}$ and $\Sh'_{n_1}$ have proper smooth reduction at every place $w$ of $F$ above $\Sigma^+_{p_1}\setminus\{\fp_1\}$. See Remark \ref{re:ns_smooth}.}

Now we deduce a contradiction. Replace $s_2$ by its image in $\rH^1_f(F,\bar\rR^{(m)})$. We also identify $\bar\rR^{(m)\tc}$ with $(\bar\rR^{(m)})^*$ via the polarization $\Xi$. Now we compute the local Tate pairing $\langle s_2,c_2\rangle_w$ \eqref{eq:local_tate} for every nonarchimedean place $w$ of $F$.
\begin{itemize}[label={\ding{109}}]
  \item Suppose that $w$ is above $\Sigma^+_\mnm$. Then we have $\loc_w(s_2)=0$ by our choice of $S$. Thus, $\langle s_2,c_2\rangle_w=0$.

  \item Suppose that $w$ is above $\Sigma^+_\ell$. Then by (L2), $\rR_\dQ$ is crystalline with Hodge--Tate weights in $[1-n,n]$. Thus, we have $\loc_w(s_2)\in\rH^1_\ns(F_w,\bar\rR^{(m)})$ by Lemma \ref{le:global_local}(2) and (L1). By \eqref{eq:selmer1_6}, Lemma \ref{le:crystalline_tate} and (L1), we have $\lambda^{m_{\r{dif}}}\langle s_2,c_2\rangle_w=0$ where $\fd_\lambda=\lambda^{m_{\r{dif}}}\subseteq O_\lambda$ is the different ideal of $E_\lambda$ over $\dQ_\ell$.

  \item Suppose that $w$ is not above $\Sigma^+_\mnm\cup\Sigma^+_\ell\cup\{\fp_1,\fp_2\}$. Then by (L2), $\rR$ is unramified. Thus, we have $\loc_w(s_2)\in\rH^1_\ns(F_w,\bar\rR^{(m)})$ by Lemma \ref{le:global_local}(1). By \eqref{eq:selmer1_6} and Lemma \ref{le:tate}, we have $\langle s_2,c_2\rangle_w=0$.

  \item Suppose that $w$ is the unique place above $\fp_1$. Then we have $\loc_w(s_2)=0$ by Proposition \ref{co:selmer_localization}. Thus, we have $\langle s_2,c_2\rangle_w=0$.

  \item Suppose that $w$ is the unique place above $\fp_2$. Then by Proposition \ref{co:selmer_localization}, we have
      \[
      \exp_\lambda\(\loc_w(s_2),\rH^1_\ns(F_w,\bar\rR^{(m)})\)\geq m-m_\Sigma-4\fr_\rR.
      \]
      By \eqref{eq:selmer1_4} and Lemma \ref{le:tate} again, we have
      \[
      \exp_\lambda\(\langle s_2,c_2\rangle_w,O_\lambda/\lambda^m\)\geq m-m_{\r{per}}-m_{\r{lat}}-m_\Sigma-8\fr_\rR.
      \]
\end{itemize}
Therefore, as long as we take $m$ such that $m>m_{\r{per}}+m_{\r{lat}}+m_\Sigma+8\fr_\rR+m_{\r{dif}}$, we will have a contradiction to the relation
\[
\sum_{w}\langle s_2,c_2\rangle_w=0,
\]
where the sum is taken over all nonarchimedean places $w$ of $F$. The theorem is proved.
\end{proof}

We also have an analogue of Corollary \ref{co:abstract} in the rank $1$ case, which we leave to the readers to formulate.

\if false

\section{Proof of main theorems}
\label{ss:8}

In the section, we prove our main theorems on bounding Selmer groups. In \S\ref{ss:admissible_prime}, we introduce the notation of (weakly) admissible primes for the coefficient field, and make some additional preparation for the main theorem. In \S\ref{ss:main_0} and \S\ref{ss:main_1}, we prove our main theorems in the (Selmer) rank $0$ and $1$ cases, respectively.

\subsection{Admissible primes for coefficient fields}
\label{ss:admissible_prime}

We keep the setup in \S\ref{ss:setup}.

\begin{definition}\label{de:lambda}
We introduce the following assumptions on a prime $\lambda$ of $E$ with the underlying rational prime $\ell$ (and the ring of integers $O_\lambda$ of $E_\lambda$):
\begin{description}
  \item[(L1)] $\ell>4n$ and is unramified in $F$;

  \item[(L2)] $\Sigma^+_\mnm$ does not contain $\ell$-adic places;

  \item[(L3)] the Galois representation $\rho_{\Pi_0,\lambda}\otimes_{E_\lambda}\rho_{\Pi_1,\lambda}$ is absolutely irreducible;

  \item[(L4)] Assumption \ref{as:first_irreducible} is satisfied, that is, both $\rho_{\Pi_0,\lambda}$ and $\rho_{\Pi_1,\lambda}$ are residually absolutely irreducible;

  \item[(L5)] under (L4), for $\alpha=0,1$, we have a $\Gamma_F$-stable $O_\lambda$-lattice $\rR_\alpha$ in $\rho_{\Pi_\alpha,\lambda}(r_\alpha)$, unique up to homothety, that is $(1-\alpha)$-polarizable, for which we choose a $(1-\alpha)$-polarization $\Xi_\alpha\colon\rR_\alpha^\tc\xrightarrow{\sim}\rR_\alpha^\vee(1-\alpha)$ and an isomorphism $\rR_\alpha\simeq O_\lambda^{\oplus n_\alpha}$ of $O_\lambda$-modules.\footnote{In fact, (L5) does not depend on the choice of $\Xi_\alpha$ and the basis, since $\Xi_\alpha$ is unique up to units in $O_\lambda$ and the basis is unique up to conjugation in $\GL_{n_\alpha}(O_\lambda)$.} After adopting the notation in \S\ref{ss:rankin_selberg}, we have
      \begin{description}
        \item[(L5-1)] either one of the two assumptions in Lemma \ref{le:restriction} is satisfied;

        \item[(L5-2)] $(\r{GI}_{F',\sP}^1)$ from Lemma \ref{le:galois_element} holds with $F'=F^+_{\r{rflx}}$ (Definition \ref{de:reflexive}) and $\sP(T)=T^2-1$ (see Remark \ref{re:lambda} below for a more explicit description);
      \end{description}

  \item[(L6)] under (L4), the homomorphism $\bar\rho_{\Pi_0,\lambda,+}$ (Remark \ref{re:single_irreducible}) is rigid for $(\Sigma^+_\mnm,\emptyset)$ (Definition \ref{de:rigid}), and $\bar\rho_{\Pi_0,\lambda}\res_{\Gal(\ol{F}/F(\zeta_\ell))}$ is absolutely irreducible;

  \item[(L7)] for $\alpha=0,1$, the composite homomorphism $\dT^{\Sigma^+_\mnm}_{n_\alpha}\xrightarrow{\phi_{\Pi_\alpha}}O_E\to O_E/\lambda$ is cohomologically generic (Definition \ref{de:generic}).
\end{description}
Finally, we say that
\begin{enumerate}
  \item $\lambda$ is \emph{weakly admissible} (with respect to $(\Pi_0,\Pi_1)$) if (L1--L5) are satisfied;

  \item $\lambda$ is \emph{admissible} (with respect to $(\Pi_0,\Pi_1)$) if (L1--L7) are satisfied.
\end{enumerate}
\end{definition}

\begin{remark}\label{re:lambda}
In Definition \ref{de:lambda}, (L5-2) is equivalent to the following assertion: the image of the restriction of the homomorphism
\[
(\bar\rho_{0+},\bar\rho_{1+},\bar\epsilon_\ell)\colon\Gamma_{F^+}\to
\sG_{n_0}(O_\lambda/\lambda)\times\sG_{n_1}(O_\lambda/\lambda)\times(O_\lambda/\lambda)^\times
\]
(see Notation \ref{no:reduction} for the notation) to $\Gal(\ol{F}/F^+_{\r{rflx}})$ contains an element $(\gamma_0,\gamma_1,\xi)$ satisfying
\begin{enumerate}[label=(\alph*)]
  \item $\xi^2-1\neq 0$;

  \item for $\alpha=0,1$, $\gamma_\alpha$ belongs to $(\GL_{n_\alpha}(O_\lambda/\lambda)\times(O_\lambda/\lambda)^\times)\fc$ with order coprime to $\ell$;

  \item $1$ appears in the eigenvalues of each of $h_{\gamma_0}$, $h_{\gamma_1}$, and $h_{\gamma_0}\otimes h_{\gamma_1}$ (Notation \ref{no:h_gamma}) with multiplicity one;

  \item $h_{\gamma_0}$ does not have an eigenvalue that is equal to $-1$ in $O_\lambda/\lambda$;

  \item $h_{\gamma_1}$ does not have an eigenvalue that is equal to $-\xi$ in $O_\lambda/\lambda$.
\end{enumerate}
\end{remark}

\begin{lem}\label{pr:elliptic_admissible}
Suppose that $E=\dQ$ and that there are two elliptic curves $A_0$ and $A_1$ over $F^+$ such that for every rational prime $\ell$ of $E$ and $\alpha=0,1$, we have $\rho_{\Pi_\alpha,\ell}\simeq\Sym^{n_\alpha-1}\rH^1_\et({A_\alpha}_{\ol{F}},\dQ_\ell)\res_{\Gamma_F}$. If ${A_0}_{\ol{F}}$ and ${A_1}_{\ol{F}}$ are not isogenous to each other and $\End({A_0}_{\ol{F}})=\End({A_1}_{\ol{F}})=\dZ$, then all but finitely many rational primes $\ell$ are weakly admissible; and when $F^+\neq\dQ$, all but finitely many rational primes $\ell$ are admissible.
\end{lem}

\begin{proof}
We need to show that every condition in Definition \ref{de:lambda} excludes only finitely many $\ell$ (for (L7) we assume $F^+\neq\dQ$). By \cite{Ser72}*{Th\'{e}or\`{e}me~6}, for sufficiently large $\ell$, the homomorphisms
\[
\Gamma_{F^+}\to\GL(\rH^1_\et({A_\alpha}_{\ol{F}},\dF_\ell))\simeq\GL_2(\dF_\ell)
\]
are both surjective for $\alpha=0,1$. Thus, we may assume that this is the case.

For (L1) and (L2), this is trivial.

For (L3), (L4), and (L5), this has been proved in Proposition \ref{pr:elliptic_curve}.

For (L6), by \cite{LTXZZ}*{Corollary~4.1.2}, the condition that $\bar\rho_{\Pi_0,\lambda,+}$ is rigid for $(\Sigma^+_\mnm,\emptyset)$ excludes only finitely many $\ell$. It is clear that the remaining two conditions also exclude only finitely many $\ell$.

For (L7), this follows from Corollary \ref{co:generic}.
\end{proof}

\begin{lem}\label{le:abstract}
Keep the setup in \S\ref{ss:setup}. Suppose that
\begin{enumerate}[label=(\alph*)]
  \item there exists a very special inert prime $\fp$ of $F^+$ (Definition \ref{de:special_inert}) such that $\Pi_{0,\fp}$ is Steinberg, and $\Pi_{1,\fp}$ is unramified whose Satake parameter contains $1$ exactly once; and

  \item for $\alpha=0,1$, there exists a nonarchimedean place $w_\alpha$ of $F$ such that $\Pi_{\alpha,w_\alpha}$ is supercuspidal.
\end{enumerate}
Then all but finitely many primes $\lambda$ of $E$ are weakly admissible; and when $F^+\neq\dQ$, all but finitely many primes $\lambda$ of $E$ are admissible.
\end{lem}

\begin{proof}
It suffices to show that all but finitely many primes $\lambda$ of $E$ are admissible (or weakly admissible if $n=2)$. It suffices to show that each of conditions (L1--L6) in Definition \ref{de:lambda} excludes only finitely many $\lambda$, and also for (L7) if $F^+\neq\dQ$.

For (L1) and (L2), this is trivial.

For (L4), this follows from \cite{LTXZZ}*{Proposition~4.2.3(1)} by (b).

For (L3), this follows from Lemma \ref{le:jordan} below by (L4) and (a).

For (L6), this follows from \cite{LTXZZ}*{Theorem~4.2.6} by (b).

For (L7), this follows from Corollary \ref{co:generic} if $F^+\neq\dQ$.

For (L5-1), let $\lambda$ be a prime of $E$ satisfying (L4) and (L6), whose underlying rational prime is at least $2n(n+1)-1$. Then by (a), $\bar\rho_{\Pi_0,\lambda}$ and $\bar\rho_{\Pi_1,\lambda}$ satisfy the assumptions in Lemma \ref{le:jordan} below, with $k=O_\lambda/\lambda$ and $\Gamma=\Gamma_F$. Thus, by Lemma \ref{le:jordan}(2), assumption (b) of Lemma \ref{le:restriction} and hence (L5-1) hold.

For (L5-2), take an arithmetic Frobenius element $\phi_\fp\in\Gamma_{F^+_\fp}$. By Definition \ref{de:special_inert}, $\phi_\fp$ belongs to $\Gal(\ol{F}/F^+_{\r{rflx}})$. For $\alpha=0,1$, put $r_\alpha\coloneqq\lfloor\tfrac{n_\alpha}{2}\rfloor$ as always. By (a), the Satake parameter of $\Pi_{0,\fp}$ is $\{p^{\pm{1}},\dots,p^{\pm(2r_0-1)}\}$; and we may write the Satake parameter of $\Pi_{1,\fp}$ as $\{1,\alpha_1^{\pm{1}},\dots,\alpha_{r_1}^{\pm{1}}\}$ in which $\alpha_i$ is an algebraic number other than $1$ for $1\leq i\leq r_1$. For our purpose, we may replace $E$ by a finite extension in $\dC$ so that $\alpha_i\in E$ for $1\leq i\leq r_1$. By Proposition \ref{pr:galois}(1), we have $|\alpha_i|=1$ for $1\leq i\leq r_1$. Therefore, for all but finitely many prime $\lambda$ of $E$, we have
\begin{itemize}[label={\ding{109}}]
  \item $\{p,\alpha_1,\dots,\alpha_{r_1}\}$ is contained in $O_\lambda^\times$;

  \item $\{p^{\pm{1}}\modulo\lambda,\dots,p^{\pm(2r_0-1)}\modulo\lambda\}$ consists of distinct elements and does not contain $-1$;

  \item $\{\alpha_i\modulo\lambda\res 1\leq i\leq r_1\}$ is disjoint from $\{1,-p,-p^{-1}\}$;

  \item $\{p^{\pm{1}}\alpha_i\modulo\lambda,\dots,p^{\pm(2r_0-1)}\alpha_i\modulo\lambda\res 1\leq i\leq r_1\}$ is disjoint from $\{p,p^{-1}\}$.
\end{itemize}
Then for every prime $\lambda$ satisfying (L4) and the above properties, (L5-2) (that is, $(\r{GI}_{F',\sP}^1)$ from Lemma \ref{le:galois_element}) is satisfied by taking the element $(\bar\rho_{0+},\bar\rho_{1+},\bar\epsilon_\ell)(\phi_\fp)$.

The lemma is proved.
\end{proof}

For every integer $m\geq 1$, we denote by $J_m$ the standard upper triangular nilpotent Jordan block
\[
\begin{pmatrix}
0 & 1 & 0 & \cdots & 0 \\
& 0 & 1 & \cdots & 0 \\
&  & \ddots & \ddots & \vdots \\
&&& 0 & 1\\
&&&& 0
\end{pmatrix}
\]
or size $m$.

\begin{lem}\label{le:jordan}
Let $\Gamma$ be a group, and $k$ a field of characteristic either zero or at least $2n(n+1)-1$. Let $\rho_0\colon\Gamma\to\GL_{n_0}(k)$ and $\rho_1\colon\Gamma\to\GL_{n_1}(k)$ be two homomorphisms that are absolutely irreducible. Suppose that there exists an element $t\in\Gamma$ such that $\rho_0(t)=1+J_{n_0}$ and $\rho_1(t)=1$. Then we have
\begin{enumerate}
  \item $\rho_0\otimes\rho_1$ is absolutely irreducible;

  \item $\rho_0\otimes\rho_1$ is not a subquotient of $\ad(\rho_0\otimes\rho_1)$.
\end{enumerate}
\end{lem}

\begin{proof}
We may assume that $k$ is algebraically closed. For $\alpha=0,1$, let $V_i=k^{\oplus n_i}$ be the space which $\Gamma$ acts on through $\rho_\alpha$. By \cite{Ser94}*{Corollaire~1}, we know that both $\rho_0\otimes\rho_1$ and $\ad(\rho_0\otimes\rho_1)$ are semisimple.

For (1), we fix an element $e\in V_0$ such that the $t$-invariant subspace of $V_0$ is spanned by $e$. Then it is clear that the $t$-invariant subspace of $V_0\otimes_kV_1$ is $k.e\otimes_k V_1$. Now suppose that $W$ is a nonzero direct summand of the $k[\Gamma]$-module $V_0\otimes_kV_1$. Let $V'_1\subseteq V_1$ be the subspace such that $k.e\otimes_k V'_1$ is the $t$-invariant subspace of $W$. Then it is easy to see that $V'_1$ is closed under the action of $\Gamma$, which forces $V'_1=V_1$ since $\rho_1$ is irreducible. This further implies that $W=V_0\otimes_kV_1$ by looking at the Jordan decomposition of $t$ on $W$, hence $\rho_0\otimes\rho_1$ is irreducible.

For (2), note that $(\rho_0\otimes\rho_1)(t)$ is conjugate to $(1+J_{n_0})^{\oplus n_1}$. On the other hand, $\ad(\rho_0\otimes\rho_1)(t)$ is conjugate to
\[
\bigoplus_{i=1}^{n_0}(1+J_{2i-1})^{\oplus n_1^2}.
\]
Since $n_0$ is even and $1,3,\dots,2n_0-1$ are odd, $\rho_0\otimes\rho_1$ is not a subquotient of $\ad(\rho_0\otimes\rho_1)$ as $\ad(\rho_0\otimes\rho_1)$ is semisimple.

The lemma is proved.
\end{proof}

The following two lemmas will be used in later subsections.

\begin{lem}\label{le:purity}
The representation $\rho_{\Pi_0,\lambda}\otimes_{E_\lambda}\rho_{\Pi_1,\lambda}(n)$ is pure of weight $-1$ at every nonarchimedean place $w$ of $F$ not above $\ell$ (Definition \ref{de:purity}).
\end{lem}

\begin{proof}
This is a consequence of Proposition \ref{pr:galois}(1) and \cite{TY07}*{Lemma~1.4(3)}.
\end{proof}

\begin{lem}\label{le:intertwining}
Assume Hypothesis \ref{hy:unitary_cohomology} for $n_1$. Let $\rV_{n_1}$ be a standard indefinite hermitian space of rank $n_1$ over $F$, $\Lambda_{n_1}$ a self-dual $\prod_{v\not\in\Sigma^+_\infty\cup\Sigma^+_\mnm}O_{F_v}$-lattice in $\rV_{n_1}\otimes_F\dA_F^{\Sigma^+_\infty\cup\Sigma^+_\mnm}$, and $\lambda$ a prime of $E$. Consider a finite set $\fP$ of special inert primes of $F^+$ whose underlying rational primes are distinct and coprime to $\Sigma^+_\mnm$, and an object $\rK_{n_1}\in\fK(\rV_{n_1})$ of the form $(\rK_{n_1})_{\Sigma^+_\mnm}\times\prod_{v\not\in\Sigma^+_\infty\cup\Sigma^+_\mnm}\rU(\Lambda_{n_1})(O_{F^+_v})$. Put
\[
\fm_1\coloneqq\dT^{\Sigma^+_\mnm\cup\Sigma^+_\fP}_{n_1}\cap\Ker\(\dT^{\Sigma^+_\mnm}_{n_1}\xrightarrow{\phi_{\Pi_1}}O_E\to O_E/\lambda\)
\]
where $\Sigma^+_\fP$ is the union of $\Sigma^+_p$ for all underlying rational primes $p$ of $\fP$. Suppose that $P_{\balpha(\Pi_{1,\fp})}\modulo\lambda$ is intertwining generic (Definition \ref{de:satake_condition}) for every $\fp\in\fP$, and that
\begin{enumerate}[label=(\alph*)]
  \item either the composite homomorphism $\dT^{\Sigma^+_\mnm}_{n_1}\xrightarrow{\phi_{\Pi_1}}O_E\to O_E/\lambda$ is cohomologically generic;

  \item or $n_1=3$ and $\rH^i_\et(\Sh(\rV_3,\rK_3)_{\ol{F}},O_E/\lambda)_{\fm_1}=0$ if $i\neq 2$.
\end{enumerate}
Then for every special maximal subgroup $\rK'_{n_1,\fP}$ of $\prod_{\fp\in\fP}\rU(\rV_{n_1})(F^+_{\fp})$ and every $i\in\dZ$, we have an isomorphism
\begin{align*}
\rH^i_\et(\Sh(\rV_{n_1},\rK_{n_1})_{\ol{F}},O_\lambda)_{\fm_1}
\simeq\rH^i_\et(\Sh(\rV_{n_1},\rK^\fP_{n_1}\rK'_{n_1,\fP})_{\ol{F}},O_\lambda)_{\fm_1}
\end{align*}
of $O_\lambda[\Gamma_F]$-modules.
\end{lem}

Note that (a) is stronger than (b) when $n_1=3$.

\begin{proof}
We first note that for every $\fp\in\fP$, $\rU(\rV_{n_1})(F^+_\fp)$ has two special maximal subgroups up to conjugation, exact one of which is hyperspecial maximal.

Take an element $\fp\in\fP$, a special maximal subgroup $\rK'^\fp_{n_1,\fP}$ of $\prod_{\fp'\in\fP\setminus\{\fp\}}\rU(\rV_{n_1})(F^+_{\fp'})$, a hyperspecial maximal subgroup $\rK^\circ_{n_1,\fp}$ of $\rU(\rV_{n_1})(F^+_\fp)$, and a non-hyperspecial special maximal subgroup $\rK^\bullet_{n_1,\fp}$ of $\rU(\rV_{n_1})(F^+_\fp)$. We claim that if $\rH^i_\et(\Sh(\rV_{n_1},\rK_{n_1}^\fP\rK'^\fp_{n_1,\fP}\rK^\circ_{n_1,\fp})_{\ol{F}},O_E/\lambda)_{\fm_1}=0$ for $i\neq 2r_1$, then there is an isomorphism
\begin{align*}
\rH^i_\et(\Sh(\rV_{n_1},\rK_{n_1}^\fP\rK'^\fp_{n_1,\fP}\rK^\circ_{n_1,\fp})_{\ol{F}},O_\lambda)_{\fm_1}
\simeq
\rH^i_\et(\Sh(\rV_{n_1},\rK_{n_1}^\fP\rK'^\fp_{n_1,\fP}\rK^\bullet_{n_1,\fp})_{\ol{F}},O_\lambda)_{\fm_1}
\end{align*}
of $O_\lambda[\Gamma_F]$-modules for every $i\in\dZ$.

Fix an embedding $E_\lambda\hookrightarrow\ol\dQ_\ell$. We first show that there is an isomorphism
\begin{align}\label{eq:intertwining}
\rH^i_\et(\Sh(\rV_{n_1},\rK_{n_1}^\fP\rK'^\fp_{n_1,\fP}\rK^\circ_{n_1,\fp})_{\ol{F}},O_\lambda)_{\fm_1}\otimes_{O_\lambda}\ol\dQ_\ell
\simeq
\rH^i_\et(\Sh(\rV_{n_1},\rK_{n_1}^\fP\rK'^\fp_{n_1,\fP}\rK^\bullet_{n_1,\fp})_{\ol{F}},O_\lambda)_{\fm_1}\otimes_{O_\lambda}\ol\dQ_\ell
\end{align}
of $\ol\dQ_\ell[\Gamma_F]$-modules for every $i\in\dZ$. Let $\Lambda^\circ_{n_1,\fp}$ be the self-dual $O_{F_\fp}$-lattice in $\rV_{n_1}\otimes_FF_\fp$ whose stabilizer is $\rK^\circ_{n_1,\fp}$. Without loss of generality, we may assume that $\rK^\bullet_{n_1,\fp}$ is the stabilizer of a lattice $\Lambda^\bullet_{n_1,\fp}$ satisfying $\Lambda^\circ_{n_1,\fp}\subseteq\Lambda^\bullet_{n_1,\fp}$ and $(\Lambda^\bullet_{n_1,\fp})^\vee/\fp\Lambda^\bullet_{n_1,\fp}\simeq\dF_{p^2}$. To show \eqref{eq:intertwining}, it suffices to show that for every (necessarily cuspidal) automorphic representation $\pi_1$ of $\rU(\rV_{n_1})(\dA_{F^+})$ that appears in either side of \eqref{eq:intertwining}, the maps
\begin{align}\label{eq:intertwining_1}
\tT^{\bullet\circ}_{n_1,\fp}\colon\pi_{1,\fp}^{\rK^\circ_{n_1,\fp}}\to\pi_{1,\fp}^{\rK^\bullet_{n_1,\fp}},\quad
\tT^{\circ\bullet}_{n_1,\fp}\colon\pi_{1,\fp}^{\rK^\bullet_{n_1,\fp}}\to\pi_{1,\fp}^{\rK^\circ_{n_1,\fp}}
\end{align}
are both isomorphisms. Here, $\tT^{\bullet\circ}_{n_1,\fp}$ and $\tT^{\circ\bullet}_{n_1,\fp}$ are introduced in Definition \ref{de:ns_hecke}. Hypothesis \ref{hy:unitary_cohomology} and the Chebotarev density theorem imply that $\rho_{\BC(\pi_1),\iota_\ell}$ and $\rho_{\Pi_1,\lambda}\otimes_{E_\lambda}\ol\dQ_\ell$ have the isomorphic (irreducible) residual representations. In particular, the Satake parameter of $\BC(\pi_1)_\fp$ does not contain $\{-p,-p^{-1}\}$ by Proposition \ref{pr:galois}(2) and the assumption that $P_{\balpha(\Pi_{1,\fp})}\modulo\lambda$ is intertwining generic. Let $\tilde\pi$ be an (unramified) principal series representation of $\rU(\rV_{n_1})(F^+_\fp)$ that has $\pi_{1,\fp}$ as a constituent. By Proposition \ref{pr:enumeration_odd}(1) and the definition of the intertwining Hecke operator $\rI^\circ_{n_1,\fp}\coloneqq\tT^{\circ\bullet}_{n_1,\fp}\circ\tT^{\bullet\circ}_{n_1,\fp}$ from Definition \ref{de:ns_hecke} or Definition \ref{de:lattice}, the composite map $\tT^{\circ\bullet}_{n_1,\fp}\circ\tT^{\bullet\circ}_{n_1,\fp}\colon\tilde\pi^{\rK^\circ_{n_1,\fp}}\to\tilde\pi^{\rK^\circ_{n_1,\fp}}$ is an isomorphism. Since both $\rK^\circ_{n_1,\fp}$ and $\rK^\bullet_{n_1,\fp}$ are special maximal subgroups of $\rU(\rV_{n_1})(F^+_\fp)$, both $\tilde\pi^{\rK^\circ_{n_1,\fp}}$ and $\tilde\pi^{\rK^\bullet_{n_1,\fp}}$ are one-dimensional. It follows that the constituent of $\tilde\pi$ that has nonzero $\rK^\circ_{n_1,\fp}$-invariants is the same as the constituent that has nonzero $\rK^\bullet_{n_1,\fp}$-invariants, which further implies that the two maps in \eqref{eq:intertwining_1} are both isomorphisms. Thus, we obtain the isomorphism \eqref{eq:intertwining}.

To prove the claim it suffices to show that $\rH^i_\et(\Sh(\rV_{n_1},\rK_{n_1}^\fP\rK'^\fp_{n_1,\fP}\rK^\bullet_{n_1,\fp})_{\ol{F}},O_\lambda)_{\fm_1}$ is a free $O_\lambda$-module for every $i\in\dZ$. If we assume (a), then this follows immediately. If we assume (b), then this follows from Proposition \ref{pr:three}, Lemma \ref{le:tate_vanishing}, and a straightforward computation on the spectral sequence in Lemma \ref{le:ns_weight_odd}.

The lemma follows immediately from the above claim by induction on the number of primes $\fp\in\fP$ for which $\rK'_{n_1,\fP}$ is not hyperspecial maximal at $\fp$. Note that the initial induction hypothesis is satisfied by either (a) or (b).
\end{proof}

\subsection{Main theorem in the Selmer rank 0 case}
\label{ss:main_0}

The following lemma is a key ingredient in the proof of Theorem \ref{th:selmer0}, which is essentially the solution of the Gan--Gross--Prasad conjecture for $\Pi_0\times\Pi_1$.

\begin{lem}\label{le:ggp}
Keep the setup in \S\ref{ss:setup}. If $L(\frac{1}{2},\Pi_0\times\Pi_1)\neq 0$, then there exist
\begin{itemize}[label={\ding{109}}]
  \item a standard definite hermitian space $\rV^\circ_n$ of rank $n$ over $F$, together with a self-dual $\prod_{v\not\in\Sigma^+_\infty\cup\Sigma^+_\mnm}O_{F_v}$-lattice $\Lambda^\circ_n$ in $\rV^\circ_n\otimes_F\dA_F^{\Sigma^+_\infty\cup\Sigma^+_\mnm}$ (and put $\rV^\circ_{n+1}\coloneqq(\rV^\circ_n)_\sharp$ and $\Lambda^\circ_{n+1}\coloneqq(\Lambda^\circ_n)_\sharp$),

  \item an object $(\rK^\circ_n,\rK^\circ_{n+1})\in\fK(\rV^\circ_n)_\sp$ in which $\rK^\circ_{n_\alpha}$ is of the form
      \[
      \rK^\circ_{n_\alpha}=\prod_{v\in\Sigma^+_\mnm}(\rK^\circ_{n_\alpha})_v\times
      \prod_{v\not\in\Sigma^+_\infty\cup\Sigma^+_\mnm}\rU(\Lambda^\circ_{n_\alpha})(O_{F^+_v})
      \]
      for $\alpha=0,1$,
\end{itemize}
such that
\[
\sum_{s\in\Sh(\rV^\circ_n,\rK^\circ_n)}f(s,\sh_\uparrow(s))\neq 0
\]
for some element $f\in O_E[\Sh(\rV^\circ_{n_0},\rK^\circ_{n_0})][\Ker\phi_{\Pi_0}]\otimes_{O_E}O_E[\Sh(\rV^\circ_{n_1},\rK^\circ_{n_1})][\Ker\phi_{\Pi_1}]$.
\end{lem}

\begin{proof}
This follows from the direction (1)$\Rightarrow$(2) of \cite{BPLZZ}*{Theorem~1.8}, together with \cite{BPLZZ}*{Remark~4.17}. Note that since our $\Pi_0$ and $\Pi_1$ are relevant representations of $\GL_{n_0}(\dA_F)$ and $\GL_{n_1}(\dA_F)$, respectively, both members in the pair of hermitian spaces in (2) of \cite{BPLZZ}*{Theorem~1.8} have to be standard definite.
\end{proof}

\begin{theorem}\label{th:selmer0}
Keep the setup in \S\ref{ss:setup}. Assume Hypothesis \ref{hy:unitary_cohomology} for both $n$ and $n+1$. If $L(\frac{1}{2},\Pi_0\times\Pi_1)\neq 0$, then for all admissible primes $\lambda$ of $E$, and for all but finitely many weakly admissible primes $\lambda$ of $E$ when $n=2$, we have
\[
\rH^1_f(F,\rho_{\Pi_0,\lambda}\otimes_{E_\lambda}\rho_{\Pi_1,\lambda}(n))=0.
\]
\end{theorem}

\begin{proof}
By Lemma \ref{le:ggp}, we may fix the choices of $\rV^\circ_n$, $\Lambda^\circ_n$, $(\rK^\circ_n,\rK^\circ_{n+1})$ in that lemma such that
\[
\sum_{s\in\Sh(\rV^\circ_n,\rK^\circ_n)}f(s,\sh_\uparrow(s))\neq 0
\]
for some $f\in O_E[\Sh(\rV^\circ_{n_0},\rK^\circ_{n_0})][\Ker\phi_{\Pi_0}]\otimes_{O_E}O_E[\Sh(\rV^\circ_{n_1},\rK^\circ_{n_1})][\Ker\phi_{\Pi_1}]$. Moreover, by Lemma \ref{le:transferable}(3), we may assume that $(\rK^\circ_{n_0})_v$ is transferable (Definition \ref{de:transferable}) for $v\in\Sigma^+_\mnm$.

We take a prime $\lambda$ of $E$ with the underlying rational prime $\ell$. We adopt notation in \S\ref{ss:rankin_selberg} with the initial data in Definition \ref{de:lambda}. Define two nonnegative integers $m_{\r{per}}$ and $m_{\r{lat}}$ as follows.
\begin{enumerate}
  \item Let $m_{\r{per}}$ be the largest (nonnegative) integer such that
     \[
     \sum_{s\in\Sh(\rV^\circ_n,\rK^\circ_n)}f(s,\sh_\uparrow(s))\in \lambda^{m_{\r{per}}}O_E
     \]
     for every $f\in O_E[\Sh(\rV^\circ_{n_0},\rK^\circ_{n_0})][\Ker\phi_{\Pi_0}]\otimes_{O_E}O_E[\Sh(\rV^\circ_{n_1},\rK^\circ_{n_1})][\Ker\phi_{\Pi_1}]$.

  \item We choose a standard \emph{indefinite} hermitian space $\rV_{n_1}$ over $F$ of rank $n_1$, together with an identification $\rU((\rV^\circ_{n_1})^\infty)\simeq\rU(\rV_{n_1}^\infty)$ of reductive groups over $\dA^\infty_{F^+}$.\footnote{There are many choices of such $\rV_{n_1}$ and the isomorphism. We choose one only to get some control on the discrepancy of the integral cohomology of Shimura varieties and the lattice coming from Galois representations.} In particular, we have the Shimura variety $\Sh(\rV_{n_1},\rK^\circ_{n_1})$. By Hypothesis \ref{hy:unitary_cohomology}, we have an isomorphism
      \[
      \rH^{2r_1}_\et(\Sh(\rV_{n_1},\rK^\circ_{n_1})_{\ol{F}},E_\lambda(r_1))/\Ker\phi_{\Pi_1}
      \simeq(\rR_1^\tc\otimes_{O_\lambda}E_\lambda)^{\oplus\mu_1}
      \]
      of $E_\lambda[\Gamma_F]$-modules for some integer $\mu_1>0$. We fix a map
      \[
      \rH^{2r_1}_\et(\Sh(\rV_{n_1},\rK^\circ_{n_1})_{\ol{F}},O_\lambda(r_1))/\Ker\phi_{\Pi_1}
      \to(\rR_1^\tc)^{\oplus\mu_1}
      \]
      of $O_\lambda[\Gamma_F]$-modules whose kernel and cokernel are both $O_\lambda$-torsion. Then we let $m_{\r{lat}}$ be the smallest nonnegative integer such that both the kernel and the cokernel are annihilated by $\lambda^{m_{\r{lat}}}$.
\end{enumerate}
Now we assume that either $\lambda$ is admissible, or $n=2$ and $\lambda$ is weakly admissible and satisfies $\rH^i_\et(\Sh(\rV_3,\rK^\circ_3)_{\ol{F}},O_E/\lambda)/\Ker\phi_{\Pi_1}=0$ for $i\neq 2$ (which only excludes finitely many primes). Note that in both cases we have $\rH^i_\et(\Sh(\rV_{n_1},\rK^\circ_{n_1})_{\ol{F}},O_E/\lambda)/\Ker\phi_{\Pi_1}=0$ if $i\neq 2r_1$ by Proposition \ref{th:generic}.

We start to prove the theorem by contradiction, hence assume
\[
\dim_{E_\lambda}\rH^1_f(F,\rho_{\Pi_0,\lambda}\otimes_{E_\lambda}\rho_{\Pi_1,\lambda}(n))\geq 1.
\]

Take a sufficiently large positive integer $m$ which will be determined later. By Lemma \ref{le:purity}, we may apply Proposition \ref{pr:selmer_reduction} by taking $\Sigma$ to be the set of places of $F$ above $\Sigma^+_\mnm\cup\Sigma^+_\ell$. Then we obtain a submodule $S$ of $\rH^1_{f,\rR}(F,\bar\rR^{(m)})$ that is free of rank $1$ over $O_\lambda/\lambda^{m-m_\Sigma}$ such that $\loc_w\res_S=0$ for every nonarchimedean place $w\in\Sigma$ not above $\ell$. Now we apply the discussion in \S\ref{ss:galois_lemma} to the submodule $S\subseteq\rH^1(F,\bar\rR^{(m)})$. By (L5-1) and Lemma \ref{le:image}, we obtain an injective map
\[
\theta_S\colon\Gal(F_S/F_{\bar\rho^{(m)}})\to\Hom_{O_\lambda}(S,\bar\rR^{(m)})
\]
whose image generates an $O_\lambda$-submodule containing $\lambda^{\fr_{\bar\rR^{(m)}}}\Hom_{O_\lambda}(S,\bar\rR^{(m)})$, which further contains $\lambda^{\fr_\rR}\Hom_{O_\lambda}(S,\bar\rR^{(m)})$ by Lemma \ref{le:reducibility} and (L3). By (L5-2) and Lemma \ref{le:galois_element}, we may choose an element $(\gamma_1,\gamma_2,\xi)$ in the image of $(\bar\rho^{(m)}_{1+},\bar\rho^{(m)}_{2+},\bar\epsilon_\ell^{(m)})\res_{\Gal(\ol{F}/F^+_{\r{rflx}})}$ satisfying (a--e) in Lemma \ref{le:galois_element}. It then gives rise to an element $\gamma\in(\GL_{n_0n_1}(O_\lambda/\lambda^m)\times(O_\lambda/\lambda^m)^\times)\fc$ as in Notation \ref{no:h_gamma} such that $(\bar\rR^{(m)})^{h_\gamma}$ is a free $O_\lambda/\lambda^m$-module of rank $1$. Now we apply the discussion in \S\ref{ss:localization}. By Proposition \ref{pr:selmer_localization} (with $m_0=m_\Sigma$ and $r_S=1$), we may fix an $(S,\gamma)$-abundant element $\Psi\in G_{S,\gamma}$ (Definition \ref{de:abundant}).

We apply the discussion and notation in \S\ref{ss:first_reciprocity} to our situation with $\lambda$, $m$, $\Sigma^+_{\lr,\rI}=\emptyset$, $\Sigma^+_\rI=\Sigma^+_\mnm$, $(\rV^\circ_n,\Lambda^\circ_n)$, $\rK^\circ_n$ and $(\rK^\circ_n,\rK^\circ_{n+1})$. By the Chebotarev density theorem, we can choose a $\gamma$-associated place (Definition \ref{de:associated}) $w^{(m)}_+$ of $F^{(m)}_+$ satisfying $\Psi_{w^{(m)}}=\Psi$ and whose underlying prime $\fp$ of $F^+$ (and the underlying rational prime $p$) is a special inert prime satisfying (PI1)--(PI7) and
\begin{description}
  \item[(PI8)] the natural map
     \[
     \rH^i_\et(\Sh(\rV_{n_1},\rK^\circ_{n_1})_{\ol{F}},O_\lambda(r_1))/(\dT^{\Sigma^+_\rI\cup\Sigma^+_p}_{n_1}\cap\Ker\phi_{\Pi_1})
     \to\rH^i_\et(\Sh(\rV_{n_1},\rK^\circ_{n_1})_{\ol{F}},O_\lambda(r_1))/\Ker\phi_{\Pi_1}
     \]
     is an isomorphism for every integer $i$.
\end{description}
We also choose remaining data in \S\ref{ss:ns_initial} with $\dQ_p^\Phi=\dQ_{p^2}$, data as in Notation \ref{no:ns_uniformization}, and an indefinite uniformization datum as in Notation \ref{no:ns_uniformization_indefinite}. By the definition of $m_{\r{per}}$, we have
\begin{align}\label{eq:selmer0_5}
\exp_\lambda\(\CF_{\Sh(\rV^\circ_n,\rK^\circ_\sp)},
O_E[\Sh(\rV^\circ_{n_0},\rK^\circ_{n_0})\times\Sh(\rV^\circ_{n_1},\rK^\circ_{n_1})]/(\fn_0,\fn_1)\)\geq m-m_{\r{per}},
\end{align}
where we recall that
\[
\fn_\alpha=\dT^{\Sigma^+_\rI\cup\Sigma^+_p}_{n_\alpha}\cap\Ker\(\dT^{\Sigma^+_\mnm}_{n_\alpha}\xrightarrow{\phi_{\Pi_\alpha}}O_E\to O_E/\lambda^m\)
\]
for $\alpha=0,1$. Here, $\CF_{\Sh(\rV^\circ_n,\rK^\circ_\sp)}$ is nothing but the characteristic function of the graph $\graph\Sh(\rV^\circ_n,\rK^\circ_n)$ of the map $\Sh(\rV^\circ_n,\rK^\circ_n)\to\Sh(\rV^\circ_{n+1},\rK^\circ_{n+1})$.

We claim that there exists an element $c_1\in\rH^1(F,\bar\rR^{(m_1)\tc})$ for some positive integer $m_1\leq m$ satisfying \begin{align}\label{eq:selmer0_4}
\exp_\lambda\(\partial_\fp\loc_\fp(c_1),\rH^1_\sing(F_\fp,\bar\rR^{(m_1)\tc})\)\geq m-m_{\r{per}}-m_{\r{lat}};
\end{align}
and such that for every nonarchimedean place $w$ of $F$ not above $\Sigma^+\cup\{\fp\}$,
\begin{align}\label{eq:selmer0_6}
\loc_w(c_1)\in\rH^1_\ns(F_w,\bar\rR^{(m_1)\tc})
\end{align}
holds.

We first prove the theorem assuming the existence of such $c_1$. Fix a generator of the submodule $S\subseteq\rH^1_{f,\rR}(F,\bar\rR^{(m)})$ and denote by its image in $\rH^1(F,\bar\rR^{(m_1)})$ by $s_1$. We also identify $\bar\rR^{(m_1)\tc}$ with $(\bar\rR^{(m_1)})^*$ via the polarization $\Xi$. Now we compute the local Tate pairing $\langle s_1,c_1\rangle_w$ \eqref{eq:local_tate} for every nonarchimedean place $w$ of $F$.
\begin{itemize}[label={\ding{109}}]
  \item Suppose that $w$ is above $\Sigma^+_\mnm$. Then we have $\loc_w(s_1)=0$ by our choice of $S$. Thus, $\langle s_1,c_1\rangle_w=0$.

  \item Suppose that $w$ is above $\Sigma^+_\ell$. Then by (L2), $\rR_\dQ$ is crystalline with Hodge--Tate weights in $[-n,n-1]$. Thus, we have $\loc_w(s_1)\in\rH^1_\ns(F_w,\bar\rR^{(m_1)})$ by Lemma \ref{le:global_local}(2) and (L1). By \eqref{eq:selmer0_6}, Lemma \ref{le:crystalline_tate} and (L1), we have $\lambda^{m_{\r{dif}}}\langle s_1,c_1\rangle_w=0$ where $\fd_\lambda=\lambda^{m_{\r{dif}}}\subseteq O_\lambda$ is the different ideal of $E_\lambda$ over $\dQ_\ell$.

  \item Suppose that $w$ is not above $\Sigma^+_\mnm\cup\Sigma^+_\ell\cup\{\fp\}$. Then by (L2), $\rR$ is unramified. Thus, we have $\loc_w(s_1)\in\rH^1_\ns(F_w,\bar\rR^{(m_1)})$ by Lemma \ref{le:global_local}(1). By \eqref{eq:selmer0_6} and Lemma \ref{le:tate}, we have $\langle s_1,c_1\rangle_w=0$.

  \item Suppose that $w$ is the unique place above $\fp$. By Proposition \ref{co:selmer_localization}, we have
      \[
      \exp_\lambda\(\loc_w(s_1),\rH^1_\ns(F_w,\bar\rR^{(m_1)})\)\geq m_1-m_\Sigma-\fr_\rR.
      \]
      By \eqref{eq:selmer0_4} and Lemma \ref{le:tate} again, we have
      \[
      \exp_\lambda\(\langle s_1,c_1\rangle_w,O_\lambda/\lambda^{m_1}\)\geq m-m_{\r{per}}-m_{\r{lat}}-m_\Sigma-\fr_\rR.
      \]
\end{itemize}
Therefore, as long as we take $m$ such that $m>m_{\r{per}}+m_{\r{lat}}+m_\Sigma+\fr_\rR+m_{\r{dif}}$, we will have a contradiction to the relation
\[
\sum_{w}\langle s_1,c_1\rangle_w=0,
\]
where the sum is taken over all nonarchimedean places $w$ of $F$. The theorem is proved.

Now we consider the claim on the existence of $c_1$. By (L4), (L6), and Theorem \ref{th:raising}(5), we have an isomorphism
\begin{align}\label{eq:selmer0_0}
\rH^{2r_0-1}_\et((\Sh(\rV'_{n_0},\tj_{n_0}\rK^{p\circ}_{n_0}\rK'_{n_0,p})_{\ol{F}},O_\lambda(r_0))/\fn_0
\xrightarrow{\sim}\bigoplus_{i=1}^{\mu_0}\bar\rR_0^{(m_i)\tc}
\end{align}
of $O_\lambda[\Gamma_F]$-modules, for finitely many positive integers $m_1,\dots,m_{\mu_0}$ at most $m$. Assumption \ref{as:first_vanishing} for $\alpha=0$ is satisfied by Lemma \ref{le:ns_pbc} and (L7) when $n\geq 3$, and by Lemma \ref{le:ns_pbc} and (L4) when $n=2$.

By Lemma \ref{le:intertwining}, we have an isomorphism
\[
\rH^i_\et(\Sh(\rV_{n_1},\rK^\circ_{n_1})_{\ol{F}},O_\lambda)_{\fm_1}\simeq
\rH^i_\et(\Sh(\rV'_{n_1},\tj_{n_1}\rK^{p\circ}_{n_1}\rK'_{n_1,p})_{\ol{F}},O_\lambda)_{\fm_1}
\]
of $O_\lambda[\Gamma_F]$-modules. Then Assumption \ref{as:first_vanishing} for $\alpha=1$ is satisfied by Lemma \ref{le:ns_pbc}, (PI8), and the fact that $\rH^i_\et(\Sh(\rV_{n_1},\rK^\circ_{n_1})_{\ol{F}},O_E/\lambda)/\Ker\phi_{\Pi_1}=0$ if $i\neq 2r_1$. Moreover, by (PI8), we may fix a map
\[
\rH^{2r_1}_\et(\Sh(\rV'_{n_1},\tj_{n_1}\rK^{p\circ}_{n_1}\rK'_{n_1,p})_{\ol{F}},O_\lambda(r_1))/
(\dT^{\Sigma^+_\rI\cup\Sigma^+_p}_{n_1}\cap\Ker\phi_{\Pi_1})
\to(\rR_1^\tc)^{\oplus\mu_1}
\]
of $O_\lambda[\Gamma_F]$-modules whose kernel and cokernel are both annihilated by $\lambda^{m_{\r{lat}}}$. Taking quotient by $\lambda^m$, we obtain a map
\begin{align}\label{eq:selmer0_1}
\rH^{2r_1}_\et(\Sh(\rV'_{n_1},\tj_{n_1}\rK^{p\circ}_{n_1}\rK'_{n_1,p})_{\ol{F}},O_\lambda(r_1))/\fn_1
\to(\bar\rR_1^{(m)\tc})^{\oplus\mu_1}
\end{align}
of $O_\lambda[\Gamma_F]$-modules whose kernel and cokernel are both annihilated by $\lambda^{m_{\r{lat}}}$.

To continue, we adopt the notational abbreviation prior to Corollary \ref{co:first}. By Lemma \ref{le:single_compact} and the K\"{u}nneth formula, we obtain a map
\begin{align}\label{eq:selmer0_2}
\Upsilon\colon\rH^{2n-1}_\et((\Sh'_{n_0}\times_{\Spec{F}}\Sh'_{n_1})_{\ol{F}},O_\lambda(n))/(\fn_0,\fn_1)\to\bigoplus_{i=1}^\mu\bar\rR^{(m_i)\tc}
\end{align}
of $O_\lambda[\Gamma_F]$-modules whose kernel and cokernel are both annihilated by $\lambda^{m_{\r{lat}}}$, from \eqref{eq:selmer0_0} and \eqref{eq:selmer0_1}. Here, we have re-indexed $\mu_1$ copies of $\{m_1,\dots,m_{\mu_0}\}$ into $\mu\coloneqq\mu_0\mu_1$ positive integers at most $m$. Recall that we have a class
\[
\AJ(\Sh'_\sp)\in\rH^1(F,\rH^{2n-1}_\et((\Sh'_{n_0}\times_{\Spec{F}}\Sh'_{n_1})_{\ol{F}},O_\lambda(n))/(\fn_0,\fn_1)),
\]
where $\Sh'_\sp$ is nothing but the graph of the morphism $\Sh'_n\to\Sh'_{n+1}$. By Corollary \ref{co:first} and \eqref{eq:selmer0_5}, we have
\begin{align}\label{eq:selmer0_3}
\exp_\lambda\(\partial_\fp\loc_\fp\AJ(\Sh'_\sp),
\rH^1_\sing(F_\fp,\rH^{2n-1}_\et((\Sh'_{n_0}\times_{\Spec{F}}\Sh'_{n_1})_{\ol{F}},O_\lambda(n))/(\fn_0,\fn_1))\)\geq m-m_{\r{per}}.
\end{align}
For every $1\leq i\leq\mu$, let
\[
\Upsilon_i\colon\rH^{2n-1}_\et((\Sh'_{n_0}\times_{\Spec{F}}\Sh'_{n_1})_{\ol{F}},O_\lambda(n))/(\fn_0,\fn_1)\to\bar\rR^{(m_i)\tc}
\]
be the composition of $\Upsilon$ \eqref{eq:selmer0_2} with the projection to the $i$-th factor; and put
\[
c_i\coloneqq\rH^1(F,\Upsilon_i)(\AJ(\Sh'_\sp))\in\rH^1(F,\bar\rR^{(m_i)\tc}).
\]
Then \eqref{eq:selmer0_3} implies
\[
\max_{1\leq i\leq\mu}\exp_\lambda\(\partial_\fp\loc_\fp(c_i),\rH^1_\sing(F_\fp,\bar\rR^{(m_i)\tc})\)\geq m-m_{\r{per}}-m_{\r{lat}}.
\]
Without loss of generality, we obtain \eqref{eq:selmer0_4}. On the other hand, as both $\Sh'_n$ and $\Sh'_{n+1}$ have smooth models over $O_{F_w}$ for which (an analogue of) Lemma \ref{le:qs_pbc} holds, we obtain \eqref{eq:selmer0_6}.
\end{proof}

Now we deduce two concrete consequences from Theorem \ref{th:selmer0}.

\begin{corollary}\label{co:sym}
Let $n\geq 2$ be an integer and denote by $n_0$ and $n_1$ the unique even and odd numbers in $\{n,n+1\}$, respectively. Let $A_0$ and $A_1$ be two modular elliptic curves over $F^+$ such that $\End({A_0}_{\ol{F}})=\End({A_1}_{\ol{F}})=\dZ$. Suppose that
\begin{enumerate}[label=(\alph*)]
  \item ${A_0}_{\ol{F}}$ and ${A_1}_{\ol{F}}$ are not isogenous to each other;

  \item both $\Sym^{n_0-1}A_0$ and $\Sym^{n_1-1}A_1$ are modular; and

  \item $F^+\neq\dQ$ if $n\geq 3$.
\end{enumerate}
If the (central critical) $L$-value $L(n,\Sym^{n_0-1}{A_0}_F\times\Sym^{n_1-1}{A_1}_F)$ does not vanish, then we have
\[
\rH^1_f(F,\Sym^{n_0-1}\rH^1_\et({A_0}_{\ol{F}},\dQ_\ell)\otimes_{\dQ_\ell}\Sym^{n_1-1}\rH^1_\et({A_1}_{\ol{F}},\dQ_\ell)(n))=0
\]
for all but finitely many rational primes $\ell$.
\end{corollary}

\begin{proof}
By (b) and \cite{AC89}, both $\Sym^{n_0-1}{A_0}_F$ and $\Sym^{n_1-1}{A_1}_F$ are modular. Thus, we may let $\Pi_\alpha$ be the (cuspidal) automorphic representation of $\GL_{n_\alpha}(\dA_F)$ associated to $\Sym^{n_\alpha-1}{A_\alpha}_F$ for $\alpha=0,1$, which is a relevant representation (Definition \ref{de:relevant}). We also have the identity
\[
L(n+s,\Sym^{n_0-1}{A_0}_F\times\Sym^{n_1-1}{A_1}_F)=L(\tfrac{1}{2}+s,\Pi_0\times\Pi_1)
\]
of $L$-functions, and that the representation of $\Gamma_F$ on $\Sym^{n_\alpha-1}\rH^1_\et({A_\alpha}_{\ol{F}},\dQ_\ell)$ is isomorphic to $\rho_{\Pi_\alpha,\ell}$ for $\alpha=0,1$. By Proposition \ref{pr:unitary_cohomology} and (c), Hypothesis \ref{hy:unitary_cohomology} is known in this case. Then the corollary follows immediately from Theorem \ref{th:selmer0} and Lemma \ref{pr:elliptic_admissible} (where we use (a) and (c)) with $E=\dQ$.
\end{proof}

\begin{remark}\label{re:sym}
In this remark, we summarize the current knowledge on the modularity of symmetric powers of elliptic curves, namely, condition (a) in Corollary \ref{co:sym}. Let $A$ be a modular elliptic curve over $F^+$ such that $\End(A_{\ol{F}})=\dZ$. We have
\begin{itemize}[label={\ding{109}}]
  \item $\Sym^2A$ is modular by \cite{GJ76};

  \item $\Sym^3A$ is modular by \cite{KS02};

  \item $\Sym^4A$ is modular by \cite{Kim03};

  \item $\Sym^5A$ and $\Sym^6A$ are modular if $F^+$ is linearly disjoint from $\dQ(\zeta_5)$ over $\dQ$;

  \item $\Sym^7A$ is modular if $F^+$ is linearly disjoint from $\dQ(\zeta_{35})$ over $\dQ$;

  \item $\Sym^8A$ is modular if $F^+$ is linearly disjoint from $\dQ(\zeta_7)$ over $\dQ$;
\end{itemize}
in which the last three cases are obtained in a series of recent work \cites{CT1,CT2,CT3} of Clozel and Thorne.

After we completed this article, we have learnt the groundbreaking result of Newton--Thorne \cites{NT1,NT2} where they prove the modularity of all symmetric powers of elliptic curves over $\dQ$ without complex multiplication. In particular, it follows that $\Sym^n A$ is modular if $F^+/\dQ$ is solvable and $A$ is the base change of an elliptic curve over $\dQ$.
\end{remark}

\begin{corollary}\label{co:abstract}
Keep the setup in \S\ref{ss:setup}. Suppose that
\begin{enumerate}[label=(\alph*)]
  \item there exists a very special inert prime $\fp$ of $F^+$ (Definition \ref{de:special_inert}) such that $\Pi_{0,\fp}$ is Steinberg, and $\Pi_{1,\fp}$ is unramified whose Satake parameter contains $1$ exactly once;

  \item for $\alpha=0,1$, there exists a nonarchimedean place $w_\alpha$ of $F$ such that $\Pi_{\alpha,w_\alpha}$ is supercuspidal; and

  \item $F^+\neq\dQ$ if $n\geq 3$.
\end{enumerate}
If $L(\frac{1}{2},\Pi_0\times\Pi_1)\neq 0$, then for all but finitely many primes $\lambda$ of $E$, we have
\[
\rH^1_f(F,\rho_{\Pi_0,\lambda}\otimes_{E_\lambda}\rho_{\Pi_1,\lambda}(n))=0.
\]
\end{corollary}

\begin{proof}
This follows from Theorem \ref{th:selmer0} and Lemma \ref{le:abstract}.
\end{proof}

\subsection{Main theorem in the Selmer rank 1 case}
\label{ss:main_1}

We state the following weak version of the arithmetic Gan--Gross--Prasad conjecture.

\begin{conjecture}\label{co:aggp}
Suppose that $L(\frac{1}{2},\Pi_0\times\Pi_1)=0$ but $L'(\frac{1}{2},\Pi_0\times\Pi_1)\neq 0$. Then there exist
\begin{itemize}[label={\ding{109}}]
  \item a standard indefinite hermitian space $\rV_n$ of rank $n$ over $F$, together with a self-dual $\prod_{v\not\in\Sigma^+_\infty\cup\Sigma^+_\mnm}O_{F_v}$-lattice $\Lambda_n$ in $\rV_n\otimes_F\dA_F^{\Sigma^+_\infty\cup\Sigma^+_\mnm}$ (and put $\rV_{n+1}\coloneqq(\rV_n)_\sharp$ and $\Lambda_{n+1}\coloneqq(\Lambda_n)_\sharp$),

  \item an object $(\rK_n,\rK_{n+1})\in\fK(\rV_n)_\sp$ in which $\rK_{n_\alpha}$ is of the form
      \[
      \rK_{n_\alpha}=\prod_{v\in\Sigma^+_\mnm}(\rK_{n_\alpha})_v\times
      \prod_{v\not\in\Sigma^+_\infty\cup\Sigma^+_\mnm}\rU(\Lambda_{n_\alpha})(O_{F^+_v})
      \]
      for $\alpha=0,1$,
\end{itemize}
such that for every prime $\lambda$ of $E$, the graph $\graph\Sh(\rV_n,\rK_n)$ of the morphism $\sh_\uparrow\colon\Sh(\rV_n,\rK_n)\to\Sh(\rV_{n+1},\rK_{n+1})$ \eqref{eq:qs_functoriality_shimura} is nonvanishing in the quotient Chow group
\[
\CH^n(\Sh(\rV_{n_0},\rK_{n_0})\times_{\Spec{F}}\Sh(\rV_{n_1},\rK_{n_1}))_E/(\Ker\phi_{\Pi_0},\Ker\phi_{\Pi_1}).
\]
\end{conjecture}

In the situation of the above conjecture, since both $\Pi_0$ and $\Pi_1$ are cuspidal, we have
\[
\rH^i_\et((\Sh(\rV_{n_0},\rK_{n_0})\times_{\Spec{F}}\Sh(\rV_{n_1},\rK_{n_1}))_{\ol{F}},E_\lambda)/(\Ker\phi_{\Pi_0},\Ker\phi_{\Pi_1})=0
\]
if $i\neq 2n-1$. In particular, the Hochschild--Serre spectral sequence gives rise to a coboundary map
\begin{align*}
\AJ^{\Pi_0,\Pi_1}_\lambda&\colon\rZ^n(\Sh(\rV_{n_0},\rK_{n_0})\times_{\Spec{F}}\Sh(\rV_{n_1},\rK_{n_1})) \\
&\to\rH^1(F,\rH^{2n-1}_\et((\Sh(\rV_{n_0},\rK_{n_0})\times_{\Spec{F}}\Sh(\rV_{n_1},\rK_{n_1}))_{\ol{F}},E_\lambda(n))
/(\Ker\phi_{\Pi_0},\Ker\phi_{\Pi_1})).
\end{align*}

\begin{theorem}\label{th:selmer1}
Keep the setup in \S\ref{ss:setup}. Assume Hypothesis \ref{hy:unitary_cohomology} for both $n$ and $n+1$. Let $\lambda$ be a prime of $E$ for which there exist
\begin{itemize}[label={\ding{109}}]
  \item a standard indefinite hermitian space $\rV_n$ of rank $n$ over $F$, together with a self-dual $\prod_{v\not\in\Sigma^+_\infty\cup\Sigma^+_\mnm}O_{F_v}$-lattice $\Lambda_n$ in $\rV_n\otimes_F\dA_F^{\Sigma^+_\infty\cup\Sigma^+_\mnm}$ (and put $\rV_{n+1}\coloneqq(\rV_n)_\sharp$ and $\Lambda_{n+1}\coloneqq(\Lambda_n)_\sharp$),

  \item an object $(\rK_n,\rK_{n+1})\in\fK(\rV_n)_\sp$ in which $\rK_{n_\alpha}$ is of the form
      \[
      \rK_{n_\alpha}=\prod_{v\in\Sigma^+_\mnm}(\rK_{n_\alpha})_v\times
      \prod_{v\not\in\Sigma^+_\infty\cup\Sigma^+_\mnm}\rU(\Lambda_{n_\alpha})(O_{F^+_v})
      \]
      for $\alpha=0,1$, satisfying that $(\rK_{n_0})_v$ is a transferable open compact subgroup (Definition \ref{de:transferable}) of $\rU(\rV^\circ_{n_0})(F^+_v)$ for $v\in\Sigma^+_\mnm$,
\end{itemize}
such that
\begin{align}\label{eq:selmer1}
\AJ^{\Pi_0,\Pi_1}_\lambda(\graph\Sh(\rV_n,\rK_n))\neq 0.
\end{align}
If we further assume that either $\lambda$ is admissible, or $n=2$ and $\lambda$ is weakly admissible and satisfies $\rH^i_\et(\Sh(\rV_3,\rK_3)_{\ol{F}},O_E/\lambda)/\Ker\phi_{\Pi_1}=0$ for $i\neq 2$, then we have
\[
\dim_{E_\lambda}\rH^1_f(F,\rho_{\Pi_0,\lambda}\otimes_{E_\lambda}\rho_{\Pi_1,\lambda}(n))=1.
\]
\end{theorem}

\begin{remark}
In fact, \eqref{eq:selmer1} already implies that the global epsilon factor of $\Pi_0\times\Pi_1$ is $-1$.
\end{remark}

\begin{proof}[Proof of Theorem \ref{th:selmer1}]
We take a prime $\lambda$ of $E$ for which we may choose data $\rV_n$, $\Lambda_n$, $(\rK_n,\rK_{n+1})$ as in the statement of the theorem such that $\AJ^{\Pi_0,\Pi_1}_\lambda(\graph\Sh(\rV_n,\rK_n))\neq 0$. We assume that $\lambda$ satisfies either (a) or (b) of Lemma \ref{le:intertwining}. Lemma \ref{le:purity} and (L2) imply that $\AJ^{\Pi_0,\Pi_1}_\lambda(\graph\Sh(\rV_n,\rK_n))$ belongs to the subspace
\[
\rH^1_f(F,\rH^{2n-1}_\et((\Sh(\rV_{n_0},\rK_{n_0})
\times_{\Spec{F}}\Sh(\rV_{n_1},\rK_{n_1}))_{\ol{F}},E_\lambda(n))/(\Ker\phi_{\Pi_0},\Ker\phi_{\Pi_1}))
\]
and hence to the submodule
\[
\rH^1_f(F,\rH^{2n-1}_\et((\Sh(\rV_{n_0},\rK_{n_0})
\times_{\Spec{F}}\Sh(\rV_{n_1},\rK_{n_1}))_{\ol{F}},O_\lambda(n))/(\Ker\phi_{\Pi_0},\Ker\phi_{\Pi_1}))
\]
by Definition \ref{de:bk_integral}. We also note that $\rH^i_\et(\Sh(\rV_{n_1},\rK_{n_1})_{\ol{F}},O_E/\lambda)/\Ker\phi_{\Pi_1}=0$ if $i\neq 2r_1$ by Proposition \ref{th:generic} and the same condition for $r_1=1$.

We adopt notation in \S\ref{ss:rankin_selberg} with the initial data in Definition \ref{de:lambda}. Define two nonnegative integers $m_{\r{per}}$ and $m_{\r{lat}}$ as follows.
\begin{enumerate}
  \item By Hypothesis \ref{hy:unitary_cohomology}, we may choose a map
     \[
     \rH^{2n-1}_\et((\Sh(\rV_{n_0},\rK_{n_0})
     \times_{\Spec{F}}\Sh(\rV_{n_1},\rK_{n_1}))_{\ol{F}},O_\lambda(n))/(\Ker\phi_{\Pi_0},\Ker\phi_{\Pi_1})\to\rR^\tc
     \]
     of $O_\lambda[\Gamma_F]$-modules such that the induced image of $\AJ^{\Pi_0,\Pi_1}_\lambda(\graph\Sh(\rV_n,\rK_n))$ in $\rH^1_f(F,\rR^\tc)$, denoted by $s^\tc$, is non-torsion. Let $s\in\rH^1_f(F,\rR)$ be the element corresponding to $s^\tc$ under the isomorphism in Lemma \ref{le:conjugate}. We put
     \[
     m_{\r{per}}\coloneqq\ord_\lambda\(s,\rH^1_f(F,\rR)/\rH^1_f(F,\rR)_\tor\)
     \]
     (Definition \ref{de:divisibility}), which is a nonnegative integer.

  \item By Hypothesis \ref{hy:unitary_cohomology}, we have an isomorphism
      \[
      \rH^{2r_1}_\et(\Sh(\rV_{n_1},\rK_{n_1})_{\ol{F}},E_\lambda(r_1))/\Ker\phi_{\Pi_1}
      \simeq(\rR_1^\tc\otimes_{O_\lambda}E_\lambda)^{\oplus\mu_1}
      \]
      of $E_\lambda[\Gamma_F]$-modules for some integer $\mu_1>0$. We fix a map
      \[
      \rH^{2r_1}_\et(\Sh(\rV_{n_1},\rK_{n_1})_{\ol{F}},O_\lambda(r_1))/\Ker\phi_{\Pi_1}\to(\rR_1^\tc)^{\oplus\mu_1}
      \]
      of $O_\lambda[\Gamma_F]$-modules whose kernel and cokernel are both $O_\lambda$-torsion. Then we let $m_{\r{lat}}$ be the smallest nonnegative integer such that both the kernel and the cokernel are annihilated by $\lambda^{m_{\r{lat}}}$.
\end{enumerate}
Note that in (1), we obtain an element $s\in\rH^1_f(F,\rR)_\dQ=\rH^1_f(F,\rR_\dQ)=\rH^1_f(F,\rho_{\Pi_0,\lambda}\otimes_{E_\lambda}\rho_{\Pi_1,\lambda}(n))$ that is nonzero. In particular, we have $\dim_{E_\lambda}\rH^1_f(F,\rho_{\Pi_0,\lambda}\otimes_{E_\lambda}\rho_{\Pi_1,\lambda}(n))\geq 1$.

We start to prove the theorem by contradiction, hence assume
\[
\dim_{E_\lambda}\rH^1_f(F,\rho_{\Pi_0,\lambda}\otimes_{E_\lambda}\rho_{\Pi_1,\lambda}(n))\geq 2.
\]

Take a sufficiently large positive integer $m$ which will be determined later. We fix a uniformizer $\lambda_0$ of $E_\lambda$. By Lemma \ref{le:purity}, we may apply Proposition \ref{pr:selmer_reduction} by taking $\Sigma$ to be the set of places of $F$ above $\Sigma^+_\mnm\cup\Sigma^+_\ell$. Then we obtain a submodule $S$ of $\rH^1_{f,\rR}(F,\bar\rR^{(m)})$ containing (the image of) $\lambda_0^{m_\Sigma-m_{\r{per}}}s$ of order $0$,\footnote{Here, $\lambda_0^{-m_{\r{per}}}s$ is any element in $\rH^1_f(F,\rR)$ satisfying $\lambda_0^{m_{\r{per}}}(\lambda_0^{-m_{\r{per}}}s)=s$.} that is free of rank $2$ over $O_\lambda/\lambda^{m-m_\Sigma}$, and such that $\loc_w\res_S=0$ for every nonarchimedean place $w\in\Sigma$ not above $\ell$. Now we apply the discussion in \S\ref{ss:galois_lemma} to the submodule $S\subseteq\rH^1(F,\bar\rR^{(m)})$. By (L5-1) and Lemma \ref{le:image}, we obtain an injective map
\[
\theta_S\colon\Gal(F_S/F_{\bar\rho^{(m)}})\to\Hom_{O_\lambda}(S,\bar\rR^{(m)})
\]
whose image generates an $O_\lambda$-submodule containing $\lambda^{4\fr_{\bar\rR^{(m)}}}\Hom_{O_\lambda}(S,\bar\rR^{(m)})$, which further contains $\lambda^{4\fr_\rR}\Hom_{O_\lambda}(S,\bar\rR^{(m)})$ by Lemma \ref{le:reducibility} and (L3). By (L5-2) and Lemma \ref{le:galois_element}, we may choose an element $(\gamma_1,\gamma_2,\xi)$ in the image of $(\bar\rho^{(m)}_{1+},\bar\rho^{(m)}_{2+},\bar\epsilon_\ell^{(m)})\res_{\Gal(\ol{F}/F^+_{\r{rflx}})}$ satisfying (a--e) in Lemma \ref{le:galois_element}. It then gives rise to an element $\gamma\in(\GL_{n_0n_1}(O_\lambda/\lambda^m)\times(O_\lambda/\lambda^m)^\times)\fc$ as in Notation \ref{no:h_gamma} such that $(\bar\rR^{(m)})^{h_\gamma}$ is a free $O_\lambda/\lambda^m$-module of rank $1$. Now we apply the discussion in \S\ref{ss:localization}. By Proposition \ref{pr:selmer_localization} (with $m_0=m_\Sigma$ and $r_S=2$), we may fix an $(S,\gamma)$-abundant pair $(\Psi_1,\Psi_2)\in G_{S,\gamma}^2$ (Definition \ref{de:abundant}). By Proposition \ref{co:selmer_localization}, we may choose a basis $\{s_1,s_2\}$ of $S$ such that $\theta_S(\Psi_1)(s_2)=\theta_S(\Psi_2)(s_1)=0$, and
\begin{align}\label{eq:selmer1_1}
\exp_\lambda\(\theta_S(\Psi_j)(s_j),(\bar\rR^{(m)})^{h_\gamma}\)\geq m-m_\Sigma-4\fr_\rR
\end{align}
for $j=1,2$. Moreover, without loss of generality, we may assume $\lambda_0^{m_\Sigma-m_{\r{per}}}s=a_1s_1+a_2s_2$ in which $a_1\in O_\lambda^\times$.

First, we apply the discussion and notation in \S\ref{ss:second_reciprocity} to our situation with $\lambda$, $m$, $\Sigma^+_{\lr,\r{II}}=\emptyset$, $\Sigma^+_{\r{II}}=\Sigma^+_\r{min}$, $(\rV_n,\Lambda_n)$, $\rK_n$ and $(\rK_n,\rK_{n+1})$. By the Chebotarev density theorem, we can choose a $\gamma$-associated place (Definition \ref{de:associated}) $w^{(m)}_{1+}$ of $F^{(m)}_+$ satisfying $\Psi_{w^{(m)}_1}=\Psi_1$ and whose underlying prime $\fp_1$ of $F^+$ (and the underlying rational prime $p_1$) is a special inert prime satisfying (PII1)--(PII7) and
\begin{description}
  \item[(PII8)] the natural map
     \[
     \rH^i_\et(\Sh(\rV_{n_1},\rK_{n_1})_{\ol{F}},O_\lambda(r_1))/(\dT^{\Sigma^+_{\r{II}}\cup\Sigma^+_{p_1}}_{n_1}\cap\Ker\phi_{\Pi_1})\to
     \rH^i_\et(\Sh(\rV_{n_1},\rK_{n_1})_{\ol{F}},O_\lambda(r_1))/\Ker\phi_{\Pi_1}
     \]
     is an isomorphism for every integer $i$.
\end{description}
We also choose remaining data in \S\ref{ss:qs_initial} with $\dQ_{p_1}^\Phi=\dQ_{p_1^2}$, a definite uniformization datum $(\rV^\star_{n_\alpha},\ti_{n_\alpha},\{\Lambda^\star_{n_\alpha,\fq}\}_{\fq\mid p_1})$ for $\alpha=0,1$ as in Notation \ref{no:qs_uniformization}. By \eqref{eq:selmer1_1} and our choice of $S$, we have
\[
\exp_\lambda\(s,\rH^1_\ns(F_{w_1},\bar\rR^{(m)})\)\geq m-m_{\r{per}}-4\fr_\rR,
\]
which implies that
\[
\resizebox{\hsize}{!}{
\xymatrix{
\exp_\lambda\(\loc_{\fp_1}([\graph\Sh(\rV_n,\rK_n)]),
\rH^{2n}_\et((\Sh(\rV_{n_0},\rK_{n_0})\times_{\Spec{F}}\Sh(\rV_{n_1},\rK_{n_1}))_{F_{\fp_1}},L(n))/(\fn_0,\fn_1)\)\geq m-m_{\r{per}}-4\fr_\rR.
}
}
\]
Here, we recall that
\[
\fn_\alpha=\dT^{\Sigma^+_{\r{II}}\cup\Sigma^+_{p_1}}_{n_\alpha}\cap\Ker\(\dT^{\Sigma^+_\mnm}_{n_\alpha}\xrightarrow{\phi_{\Pi_\alpha}}O_E\to O_E/\lambda^m\)
\]
for $\alpha=0,1$. Assumption \ref{as:first_vanishing} for $\alpha=0$ is satisfied by Lemma \ref{le:ns_pbc} and (L7) when $n\geq 3$, and by Lemma \ref{le:ns_pbc} and (L4) when $n=2$. Assumption \ref{as:first_vanishing} for $\alpha=1$ is satisfied by Lemma \ref{le:ns_pbc}, (PII8), and the fact that $\rH^i_\et(\Sh(\rV_{n_1},\rK_{n_1})_{\ol{F}},O_E/\lambda)/\Ker\phi_{\Pi_1}=0$ if $i\neq 2r_1$. Thus, we may apply Theorem \ref{th:second}, hence obtain
\begin{align}\label{eq:selmer1_0}
\exp_\lambda\(\CF_{\Sh(\rV^\star_n,\rK^\star_\sp)},
O_E[\Sh(\rV^\star_{n_0},\rK^\star_{n_0})\times\Sh(\rV^\star_{n_1},\rK^\star_{n_1})]/(\fn_0,\fn_1)\)\geq m-m_{\r{per}}-4\fr_\rR.
\end{align}

Second, we apply the discussion and notation in \S\ref{ss:first_reciprocity} to our situation with $\lambda$, $m$, $\Sigma^+_{\lr,\rI}=\{\fp_1\}$, $\Sigma^+_\rI=\Sigma^+_\mnm\cup\Sigma^+_{p_1}$, $\rV^\circ_n=\rV^\star_n$, $\rK^\circ_n=\rK^\star_n$ and $(\rK^\circ_\sp,\rK^\circ_{n+1})=(\rK^\star_\sp,\rK^\star_{n+1})$. By the Chebotarev density theorem, we can choose a $\gamma$-associated place $w^{(m)}_{2+}$ of $F^{(m)}_+$ satisfying $\Psi_{w^{(m)}_2}=\Psi_2$ and whose underlying prime $\fp_2$ of $F^+$ (and the underlying rational prime $p_2$) is a special inert prime satisfying (PI1)--(PI7), $p_2\neq p_1$, and
\begin{description}
  \item[(PI8)] the natural map
     \[
     \rH^{2r_1}_\et(\Sh(\rV_{n_1},\rK_{n_1})_{\ol{F}},O_\lambda(r_1))
     /(\dT^{\Sigma^+_\rI\cup\Sigma^+_{p_2}}_{n_1}\cap\Ker\phi_{\Pi_1})\to
     \rH^{2r_1}_\et(\Sh(\rV_{n_1},\rK_{n_1})_{\ol{F}},O_\lambda(r_1))/\Ker\phi_{\Pi_1}
     \]
     is an isomorphism.
\end{description}

We claim that there exists an element $c_2\in\rH^1(F,\bar\rR^{(m_2)\tc})$ for some positive integer $m_2\leq m$ satisfying \begin{align}\label{eq:selmer1_4}
\exp_\lambda\(\partial_{\fp_2}\loc_{\fp_2}(c_2),\rH^1_\sing(F_{\fp_2},\bar\rR^{(m_2)\tc})\)\geq m-m_{\r{per}}-4\fr_\rR-m_{\r{lat}};
\end{align}
and such that for every nonarchimedean place $w$ of $F$ not above $\Sigma^+\cup\{\fp_1,\fp_2\}$,
\begin{align}\label{eq:selmer1_6}
\loc_w(c_2)\in\rH^1_\ns(F_w,\bar\rR^{(m_2)\tc})
\end{align}
holds.

By Remark \ref{re:qs_tate_hecke} and Remark \ref{re:qs_uniformization}, we know that there exists an isomorphism $\rU((\rV^\circ_{n_1})^\infty)\simeq\rU(\rV_{n_1}^\infty)$ sending $\rK^\circ_{n_1}$ to $\rK_{n_1}$. Then the claim can be proved by the exactly same argument for the similar claim in the proof of Theorem \ref{th:selmer0}, using \eqref{eq:selmer1_0} and the fact that $\bar\rho_{\Pi_0,\lambda,+}$ is rigid for $(\Sigma^+_\mnm,\Sigma^+_{\lr,\rI})$.\footnote{In fact, one needs to use the additional fact that when $F^+\neq\dQ$, both Shimura varieties $\Sh'_{n_0}$ and $\Sh'_{n_1}$ have proper smooth reduction at every place $w$ of $F$ above $\Sigma^+_{p_1}\setminus\{\fp_1\}$. See Remark \ref{re:ns_smooth}.}

Now we deduce a contradiction. Replace $s_2$ by its image in $\rH^1_f(F,\bar\rR^{(m_2)})$. We also identify $\bar\rR^{(m_2)\tc}$ with $(\bar\rR^{(m_2)})^*$ via the polarization $\Xi$. Now we compute the local Tate pairing $\langle s_2,c_2\rangle_w$ \eqref{eq:local_tate} for every nonarchimedean place $w$ of $F$.
\begin{itemize}[label={\ding{109}}]
  \item Suppose that $w$ is above $\Sigma^+_\mnm$. Then we have $\loc_w(s_2)=0$ by our choice of $S$. Thus, $\langle s_2,c_2\rangle_w=0$.

  \item Suppose that $w$ is above $\Sigma^+_\ell$. Then by (L2), $\rR_\dQ$ is crystalline with Hodge--Tate weights in $[1-n,n]$. Thus, we have $\loc_w(s_2)\in\rH^1_\ns(F_w,\bar\rR^{(m_2)})$ by Lemma \ref{le:global_local}(2) and (L1). By \eqref{eq:selmer1_6}, Lemma \ref{le:crystalline_tate} and (L1), we have $\lambda^{m_{\r{dif}}}\langle s_2,c_2\rangle_w=0$ where $\fd_\lambda=\lambda^{m_{\r{dif}}}\subseteq O_\lambda$ is the different ideal of $E_\lambda$ over $\dQ_\ell$.

  \item Suppose that $w$ is not above $\Sigma^+_\mnm\cup\Sigma^+_\ell\cup\{\fp_1,\fp_2\}$. Then by (L2), $\rR$ is unramified. Thus, we have $\loc_w(s_2)\in\rH^1_\ns(F_w,\bar\rR^{(m_2)})$ by Lemma \ref{le:global_local}(1). By \eqref{eq:selmer1_6} and Lemma \ref{le:tate}, we have $\langle s_2,c_2\rangle_w=0$.

  \item Suppose that $w$ is the unique place above $\fp_1$. Then we have $\loc_w(s_2)=0$ by Proposition \ref{co:selmer_localization}. Thus, we have $\langle s_2,c_2\rangle_w=0$.

  \item Suppose that $w$ is the unique place above $\fp_2$. Then by Proposition \ref{co:selmer_localization}, we have
      \[
      \exp_\lambda\(\loc_w(s_2),\rH^1_\ns(F_w,\bar\rR^{(m_2)})\)\geq m_2-m_\Sigma-4\fr_\rR.
      \]
      By \eqref{eq:selmer1_4} and Lemma \ref{le:tate} again, we have
      \[
      \exp_\lambda\(\langle s_2,c_2\rangle_w,O_\lambda/\lambda^{m_2}\)\geq m-m_{\r{per}}-m_{\r{lat}}-m_\Sigma-8\fr_\rR.
      \]
\end{itemize}
Therefore, as long as we take $m$ such that $m>m_{\r{per}}+m_{\r{lat}}+m_\Sigma+8\fr_\rR+m_{\r{dif}}$, we will have a contradiction to the relation
\[
\sum_{w}\langle s_2,c_2\rangle_w=0,
\]
where the sum is taken over all nonarchimedean places $w$ of $F$. The theorem is proved.
\end{proof}

We also have an analogue of Corollary \ref{co:abstract} in the rank $1$ case, which we leave to the readers to formulate.

\fi

\appendix

\section{Unitary Deligne--Lusztig varieties}
\label{ss:dl}

In this appendix, we study some unitary Deligne--Lusztig varieties in \S\ref{ss:dl_smooth} and \S\ref{ss:dl_semistable} for those used in \S\ref{ss:qs} and \S\ref{ss:ns}, respectively.

We fix a rational prime $p$. Let $\kappa$ be a field containing $\dF_{p^2}$. Recall from \S\ref{ss:notation} that we denote by $\sigma\colon S\to S$ the absolute $p$-power Frobenius morphism for schemes $S$ in characteristic $p$.

\subsection{Unitary Deligne--Lusztig varieties in the smooth case}
\label{ss:dl_smooth}

In this subsection, we introduce certain Deligne--Lusztig varieties that appear in the special fiber of the smooth integral model studied in \S\ref{ss:qs}.

Consider a pair $(\sV,\{\;,\;\})$ in which $\sV$ is a finite dimensional $\kappa$-linear space, and $\{\;,\;\}\colon\sV\times\sV\to\kappa$ is a (not necessarily non-degenerate) pairing that is $(\kappa,\sigma)$-linear in the first variable and $\kappa$-linear in the second variable. For every $\kappa$-scheme $S$, put $\sV_S\coloneqq\sV\otimes_\kappa\cO_S$. Then there is a unique pairing $\{\;,\;\}_S\colon\sV_S\times\sV_S\to\cO_S$ extending $\{\;,\;\}$ that is $(\cO_S,\sigma)$-linear in the first variable and $\cO_S$-linear in the second variable. For a subbundle $H\subseteq\sV_S$, we denote by $H^\dashv\subseteq\sV_S$ its \emph{right} orthogonal complement under $\{\;,\;\}_S$.

\begin{definition}\label{de:dl_admissible}
We say that a pair $(\sV,\{\;,\;\})$ is \emph{admissible} if there exists an $\dF_{p^2}$-linear subspace $\sV_0\subseteq\sV_{\ol\kappa}$ such that the induced map $\sV_0\otimes_{\dF_{p^2}}\ol\kappa\to\sV_{\ol\kappa}$ is an isomorphism, and $\{x,y\}=-\{y,x\}^\sigma$ for every $x,y\in\sV_0$.
\end{definition}

\begin{definition}\label{de:dl}
For a pair $(\sV,\{\;,\;\})$ and an integer $h$, we define a presheaf
\[
\DL(\sV,\{\;,\;\},h)
\]
on $\Sch_{/\kappa}$ such that for every $S\in\Sch_{/\kappa}$, $\DL(\sV,\{\;,\;\},h)(S)$ is the set of subbundles $H$ of $\sV_S$ of rank $h$ such that $H^\dashv\subseteq H$. We call $\DL(\sV,\{\;,\;\},h)$ the \emph{(unitary) Deligne--Lusztig variety} (see Proposition \ref{pr:dl} below) attached to $(\sV,\{\;,\;\})$ of rank $h$.
\end{definition}

\begin{proposition}\label{pr:dl}
Consider an admissible pair $(\sV,\{\;,\;\})$. Put $N\coloneqq\dim_\kappa\sV$ and $d\coloneqq\dim_\kappa\sV^\dashv$.
\begin{enumerate}
  \item If $2h<N+d$ or $h>N$, then $\DL(\sV,\{\;,\;\},h)$ is empty.

  \item If $N+d\leq 2h\leq 2N$, then $\DL(\sV,\{\;,\;\},h)$ is represented by a projective smooth scheme over $\kappa$ of dimension $(2h-N-d)(N-h)$ with a canonical isomorphism for its tangent sheaf
      \[
      \cT_{\DL(\sV,\{\;,\;\},h)/\kappa}\simeq\HOM\(\cH/\cH^\dashv,\sV_{\DL(\sV,\{\;,\;\},h)}/\cH\)
      \]
      where $\cH\subseteq\sV_{\DL(\sV,\{\;,\;\},h)}$ is the universal subbundle.

  \item If $N+d< 2h\leq 2N$, then $\DL(\sV,\{\;,\;\},h)$ is geometrically irreducible.
\end{enumerate}
\end{proposition}

\begin{proof}
Part (1) is obvious from the definitions.

For (2), $\DL(\sV,\{\;,\;\},h)$ is a closed sub-presheaf of the Grassmannian scheme $\r{Gr}(\sV,h)$ classifying subbundles of $\sV$ of rank $h$, hence is represented by a projective scheme over $\kappa$. Now we compute the tangent sheaf. Consider a closed immersion $S\hookrightarrow\hat{S}$ in $\Sch_{/\kappa}$ defined by an ideal sheaf $\cI$ with $\cI^2=0$. Take an object $H\subseteq\sV_S$ in $\DL(\sV,\{\;,\;\},h)(S)$. Let $D_H$ and $G_H$ be the subset of $\DL(\sV,\{\;,\;\},h)(\hat{S})$ and $\r{Gr}(\sV,h)(\hat{S})$ of elements that reduce to $H$, respectively. It is well-known that $G_H$ is canonically a torsor over $\Hom_{\cO_S}(H,(\sV_S/H)\otimes_{\cO_S}\cI)$. Since $\cI^p=0$, the right orthogonal complement $\hat{H}^\dashv$ depends only on $H$ for every $\hat{H}\in G_H$. In particular, the subset $D_H$ is canonically a torsor over the subgroup $\Hom_{\cO_S}(H/H^\dashv,(\sV_S/H)\otimes_{\cO_S}\cI)$ of $\Hom_{\cO_S}(H,(\sV_S/H)\otimes_{\cO_S}\cI)$. Thus, $\DL(\sV,\{\;,\;\},h)$ is smooth; and we have a canonical isomorphism for the tangent sheaf
\[
\cT_{\DL(\sV,\{\;,\;\},h)/\kappa}\simeq\HOM\(\cH/\cH^\dashv,\sV_{\DL(\sV,\{\;,\;\},h)}/\cH\)
\]
where $\cH$ is the universal subbundle. Note that this is a locally free $\cO_{\DL(\sV,\{\;,\;\},h)}$-module of rank $(2h-N-d)(N-h)$.

For (3), we may assume that $\kappa$ is algebraically closed. By Definitions \ref{de:dl_admissible} and \ref{de:dl}, we have a canonical isomorphism $\DL(\sV,\{\;,\;\},h)\simeq\DL(\sV_0,\{\;,\;\}_0,h)\otimes_{\dF_{p^2}}\kappa$, where $\{\;,\;\}_0$ denotes the restriction of $\{\;,\;\}$ to $\sV_0$. Suppose that $d=0$. Then $\{\;,\;\}_0$ is non-degenerate. By \cite{BR06}*{Theorem~1}, we know that $\DL(\sV_0,\{\;,\;\}_0,h)$ is geometrically irreducible. In general, we consider $\sV'_0\coloneqq\sV_0/\sV_0^\dashv$ equipped with a pairing $\{\;,\;\}'_0$ induced from $\{\;,\;\}_0$. Then it is clear that the morphism $\DL(\sV_0,\{\;,\;\}_0,h)\to\DL(\sV'_0,\{\;,\;\}'_0,h)$ sending a point $H\in\DL(\sV_0,\{\;,\;\}_0,h)(S)$ to $H/\sV_{0S}^\dashv$ is an isomorphism. Thus, $\DL(\sV_0,\{\;,\;\}_0,h)$ is geometrically irreducible by the previous case. The proposition is proved.
\end{proof}

\begin{lem}\label{le:dl_cohomology}
Consider a pair $(\sV,\{\;,\;\})$ with $\dim_\kappa\sV=N\geq 2$ and $\dim_\kappa\sV^\dashv=0$, and a $p$-coprime coefficient ring $L$. Suppose that $p+1$ is invertible in $L$.
\begin{enumerate}
  \item The subscheme $\DL(\sV,\{\;,\;\},N-1)$ is a hypersurface in $\dP(\sV)$ of degree $p+1$.

  \item The restriction map
     \[
     \rH^i_\et(\dP(\sV)_{\ol\kappa},L)\to\rH^i_\et(\DL(\sV,\{\;,\;\},N-1)_{\ol\kappa},L)
     \]
     induced by the obvious inclusion $\DL(\sV,\{\;,\;\},N-1)\to\dP(\sV)$ is an isomorphism for $i\not\in\{N-2,2N-2\}$.

  \item For every $i\in\dZ$, $\rH^i_\et(\DL(\sV,\{\;,\;\},N-1)_{\ol\kappa},L)$ is a free $L$-module.

  \item When $N$ is even, the action of $\Gal(\ol\kappa/\kappa)$ on $\rH^{N-2}_\et(\DL(\sV,\{\;,\;\},N-1)_{\ol\kappa},L(\tfrac{N-2}{2}))$ is trivial.
\end{enumerate}
\end{lem}

\begin{proof}
The lemma is trivial if $N=2$. Now we assume $N\geq 3$. Then $S\coloneqq\DL(\sV,\{\;,\;\},N-1)$ is a geometrically connected smooth hypersurface in $\dP(\sV)$ by Proposition \ref{pr:dl}.

Part (1) follows since $S$ is defined by a homogenous polynomial of degree $p+1$, by its definition.

For (2), by the Lefschetz hyperplane theorem, the restriction map $\rH^i_\et(\dP(\sV)_{\ol\kappa},L)\to\rH^i_\et(S_{\ol\kappa},L)$ is an isomorphism for $0\leq i\leq N-3$; and the Gysin map $\rH^i_\et(S_{\ol\kappa},L)\to\rH^{i+2}_\et(\dP(\sV)_{\ol\kappa},L(1))$ is an isomorphism for $N-1\leq i\leq 2(N-2)$. By (1), the composite map
\[
\rH^i_\et(\dP(\sV)_{\ol\kappa},L)\to\rH^i_\et(S_{\ol\kappa},L)\to\rH^{i+2}_\et(\dP(\sV)_{\ol\kappa},L(1))
\]
is given by the cup product with $c_1(\cO_{\dP(\sV)_{\ol\kappa}}(p+1))$, which is an isomorphism for $i\neq 2N-2$ since $p+1$ is invertible in $L$. Thus, (2) follows.

Part (3) is an immediate consequence of (2).

For (4), it suffices to consider the case where $L=\dQ_\ell$ for some $\ell\neq p$ by (3). Then it is well-known that $\rH^{N-2}_\et(\DL(\sV,\{\;,\;\},N-1)_{\ol\kappa},\dQ_\ell(\tfrac{N-2}{2}))$ is spanned by Tate cycles over $\kappa$ (see, for example, \cite{HM78}). In particular, (4) follows.
\end{proof}

\begin{proposition}\label{pr:dl_excess_pre}
Suppose that $\kappa$ is algebraically closed. Consider an admissible pair $(\sV,\{\;,\;\})$ over $\kappa$ with $\dim_\kappa\sV=2r+1$ for some integer $r\geq 1$ and $\dim_\kappa\sV^\dashv=0$. Let $\cH$ be the universal object over $\DL(\sV,\{\;,\;\},r+1)$. Then we have
\[
\int_{\DL(\sV,\{\;,\;\},r+1)}c_r\(\(\sigma^*\cH^{\vdash}\)\otimes_{\DL(\sV,\{\;,\;\},r+1)}\(\cH/\cH^{\vdash}\)\)
=\td_{r,p},
\]
where $\td_{r,p}$ is the number introduced in Notation \ref{no:numerical}.
\end{proposition}

\begin{proof}
This is \cite{XZ17}*{Proposition~9.3.10}.
\end{proof}

Now we construct the special morphisms between Deligne--Lusztig varieties when rank increases.

\begin{construction}\label{cs:dl_special_morphism}
Let $(\sV,\{\;,\;\})$ be an admissible pair with $\dim_\kappa\sV=n\geq 1$ satisfying $\dim\sV^\dashv=n+1-2\floor{\tfrac{n+1}{2}}$. We put $\sV_\sharp\coloneqq\sV\oplus\kappa1$ and extend $\{\;,\;\}$ to a pairing $\{\;,\;\}_\sharp$ on $\sV_\sharp$ with $\{1,1\}_\sharp=0$. Suppose that we have another admissible pair $(\sV_\natural,\{\;,\;\}_\natural)$ with $\dim_\kappa\sV_\natural=n+1$ satisfying $\dim\sV_\natural^\dashv=n-2\floor{\tfrac{n}{2}}$, together with a $\kappa$-linear map $\delta\colon\sV_\sharp\to\sV_\natural$ of corank $\dim\sV^\dashv$ such that $\{\delta(x),\delta(y)\}_\natural=\{x,y\}_\sharp$ for every $x,y\in\sV_\sharp$. We construct a morphism
\[
\delta_\uparrow\colon\DL(\sV,\{\;,\;\},\ceil{\tfrac{n+1}{2}})\to\DL(\sV_\natural,\{\;,\;\}_\natural,\ceil{\tfrac{n+2}{2}})
\]
by sending $H\in\DL(\sV,\{\;,\;\},\ceil{\tfrac{n+1}{2}})(S)$ to $\delta(H\oplus\cO_S1)$. We call $\delta_\uparrow$ a \emph{special morphism}.
\end{construction}

\begin{proposition}\label{pr:dl_special_morphism}
The morphism $\delta_\uparrow$ is well-defined, and is a regular embedding.
\end{proposition}

\begin{proof}
When $n$ is odd, $\delta$ is an isomorphism, which implies that $\delta_\uparrow$ is well-defined an is an isomorphism.

When $n$ is even, $\delta$ is of corank $1$. The identity $\{\delta(x),\delta(y)\}_\natural=\{x,y\}_\sharp$ for every $x,y\in\sV_\sharp$ implies $\Ker\delta\subset\sV_{\sharp}^{\dashv}=\sV^{\dashv}\oplus\kappa1$. Take $S\in\Sch_{/\kappa}$. For $H\in\DL(\sV,\{\;,\;\},\ceil{\tfrac{n+1}{2}})(S)$, $H\oplus\cO_S1$ must contain $\sV_\sharp^\dashv$ and hence $(\Ker\delta)_S$. It follows that $\delta(H\oplus\cO_S1)$ has the same rank as $H$, which is $\ceil{\tfrac{n+1}{2}}=\ceil{\tfrac{n+2}{2}}$. The identity $\{\delta(x),\delta(y)\}_\natural=\{x,y\}_\sharp$ for every $x,y\in\sV_\sharp$ also implies $\delta(H^{\dashv}\oplus\cO_S1)\subseteq(\delta(H\oplus\cO_S1))^\dashv$, which forces $\delta(H^{\dashv}\oplus\cO_S1)=(\delta(H\oplus\cO_S1))^\dashv$ as both sides have the same rank $\tfrac{n}{2}$. It follows that $(\delta(H\oplus\cO_S1))^\dashv\subseteq\delta(H\oplus\cO_S1)$ as $H^{\dashv}\subseteq H$. In other words, $\delta_\uparrow$ is well-defined. On the other hand, for $H_\natural\in\DL(\sV_\natural,\{\;,\;\}_\natural,\ceil{\tfrac{n+2}{2}})(S)$, whether $(\delta\kappa1)_S\subseteq H\subseteq(\delta\sV_\sharp)_S$ holds is a closed condition; and once it does, there is a unique element $H\in\DL(\sV,\{\;,\;\},\ceil{\tfrac{n+1}{2}})(S)$ such that $H_\natural=\delta(H\oplus\cO_S1)$. Thus, $\delta_\uparrow$ is a regular embedding by Proposition \ref{pr:dl}(2).

The proposition is proved.
\end{proof}

\subsection{Unitary Deligne--Lusztig varieties in the semistable case}
\label{ss:dl_semistable}

In this subsection, we introduce certain Deligne--Lusztig varieties that appear in the special fiber of the semistable integral model studied in \S\ref{ss:ns}. We keep the notation from the previous subsection.

\begin{definition}\label{de:dl_bullet}
For a pair $(\sV,\{\;,\;\})$ with $\dim_\kappa\sV=N$, we define a presheaf
\[
\DL^\bullet(\sV,\{\;,\;\})
\]
on $\Sch_{/\kappa}$ such that for every $S\in\Sch_{/\kappa}$, $\DL^\bullet(\sV,\{\;,\;\})(S)$ is the set of pairs $(H_1,H_2)$ of subbundles of $\sV_S$ of ranks $\ceil{\tfrac{N}{2}}$ and $\ceil{\tfrac{N}{2}}-1$, respectively, satisfying the following inclusion relations
\[
\xymatrix{
&& H_1  \ar@{}[dr]|-*[@]{\subset} \\
\sV_S^\dashv \ar@{}[r]|-*[@]{\subset}  & H_2 \ar@{}[ur]|-*[@]{\subset}\ar@{}[dr]|-*[@]{\subset}  && H_2^\dashv \\
&& H_1^\dashv \ar@{}[ur]|-*[@]{\subset}
}
\]
of subbundles of $\sV_S$.
\end{definition}

\begin{proposition}\label{pr:dl_bullet}
Consider an admissible pair $(\sV,\{\;,\;\})$. Put $N\coloneqq\dim_\kappa\sV$ and $d\coloneqq\dim_\kappa\sV^\dashv$.
\begin{enumerate}
  \item If $d\geq \ceil{\tfrac{N}{2}}$, then $\DL^\bullet(\sV,\{\;,\;\})$ is empty.

  \item If $d\leq \ceil{\tfrac{N}{2}}-1$, then $\DL^\bullet(\sV,\{\;,\;\})$ is represented by a projective smooth scheme over $\kappa$, whose tangent sheaf fits canonically into an exact sequence
      \[
      0\to\HOM\(\cH_1/\cH_2,\cH_2^\dashv/\cH_1\)\to\cT_{\DL^\bullet(\sV,\{\;,\;\})/\kappa}\to
      \HOM(\cH_2/\sV^\dashv_{\DL^\bullet(\sV,\{\;,\;\})},\cH_1^\dashv/\cH_2)\to 0
      \]
      where $\sV^\dashv_{\DL^\bullet(\sV,\{\;,\;\})}\subseteq\cH_2\subseteq\cH_1\subseteq\sV_{\DL^\bullet(\sV,\{\;,\;\})}$ are the universal subbundles.

  \item If $N\geq 2$ and $d=N-2\floor{\tfrac{N}{2}}$, then $\DL^\bullet(\sV,\{\;,\;\})$ is geometrically irreducible of dimension $\floor{\tfrac{N}{2}}$.
\end{enumerate}
\end{proposition}

\begin{proof}
Part (1) is obvious from the definitions.

For (2), let $\r{Gr}(\sV, r)$ denote by the Grassmannian variety that classifies subspaces of $\sV$ of dimension $r$. Then $\DL^\bullet(\sV,\{\;,\;\})$ is a closed sub-presheaf of $\r{Gr}(\sV,\ceil{\tfrac{N}{2}})\times\r{Gr}(\sV,\ceil{\tfrac{N}{2}}-1)$, hence it is represented by a projective scheme over $\kappa$. Now we prove that $\DL^\bullet(\sV,\{\;,\;\})$ is smooth and compute its tangent sheaf. Consider a closed immersion $S\hookrightarrow\hat{S}$ in $\Sch_{/\kappa}$ defined by an ideal sheaf $\cI$ with $\cI^2=0$. Take an object $\sV_S^\dashv\subseteq H_2\subseteq H_1\subseteq\sV_S$ in $\DL^\bullet(\sV,\{\;,\;\})(S)$. To lift $(H_1,H_2)$ to a pair $(\hat{H}_1,\hat{H}_2)\in\DL^\bullet(\sV,\{\;,\;\})(\hat{S})$, we first lift $H_2$, where the set of all possible lifts canonically form a torsor under the group $\Hom_{\cO_S}(H_2/\sV_S^\dashv,(H_1^\dashv/H_2)\otimes_{\cO_S}\cI)$ as $\hat{H}_1^\dashv$ depends only on $H_1^\dashv$. Once such a lift $\hat{H}_2$ is given, the possible lifts of $H_1$ form a torsor under the group $\Hom_{\cO_S}(H_1/H_2,(H_2^\dashv/H_1)\otimes_{\cO_S}\cI)$.  In particular, Zariski locally, there is no obstruction to lifting $(H_1, H_2)$, hence $\DL^\bullet(\sV,\{\;,\;\})$ is smooth. The statement on the tangent bundle of $\DL^\bullet(\sV,\{\;,\;\})$ follows immediately from the above discussion applied to the universal object on $\DL^\bullet(\sV,\{\;,\;\})$.

For (3), similar to the argument for Proposition \ref{pr:dl}(3), we may assume that $N$ is even this time. Then the statement follows again by \cite{BR06}*{Theorem~1}.
\end{proof}

\begin{construction}\label{cs:dl_bullet_special_morphism}
Let $(\sV,\{\;,\;\})$ be an admissible pair with $\dim_\kappa\sV=n\geq 2$ satisfying $\dim_\kappa\sV^\dashv=n-2\floor{\tfrac{n}{2}}$. We put $\sV_\sharp\coloneqq\sV\oplus\kappa1$ and extend $\{\;,\;\}$ to a pairing $\{\;,\;\}_\sharp$ on $\sV_\sharp$ with $\{1,1\}_\sharp=0$. Suppose that we have another admissible pair $(\sV_\natural,\{\;,\;\}_\natural)$ with $\dim_\kappa\sV_\natural=n+1$ satisfying $\dim\sV_\natural^\dashv=n+1-2\floor{\tfrac{n+1}{2}}$, together with a $\kappa$-linear map $\delta\colon\sV_\sharp\to\sV_\natural$ of corank $\dim\sV^\dashv$ such that $\{\delta(x),\delta(y)\}_\natural=\{x,y\}_\sharp$ for every $x,y\in\sV_\sharp$. Then similar to Construction \ref{cs:dl_special_morphism} and Proposition \ref{pr:dl_special_morphism}, we have a morphism
\[
\delta_\uparrow\colon\DL^\bullet(\sV,\{\;,\;\})\to\DL^\bullet(\sV_\natural,\{\;,\;\}_\natural)
\]
by sending $(H_1,H_2)\in\DL^\bullet(\sV,\{\;,\;\})(S)$ to $(\delta(H_1\oplus\cO_S1),\delta(H_2\oplus\cO_S1))\in\DL^\bullet(\sV_\natural,\{\;,\;\}_\natural)(S)$, which is a regular embedding.
\end{construction}

\begin{proposition}\label{pr:dl_excess}
Suppose that $\kappa$ is algebraically closed. Consider an admissible pair $(\sV,\{\;,\;\})$ over $\kappa$. Let $(\cH_1,\cH_2)$ be the universal object over $\DL^\bullet(\sV,\{\;,\;\})$.
\begin{enumerate}
  \item Suppose that $\dim_\kappa\sV=2r+1$ for some integer $r\geq 1$ and $\dim_\kappa\sV^\dashv=1$. Then we have
    \[
    \int_{\DL^\bullet(\sV,\{\;,\;\})}c_r\(\(\sigma^*\cH_2\)\otimes_{\cO_{\DL^\bullet(\sV,\{\;,\;\})}}\(\cH_1^\dashv/\cH_2\)\)
    =\td^\bullet_{r,p}.
    \]

  \item Suppose that $\dim_\kappa\sV=2r$ for some integer $r\geq 1$ and $\dim_\kappa\sV^\dashv=0$. Then we have
    \[
    \int_{\DL^\bullet(\sV,\{\;,\;\})}c_{r-1}\(\(\sigma^*\cH_2\)\otimes_{\cO_{\DL^\bullet(\sV,\{\;,\;\})}}\(\cH_1^\dashv/\cH_2\)\)\cdot c_1\(\cH_1^\dashv/\cH_2\)=\td^\bullet_{r,p}.
    \]
\end{enumerate}
Here, $\td^\bullet_{r,p}$ is the number introduced in Notation \ref{no:numerical}.
\end{proposition}

Note that $\DL^\bullet(\sV,\{\;,\;\})$ is irreducible of dimension $r$, by Proposition~\ref{pr:dl_bullet}.

\begin{proof}
For (1), we let $\bar\sV$ be the quotient space $\sV/\sV^\dashv$, equipped with the induced pairing, which we still denote by $\{\;,\;\}$. Then we have a canonical isomorphism $\DL^\bullet(\sV,\{\;,\;\})\xrightarrow{\sim}\DL^\bullet(\bar\sV,\{\;,\;\})$ by sending a pair $(H_1, H_2)$ to $(H_1/\sV^{\dashv}, H_2/\sV^{\dashv})$. If we denote by $(\bar\cH_1,\bar\cH_2)$ the universal object over $\DL^\bullet(\bar\sV,\{\;,\;\})$. Then we have
\[
c_r\(\(\sigma^*\cH_2\)\otimes_{\cO_{\DL^\bullet(\sV,\{\;,\;\})}}\(\cH_1^\dashv/\cH_2\)\)
=c_{r-1}\(\(\sigma^*\bar\cH_2\)\otimes_{\cO_{\DL^\bullet(\bar\sV,\{\;,\;\})}}\(\bar\cH_1^\dashv/\bar\cH_2\)\)\cdot c_1\(\bar\cH_1^\dashv/\bar\cH_2\)
\]
under the above isomorphism. Therefore, (1) follows from (2).

For (2), consider $\sV_\sharp\coloneqq\sV\oplus\kappa1$ and extend $\{\;,\;\}$ to a pairing $\{\;,\;\}_\sharp$ on $\sV_\sharp$ with $\{1,1\}_\sharp=1$. Then we have Deligne--Lusztig varieties $\DL(\sV_\sharp,\{\;,\;\}_\sharp,h)$. In what follows, we only need to study the one with $h=r+1$, and will simply write $\DL(\sV_\sharp)$ for $\DL(\sV_\sharp,\{\;,\;\}_\sharp,r+1)$. Since we will work with two spaces, we will denote by $(\vdash,\dashv)$ for the (left,right) orthogonal complement for $\sV$, and $(\vDash,\Dashv)$ for the (left,right) orthogonal complement for $\sV_\sharp$.

We now define a correspondence
\[
\DL(\sV_\sharp) \xleftarrow{\pi} \wt\DL(\sV) \xrightarrow{\pi^\bullet} \DL^\bullet(\sV)
\]
of schemes over $\kappa$. For every $\kappa$-scheme $S$,
\begin{itemize}[label={\ding{109}}]
  \item $\wt\DL(\sV)(S)$ is the set of pairs $(H,H_2)$ where $H$ is an element in $\DL(\sV_\sharp)(S)$ and $H_2$ is a subbundle of $H^{\vDash}$ of rank $r-1$ that is contained in $\sV_S$;

  \item $\pi$ sends $(H,H_2)\in\wt\DL(\sV)(S)$ to $H\in\wt\DL(\sV)(S)$; and

  \item $\pi^\bullet$ sends $(H,H_2)\in\wt\DL(\sV)(S)$ to $(H_1,H_2)\in\DL^\bullet(\sV)(S)$ where $H_1\coloneqq (H\cap\sV_S)^\vdash$.
\end{itemize}
It needs to show that $\pi^\bullet$ is well-defined, which amounts to the following four statements:
\begin{itemize}[label={\ding{109}}]
  \item $H_1$ is a subbundle of $\sV_S$ of rank $rp$: It suffices to show that the composite map $H\to\sV_{\sharp S}\to\cO_S1$ is surjective, where the latter map is induced by the projection $\sV_\sharp\to\kappa1$. If not, then there exists a geometric point $s$ of $S$ such that $H_s$ is contained in $\sV_s$, which contradicts the inclusion $H_s^{\Dashv}\subseteq H_s$.

  \item $H_2\subseteq H_1$: As $H^{\Dashv}\subseteq H$ by the definition of $\DL(\sV_\sharp)$, we have $H^{\vDash}\subseteq H$ and  $\{H^{\vDash},H\}_\sharp=0$. Thus, $\{H^{\vDash}\cap\sV_S,H\cap\sV_S\}=0$, which implies $H_2\subseteq H^{\vDash}\cap\sV_S\subseteq(H\cap\sV_S)^\vdash=H_1$.

  \item $H_1\subseteq H_2^\vdash$: As $H^{\vDash}\subseteq H$, we have that $H_1^\dashv=H\cap\sV_S$ contains $H_2$, which implies $H_1=(H_1^\dashv)^\vdash\subseteq H_2^\vdash$.

  \item $H_1\subseteq H_2^\dashv$: As $H^{\Dashv}\subseteq H$, we have
      $(H^{\vDash})^{\Dashv\Dashv}\cap\sV_S\subseteq H\cap\sV_S$, which is equivalent to $(H^{\vDash}\cap\sV_S)^{\dashv\dashv}\subseteq H\cap\sV_S$. As $H_2$ is contained in $H^{\vDash}\cap\sV_S$, we have $H_2^{\dashv\dashv}\subseteq H\cap\sV_S=H_1^\dashv$, which implies $H_1\subseteq H_2^\dashv$.
\end{itemize}
We denote by $\cH$, $(\wt\cH,\wt\cH_2)$, and $(\cH_1,\cH_2)$ the universal objects over $\DL(\sV_\sharp)$, $\wt\DL(\sV)$, and $\DL^\bullet(\sV)$, respectively. By definition, we have $\wt\cH=\pi^*\cH$ and $\wt\cH_2=\pi^{\bullet*}\cH_2$.

We first study the morphism $\pi$. We say that a point $s\in\DL(\sV_\sharp)(\kappa)$ represented by $H_s$ is \emph{special} if $H_s^\vDash$ is a maximal isotropic subspace of $\sV$ satisfying $H_s^{\Dashv}=H_s^\vDash$. Then there are exactly $(p+1)(p^3+1)\cdots(p^{2r-1}+1)$ special points. Let $\DL(\sV_\sharp)'$ be the locus of special points. It is clear that for every morphism $S\to\DL(\sV_\sharp)\setminus\DL(\sV_\sharp)'$, $\pi^{-1}(S)$ is a singleton; and for a special point $s$, we have $\pi^{-1}(s)=\dP(H_s^\vDash)\simeq\dP^{r-1}_\kappa$. In particular, $\pi$ is a blow-up along $\DL(\sV_\sharp)'$, for which we denote by $E\subseteq\wt\DL(\sV)$ the exceptional divisor. In particular, $\pi$ is projective. Moreover, $E$ is exactly the zero locus of the canonical projection map
\[
\wt\cH^\vDash/\wt\cH_2\to\cO_{\wt\DL(\sV)}1\subseteq\cO_{\wt\DL(\sV)}\otimes_\kappa\sV_{\sharp},
\]
which implies
\begin{align}\label{eq:dl_excess_2}
\wt\cH^\vDash/\wt\cH_2\simeq\cO_{\wt\DL(\sV)}(-E).
\end{align}

Next we study the morphism $\pi^\bullet$. We claim that $\pi^\bullet$ is generically finite of degree $p+1$. Take a point $s\in\DL^\bullet(\sV)(\kappa)$ represented by $(H_{1s},H_{2s})$. Then by construction, for every scheme $S$ over $\{s\}\times_{\DL^\bullet(\sV)}\wt\DL(\sV)$, $\wt\DL(\sV)(S)$ consists of subbundles $H\subseteq\sV_\sharp\otimes_\kappa\cO_S$ satisfying $H_{2s}\otimes_\kappa\cO_S\subseteq H^{\vDash} \subseteq H_{1s}\otimes_\kappa\cO_S\oplus\cO_S1$ and $H^{\vDash}\subseteq H$. Note that we have an induced pairing
\[
\{\;,\;\}_s\colon\frac{H_{1s}\oplus\kappa 1}{H_{2s}}\times\frac{H_{1s}\oplus\kappa 1}{H_{2s}}\to\kappa
\]
that is $\sigma$-linear in the first variable and linear in the second variable. Then it is clear that when $\{\;,\;\}_s$ is perfect, $\{s\}\times_{\DL^\bullet(\sV)}\wt\DL(\sV)$ is isomorphic to the union of $p+1$ copies of $\Spec\kappa$. However, $\{\;,\;\}_s$ fails to be perfect if and only if $H_1^\dashv=H_1$. Thus, the locus where $\{\;,\;\}_s$ fails to be perfect is a finite union of $\dP^{r-1}_\kappa$. Therefore, $\pi^\bullet$ is generically finite of degree $p+1$.

To proceed, we introduce two more bundles
\[
\cE\coloneqq\(\sigma^*\cH^{\vDash}\)\otimes_{\DL(\sV_\sharp)}\(\cH/\cH^{\vDash}\),\quad
\cE^\bullet\coloneqq\(\sigma^*\cH_2\)\otimes_{\DL^\bullet(\sV)}\(\cH_1^\dashv/\cH_2\)
\]
on $\DL(\sV_\sharp)$ and $\DL^\bullet(\sV)$ of ranks $r$ and $r-1$, respectively.

We claim that
\begin{align}\label{eq:dl_excess_3}
\cL\coloneqq\pi^{\bullet*}\(\cH_1^\dashv/\cH_2\)\simeq\cO_{\wt\DL(\sV)}(-E)\otimes_{\cO_{\wt\DL(\sV)}}
\(\wt\cH/\wt\cH^{\vDash}\).
\end{align}
In fact, we have
\[
\cL=\(\wt\cH\cap\sV_{\wt\DL(\sV)}\)/\wt\cH_2
\]
by definition. Thus, the claim follows from the following injective map
\[
\xymatrix{
0 \ar[r] & \wt\cH_2 \ar[r]\ar[d] & \wt\cH^\vDash \ar[r]\ar[d] & \cO_{\wt\DL(\sV)}(-E) \ar[r]\ar[d] & 0 \\
0 \ar[r] & \wt\cH\cap\sV_{\wt\DL(\sV)} \ar[r] & \wt\cH \ar[r] & \cO_{\wt\DL(\sV)}1 \ar[r] & 0
}
\]
of short exact sequences of coherent sheaves on $\wt\DL(\sV)$ by \eqref{eq:dl_excess_2} and the Snake Lemma.

By \eqref{eq:dl_excess_2} and \eqref{eq:dl_excess_3}, we have
\begin{align*}
&\quad\pi^*\(c_r(\cE)\) \\
&=c_r\(\pi^*\cE\) \\
&=c_{r-1}\(\(\sigma^*\wt\cH_2\)\otimes_{\cO_{\wt\DL(\sV)}}\(\cH/\cH^{\vDash}\)\)
\cdot c_1\(\cO_{\wt\DL(\sV)}(-pE)\otimes_{\cO_{\wt\DL(\sV)}}\(\cH/\cH^{\vDash}\)\)  \\
&=c_{r-1}\(\(\sigma^*\wt\cH_2\)\otimes_{\cO_{\wt\DL(\sV)}}\cL(E)\)\cdot c_1(\cL((1-p)E)) \\
&=c_{r-1}\(\pi^{\bullet*}\cE^\bullet\otimes_{\cO_{\wt\DL(\sV)}}\cO_{\wt\DL(\sV)}(E)\)\cdot c_1(\cL((1-p)E)) \\
&=\(c_{r-1}\(\pi^{\bullet*}\cE^\bullet\)+\sum_{i=1}^{r-1}c_1(E)^ic_{r-i-1}\(\pi^{\bullet*}\cE^\bullet\)\)
\cdot\(c_1(\cL)+(1-p)c_1(E)\) \\
&=c_{r-1}\(\pi^{\bullet*}\cE^\bullet\)\cdot c_1(\cL)+\sum_{i=1}^{r-1}c_1(E)^ic_1(\cL)c_{r-i-1}\(\pi^{\bullet*}\cE^\bullet\)
+(1-p)\sum_{i=1}^rc_1(E)^ic_{r-i}\(\pi^{\bullet*}\cE^\bullet\) \\
&=\pi^{\bullet*}\(c_{r-1}(\cE^\bullet)\cdot c_1\(\cH_1^\dashv/\cH_2\)\)
+\sum_{i=1}^{r-1}c_1(E)^ic_1(\cL)c_{r-i-1}\(\pi^{\bullet*}\cE^\bullet\)
+(1-p)\sum_{i=1}^rc_1(E)^ic_{r-i}\(\pi^{\bullet*}\cE^\bullet\).
\end{align*}
Since $\pi$ and $\pi^\bullet$ are generically finite of degrees $1$ and $p+1$, respectively, it follows that
\begin{align}\label{eq:dl_excess_4}
&\quad (p+1)\int_{\DL^\bullet(\sV)}c_{r-1}(\cE^\bullet)\cdot c_1\(\cH_1^\dashv/\cH_2\)
-\int_{\DL(\sV_\sharp)}c_r(\cE) \\
&=(p-1)\sum_{i=1}^r\int_{\wt\DL(\sV)}c_1(E)^ic_{r-i}\(\pi^{\bullet*}\cE^\bullet\)
-\sum_{i=1}^{r-1}\int_{\wt\DL(\sV)}c_1(E)^ic_1(\cL)c_{r-i-1}\(\pi^{\bullet*}\cE^\bullet\) \notag \\
&=(p-1)\sum_{i=0}^{r-1}\int_E(-\eta)^ic_{r-i-1}\(\pi^{\bullet*}\cE^\bullet\res_E\)
-\sum_{i=0}^{r-2}\int_E(-\eta)^ic_1(\cL\res_E)c_{r-i-2}\(\pi^{\bullet*}\cE^\bullet\res_E\) \notag
\end{align}
where $\eta\coloneqq c_1(\cO_E(1))$. As $\wt\cH/\wt\cH^{\vDash}=\pi^*\(\cH/\cH^{\vDash}\)$, we have $\cL\res_E\simeq\cO_E(-E)=\cO_E(1)$. On the other hand, $\wt\cH_2\res_E$ is the tautological subbundle (of rank $r-1$), which satisfies the short exact sequence
\[
0 \to \wt\cH_2\res_E \to \cO_E^{\oplus r} \to \cO_E(1) \to 0.
\]
Thus, $\cF\coloneqq\pi^{\bullet*}\cE^\bullet\res_E$, which equals $(\sigma^*\wt\cH_2\res_E)\otimes_{\cO_E}(\cL\res_E)$, satisfies the short exact sequence
\[
0 \to \cF \to \cO_E(1)^{\oplus r} \to \cO_E(p+1) \to 0.
\]
Therefore, we have
\begin{align}\label{eq:dl_excess_5}
\eqref{eq:dl_excess_4}&=p\sum_{i=0}^{r-1}\int_E(-\eta)^ic_{r-i-1}(\cF)-\int_Ec_{r-1}(\cF) \\
&=p\int_E c_{r-1}(\cF(-1))-\int_Ec_{r-1}(\cF) \notag \\
&=p\int_E (-p)^{r-1}\eta^{r-1}-\int_E\frac{1-(-p)^r}{p+1}\eta^{r-1} \notag\\
&=\frac{(-p)^{r+1}-1}{p+1}\int_E\eta^{r-1} \notag\\
&=\frac{(-p)^{r+1}-1}{p+1}\cdot|\DL(\sV_\sharp)'(\kappa)| \notag \\
&=\frac{(-p)^{r+1}-1}{p+1}(p+1)(p^3+1)\cdots(p^{2r-1}+1). \notag
\end{align}
By Proposition \ref{pr:dl_excess_pre}, we have
\begin{align}\label{eq:dl_excess_6}
\int_{\DL(\sV_\sharp)}c_r(\cE)=\td_{r,p}.
\end{align}
Thus, (2) follows from \eqref{eq:dl_excess_4}, \eqref{eq:dl_excess_5} and \eqref{eq:dl_excess_6}. The proposition is proved.
\end{proof}

\section{Computation in Hecke algebras}
\label{ss:a}

In this appendix, we compute several explicit formulae on the evaluation of certain Hecke elements. In \S\ref{ss:characters}, we prove some combinatorial formulae on characters of the dual group (of a unitary group). In \S\ref{ss:two_hecke}, we introduce the two unitary Hecke algebras and prove a formula for an intertwining operator between the two Hecke algebras. In \S\ref{ss:enumeration_even} and \S\ref{ss:enumeration_odd}, we evaluate certain Hecke operators under a Satake parameter in the even and odd rank cases, respectively.

\subsection{Characters of the dual group}
\label{ss:characters}

Let $N\geq 1$ be an integer with $r\coloneqq\lfloor\tfrac{N}{2}\rfloor$. We let $\GL_N$ be the group of automorphism of the $\dZ$-module $\dZ^{\oplus N}$, which is a group scheme over $\dZ$. Let $\rT_N\subseteq\GL_N$ be the subgroup of diagonal matrices. The group of homomorphisms from $\rT_N$ to $\dG_m$, denoted by $\dX^*_N$, is a free abelian group generated by $\{\mu_1,\dots,\mu_N\}$ where $\mu_i$ is the projection to the $i$-th factor. For $\mu\in\dX^*_N$, we denote by $[\mu]$ the corresponding element in $\dZ[\dX^*_N]$. For $1\leq i\leq r$, we put
\[
\bmu_i\coloneqq[\mu_i-\mu_{N+1-i}]+[\mu_{N+1-i}-\mu_i]\in\dZ[\dX^*_N].
\]
For $0\leq\delta\leq r$, let $\fs_\delta\in\dZ[\dX^*_N]$ be the elementary symmetric polynomial in $\bmu_1,\dots,\bmu_r$ of degree $\delta$. Finally, we denote by $\dZ[\dX^*_N]^\sym$ the subring of $\dZ[\dX^*_N]$ generated by $\{\fs_1,\dots,\fs_r\}$ over $\dZ$.

Now we consider $\GL^\ext_N\coloneqq\GL_N\rtimes\{1,\sigma\}$ in which the involution $\sigma$ sends $A\in\GL_N$ to
\[
\begin{pmatrix}
& & & &  & 1\\
& &  & -1 & \\
& & \iddots & & \\
& (-1)^{N-2} &  &  & \\
(-1)^{N-1} & & & &
\end{pmatrix}
\tp{A}^{-1}
\begin{pmatrix}
& & & &  & 1\\
& &  & -1 & \\
& & \iddots & & \\
& (-1)^{N-2} &  &  & \\
(-1)^{N-1} & & & &
\end{pmatrix}^{-1}.
\]

For every algebraic representation $\rho$ of $\GL^\ext_N$ (over $\dZ$), we denote by $\chi(\rho)$ the restriction of the character of $\rho$ to $\rT_N\sigma$, regarded as an element in $\dZ[\dX^*_N]$. Let $\rho_{N,\std}$ be the standard representation of $\GL_N$ and $\rho_{N,\std}^\vee$ its dual. We let $\{\varepsilon_1,\dots,\varepsilon_N\}$ be the standard basis of $\rho_{N,\std}$ and $\{\varepsilon^\vee_1,\dots,\varepsilon^\vee_N\}$ the dual basis of $\rho_{N,\std}^\vee$. For a subset $I\subseteq\{1,\dots,N\}$, we put $\langle I\rangle\coloneqq\sum_{i\in I}i$, $I^\vee\coloneqq\{N+1-i\res i\in I\}$, $\varepsilon_I\coloneqq\wedge_{i\in I}\varepsilon_i$ and $\varepsilon^\vee_I\coloneqq\wedge_{i\in I}\varepsilon^\vee_i$ (in the increasing order of the indices). For $0\leq \delta\leq r$, put
\[
\rho_{N;\delta}\coloneqq\(\bigwedge^\delta\rho_{N,\std}\)\otimes\(\bigwedge^\delta\rho_{N,\std}^\vee\),
\]
which extends uniquely to a representation of $\GL^\ext_N$ such that $\sigma$ sends $\varepsilon_I\otimes\varepsilon^\vee_{J^\vee}$ to $(-1)^{\langle I\rangle+\langle J\rangle}\varepsilon_{J}\otimes\varepsilon^\vee_{I^\vee}$.

\begin{remark}
In the next subsection, we will study the unramified unitary group $\rU(\rV_N)$ over nonarchimedean local fields. Then $\GL^\ext_N(\dC)$ is simply the Langlands dual group of $\rU(\rV_N)$, and we have $\dZ[\dX^*_N]^\sym\simeq\dZ[\dX^*(\widehat{\rU(\rV_N)})^\sigma]^{\rW_N}$.
\end{remark}

\begin{lem}\label{le:character1}
We have
\[
\chi(\rho_{N;\delta})=
\begin{dcases}
  \sum_{i=0}^\delta\binom{r-\delta+i}{\lfloor\tfrac{i}{2}\rfloor}\cdot\fs_{\delta-i}, &\text{if $N$ is odd};\\
  \sum_{j=0}^{\lfloor\tfrac{\delta}{2}\rfloor}\binom{r-\delta+2j}{j}\cdot\fs_{\delta-2j}, &\text{if $N$ is even}.
\end{dcases}
\]
In particular, $\chi(\rho_{N;\delta})$ belongs to $\dZ[\dX^*_N]^\sym$.
\end{lem}

\begin{proof}
Note that for every $t\in\rT_N$, $t\sigma$ sends $\varepsilon_I\otimes\varepsilon^\vee_{J^\vee}$ to
\[
(-1)^{\langle I\rangle+\langle J\rangle}\prod_{i\in I^\vee}\mu_i(t)^{-1}\prod_{j\in J}\mu_j(t)\cdot\varepsilon_{J}\otimes\varepsilon^\vee_{I^\vee}.
\]
In particular, such term contributes to $\chi(\rho_{N,\delta})(t\sigma)$ exactly when $I=J$. It follows that
\begin{align*}
\chi(\rho_{N,\delta})(t\sigma)
&=\sum_{I\subseteq\{1,\dots,N\},|I|=\delta}\prod_{i\in I^\vee}\mu_i(t)^{-1}\prod_{i\in I}\mu_i(t) \\
&=\sum_{I\subseteq\{1,\dots,N\},|I|=\delta}\prod_{i\in I}\mu_i(t)\mu_{N+1-i}(t)^{-1}.
\end{align*}
To evaluate the above sum, we consider $i\coloneqq|I\cap I^\vee|$, which has to be even when $N$ is even. It is easy to see that for fixed $0\leq i\leq \delta$ (that is even if $N$ is even), the contribution from those subsets $I$ to the above sum is
\[
\binom{r-\delta+i}{\lfloor\tfrac{i}{2}\rfloor}\cdot\fs_{\delta-i}(t).
\]
Thus, the lemma follows.
\end{proof}

\begin{lem}\label{le:character2}
Suppose that $N=2r$ is even.
\begin{enumerate}
  \item We have
     \[
     \prod_{i=1}^r\(\lambda+\lambda^{-1}+\bmu_i\)=
     \chi(\rho_{N;r})+\sum_{\delta=1}^r\chi(\rho_{N;r-\delta})(\lambda^\delta+\lambda^{-\delta})
     \]
     in $\dZ[\dX^*_N]^\sym\otimes\dZ[\lambda,\lambda^{-1}]$.

  \item We have
     \[
     \sum_{j=1}^r\prod_{\substack{i=1\\i\neq j}}^r\(\lambda+\lambda^{-1}+\bmu_i\)=\sum_{\delta=1}^r\delta\cdot
     \chi(\rho_{N;r-\delta})\frac{\lambda^\delta-\lambda^{-\delta}}{\lambda-\lambda^{-1}}
     \]
     in $\dZ[\dX^*_N]^\sym\otimes\dZ[\lambda,\lambda^{-1}]$.
\end{enumerate}
\end{lem}

\begin{proof}
Part (1) is follows from Lemma \ref{le:character1} by comparing coefficients of powers of $\lambda$. Part (2) follows from (1) by taking derivative with respect to $\lambda$ and dividing both sides of the resulted equality by $1-\lambda^{-2}$.
\end{proof}

\begin{lem}\label{le:character3}
Suppose that $N=2r+1$ is odd. We have
\[
\prod_{i=1}^r\(\lambda+\lambda^{-1}+\bmu_i\)=
\sum_{\delta=0}^r\chi(\rho_{N;r-\delta})\frac{\lambda^{\delta+1}+\lambda^{-\delta}}{\lambda+1}
\]
in $\dZ[\dX^*_N]^\sym\otimes\dZ[\lambda,\lambda^{-1}]$.
\end{lem}

\begin{proof}
By Lemma \ref{le:character1}, the right-hand side of the desired identity equals
\[
\sum_{\delta=0}^r\frac{\lambda^{\delta+1}+\lambda^{-\delta}}{\lambda+1}
\sum_{i=0}^{r-\delta}\binom{\delta+i}{\lfloor\tfrac{i}{2}\rfloor}\cdot\fs_{r-\delta-i},
\]
which coincides with
\[
\sum_{i=0}^r\(\sum_{\delta=0}^{r-i}\frac{\lambda^{\delta+1}+\lambda^{-\delta}}{\lambda+1}\binom{r-i}{\lfloor\tfrac{r-i-\delta}{2}\rfloor}\)\fs_i
\]
by substituting $i$ by $r-\delta-i$. Thus, it remains to show that
\[
\sum_{\delta=0}^k\frac{\lambda^{\delta+1}+\lambda^{-\delta}}{\lambda+1}\binom{k}{\lfloor\tfrac{k-\delta}{2}\rfloor}
=(\lambda+\lambda^{-1})^k
\]
for $0\leq k\leq r$. However, we have
\begin{align*}
&\quad\sum_{\delta=0}^k\frac{\lambda^{\delta+1}+\lambda^{-\delta}}{\lambda+1}\binom{k}{\lfloor\tfrac{k-\delta}{2}\rfloor} \\
&=\binom{k}{0}\(\frac{\lambda^{k+1}+\lambda^{-k}}{\lambda+1}+\frac{\lambda^{k}+\lambda^{-(k-1)}}{\lambda+1}\)
+\binom{k}{1}\(\frac{\lambda^{k-1}+\lambda^{-(k-2)}}{\lambda+1}+\frac{\lambda^{k-2}+\lambda^{-(k-3)}}{\lambda+1}\)
+\cdots \\
&=\binom{k}{0}(\lambda^k+\lambda^{-k})+\binom{k}{1}(\lambda^{k-1}+\lambda^{-(k-1)})+\cdots \\
&=(\lambda+\lambda^{-1})^k.
\end{align*}
The lemma follows.
\end{proof}

\subsection{Two Hecke algebras}
\label{ss:two_hecke}

From now to the end of this section, we fix an unramified quadratic extension $F/F^+$ of nonarchimedean local fields. Let $q$ be the residue cardinality of $F^+$ and $\fp$ the maximal ideal of $O_F$.

Let $N\geq 1$ be an integer with $r\coloneqq\lfloor\tfrac{N}{2}\rfloor$. Consider a hermitian space $\rV_N$ over $F$ (with respect to $F/F^+$) of rank $N$ together with a basis $\{e_{-r},\dots,e_r\}$ (with $e_0$ omitted if $N$ is even) such that $(e_{-i},e_j)_{\rV_N}=\delta_{ij}$ for $0\leq i,j\leq r$. Via this basis, we identify $\rU(\rV_N)$ as a closed subgroup of $\Res_{F/F^+}\GL_N$. We study two lattices
\begin{align}\label{eq:lattice}
\Lambda^\circ_N=O_Fe_{-r}\oplus\cdots\oplus O_F e_r,\quad
\Lambda^\bullet_N=\fp^{-1}e_{-r}\oplus\cdots\oplus \fp^{-1}e_{-1}\oplus O_F e_0 \oplus\cdots\oplus O_F e_r
\end{align}
of $\rV_N$. We have $(\Lambda^\circ_N)^\vee=\Lambda^\circ_N$, $\fp\Lambda^\bullet_N\subseteq(\Lambda^\bullet_N)^\vee$, and that the $O_F$-module $(\Lambda^\bullet_N)^\vee/\fp\Lambda^\bullet_N$ has length $N-2r$. Let $\rK^\circ_N$ and $\rK^\bullet_N$ be the stabilizers of $\Lambda^\circ_N$ and $\Lambda^\bullet_N$, respectively, which are subgroups of $\rU(\rV_N)(F^+)$. It is clear that $\rK^\circ_N$ is hyperspecial maximal; $\rK^\bullet_N$ is special maximal and is hyperspecial if and only if $N$ is even. We have two \emph{commutative} Hecke algebras
\[
\dT_N^\circ\coloneqq\dZ[\rK^\circ_N\backslash\rU(\rV_N)(F^+)/\rK^\circ_N],\quad
\dT_N^\bullet\coloneqq\dZ[\rK^\bullet_N\backslash\rU(\rV_N)(F^+)/\rK^\bullet_N].
\]
Recall that by our convention in \S\ref{ss:notation}, the units in $\dT_N^\circ$ and $\dT_N^\bullet$ are $\CF_{\rK^\circ_N}$ and $\CF_{\rK^\bullet_N}$, respectively. Let $\rA_N(F^+)$ (resp.\ $\rA_N(O_{F^+})$) be the subgroup of $\rU(\rV_N)(F^+)$ that acts on $e_i$ by a scalar in $F^+$ (resp.\ $O_{F^+}$) for every $-r\leq i\leq r$.

\begin{notation}\label{no:hecke}
For each element $\bbt=(t_1,\dots,t_N)\in\dZ^N$ satisfying $t_i+t_{N+1-i}=0$ and $a\in F^\times$, we have an element
$a^{\bbt}\in\rA_N(F^+)$ such that $a^{\bbt}\cdot e_{-i}=a^{t_{r+1-i}}e_{-i}$ for $0\leq i\leq r$. For $0\leq\delta\leq r$, put $\bbt_\delta\coloneqq(1^\delta,0^{N-2\delta},(-1)^\delta)$. We let $\tT^\circ_{N;\bbt}$ (resp.\ $\tT^\bullet_{N;\bbt}$) be the element in $\dT_N^\circ$ (resp.\ $\dT_N^\bullet$) corresponding to the double coset $\rK^\circ_N\varpi^{\bbt} \rK^\circ_N$ (resp.\ $\rK^\bullet_N\varpi^{\bbt} \rK^\bullet_N$) for some uniformizer $\varpi$ of $F$; and simply write $\tT^\circ_{N;\delta}$ (resp.\ $\tT^\bullet_{N;\delta}$) for $\tT^\circ_{N;\bbt_\delta}$ (resp.\ $\tT^\bullet_{N;\bbt_\delta}$).
\end{notation}

\begin{remark}
The elements $\tT^\circ_{N;\bbt}\in\dT^\circ_N$ and $\tT^\bullet_{N;\bbt}\in\dT^\bullet_N$ do not depend on the choice of the basis $\{e_{-r},\dots,e_r\}$ satisfying \eqref{eq:lattice}.
\end{remark}

\begin{definition}\label{de:lattice}
We denote
\begin{itemize}[label={\ding{109}}]
  \item $\Lat^\circ_N$ the set of all self-dual lattices in $\rV_N$;

  \item $\Lat^\bullet_N$ the set of all lattices $\rL$ in $\rV_N$ satisfying $\fp\rL\subseteq\rL^\vee$ and that $\rL^\vee/\fp\rL$ has length $N-2\lfloor\frac{N}{2}\rfloor$;

  \item $\tT^{\bullet\circ}_N\in\dZ[\rK^\bullet_N\backslash\rU(\rV_N)(F^+)/\rK^\circ_N]$ the characteristic function of $\rK^\bullet_N\rK^\circ_N$; and

  \item $\tT^{\circ\bullet}_N\in\dZ[\rK^\circ_N\backslash\rU(\rV_N)(F^+)/\rK^\bullet_N]$ the characteristic function of $\rK^\circ_N\rK^\bullet_N$.
\end{itemize}
Moreover, we define the \emph{intertwining Hecke operator}
\[
\tI^\circ_N\coloneqq\tT^{\circ\bullet}_N\circ\tT^{\bullet\circ}_N\in\dT^\circ_N
\]
where the composition is taken as composition of cosets.
\end{definition}

Note that we have canonical injective homomorphisms
\[
\dT^\circ_N\to\End_\dZ(\dZ[\Lat^\circ_N]),\quad
\dT^\bullet_N\to\End_\dZ(\dZ[\Lat^\bullet_N])
\]
sending $\tT^?_{N;\bbt}$ to the endomorphism that takes $f\in\dZ[\Lat^?_N]$ to the function $\tT^?_{N;\bbt}f$ satisfying $(\tT^?_{N;\bbt}f)(\rL)=\sum f(\rL')$ where the sum is taken over all $\rL'\in\Lat^?_N$ such that $\rL'$ and $\rL$ have relative position $\varpi^{\bbt}$ for $?=\circ,\bullet$.

\begin{lem}\label{le:hecke_circ}
We have the identity
\begin{align*}
\tI^\circ_N=
\begin{dcases}
\tT^\circ_{N;r}+(q+1)\tT^\circ_{N;r-1}+(q+1)(q^3+1)\tT^\circ_{N;r-2}+\cdots+\prod_{i=1}^r(q^{2i-1}+1)\tT^\circ_{N;0},&\text{if $N=2r$;}\\
\tT^\circ_{N;r}+(q^3+1)\tT^\circ_{N;r-1}+(q^3+1)(q^5+1)\tT^\circ_{N;r-2}+\cdots+\prod_{i=1}^r(q^{2i+1}+1)\tT^\circ_{N;0},&\text{if $N=2r+1$}
\end{dcases}
\end{align*}
in $\dT^\circ_N$.
\end{lem}

\begin{proof}
For a pair $(\rL^\circ_1,\rL^\circ_2)\in(\Lat^\circ_N)^2$, we denote by $\Disc(\rL^\circ_1,\rL^\circ_2)$ the sum of the lengths of $\rL^\circ_1/(\rL^\circ_1\cap\rL^\circ_2)$ and $\rL^\circ_2/(\rL^\circ_1\cap\rL^\circ_2)$.

To compute $\tI^\circ_N$, it suffices to compute its induced endomorphism on $\dZ[\Lat^\circ_N]$. Now we take an element $f\in\dZ[\Lat^\circ_N]$. Then
\begin{align*}
(\tT^{\circ\bullet}_N(\tT^{\bullet\circ}_N f))(\rL^\circ_1)
=\sum_{\substack{\rL^\bullet\in\Lat^\bullet_N\\\rL^\circ_1\subseteq \rL^\bullet\subseteq \fp^{-1}\rL^\circ_1}}(\tT^{\bullet\circ}_N f)(\rL^\bullet)
=\sum_{\substack{\rL^\bullet\in\Lat^\bullet_N\\\rL^\circ_1\subseteq \rL^\bullet\subseteq \fp^{-1}\rL^\circ_1}}
\sum_{\substack{\rL^\circ_2\in\Lat^\circ_N\\\rL^\circ_2\subseteq \rL^\bullet\subseteq\fp^{-1}\rL^\circ_2}}f(\rL^\circ_2)
\end{align*}
for every $\rL^\circ_1\in\Lat^\circ_N$. Note that for pairs $(\rL^\circ_1,\rL^\circ_2)\in(\Lat^\circ_N)^2$ appearing in the formula above, we have $\fp \rL^\circ_2\subseteq \rL^\circ_1\subset\fp^{-1}\rL^\circ_2$ and $\Disc(\rL^\circ_1,\rL^\circ_2)\in\{0,2,\dots,2r\}$.

Now for a pair $(\rL^\circ_1,\rL^\circ_2)\in(\Lat^\circ_N)^2$ satisfying $\fp\rL^\circ_2\subseteq \rL^\circ_1\subset\fp^{-1}\rL^\circ_2$, we consider the set
\[
\Lat^\bullet_N(\rL^\circ_1,\rL^\circ_2)\coloneqq\{\rL^\bullet\in\Lat^\bullet_N\res \rL^\circ_1\subseteq \rL^\bullet\subseteq\fp^{-1}\rL^\circ_1,
\rL^\circ_2\subseteq \rL^\bullet\subseteq\fp^{-1}\rL^\circ_2\}.
\]
It is easy to see that the cardinality of $\Lat^\bullet_N(\rL^\circ_1,\rL^\circ_2)$ depends only on $\Disc(\rL^\circ_1,\rL^\circ_2)$. For $0\leq \delta\leq r$, we denote by $c_{N,\delta}$ the cardinality of $\Lat^\bullet_N(\rL^\circ_1,\rL^\circ_2)$ with $\Disc(\rL^\circ_1,\rL^\circ_2)=2\delta$. Then the lemma is equivalent to showing that $c_{N,r}=1$ and
\[
c_{N,\delta}=
\begin{dcases}
\prod_{i=1}^{r-\delta}(q^{2i-1}+1),\quad 0\leq \delta <r,&\text{when $N=2r$;}\\
\prod_{i=1}^{r-\delta}(q^{2i+1}+1),\quad 0\leq \delta <r,&\text{when $N=2r+1$.}
\end{dcases}
\]
Without loss of generality, we may assume $\rL^\circ_1=\Lambda^\circ_N$ and
\[
\rL^\circ_2=\fp^{-1}e_{-r}\oplus\cdots\oplus\fp^{-1}e_{-r+\delta-1} \oplus O_F e_{-r+\delta} \oplus\cdots\oplus O_F e_{r-\delta}
\oplus\fp O_F e_{r-\delta+1}\oplus\cdots\oplus\fp O_F e_r.
\]
When $\delta=r$, $\Lambda^\bullet_N$ is the only element in $\Lat^\bullet_N(\rL^\circ_1,\rL^\circ_2)$. Thus, we have $c_{N,r}=1$. For $0\leq\delta<r$, we have $c_{N,\delta}=c_{N-2\delta,0}$. Thus, it suffices to show
\[
c_{N,0}=
\begin{dcases}
\prod_{i=1}^r(q^{2i-1}+1)=(q+1)\cdots(q^{2r-1}+1),&\text{when $N=2r$;}\\
\prod_{i=1}^r(q^{2i+1}+1)=(q^3+1)\cdots(q^{2r+1}+1),&\text{when $N=2r+1$.}
\end{dcases}
\]
However, $c_{N,0}$ is nothing but the number of maximal isotropic subspaces of the hermitian space $\Lambda^\circ_N\otimes_{O_F}O_F/\fp$ over $O_F/\fp$ of dimension $N$, which is given by the above formula. Thus, the lemma is proved.
\end{proof}

Now we recall Satake transforms. Denote by $\rW_N$ the Weyl group of $\rA_N(F^+)$ in $\rU(\rV_N)(F^+)$, which preserves $\rA_N(O_{F^+})$; and we have the two Satake transforms
\begin{align*}
\Sat^\circ_N&\colon\dT^\circ_N\to\dZ[q^{-1}][\rA_N(F^+)/\rA_N(O_{F^+})]^{\rW_N},\\
\Sat^\bullet_N&\colon\dT^\bullet_N\to\dZ[q^{-1}][\rA_N(F^+)/\rA_N(O_{F^+})]^{\rW_N}.
\end{align*}
In addition, we have an isomorphism
\[
\dZ[q^{-1}][\rA_N(F^+)/\rA_N(O_{F^+})]^{\rW_N}\simeq\dZ[q^{-1}][\dX^*_N]^\sym
\]
of $\dZ[q^{-1}]$-rings under which $\fs_\delta$ corresponds to the sum of elements in the $\rW_N$-orbit of $\varpi^{\bbt_\delta}\rA_N(O_{F^+})$ for every $0\leq \delta\leq r$. In what follows, we will regard $\dZ[q^{-1}][\dX^*_N]^\sym$ as the target of both Satake transforms $\Sat^\circ_N$ and $\Sat^\bullet_N$.

\begin{notation}\label{no:parameter}
Let $\dZ[q^{-1}][\dX^*_N]'$ be the $\dZ[q^{-1}]$-subring of $\dZ[q^{-1}][\dX^*_N]$ generated by the subset $\{\bmu_1,\dots,\bmu_r\}$. For every $\dZ[q^{-1}]$-ring $L$ and every tuple $\balpha=(\alpha_1,\dots,\alpha_N)\in L^N$ satisfying $\alpha_i\alpha_{N+1-i}=1$, we have a homomorphism $\phi'_{\balpha}\colon\dZ[q^{-1}][\dX^*_N]'\to L$ sending $\bmu_i$ to $\alpha_i+\alpha_i^{-1}$ for $1\leq i\leq r$, similar to Construction \ref{cs:satake_hecke_pre}, and denote by
\begin{align*}
\phi^\circ_{\balpha}&\colon\dT^\circ_N\xrightarrow{\Sat^\circ_N}\dZ[q^{-1}][\dX^*_N]^\sym
\subseteq\dZ[q^{-1}][\dX^*_N]'\xrightarrow{\phi'_{\balpha}}L,\\
\phi^\bullet_{\balpha}&\colon\dT^\bullet_N\xrightarrow{\Sat^\bullet_N}\dZ[q^{-1}][\dX^*_N]^\sym
\subseteq\dZ[q^{-1}][\dX^*_N]'\xrightarrow{\phi'_{\balpha}}L,
\end{align*}
the composite homomorphisms.
\end{notation}

The following three lemmas will be used in later computation.

\begin{lem}\label{le:enumeration1}
We have the identity
\[
q^{\delta(N-\delta)}\chi(\rho_{N,\delta})=\sum_{i=0}^\delta\qbinom{N-2i}{\delta-i}_{-q}\Sat^\circ_N(\tT^\circ_{N;i})
\]
in $\dZ[q^{-1}][\dX^*_N]^\sym$ for $0\leq\delta\leq r$.
\end{lem}

\begin{proof}
This is \cite{XZ17}*{Lemma~9.2.4}.
\end{proof}

\begin{lem}\label{le:enumeration2}
For every integer $k\geq 1$, we have
\[
\sum_{\delta=-k}^k q^{\delta^2}\qbinom{2k}{k-\delta}_{-q}=(q+1)(q^3+1)\cdots(q^{2k-1}+1).
\]
\end{lem}

\begin{proof}
For every integer $k\geq 1$, we have the Gauss polynomial identity
\[
\sum_{\delta=0}^{2k}(-1)^\delta\qbinom{2k}{\delta}_\lambda=(1-\lambda)(1-\lambda^3)\cdots(1-\lambda^{2k-1})
\]
in $\dZ[\lambda]$.\footnote{A proof can be found at \url{http://mathworld.wolfram.com/GausssPolynomialIdentity.html}.} Now we specialize the identity to $\lambda=-q^{-1}$. Then we get
\[
\sum_{\delta=0}^{2k}(-1)^\delta(-q)^{-(2k-1)-(2k-3)-\cdots-(2k-2\delta+1)}\qbinom{2k}{\delta}_{-q}
=q^{-k^2}(q+1)(q^3+1)\cdots(q^{2k-1}+1).
\]
The lemma then follows by changing $\delta$ to $k-\delta$.
\end{proof}

\begin{lem}\label{le:enumeration3}
For every integer $k\geq 1$, we have
\[
\sum_{\delta=-k-1}^k(-1)^\delta\delta q^{\delta^2+\delta}\qbinom{2k+1}{k-\delta}_{-q}
-\sum_{\delta=-k}^k(-1)^\delta\delta q^{\delta^2+\delta}\qbinom{2k}{k-\delta}_{-q}=(-q)^k(q+1)(q^3+1)\cdots(q^{2k-1}+1).
\]
\end{lem}

\begin{proof}
In fact, we have
\begin{align*}
&\quad\sum_{\delta=-k-1}^k(-1)^\delta\delta q^{\delta^2+\delta}\qbinom{2k+1}{k-\delta}_{-q}
-\sum_{\delta=-k}^k(-1)^\delta\delta q^{\delta^2+\delta}\qbinom{2k}{k-\delta}_{-q} \\
&=\sum_{\delta=-k-1}^k(-1)^\delta\delta q^{\delta^2+\delta}(-q)^{k+\delta+1}\qbinom{2k}{k-\delta-1}_{-q} \\
&=(-1)^{k+1}q^k\sum_{\delta=-k}^k(\delta-1)q^{\delta^2}\qbinom{2k}{k-\delta}_{-q}
\end{align*}
which, by Lemma \ref{le:enumeration2}, equals
\[
(-q)^k(q+1)(q^3+1)\cdots(q^{2k-1}+1)+(-1)^{k+1}q^k\sum_{\delta=-k}^k\delta q^{\delta^2}\qbinom{2k}{k-\delta}_{-q}.
\]
The lemma follows since
\[
\sum_{\delta=-k}^k\delta q^{\delta^2}\qbinom{2k}{k-\delta}_{-q}=0.
\]
\end{proof}

\subsection{Enumeration of Hecke operators in the even rank case}
\label{ss:enumeration_even}

In this subsection, we assume that $N=2r$ is even.

\begin{lem}\label{le:enumeration_even_1}
We have the identity
\[
q^{r^2}\prod_{i=1}^r\(\bmu_i+2\)=\Sat^\circ_N(\tT^\circ_{N;r})
+\sum_{\delta=1}^r(q+1)(q^3+1)\cdots(q^{2\delta-1}+1)\cdot\Sat^\circ_N(\tT^\circ_{N;r-\delta})
\]
in $\dZ[q^{-1}][\dX^*_N]^\sym$.
\end{lem}

\begin{proof}
By Lemma \ref{le:character2}(1) and Lemma \ref{le:enumeration1}, we have
\begin{align*}
q^{r^2}\prod_{i=1}^r\(\bmu_i+2\)
&=q^{r^2}\chi(\rho_{N;r})+q^{r^2}\sum_{\delta=1}^r2\chi(\rho_{N;r-\delta}) \\
&=\sum_{i=0}^r\qbinom{2r-2i}{r-i}_{-q}\Sat^\circ_N(\tT^\circ_{N;i})+
\sum_{\delta=1}^r2q^{\delta^2}\sum_{i=0}^{r-\delta}\qbinom{2r-2i}{r-\delta-i}_{-q}\Sat^\circ_N(\tT^\circ_{N;i}) \\
&=\sum_{i=0}^r\(\sum_{\delta=-(r-i)}^{r-i}q^{\delta^2}\qbinom{2r-2i}{r-\delta-i}_{-q}\)\Sat^\circ_N(\tT^\circ_{N;i}),
\end{align*}
which equals
\[
\Sat^\circ_N(\tT^\circ_{N;r})
+\sum_{\delta=1}^r(q+1)(q^3+1)\cdots(q^{2\delta-1}+1)\cdot\Sat^\circ_N(\tT^\circ_{N;r-\delta})
\]
by Lemma \ref{le:enumeration2}. The lemma is proved.
\end{proof}

\begin{lem}\label{le:enumeration_even_2}
We have the identity
\[
q^{r^2}\prod_{i=1}^r\(\bmu_i-q-q^{-1}\)=\Sat^\circ_N(\tT^\circ_{N;r})
+\sum_{\delta=1}^r(-q)^\delta(q+1)(q^3+1)\cdots(q^{2\delta-1}+1)\cdot\Sat^\circ_N(\tT^\circ_{N;r-\delta})
\]
in $\dZ[q^{-1}][\dX^*_N]^\sym$.
\end{lem}

\begin{proof}
By Lemma \ref{le:character2}(1) and Lemma \ref{le:enumeration1}, we have
\begin{align*}
&\quad q^{r^2}\prod_{i=1}^r\(\bmu_i-q-q^{-1}\) \\
&=q^{r^2}\chi(\rho_{N;r})+q^{r^2}\sum_{\delta=1}^r((-q)^\delta+(-q)^{-\delta})\chi(\rho_{N;r-\delta}) \\
&=\sum_{i=0}^r\qbinom{2r-2i}{r-i}_{-q}\Sat^\circ_N(\tT^\circ_{N;i})+
\sum_{\delta=1}^r\sum_{i=0}^{r-\delta}q^{\delta^2}((-q)^\delta+(-q)^{-\delta})\qbinom{2r-2i}{r-\delta-i}_{-q}\Sat^\circ_N(\tT^\circ_{N;i}) \\
&=\sum_{i=0}^r\(\qbinom{2r-2i}{r-i}_{-q}+\sum_{\delta=1}^{r-i}
(-1)^\delta\(q^{\delta^2+\delta}+q^{\delta^2-\delta}\)\qbinom{2r-2i}{r-\delta-i}_{-q}\)\Sat^\circ_N(\tT^\circ_{N;i}) \\
&=\sum_{i=0}^r\(\sum_{\delta=-(r-i)}^{r-i}
(-1)^\delta q^{\delta^2+\delta}\qbinom{2r-2i}{r-\delta-i}_{-q}\)\Sat^\circ_N(\tT^\circ_{N;i}).
\end{align*}
Thus, the lemma follows from Lemma \ref{le:enumeration_even_3} below by comparing coefficients.
\end{proof}

\begin{lem}\label{le:enumeration_even_3}
For every integer $k\geq 1$, we have
\[
\sum_{\delta=-k}^{k}(-1)^\delta q^{\delta^2+\delta}\qbinom{2k}{k-\delta}_{-q}
=(-q)^k(q+1)(q^3+1)\cdots(q^{2k-1}+1).
\]
\end{lem}

\begin{proof}
By Lemma \ref{le:enumeration2}, the lemma is equivalent to the identity
\[
(-q)^k\sum_{\delta=-k}^k q^{\delta^2}\qbinom{2k}{k-\delta}_{-q}=\sum_{\delta=-k}^{k}(-1)^\delta q^{\delta^2+\delta}\qbinom{2k}{k-\delta}_{-q}.
\]
However, we have
\begin{align*}
&\quad(-q)^k\sum_{\delta=-k}^k q^{\delta^2}\qbinom{2k}{k-\delta}_{-q}-\sum_{\delta=-k}^{k}(-1)^\delta q^{\delta^2+\delta}\qbinom{2k}{k-\delta}_{-q} \\
&=\sum_{\delta=-k}^k(-1)^\delta q^{\delta^2+\delta}\((-q)^{k-\delta}-1\)\qbinom{2k}{k-\delta}_{-q} \\
&=\sum_{\delta=-k}^k(-1)^\delta q^{\delta^2+\delta}\((-q)^{2k}-1\)\qbinom{2k-1}{k-\delta-1}_{-q} \\
&=\((-q)^{2k}-1\)\sum_{\delta=-k}^k(-1)^\delta q^{\delta^2+\delta}\qbinom{2k-1}{k-\delta-1}_{-q}.
\end{align*}
Note that in the last summation, the term of $\delta$ and the term of $-\delta-1$ cancel with each other for $-k\leq\delta\leq k-1$; and the term with $\delta=k$ vanishes. Thus, the above summation is zero; and the lemma follows.
\end{proof}

\begin{lem}\label{le:enumeration_even_4}
We have the identity
\begin{align*}
&\quad\(q^{r^2+1}-q^{r^2-1}\)\sum_{j=1}^r\prod_{\substack{i=1\\i\neq j}}^r\(\bmu_i-q-q^{-1}\) \\
&=\sum_{\delta=1}^r\((-q)^\delta(q+1)(q^3+1)\cdots(q^{2\delta-1}+1)-\sum_{i=0}^\delta(-1)^i(2i+1)q^{i^2+i}\qbinom{2\delta+1}{\delta-i}_{-q}\)
\Sat^\circ_N(\tT^\circ_{N;r-\delta})
\end{align*}
in $\dZ[q^{-1}][\dX^*_N]^\sym$.
\end{lem}

\begin{proof}
By Lemma \ref{le:character2}(2) and Lemma \ref{le:enumeration1}, we have
\begin{align*}
&\quad\(q^{r^2+1}-q^{r^2-1}\)\sum_{j=1}^r\prod_{i\neq j}\(\bmu_i-q-q^{-1}\) \\
&=q^{r^2}\sum_{\delta=1}^r(-1)^{\delta-1}\delta(q^\delta-q^{-\delta})\cdot\chi(\rho_{N;r-\delta}) \\
&=\sum_{\delta=1}^r(-1)^{\delta-1}q^{\delta^2}(\delta q^\delta-\delta q^{-\delta})\sum_{i=0}^{r-\delta}\qbinom{2r-2i}{r-\delta-i}_{-q}\Sat^\circ_N(\tT^\circ_{N;i}) \\
&=\sum_{i=0}^{r-1}\(\sum_{\delta=1}^{r-i}(-1)^{\delta-1}
q^{\delta^2}(\delta q^\delta-\delta q^{-\delta})\qbinom{2r-2i}{r-\delta-i}_{-q}\)\Sat^\circ_N(\tT^\circ_{N;i}).
\end{align*}
Thus the lemma is equivalent to the identity
\begin{align*}
&\quad\sum_{\delta=0}^k(-1)^\delta(2\delta+1)q^{\delta^2+\delta}\qbinom{2k+1}{k-\delta}_{-q}-\sum_{\delta=1}^k(-1)^\delta
q^{\delta^2}(\delta q^\delta-\delta q^{-\delta})\qbinom{2k}{k-\delta}_{-q} \\
&=(-q)^k(q+1)(q^3+1)\cdots(q^{2k-1}+1)
\end{align*}
for every integer $k\geq 1$. In fact, we have
\begin{align*}
&\quad\sum_{\delta=0}^k(-1)^\delta(2\delta+1)q^{\delta^2+\delta}\qbinom{2k+1}{k-\delta}_{-q}
-\sum_{\delta=1}^k(-1)^\delta q^{\delta^2}(\delta q^\delta-\delta q^{-\delta})\qbinom{2k}{k-\delta}_{-q} \\
&=\sum_{\delta=-k-1}^k(-1)^\delta\delta q^{\delta^2+\delta}\qbinom{2k+1}{k-\delta}_{-q}
-\sum_{\delta=-k}^k(-1)^\delta q^{\delta^2}\delta q^\delta\qbinom{2k}{k-\delta}_{-q} \\
&=(-q)^k(q+1)(q^3+1)\cdots(q^{2k-1}+1)
\end{align*}
by Lemma \ref{le:enumeration3}. The lemma follows.
\end{proof}

\begin{proposition}\label{pr:enumeration_even}
Let $L$ be a $\dZ[q^{-1}]$-ring. Consider an $N$-tuple $\balpha=(\alpha_1,\dots,\alpha_N)\in L^N$ satisfying $\alpha_i\alpha_{N+1-i}=1$, which determines a homomorphism $\phi^\circ_{\balpha}\colon\dT^\circ_N\to L$ as in Notation \ref{no:parameter}.
\begin{enumerate}
  \item We have
     \[
     \phi^\circ_{\balpha}(\tI^\circ_N)=q^{r^2}\prod_{i=1}^r\(\alpha_i+\frac{1}{\alpha_i}+2\).
     \]

  \item We have
     \[
     \phi^\circ_{\balpha}\((q+1)\tR^\circ_N-\tI^\circ_N\)=-q^{r^2}\prod_{i=1}^r\(\alpha_i+\frac{1}{\alpha_i}-q-\frac{1}{q}\)
     \]
     where
     \[
     \tR^\circ_N\coloneqq\sum_{\delta=0}^{r-1}
     \frac{1-(-q)^{r-\delta}}{q+1} (q+1)(q+3)\cdots(q^{2(r-\delta)-1}+1)\cdot\tT^\circ_{N;\delta}.
     \]

  \item We have
     \[
     \phi^\circ_{\balpha}\(\tR^\circ_N+(q+1)\tT^\circ_N\)=
     -\(q^{r^2+1}-q^{r^2-1}\)\sum_{j=1}^r\prod_{\substack{i=1\\i\neq j}}^r\(\alpha_i+\frac{1}{\alpha_i}-q-\frac{1}{q}\)
     \]
     where
     \[
     \tT^\circ_N\coloneqq\sum_{\delta=0}^{r-1}\rd^\bullet_{r-\delta,q}\cdot\tT^\circ_{N;\delta}
     \]
     in which the numbers $\td^\bullet_{r-\delta,q}$ are introduced in Notation \ref{no:numerical}.
\end{enumerate}
\end{proposition}

\begin{proof}
Part (1) follows from Lemma \ref{le:hecke_circ} and Lemma \ref{le:enumeration_even_1}. Part (2) follows from Lemma \ref{le:hecke_circ} and Lemma \ref{le:enumeration_even_2}. Part (3) follows from Lemma \ref{le:enumeration_even_4}.
\end{proof}

\begin{lem}\label{le:enumeration_even_0}
We have
\[
\tT^{\bullet\circ}_N\circ\tR^\circ_N=\tR^\bullet_N\circ\tT^{\bullet\circ}_N,\quad
\tT^{\bullet\circ}_N\circ\tT^\circ_N=\tT^\bullet_N\circ\tT^{\bullet\circ}_N
\]
in $\dZ[\rK^\bullet_N\backslash\rU(\rV_N)(F^+)/\rK^\circ_N]$, where $\tR^\circ_N$ and $\tT^\circ_N$ are defined in Proposition \ref{pr:enumeration_even} (2) and (3), respectively, and
\[
\begin{dcases}
\tR^\bullet_N\coloneqq\sum_{\delta=0}^{r-1}
\frac{1-(-q)^{r-\delta}}{q+1} (q+1)(q+3)\cdots(q^{2(r-\delta)-1}+1)\cdot\tT^\bullet_{N;\delta},\\
\tT^\bullet_N\coloneqq\sum_{\delta=0}^{r-1}\rd^\bullet_{r-\delta,q}\cdot\tT^\bullet_{N;\delta}.
\end{dcases}
\]
\end{lem}

\begin{proof}
In fact, by the same lattice counting argument as for Lemma \ref{le:hecke_circ}, we have
\[
\tT^{\bullet\circ}_N\circ\tT^\circ_{N;\delta}=\tT^\bullet_{N;\delta}\circ\tT^{\bullet\circ}_N
\]
for every $0\leq\delta\leq r$. Then the lemma follows immediately.
\end{proof}

\subsection{Enumeration of Hecke operators in the odd rank case}
\label{ss:enumeration_odd}

In this subsection, we assume that $N=2r+1$ is odd.

\begin{lem}\label{le:enumeration_odd_1}
We have the identity
\[
q^{r^2+r}\prod_{i=1}^r\(\bmu_i+q+q^{-1}\)=\Sat^\circ_N(\tT^\circ_{N;r})
+\sum_{\delta=1}^r(q^3+1)(q^5+1)\cdots(q^{2\delta+1}+1)\cdot\Sat^\circ_N(\tT^\circ_{N;r-\delta})
\]
in $\dZ[q^{-1}][\dX^*_N]^\sym$.
\end{lem}

\begin{proof}
By Lemma \ref{le:character3} and Lemma \ref{le:enumeration1}, we have
\begin{align*}
q^{r^2+r}\prod_{i=1}^r\(\bmu_i+q+q^{-1}\)
&=q^{r^2+r}\sum_{\delta=0}^r\frac{q^{\delta+1}+q^{-\delta}}{q+1}\cdot\chi(\rho_{N;r-\delta}) \\
&=q^{r^2+r}\sum_{\delta=0}^r\frac{q^{\delta+1}+q^{-\delta}}{q+1}\cdot q^{-(r-\delta)(r+1+\delta)}
\sum_{i=0}^{r-\delta}\qbinom{2r+1-2i}{r-\delta-i}_{-q}\Sat^\circ_N(\tT^\circ_{N;i}) \\
&=\frac{1}{q+1}\sum_{i=0}^r\(\sum_{\delta=0}^{r-i}(q^{2\delta+1}+1)q^{\delta^2}\qbinom{2(r-i)+1}{r-i-\delta}_{-q}\)\Sat^\circ_N(\tT^\circ_{N;i}) \\
&=\frac{1}{q+1}\sum_{i=0}^r\(\sum_{\delta=-(r-i)-1}^{r-i}q^{\delta^2}\qbinom{2(r-i)+1}{r-i-\delta}_{-q}\)\Sat^\circ_N(\tT^\circ_{N;i}).
\end{align*}
Thus the lemma is equivalent to the identity
\begin{align*}
\sum_{\delta=-k-1}^k q^{\delta^2}\qbinom{2k+1}{k-\delta}_{-q}=(q+1)(q^3+1)\cdots(q^{2k+1}+1)
\end{align*}
for every integer $k\geq 0$. By Lemma \ref{le:enumeration2}, we have
\[
\sum_{\delta=-k-1}^{k+1} q^{\delta^2}\qbinom{2k+2}{k+1-\delta}_{-q}=(q+1)(q^3+1)\cdots(q^{2k+1}+1).
\]
Thus, it remains to show
\[
\sum_{\delta=-k-1}^{k+1} q^{\delta^2}\qbinom{2k+2}{k+1-\delta}_{-q}=\sum_{\delta=-k-1}^k q^{\delta^2}\qbinom{2k+1}{k-\delta}_{-q}.
\]
However, the difference equals
\begin{align*}
\sum_{\delta=-k-1}^{k+1} q^{\delta^2}\(\qbinom{2k+2}{k+1-\delta}_{-q}-\qbinom{2k+1}{k-\delta}_{-q}\)
&=\sum_{\delta=-k-1}^{k+1} q^{\delta^2}(-q)^{k+1-\delta}\qbinom{2k+1}{k+1-\delta}_{-q} \\
&=(-q)^{k+1}\sum_{\delta=-k-1}^{k+1} (-1)^\delta q^{\delta^2-\delta}\qbinom{2k+1}{k+1-\delta}_{-q}
\end{align*}
which equals zero as the term of $\delta$ and the term of $-\delta+1$ cancel each other for $-k\leq\delta\leq k+1$ and the term with $\delta=-k-1$ vanishes. The lemma follows.
\end{proof}

\begin{lem}\label{le:enumeration_odd_2}
We have the identity
\[
q^{r^2+r}\prod_{i=1}^r\(\bmu_i-2\)=\sum_{\delta=0}^r\td_{\delta,q}\cdot\Sat^\circ_N(\tT^\circ_{N;r-\delta})
\]
in $\dZ[q^{-1}][\dX^*_N]^\sym$, in which the numbers $\td_{\delta,q}$ are introduced in Notation \ref{no:numerical}.
\end{lem}

\begin{proof}
By Lemma \ref{le:character3} and Lemma \ref{le:enumeration1}, we have
\begin{align*}
q^{r^2+r}\prod_{i=1}^r\(\bmu_i-2\)
&=q^{r^2+r}\sum_{\delta=0}^r(-1)^\delta(2\delta+1)\cdot\chi(\rho_{N;r-\delta}) \\
&=q^{r^2+r}\sum_{\delta=0}^r(-1)^\delta(2\delta+1)\cdot q^{-(r-\delta)(r+1+\delta)}\sum_{i=0}^{r-\delta}\qbinom{2r+1-2i}{r-\delta-i}_{-q}\Sat^\circ_N(\tT^\circ_{N;i}) \\
&=\sum_{i=0}^r\(\sum_{\delta=0}^{r-i}(-1)^\delta(2\delta+1)q^{\delta(\delta+1)}
\qbinom{2(r-i)+1}{r-i-\delta}_{-q}\)\Sat^\circ_N(\tT^\circ_{N;i}) \\
&=\sum_{\delta=0}^r\td_{\delta,q}\cdot\Sat^\circ_N(\tT^\circ_{N;r-\delta}).
\end{align*}
The lemma is proved.
\end{proof}

\begin{proposition}\label{pr:enumeration_odd}
Let $L$ be a $\dZ[q^{-1}]$-ring. Consider an $N$-tuple $\balpha=(\alpha_1,\dots,\alpha_N)\in L^N$ satisfying $\alpha_i\alpha_{N+1-i}=1$, which determines a homomorphism $\phi^\circ_{\balpha}\colon\dT^\circ_N\to L$ as in Notation \ref{no:parameter}.
\begin{enumerate}
  \item We have
     \[
     \phi^\circ_{\balpha}(\tI^\circ_N)=q^{r^2+r}\prod_{i=1}^r\(\alpha_i+\frac{1}{\alpha_i}+q+\frac{1}{q}\).
     \]

  \item We have
     \[
     \phi^\circ_{\balpha}(\tT^\circ_N)=q^{r^2+r}\prod_{i=1}^r\(\alpha_i+\frac{1}{\alpha_i}-2\),
     \]
     where
     \[
     \tT^\circ_N\coloneqq\sum_{\delta=0}^r\td_{r-\delta,q}\cdot\tT^\circ_{N;\delta}
     \]
     in which the numbers $\td_{r-\delta,q}$ are introduced in Notation \ref{no:numerical}.
\end{enumerate}
\end{proposition}

\begin{proof}
Part (1) follows from Lemma \ref{le:hecke_circ} and Lemma \ref{le:enumeration_odd_1}. Part (2) follows from Lemma \ref{le:enumeration_odd_2}.
\end{proof}

\begin{lem}\label{le:enumeration_odd_0}
We have
\[
\tT^{\bullet\circ}_N\circ\tT^\circ_N=\((q+1)^2\tT^\bullet_N+\tT^{\bullet\circ}_N\circ\tT^{\circ\bullet}_N\)\circ\tT^{\bullet\circ}_N
\]
in $\dZ[\rK^\bullet_N\backslash\rU(\rV_N)(F^+)/\rK^\circ_N]$, where $\tT^\circ_N$ is defined in Proposition \ref{pr:enumeration_odd}(2), and
\[
\tT^\bullet_N\coloneqq\sum_{\delta=0}^{r-1}\td^\bullet_{r-\delta,q}\cdot\tT^\bullet_{N;\delta}.
\]
\end{lem}

This lemma is a hard exercise in combinatorics. In fact, our proof below is by brutal force; it would be interesting to find a conceptual proof.

\begin{proof}
It suffices to show that for every element $f\in\dZ[\Lat^\circ_N]$, we have
\begin{align}\label{eq:enumeration_odd_0}
\((q+1)^2\tT^\bullet_N+\tT^{\bullet\circ}_N\circ\tT^{\circ\bullet}_N\)(\tT^{\bullet\circ}_N(f))
=\tT^{\bullet\circ}_N(\tT^\circ_N(f))
\end{align}
in $\dZ[\Lat^\bullet_N]$. Without loss of generality, we may just consider their values on $\Lambda^\bullet_N$.

For every $\rL\in\Lat^\circ_N$ and $0\leq\delta\leq r$, we denote
\begin{itemize}[label={\ding{109}}]
  \item $c^\bullet_\delta(\rL)$ the number of $\rL^\bullet\in\Lat^\bullet_N$ satisfying $\rL\subseteq \rL^\bullet$ and $(\rL^\bullet+\Lambda^\bullet_N)/\Lambda^\bullet_N\simeq(O_F/\fp)^{\oplus\delta}$; and

  \item $c^\circ_\delta(\rL)$ the number of $\rL^\circ\in\Lat^\circ_N$ satisfying $\rL^\circ\subseteq\Lambda^\bullet_N$ and $\rL/(\rL\cap \rL^\circ)\simeq(O_F/\fp)^{\oplus\delta}$.
\end{itemize}
We then have
\begin{align*}
(\tT^\bullet_{N;\delta}(\tT^{\bullet\circ}_N(f)))(\Lambda^\bullet_N)&=\sum_{\rL\in\Lat^\circ_N}c^\bullet_\delta(\rL)\cdot f(\rL),\\
(\tT^{\bullet\circ}_N(\tT^\circ_{N;\delta}(f)))(\Lambda^\bullet_N)&=\sum_{\rL\in\Lat^\circ_N}c^\circ_\delta(\rL)\cdot f(\rL).
\end{align*}
We claim the following identities
\begin{align}\label{eq:enumeration_odd_1}
c^\bullet_\delta(\rL)=
\begin{dcases}
q^{(\delta-\gamma)(\delta-\gamma+2)}\qbinom{r-\gamma}{\delta-\gamma}_{q^2},
&\text{if $(\rL+\Lambda^\bullet_N)/\Lambda^\bullet_N\simeq(O_F/\fp)^{\oplus\gamma}$ for some $0\leq\gamma\leq\delta$;} \\
0, &\text{otherwise;}
\end{dcases}
\end{align}
\begin{align}\label{eq:enumeration_odd_2}
c^\circ_\delta(\rL)=
\begin{dcases}
q^{(\delta-\gamma)^2}\qbinom{r-\gamma}{\delta-\gamma}_{q^2},
&\text{if $(\rL+\Lambda^\bullet_N)/\Lambda^\bullet_N\simeq(O_F/\fp)^{\oplus\gamma}$ for some $0\leq\gamma\leq\delta$;} \\
0, &\text{otherwise.}
\end{dcases}
\end{align}

For \eqref{eq:enumeration_odd_1}, we must have $(\rL+\Lambda^\bullet_N)/\Lambda^\bullet_N\subseteq(\rL^\bullet+\Lambda^\bullet_N)/\Lambda^\bullet_N\simeq(O_F/\fp)^{\oplus\delta}$. Thus, the otherwise case is confirmed. Suppose that $(\rL+\Lambda^\bullet_N)/\Lambda^\bullet_N\simeq(O_F/\fp)^{\oplus\gamma}$ for some $0\leq\gamma\leq\delta$. Then $(\fp\Lambda^\bullet_N+\rL)/\rL$ is an isotropic subspace of $\fp^{-1}\rL/\rL$ of dimension $\gamma$. Moreover, $c^\bullet_\delta(\rL)$ is the same as the number of maximal isotropic subspaces of $((\fp\Lambda^\bullet_N+\rL)/\rL)^\perp/((\fp\Lambda^\bullet_N+\rL)/\rL)$ whose intersection with (the image of) $(\fp^{-1}\rL\cap\Lambda^\bullet_N+\rL)/\rL$, which itself is a maximal isotropic subspace, has dimension $r-\delta$. Thus, we obtain \eqref{eq:enumeration_odd_1} by Lemma \ref{le:enumeration_odd_3} below since $((\fp\Lambda^\bullet_N+\rL)/\rL)^\perp/((\fp\Lambda^\bullet_N+\rL)/\rL)$ has dimension $2r+1-2\gamma$.

For \eqref{eq:enumeration_odd_2}, we must have $(\rL+\Lambda^\bullet_N)/\Lambda^\bullet_N\simeq\rL/(\rL\cap\Lambda^\bullet_N)$ which is a quotient of $\rL/(\rL\cap \rL^\circ)\simeq(O_F/\fp)^{\oplus\delta}$. Thus, the otherwise case is confirmed. Suppose that $(\rL+\Lambda^\bullet_N)/\Lambda^\bullet_N\simeq(O_F/\fp)^{\oplus\gamma}$ for some $0\leq\gamma\leq\delta$. Then $(\rL+\Lambda^\bullet_N)/\Lambda^\bullet_N$ is an isotropic subspace of $\fp^{-1}\Lambda^\bullet_N/\Lambda^\bullet_N$ of dimension $\gamma$. Moreover, $c^\circ_\delta(\rL)$ is the same as the number of maximal isotropic subspaces of $((\rL+\Lambda^\bullet_N)/\Lambda^\bullet_N)^\perp/((\rL+\Lambda^\bullet_N)/\Lambda^\bullet_N)$ whose intersection with (the image of) $(\fp^{-1}\Lambda^\bullet_N\cap\fp^{-1}\rL+\Lambda^\bullet_N)/\Lambda^\bullet_N$, which itself is a maximal isotropic subspace, has dimension $r-\delta$. Thus, we obtain \eqref{eq:enumeration_odd_2} by Lemma \ref{le:enumeration_odd_3} since $((\rL+\Lambda^\bullet_N)/\Lambda^\bullet_N)^\perp/((\rL+\Lambda^\bullet_N)/\Lambda^\bullet_N)$ has dimension $2r-2\gamma$.

Now we come back to the values of \eqref{eq:enumeration_odd_0} on $\Lambda^\bullet_N$. By a similar proof of Lemma \ref{le:hecke_circ}, we have
\[
\tT^{\bullet\circ}_N\circ\tT^{\circ\bullet}_N=
\tT^\bullet_{N;r}+(q+1)\tT^\bullet_{N;r-1}+(q+1)(q^3+1)\tT^\bullet_{N;r-2}+\cdots+\prod_{i=1}^r(q^{2i-1}+1)\tT^\bullet_{N;0}
\]
in $\dT^\bullet_N$. Then under Notation \ref{no:numerical}, we have
\begin{align}\label{eq:enumeration_odd_5}
&\quad\((q+1)^2\tT^\bullet_N+\tT^{\bullet\circ}_N\circ\tT^{\circ\bullet}_N\)\circ\tT^{\bullet\circ}_N \\
&=\tT^\bullet_{N;r}\circ\tT^{\bullet\circ}_N+\sum_{\delta=0}^{r-1}
\((q+1)\td_{r-\delta,q}+(-q)^{r-\delta+1}(q+1)(q^3+1)\cdots(q^{2(r-\delta)-1}+1)\)\tT^\bullet_{N;\delta}\circ\tT^{\bullet\circ}_N. \notag
\end{align}
By \eqref{eq:enumeration_odd_1}, \eqref{eq:enumeration_odd_2} and \eqref{eq:enumeration_odd_5}, the lemma is equivalent to that for every integer $k\geq 0$, we have
\begin{align*}
\resizebox{\hsize}{!}{
\xymatrix{
\displaystyle
\sum_{\delta=0}^k\td_{k-\delta,q}q^{\delta^2}\qbinom{k}{\delta}_{q^2}
=q^{k(k+2)}+\sum_{\delta=0}^{k-1}
\((q+1)\td_{k-\delta,q}+(-q)^{k-\delta+1}(q+1)(q^3+1)\cdots(q^{2(k-\delta)-1})\)q^{\delta(\delta+2)}\qbinom{k}{\delta}_{q^2},
}
}
\end{align*}
or equivalently,
\begin{align}\label{eq:enumeration_odd_3}
\resizebox{\hsize}{!}{
\xymatrix{
\displaystyle
\sum_{\delta=0}^k\td_{\delta,q}q^{(k-\delta)^2}\qbinom{k}{\delta}_{q^2}
=q^{k(k+2)}+\sum_{\delta=1}^{k}
\((q+1)\td_{\delta,q}+(-q)^{\delta+1}(q+1)(q^3+1)\cdots(q^{2\delta-1}+1)\)q^{(k-\delta)(k-\delta+2)}\qbinom{k}{\delta}_{q^2}.
}
}
\end{align}
By Lemma \ref{le:enumeration3}, we have
\begin{align*}
&\quad (-q)^{\delta+1}(q+1)(q^3+1)\cdots(q^{2\delta-1}+1) \\
&=-q\sum_{j=-\delta-1}^\delta(-1)^j j q^{j^2+j}\qbinom{2\delta+1}{\delta-j}_{-q}
+q\sum_{j=-\delta}^\delta(-1)^j j q^{j^2+j}\qbinom{2\delta}{\delta-j}_{-q} \\
&=-q\td_{\delta,q}+q\sum_{j=-\delta}^\delta(-1)^j j q^{j^2+j}\qbinom{2\delta}{\delta-j}_{-q}.
\end{align*}
Thus, \eqref{eq:enumeration_odd_3} is equivalent to
\begin{align*}
\sum_{\delta=0}^k\td_{\delta,q}q^{(k-\delta)^2}\qbinom{k}{\delta}_{q^2}
=\sum_{\delta=0}^{k}\(\td_{\delta,q}+q\sum_{j=-\delta}^\delta(-1)^j j q^{j^2+j}\qbinom{2\delta}{\delta-j}_{-q}\)q^{(k-\delta)(k-\delta+2)}\qbinom{k}{\delta}_{q^2},
\end{align*}
or equivalently,
\begin{align}\label{eq:enumeration_odd_4}
\sum_{\delta=0}^k\td_{\delta,q}q^{(k-\delta)^2}(q^{2(k-\delta)}-1)\qbinom{k}{\delta}_{q^2}=
-\sum_{\delta=0}^{k}\sum_{j=-\delta}^\delta(-1)^j j q^{j^2+j}\qbinom{2\delta}{\delta-j}_{-q}q^{(k-\delta+1)^2}\qbinom{k}{\delta}_{q^2}.
\end{align}
However, we have
\begin{align*}
&\quad\sum_{\delta=0}^k\td_{\delta,q}q^{(k-\delta)^2}(q^{2(k-\delta)}-1)\qbinom{k}{\delta}_{q^2} \\
&=\sum_{\delta=0}^{k-1}\td_{\delta,q}q^{(k-\delta)^2}(q^{2(k-\delta)}-1)\qbinom{k}{\delta}_{q^2} \\
&=\sum_{\delta=0}^{k-1}\sum_{j=-\delta-1}^\delta(-1)^j j q^{j^2+j}\qbinom{2\delta+1}{\delta-j}_{-q}q^{(k-\delta)^2}(q^{2(k-\delta)}-1)\qbinom{k}{\delta}_{q^2} \\
&=\sum_{\delta=0}^{k-1}\sum_{j=-\delta-1}^\delta(-1)^j j q^{j^2+j}\qbinom{2\delta+1}{\delta-j}_{-q}q^{(k-\delta)^2}(q^{2\delta+2}-1)\qbinom{k}{\delta+1}_{q^2} \\
&=\sum_{\delta=0}^{k-1}\sum_{j=-\delta-1}^\delta(-1)^j j q^{(k-\delta)^2+j^2+j}((-q)^{2\delta+2}-1)\qbinom{2\delta+1}{\delta-j}_{-q}\qbinom{k}{\delta+1}_{q^2} \\
&=\sum_{\delta=0}^{k-1}\sum_{j=-\delta-1}^\delta(-1)^j j q^{(k-\delta)^2+j^2+j}((-q)^{\delta-j+1}-1)\qbinom{2\delta+2}{\delta-j+1}_{-q}\qbinom{k}{\delta+1}_{q^2} \\
&=\sum_{\delta=1}^k\sum_{j=-\delta}^{\delta-1}(-1)^j j q^{(k+1-\delta)^2+j^2+j}((-q)^{\delta-j}-1)\qbinom{2\delta}{\delta-j}_{-q}\qbinom{k}{\delta}_{q^2} \\
&=\sum_{\delta=0}^k\sum_{j=-\delta}^{\delta}(-1)^j j q^{(k+1-\delta)^2+j^2+j}((-q)^{\delta-j}-1)\qbinom{2\delta}{\delta-j}_{-q}\qbinom{k}{\delta}_{q^2}.
\end{align*}
Thus, \eqref{eq:enumeration_odd_4} is equivalent to
\[
\sum_{\delta=0}^k\sum_{j=-\delta}^{\delta}(-1)^j j q^{(k+1-\delta)^2+j^2+j}(-q)^{\delta-j}\qbinom{2\delta}{\delta-j}_{-q}\qbinom{k}{\delta}_{q^2}=0,
\]
which is obvious since
\[
\sum_{j=-\delta}^{\delta}j q^{j^2}\qbinom{2\delta}{\delta-j}_{-q}=0.
\]
The lemma is finally proved.
\end{proof}

\begin{lem}\label{le:enumeration_odd_3}
Let $V$ be a (nondegenerate) hermitian space over $O_F/\fp$ of dimension $m\geq 1$ with $r=\lfloor\tfrac{m}{2}\rfloor$, and $Y_0\subseteq V$ a maximal isotropic subspace. Then the number of maximal isotropic subspaces $Y\subseteq V$ satisfying $\dim_{O_F/\fp}(Y\cap Y_0)=r-s$ with $0\leq s\leq r$ is given by
\[
\begin{dcases}
q^{s(s+2)}\qbinom{r}{s}_{q^2}, &\text{if $m=2r+1$;} \\
q^{s^2}\qbinom{r}{s}_{q^2}, &\text{if $m=2r$.}
\end{dcases}
\]
\end{lem}

\begin{proof}
We will prove the case for $m$ odd and leave the similar case for $m$ even to the readers. We fix an integer $0\leq s\leq r$. It is easy to see that the number of choices of the intersection $Y\cap Y_0$ (of dimension $r-s$) is
\[
\frac{(q^{2r}-1)(q^{2(r-1)}-1)\cdots(q^{2(r-s+1)}-1)}{(q^{2s}-1)(q^{2(s-1)}-1)\cdots(q^2-1)}
=\qbinom{r}{s}_{q^2}.
\]
Then we count the number of $Y$ with $Y\cap Y_0$ fixed. We take a basis $\{e_{-r},\dots,e_r\}$ of $V$ such that $(e_{-i},e_j)_V=\delta_{i,j}$ for $0\leq i,j\leq r$; $Y_0$ is spanned by $\{e_{-r},\dots,e_{-1}\}$; and $Y\cap Y_0$ is spanned by $\{e_{-r},\dots,e_{-s-1}\}$. Let $\{f_1,\dots,f_s\}$ be an element in $Y^s$ such that $\{e_{-r},\dots,e_{-s-1},f_1,\dots,f_s\}$ form a basis of $Y$. Then since $Y$ is isotropic, the coefficients on $\{e_{s+1},\dots,e_r\}$ of each $f_i$ have to be zero. In particular, there is unique such element $\{f_1,\dots,f_s\}\in Y^s$ that is of the form
\[
(f_1,\dots,f_s)=(e_1,\dots,e_s)+(e_{-s},\dots,e_{-1},e_0)
\begin{pmatrix}
A \\
v
\end{pmatrix}
\]
with (uniquely determined) $A\in\rM_{s,s}(O_F/\fp)$ and $v\in\rM_{1,s}(O_F/\fp)$. Moreover, the isotropic condition on $Y$ is equivalent to that $\tp{A}^\tc+A+\tp{v}^\tc\cdot v=0$, where $\tc$ denotes the Galois involution of $F/F^+$. It follows that the number for such $Y$ with given $Y\cap Y_0$ (of dimension $r-s$) is $q^{s(s+2)}$. Thus, the lemma follows.
\end{proof}

\section{Some representation theory for unitary groups}
\label{ss:b}

In this section, we prove several results for representations of unitary groups. Unless specified otherwise, all representations will have coefficients in $\dC$. In \S\ref{ss:local_bc}, we recall some general facts about the local base change for unitary groups. In \S\ref{ss:omega}, we study the representation appearing in the cohomology of Fermat hypersurfaces, and also compute the local base change of some admissible representations with nonzero Iwahori fixed vectors. In \S\ref{ss:endoscopic}, we collect everything we need from the endoscopic classification for unitary groups in Proposition \ref{pr:arthur} and derive two corollaries from it.

\subsection{Local base change for unitary groups}
\label{ss:local_bc}

In this subsection, we fix an unramified quadratic extension $F/F^+$ of nonarchimedean local fields. For every element $\alpha\in\dC^\times$, we denote by $\ul\alpha\colon F^\times\to\dC^\times$ the unramified character that sends every uniformizer to $\alpha$.

Consider a hermitian space $\rV$ over $F$ (with respect to $F/F^+$) of rank $N$. Put $G\coloneqq\rU(\rV)$. For an irreducible admissible representation $\pi$ of $G(F^+)$, we denote by $\BC(\pi)$ its base change, which is an irreducible admissible representation of $\GL_N(F)$. Such local base change is defined by \cite{Rog90} when $N\leq 3$ and by \cites{Mok15,KMSW} for general $N$.

We review the construction of $\BC(\pi)$ in certain special cases. For a parabolic subgroup $P$ of $G$ and an admissible representation $\sigma$ of $P(F^+)$, we denote by $\rI^G_P(\sigma)$ the normalized parabolic induction, which is an admissible representation of $G(F^+)$. Fix a minimal parabolic subgroup $P_\mnm$ of $G$.

We first review Langlands classification of irreducible admissible representations of $G(F^+)$ (see, for example, \cite{Kon03}).  For an irreducible admissible representation $\pi$ of $G(F^+)$, there is a unique parabolic subgroup $P$ of $G$ containing $P_\mnm$ with Levi quotient $M_P$, a unique tempered representation $\tau$ of $M_P(F^+)$, and a unique strictly positive (unramified) character $\chi$ of $P_\pi(F^+)$, such that $\pi$ is isomorphic to the unique irreducible quotient of $\rI^G_P(\tau\chi)$, which we denote by $\rJ^G_P(\tau\chi)$, known as the \emph{Langlands quotient}. Suppose that $\pi\simeq\rJ^G_P(\tau\chi)$ is a Langlands quotient. Then we may write
\[
M_P=G_0\times\Res_{F/F^+}\GL_{r_1}\times\cdots\times\Res_{F/F^+}\GL_{r_t}
\]
with $G_0$ the unitary factor, under which
\[
\chi=\b{1}\boxtimes\(\ul{\alpha_1}\circ\r{det}_{r_1}\)\boxtimes\cdots\boxtimes\(\ul{\alpha_t}\circ\r{det}_{r_t}\)
\]
for unique real numbers $1<\alpha_1<\cdots<\alpha_t$, where $\r{det}_r$ denotes the determinant on $\GL_r(F)$. Suppose that $\tau=\tau_0\boxtimes\tau_1\boxtimes\cdots\boxtimes\tau_t$ under the above decomposition. Consider a standard parabolic subgroup $P'$ of $\GL_N$ whose Levi is $\GL_{r_t}\times\cdots\times\GL_{r_1}\times\GL_{N_0}\times\GL_{r_1}\times\cdots\times\GL_{r_t}$. Then $\BC(\pi)$ is isomorphic to
\[
\rJ^{\GL_N}_{P'}
\(\tau_t^{\vee\tc}\(\ul{\alpha_t^{-1}}\circ\r{det}_{r_t}\)\boxtimes\cdots\boxtimes\tau_1^{\vee\tc}\(\ul{\alpha_1^{-1}}\circ\r{det}_{r_1}\)
\boxtimes\BC(\tau_0)\boxtimes\tau_1\(\ul{\alpha_1}\circ\r{det}_{r_1}\)\boxtimes\cdots\boxtimes\tau_t\(\ul{\alpha_t}\circ\r{det}_{r_t}\)\)
\]
which is a Langlands quotient of $\GL_N(F)$. Here, $\tau^\tc$ stands for $\tau\circ\tc$.

We then review the construction of tempered representations from discrete series representations (see, for example, \cite{Jan14}). Let $\tau$ be an irreducible admissible tempered representation of $G(F^+)$. Then there is a unique parabolic subgroup $P$ of $G$ containing $P_\mnm$, and a discrete series representation $\sigma$ of $M_P(F^+)$ such that $\tau$ is a direct summand of $\rI^G_P(\sigma)$. In fact, $\rI^G_P(\sigma)$ is a direct sum of finitely many tempered representations of multiplicity one. Write $\sigma=\sigma_0\boxtimes\sigma_1\boxtimes\cdots\boxtimes\sigma_t$, similar to the previous case. Then under the same notation, we have
\[
\BC(\tau)\simeq\rI^{\GL_N}_{P'}\(\sigma_t^{\vee\tc}\boxtimes\cdots\boxtimes\sigma_1^{\vee\tc}
\boxtimes\BC(\sigma_0)\boxtimes\sigma_1\boxtimes\cdots\boxtimes\sigma_t\)
\]
which is an irreducible admissible representation of $\GL_N(F)$.

Finally, if $\pi$ is an irreducible admissible representation of $G(F^+)$ that is a constituent of an unramified principal series, then $\BC(\pi)$ is a constituent of an unramified principal series of $\GL_N(F)$. Thus, it makes sense to talk about the Satake parameter of $\BC(\pi)$, denoted by $\balpha(\BC(\pi))$.

In what follows, we will suppress the parabolic subgroup $P'$ of $\GL_N$ when it is clear. We denote by $\St_N$ the Steinberg representation of $\GL_N(F)$.

\subsection{Tate--Thompson representations}
\label{ss:omega}

In this subsection, let $F/F^+$ be as in the previous subsection, with residue field extension $\kappa/\kappa^+$. Let $q$ be the residue cardinality of $F^+$ and $\fp$ the maximal ideal of $O_F$.

Let $N\geq 2$ be an integer with $r\coloneqq\lfloor\tfrac{N}{2}\rfloor$. Consider a hermitian space $\rV_N$ over $F$ of rank $N$ together with a self-dual lattice $\Lambda_N$. Put $\rU_N\coloneqq\rU(\rV_N)$, and let $\rK_N$ be the stabilizer of $\Lambda_N$ which is a hyperspecial maximal subgroup of $\rU_N(F^+)$. Put $\bar\Lambda_N\coloneqq\Lambda_N\otimes_{O_{F^+}}\kappa^+$ and $\bar\rU_N\coloneqq\rU(\bar\Lambda_N)$. Then we have the reduction homomorphism $\rK_N\to\bar\rU_N(\kappa^+)$.

Let $\Iso(\bar\Lambda_N)\subseteq\dP(\bar\Lambda_N)$ be the isotropic locus, that is, it parameterizes hyperplanes $H$ of $\bar\Lambda_N$ satisfying $H^\perp\subseteq H$. Then $\Iso(\bar\Lambda_N)$ is a smooth hypersurface in $\dP(\bar\Lambda_N)$, known as the \emph{Fermat hypersurface}. In particular, $\Iso(\bar\Lambda_N)$ has dimension $N-2$ and admits a natural action by $\bar\rU_N(\kappa^+)$. For a rational prime $\ell$ that is invertible in $\kappa$, put
\[
\rH^\prim(\Iso(\bar\Lambda_N)_{\ol\kappa},\ol\dQ_\ell)\coloneqq
\Ker\(\cup c_1(\cO_{\dP(\bar\Lambda_N)}(1))\colon
\rH^{N-2}_\et(\Iso(\bar\Lambda_N)_{\ol\kappa},\ol\dQ_\ell)\to\rH^N_\et(\Iso(\bar\Lambda_N)_{\ol\kappa},\ol\dQ_\ell(1))\).
\]
It is well-known by Tate--Thompson that (see, for example, \cite{HM78}) there is a unique irreducible representation $\Omega_N$ of $\bar\rU_N(\kappa^+)$ such that $\Omega_N$ is isomorphic to $\iota_\ell^{-1}\rH^\prim(\Iso(\bar\Lambda_N)_{\ol\kappa},\ol\dQ_\ell)$ as representations of $\bar\rU_N(\kappa^+)$ for every isomorphism $\iota_\ell\colon\dC\xrightarrow{\sim}\ol\dQ_\ell$. We call $\Omega_N$ the \emph{Tate--Thompson representation}. We often regard $\Omega_N$ as a representation of $\rK_N$ by inflation according to the context.

To describe $\Omega_N$, we first recall some notation from parabolic induction of finite reductive groups. For every $N$, we fix a Borel subgroup $\rP_N$ of $\bar\rU_N$. For positive integers $r_1,\dots,r_t$ satisfying $r_1+\cdots+r_t\leq r$, we obtain a parabolic subgroup $\rP_N^{(r_1,\dots,r_t)}$ of $\bar\rU_N$ containing $\rP_N$, whose Levi quotient $\rM_N^{(r_1,\dots,r_t)}$ is isomorphic to
$\bar\rU_{N-2(r_1+\cdots+r_t)}\times\Res_{\kappa/\kappa^+}\GL_{r_1}\times\cdots\times\Res_{\kappa/\kappa^+}\GL_{r_t}$. For example, we have $\rP_N^{(1^r)}=\rP_N$. Given a representation $\sigma$ of $\rM_N^{(r_1,\dots,r_t)}(\kappa^+)$, we denote by $\Ind_{\rP_N^{(r_1,\dots,r_t)}}^{\bar\rU_N}\sigma$ the parabolic induction, which is a representation of $\bar\rU_N(\kappa^+)$.

Now we suppose that $N=2r$ is even. The irreducible constituents of $\Ind_{\rP_N}^{\bar\rU_N}\mathbf{1}$ are parameterized by irreducible representations of the Weyl group $\rW_N\simeq\{\pm1\}^r\rtimes\fS_r$. For every irreducible representation $\epsilon$ of $\rW_N$, we denote by $\r{PS}(\epsilon)$ the corresponding irreducible representation of $\bar\rU_N(\kappa^+)$. We now specify a character $\epsilon^{\r{TT}}_N\colon\rW_N\to\{\pm1\}$ as the extension of the product homomorphism $\{\pm1\}^r\to\{\pm1\}$, which is invariant under the $\fS_r$-action, to $\rW_N$ that is trivial on $\{+1\}^r\rtimes\fS_r$.

\begin{proposition}\label{pr:omega}
We have
\begin{enumerate}
  \item When $N=2r$ is even, the representation $\Omega_N$ is isomorphic to $\r{PS}(\epsilon^{\r{TT}}_N)$.

  \item When $N=2r$ is even, $\Omega_N$ is the unique nontrivial irreducible representation of $\bar\rU_N(\kappa^+)$ satisfying $\dim\Omega_N^{\rP_N(\kappa^+)}=\dim\Omega_N^{\rP_N^{(r)}(\kappa^+)}=1$.

  \item The representation $\Omega_3$ is the (unique) cuspidal unipotent representation of $\bar\rU_3(\kappa^+)$.

  \item When $N=2r+1$ is odd with $r\geq 1$, the representation $\Omega_N$ is a multiplicity free constituent of $\Ind_{\rP_N^{(1^{r-1})}}^{\bar\rU_N}\Omega_3\boxtimes\mathbf{1}^{\boxtimes r-1}$.
\end{enumerate}
\end{proposition}

\begin{proof}
We recall some notion of Deligne--Lusztig characters. Let $\fS_N$ be the group of $N$-permutations, and $\fP_N$ its conjugacy classes which is canonically identified with the set of partitions of $N$. For every $\pi\in\fP_N$, we let $R_\pi$ be the Deligne--Lusztig character (of $\bar\rU_N(\kappa^+)$) \cite{DL76}*{Corollary~4.3} associated to the trivial representation of the maximal torus corresponding to $\pi$. Let $R_N$ be the character of the representation $\Omega_N$. Then by \cite{HM78}*{Theorem~1}, we have
\begin{align}\label{eq:omega1}
R_N=(-1)^{N+1}\sum_{\pi\in\fP_N}\frac{\chi_N(\pi)}{z_\pi}R_\pi
\end{align}
where $\chi_N$ is the character function (on $\fP_N$) of the unique nontrivial subrepresentation of the standard representation of $\fS_N$; and $N!/z_\pi$ is the cardinality of the conjugacy class $\pi$. By \cite{DL76}*{Theorem~6.8}, we have the following orthogonality relation
\begin{align}\label{eq:omega2}
\langle R_\pi, R_{\pi'}\rangle=
\begin{dcases}
0, &\text{if $\pi\neq\pi'$;} \\
z_\pi, &\text{if $\pi=\pi'$.}
\end{dcases}
\end{align}
We are ready to prove the proposition. In what follows, we write $(s^r)$ for the $r$-tuple $(s,\dots,s)$.

For (1), note that $\epsilon^{\r{TT}}_N$ is the unique nontrivial character of $\rW_N$ that is trivial on $\{+1\}^r\rtimes\fS_r$. Thus, (1) follows from (2) by \cite{Cur79}*{Theorem~4.4.5}.

For (2), we first show the uniqueness of $\Omega_N$. The condition $\dim\Omega_N^{\rP_N(\kappa^+)}=1$ implies that $\Omega_N$ is a constituent of $\Ind_{\rP_N}^{\bar\rU_N}\mathbf{1}$ corresponding to a character of $\rW_N$. However, there are only four characters of $\rW_N$, among which only the trivial character and $\epsilon^{\r{TT}}_N$ will give constituents with nonzero $\rP_N^{(r)}(\kappa^+)$-invariants. Thus, the uniqueness follows. For the identity $\dim\Omega_N^{\rP_N(\kappa^+)}=\dim\Omega_N^{\rP_N^{(r)}(\kappa^+)}=1$, it suffices to show that $\dim\Omega_N^{\rP_N(\kappa^+)}=1$ and $\Omega_N^{\rP_N^{(r)}(\kappa^+)}\neq 0$. Let $R'_{2r}$ be the character of $\Ind_{\rP_{2r}}^{\bar\rU_{2r}}\mathbf{1}$. Then by \cite{DL76}*{Proposition~8.2}, we have $R'_{2r}=R_{(2^r)}$. By \eqref{eq:omega1} and \eqref{eq:omega2}, we have
\begin{align*}
\langle R_{2r},R'_{2r}\rangle=\left\langle-\sum_{\pi\in\fP_{2r}}\frac{\chi_{2r}(\pi)}{z_\pi}R_\pi,R_{(2^r)}\right\rangle
=-\chi_{2r}((2^r))=-(-1)=1,
\end{align*}
which implies $\dim\Omega_N^{\rP_N(\kappa^+)}=1$. Let $\rY_N\subseteq\bar\Lambda_N$ be the maximal isotropic subspace stabilized by $\rP_N^{(r)}$. Then $\dP(\rY_N)$ is contained in $\Iso(\bar\Lambda_N)$, which gives rise to an element in $\CH^{r-1}(\Iso(\bar\Lambda_N))$. It is well-known that its cohomology class subtracted by $c_1(\cO_{\dP(\bar\Lambda_N)_{\ol\kappa}}(1))$ is a nonzero element in $\rH^\prim(\Iso(\bar\Lambda_N)_{\ol\kappa},\ol\dQ_\ell)(r-1)$, which is fixed by $\rP_N^{(r)}(\kappa^+)$ by construction. Thus, we have $\Omega_N^{\rP_N^{(r)}(\kappa^+)}\neq 0$; and (2) follows.

For (3), we have $R_3=\frac{1}{3}(R_{(1^3)}-R_{(3)})$ by \eqref{eq:omega1}. Then as computed in \cite{Pra}*{Example~6.2}, $\Omega_3$ is the unique cuspidal unipotent representation of $\bar\rU_3(\kappa^+)$.

For (4), let $R'_{2r+1}$ be the character of $\Ind_{\rP_{2r+1}^{(1^{r-1})}}^{\bar\rU_{2r+1}}\(\Omega_3\boxtimes\mathbf{1}^{\boxtimes r-1}\)$. Then by \cite{DL76}*{Proposition~8.2}, we have
\[
R'_{2r+1}=\frac{1}{3}\(R_{(2^{r-1},1^3)}-R_{(2^{r-1},3)}\).
\]
By \eqref{eq:omega1} and \eqref{eq:omega2}, we have
\begin{align*}
\langle R_{2r+1},R'_{2r+1}\rangle&=
\left\langle\sum_{\pi\in\fP_{2r+1}}\frac{\chi_{2r+1}(\pi)}{z_\pi}R_\pi,\frac{1}{3}\(R_{(2^{r-1},1^3)}-R_{(2^{r-1},3)}\)\right\rangle \\
&=\frac{1}{3}\(\chi_{2r+1}((2^{r-1},1^3))-\chi_{2r+1}((2^{r-1},3))\)=\frac{1}{3}(2-(-1))=1.
\end{align*}
Thus, (4) follows.
\end{proof}

From now on, we assume that $N=2r$ is even.

\begin{lem}\label{le:base_change_0}
Let $\pi$ be an irreducible admissible representation of $\rU_{2r}(F^+)$ such that $\pi\res_{\rK_{2r}}$ contains $\Omega_{2r}$ (hence is a constituent of an unramified principal series).
\begin{enumerate}
  \item If the Satake parameter of $\BC(\pi)$ contains neither $\{q,q^{-1}\}$ nor $\{-1,-1\}$, then $\pi\res_{\rK_{2r}}$ contains the trivial representation.

  \item If the Satake parameter of $\BC(\pi)$ contains $\{q,q^{-1}\}$, then there exists an element $(\alpha_2,\dots,\alpha_r)\in(\dC^\times)^{r-1}$ satisfying $1\leq|\alpha_2|\leq\cdots\leq|\alpha_r|$, unique up to permutation, such that $\BC(\pi)$ is isomorphic to the unique irreducible quotient of
      \[
      \rI^{\GL_{2r}}
      \(\ul{\alpha_r^{-1}}\boxtimes\cdots\boxtimes\ul{\alpha_2^{-1}}\boxtimes\St_2\boxtimes\ul{\alpha_2}\boxtimes\cdots\boxtimes\ul{\alpha_r}\).
      \]
\end{enumerate}
\end{lem}

\begin{proof}
We fix a decomposition
\[
\Lambda_{2r}=O_F e_{-r}\oplus\cdots\oplus O_F e_{-1} \oplus O_F e_1\oplus \cdots\oplus O_F e_r,
\]
in which $(e_{-i},e_j)=\delta_{ij}$ for $1\leq i,j\leq r$. For $0\leq i\leq r$, put
\[
\rV_{2i}\coloneqq F e_{-i}\oplus\cdots\oplus F e_{-1} \oplus F e_1\oplus \cdots\oplus F e_i,
\]
which is a hermitian subspace of $\rV_{2r}$. We take the minimal parabolic (Borel) subgroup $P_\mnm$ of $G\coloneqq\rU_{2r}$ to be the stabilizer of the flag $Fe_{-r}\subseteq\cdots\subseteq F e_{-r}\oplus\cdots\oplus F e_{-1}$. We also fix a Levi subgroup of $P_\mnm$ to be $\Res_{F/F^+}\GL(F e_1)\times\cdots\times\Res_{F/F^+}\GL(F e_r)$.

Put $K\coloneqq\rK_{2r}$, which is a hyperspecial maximal subgroup of $G(F^+)$. Let $I$ be the subgroup of $K$ of elements whose reduction modulo $\fp$ stabilizes the flag $\kappa e_{-r}\subseteq\cdots\subseteq\kappa e_{-r}\oplus\cdots\oplus\kappa e_{-1}$, which is an Iwahori subgroup of $G(F^+)$. Let $J$ be the subgroup of $K$ of elements whose reduction modulo $\fp$ stabilizes the subspace $\kappa e_{-r}\oplus\cdots\oplus\kappa e_{-1}$, which is a parahoric subgroup of $G(F^+)$. We clearly have $I\subseteq J\subseteq K$. Now we realize the Weyl group $\rW_{2r}\simeq\{\pm1\}^r\rtimes\fS_r$ explicitly as a subgroup of $K$. For $1\leq i\leq r$, we let $i$-th $-1$ in $\rW_{2r}$ correspond to the element that only switches $e_{-i}$ and $e_i$, denoted by $w_i$. For every $\sigma\in\fS_r$, we let $(1^r,\sigma)\in\rW_{2r}$ correspond to the element that sends $e_{\pm i}$ to $e_{\pm\sigma(i)}$, denoted by $w'_\sigma\in J$. Then $\{w_1,w'_{(1,2)},\dots,w'_{(r-1,r)}\}$ is a set of distinguished generators of $\rW_{2r}$. We recall the Bruhat decompositions
\[
K=\coprod_{w\in\rW_{2r}}I w I,\quad
K=\coprod_{i=0}^r J w_1\cdots w_i J.
\]
For $w\in\rW_{2r}$, we let $0\leq i(w)\leq r$ be the unique integer such that $w\in Jw_1\cdots w_{i(w)}J$.

By Proposition \ref{pr:omega}(2), we have a $K$-equivariant embedding $\Omega_{2r}\hookrightarrow\dC[I\backslash K]$, unique up to scalar, hence obtain a distinguished subspace $\Omega_{2r}^I\subseteq\dC[I\backslash K/I]$ of dimension one. We would like to find a generator of $\Omega_{2r}^I$. Now we compute the character of the $\dC[I\backslash K/I]$-module $\Omega_{2r}^I$. By Proposition \ref{pr:omega}(2), $\Omega_{2r}^I$ is contained in $\dC[J\backslash K/J]$. It follows that the element $\CF_{I w_1I}$ acts on $\Omega_{2r}^I$ by either $q$ or $-1$, in which the former case corresponds to the $K$-spherical one, which is not our case by Proposition \ref{pr:omega}(1). Thus, $\Omega_{2r}^I$ is spanned by the following function:
\[
f\coloneqq\sum_{w\in\rW_{2r}}(-q)^{-i(w)}\cdot\CF_{IwI}\in\dC[I\backslash K/I].
\]
For every element $\balpha=(\alpha_1,\dots,\alpha_r)\in(\dC^\times)^r$, we have the projection map
\[
\sP_{\balpha}\colon\dC[I\backslash K/I]\to\rI^G_{P_\mnm}\(\ul{\alpha_1}\boxtimes\cdots\boxtimes\ul{\alpha_r}\)^I
\]
defined at the beginning of \cite{Cas80}*{\S2}, which is $\dC[I\backslash K/I]$-equivariant. Put $\phi_{\balpha}\coloneqq\sP_{\balpha}(f)$.

Take an irreducible admissible representation $\pi$ of $\rU_{2r}(F^+)$ such that $\pi\res_K$ contains $\Omega_{2r}$. Then $\pi$ is a constituent of an unramified principal series. Now we separate the discussion.

Suppose that we are in the situation of (1). Then there exists an element $\balpha=(\alpha_1,\dots,\alpha_r)\in(\dC^\times)^r$ satisfying $1\leq|\alpha_1|\leq\cdots\leq|\alpha_r|$ and $\alpha_i\not\in\{-1,q\}$, unique up to permutation, such that $\pi$ is a constituent of $\rI^G_{P_\mnm}\(\ul{\alpha_1}\boxtimes\cdots\boxtimes\ul{\alpha_r}\)$. There exist a unique nonnegative integer $r_0$ and unique positive integers $r_1,\dots,r_t$ satisfying $r_0+\cdots+r_t=r$, such that
\[
1=|\alpha_1|=\cdots=|\alpha_{r_0}| < |\alpha_{r_0+1}|=\cdots=|\alpha_{r_0+r_1}| < \cdots < |\alpha_{r_0+\cdots+r_{t-1}+1}|=\cdots=|\alpha_r|
\]
holds. For every $1\leq i\leq t$, put
\[
\tau_i\coloneqq\rI^{\GL_{r_i}}\(\ul{\alpha_{r_0+\cdots+r_{i-1}+1}}\boxtimes\cdots\boxtimes\ul{\alpha_{r_0+\cdots+r_i}}\)
\otimes\(\ul{|\alpha_{r_0+\cdots+r_i}^{-1}|}\circ\r{det}_{r_i}\),
\]
which is an irreducible tempered representation of $\GL_{r_i}(F)$. Put $G_0\coloneqq\rU(\rV_{2r_0})$ and $P_{0\mnm}\coloneqq G_0\cap P_\mnm$. As $\ul{\alpha_1}\boxtimes\cdots\boxtimes\ul{\alpha_{r_0}}$ is a discrete series representation of $P_{0\mnm}(F^+)$, the parabolic induction
\[
\tau_0\coloneqq\rI^{G_0}_{P_{0\mnm}}\(\ul{\alpha_1}\boxtimes\cdots\boxtimes\ul{\alpha_{r_0}}\)
\]
is a finite direct sum of irreducible tempered representations of $G_0(F^+)$. As $\{\alpha_1,\dots,\alpha_{r_0}\}$ does not contain $-1$, $\tau_0$ is irreducible by \cite{Gol95}*{Theorem~1.4 \& Theorem~3.4}. In particular, we obtain a Langlands quotient
\[
\rJ^G_P\(\tau_0\boxtimes\(\boxtimes_{i=1}^t\tau_i\(\ul{|\alpha_{r_0+\cdots+r_i}|}\circ\r{det}_{r_i}\)\)\),
\]
where $P$ is the parabolic subgroup of $G$ containing $P_0$ whose Levi quotient is isomorphic to $G_0\times\Res_{F/F^+}\GL_{r_1}\times\cdots\times\Res_{F/F^+}\GL_{r_t}$. We claim that
\begin{align}\label{eq:base_change_0}
\phi_{\balpha}\neq 0\in\rJ^G_P\(\tau_0\boxtimes\(\boxtimes_{i=1}^t\tau_i\(\ul{|\alpha_{r_0+\cdots+r_i}|}\circ\r{det}_{r_i}\)\)\).
\end{align}
Assuming this claim, then $\pi$ is isomorphic to the above Langlands quotient, which is the unique irreducible quotient of $\rI^G_{P_\mnm}\(\ul{\alpha_1}\boxtimes\cdots\boxtimes\ul{\alpha_r}\)$. In particular, $\pi\res_{\rK_{2r}}$ contains the trivial representation. Thus, (1) follows.

Now we prove \eqref{eq:base_change_0}. Let $w\in\rW_{2r}$ be the element acting trivially on $\rV_{2r_0}$ and switching $e_{-(r_0+\cdots+r_{i-1}+j)}$ with $e_{r_0+\cdots+r_i+1-j}$ for every $1\leq j\leq r_i$ and then every $1\leq i\leq t$. By \cite{Kon03}*{Corollary~3.2}, \eqref{eq:base_change_0} is equivalent to
\begin{align}\label{eq:base_change_01}
T_w\phi_{\balpha}\neq 0,
\end{align}
where $T_w$ is the intertwining operator, which, in this case, is defined by an absolutely convergent integral
\begin{align*}
(T_w\phi_{\balpha})(g)=\int_{N(F^+)}\phi_{\balpha}(w^{-1}ng)\rd n,
\end{align*}
where $N$ is the unipotent radical of $P$ and the integral is absolutely convergent (see the discussion after \cite{Kon03}*{Proposition~2.2}). Since the eigenspace for the character of $\Omega_{2r}^I$ has dimension $1$, we must have
\[
T_w\phi_{\balpha}=C(\balpha)\phi_{w\balpha}
\]
for some complex number $C(\balpha)$. By \cite{Cas80}*{Theorem~3.4} and the continuity, we have
\[
C(\balpha)=\prod_{i=r_0+1}^r\(\frac{q-\alpha_i}{q(\alpha_i-1)}\prod_{|\alpha_j|<|\alpha_i|}\frac{\alpha_i-q^{-2}\alpha_j}{\alpha_i-\alpha_j}
\prod_{j=1}^{i-1}\frac{\alpha_i\alpha_j-q^{-2}}{\alpha_i\alpha_j-1}\),
\]
which is nonzero in the situation of (1). From this we obtain \eqref{eq:base_change_01}, hence \eqref{eq:base_change_0}.

Suppose that we are in the situation of (2). Then there exists an element $\balpha=(q,\alpha_2,\dots,\alpha_r)\in(\dC^\times)^r$ satisfying $1\leq|\alpha_2|\leq\cdots\leq|\alpha_r|$, unique up to permutation, such that $\pi$ is a constituent of
\[
\rI^G_{P_\mnm}\(\ul{q}\boxtimes\ul{\alpha_2}\boxtimes\cdots\boxtimes\ul{\alpha_r}\).
\]
Let $Q$ be the parabolic subgroup of $G$ stabilizing the flag $Fe_{-r}\subseteq\cdots\subseteq F e_{-r}\oplus\cdots\oplus F e_{-2}$, whose Levi quotient is $\rU(\rV_2)\times\Res_{F/F^+}\GL(F e_2)\times\cdots\times\Res_{F/F^+}\GL(F e_r)$. Then we have a canonical inclusion
\[
\rI^G_Q\(\Sp_2\boxtimes\ul{\alpha_2}\boxtimes\cdots\boxtimes\ul{\alpha_r}\)\subseteq
\rI^G_{P_\mnm}\(\ul{q}\boxtimes\ul{\alpha_2}\boxtimes\cdots\boxtimes\ul{\alpha_r}\)
\]
where $\Sp_2$ denotes the Steinberg representation of $\rU(\rV_2)(F^+)$. As $\CF_{Iw_1I}$ acts by $-1$ on $\phi_{\balpha}$, we have
\[
\phi_{\balpha}\in\rI^G_Q\(\Sp_2\boxtimes\ul{\alpha_2}\boxtimes\cdots\boxtimes\ul{\alpha_r}\).
\]
In particular, it follows that $\pi$ is a constituent of $\rI^G_Q\(\Sp_2\boxtimes\ul{\alpha_2}\boxtimes\cdots\boxtimes\ul{\alpha_r}\)$. To proceed, there exist unique positive integers $r_0,\dots,r_t$ satisfying $r_0+\cdots+r_t=r$, such that
\[
1=|\alpha_2|=\cdots=|\alpha_{r_0}| < |\alpha_{r_0+1}|=\cdots=|\alpha_{r_0+r_1}| < \cdots < |\alpha_{r_0+\cdots+r_{t-1}+1}|=\cdots=|\alpha_r|
\]
holds. For every $1\leq i\leq t$, put
\[
\tau_i\coloneqq\rI^{\GL_{r_i}}\(\ul{\alpha_{r_0+\cdots+r_{i-1}+1}}\boxtimes\cdots\boxtimes\ul{\alpha_{r_0+\cdots+r_i}}\)
\otimes\(\ul{|\alpha_{r_0+\cdots+r_i}^{-1}|}\circ\r{det}_{r_i}\),
\]
which is an irreducible tempered representation of $\GL_{r_i}(F)$. Put $G_0\coloneqq\rU(\rV_{2r_0})$ and $Q_0\coloneqq G_0\cap Q$. As $\Sp_2\boxtimes\ul{\alpha_2}\boxtimes\cdots\boxtimes\ul{\alpha_{r_0}}$ is a discrete series representation of $Q_0(F^+)$, the parabolic induction
\[
\rI^{G_0}_{Q_0}\(\Sp_2\boxtimes\ul{\alpha_2}\boxtimes\cdots\boxtimes\ul{\alpha_{r_0}}\)
\]
is a finite direct sum of irreducible tempered representations of $G_0(F^+)$. Let $\tau_0$ be the unique direct summand such that $\phi_{\balpha}$ is contained in the subspace
\[
\rI^G_P\(\tau_0\boxtimes\(\boxtimes_{i=1}^t\tau_i\(\ul{|\alpha_{r_0+\cdots+r_i}|}\circ\r{det}_{r_i}\)\)\)
\subseteq\rI^G_Q\(\Sp_2\boxtimes\ul{\alpha_2}\boxtimes\cdots\boxtimes\ul{\alpha_r}\),
\]
where $P$ is the parabolic subgroup of $G$ containing $P_0$ whose Levi quotient is isomorphic to $G_0\times\Res_{F/F^+}\GL_{r_1}\times\cdots\times\Res_{F/F^+}\GL_{r_t}$. In particular, we obtain a Langlands quotient
\[
\rJ^G_P\(\tau_0\boxtimes\(\boxtimes_{i=1}^t\tau_i\(\ul{|\alpha_{r_0+\cdots+r_i}|}\circ\r{det}_{r_i}\)\)\).
\]
By the same proof of \eqref{eq:base_change_0}, we obtain
\[
\phi_{\balpha}\neq 0\in\rJ^G_P\(\tau_0\boxtimes\(\boxtimes_{i=1}^t\tau_i\(\ul{|\alpha_{r_0+\cdots+r_i}|}\circ\r{det}_{r_i}\)\)\).
\]
In fact, in this case, we have the formula
\[
C(\balpha)=\prod_{i=r_0+1}^r\(\frac{q-\alpha_i}{q(\alpha_i-1)}\frac{\alpha_i-q^{-1}}{\alpha_i-q}
\prod_{\substack{j>1\\ |\alpha_j|<|\alpha_i|}}\frac{\alpha_i-q^{-2}\alpha_j}{\alpha_i-\alpha_j}
\prod_{j=1}^{i-1}\frac{\alpha_i\alpha_j-q^{-2}}{\alpha_i\alpha_j-1}\).
\]

Then $\BC(\pi)$ is isomorphic to the unique irreducible quotient of
\[
\rI^{\GL_{2r}}\(\(\boxtimes_{i=t}^1\tau_i^{\vee\tc}\(\ul{|\alpha_{r_0+\cdots+r_i}^{-1}|}\circ\r{det}_{r_i}\)\)
\boxtimes\BC(\tau_0)\boxtimes\(\boxtimes_{i=1}^t\tau_i\(\ul{|\alpha_{r_0+\cdots+r_i}|}\circ\r{det}_{r_i}\)\)\).
\]
However, $\BC(\tau_0)$ is isomorphic to
\begin{align*}
&\quad\rI^{\GL_{2r_0}}\(\ul{\alpha_{r_0}^{-1}}\boxtimes\cdots\boxtimes\ul{\alpha_2^{-1}}\boxtimes\BC(\Sp_2)
\boxtimes\ul{\alpha_2}\boxtimes\cdots\boxtimes\ul{\alpha_{r_0}}\) \\
&\simeq\rI^{\GL_{2r_0}}\(\ul{\alpha_{r_0}^{-1}}\boxtimes\cdots\boxtimes\ul{\alpha_2^{-1}}\boxtimes\St_2
\boxtimes\ul{\alpha_2}\boxtimes\cdots\boxtimes\ul{\alpha_{r_0}}\)
\end{align*}
which is irreducible. Thus, (2) follows.

The lemma is proved.
\end{proof}

\begin{remark}
In the situation of Lemma \ref{le:base_change_0}, the proof actually shows that if the Satake parameter of $\BC(\pi)$ does not contain $\{q,q^{-1}\}$ but possibly contains $\{-1,-1\}$, then $\pi$ is unramified with respect to either $\rK_{2r}$ or the other (conjugacy class of) hyperspecial maximal subgroup that is not conjugate to $\rK_{2r}$ in $\rU_{2r}(F^+)$.
\end{remark}

Let $\rV'_{2r}$ be another hermitian space over $F$ together with a lattice $\Lambda'_{2r}$ satisfying $\Lambda'_{2r}\subseteq (\Lambda'_{2r})^\vee$ and $(\Lambda'_{2r})^\vee/\Lambda'_{2r}\simeq\kappa$. Put $\rU'_{2r}\coloneqq\rU(\rV'_{2r})$, and let $\rK'_{2r}$ be the stabilizer of $\Lambda'_{2r}$ which is a special maximal subgroup of $\rU'_{2r}(F^+)$.

\begin{lem}\label{le:base_change_1}
Let $\pi'$ be an irreducible admissible representation of $\rU'_{2r}(F^+)$ such that $(\pi')^{\rK'_{2r}}\neq\{0\}$. Then there exists an element $(\alpha_2,\dots,\alpha_r)\in(\dC^\times)^{r-1}$ satisfying $1\leq|\alpha_2|\leq\cdots\leq|\alpha_r|$, unique up to permutation, such that  $\BC(\pi')$ is isomorphic to the unique irreducible quotient of
\[
\rI^{\GL_{2r}}\(\ul{\alpha_r^{-1}}\boxtimes\cdots\boxtimes\ul{\alpha_2^{-1}}
\boxtimes\St_2\boxtimes\ul{\alpha_2}\boxtimes\cdots\boxtimes\ul{\alpha_r}\).
\]
\end{lem}

\begin{proof}
We fix a decomposition
\[
\Lambda'_{2r}=O_F e_{-r}\oplus\cdots\oplus O_F e_{-2} \oplus \Lambda'_2\oplus O_F e_2\oplus \cdots\oplus O_F e_r,
\]
in which $(e_{-i},e_j)=\delta_{ij}$ for $2\leq i,j\leq r$. For $1\leq i\leq r$, put
\[
\rV'_{2i}\coloneqq F e_{-i}\oplus\cdots\oplus F e_{-2} \oplus \Lambda'_2\otimes_{O_F}F \oplus F e_2\oplus \cdots\oplus F e_i,
\]
which is a hermitian subspace of $\rV'_{2r}$. We take the minimal parabolic subgroup $P_\mnm$ of $G\coloneqq\rU'_{2r}$ to be the stabilizer of the flag $Fe_{-r}\subseteq\cdots\subseteq F e_{-r}\oplus\cdots\oplus F e_{-2}$. We also fix a Levi subgroup of $P_\mnm$ to be $\rU(\rV'_2)\times\Res_{F/F^+}\GL(F e_2)\times\cdots\times\Res_{F/F^+}\GL(F e_r)$. We have a similar embedding $\rW'_{2r}\hookrightarrow\rK'_{2r}$ of the Weyl group $\rW'_{2r}\simeq\rW_{2r-2}$. For every element $\balpha=(\alpha_2,\dots,\alpha_r)\in(\dC^\times)^{r-1}$, we let $\phi'_{\balpha}$ be the element in $\rI^G_{P_\mnm}\(\b{1}'_2\boxtimes\ul{\alpha_2}\boxtimes\cdots\boxtimes\ul{\alpha_r}\)$ that takes value $1$ on $\rK'_{2r}$, where $\b{1}'_2$ denotes the trivial representation of $\rU(\rV'_2)(F^+)$.

Take an irreducible admissible representation $\pi'$ of $G(F^+)$ such that $(\pi')^{\rK'_{2r}}\neq 0$. Then it is a constituent of an unramified principal series. In other words, there exists an element $\balpha=\{\alpha_2,\dots,\alpha_r\}\in(\dC^\times)^{r-1}$ satisfying $1\leq|\alpha_2|\leq\cdots\leq|\alpha_r|$, unique up to permutation, such that $\pi'$ is a constituent of
\[
\rI^G_{P_\mnm}\(\b{1}'_2\boxtimes\ul{\alpha_2}\boxtimes\cdots\boxtimes\ul{\alpha_r}\).
\]
To proceed, there exist unique positive integers $r_0,\dots,r_t$ satisfying $r_0+\cdots+r_t=r$, such that
\[
1=|\alpha_2|=\cdots=|\alpha_{r_0}| < |\alpha_{r_0+1}|=\cdots=|\alpha_{r_0+r_1}| < \cdots < |\alpha_{r_0+\cdots+r_{t-1}+1}|=\cdots=|\alpha_r|
\]
holds. For every $1\leq i\leq t$, put
\[
\tau_i\coloneqq\rI^{\GL_{r_i}}\(\ul{\alpha_{r_0+\cdots+r_{i-1}+1}}\boxtimes\cdots\boxtimes\ul{\alpha_{r_0+\cdots+r_i}}\)
\otimes\(\ul{|\alpha_{r_0+\cdots+r_i}^{-1}|}\circ\r{det}_{r_i}\),
\]
which is an irreducible tempered representation of $\GL_{r_i}(F)$. Put $G_0\coloneqq\rU(\rV'_{2r_0})$ and $P_{0\mnm}\coloneqq G_0\cap P_\mnm$. As $\b{1}'_2\boxtimes\ul{\alpha_2}\boxtimes\cdots\boxtimes\ul{\alpha_{r_0}}$ is a discrete series representation of $P_{0\mnm}(F^+)$, the parabolic induction
\[
\rI^{G_0}_{P_{0\mnm}}\(\b{1}'_2\boxtimes\ul{\alpha_2}\boxtimes\cdots\boxtimes\ul{\alpha_{r_0}}\)
\]
is a finite direct sum of irreducible tempered representations of $G_0(F^+)$. Let $\tau_0$ be the unique direct summand with nonzero invariants under $\rK'_{2r}\cap G_0(F^+)$. In particular, we obtain a Langlands quotient
\[
\rJ^G_P\(\tau_0\boxtimes\(\boxtimes_{i=1}^t\tau_i\(\ul{|\alpha_{r_0+\cdots+r_i}|}\circ\r{det}_{r_i}\)\)\),
\]
where $P$ is the parabolic subgroup of $G$ containing $P_0$ whose Levi quotient is isomorphic to $G_0\times\Res_{F/F^+}\GL_{r_1}\times\cdots\times\Res_{F/F^+}\GL_{r_t}$. We claim
\begin{align}\label{eq:base_change_1}
\rJ^G_P\(\tau_0\boxtimes\(\boxtimes_{i=1}^t\tau_i\(\ul{|\alpha_{r_0+\cdots+r_i}|}\circ\r{det}_{r_i}\)\)\)^{\rK'_{2r}}\neq\{0\}.
\end{align}
Assuming this claim, then $\BC(\pi')$ is isomorphic to the unique irreducible quotient of
\[
\rI^{\GL_{2r}}\(\(\boxtimes_{i=t}^1\tau_i^{\vee\tc}\(\ul{|\alpha_{r_0+\cdots+r_i}^{-1}|}\circ\r{det}_{r_i}\)\)
\boxtimes\BC(\tau_0)\boxtimes\(\boxtimes_{i=1}^t\tau_i\(\ul{|\alpha_{r_0+\cdots+r_i}|}\circ\r{det}_{r_i}\)\)\).
\]
However, $\BC(\tau_0)$ is isomorphic to
\begin{align*}
&\quad\rI^{\GL_{2r_0}}\(\ul{\alpha_{r_0}^{-1}}\boxtimes\cdots\boxtimes\ul{\alpha_2^{-1}}\boxtimes\BC(\b{1}')
\boxtimes\ul{\alpha_2}\boxtimes\cdots\boxtimes\ul{\alpha_{r_0}}\) \\
&\simeq\rI^{\GL_{2r_0}}\(\ul{\alpha_{r_0}^{-1}}\boxtimes\cdots\boxtimes\ul{\alpha_2^{-1}}\boxtimes\St_2
\boxtimes\ul{\alpha_2}\boxtimes\cdots\boxtimes\ul{\alpha_{r_0}}\)
\end{align*}
which is irreducible. The lemma follows.

Now we prove \eqref{eq:base_change_1}. Note that we have a canonical $G(F^+)$-equivariant inclusion
\[
\rI^G_P\(\tau_0\boxtimes\(\boxtimes_{i=1}^t\tau_i\(\ul{|\alpha_{r_0+\cdots+r_i}|}\circ\r{det}_{r_i}\)\)\)\subseteq
\rI^G_{P_\mnm}\(\b{1}'_2\boxtimes\ul{\alpha_2}\boxtimes\cdots\boxtimes\ul{\alpha_r}\),
\]
under which $\phi'_{\balpha}$ belongs to the former space by our choice of $\tau_0$. Let $w\in\rW'_{2r}$ be the element acting trivially on $\rV'_{2r_0}$ and switching $e_{-(r_0+\cdots+r_{i-1}+j)}$ with $e_{r_0+\cdots+r_i+1-j}$ for every $1\leq j\leq r_i$ and then every $1\leq i\leq t$. By \cite{Kon03}*{Corollary~3.2}, \eqref{eq:base_change_1} is equivalent to
\begin{align}\label{eq:base_change_11}
T_w\phi'_{\balpha}\neq 0.
\end{align}
By \cite{Cas80}*{Theorem~3.1} and the continuity, we have $T_w\phi'_{\balpha}=C(\balpha)\phi'_{w\balpha}$, where
\[
C(\balpha)=\prod_{i=r_0+1}^r\(\frac{\alpha_i-q^{-1}}{\alpha_i-1}\prod_{|\alpha_j|<|\alpha_i|}\frac{\alpha_i-q^{-2}\alpha_j}{\alpha_i-\alpha_j}
\prod_{j=1}^{i-1}\frac{\alpha_i\alpha_j-q^{-2}}{\alpha_i\alpha_j-1}\),
\]
which is nonzero. From this we obtain \eqref{eq:base_change_11}, hence \eqref{eq:base_change_1}.
\end{proof}

The following proposition exhibits an example of the local Jacquet--Langlands correspondence.

\begin{proposition}\label{pr:packet}
Define
\begin{itemize}[label={\ding{109}}]
  \item $\cS$ to be the set of isomorphism classes of irreducible admissible representations $\pi$ of $\rU_{2r}(F^+)$ such that $\pi\res_{\rK_{2r}}$ contains $\Omega_{2r}$ and that the Satake parameter of $\BC(\pi)$ contains $\{q,q^{-1}\}$ (Remark \ref{re:satake_condition});

  \item $\cS'$ to be the set of isomorphism classes of irreducible admissible representations $\pi'$ of $\rU'_{2r}(F^+)$ such that $\pi'\res_{\rK'_{2r}}$ contains the trivial representation.
\end{itemize}
Then there is a unique bijection between $\cS$ and $\cS'$ such that $\pi$ and $\pi'$ correspond if and only if $\BC(\pi)\simeq\BC(\pi')$.
\end{proposition}

\begin{proof}
We first note that both $\BC(\pi)$ and $\BC(\pi')$ are constituents of unramified principal series. We define a correspondence between $\cS$ and $\cS'$ via the condition that the two Satake parameters $\balpha(\BC(\pi))$ and $\balpha(\BC(\pi'))$ coincide. By Lemma \ref{le:base_change_0} and Lemma \ref{le:base_change_1}, the previous correspondence is a bijection, and we have $\BC(\pi)\simeq\BC(\pi')$ if $\pi$ and $\pi'$ correspond. The proposition is proved.
\end{proof}

\begin{remark}
In fact, for $\pi\in\cS$ and $\pi'\in\cS'$ in Proposition \ref{pr:packet} that correspond to each other, they should also correspond under the local theta correspondence with respect to the trivial splitting character. When $q$ is odd, this has been verified in \cite{Liu4}.
\end{remark}

\subsection{Results from the endoscopic classification}
\label{ss:endoscopic}

Now $F/F^+$ will stand for a totally imaginary quadratic extension of a totally real number field as in the main text. We state the following proposition, which summarises all we need from the endoscopic classification for unitary groups in this article. In particular, we will use the notion of local base change for unitary groups defined over $F^+_v$ for every place $v$ of $F^+$, denoted by $\BC$ as well, for which we have discussed some special cases when $v$ is inert in $F$ in \S\ref{ss:local_bc}.

\begin{proposition}\label{pr:arthur}
Take a relevant representation (Definition \ref{de:relevant}) $\Pi$ of $\GL_N(\dA_F)$. Let $\rV$ be a standard definite or indefinite hermitian space over $F$ of rank $N$ and $\pi=\otimes_v\pi_v$ an irreducible admissible representation of $\rU(\rV)(\dA_{F^+})$. We have
\begin{enumerate}
  \item If $\BC(\pi_v)\simeq\Pi_v$ for every place $v$ of $F^+$, then the discrete automorphic multiplicity of $\pi$ is $1$.

  \item If $\pi$ is automorphic and $\Pi$ is its automorphic base change (Definition \ref{de:bc_global}), then $\BC(\pi_v)\simeq\Pi_v$ holds for every place $v$ of $F^+$. In particular, the discrete automorphic multiplicity of $\pi$ is $1$ by (1).

  \item If $v$ is archimedean but not $\ul\tau_\infty$, then $\BC(\pi_v)\simeq\Pi_v$ if and only if $\pi_v$ is the trivial representation.

  \item If $v=\ul\tau_\infty$, then $\BC(\pi_v)\simeq\Pi_v$ if and only if $\pi_v$ is the trivial representation (resp.\ is one of the $N$ discrete series representations with the Harish-Chandra parameter $\{\frac{1-N}{2},\frac{3-N}{2},\dots,\frac{N-3}{2},\frac{N-1}{2}\}$) when $\rV$ is definite (resp.\ indefinite).
\end{enumerate}
\end{proposition}

\begin{proof}
Parts (1) and (2) are consequences of \cite{KMSW}*{Theorem~1.7.1} for generic packets. Parts (3) and (4) follow from (1), (2), and the definition of relevant representations.
\end{proof}

The above proposition has the following two immediate corollaries as two examples of the global Jacquet--Langlands correspondence.

\begin{corollary}\label{pr:jacquet_langlands_0}
Take a prime $\fp$ of $F^+$ inert in $F$. Let $\rV$ and $\rV'$ be a standard definite and a standard indefinite hermitian space over $F$, respectively, of even rank $N=2r$, satisfying $\rV_v\simeq\rV'_v$ (for which we fix) for every place $v$ of $F^+$ other than $\ul\tau_\infty$ and $\fp$. Let $\pi$ be an automorphic representation of $\rU(\rV)(\dA_{F^+})$ such that $\pi_\infty$ is trivial, that $\BC(\pi)$ (Definition \ref{de:bc_global}, which exists by Proposition \ref{pr:bc_global}) is cuspidal, and that $\pi_\fp$ belongs to the set $\cS$ in Proposition \ref{pr:packet} (in particular, $\rV\otimes_{F^+}F^+_\fp$ admits a self-dual lattice). Consider the representation $\pi'\coloneqq\pi'_{\ul\tau_\infty}\otimes\pi'_\fp\otimes\pi^{\ul\tau_\infty,\fp}$ of $\rU(\rV')(\dA_{F^+})$ where
\begin{itemize}[label={\ding{109}}]
  \item $\pi'_{\ul\tau_\infty}$ is a discrete series representation of $\rU(\rV')(F^+_{\ul\tau_\infty})$ with the Harish-Chandra parameter $\{\tfrac{1}{2}-r,\tfrac{3}{2}-r,\dots,r-\tfrac{3}{2},r-\tfrac{1}{2}\}$; and

  \item $\pi'_\fp\in\cS'$ is the representation of $\rU(\rV')(F^+_\fp)$ corresponding to $\pi_\fp$ as in Proposition \ref{pr:packet}.
\end{itemize}
Then the discrete automorphic multiplicity of $\pi'$ is $1$.
\end{corollary}

\begin{proof}
Put $\Pi\coloneqq\BC(\pi)$. By Proposition \ref{pr:arthur} and Proposition \ref{pr:packet}, we have $\BC(\pi'_v)\simeq\Pi_v$ for every place $v$ of $F^+$. The corollary follows by Proposition \ref{pr:arthur}(1).
\end{proof}

\begin{corollary}\label{pr:jacquet_langlands_1}
Take a prime $\fp$ of $F^+$ inert in $F$. Let $\rV$ and $\rV'$ be a standard definite and a standard indefinite hermitian space over $F$, respectively, of odd rank $N=2r+1$, satisfying $\rV_v\simeq\rV'_v$ (for which we fix) for every place $v$ of $F^+$ other than $\ul\tau_\infty$ and $\fp$. Let $\pi'$ be an automorphic representation of $\rU(\rV')(\dA_{F^+})$ such that $\BC(\pi')$ exists and is cuspidal, that $\pi'_{\ul\tau_\infty}$ is a discrete series representation of $\rU(\rV')(F^+_{\ul\tau_\infty})$ (Definition \ref{de:bc_global}) with the Harish-Chandra parameter $\{-r,1-r,\dots,r-1,r\}$, that $\pi'_{\ul\tau}$ is trivial for every archimedean place $\ul\tau\neq\ul\tau_\infty$, and that $\pi'_\fp$ is unramified. Consider the representation $\pi\coloneqq\pi_{\ul\tau_\infty}\otimes\pi_\fp\otimes(\pi')^{\ul\tau_\infty,\fp}$ of $\rU(\rV)(\dA_{F^+})$ where
\begin{itemize}[label={\ding{109}}]
  \item $\pi_{\ul\tau_\infty}$ is trivial; and

  \item $\pi_\fp$ is unramified satisfying $\BC(\pi_\fp)\simeq\BC(\pi'_\fp)$.
\end{itemize}
Then the discrete automorphic multiplicity of $\pi$ is $1$.
\end{corollary}

\begin{proof}
Put $\Pi'\coloneqq\BC(\pi')$. By Proposition \ref{pr:arthur} and Proposition \ref{pr:packet}, we have $\BC(\pi_v)\simeq\Pi'_v$ for every place $v$ of $F^+$. The corollary follows by Proposition \ref{pr:arthur}(1).
\end{proof}

\section{Some trace formulae argument}
\label{ss:c}

This appendix has two goals. In \S\ref{ss:vanishing}, we remove some conditions in a theorem of Caraiani and Scholze \cite{CS17}. In \S\ref{ss:dimension}, we prove a formula computing the dimension of old forms in an $L$-packet for unitary groups. These two subsections are independent on a logical level; we collect them together in one appendix mainly because the argument we use are similar, namely, trace formulae.

We keep the setup in Section \ref{ss:3}.

\subsection{Vanishing of cohomology off middle degree}
\label{ss:vanishing}

\begin{definition}\label{de:generic}
Let $N\geq 1$ be an integer, and $\Sigma^+$ a finite set of nonarchimedean places of $F^+$ containing $\Sigma^+_\bad$. Consider a homomorphism $\phi\colon\dT_N^{\Sigma^+}\to\kappa$ with $\kappa$ a field. We say that $\phi$ is \emph{cohomologically generic} if
\[
\rH^i_\et(\Sh(\rV,\rK)_{\ol{F}},\kappa)_{\dT_N^{\Sigma^{+\prime}}\cap\Ker\phi}=0
\]
holds for
\begin{itemize}[label={\ding{109}}]
  \item every finite set $\Sigma^{+\prime}$ of nonarchimedean places of $F^+$ containing $\Sigma^+$,

  \item every integer $i\neq N-1$, and

  \item every standard indefinite hermitian space $\rV$ over $F$ of dimension $N$ and every object $\rK\in\fK(\rV)$ of the form $\rK_{\Sigma^{+\prime}}\times\prod_{v\not\in\Sigma^+_\infty\cup\Sigma^{+\prime}}\rU(\Lambda)(O_{F^+_v})$ for a self-dual $\prod_{v\not\in\Sigma^+_\infty\cup\Sigma^{+\prime}}O_{F_v}$-lattice $\Lambda$ in $\rV\otimes_F\dA_F^{\Sigma^+_\infty\cup\Sigma^{+\prime}}$.
\end{itemize}
\end{definition}

The following definition is essentially \cite{CS17}*{Definition~1.9}.

\begin{definition}\label{de:decomposed_generic}
Let $\phi\colon\dT_N^{\Sigma^+}\to\kappa$ be a homomorphism with $\kappa$ a field. For a place $w$ of $F^+$ not in $\Sigma^+$ that splits in $F$, we say that $\phi$ is \emph{decomposed generic at $w$} if $\phi(H_w)\in\kappa[T]$ has distinct (nonzero) roots in which there is no pair with ratio equal to $\|w\|$.\footnote{In fact, as pointed out in \cite{CS19}*{Remark~1.4}, there is no need to assume that the roots are distinct.} Here, $H_w\in\dT_{N,w}[T]$ is the Hecke polynomial.
\end{definition}

\begin{proposition}\label{th:generic}
Let $N\geq 1$ be an integer, and $\Sigma^+$ a finite set of nonarchimedean places of $F^+$ containing $\Sigma^+_\bad$. Let $\rV$ be a standard indefinite hermitian space over $F$ of dimension $N$ such that $\rV_v$ is split for $v\not\in\Sigma^+_\infty\cup\Sigma^+$. Let $\phi\colon\dT_N^{\Sigma^+}\to\ol\dF_\ell$ be a homomorphism. Suppose that $F^+\neq\dQ$. Suppose that there exists a place $w$ of $F^+$ not in $\Sigma^+\cup\Sigma^+_\ell$ that splits in $F$, such that $\phi$ is decomposed generic at $w$. Then we have
\[
\rH^i_\et(\Sh(\rV,\rK)_{\ol{F}},\ol\dF_\ell)_{\Ker\phi}=0
\]
for every integer $i\neq N-1$, and every object $\rK\in\fK(\rV)$ of the form $\rK_{\Sigma^+}\times\prod_{v\not\in\Sigma^+_\infty\cup\Sigma^+}\rU(\Lambda)(O_{F^+_v})$ for a self-dual $\prod_{v\not\in\Sigma^+_\infty\cup\Sigma^+}O_{F_v}$-lattice $\Lambda$ in $\rV\otimes_F\dA_F^{\Sigma^+_\infty\cup\Sigma^+}$.
\end{proposition}

\begin{proof}
When $F$ contains an imaginary quadratic field and every place in $\Sigma^+$ splits in $F$ (which implies $F^+\neq\dQ$), the proposition can be deduced from the analogous statement for the unitary similitude group, namely Case 2 of \cite{CS17}*{Theorem~6.3.1(2)}. We now explain how to remove these restrictions.

In the statement of the proposition, let $w_0$ be the underlying rational prime of $w$. We fix an isomorphism $\dC\simeq\ol\dQ_{w_0}$ that induces the place $w$ of $F^+$. Put $\rG\coloneqq\Res_{F^+/\dQ}\rU(\rV)$. We have the Deligne homomorphism $\rh\colon\Res_{\dC/\dR}\bG_m\to\rG\otimes_\dQ\dR$ as in \S\ref{ss:unitary_shimura}. Put $\rK_{w_0,0}\coloneqq\prod_{v\mid w_0}\rU(\Lambda)(O_{F^+_v})$, which is a hyperspecial maximal subgroup of $\rG(\dQ_{w_0})$. We fix a character $\varpi\colon F^\times\backslash\dA_F^\times\to\dC^\times$ that is unramified outside $\Sigma^+$ such that $\varpi\res_{\dA_{F^+}^\times}$ is the quadratic character $\eta_{F/F^+}$ associated to $F/F^+$. Put $\Sigma\coloneqq\{p\res\Sigma^+_p\cap\Sigma^+\neq\emptyset\}$.

We define a subtorus $\rT\subseteq\Res_{F/\dQ}\bG_m$ such that for every $\dQ$-ring $R$,
\[
\rT(R)=\{a\in F\otimes_\dQ R\res \Nm_{F/{F^+}}a\in R^\times\}.
\]
We fix a CM type $\Phi$ containing $\tau_\infty$ satisfying that all elements in $\Phi$ inducing the place $w$ of $F^+$ induce the same place of $F$, and a sufficiently small open compact subgroup $\rK_\rT\subseteq\rT(\dA^\infty)$ such that $(\rK_\rT)_p$ is maximal for every $p\not\in\Sigma$. Then $\Phi$ induces a Deligne homomorphism $\rh_\Phi\colon\Res_{\dC/\dR}\bG_m\to\rT\otimes_\dQ\dR$. We also put $\fT\coloneqq\rT(\dA^{\infty,w_0})/\rT(\dZ_{(w_0)})\rK_\rT^{w_0}$ similar to Definition \ref{de:groupoid}.

Put $\tilde\rG\coloneqq\rG\times\rT$ and $\tilde\rh\coloneqq\rh\times\rh_\Phi$. Then we have the Shimura datum $(\tilde\rG,\tilde\rh)$, which is of Hodge type. Its reflex field is the composition $F.F_\Phi\subseteq\dC$. Therefore, for every sufficiently small open compact subgroup $\rK\subseteq\rG(\dA^\infty)$, we have the Shimura variety $\Sh(\tilde\rG,\tilde\rh)_{\rK\times\rK_\rT}$, which is smooth projective (as $F^+\neq\dQ$) over $F.F_\Phi$ of dimension $N-1$. When $\rK$ is of the form $\rK^{w_0}\rK_{w_0,0}$, it has a canonical smooth projective model $\sS(\tilde\rG,\tilde\rh)_{\rK^{w_0}}$ over $W(\ol\dF_{w_0})$ which admits a moduli interpretation similar to the one introduced in \S\ref{ss:qs_moduli_scheme}. Note that $F.F_\Phi$ is contained in $W(\ol\dF_{w_0})_\dQ$ under the isomorphism $\dC\simeq\ol\dQ_{w_0}$.

The discussion in \cite{CS17}, except in \S5, is valid for all proper Shimura varieties of Hodge type including the above one. Thus, we need to modify the argument in \cite{CS17}*{\S5} for our case.

Let $\mu$ and $\tilde\mu$ be the Hodge cocharacters corresponding to $\rh$ and $\tilde\rh$, respectively. We have the natural projection map $B(\tilde\rG,\tilde\mu)\to B(\rG,\mu)$ of Kottwitz sets, which is a bijection. For every $b\in B(\rG,\mu)$, we have the corresponding Kottwitz groups $\tilde{J}_b$ and $J_b$, with a canonical isomorphism $\tilde{J}_b\simeq J_b\times\rT$. For every (sufficiently small) open compact subgroup $\rK^{w_0}\subseteq\rG(\dA^{\infty,w_0})$ and positive integer $m$, we have the Igusa variety $\sI^b_{\r{Mant},\rK^{w_0},m}$ for the integral model $\sS(\tilde\rG,\tilde\rh)_{\rK^{w_0}}$, which is a $\fT$-scheme over $\ol\dF_{w_0}$. Define
\[
[\rH_{\fT,c}(\sI^b_{\r{Mant}},\ol\dQ_\ell)]\coloneqq\bigoplus_i(-1)^i\varinjlim_{\rK^{w_0},m}
\rH^i_{\fT,c}(\sI^b_{\r{Mant},\rK^{w_0},m},\ol\dQ_\ell),
\]
which is a virtual representation of $\rG(\dA^{\infty,w_0})\times J_b(\dQ_{w_0})$. The crucial point is that our $\rG$ is the honest unitary group, rather than the unitary similitude group. Then \cite{CS17}*{Theorem~5.2.3} is modified as
\[
\tr\(\phi\res\rH_{\fT,c}(\sI^b_{\r{Mant}},\ol\dQ_\ell)\)=\sum_{(\rH,s,\eta)}\iota(\rG,\rH)\r{ST}^\rH_e(\phi^\rH)
\]
where the sum is taken over equivalent classes of elliptic endoscopic triples $(\rH,s,\eta)$ of $\rG$; and we use the character $\varpi$ for the Langlands--Shelstad transfer. This formula can be proved in the same way as for \cite{Shi10}*{Theorem~7.2} since our Shimura variety has a similar moduli interpretation as seen in \S\ref{ss:qs_moduli_scheme}, although the Shimura datum $(\tilde\rG,\tilde\rh)$ is not of PEL type in the sense of Kottwitz. We can fix the representatives of the triples $(\rH,s,\eta)$ as in \cite{CS17}*{Page~734} but without the similitude factor. In particular, \cite{CS17}*{Corollary~5.2.5} is modified as
\[
\tr\(\phi\res\rH_{\fT,c}(\sI^b_{\r{Mant}},\ol\dQ_\ell)\)=\sum_{\rG_{\vec{n}}}\iota(\rG,\rG_{\vec{n}})\r{ST}^{\rG_{\vec{n}}}_e(\phi^{\vec{n}}).
\]

The next statement \cite{CS17}*{Proposition~5.3.1} or rather \cite{Shi11}*{Corollary~4.7}, namely,
\[
I^{\dG_{\vec{n}}\theta}_{\r{geom}}(f^{\vec{n}}\theta)=\tau(\rG_{\vec{n}})^{-1}\r{ST}^{\rG_{\vec{n}}}_e(\phi^{\vec{n}})
\]
holds as long as $f^{\vec{n}}$ and $\phi^{\vec{n}}$ are associated in the sense of \cite{Lab99}*{3.2}. Here, $\dG_{\vec{n}}$ is the group $\Res_{F/\dQ}\GL_{\vec{n}}\rtimes\{1,\theta\}$. Note that, for rational primes in $\Sigma$, we do not have explicit local base change transfer. However, we will see shortly that there are enough associated pairs at these primes to make the remaining argument work, following an idea in \cite{Shi}.

For the test function $\phi\in C^\infty_c(\rG(\dA^{\infty,w_0})\times J_b(\dQ_{w_0}))$ in \cite{CS17}*{Theorem~5.3.2}, if we assume $\phi=\phi_\Sigma\otimes\phi^\Sigma$ in which $\phi_\Sigma$ is the characteristic function of some open compact subgroup $\rK_\Sigma\subseteq\rG(\dQ_\Sigma)$, then for every $\rG_{\vec{n}}$, $\phi^{\vec{n}}$ is associated to some function $f^{\vec{n}}$ in the sense above. This is shown in the claim in the proof of \cite{Shi}*{Proposition~1.4}. In particular, for such $\phi$, we have
\[
\tr\(\phi\res\rH_{\fT,c}(\sI^b_{\r{Mant}},\ol\dQ_\ell)\)=\sum_{\rG_{\vec{n}}}\iota(\rG,\rG_{\vec{n}})
I^{\dG_{\vec{n}}\theta}_{\r{spec}}(f^{\vec{n}}\theta)
\]
in view of the above identities and \cite{CS17}*{(5.3.2)}. The remaining argument toward \cite{CS17}*{Theorem~5.5.7} is same as it is on the $\GL$-side, for which it suffices to use the above test functions $\phi$. In fact, our case is slightly easier as we do not have the similitude factor.

The argument towards Proposition \ref{th:generic} or \cite{CS17}*{Theorem~6.3.1(2)} only uses \cite{CS17}*{Theorem~5.5.7}. Therefore, the proposition holds.
\end{proof}

\begin{corollary}\label{co:generic}
Let the situation be as in \S\ref{ss:preliminaries}. Suppose that $F^+\neq\dQ$. Then for all but finitely many primes $\lambda$ of $E$, the composite homomorphism
\begin{align}\label{eq:generic}
\dT^{\Sigma^+}_N\xrightarrow{\phi_\Pi}O_E\to O_E/\lambda
\end{align}
is cohomologically generic (Definition \ref{de:generic}).
\end{corollary}

\begin{proof}
As pointed out in the proof of \cite{CH13}*{Proposition~3.2.5}, we can choose a nonarchimedean place $w$ of $F$ such that $\Pi_w$ is unramified whose Satake parameter contains distinct elements $\alpha_1,\dots,\alpha_N$, which are nonzero algebraic numbers. Since $\Pi_w$ is generic, we have $\alpha_i/\alpha_j\not\in\{1,\|w\|\}$ for $i\neq j$. Thus, for every sufficiently large rational prime $\ell$, we have $\alpha_i/\alpha_j\not\in\{1,\|w\|\}$ for $i\neq j$ even in $\ol\dF_\ell$. Let $\lambda$ be a prime of $E$ above such a rational prime $\ell$. Applying the Chebotarev density theorem to any residual Galois representation $\bar\rho_{\Pi,\lambda}$ of $\rho_{\Pi,\lambda}$, we conclude that there are infinitely many nonarchimedean places $w$ of $F^+$ not in $\Sigma^+\cup\Sigma^+_\ell$ that splits in $F$, such that \eqref{eq:generic} is decomposed generic at $w$ (Definition \ref{de:decomposed_generic}). Thus, \eqref{eq:generic} is cohomologically generic by Proposition \ref{th:generic}. The corollary follows.
\end{proof}

\subsection{Dimension of old forms}
\label{ss:dimension}

Let $N=2r$ be an even positive integer. We consider

\begin{itemize}[label={\ding{109}}]
  \item a relevant representation $\Pi$ of $\GL_N(\dA_F)$,

  \item two disjoint finite sets $\Sigma^+_\mnm$ and $\Sigma^+_\lr$ of nonarchimedean places of $F^+$ such that $\Sigma^+_\mnm$ contains $\Sigma^+_\bad$; $\Sigma^+_\mnm\cup\Sigma^+_\lr$ contains $\Sigma^+_\Pi$ (Notation \ref{no:satake}); and every place in $\Sigma^+_\lr$ is inert in $F$,

  \item a finite set $\Sigma^+$ of nonarchimedean places of $F^+$ containing $\Sigma^+_\mnm\cup\Sigma^+_\lr$,

  \item a standard definite or indefinite hermitian space $\rV$ over $F$ of rank $N$ such that $\rV_v$ is not split for $v\in\Sigma^+_\lr$,

  \item a self-dual $\prod_{v\not\in\Sigma^+_\infty\cup\Sigma^+_\mnm\cup\Sigma^+_\lr}O_{F_v}$-lattice $\Lambda$ in $\rV\otimes_F\dA_F^{\Sigma^+_\infty\cup\Sigma^+_\mnm\cup\Sigma^+_\lr}$,

  \item an object $\rK\in\fK(\rV)$ of the form
      \[
      \rK=\prod_{v\in\Sigma^+_\mnm\cup\Sigma^+_\lr}\rK_v\times
      \prod_{v\not\in\Sigma^+_\infty\cup\Sigma^+_\mnm\cup\Sigma^+_\lr}\rU(\Lambda)(O_{F^+_v}),
      \]
      satisfying that $\rK_v$ is special maximal for $v\in\Sigma^+_\lr$.
\end{itemize}
We have the homomorphism
\[
\phi_\Pi\colon\dT^{\Sigma^+}_N\to\ol\dQ_\ell
\]
given by $\Pi$. Fix an isomorphism $\iota_\ell\colon\dC\xrightarrow{\sim}\ol\dQ_\ell$.

\begin{definition}\label{de:transferable}
Let $v$ be a nonarchimedean place of $F^+$. We say that an open compact subgroup $\rK_v$ of $\rU(\rV)(F^+_v)$ is \emph{transferable} if the following two conditions are satisfied.
\begin{enumerate}
  \item For every endoscopic group $\rH$ of $\rU(\rV_v)$, there exist an endoscopic transfer $f_{\rK_v}^\rH$ of $\CF_{\rK_v}$ to $\rH$ and a compactly supported smooth function $\phi_{\rK_v}^\rH$ on $\rH(F_v)$ such that $f_{\rK_v}^\rH$ and $\phi_{\rK_v}^\rH$ are associated in the sense of \cite{Lab99}*{\S3.2}.

  \item When $\rH$ is the quasi-split unitary group of rank $N$, we can take $\phi_{\rK_v}^\rH$ to be supported on a maximal open compact subgroup of $\rH(F_v)$ (which is isomorphic to $\GL_N(F_v)$).\footnote{In fact, this restriction is not necessary for Proposition \ref{pr:old} below; it is only used in the application of this proposition, namely, Proposition \ref{pr:old_same}.}
\end{enumerate}
We call the function $\phi_{\rK_v}^\rH$ in (2) an \emph{inertial transfer} of $\rK_v$ if $\rK_v$ is transferable, and will drop the superscript $\rH$ in practice.
\end{definition}

\begin{lem}\label{le:transferable}
Let $v$ be a nonarchimedean place of $F^+$.
\begin{enumerate}
  \item If $v$ splits in $F$, then every open compact subgroup $\rK_v$ is transferable.

  \item If $v$ is not in $\Sigma^+_\infty\cup\Sigma^+_\mnm\cup\Sigma^+_\lr$, then the characteristic function of the hyperspecial maximal subgroup $\rU(\Lambda)(O_{F^+_v})$ is transferable and admits $\CF_{\GL_N(O_{F_v})}$ as an inertial transfer.

  \item If $v$ is in $\Sigma^+_\mnm\cup\Sigma^+_\lr$, then every sufficiently small open compact subgroup $\rK_v$ is transferable.
\end{enumerate}
\end{lem}

\begin{proof}
Part (1) is trivial. Part (2) is the combination of the endoscopic fundamental lemma \cite{LN08} and the base change fundamental lemma \cite{Lab99}.

For (3), for sufficiently small $\rK_v$, condition (1) in Definition \ref{de:transferable} is proved in \cite{Mor10}*{Lemma~8.4.1(1)}; and condition (2) can be achieved by \cite{Lab99}*{Proposition~3.1.7(2)} (see the proof of \cite{Lab99}*{Proposition~3.3.2}).
\end{proof}

\begin{proposition}\label{pr:old}
Suppose that $\rK_v$ is transferable for $v\in\Sigma^+_\mnm$. For every $v\in\Sigma^+_\lr$, let $c_v$ be equal to $1$ (resp.\ $0$) if one can (resp.\ cannot) find complex numbers $\alpha_2,\dots,\alpha_r$ of norm one such that $\Pi_v$ is isomorphic to the induction
\[
\rI^{\GL_{2r}}
\(\ul{\alpha_r^{-1}}\boxtimes\cdots\boxtimes\ul{\alpha_2^{-1}}\boxtimes\St_2\boxtimes\ul{\alpha_2}\boxtimes\cdots\boxtimes\ul{\alpha_r}\)
\]
(see \S\ref{ss:local_bc} for the notation of induced representations). Then we have the identities
\begin{align*}
\dim\ol\dQ_\ell[\Sh(\rV,\rK)][\iota_\ell\phi_\Pi]
&=\left|\prod_{v\in\Sigma^+_\mnm}\tr(\Pi_v(\phi_{\rK_v})\circ A_{\Pi_v})\prod_{v\in\Sigma^+_\lr}c_v\right|, \\
\dim\rH^{N-1}_\et(\Sh(\rV,\rK)_{\ol{F}},\ol\dQ_\ell)[\iota_\ell\phi_\Pi]
&=N\left|\prod_{v\in\Sigma^+_\mnm}\tr(\Pi_v(\phi_{\rK_v})\circ A_{\Pi_v})\prod_{v\in\Sigma^+_\lr}c_v\right|,
\end{align*}
when $\rV$ is definite and indefinite, respectively, for any inertial transfer $\phi_{\rK_v}$ for $\rK_v$ and any normalized intertwining operator $A_{\Pi_v}$ for $\Pi_v$ \cite{Shi11}*{\S4.1}, for $v\in\Sigma^+_\mnm$.
\end{proposition}

\begin{proof}
We only prove the case where $\rV$ is indefinite, and leave the case where $\rV$ is definite (which is slightly easier) to the readers.

By Proposition \ref{pr:galois}(1), we know that $\Pi$ is tempered everywhere. Moreover, every discrete automorphic representation of $\rU(\rV)(\dA_{F^+})$ whose automorphic base change is isomorphic to $\Pi$ has to be cuspidal as well. Thus, we have $\rH^i_\et(\Sh(\rV,\rK)_{\ol{F}},\ol\dQ_\ell)[\iota_\ell\phi_\Pi]=0$ for $i\neq N-1$.

If there exists $v\in\Sigma^+_\lr$ such that $c_v=0$, then by Lemma \ref{le:base_change_1} and the above fact that $\Pi_v$ is tempered, we have $\rH^{N-1}_\et(\Sh(\rV,\rK)_{\ol{F}},\ol\dQ_\ell)[\iota_\ell\phi_\Pi]=0$. Thus, the proposition follows. In what follows, we assume $c_v=1$ for every $v\in\Sigma^+_\lr$.

By Proposition \ref{pr:arthur} and Lemma \ref{le:base_change_1}, we have
\[
\dim\rH^{N-1}_\et(\Sh(\rV,\rK)_{\ol{F}},\ol\dQ_\ell)[\iota_\ell\phi_\Pi]=N\prod_{v\in\Sigma^+_\mnm}
\sum_{\BC(\pi_v)\simeq\Pi_v}\dim(\pi_v)^{\rK_v},
\]
where the sum is taken over isomorphism classes of irreducible admissible representations $\pi_v$ of $\rU(\rV)(F^+_v)$ such that $\BC(\pi_v)\simeq\Pi_v$ (for $v\in\Sigma^+_\mnm$). Thus, our goal is to show
\begin{align}\label{eq:old}
\prod_{v\in\Sigma^+_\mnm}
\sum_{\BC(\pi_v)\simeq\Pi_v}\dim(\pi_v)^{\rK_v}=
\left|\prod_{v\in\Sigma^+_\mnm}\tr(\Pi_v(\phi_{\rK_v})\circ A_{\Pi_v})\right|.
\end{align}

We choose a quadratic totally real extension $\breve{F}^+/F^+$ in $\dC$ satisfying
\begin{itemize}[label={\ding{109}}]
  \item every prime in $\Sigma^+_\mnm$ splits in $\breve{F}^+$;

  \item every prime in $\Sigma^+_\lr$ is inert in $\breve{F}^+$;

  \item the quadratic base change of $\Pi$ to $\breve{F}\coloneqq F.\breve{F}^+$, denoted by $\breve\Pi$, remains cuspidal (hence relevant).
\end{itemize}
By the same proof of \cite{Shi11}*{Proposition~7.4}, we know that such $\breve{F}^+$ exists. Let $\breve\rV$ be the standard definite hermitian space over $\breve{F}$ of rank $N$ that is split at all primes not above $\Sigma^+_\mnm$ and such that $\breve\rV_{\breve{v}}\simeq\rV_{v}$ for every $v\in\Sigma^+_\mnm$ and every prime $\breve{v}$ of $\breve{F}^+$ above $v$, which exists as $[\breve{F}^+:F^+]=2$. Let $\breve\Sigma^+_\mnm$ be the set of primes of $\breve{F}^+$ above $\Sigma^+_\mnm$. Take a finite set $\breve\Sigma^+$ of primes of $\breve{F}^+$ satisfying
\begin{itemize}[label={\ding{109}}]
  \item $\breve\Sigma^+$ contains $\breve\Sigma^+_\mnm$;

  \item $\breve\Pi_{\breve{v}}$ is unramified for every prime of $\breve{F}^+$ not in $\breve\Sigma^+$;

  \item every prime in $\breve\Sigma^+\setminus\breve\Sigma^+_\mnm$ splits in $\breve{F}$.
\end{itemize}
By our choice of $\breve{F}^+$, such $\breve\Sigma^+$ exists. Take an object $\breve\rK\in\fK(\breve\rV)$ of the form $\breve\rK=\prod\breve\rK_{\breve{v}}$ satisfying
\begin{itemize}[label={\ding{109}}]
  \item $\breve\rK_{\breve{v}}$ is hyperspecial maximal if $\breve{v}\not\in\breve\Sigma^+$;

  \item $\breve\rK_{\breve{v}}=\rK_v$ under a chosen isomorphism $\breve\rV_{\breve{v}}\simeq\rV_{v}$ if $\breve{v}$ is above a prime $v\in\Sigma^+_\mnm$;

  \item $\breve\Pi_{\breve{v}}$ has nonzero $\breve\rK_{\breve{v}}\times\breve\rK_{\breve{v}}$ invariants if $\breve{v}\in\breve\Sigma^+\setminus\breve\Sigma^+_\mnm$.
\end{itemize}

Then we have
\begin{align}\label{eq:old2}
\dim\ol\dQ_\ell[\Sh(\breve\rV,\breve\rK)][\iota_\ell\phi_{\breve\Pi}]=
\prod_{\breve{v}\in\breve\Sigma^+}\sum_{\BC(\breve\pi_{\breve{v}})\simeq\breve\Pi_{\breve{v}}}
\dim(\breve\pi_{\breve{v}})^{\breve\rK_{\breve{v}}}.
\end{align}
On the other hand, by \cite{Shi}*{(1.8) \& (1.9)}, we have
\begin{align}\label{eq:old3}
\dim\ol\dQ_\ell[\Sh(\breve\rV,\breve\rK)][\iota_\ell\phi_{\breve\Pi}]=
\left|\prod_{\breve{v}\in\breve\Sigma^+}\tr(\breve\Pi_{\breve{v}}(\phi_{\breve\rK_{\breve{v}}})\circ A_{\breve\Pi_{\breve{v}}})\right|.
\end{align}
Here, for $\breve{v}\in\breve\Sigma^+\setminus\breve\Sigma^+_\mnm$, we take $\phi_{\breve\rK_{\breve{v}}}$ to be $\CF_{\breve\rK_{\breve{v}}}\otimes\CF_{\breve\rK_{\breve{v}}}$; and it is easy to see that
\begin{align}\label{eq:old4}
\left|\tr(\breve\Pi_{\breve{v}}(\phi_{\breve\rK_{\breve{v}}})\circ A_{\breve\Pi_{\breve{v}}})\right|=
\sum_{\BC(\breve\pi_{\breve{v}})\simeq\breve\Pi_{\breve{v}}}\dim(\breve\pi_{\breve{v}})^{\breve\rK_{\breve{v}}}\geq 1
\end{align}
(in fact, the sum is taken over a singleton). Combining \eqref{eq:old2}, \eqref{eq:old3}, and \eqref{eq:old4}, we obtain
\[
\prod_{\breve{v}\in\breve\Sigma^+_\mnm}\sum_{\BC(\breve\pi_{\breve{v}})\simeq\breve\Pi_{\breve{v}}}
\dim(\breve\pi_{\breve{v}})^{\breve\rK_{\breve{v}}}=
\left|\prod_{\breve{v}\in\breve\Sigma^+_\mnm}\tr(\breve\Pi_{\breve{v}}(\phi_{\breve\rK_{\breve{v}}})\circ A_{\breve\Pi_{\breve{v}}})\right|,
\]
which is nothing but
\[
\(\prod_{v\in\Sigma^+_\mnm}\sum_{\BC(\pi_v)\simeq\Pi_v}\dim(\pi_v)^{\rK_v}\)^2=
\left|\prod_{v\in\Sigma^+_\mnm}\tr(\Pi_v(\phi_{\rK_v})\circ A_{\Pi_v})\right|^2.
\]
Thus, \eqref{eq:old} follows. The proposition is proved.
\end{proof}

\end{document}